\documentclass{svmult}

\usepackage{amsmath}
\usepackage{amssymb}

\usepackage{graphics} 
\usepackage{color}
\usepackage{makeidx}         
\usepackage{graphicx}        
\usepackage{multicol}        
\newtheorem{dfn}{Definition}[section]
\newtheorem{rem}[dfn]{Remark}

\newtheorem{ex}[dfn]{Example}
\newtheorem{lem}[dfn]{Lemma}
\newtheorem{prop}[dfn]{Proposition}
\newtheorem{thm}[dfn]{Theorem}
\newtheorem{cor}[dfn]{Corollary}
\newcommand{\bn}{{\bold n}}
\newcommand{\ba}{{\bold a}}
\newcommand{\bb}{{\bold b}}
\newcommand{\bd}{{\bold d}}
\newcommand{\bu}{{\bold u}}
\newcommand{\bv}{{\bold v}}
\newcommand{\bw}{{\bold w}}
\newcommand{\bff}{{\bold f}}
\newcommand{\be}{{\bold e}}
\newcommand{\bg}{{\bold g}}
\newcommand{\bh}{{\bold h}}
\newcommand{\bk}{{\bold k}}
\newcommand{\bp}{{\bold p}}

\newcommand{\bA}{{\bold A}}
\newcommand{\bB}{{\bold B}}

\newcommand{\bD}{{\bold D}}

\newcommand{\bF}{{\bold F}}

\newcommand{\bH}{{\bold H}}
\newcommand{\bI}{{\bold I}}
\newcommand{\bJ}{{\bold J}}
\newcommand{\bK}{{\bold K}}

\newcommand{\bM}{{\bold M}}

\newcommand{\bP}{{\bold P}}
\newcommand{\bQ}{{\bold Q}}
\newcommand{\bS}{{\bold S}}
\newcommand{\bT}{{\bold T}}
\newcommand{\bU}{{\bold U}}
\newcommand{\bV}{{\bold V}}
\newcommand{\bW}{{\bold W}}
\newcommand{\BPi}{{\bold \Pi}}
\newcommand{\DV}{{\rm Div}\,}
\newcommand{\dv}{{\rm div}\,}

\newcommand{\BB}{{\Bbb B}}
\newcommand{\BR}{{\Bbb R}}
\newcommand{\BC}{{\Bbb C}}

\newcommand{\BN}{{\Bbb N}}

\newcommand{\BZ}{{\Bbb Z}}

\newcommand{\tzeta}{\tilde{\zeta}}

\newcommand{\CA}{{\mathcal A}}
\newcommand{\CB}{{\mathcal B}}
\newcommand{\CC}{{\mathcal C}}
\newcommand{\CD}{{\mathcal D}}
\newcommand{\CE}{{\mathcal E}}
\newcommand{\CF}{{\mathcal F}}
\newcommand{\CI}{{\mathcal I}}
\newcommand{\CJ}{{\mathcal J}}
\newcommand{\CK}{{\mathcal K}}
\newcommand{\CL}{{\mathcal L}}
\newcommand{\CM}{{\mathcal M}}
\newcommand{\CN}{{\mathcal N}}

\newcommand{\CR}{{\mathcal R}}
\newcommand{\CS}{{\mathcal S}}
\newcommand{\CT}{{\mathcal T}}
\newcommand{\CH}{{\mathcal H}}
\newcommand{\CO}{{\mathcal O}}
\newcommand{\CP}{{\mathcal P}}

\newcommand{\CU}{{\mathcal U}}
\newcommand{\CV}{{\mathcal V}}
\newcommand{\CW}{{\mathcal W}}
\newcommand{\CX}{{\mathcal X}}
\newcommand{\CY}{{\mathcal Y}}
\newcommand{\CZ}{{\mathcal Z}}

\newcommand{\fg}{{\frak g}}
\newcommand{\fp}{{\frak p}}
\newcommand{\fq}{{\frak q}}
\newcommand{\fr}{{\frak r}}

\newcommand{\pd}{\partial}

\newcommand{\Hol}{{\rm Hol}\,}
\newcommand{\HS}{{\BR^N_+}}
\newcommand{\pf}{\noindent {\bf Proof.}\quad}
\newenvironment{cases*}%
{%
\left\{
\begin{array}{@{}r@{\;}l@{\quad}l@{}}
}%
{\end{array}\right.}
\makeindex             

\begin{document}

\title*{$\CR$ Boundedness, Maximal Regularity and  
  Free Boundary Problems
for the Navier Stokes Equations}
\titlerunning{$\CR$-boundedness, Maximal Regularity and Free Boundary Problem}
\author{Yoshihiro Shibata 
\thanks{Partially supported by JSPS
Grant-in-aid for Scientific Research  (A) 17H01097 and 
Top Global University Project.
\endgraf
Adjunct faculty member in the Department of Mechanical Engineering
and Materials Science, University of Pittsburgh.}
}
\institute{Department of Mathematics and Research Institute of 
Science and Engineering, Waseda University, 
3-4-1 Ohkubo Shinjuku-ku, 169-8555 Tokyo, Japan. 
\texttt{yshibata325@gmail.com}}
%
%

\maketitle
\begin{abstract}
In this lecture note, we study  free bounary  problems for the 
Navier-Stokes equations with and without surface tension. 
The local wellposedness, the global wellposedness,
and asymptotics of solutions as time goes to infinity are studied 
in the $L_p$ in time
and $L_q$ in space framework. The tool in proving the local well-posedness
is the maximal $L_p$-$L_q$ regularity for the Stokes equations with
non-homogeneous free
boundary conditions. The approach here of proving the maximal $L_p$-$L_q$ 
regularity is based on the $\CR$ bounded solution operators of the 
generalized resolvent problem for the Stokes equations with 
non-homogeneous free
boundary conditions and the Weis operator valued Fourier
multiplier.  
 
The key issue of proving the global well-posedness for the strong solutions 
is the decay properties of  Stokes semigroup, which are derived 
by  spectral analysis of the Stokes operator in the bulk space
 and the Laplace-Beltrami operator on the boundary.  In this 
lecture note, we study the following two cases: \thetag1 a bounded domain with
surface tension and \thetag2 an exterior domain without surface tension.
In particular, in studying the exterior domain case, it is essential
to choose different exponents $p$ and $q$.  Because, in the unbounded
domain case, we can obtain only polynomial decay in suitable 
$L_q$ norms in space, and so to guarantee the integrability of 
$L_p$ norm of solutions in time, it is necessary to have 
freedom to choose an exponent with respect to time variable. 

\end{abstract}


\section{Introduction} \label{sec:1}

\subsection{Free boundary problem for the Navier-Stokes equations}
In this chapter, we study free boundary problem for 
the Navier-Stokes equations and $\CR$ bounded solution operators
and $L_p$-$L_q$ maximal regularity theorem for the Stokes equations
with free boundary conditions. 
Typical problems  are \thetag{P1} the motion of an isolated
liquid mass and \thetag{P2}
 the motion of a viscous incompressible fluid contained in an ocean of infinite extent.

The mathematical problem is to find a time dependent domain $\Omega_t$, 
$t$ being time variable, in the $N$-dimensional Euclidean
space $\BR^N$, the velocity vector field, $\bv(x, t)={}^\top(v_1(x, t),
\ldots, v_N(x, t))$, where ${}^\top M$ denotes the transposed $M$, 
and the pressure field $\fp = \fp(x, t)$ satisfying the 
Navier-Stokes equations in $\Omega_t$:
\begin{equation}\label{navier:1}\left\{\begin{aligned}
\pd_t\bv + (\bv\cdot\nabla)\bv - 
\DV(\mu\bD(\bv) - \fp\bI)  =0& 
&\quad&\text{in $\bigcup_{0 < t < T}\Omega_t\times\{t\}$}, \\
\dv\bv  = 0& &\quad&\text{in $\bigcup_{0 < t < T}\Omega_t\times\{t\}$}, \\
(\mu\bD(\bv) - \fp\bI)\bn_t  = \sigma H(\Gamma_t)\bn_t&
&\quad&\text{on $\bigcup_{0 < t < T}\Gamma_t\times\{t\}$}, \\
V_n  = \bv\cdot\bn_t& 
&\quad&\text{on $\bigcup_{0 < t < T}\Gamma_t\times\{t\}$}, \\
\bv|_{t=0}  = \bv_0 \quad\text{in $\Omega_0$}, 
\quad \Omega_t|_{t=0} = \Omega_0&,
\end{aligned}\right.\end{equation}
where $\Gamma_t$ is the boundary of $\Omega_t$ and $\bn_t$ is the unit 
outer normal to $\Gamma_t$.   As for the remaining  notation
in \eqref{navier:1}, $\pd_t= \pd/\pd t$, $\bv_0= {}^\top(v_{01},
\ldots, v_{0N})$ is a given initial velocity field, $\bD(\bv)
= \nabla \bv + {}^\top\nabla\bv$ the doubled deformation tensor
with $(i, j)^{\rm th}$ component $D_{ij}(\bv) = \pd_iv_j + \pd_jv_i$,
$\pd_i = \pd/\pd x_i$, $\bI$ the $N\times N$ identity matrix,
$H(\Gamma_t)$ the $N-1$ times mean curvature of $\Gamma_t$
given by $H(\Gamma_t)\bn_t = \Delta_{\Gamma_t}x$ ($x \in \Gamma_t$),
where $\Delta_{\Gamma_t}$ is the Laplace-Beltrami operator on 
$\Gamma_t$, $V_n$ the velocity of the evolution of free surface $\Gamma_t$
in the direction of $\bn_t$, and $\mu$ and $\sigma$
are positive constants representing 
the viscous coefficient and the 
coefficient of the surface tension, respectively. Moreover,
for any matrix field $\bK$ with $(i, j)^{\rm th}$ component 
$K_{ij}$, the quantity $\DV \bK$ is an $N$-vector with
$i^{\rm th}$ component $\sum_{j=1}^N\pd_jK_{ij}$, and for 
any vector of functions $\bw = {}^\top(w_1, \ldots, w_N)$, 
we set $\dv\bw = \sum_{j=1}^N\pd_jw_j$ and 
$(\bw\cdot\nabla)\bw$ is an $N$-vector with $i^{\rm th}$ component
$\sum_{j=1}^Nw_j\pd_jw_i$. The domain $\Omega_0=\Omega$ is given,  
$\Gamma$ denotes the boundary of $\Omega$, and $\bn$ the unit outer
normal to $\Gamma$. For simplicity, we assume that the mass 
density equals one in this chapter. 

Since $\Omega_t$ is unknown, we transform $\Omega_t$ to some fixed domain
$\Omega$, and the system of linearized equations of the nonlinear equations on 
$\Omega$ has the following forms:
\begin{equation}\label{linear:1.1}\left\{\begin{aligned}
\pd_t\bu - \DV(\mu\bD(\bu) - \fp\bI) = \bff&&\quad&\text{in $\Omega^T$},
\\
\dv\bu = g = \dv\bg&&\quad&\text{in $\Omega^T$}, \\
\pd_t\eta + <\ba \mid \nabla'_\Gamma\eta> - \bn\cdot\bu =d&
&\quad&\text{on $\Gamma^T$}, \\
(\mu\bD(\bu) - \fp\bI)\bn - \sigma((\Delta_\Gamma + b)\eta)\bn
= \bh&&\quad&\text{on $\Gamma^T$}, \\
(\bu, \eta)|_{t=0} = (\bu_0, \eta_0)&&\quad&\text{in $\Omega\times\Gamma$}.
\end{aligned}\right.\end{equation}
Here, $\eta$ is a unknown scalar function obtained by linearization of
kinetic equation $V_n = \bv\cdot\bn_t$, $\Gamma$ is the boundary of
$\Omega$, $\nabla'_\Gamma\eta$ denotes the tangential derivative of $\eta$
on $\Gamma$, $\ba$ and $b$ are given funtions defined on $\Gamma$,
and $\bff$, $g$, $\bg$, $d$, $\bh$, $\bu_0$ and $\eta_0$ are
prescribed functions. Moreover, we set $\Omega^T = \Omega\times(0, T)$ and 
$\Gamma^T = \Gamma\times(0, T)$. 

The main topics for Eq. \eqref{linear:1.1} is to prove the 
maximal $L_p$-$L_q$ regularity and the decay properties of 
solutions. To prove these properties, we consider the corresponding 
resolvent problem:
\begin{equation}\label{linear:1.2}\left\{\begin{aligned}
\lambda\bu - \DV(\mu\bD(\bu) - \fp\bI) = \bff&&\quad&\text{in $\Omega$},
\\
\dv\bu = g = \dv\bg&&\quad&\text{in $\Omega$}, \\
\lambda\eta + <\ba \mid \nabla'_\Gamma\eta> - \bn\cdot\bu =d&
&\quad&\text{on $\Gamma$}, \\
(\mu\bD(\bu) - \fp\bI)\bn - \sigma((\Delta_\Gamma + b)\eta)\bn
= \bh&&\quad&\text{on $\Gamma$}.
\end{aligned}\right.\end{equation}
The main issue of Eq. \eqref{linear:1.2} is to prove the existence of
$\CR$ bounded solution operators, 
which, combined with the Weis operator valued Fourier multiplier
theorem \cite{Weis}, yields the maximal $L_p$-$L_q$ regularity for Eq.
\eqref{navier:1}. Moreover, using some spectral analysis of solutions
to Eq. \eqref{linear:1.2}, we derive decay properties of 
solutions to Eq. \eqref{linear:1.1}. 

This chapter is organized as follows. In the rest of Sect. \ref{sec:1}
we derive the boundary condition in Eq. \eqref{navier:1} in view of 
the conservation of mass and momentum. Moreover, we derive
the consevation of angular momentum. To end Sect. \ref{sec:1},
we give a short history of mathematical studies of Eq. \eqref{navier:1},
and further notation used throughout this chapter. In Sect. \ref{sec:p}, we 
give the definition of uniformly $C^k$ domains ($k \geq 2$), 
the weak Dirichlet problem, and the Laplace-Beltrami operators. 
The weak Dirichlet problem is used to define the reduced Stokes 
equations obtained by eliminating the pressure term and the
Laplace-Beltrami operator plays an essential role to describe the 
surface tension.  In Sect. \ref{sec:2}, we explain some 
transformation of $\Omega_t$ to a fixed domain and derive the system
of nonlinear equations on this fixed domain. In Sect. \ref{sec:3}, 
we state the maximal $L_p$-$L_q$ regularity theorems for the linearized
equations \eqref{linear:1.1} and also the existence theorem of $\CR$-bounded
solution operator for the resolvent problem \eqref{linear:1.2}. Moreover,
using the $\CR$ bounded solution operator and the Weis operator valued
Fourier multiplier theorem, we prove the maximal $L_p$-$L_q$ regularity 
theorem. In Sect. \ref{sec:5.0}, 
we prove the existence of $\CR$-bounded solution
operators by dividing the studies into the following subsections:
 Subsec. \ref{subsec:5.1} is devoted to a model problem in $\BR^N$, 
Subsec. \ref{subsec:3.2} perturbed problem in $\BR^N$, 
Subsec. \ref{subsec:4} a model problem in $\BR^N_+$, and 
Subsec. \ref{sec:5} a problem in a bent half space.
And then, in the rest of subsections, putting together the results 
obtained in previous subsections and using the partition of 
unity, we prove the existence of $\CR$-bounded solution operators
and also we prove the uniqueness of solutions. Since
the pressure term is non-local, it does not fit the  usual 
localization argument using in the parameter elliptic problem,
and so according to an idea due to Grubb-Solonnikov \cite{GS},
we consider the reduced Stokes equation which is defined in
Subsec \ref{subsec:3.3}. In Sect. \ref{sec:loc6}, we prove the 
local well-posedness of Eq. \eqref{navier:1}.
Where, we use the Hanzawa transform to transform Eq. \eqref{navier:1} to
nonlinear equations on a fixed domain. Thus, one of main difficulties is 
to treat the nonlinear term of the form 
$<\bu \mid \nabla'_\Gamma\rho>$ on the boundary which arises
 from the kinetic equation $V_n=\bv\cdot\bn_t$. If we assume that 
the initial data of both  $\bu$ and $\rho$ are small, then this term
is harmless.  But, if we want to avoid the smallness condition on
 the initial
velocity field $\bu_0=\bu|_{t=0}$, we approximate
initial data $\bu_0$ by a family of functions $\{\bu_\kappa\}$ such that 
$\|\bu_\kappa-\bu_0\|_{L_q} \sim \kappa^a$ and $\|\bu_\kappa\|_{H^2_q} 
\sim \kappa^{-b}$ as $\kappa\to0$ with some positive constants
$a$ and $b$. In Sect. \ref{sec:6}, we prove the global well-posedness of
Eq. \eqref{navier:1} under the assumption that the initial domain 
$\Omega_t|_{t=0} = \Omega$ is closed to a ball and initial data are small.
We use the Hanzawa transform whose center is the barycenter of $\Omega_t$
which is crucial and the result is, roughly speaking, that 
the barycenter converges to some point $\xi_\infty$ and 
the shape of domain becomes a ball with the center at $\xi_\infty$ when
time goes infinity. In the final section, Sect. 8, we consider Eq.
\eqref{navier:1} in the case where $\sigma=0$, that is the surface
tension is not taken into account, and we prove the global
well-posedness of Eq. \eqref{navier:1} in the case 
where $\Omega= \Omega_t|_{t=0}$ is an exterior domain in 
$\BR^N$ ($N \geq 3$).  Since we consider the without surface tension case,
we can not represent the free surface by using some representing function
like Hanzawa transform, because we can not obtain enough regularity of such 
functions. Thus, we use the local Lagrange transform to transform
$\Gamma_t$ to $\Gamma_t|_{t=0} = \Gamma$, which is identity on the outside 
of some large ball.

\subsection{Derivation of boundary conditions.}\label{subsec:1.2}
We drive the boundary 
conditions which guarantee the conservation of mass and the conservation of 
momentum.  For a while,  the mass density $\rho$ is assumed to be a function of
$(x, t)$. 
Let $\phi_t : \Omega \to \BR^N$ be a smooth injective map 
with suitable regularity for 
each time $t \geq 0$ such that $\phi_t(y)|_{t=0} = y$, and let    
$$\Omega_t = \{x = \phi_t(y) \mid y \in \Omega\}.$$ 
Let $\rho$ and $\bv$ satisfy
 the Navier-Stokes equations: 
\begin{alignat}2
\pd_t\rho + \dv(\rho\bv)  &= 0  &\quad&\text{in $\bigcup_{0 < t < T}
\Omega_t\times\{t\}$, } \label{eq:mass}\\
\rho(\pd_t\bv + \bv\cdot\nabla\bv) - \DV\bS(\bv, \fp)  &= 0
&\quad&\text{in $\bigcup_{0 < t < T}
\Omega_t\times\{t\}$}.
\label{eq:momentum}
\end{alignat}
We say that  Eq. \eqref{eq:mass} is the equation of mass
and Eq. \eqref{eq:momentum} the equation of momentum. 
Let 
\begin{equation}\label{stress:1.1}
\bS(\bv, \fp) = \mu\bD(\bv) + (\nu-\mu)\dv\bv\bI - \fp \bI,
\end{equation}
which denotes the stress tensor, 
where $\mu$ and $\nu$ are positive constants describing the 
first and second viscosity coefficients, respectively. 
If we assume that $\rho$ is a positive 
constant, then  by \eqref{eq:mass} $\dv \bv=0$, which, combined with
\eqref{stress:1.1}, leads to $\bS(\bv, \fp) = \mu \bD(\bv) - \fp \bI$.
This is the situation of Eq. \eqref{navier:1}.

 Let  
$$\bS(\bv, \fp) 
=\left( \begin{matrix}
\bS_1(\bv, \fp)\\ \vdots \\ \bS_N(\bv, \fp)
\end{matrix}\right).$$
Where $\bS_i(\bv, \fp)$ are  $N$ row vectors of functions whose
$j$-th component is $\mu D_{ij}(\bv) + (\nu-\mu)\dv\bv\delta_{ij}-
\fp \delta_{ij}$, $\delta_{ij}$ being the
Kronecker delta symbols, that is, $\delta_{ii} = 1$ and 
$\delta_{ij} = 0$ ($i\not=j$).  
Then, Eq. \eqref{eq:momentum} 
is written componentwise as 
\begin{align*}
&\rho(\pd_t v_i + \sum_{j=1}^Nv_j\pd_jv_i)
-\dv\bS_i(\bv, \fp) \\
&= 
\rho(\pd_t v_i + \sum_{j=1}^Nv_j\pd_jv_i)
-\mu\Delta v_i - \nu\pd_i\dv\bv + \pd_i\fp = 0
\quad(i=1, \ldots, N),
\end{align*}
where $\displaystyle{
\Delta v_i = \sum_{j=1}^N\frac{\pd^2v_i}{\pd x_j^2}}$.

 To drive boundary conditions
that preserve the conservation of mass and the conservation
of momentum, we use the following theorem.
\begin{thm}[Reynolds transport theorem]\label{thm:rey}
Let $y = \phi_t^{-1}(x)$ be the inverse map of 
$x =\phi_t(y)$.  Let $J=J(y, t)$ be the Jacobian of the transformation:
$x = \phi_t(y)$ and $\bw(x, t)= (\pd_t\phi_t)(\phi_t^{-1}(x))$.   Then, 
$$\frac{\pd }{\pd t}J(y, t) = (\dv_x\bw(x,t))J(y,t).$$
\end{thm}
\begin{remark} Reynolds transport theorem will be proved in \eqref{eq:ray}
of Subsec. \ref{sec:para} below.
\end{remark}
Let $f = f(x, t)$ be a function defined on $\overline{\Omega_t}$,
and then
\begin{equation}\label{rey:1}
\frac{d}{dt}\int_{\Omega_t} f(x, t)\,dx 
= \int_{\Omega_t} \pd_t f(x, t) + \dv(f(x, t)\bw(x, t))\,dx.
\end{equation}
In fact,  by Reynolds transport theorem, 
\begin{align*}
&\frac{d}{dt}\int_{\Omega_t}f(x, t)\,dx
= \frac{d}{dt}\int_\Omega f(\phi_t(y), t)J(y, t)\,dy \\
&= \int_\Omega \{(\pd_tf + (\nabla f)\cdot \pd_t\phi_t)J
+f \pd_tJ\}\,dy 
= \int_{\Omega_t}( \pd_tf +\dv(f\bw))\,dx.
\end{align*}

We first consider the conservation of mass: 
\begin{equation}\label{mass:1}
\frac{d}{dt}\int_{\Omega_t} \rho\,dx = 0.
\end{equation}
By \eqref{rey:1} and the equation of mass, Eq. \eqref{eq:mass}, 
\begin{align*}
\frac{d}{dt}\int_{\Omega_t} \rho\,dx
&= \int_{\Omega_t} \pd_t \rho + \dv(\rho\bw)\,dx
= \int_{\Omega_t} \dv(\rho(\bw-\bv))\,dx \\
&= \int_{\Gamma_t} \rho(\bw - \bv)\cdot\bn_t\,d\omega,
\end{align*}
where $d\omega$ is the area element on $\Gamma_t$. 
Thus, in order that the mass conservation \eqref{mass:1}
holds, it is sufficient to assume that 
\begin{equation}\label{mass:2}
\rho(\bw - \bv)\cdot\bn_t = 0.
\end{equation}
Since $\rho$ is the mass density, we may assume that 
$\rho > 0$, and so,  we have
\begin{equation}\label{mass:3}
\bw\cdot\bn_t = \bv\cdot\bn_t 
\quad\text{on $\Gamma_t$}.
\end{equation}
Let $V_{\Gamma_t} = \bw\cdot\bn_t$, which represents the velocity of
the evolution of the free surface $\Gamma_t$ in the $\bn_t$ direction, 
and then \eqref{mass:3} is written as  
\begin{equation}\label{mass:4}
V_{\Gamma_t} = \bv\cdot\bn_t \quad\text{on $\Gamma_t$}
\end{equation}
for $t \in (0, T)$.  This is called a {\bf kinematic condition},
or non-slip condition. 
Let us assume that $\Gamma_t$ is
represented by $F(x, t) = 0$ for $x \in \Gamma_t$ locally. Since 
$\Gamma_t = \{x = \phi_t(y) \mid y\in \Gamma\}$ for $t \in (0, T)$,
we have $F(\phi_t(y), t) = 0$ for $y \in \Gamma$. 
Differentiation of this formula with respect to $t$ yields that
$$\pd_tF + (\pd_t\phi_t)\cdot \nabla_xF = \pd_tF + \bw\cdot\nabla F = 0.
$$
Since $\bn_t = \nabla_xF/|\nabla_xF|$, it follows from 
\eqref{mass:3} that $\bw\cdot\nabla_xF= \bv\cdot\nabla_xF$ on 
$\Gamma_t$, and so 
\begin{equation}\label{mass:5}
\pd_t F + \bv\cdot\nabla F
=0
\quad\text{on $\bigcup_{0 < t < T}\Gamma_t\times\{t\}$}.
\end{equation}
This expresses the fact that the free surface $\Gamma_t$ consists for all $t > 0$ of the 
same fluid particles, which do not leave it and are not incident on it from
inside $\Omega_t$. 

Secondly, we consider the conservation of momentum:
\begin{equation}\label{mom:1}
\frac{d}{dt}\int_{\Omega_t} \rho v_i\,dx = 0 
\quad(i=1, \ldots, N).
\end{equation} 
As $\pd_t\rho + \dv(\rho\bv) = 0$ (cf. \eqref{eq:mass}), 
 Eq. \eqref{eq:momentum} is 
rewritten as 
$$\pd_t(\rho v_i) + \dv(\rho v_i\bv)-
\dv \bS_i(\bv, \fp)= 0.
$$
And then, by \eqref{rey:1} 
\begin{align*}
\frac{d}{dt}\int_{\Omega_t}\rho v_i\,dx
&= \int_{\Omega_t}\pd_t(\rho v_i) + \dv(\rho v_i\bw) \,dx \\
&= \int_{\Omega_t}\dv(\rho v_i(\bw - \bv))
+\dv \bS_i(\bv, \fp)\,dx\\
& = \int_{\Gamma_t} \rho v_i(\bw-\bv)\cdot\bn_t\,d\omega
+ \int_{\Gamma_t} \bS_i(\bv, \fp)\cdot\bn_t\,d\omega.
\end{align*}
Thus, by \eqref{mass:2}, 
$$\frac{d}{dt}\int_{\Omega_t} \rho v_i\,dx = 
\int_{\Gamma_t} \bS_i(\bv, \fp)\cdot\bn_t\,d\omega.
$$
To obtain \eqref{mom:1}, it suffices to assume that 
$$\int_{\Gamma_t} \bS_i(\bv, \fp)\cdot\bn_t \,d\omega = 0.
 \quad(i = 1, \ldots, N).
$$

 Let 
$H(\Gamma_t)$ and $\Delta_{\Gamma_t}$ be  the $N-1$ times  mean curvature
of $\Gamma_t$ and the Laplace-Beltrami operator on $\Gamma_t$.
We know that 
\begin{equation}\label{rap:1}
 H(\Gamma_t)\bn_t  =  
\Delta_{\Gamma_t}x
\quad(x \in \Gamma_t).
\end{equation}
In this lecture note, our boundary conditions are
\begin{equation}\label{bound:1}
\bS(\bv, \fp)\bn_t = \sigma H(\Gamma_t)\bn_t,
\quad V_{\Gamma_t} = \bv\cdot\bn_t
\quad\text{on $\Gamma_t$}
\end{equation}
for $t \in (0, T)$.  Where, $\bS(\bv, \fp)\bn_t
= {}^\top(\bS_1(\bv, \fp)\cdot\bn_t, \ldots, 
\bS_N(\bv, \fp)\cdot\bn_t)$, 
and by \eqref{bound:1} 
$\bS_i(\bv, \fp)\cdot\bn_t =  \sigma \Delta_{\Gamma_t}x_i$
for $i=1, \ldots, N$. In particular, 
$$\int_{\Gamma_t} \bS_i(\bv, \fp)\cdot\bn_t\,d\omega
= \int_{\Gamma_t} \Delta_{\Gamma_t}x_i\,d\omega =0.
$$
Thus, the conservation of momentum \eqref{mom:1} holds.

We will show two more identities. 
By \eqref{eq:mass}, \eqref{eq:momentum}, and \eqref{bound:1}, 
we have 
\begin{gather}\label{mom:2}
\frac{d}{dt}\int_{\Omega_t}\rho(x_iv_j -x_jv_i)\,dx = 0, \\
\frac{d}{dt}\int_{\Omega_t}\rho x_i\,dx = \int_{\Omega_t}
\rho v_i\,dx.
\label{mom:3}
\end{gather}
The formula \eqref{mom:2} is called the conservation of angular momentum. 
In fact, by \eqref{rey:1} 
\allowdisplaybreaks{
\begin{align*}
&\frac{d}{dt}\int_{\Omega_t}\rho(x_iv_j - x_jv_i)\,dx \\
& = \int_{\Omega_t}\{\pd_t(\rho(x_iv_j-x_jv_i))
+ \dv(\rho(x_iv_j-x_jv_i)\bw)\}\,dx\\
& = \int_{\Omega_t}\{(\pd_t\rho)(x_iv_j-x_jv_i)
+ x_i\rho\pd_tv_j - x_j\rho\pd_tv_i 
+ \dv(\rho(x_iv_j-x_jv_i)\bw)\}\,dx \\
& = \int_{\Omega_t}\{-(\dv(\rho\bv))(x_iv_j-x_jv_i) 
-\rho(x_i\bv\cdot\nabla v_j - x_j\bv\cdot\nabla v_i)\\
&\qquad +x_i\dv\bS_j(\bv, \fp) - x_j\dv\bS_i(\bv, \fp)
+ \dv(\rho(x_iv_j-x_jv_i)\bw)\}\,dx \\
& = \int_{\Omega_t}[ \dv\{\rho(x_iv_j-x_jv_i)(\bw-\bv)\}
+\dv(x_i\bS_j(\bv, \fp))-\dv(x_j\bS_i(\bv, \fp))\\
&\qquad +\sum_{k=1}^N(\delta_{ik}S_{jk}(\bv, \fp)
-\delta_{jk}S_{ik}(\bv, \fp)]\,dx,
\end{align*}
}
where $S_{jk}(\bv, \fp) = D_{jk}(\bv) + ((\nu-\mu)\dv\bv - \fp)\delta_{jk}$.
By \eqref{mass:3}, we have
$$\int_{\Omega_t}\dv\{\rho(x_iv_j-x_jv_i)(\bw-\bv)\}\,dx
= \int_{\Gamma_t} \rho(x_iv_j-x_jv_i)(\bw-\bv)\cdot\bn_t\,d\omega = 0.
$$
By  \eqref{rap:1} and \eqref{bound:1}, we have 
\begin{align*}
&\int_{\Omega_t}\dv(x_i\bS_j(\bv, \fp))-\dv(x_j\bS_i(\bv, \fp))\,dx \\
&= \int_{\Gamma_t} \{x_i\bS_j(\bv, \fp)\cdot\bn_t - x_j \bS_i(\bv, \fp)
\cdot\bn_t\} \,d\omega\\
& = \sigma\int_{\Gamma_t}\{x_i\Delta_{\Gamma_t}x_j - 
x_j\Delta_{\Gamma_t}x_i\}\,d\omega \\
&= -\sigma\int_{\Gamma_t} \nabla_{\Gamma_t}x_i \cdot \nabla_{\Gamma_t}x_j
-\nabla_{\Gamma_t}x_j \cdot \nabla_{\Gamma_t}x_i \}\,d\omega
= 0.
\end{align*}
Since $S_{ij}(\bv, \fp) = S_{ji}(\bv, \fp)$, we have 
$$\int_{\Omega_t}\{
\sum_{k=1}^N(\delta_{ik}S_{jk}(\bv, \fp)
-\delta_{jk}S_{ik}(\bv, \fp))\,dx = 0.
$$
Thus, we have \eqref{mom:2}.

Analogously, by \eqref{rey:1}, \eqref{eq:mass},  and 
\eqref{mass:3}, we have
\begin{align*}
&\frac{d}{dt}\int_{\Omega_t}\rho x_i\,dx
= \int_{\Omega_t}\{\pd_t(\rho x_i)
+\dv(\rho x_i\bw)\}\,dx \\
&= \int_{\Omega_t}\{(\pd_t\rho)x_i + \dv(\rho x_i\bw)\}\,dx
= \int_{\Omega_t}\{-(\dv(\rho\bv))x_i
+\dv(\rho x_i\bw)\}\,dx \\
&= \int_{\Omega_t} \{\dv(x_i\rho(\bw-\bv)) + \rho v_i\}\,dx
= \int_{\Gamma_t}\rho x_i(\bw-\bv)\cdot\bn_t\,d\omega
+ \int_{\Omega_t} \rho v_i\,dx \\
&= \int_{\Omega_t} \rho v_i\,dx.
\end{align*}
Thus, we have \eqref{mom:3}. 

Finally, assuming that $\rho=1$ 
and adding  the initial conditions:
\begin{equation}\label{initial:1}
\bv|_{t=0} = \bv_0 \quad\text{in $\Omega$}, \quad 
\Omega_t|_{t=0} = \Omega,
\end{equation}
we have the set of equations in Eq. \eqref{navier:1}.

\subsection{Short history}

The problems \thetag{P1} and \thetag{P2}
have been studied by many mathematicians. 
In  case of \thetag{P1}, the initial domain, $\Omega$,  is bounded.  
In the case that $\sigma > 0$,
the local in time unique existence theorem  was proved by Solonnikov 
\cite{Sol1, Sol2, Solonnikov7, Sol3, Sol4, Sol5}, and Padula and 
Solonnikov \cite{Pad-Sol2} 
in the $L_2$ Sobolev-Slobodetskii space, by
Schweizer \cite{Schweizer} in the semigroup setting, by 
Moglilevski\u\i \, and Solonnikov \cite{Mog-Sol1, Mog-Sol2, Sol4} 
in the H\"older spaces. In the case that 
$\sigma=0$, the local in time unique existence
theorem was proved by 
Solonnikov \cite{Solonnikov6} and  Mucha and Zaj\c{a}czkowski 
\cite{MZ2} in the $L_p$ Sobolev-Slobodetskii space, and by Shibata 
and Shimizu \cite{Shibata-Shimizu3, SS2}
in the $L_p$ in time and $L_q$ in space 
setting.

 The global in time unique 
existence theorem  was proved,
in the case that $\sigma=0$, 
 by Solonnikov \cite{Solonnikov6} in the 
$L_p$ Sobolev-Slobodetskii space, and by Shibata 
and Shimizu \cite{Shibata-Shimizu3, SS2}
in the $L_p$ in time and $L_q$ in space
under the assumption that the initial velocities 
are small and orthogonal to the rigid 
space,  $\{\bu \mid 
\bD(\bu) = 0\}$, 
and, in the case that $\sigma >0$,  was proved  
by Solonnikov \cite{Solonnikov5} 
in the $L_2$ Sobolev-Slobodetskii space, 
  by  Padula and Solonnikov \cite{Pad-Sol} in the H\"older spaces, and 
by Shibata \cite{S5} in the $L_p$ in time and $L_q$ in space 
setting 
under the assumptions that
 the initial domain, $\Omega$, is sufficiently close to a 
ball and the initial velocities are 
small and orthogonal to the rigid space.

In case of \thetag{P2}, 
 the initial domain, $\Omega$, is a  
perturbed layer given by $\Omega = \{ x \in \BR^N \mid 
-b < x_N < \eta(x'), \enskip x' = (x_1, \ldots, x_N)
\in \BR^{N-1}\}$.   The local in time 
unique existence theorem 
was proved by   Beale \cite{Beale1, Beale2}, Allain \cite{Allain} and 
Tani \cite{Tani} in the $L_2$ Sobolev-Slobodetskii space
when $\sigma > 0$  
and by Abels \cite{Abels} in the 
$L_p$ Sobolev-Slobodetskii space when $\sigma =0$.  
The global in time unique existence theorem for small  initial velocities 
 was  proved in the $L_2$ 
Sobolev-Slobodetskii space 
by Beale \cite{Beale1, Beale2} and Tani and  Tanaka \cite{Tani-Tanaka} 
in the case that 
$\sigma > 0$ 
and by Saito \cite{Saito0} in the $L_p$ in time and 
$L_q$ in space setting in the case that $\sigma = 0$. 
The decay rate was studied
 by Beale and Nishida \cite{Beale-Nishida} (cf. also \cite{Hataya1} for the 
detailed proof), Lynn and 
Sylvester \cite{Sy1}, Hataya \cite{Hataya} and 
Hataya-Kawashima \cite{Hataya2} in the $L_2$ framework with polynomial decay,
and by Saito \cite{Saito0} with exponential decay.

Recently, the local well-posedness in general unbounded domains was proved 
by Shibata \cite{S3} and \cite{SSu1} under the assumption that the initial 
domain is uniformly $C^3$ and 
the weak Dirichlet problem is uniquely solvable in the initial domain.
To transform Eq.\eqref{navier:1} to the problem in the reference domain,
in \cite{S3} the Lagrange transform was used and in \cite{SSu1} the Hanzawa
transform was used. In this lecture note, the latter method will be explained.
Moreover, in the case that $\Omega_t = 
\{x \in \BR^N \mid x_N < \eta(x', t)
\enskip x'=(x_1, \ldots, x_{N-1}) \in \BR^{N-1}\}$,
which is  corresponding to the ocean problem without
bottom physically,  
the global in time unique existence theorem for Eq. \eqref{navier:1}
was proved by Saito and Shibata \cite{Saito1, Saito2}, and 
in the case that $\Omega$ is an exterior domain
in $\BR^N$ ($N \geq 3$), the local and global in time 
unique existence theorems were proved by
Shibata \cite{S8, S6, S7}. In such unbounded domains case, 
we can only show that the $L_q$ space norm of
solutions of the Stokes equations with free boundary conditions
decay polynomially, 
and so  we have to choose different exponents $p$ and $q$ to guarantee
the $L_p$ integrability on time interval $(0, \infty)$.
This is one of the reasons why the maximal $L_p$-$L_q$ regularity theory 
with freedom of choice of 
$p$ and $q$ 
is necessary in the study of the free boundary problems.
In Sect. 8 below, it is explained how to prove the global
in time unique existence theorem by 
combining the decay properties of the lower order terms
and the maximal $L_p$-$L_q$ regularity of the highest order
terms.  Finally, the author is  honored to refer the readers to 
the lecture note written by  Solonnikov, who
is a pioneer of the study of free boundary problems of the Navier-
Stokes equations,  and Denisova \cite{SD}, where  
it is given 
 an overview of studies achieved  by Solonnikov and his followers
about the free boundary problems of the 
Navier-Stokes equations in a bounded domain and related topics,

We remark that 
the two phase fluid flows separated by sharp interface problem
has been studied by Abels \cite{Abels2}, Denisova and Solonnikov  
\cite{De1, De3, DeSo2, DeSo3}, Giga and Takahashi \cite{GiTa, Taka},
Nouri and Poupaud \cite{NP}, Pruess,  Simonett et al. \cite{ KPW, 
PS1, PS2, PS17, PS18},
Shibata and Shimizu \cite{Shibata-Shimizu2}, 
Shimizu \cite{Sh08, Shp1},  Tanaka \cite{Tana} 
and references therein.

\subsection{Further Notation}

This section is ended by explaining further notation used in this
lecture note.
We denote the sets of all complex numbers, real numbers, integers, 
 and natural numbers by
$\BC$, $\BR$, $\BZ$, and $\BN$, respectively. Let $\BN_0 = \BN \cup\{0\}.$
For any multi-index $\alpha = (\alpha_1, \ldots, \alpha_N) 
\in \BN_0^N$ we set 
$\partial_x^\alpha h = 
\partial_1^{\alpha_1}\cdots \partial_N^{\alpha_N} h$
with $\partial_i = \partial/\partial x_i$.  
For any scalor function $f$, we write
\begin{align*}
&\nabla f= (\pd_1f, \ldots, \pd_Nf), \quad 
\bar\nabla f=(f, \pd_1f, \ldots, \pd_Nf),\\
&\nabla^nf = (\pd_x^\alpha f \mid |\alpha|=n), \quad 
\bar\nabla^n f = (\pd_x^\alpha f \mid |\alpha| \leq n)
\quad(n \geq 2),
\end{align*}
where $\pd_x^0f = f$. For any $m$-vector of function $\bff
={}^\top(f_1, \ldots, f_m)$, we write 
\begin{align*}
&\nabla \bff= (\nabla f_1, \ldots, \nabla f_m), \quad 
\bar\nabla \bff=(\bar\nabla f_1, \ldots, \bar\nabla f_m),\\
&\nabla^n\bff = (\nabla^n f_1, \ldots, \nabla^nf_m), \quad 
\bar\nabla^n\bff = (\bar\nabla^nf_1,\ldots, \bar\nabla^nf_m). 
\end{align*}
For any $N$-vector of functions,
$\bu={}^\top(u_1, \ldots, u_N)$, sometime $\nabla\bu$ is regarded as
an $N\times N$ matrix of functions whose $(i, j)^{\rm th}$ component is
$\pd_iu_j$, that is 
$$\nabla \bu = \left(\begin{matrix} 
\pd_1u_1 & \pd_2u_1 & \cdots & \pd_Nu_1 \\
\pd_1u_2 & \pd_2u_2 & \cdots & \pd_Nu_2 \\
\vdots & \vdots & \ddots & \vdots \\
\pd_1u_N & \pd_2u_N & \cdots & \pd_Nu_N 
\end{matrix}\right).
$$
For any $m$-vector $V=(v_1, \ldots, v_m)$ and $n$-vector
$W=(w_1, \ldots, w_n)$, $V\otimes W$ denotes an $(m,n)$ matrix 
whose $(i, j)^{\rm th}$ component is $V_iW_j$.  For any 
$(mn, N)$ matrix $A=(A_{ij, k} \mid i=1, \ldots, m, j=1, \ldots, n,
k=1, \ldots, N)$, 
$AV\otimes W$ denotes an $N$ column vector whose $i^{\rm th}$ component is
the quantity: $\sum_{j=1}^m\sum_{j=1}^n A_{jk, i}v_jw_k$.  

For any $N$ vector $\ba$, $\ba_i$ denotes the $i^{\rm th}$ component
of $\ba$ and for an $N\times N$ matrix $\bA$, $\bA_{ij}$ denotes 
the $(i, j)^{\rm th}$ component of $\bA$, and moreover, 
the $N\times N$ matrix whose
$(i, j)^{\rm th}$ component is $K_{ij}$ is written as  $(K_{ij})$.
Let $\delta_{ij}$
be the Kronecker delta symbol, that is $\delta_{ii} = 1$
and $\delta_{ij} = 0$ for $i\not=j$. In particular, 
$\bI = (\delta_{ij})$ is the $N \times N$ identity matrix.  
Let $\ba\cdot \bb =<\ba, \bb>= \sum_{j=1}^N\ba_j\bb_j$
for any $N$-vectors $\ba$ and $\bb$.  
For any $N$-vector $\ba$, let $\BPi_0\ba = \ba_\tau
: = \ba - <\ba, \bn>\bn$.  For any two $N\times N$ matrices 
$\bA$ and $\bB$, we write  $\bA:\bB = {\rm tr}\bA\bB
: = \sum_{i,j=1}^N\bA_{ij}\bB_{ji}$. Given  
$1 < q < \infty$, let $q' = q/(q-1)$. 
For $L > 0$, let $B_L = \{x \in \BR^N \mid |x| < L\}$ and
$S_L = \{x \in \BR^N \mid |x| = L \}$. Moreover, for
$L < M$, we set $D_{L,M} = \{x \in \BR^N \mid L < |x| < M\}$.

For any domain $G$ in $\BR^N$, let $L_q(G)$, $H^m_q(G)$,  and 
$B^s_{q,p}(G)$  
be the standard Lebesgue, Sobolev, and 
Besov spaces on $G$, 
and let  $\|\cdot\|_{L_q(G)}$, 
$\|\cdot\|_{H^m_q(G)}$, and  
$\|\cdot\|_{B^s_{q,p}(G)}$ denote their respective norms.
We write $L_q(G)$ as $H^0_q(G)$, and $B^s_{q,q}(G)$ as simply  
$W^s_q(G)$. 
For a Banach space $X$ with norm $\|\cdot\|_X$  
and for $1 \leq p \leq \infty$, $L_p((a, b), X)$ and 
$H^m_p((a, b), X)$ denote the standard  Lebesgue and Sobolev spaces of
$X$-valued functions defined on an interval $(a, b)$, 
and $\|\cdot\|_{L_p((a, b), X)}$,  
$\|\cdot\|_{H^m_p((a, b), X)}$ denote their respective norms. 
For $\theta \in (0, 1)$, $H^\theta_p(\BR, X)$ denotes the standard
$X$-valued Bessel potential space defined by 
\begin{align*}
H^\theta_p(\BR, X) &= \{f \in L_p(\BR, X) \mid 
\|f\|_{H^\theta_p(\BR, X)} < \infty\}, \\
\|f\|_{H^\theta_p(\BR, X)} &= \Bigl(\int_{\BR}
\|\CF^{-1}[(1+\tau^2)^{\theta/2}\CF[f](\tau)](t)\|_X^p\,dt
\Bigr)^{1/p},
\end{align*}
where $\CF$ and $\CF^{-1}$ denote the Fourier transform and 
the inverse Fourier transform, respectively. 
Let $X^d = \{(f_1, \ldots, f_d) \mid f_i \in X \enskip
 (i=1, \ldots, d)\}$, 
and write the norm of $X^d$ as simply  $\|\cdot\|_X$, which is 
defined by    
$\|f\|_X = \sum_{j=1}^d \|f_j\|_X$ for 
$f = (f_1, \ldots, f_d)
\in X^d$.

For two Banach spaces $X$ and $Y$, $X+Y = \{x + y \mid x \in X, y\in Y\}$,
$\CL(X, Y)$ denotes the 
set of all bounded linear operators from $X$ into $Y$ and 
$\CL(X, X)$ is written simply as $\CL(X)$. 
For a domain $U$ in $\BC$, $\Hol(U, \CL(X, Y))$
 denotes the set of all $\CL(X, Y)$-valued holomorphic 
functions defined on $U$. 
Let 
$\CR_{\CL(X, Y)}(\{\CT(\lambda) \mid \lambda \in U\})$
be the $\CR$ norm of the operator family 
$\CT(\lambda) \in \Hol(U, \CL(X, Y))$.
Let 
\begin{gather*}
\Sigma_{\epsilon_0} = \{\lambda \in \BC\setminus\{0\} \mid
|\arg\lambda| \leq \pi-\epsilon_0\}, 
\quad \Sigma_{\epsilon_0, \lambda_0}
= \{\lambda \in \Sigma_{\epsilon_0} \mid |\lambda| \geq \lambda_0\}, \\
\BC_{+, \lambda_0} = \{\lambda \in \BC \mid 
{\rm Re}\,\lambda \geq \max(0, \lambda_0)\}.
\end{gather*} 
Let $C^\infty_0(G)$ be the set of all $C^\infty$ functions
whose supports are compact and contained in $G$. 
Let 
$(\bu, \bv)_G = \int_G\bu\cdot\bv\,dx$ and   
$(\bu,\bv)_{\partial G} = \int_{\partial G} \bu\cdot \bv\,ds$, 
where $ds$ denotes the surface element on $\partial G$ and $\pd G$ is 
the boundary of $G$.   
For $T > 0$, $G\times(0, T) = \{(x, t) \mid x \in G, 
t \in (0, T)\}$ is written simply by $G^T$. 
Let
$\BR^N_+ = \{x = (x_1, \ldots, x_N) \mid x_N > 0\}$ and 
$\BR^N_0 = \{x = (x_1, \ldots, x_N) \mid x_N = 0\}$.
The letter $C$ denotes a
generic constant and $C_{a,b,c,\cdots}$ denotes that the 
constant $C_{a,b,c,\cdots}$ depends 
on $a$, $b$, $c, \cdots$.  
The value of 
$C$ and $C_{a,b,c,\cdots}$ may change from line to line.


\section{Preliminaries} \label{sec:p}
In this section, we study a uniform $C^k$ domain, a uniorm $C^k$ domain
 whose inside has
a finite covering, 
the weak Dirichlet problem,  Besov spaces on the boundary, and the Laplace-
Beltrami operator on the boundary.
\subsection{Definitions of domains}\label{subsec:p-1}
We first introduce the definition of a uniform  $C^k$ ($k=2$ or $3$) domain. 
\begin{dfn}\label{dfn:1}
Let $k \in \BN$.  We say that $\Omega$ is a uniform $C^k$ domain,
if there exist positive constants $a_1$, $a_2$, and $A$ such 
that the following assertion holds:  For any $x_0 = (x_{01},
\ldots, x_{0N}) \in \Gamma$ there exist a coordinate number
$j$ and a $C^k$ function $h(x')$ defined on $B'_{a_1}(x_0')$ such that 
$\|h\|_{H^k_\infty(B'_{a_1}(x_0'))} \leq A$ and
\begin{align*}
\Omega\cap B_{a_2}(x_0) & = \{x \in \BR^N \mid x_j > h(x') \enskip
(x' \in B_{a_1}(x_0')) \} \cap B_{a_2}(x_0), \\
\Gamma \cap B_{a_2}(x_0) & = \{x \in \BR^N \mid x_j = h(x') \enskip
(x' \in B_{a_1}'(x_0')) \} \cap B_{a_2}(x_0).
\end{align*}
Here, we have set 
\begin{gather*}
y' = (y_1, \ldots, y_{j-1}, y_{j+1}, \ldots, y_N) \enskip (y \in \{x, x_0\}),
 \\
B'_{a_1}(x'_0)  = \{x' \in \BR^{N-1} \mid |x' - x'_0| < a_1\}, \\
B_{a_2}(x_0) = \{x \in \BR^N \mid |x-x_0| < a_2\}.
\end{gather*}
\end{dfn}
The uniform $C^k$ domains are characterized as follows.
\begin{prop}\label{prop:lap} Let $k = 2$ or $k=3$. 
Let $\Omega$ be a uniform $C^k$ domain in $\BR^N$.
Then, for any $M_1 \in (0, 1)$, 
 there exist  constants $M_2>0$ and $0 < r_0 < 1$, 
at most countably many $N$-vector of functions 
$\Phi_j \in C^k(\BR^N)^N$  and points 
$x^0_j \in \Omega$ and  $x^1_j \in \Gamma$  
such that the following assertions hold:
\begin{itemize}
\item[\thetag{i}]~The maps: $\BR^N \ni x \mapsto \Phi_j(x) \in \BR^N$
are bijections satisfying the following conditions:
$\nabla\Phi_j = \CA_j + B_j$, 
$\nabla(\Phi_j)^{-1}=\CA_{j,-}
+ B_{j,-}$, where $\CA_j$ and $\CA_{j,-}$ are $N\times N$ constant 
orthogonal matrices, and $B_j$ and $B_{j,-}$ are $N\times N$ matrices of 
$C^{k-1}(\BR^N)$ functions defined on $\BR^N$ satisfying 
the conditions:
$\|(B_j, B_{j,-})\|_{L_\infty(\BR^N)} \leq M_1$, and
$\|\nabla (B_j, B_{j,-})\|_{L_\infty(\BR^N)} \leq C_A$,
where $C_A$ is a constant depending on constants $A$, $\alpha_1$ and $\alpha_2$
appearing in Definition \ref{dfn:1}.  Moreover, if $k=3$, then  
$\|\nabla^2(B_j, B_{j-}) \|_{L_\infty(\BR^N)} \leq M_2$.
\item[\thetag{ii}]~$\Omega = \Bigl(\bigcup_{j=1}^\infty B_{r_0}(x^0_j)\Bigr)
\cup \Bigl(\bigcup_{j=1}^\infty(\Phi_j(\HS) \cap B_{r_0}(x^1_j))\Bigr)$, 
$B_{r_0}(x^0_j) \subset \Omega$, 
$\Phi_j(\BR^N_0)\cap B_{r_0}(x^1_j) = \Gamma\cap B_{r_0}(x^1_j)$.
\item[\thetag{iii}]~There exist $C^\infty$ functions $\zeta^i_j$ and  
$\tilde\zeta^i_j$  $(i=0,1,\enskip j \in \BN)$ such that 
\end{itemize}
\vskip-1pc
\begin{gather*}
0 \leq \zeta^i_j, \enskip \tilde\zeta^i_j \leq 1,\enskip 
{\rm supp}\, \zeta^i_j \subset {\rm supp}\, \tilde \zeta^i_j
\subset B_{r_0}(x^i_j),  \enskip 
\tilde\zeta^i_j = 1 \,\text{on ${\rm supp}\,\zeta^i_j$},
\\
\sum_{i=0}^1\sum_{j=1}^\infty \zeta^i_j = 1\,
\text{on $\overline{\Omega}$},
\sum_{j=1}^\infty \zeta^1_j = 1\,\text{on $\Gamma$}, 
\|\nabla\zeta^i_j\|_{H^{k-1}_\infty(\BR^N)}, 
\, \|\nabla\tilde\zeta^i_j\|_{H^{k-1}_\infty(\BR^N)} \leq M_2.
\end{gather*}
\begin{itemize}
\item[\thetag{iv}]~\textcolor{red}{Below, 
for the notational simplicity, we set 
$B^i_j = B_{r_0}(x^i_j)$. }
For each $j$, let $\ell^{ij}_k$ $(k=1, \ldots, m^i_j)$ 
be numbers  for which $B^i_j \cap B^i_{\ell^{ij}_k} \not=\emptyset$ and 
$B^i_j \cap B^i_m = \emptyset$ for $m \not\in \{\ell^{ij}_k\ \mid 
k=1, \ldots, m^i_j\}$. Then, there exists an $L \geq 2$ independent of 
$M_1$ such that $m^i_j \leq L$. 
\end{itemize}
\end{prop}
\pf Proposition \ref{prop:lap} was essentially proved by 
Enomoto-Shibata \cite[Appendix]{ES1}
 instead of $\|\nabla(B^i_j, 
B^i_{j,-})\|_{L_\infty(\BR^N)} \leq C_A$.
There, it was proved that 
$$\|\nabla(B^i_j, 
B^i_{j,-})\|_{L_\infty(\BR^N)} \leq CM_2,
$$
that is the estimate of 
$\nabla(B^i_j, B^i_{j,-})$ depends on $M_1$. Since the proof
of Proposition \ref{prop:lap} is almost the same as in
\cite[Appendix]{ES1}, we shall give an idea
how to improve this point below.  Let $x_0 =(x_0', x_{0N})\in \Gamma$ 
and we assume that
\begin{align*}
\Omega \cap B_\beta(x_0) = \{x \in \BR^N \mid
x_N > h(x') \enskip(x' \in B'_\alpha(x_0'))\} \cap B_\beta(x_0), \\
\Gamma \cap B_\beta(x_0) = \{x \in \BR^N \mid
x_N = h(x') \enskip(x' \in B'_\alpha(x_0'))\} \cap B_\beta(x_0).
\end{align*}
We only consider the case where $k=3$. In fact, by the same argument, we can 
improve the estimate in the case where $k=2$.  We assume that 
$h \in C^3(B'_\alpha(x_0'))$, $\|h\|_{H^3_\infty(B'_\alpha(x_0'))} \leq K$,
and $x_{0N} = h(x_0')$. Below, $C$ denotes a generic constant depending on
$K$, $\alpha$ and $\beta$ but independent of $\epsilon$. Let $\rho(y)$ be
a function in $C^\infty_0(\BR^N)$ such that $\rho(y) = 1$ for $|y'| \leq 1/2$
and $|y_N| \leq 1/2$ and $\rho(y) = 0$ for $|y'| \geq 1$ or  $|y_N| \geq 1$.
Let $\rho_\epsilon(y) = \rho(y/\epsilon)$.  We consider a $C^\infty$
diffeomorphism:
$$x_j = \Phi_j^\epsilon(y)  = x_{0j} + \sum_{k=1}^Nt_{j,k}y_k
+ \sum_{k,\ell=1}^Ns_{j, k\ell}y_ky_\ell \rho_\epsilon(y).
$$
Where, $t_{j,k}$ and $s_{j, k\ell}$ are some constants satisfying the
conditions \eqref{tr:1*}, \eqref{tr:3}, and \eqref{tr:4},
below.  Let 
$$G_\epsilon(y) = \Phi^\epsilon_N(y) - h(\Phi^\epsilon_1(y), \ldots, 
\Phi^\epsilon_{N-1}(y)).
$$
Notice that $G_\epsilon(0) = x_{0N} - h(x_0') = 0$.  
We choose $t_{j,k}$ and $s_{\ell, mn}$in such a way that
\begin{equation}\label{tr:3}
\begin{aligned}
\frac{\pd G_\epsilon}{\pd y_N}(0) &= t_{N,N} - \sum_{k=1}^{N-1}
\frac{\pd h}{\pd x_k}(x_0')t_{k, N} \not=0,
 \\ 
\frac{\pd G_\epsilon}{\pd y_j}(0)  &= t_{N,j} - 
\sum_{k=1}^{N-1}\frac{\pd h}{\pd x_k}(x_0')t_{k,j} = 0;
\end{aligned}
\end{equation}
\begin{equation}\label{tr:4}
\begin{aligned}
\frac{\pd^2 G_\epsilon}{\pd y_\ell\pd y_m}(0)
 &= s_{N, \ell m}+s_{N, m\ell}
-\sum_{k=1}^{N-1}\frac{\pd h}{\pd x_k}(x_0')(s_{k,\ell m}
+ s_{k, m \ell})\\
&- \sum_{j,k=1}^{N-1}\frac{\pd^2 G_\epsilon}{\pd x_j\pd x_k}(x'_0)
t_{j,\ell}t_{k,m} = 0.
\end{aligned}
\end{equation}
Moreover, setting 
$$T = \left(\begin{matrix}
t_{1,1} & t_{2,1} &\cdots & t_{N,1} \\
t_{1,2} & t_{2,2} & \cdots & t_{N,2} \\
\vdots & \vdots & \ddots & \vdots \\
t_{1,N} & t_{2,N} & \cdots & t_{N,N}
\end{matrix}\right)
$$
we assume that $T$ is an orthogonal matrix, that is
\begin{equation}\label{tr:1*}
\sum_{\ell=1}^N t_{\ell, m}t_{\ell, n}
= \delta_{mn} = \begin{cases}
1 &\quad\text{for $m=n$}, \\
0 &\quad\text{for $m\not=n$}.
\end{cases}
\end{equation}
We write $\dfrac{\pd h}{\pd x_j}(x_0')$ simply by $h_j$ and 
set 
$$H_j = \sqrt{1 + \sum_{\ell=j}^{N-1}h_\ell^2} = 
\sqrt{1 + h_j^2 + h_{j+1}^2 + \cdots+ h_{N-1}^2}.
$$
Let 
\begin{align*}
t_{N,N-j} &= \frac{h_{N-j}}{H_{N-j}H_{N+1-j}}, \quad
t_{N-k, N-j} = -\frac{h_{N-k}h_{N-j}}{H_{N-j}H_{N+1-j}} 
\intertext{for $k=1, \ldots, j-1$, and}
t_{N-j, N-j} & = \frac{H_{N+1-j}}{H_{N-1}}, \quad
t_{k, N-j} = 0, 
\intertext{for $k=1, \ldots, N-j-1$ and $j=1, \ldots, N-1$, and}
t_{i,N} & = -\frac{h_i}{H_1},\quad
t_{N,N} = \frac{1}{H_1}
\end{align*}
for $i=1, \ldots, N-1$. 
Then, we see that such $t_{j,k}$ satisfy \eqref{tr:3}
and \eqref{tr:1*}.  In particular, 
\begin{equation}\label{tr:5} \frac{\pd G_\epsilon}{\pd y_N}(0) 
= \frac{1}{H_1}.
\end{equation}
Moreover, assuming the symmetry: $s_{\ell, jk} = s_{\ell, kj}$, we have
\begin{align*}
s_{N,jk} &= \frac{1}{2H_2}\sum_{m,n=1}^{N-1}
\frac{\pd^2h}{\pd x_m\pd x_n}(x_0')t_{m,j}t_{n,k}, \\
s_{i,jk} &= -\frac{h_i}{2H_1^2}\sum_{m,n=1}^{N-1}
\frac{\pd^2h}{\pd x_m\pd x_n}(x_0')t_{m,j}t_{n,k}.
\end{align*}
By successive approximation, we see that there exists a constant  
$\epsilon_0 > 0$ such that for any $\epsilon \in (0, \epsilon_0)$
there exists a function $\psi_\epsilon \in C^3(B'_\epsilon(0))$
satisfies the following conditions:
\begin{gather}
\psi_\epsilon(0) = \pd_i\psi_\epsilon(0) = \pd_i\pd_j\psi_\epsilon(0) =0,
\nonumber \\
\|\psi_\epsilon\|_{L_\infty(B'_\epsilon(0))} \leq C\epsilon^2, 
\quad
\|\pd_i\psi_\epsilon\|_{L_\infty(B'_\epsilon(0))} \leq C\epsilon, 
\nonumber\\
\|\pd_i\pd_j\psi_\epsilon\|_{L_\infty(B'_\epsilon(0))} \leq C, 
\quad
\|\pd_i\pd_j\pd_k\psi_\epsilon\|_{L_\infty(B'_\epsilon(0))} 
\leq C\epsilon^{-1}, \nonumber \\
G_\epsilon(y', \psi_\epsilon(y')) = 0 \quad\text{for $y' \in B'_\epsilon(0)$},
\label{t:8}
\end{gather}
where $i$, $j$ and $k$ run from $1$ through $N-1$.  Notice that 
\begin{equation}\label{t:10}
\begin{aligned}
&x_N- h_\epsilon(x')  = G_\epsilon(y)\\ 
&\quad = G_\epsilon(y', \psi_\epsilon(y'))\\
&\quad
+ \int^1_0(\pd_NG_\epsilon)(y', \psi_\epsilon(y') +
\theta(y_N-\psi_\epsilon(y')))\,d\theta(y_N - \psi_\epsilon(y'))
\\
&\quad = ((\pd_NG_\epsilon)(0) + \tilde G_\epsilon(y))(y_N-\psi_\epsilon(y')),
\end{aligned}
\end{equation}
where we have used $G_\epsilon(y', \psi_\epsilon(y'))=0$ and 
\begin{align*}
&\tilde G_\epsilon(y)  =\int^1_0\int^1_0\{
\sum_{\ell=1}^{N-1}(\pd_\ell\pd_NG_\epsilon)
(\tau y', \tau(\psi_\epsilon(y') + \theta(y_N-\psi_\epsilon(y')))y_\ell\\
&+ \pd_N^2G_\epsilon
(\tau y', \tau(\psi_\epsilon(y') + \theta(y_N-\psi_\epsilon(y')))
(\psi_\epsilon(y') + \theta(y_N-\psi_\epsilon(y')))\}\,d\theta d\tau.
\end{align*}
Since $(\pd_NG_\epsilon)(0) = 1/H_1$, choosing $\epsilon_0 > 0$ so small
that $|\tilde G_\epsilon(y)| \leq 1/(2H_1)$ for $|y| \leq \epsilon_0$, 
we see that $x_N - h(x') \geq 0$ and $y_N-\psi_\epsilon(y') \geq 0$
are equivalent. 

Let $\omega$ be a function in $C^\infty_0(\BR^{N-1})$ such that 
$\omega(y') = 1$ for $|y'| \leq 1/2$ and $\omega(y') =0$ for
$|y'|\geq 1$ and set $\omega_\epsilon(y') = \psi_\epsilon(y')
\omega(y'/\epsilon)$.  Then, by \eqref{t:8} we have
\begin{equation}\label{t:9}\begin{split}
\|\omega_\epsilon\|_{L_\infty(\BR^{N-1})} &\leq C\epsilon^2, \quad
\|\pd_i\omega_\epsilon\|_{L_\infty(\BR^{N-1})} \leq C\epsilon, \\
\|\pd_i\pd_j\omega_\epsilon\|_{L_\infty(\BR^{N-1})} &\leq C,\quad
\|\pd_i\pd_j\pd_k\omega_\epsilon\|_{L_\infty(\BR^{N-1})} \leq C\epsilon^{-1}.
\end{split}\end{equation}
where $i$, $j$, and $k$ run from $1$ through $N-1$. 
Setting $\Psi^\epsilon(z) = \Phi^\epsilon(z', z_N + \omega_\epsilon(z'))$, 
that is $y_N = z_N+\omega_\epsilon(z')$, and $y_j = 
z_j$ for $j=1, \ldots, N-1$, 
we see that there exists an $\epsilon_0 > 0$ such that for any 
$\epsilon \in (0, \epsilon_0)$, the map: $z \to x=\Psi^\epsilon(z)$ is 
a diffeomorphism of $C^3$ class from $\BR^N$ onto $\BR^N$. 
Since 
$$
\frac{\pd x_m}{\pd z_k}
= t_{m,k} + b_{m,k}, 
\quad
\frac{\pd x_m}{\pd z_N}
=  t_{m,N} + b_{n,N}
$$
where we have set 
\begin{align*}
b_{m,k} & = \sum_{i,j=1}\frac{\pd}{\pd z_k}
(s_{m, ij}y_iy_j\rho_\epsilon(y))\\
&+ \{t_{m,N} + \sum_{i,j=1}^N\frac{\pd}{\pd z_k}
(s_{m,ij}y_iy_j\rho_\epsilon(y))\}
\frac{\pd\omega_\epsilon}{\pd z_k}(z'), \\
b_{n,N}
& = \sum_{i,j=1}\frac{\pd}{\pd z_N}
(s_{m, ij}y_iy_j\rho_\epsilon(y)),
\end{align*}
let $\CA$ and $B$ be the $N\times N$ matrices whose $(m,n)^{\rm th}$
components are $t_{m,n}$ and $b_{m,n}$, respectively.  Then, 
by \eqref{t:9}, $\CA$ is an orthogonal matrix and $B$ satisfies the 
estimates:
$$\|B\|_{L_\infty(\BR^N)} \leq C\epsilon, 
\quad \|\nabla B\|_{L_\infty(\BR^N)} \leq C, \quad
\|\nabla^2B\|_{L_\infty(\BR^N)} \leq C\epsilon^{-1}.
$$
Moreover, by \eqref{t:10} we have
$$x_N- h(x') = (1/H_1 + \tilde G_\epsilon(z', z_N+\omega_\epsilon(z')))
(z_N + (\omega(z'/\epsilon)-1)\psi_\epsilon(z')),
$$
which shows that when $|z'| \leq \epsilon/2$, $x_N \geq h(x')$ and 
$z_N \geq 0$ are equivalent. We can construct the sequences of 
$C^\infty_0(\BR^N)$ functions, $\{\zeta^i_j\}$, $\{\tilde\zeta^i_j\}$,
by standard manner (cf. Enomoto-Shibata \cite[Appendix]{ES1}).
This completes the proof of Proposition \ref{prop:lap}. \qed

To show the {\it a priori} estimates for the Stokes equations with
free boundary condition in a general domain, 
we need some restriction of the domains. In this lecture note, 
we adopt the following conditions. 
\begin{dfn}\label{dfn:2} Let $k=2$ or $3$ and let 
$\Omega$ be a domain in $\BR^N$.  We say that 
$\Omega$ is a uniformly $C^k$ domain whose inside has a 
finite covering if $\Omega$ is a uniformly $C^k$ domain
in the sense of Definition \ref{dfn:1} and the following
assertion hold: 
\item[\thetag{v}]~Let $\zeta^i_j$ be the partition of 
unity given in Proposition \ref{prop:lap} \thetag{iii}
and set $\psi^0 = \sum_{j=1}^\infty \zeta^0_j$. 
Let $\CO= {\rm supp}\,\nabla\psi^0 \cup \bigl(\bigcup_{j=1}^\infty 
 {\rm supp}\, \nabla\zeta^1_j\bigr)$.
Then, there exists a finite number of subdomains $\CO_j$ 
$(j = 1, \ldots, \iota)$ such that $\CO \subset
\bigcup_{j=1}^\iota \CO_j$ and each $\CO_j$
satisfies one of the following conditions:
\begin{itemize}
\item[\thetag{a}]~There exists an $R > 0$ such that
$\CO_j \subset \Omega_R$, where $\Omega_R = \{x \in \Omega \mid
|x| < R\}$, 
\item[\thetag{b}]~There exist a translation $\tau$, a rotation
$\CA$, a domain $D \subset \BR^{N-1}$, a coordinate 
function $a(x')$ defined for $x' \in D$, and 
a positive constant  $b$  such that
$0 \leq a(x') < b$ for $x \in D$,  
\begin{align*}
\CA\circ\tau(\CO_j)
\subset \{x = (x', x_N) \mid 
x' \in D, \enskip a(x') \leq x_N \leq b\}
&\subset \CA\circ\tau(\Omega), \\
\{x=(x', x_N) \in \BR^N \mid 
x' \in D, \enskip x_N=a(x') \}
&\subset \CA\circ\tau(\Gamma).
\end{align*}
Where, for any subset $E$ of $\BR^N$, $\CA(E) = \{Ax \mid x \in E\}$ 
with some orthogonal matrix $A$ and $\tau(E) = \{x+y \mid
x \in E\}$ with some $y \in \BR^N$. 
\end{itemize}
\end{dfn}
\begin{ex} Let $\Omega$ be a domain whose boundary $\Gamma$ 
is a $C^k$ hypersurface.
If $\Omega$ satisfies one of the following conditions, then $\Omega$ is a 
 uniform $C^k$ domain whose inside has a finite covering.
\begin{itemize}
\item[\thetag1]~$\Omega$ is bounded, 
or $\Omega$ is an exterior domain, that is, 
$\Omega = \BR^N\setminus\overline{\CO}$ with some bounded domain $\CO$. 
\item[\thetag2]~$\Omega = \BR^N_+$ (half space), or $\Omega$ is 
a perturbed half space, that is,  there exists an 
$R > 0$ such that $\Omega \cap B^R = \BR^N_+ \cap B^R$, where 
$B^R = \{x \in \BR^N \mid |x| >R\}$. 
\item[\thetag3]~$\Omega$ is a layer $L$ or perturbed layer, that is, 
there exists an $R>0$ such that $\Omega \cap B^R = L \cap B^R$.  Here
$L = \{x = (x', x_N) \in \BR^N \mid x' = (x_1, \ldots, x_{N-1}) \in \BR^{N-1},
\enskip a < x_N < b\}$ for some constants $a$ and $b$ for which $a < b$. 
\item[\thetag4]~$\Omega$ is a tube, that is,  there exists a bounded domain
$D$ in $\BR^{N-1}$ such that $\Omega = D\times\BR$.
\item[\thetag5]~There exist an $R > 0$ and several orthogonal transforms,
$\CR_i$ ($i=1, \ldots, M$), such that $\Gamma \cap B^R 
= \Bigl(\bigcup_{i=1}^M \CR_i\BR^N_0\Bigr) \cap B^R$.
\item[\thetag6]~There exist an $R > 0$, half tubes, $T_i$ $(i=1, \ldots, M)$,
and orthogonal transforms, $\CR_i$ $(i=1, \ldots, M)$, such that
$\Omega \cap B^R = \Bigl(\bigcup_{i=1}^M \CR_i T_i\Bigr) \cap B^R$,
where what $T_i$ is a half tube means that $T_i = D_i \times [0, \infty)$
with some bounded domain $D_i$ of $\BR^{N-1}$. 
\end{itemize}
\end{ex}
\vskip0.5pc
In the following, we  write $B_{r_0}(x^i_j)$, 
$\Phi_j(\BR^N_+)$, and $\Phi_j(\BR^N_0)$ 
simply by $B^i_j$, $\Omega_j$ and $\Gamma_j$, respectively. In view of  
Proposition \ref{prop:lap} \thetag{ii}, we have $\Omega_j \cap B^1_j
=\Omega \cap B^1_j$ and $\Gamma_j \cap B^1_j = \Gamma\cap B^1_j$.  
By the finite intersection property stated in 
 Proposition \ref{prop:lap} \thetag{iv}, for any $r \in [1, \infty)$
there exists a constant $C_{r,L}$ such that 
\begin{equation}\label{D.1}
\Bigl[\sum_{j=1}^\infty 
\|f\|_{L_r(\Omega\cap B^i_j)}^r\Bigr]^{\frac1r} 
\leq C_{r,L}\|f\|_{L_r(\Omega)}
\quad\text{for any $f \in L_r(\Omega)$}.
\end{equation}
Let $n \in \BN_0$, $f \in H^n_q(\Omega)$, and 
let $\eta^i_j$ be functions in $C^\infty_0(B^i_j)$ with
$\|\eta^i_j\|_{H^n_\infty(\BR^N)} \leq c_0$ for 
some constant $c_0$ independent of  $j \in \BN$.
Since 
$\Omega \cap B^1_j = \Omega_j \cap B^1_j$, by \eqref{D.1} 
\begin{equation}\label{D.1*}
\sum_{j=1}^\infty\|\eta^0_jf\|_{H^n_q(\BR^N)}^q
+ \sum_{j=1}^\infty\|\eta^1_jf\|_{H^n_q(\Omega_j)}^q
\leq C_q\|f\|_{H^n_q(\Omega)}^q.
\end{equation}
\subsection{Besov spaces on $\Gamma$} \label{subsec:p-3}

We now define Besov spaces on $\Gamma$. Before turning to it, 
we recall some basic facts.  Let $f$ be a 
function defined on $\Gamma$ such that ${\rm supp}\, f \subset 
\Gamma \cap B^1_j \cap B^1_k$.  Let $f_j = f\circ\Phi_j^{-1}$, and then 
by Proposition \ref{prop:lap} \thetag{iv}, we see that for any $s \in [-1, 3]$
and $j$, $k \in \BN$, 
\begin{equation}\label{def:sol.1}
\|f_j\|_{W^s_q(\BR^N_0)} \leq C_{s,q}\|f_k\|_{W^s_q(\BR^N_0)}.
\end{equation}
In fact, in the case that $s = 0, 1, 2, 3$, noting 
$W^s_q = H^s_q$, we see that  the inequality
\eqref{def:sol.1} follows from the direct calculations. When $s = -1$, 
it follows from duality argument. 
Finally, in the case that $s \not\in \BZ$, the inequality \eqref{def:sol.1}
 follows
from real interpolation. Let $W^s_q(\Gamma_j)$ and its norm
$\|\cdot\|_{W^s_q(\Gamma_j)}$ be defined by 
$$W^s_q(\Gamma_j) = \{f \mid f\circ \Phi_j \in W^s_q(\BR^N_0)\},
\quad \|f\|_{W^s_q(\Gamma_j)} = \|f\circ\Phi^{-1}_j\|_{W^s_q(\BR^N_0)}.
$$
In view of \eqref{def:sol.1}, if ${\rm supp}\,f \subset 
\Gamma \cap B^1_j \cap B^1_k$, then 
\begin{equation}\label{def:sol.2}
\|f\|_{W^s_q(\Gamma_j)} \leq C_{s,q} \|f\|_{W^s_q(\Gamma_k)}.
\end{equation}
For $s \in [-1, 3]$, we now define $W^s_q(\Gamma)$ by
\begin{multline}\label{def:sob1}
W^s_q(\Gamma) = \{f = \sum_{j=1}^\infty f_j \mid  {\rm supp}\, f_j 
\subset \Gamma \cap B^1_j, \quad 
f_j \in W^s_q(\Gamma_j), \\
\|f\|_{W^s_q(\Gamma)}
= \Bigl\{\sum_{j=1}^\infty \|f_j\|_{W^s_q(\Gamma_j)}^q
\Bigr\}^{1/q} < \infty\}.
\end{multline}
Since each $W^s_q(\Gamma_j)$ is a Banach space, so is $W^s_q(\Gamma)$.  
Given 
$f \in W^s_q(\Gamma)$, we have 
\begin{equation}\label{def:sob2}
\tilde\zeta_jf \in W^s_q(\Gamma_j), 
\quad \sum_{j=1}^\infty\|\tilde\zeta_jf\|_{W^s_q(\Gamma_j)}^q
\leq (c_0C_{s,q}3^N)^q\|f\|_{W^s_q(\Gamma)}^q.
\end{equation}
In fact, let $f = \sum_{j=1}^\infty f_j \in W^s_q(\Gamma)$, 
where  $f_j$ satisfy
\eqref{def:sob1}. For each $j$, let $\ell_k^{1j}$ $(k=1, \ldots, m_j^1)$
be the numbers given in Proposition \ref{prop:lap} \thetag{iv}
for which $B^1_j \cap B^1_{\ell^{1j}_k} \not=\emptyset$ and 
$B^1_j \cap B^1_m=\emptyset$ for $m \not\in \{\ell_k^{1j} \mid k = 1, 
\ldots, m^1_j\}$, and then 
$$\tilde\zeta_j f = \sum_{k=1, \ldots, m^1_j} \tilde\zeta_jf_{\ell^{1j}_k}.$$
Since ${\rm supp}\,\tilde\zeta_jf_{\ell^{1j}_k}
\subset \Gamma\cap B^1_{\ell^{1j}_k} \cap B^1_j$, by \eqref{def:sol.2} 
\begin{align*}
&\|\tilde\zeta_jf\|_{W^s_q(\Gamma_j)}
\leq \sum_{k=1, 2, \ldots, m^1_j} 
\|\tilde\zeta_jf_{\ell_k^{1j}}\|_{W^s_q(\Gamma_j)}
\leq C_{s,q} \sum_{k=1, 2, \ldots, m^1_j} \|\tilde\zeta_jf_{\ell_k^{1j}}
\|_{W^s_q(\Gamma_{\ell_k^j})} \\
&\quad \leq c_0C_{s,q}\sum_{k=1, 2, \ldots, m^1_j} \|f_{\ell_k^{1j}}
\|_{W^s_q(\Gamma_{\ell_k^{1j}})} < \infty,
\end{align*}
where $c_0$ is the number appearing in Proposition \ref{prop:lap}. 
Thus, we have $\tilde \zeta_jf \in W^s_q(\Gamma_j)$.  
Since $m_{1j} \leq L$ with some constant $L$ independent of 
$M_1$ as follows from Proposition \ref{prop:lap}, we have  
\begin{align*}
\sum_{j=1}^\infty\|\tilde \zeta_jf\|_{W^s_q(\Gamma_j)}^q
&\leq  \sum_{j=1}^\infty (c_0C_{s,q})^q (m^1_j)^{q-1} \sum_{k=1, \ldots, m^1_j}
\|f_{\ell^{1j}_k}\|_{W^s_q(\Gamma_{\ell^{1j}_k})}^q \\
&\leq (c_0C_{s,q}L)^q\sum_{n=1}^\infty\|f_n\|_{W^s_q(\Gamma_n)}^q,
\end{align*}
which yields \eqref{def:sob2}

\subsection{The weak Dirichlet problem}\label{subsec:p-2}

Let $\hat H^1_{q,0}(\Omega)$ be the homogeneous Sobolev space defined by
letting
\begin{equation}\label{hom:1}
\hat H^1_{q,0}(\Omega) = \{\varphi \in L_{q, {\rm loc}}(\Omega)
\mid \nabla\varphi \in L_q(\Omega)^N, \quad\varphi|_\Gamma=0\}.
\end{equation}
Let $1 < q < \infty$. The variational equation:
\begin{equation}\label{wd:1}
(\nabla u, \nabla\varphi)_\Omega = (\bff, \nabla\varphi)_\Omega
\quad\text{for all $\varphi \in \hat H^1_{q', 0}(\Omega)$}
\end{equation}
is called the weak Dirichlet problem, where $q' = q/(q-1)$.  
\begin{dfn}\label{dfn:wd} We say that the weak Dirichlet problem
\eqref{wd:1} is uniquely solvable in $\hat H^1_{q,0}(\Omega)$ if 
for any $\bff \in L_q(\Omega)^N$, problem \eqref{wd:1} admits
a unique solution $u \in \hat H^1_{q,0}(\Omega)$ possessing the 
estimate: $\|\nabla u\|_{L_q(\Omega)} \leq C\|\bff\|_{L_q(\Omega)}$.

We define an operator $\CK$ acting on $\bff \in L_q(\Omega)^N$ by 
$u = \CK(\bff)$.
\end{dfn} 
\begin{rem}\label{rem:3.2} \thetag1\quad
When $q=2$, the weak Dirichlet problem is uniquely 
solvable for any  domain $\Omega$, which is easily proved by using
the Hilbert space structure of the space $\hat H^1_{2,0}(\Omega)$.
But, for any $q \not=2$, speaking generally, 
we do not know whether the weak Dirichlet problem 
is uniquely solvable. \\
\thetag2\quad Given $\bff \in L_q(\Omega)^N$ and $g \in 
W^{1-1/q}_q(\Gamma)$, we consider the weak Dirichlet problem:
\begin{equation}\label{wd:2}
(\nabla u, \nabla \varphi)_\Omega = (\bff, \nabla\varphi)_\Omega
\quad\text{for every $\varphi \in \hat H^1_{q', 0}(\Omega)$},
\end{equation}
subject to $u = g$ on $\Gamma$.  Let $G$ be an extension of $g$ to
$\Omega$ such that $G = g$ on $\Gamma$ and  
$\|G\|_{H^1_q(\Omega)} \leq C\|g\|_{W^{1-1/q}_q(\Gamma)}$ for some
constant $C > 0$.  Let $v \in \hat H^1_{q,0}(\Omega)$
be a solution of the weak Dirichlet problem:
$$
(\nabla v, \nabla \varphi)_\Omega = (\bff-\nabla G, \nabla\varphi)_\Omega
\quad\text{for every $\varphi \in \hat H^1_{q', 0}(\Omega)$}.
$$
Then, $u = G + v \in H^1_q(\Omega) + \hat H^1_{q,0}(\Omega)$ is a unique
solution of Eq. \eqref{wd:2} possessing the estimate: 
$$\|\nabla u\|_{L_q(\Omega)}
\leq C(\|g\|_{W^{1-1/q}_q(\Gamma)} + \|\bff\|_{L_q(\Omega)})$$
for some constant $C > 0$.
\end{rem}
\begin{ex}
When $\Omega$ is a  bounded domain, an exterior domain, half space, 
a perturbed half space, layer, a perturbed layer, and a 
tube, then the weak Dirichlet problem is uniquely solvable for $q
\in (1, \infty)$.
\end{ex}
\begin{thm}\label{thm:ap-r} Let $1 < q < \infty$. Let $\Omega$ be
a uniform $C^2$ domain.  Given $\bff \in L_q(\Omega)^N$, let
$u \in \hat H^1_{q,0}(\Omega)$ be a unique solution of the weak Dirichlet
problem \eqref{wd:1} possessing the estimate:
$\|\nabla u \|_{L_q(\Omega)} \leq C\|\bff\|_{L_q(\Omega)}$.  If we assume
that $\dv\bff \in L_q(\Omega)$ in addition, then $\nabla^2 u \in L_q(\Omega)$
and 
$$\|\nabla^2u\|_{L_q(\Omega)} \leq C_{M_2, q}(\|\dv\bff\|_{L_q(\Omega)}
+ \|\bff\|_{L_q(\Omega)}).
$$
\end{thm}
\pf This theorem will be proved in Subsec. \ref{sec:ap.2}, after
some preparations for weak Laplace problem in $\BR^N$ in 
Subsec.\ref{subsec.ap.1} and weak Dirihle problem in the half space
in Subsec. \ref{subsec:ap.2} below.  \qed
\begin{rem}
From Theorem \ref{thm:ap-r}, we know the existence of solutions of the 
strong Dirichlet problem:
\begin{equation}\label{sd:1}
\Delta u = \dv\bff\quad\text{in $\Omega$}, \quad 
u|_\Gamma = 0.
\end{equation}
But, the uniqueness of solutions of Eq.\eqref{sd:1} does not hold generally.
For example, let $B_1= \{x \in \BR^N \mid |x| >1\}$ for $N \geq 2$,
and let $f(x)$ be a function defined by
$$f(x) = \begin{cases} \log|x| \quad&\text{when $N=2$}, \\
|x|^{-(N-2)}-1 \quad&\text{when $N \geq 3$}.
\end{cases}
$$
Then, $f(x)$ satisfies the homogeneous Dirichlet problem:
$\Delta f=0$ in $B_1$ and $f|_{\pd B_1} = 0$, where $\pd B_1
= \{x \in \BR^N \mid |x| = 1\}$. 
\end{rem}


\subsection{The weak Laplace problem in $\BR^N$}\label{subsec.ap.1}
In this subsection, we consider the following
weak Laplace problem in $\BR^N$:
\begin{equation}\label{ap:wd.1}
(\nabla u, \nabla\varphi)_{\BR^N} = (\bff, \nabla\varphi)_{\BR^N}
\quad\text{for any $\varphi \in \hat H^1_{q'}(\BR^N)$}.
\end{equation}
We shall prove the following theorem.
\begin{thm}\label{thm.wd.1}
Let $1 < q < \infty$.  Then, for any $\bff \in L_q(\BR^N)^N$, 
the weak Laplace problem \eqref{ap:wd.1} admits a unique
solution $u \in \hat H^1_q(\BR^N)$ possessing the estimate:
$\|\nabla u\|_{L_q(\BR^N)} \leq C\|\bff\|_{L_q(\BR^N)}$. 

Moreover, if we assume that $\dv\bff \in L_q(\BR^N)$ in addition,
then $\nabla^2 u \in L_q(\BR^N)^{N^2}$ and 
$$\|\nabla^2u\|_{L_q(\BR^N)} \leq C\|\dv\bff\|_{L_q(\BR^N)}.
$$
\end{thm}
\pf To prove the theorem, we consider the strong Laplace equation:
\begin{equation}\label{ap:1}
\Delta u = \dv \bff \quad\text{in $\BR^N$}.
\end{equation}
Let  
$$H^1_{q, {\rm div}}(D) = \{\bff \in L_q(D)^N \mid \dv \bff \in L_q(D)\},
$$
where $D$ is any domain in $\BR^N$.  Since $C^\infty_0(\BR^N)^N$ is dense
both in $L_q(\BR^N)^N$ and $H^1_{q, {\rm div}}(\BR^N)$, 
we may assume that $\bff \in C^\infty_0(\BR^N)^N$. 
Let $\CF[f] = \hat f$ and $\CF^{-1}$ denote the Fourier transfom 
$f$ and the Fourier inverse transform, respectively. We then set
$$u = -\CF^{-1}\Bigl[\frac{\CF[\dv\bff](\xi)}{|\xi^2|}\Bigr]
= -\CF^{-1}\Bigl[\frac{\sum_{j=1}^Ni\xi_j\CF[f_j](\xi)}{|\xi|^2}\Bigr]
$$
for $\bff = {}^\top(f_1, \ldots, f_N)$. 
By  the Fourier multiplier theorem, we have
\begin{equation}\label{ap:2}\begin{split}
\|\nabla u\|_{L_q(\BR^N)} &\leq C\|\bff\|_{L_q(\BR^N)}, \\
\|\nabla^2u\|_{L_q(\BR^N)}&\leq C\|\dv \bff\|_{L_q(\BR^N)}.
\end{split}\end{equation}
Of course, $u$ satisfies Eq.\eqref{ap:1}. 

We now prove that $u$ satisfies the weak Laplace equation
\eqref{ap:1}.  For this purpose, we use the following lemma.
\begin{lem}\label{lem:hardy} Let $1 < q < \infty$ and let 
$$d_q(x) = \begin{cases} (1 + |x|^2)^{1/2} &\quad\text{for $N\not=q$}, \\
(1+|x|^2)^{1/2}\log(2 + |x|^2)^{1/2} &\quad\text{for $N=q$}.
\end{cases}
$$
Then, for any $\varphi \in \hat H^1_q(\BR^N)$, there exists a constant
$c$ for which
$$
\Bigl\|\frac{\varphi-c}{d_q}\Bigr\|_{L_q(\BR^N)} \leq C\|\nabla\varphi
\|_{L_q(\BR^N)}
$$
with some constant independent of $\varphi$ and $c$.
\end{lem}
\pf For a proof,  see Galdi \cite[Chapter II]{Gal}. \qed
\vskip0.5pc
To use Lemma \ref{lem:hardy}, we use a cut-off function, 
$\psi_R$, of Sobolev's
type defined as follows: Let $\psi$ be a function in $C^\infty(\BR)$ such that
$\psi(t) = 1$ for $|t| \leq 1/2$ and $\psi(t) = 0$ for $|t| \geq 1$, and 
set 
$$\psi_R(x) = \psi\Bigl(\frac{\ln\ln|x|}{\ln\ln R}\Bigr).
$$
Notice that 
\begin{equation} \label{ap:5} |\nabla\psi_R(x)| \leq \frac{c}{\ln \ln R}
\frac{1}{|x|\ln|x|},
\quad{\rm supp}\,\nabla\psi_R \subset D_R,
\end{equation}
where we have set $D_R = \{x \in \BR^N \mid
e^{\sqrt{\ln R}} \leq |x| \leq R\}$. Noting that 
$\bff \in C^\infty_0(\BR^N)^N$,
by \eqref{ap:1} for large $R > 0$ and $\varphi \in \hat H^1_{q'}(\BR^N)$
we have 
\begin{equation}\label{ap:6}\begin{split}
(\bff, \nabla\varphi)_{\BR^N} &= (\bff, \nabla(\varphi-c))_{\BR^N}
= -(\dv \bff, \varphi -c)_{\BR^N}\\
&= -(\psi_R\dv\bff, \varphi-c)_{\BR^N}
= -(\psi_R\Delta u, \varphi-c)_{\BR^N} \\
&=((\nabla\psi_R)\cdot(\nabla u), \varphi-c)_{\BR^N}
+ (\psi_R\nabla u, \nabla\varphi)_{\BR^N},
\end{split}\end{equation}
where $c$ is a constant for which 
\begin{equation}\label{ap:4}
\Bigl\|\frac{\varphi-c}{d_{q'}}\Bigr\|_{L_{q'}(\BR^N)} 
\leq C\|\nabla\varphi\|_{L_{q'}(\BR^N)}.
\end{equation}
By \eqref{ap:5} and \eqref{ap:4}, we have
\begin{equation}\label{ap:7}\begin{split}
|((\nabla\psi_R)\cdot(\nabla u), \varphi-c)_{\BR^N}| 
&\leq \|d_{q'}(\nabla\psi_R)\cdot(\nabla u)\|_{L_q(\BR^N)}
\Bigl\|\frac{\varphi-c}{d_{q'}}\Bigr\|_{L_{q'}(\BR^N)} \\
&\leq \frac{C}{\ln \ln R}\|\nabla u\|_{L_q(D_R)}
\|\nabla \varphi\|_{L_{q'}(\BR^N)} \to 0
\end{split}\end{equation}
as $R \to \infty $. By \eqref{ap:6} and \eqref{ap:7} 
we see that $u$ satisfies the weak Dirichlet problem 
\eqref{ap:wd.1}.  The uniqueness follows from the 
existence theorem just proved for the dual problem. 
Moreover, if $\dv \bff \in L_q(\BR^N)$ in addition, then
$\nabla^2u \in L_q(\BR^N)^{N^2}$, and so by \eqref{ap:2}
we complete the proof of Theorem \ref{thm.wd.1}. \qed
\subsection{The weak Dirichlet problem in
the half space case}\label{subsec:ap.2}
In this subsection, we consider the following
weak Dirichlet problem in $\BR^N_+$: 
\begin{equation}\label{ab:wd.1}
(\nabla u, \nabla\varphi)_{\BR^N_+} = (\bff, \nabla\varphi)_{\BR^N_+}
\quad\text{for any $\varphi \in \hat H^1_{q',0}(\BR^N_+)$}.
\end{equation}
We shall prove the following theorem.
\begin{thm}\label{thm.wd.2}
Let $1 < q < \infty$.  Then, for any $\bff \in L_q(\BR^N_+)^N$, 
the weak Dirichlet problem \eqref{ab:wd.1} admits a unique
solution $u \in \hat H^1_{q,0}(\BR^N_+)$ possessing the estimate:
$\|\nabla u\|_{L_q(\BR^N_+)} \leq C\|\bff\|_{L_q(\BR^N_+)}$. 

Moreover, if we assume that $\dv\bff \in L_q(\BR^N_+)$ in addition,
then $\nabla^2 u \in L_q(\BR^N_+)^{N^2}$ and 
$$\|\nabla^2u\|_{L_q(\BR^N_+)} \leq C\|\dv\bff\|_{L_q(\BR^N_+)}.
$$
\end{thm}
\pf We may assume that 
$\bff ={}^\top(f_1, \ldots, f_N) \in C^\infty_0(\HS)^N$
in the following, because $C^\infty_0(\HS)$ is dense both in $L_q(\HS)^N$ and
$H^1_{q, {\rm div}}(\HS)$. We first consider the strong Dirichlet problem:
\begin{equation}\label{ab:2}
\Delta u = \dv\bff \quad\text{in $\HS$}, \quad 
u|_{x_N=0} = 0.
\end{equation}
For any  function, $f(x)$, defined in $\HS$, let 
$f^e$ and $f^o$ be the even extension and the odd extension of
$f$ defined by letting  
\begin{equation}\label{4.ext.1}
f^e(x) = \begin{cases} f(x', x_N) \quad&x_N > 0, 
\\ f(x', -x_N) \quad&x_N < 0, \end{cases} \quad
f^o(x) = 
\begin{cases} f(x', x_N) \quad&x_N > 0, 
\\ -f(x', -x_N) \quad&x_N < 0, \end{cases}
\end{equation}
where $x' = (x_1, \ldots, x_{N-1}) \in \BR^{N-1}$ and $x = (x', x_N) 
\in\BR^N$.

 Noting that 
$(\dv\bff)^o = \sum_{j=1}^{N-1}\pd_j(f_j)^o + \pd_N(f_N)^e$, 
we define $u$ by letting
\begin{align*}
u &= -\CF^{-1}\Bigl[\frac{\CF[(\dv\bff)^o](\xi)}{|\xi|^2}\Bigr]
\\
&= -\CF^{-1}\Bigl[\frac{\sum_{j-1}^{N-1}i\xi_j\CF[(f_j)^o](\xi)
+ i\xi_N\CF[(f_N)^e](\xi)}{|\xi|^2}\Bigr].
\end{align*}
We then have
\begin{equation}\label{ab:3}\begin{split}
\|\nabla u\|_{L_q(\BR^N)} & \leq C\|\bff\|{L_q(\HS)}, \quad 
\|\nabla^2u\|_{L_q(\BR^N)}  \leq C\|\dv\bff\|_{L_q(\HS)},
\end{split}\end{equation}
and moreover $u$ satisfies Eq.\eqref{ab:2}. 

We next prove that $u$ satisfies the weak Dirichlet problem Eq.\eqref{ab:wd.1}.
For this purpose, instead of lemma \ref{lem:hardy}, we use the 
Hardy type inequality:
\begin{equation}\label{ab:4}
\Bigl(\int^\infty_0\Bigl(\int^x_0f(y)\,dy\Bigr)^px^{-r-1}
\,dy\Bigr)^{1/p}
\leq (p/r)\Bigl(\int^\infty_0(yf(y))^p
y^{-r-1}\,dy\Bigr)^{1/p},
\end{equation}
where $f \geq 0$, $p \geq 1$ and $r > 0$ (cf. Stein \cite[A.4 p.272]{Stein}).  
Of course, using zero extension of
f suitably, we can replace the interval $(0, \infty)$ by $(a, b)$ for any
$0\leq a < b < \infty$ in \eqref{ab:4}. 
Let $D_{R, 2R} = \{x \in \BR^N \mid R \leq |x| \leq 2R\}$.
Using \eqref{ab:4}, 
we see that for any $\varphi \in \hat H^1_{q',0}(\HS)$ 
\begin{equation}\label{ab:5} 
\lim_{R\to\infty} R^{-1}\|\varphi\|_{L_{q'}(D_{R, 2R})} =0.
\end{equation}
In fact, using $\varphi|_{x_N=0}=0$, we write 
$\displaystyle{\varphi(x', x_N) = \int^{x_N}_0 (\pd_s\varphi)(x', s)\,ds}$. 
Thus, by \eqref{ab:4} we have
$$\int^b_a|\varphi(x', x_N)|^{q'}\,dx_N \leq \Bigl(\frac{bq'}{q'-1}\Bigr)^{q'}
\int^b_a|(\pd_N\varphi)(x', x_N)|^{q'}\,dx_N
$$
for any $0 < a < b$.  Let 
\begin{align*}
E^1_R & =\{x \in \BR^N \mid |x'| \leq 2R, \enskip 
R/2 \leq x_N < 2R\}, \\
E^2_R & =\{x \in \BR^N \mid 0 \leq x_N  \leq 2R, \enskip 
R/2 \leq |x'|\leq 2R\},
\end{align*}
and then $D_{R,2R} \subset E^1_R \cup E^2_R$. 
Thus, by \eqref{ab:4}, 
\begin{align*}
\Bigl(\int_{D_{R, 2R}}|\varphi(x)|^{q'}\,dx\Bigr)^{1/{q'}} 
\leq &\Bigl(\frac{Rq'}{q'-1}\Bigr)\Bigr\{\underset{|x'|\leq R}
\int
\int^{2R}_{R/2}|\pd_N\varphi(x)|^{q'}\,dx_Ndx'\\
& + 
\underset{R/2
\leq |x'| \leq 2R}
\int \int^{2R}_0|\pd_N\varphi(x)|^{q'}\,dx_Ndx'
\Bigr\}^{1/{q'}},
\end{align*}
which leads to \eqref{ab:5}.

Let $\omega$ be a function in $C^\infty_0(\BR^N)$ such that
$\omega(x) = 1$ for $|x| \leq 1$ and $\varphi(x) = 0$ for $|x| \geq 2$,
and we set $\omega_R(x) = \omega(x/R)$. 
For any $\varphi \in \hat H^1_{q', 0}(\HS)$ and for large $R > 0$, we have
\begin{equation}\label{ab:6}\begin{split}
(\dv\bff, \varphi)_{\HS} & = (\omega_R\dv\bff, \varphi)_{\HS} 
= (\omega_R\Delta u, \varphi)_{\HS} \\
&= 
-((\nabla\omega_R)\cdot(\nabla u), \varphi)_{\HS}
-(\omega_R\nabla u, \nabla\varphi)_{\HS}.
\end{split}\end{equation}
By \eqref{ab:5}
$$|((\nabla\omega_R)\cdot(\nabla u), \varphi)_{\HS}|
\leq R^{-1}\|\nabla u\|_{L_q(D_{R,2R})}\|\varphi\|_{L_{q'}(D_{R,2R})}\to 0
$$
as $R \to \infty$.  On the other hand, $(\dv\bff, \varphi)_{\HS} = 
-(\bff, \nabla\varphi)_{\HS}$,
where we have used $\bff \in C^\infty_0(\HS)^N$. Thus, by \eqref{ab:6}
we have
$$(\nabla u, \nabla\varphi)_{\HS} = (\bff, \nabla\varphi)_{\HS}$$
for any $\varphi \in \hat H^1_{q', 0}(\HS)$.  This shows that $u$
is a solution of the weak Dirichlet problem.  The uniqueness follows
from the existence of solutions for the dual problem,
which completes the proof of Theorem \ref{thm.wd.2}.\qed
\subsection{Regularity of the weak Dirichlet problem}\label{sec:ap.2}

In this subsection, we shall prove  Theorem \ref{thm:ap-r} in Subsec.
\ref{subsec:p-2}.
Let $\zeta^i_j$ $(i=0,1, \, j \in \BN)$ be cut-off functions
given in Proposition \ref{prop:lap}. We first consider the regularity
of $\zeta^0_ju$. For this purpose, 
we use the following lemma. 
\begin{lem}\label{lem:5.7}
Let $\Omega$ be a uniformly $C^2$ domain in $\BR^N$. Then,
there exists a constant $c_1 > 0$ independent of $j \in \BN$
such that 
\begin{alignat*}2
\|\varphi\|_{H^1_q(\Omega_j\cap B^1_j)} 
&\leq c_1\|\nabla \varphi\|_{L_q(\Omega_j\cap B^1_j)}
&\quad&\text{for any $\varphi \in \hat H^1_{q,0}(\Omega_j)$}, \\
\|\psi\|_{H^1_q(\Omega\cap B^1_j)} 
&\leq c_1\|\nabla \psi\|_{L_q(\Omega\cap B^1_j)}
&\quad&\text{for any $\psi \in \hat H^1_{q, 0}(\Omega)$}, \\
\|\varphi-c^0_j(\varphi)\|_{H^1_q(B^0_j)} 
&\leq c_1\|\nabla \varphi\|_{L_q(B^0_j)}
&\quad&\text{for any $\varphi \in \hat H^1_q(\BR^N)$}, \\
\|\psi-c^0_j(\psi)\|_{H^1_q(B^0_j)} 
&\leq c_1\|\nabla \psi\|_{L_q(B^0_j)}
&\quad&\text{for any $\psi \in \hat H^1_{q, 0}(\Omega)$}. 
\end{alignat*}
Here, $c^0_j(\varphi)$ and $c^0_j(\psi)$ are suitable 
constants depending 
on $\varphi$ and $\psi$, respectively. 
\end{lem}
\pf For a proof, see  Shibata \cite[Lemma 3.4, Lemma  3.5]{S0}. \qed
\vskip0.5pc\noindent
{\bf Continuation of Proof of Theorem \ref{thm:ap-r}}.~ 
 Let $c^0_j=c^0_j(\varphi)$ be a constant 
in Lemma \ref{lem:5.7} such that 
\begin{equation}\label{ad:1} \|u-c^0_j\|_{L_q(B^0_j)} 
\leq c_1\|\nabla u\|_{L_q(B^0_j)}.
\end{equation}
For any $\varphi \in \hat H^1_{q'}(\BR^N)$, we have
\begin{align*}
&(\nabla(\zeta^0_j(u-c^0_j)), \nabla\varphi)_{\BR^N} \\
& = ((\nabla\zeta^0_j)(u-c^0_j), \nabla\varphi)_{\BR^N}
+ (\nabla u, \nabla(\zeta^0_j\varphi))_{\BR^N} -
((\nabla u)\cdot(\nabla\zeta^0_j), \varphi)_{\BR^N} \\
&= -((\Delta \zeta^0_j)(u-c^0_j) + 2(\nabla\zeta^0_j)
\cdot(\nabla u)
+ \zeta^0_j\dv\bff, \varphi)_{\BR^N},
\end{align*}
where we have used 
$(\nabla u, \nabla(\zeta^0_j\varphi))_{\BR^N}
 = (\bff, \nabla(\zeta^0_j\varphi))_{\BR^N}
=-(\zeta^0_j\dv\bff, \varphi)_{\BR^N}.$
Let $f = (\Delta \zeta^0_j)(u-c^0_j) + 2(\nabla\zeta^0_j)
\cdot(\nabla u)
+ \zeta^0_j\dv\bff$.  Since $C^\infty_0(\BR^N) \subset \hat H^1_q(\BR^N)$,
for any $\varphi \in C^\infty_0(\BR^N)$ we have
$$(\Delta(\zeta^0_j(u-c^0_j)), \varphi)_{\BR^N} = (f, \varphi)_{\BR^N},
$$
which yields that 
\begin{equation}\label{ad:1.1}
\Delta(\zeta^0_j(u-c^0_j)) = f \quad\text{in $\BR^N$}
\end{equation}
in the sense of distribution. By Lemma \ref{lem:5.7}, $f \in L_q(\BR^N)$ and 
\begin{equation}\label{ad:2}
\|f\|_{L_q(\BR^N)} \leq C(\|\dv\bff\|_{L_q(B^0_j)} 
+ \|\nabla u\|_{L_q(B^0_j)}). 
\end{equation}
From \eqref{ad:1.1} it follows that
$$\pd_k\pd_\ell\Delta(\zeta^0_j(u-c^0_j)) = \pd_k\pd_\ell f$$
for any $k$, $\ell \in \BN$.  Since both sides are compactly
supported distributions, we can apply the Fourier transform and
the inverse Fourier transform. We then have
$$\pd_k\pd_\ell(\zeta^0_j(u-c^0_j))
= \CF^{-1}\Bigl[\frac{\xi_k\xi_\ell}{|\xi|^2}
\CF[f](\xi)\Bigr].
$$
By the Fourier multiplier theorem, we have
$$\|\pd_k\pd_\ell(\zeta^0_j(u-c^0_j))\|_{L_q(\BR^N)}
\leq C\|f\|_{L_q(\BR^N)}.
$$
Since $\pd_k\pd_\ell(\zeta^0_j(u-c^0_j))
= \zeta^0_j\pd_k\pd_\ell u + (\pd_k\zeta^0_j)\pd_\ell u
+(\pd_\ell\zeta^0_j)\pd_ku + (\pd_k\pd_\ell\zeta^0_j)(u-c^0_j)$,
by \eqref{ad:1} and \eqref{ad:2} we have
$\zeta^0_j\nabla^2u \in L_q(\BR^N)^N$ and
\begin{equation}\label{ad:3}\|\zeta^0_j\nabla^2u
\|_{L_q(\Omega)} \leq C_{M_2}(\|\dv\bff\|_{L_q(B^0_j)}
+ \|\nabla u\|_{L_q(B^0_j)}).
\end{equation}

We next consider $\zeta^1_ju$.  
For any $\varphi \in \hat H^1_{q',0}(\Omega_j)$, we have
\begin{equation}\label{ad:4}
(\nabla(\zeta^1_ju), \nabla\varphi)_{\Omega_J} = (g, \varphi)_{\Omega_j},
\end{equation}
where we have set 
$g = -(\zeta^1_j\dv\bff + 2(\nabla u)\cdot(\nabla\zeta^1_j)
+ (\Delta\zeta^1_j)u)$. By Lemma \ref{lem:5.7}, $g \in L_q(\Omega_j)$ and 
\begin{equation}\label{ad:5}
\|g\|_{L_q(\Omega_j)} \leq C(\|\dv\bff\|_{L_q(\Omega\cap B^1_j)}
+ \|\nabla u\|_{L_q(\Omega \cap B^1_j)}).
\end{equation}
We use the symbols given in Proposition \ref{prop:lap}.  
Let $a_{k\ell}$ and $b_{k\ell}$ be the $(k, \ell)^{\rm th}$
component of  $N\times N$ matrices $\CA_j$ and $B_j$ given in
Proposition \ref{prop:lap}.  By the change of variables:
$y = \Phi_j(x)$, the variational equation \eqref{ad:4}
is transformed to
\begin{equation}\label{ad:6}
\sum_{k,\ell=1}^N((\delta_{k\ell} + A_{k\ell})\pd_kv, \pd_\ell\varphi)_{\HS}
= (h, \varphi)_{\HS}.
\end{equation}
Where, we have set 
\begin{align*}
v & = \zeta^1_ju\circ\Phi_j, \quad h = g\circ\Phi_j, \quad
J = \det(\CA_j + B_j) = 1 + J^0, \\
A_{k\ell} & = \sum_{m=1}^N\{a_{\ell m}b_{km} 
+ a_{\ell m}J_0(a_{k\ell}+b_{k\ell})
+ b_{\ell m}J(a_{km} + b_{km})\}.
\end{align*}
By Proposition \ref{prop:lap} and \eqref{ad:5}, we have
\begin{equation}\label{ad:7}\begin{split}
\|A_{k\ell}\|_{L_\infty(\BR^N)} & \leq CM_1, \quad 
\|\nabla A_{k\ell}\|_{L_\infty(\BR^N)} \leq C_A, \\
\|h\|_{L_q(\HS)} & \leq C(\|\dv\bff\|_{L_q(B^1_j\cap\Omega)} + 
\|\nabla u\|_{L_q(B^1_j\cap\Omega)}),
\end{split}\end{equation}
where $C_A$ is a constant depending $a_1$, $a_2$ and $A$ appearing
in Definition \ref{dfn:1}. Since ${\rm supp}\,g \subset
\Phi^{-1}(B^1_j) \cap \Omega$, by Lemma \ref{lem:5.7} 
$$|(h, \varphi)_{\HS}| \leq \|h\|_{L_q(B^1_j \cap \HS)}
\|\varphi\|_{L_{q'}(B^1_j\cap \HS)} 
\leq C\|g\|_{L_q(B^1_j\cap \Omega)} \|\nabla\varphi\|_{L_{q'}(\HS)}
$$
for any $\varphi \in \hat H^1_{q',0}(\HS)$, where $C$ is a constant independent
of $j \in \BN$. Thus, by the Hahn-Banach theorem, there exists a $\bh 
\in L_q(\HS)^N$ such that $\|\bh\|_{L_q(\HS)}
 \leq C\|h\|_{L_q(B^1_j\cap\HS)}$
and $(\bh, \nabla\varphi)_\HS = (h, \varphi)_\HS$ for any $\varphi \in
\hat H^1_{q', 0}(\HS)$.  In particular, $\dv \bh = -h \in L_q(\HS)$.
Thus, the variational problem \eqref{ad:6} reads
$$\sum_{k,\ell=1}^N((\delta_{k\ell} + A_{k\ell})\pd_kv, 
\pd_\ell\varphi)_\HS = (\bh, \nabla\varphi)_\HS
\quad\text{for any $\varphi \in \hat H^1_{q', 0}(\HS)$}.
$$

We now prove that if $M_1 \in (0, 1)$ is small enough,
then for any $\bg \in L_q(\HS)^N$, 
there exists a  unique solution $w \in 
\hat H^1_{q,0}(\HS)$ of the variational problem:
\begin{equation}\label{ac:2.5}
\sum_{k,\ell=1}^N((\delta_{k\ell} + A_{k\ell})\pd_kw, 
\pd_\ell\varphi)_\HS = (\bg, \nabla\varphi)_\HS
\quad\text{for any $\varphi \in \hat H^1_{q', 0}(\HS)$},
\end{equation}
possessing the estimate: 
\begin{equation}\label{ac:2.4}
\|\nabla w\|_{L_q(\HS)} \leq C\|\bg\|_{L_q(\HS)}.
\end{equation} 
Morevoer, if $\dv \bg \in L_q(\HS)$, then $\nabla w \in H^1_q(\HS)^N$ and 
\begin{equation}\label{ac:2.3}
\|\nabla^2w\|_{L_q(\HS)} \leq C\|\dv\bg\|_{L_q(\HS)} + 
C_A\|\bg\|_{L_q(\HS)}.
\end{equation}

In fact, we prove the existence of $w$ by the successive approximation.
Let $w_1 \in \hat H^1_{q,0}(\HS)$ be a solution of the weak 
Dirichlet problem:
\begin{equation}\label{ac:2.6}
(\nabla w_1, \nabla\varphi)_{\HS} = (\bg, \nabla\varphi)_{\HS}
\quad\text{for any $\varphi \in \hat H^1_{q', 0}(\HS)$}.
\end{equation}
By Theorem \ref{thm.wd.2}, $w_1$ uniquely exists and satisfies
the estimate:
\begin{equation}\label{ac:2.7}
\|\nabla w_1\|_{L_q(\HS)} \leq C\|\bg\|_{L_q(\HS)}.
\end{equation}
Moreover, if we assume that $\dv\bg  \in L_q(\HS)$ additionally,
then $\nabla^2 w_1 \in H^1_q(\HS)^N$ and 
\begin{equation}\label{ac:2.8}
\|\nabla^2w_1\|_{L_q(\HS)} \leq C\|\dv \bg\|_{L_q(\HS)}.
\end{equation}
Given $w_j \in \hat H^1_{q, 0}(\HS)$, let $w_{j+1} \in \hat H^1_{q,0}(\HS)$
be a solution of the weak Dirichlet problem:
\begin{equation}\label{ac.2.9}
(\nabla w_{j+1}, \nabla\varphi)_{\HS} 
= (\bg, \nabla\varphi)_{\HS} - \sum_{k, \ell=1}^N(A_{k\ell}\pd_kw_j,
\pd_\ell \varphi)_{\HS}
\end{equation}
for any $\varphi \in \hat H^1_{q', 0}(\HS)$.
By Theorem \ref{thm.wd.2} and \eqref{ad:7}, $w_{j+1}$ exists and 
satisfies the estimate:
\begin{equation}\label{ac:2.10}
\|\nabla w_{j+1}\|_{L_q(\HS)} \leq C(\|\bg\|_{L_q(\Omega_+)} 
+ M_1\|\nabla w_j\|_{L_q(\HS)}).
\end{equation}
Applying Theorem \ref{thm.wd.2} and \eqref{ad:7} to the 
difference $w_{j+1}-w_j$, we have
\begin{equation}\label{ac:2.11}
\|\nabla (w_{j+1} - w_j)\|_{L_q(\HS)} 
\leq CM_1\|\nabla (w_j-w_{j-1})\|_{L_q(\HS)}.
\end{equation}
Choosing $CM_1 \leq 1/2$ in \eqref{ac:2.11}, 
we see  that $\{w_j\}_{j=1}^\infty$ is a Cauchy sequence in 
$\hat H^1_{q,0}(\HS)$, and so the limit $w \in H^1_{q, 0}(\HS)$
exists and 
satisfies the weak Dirichlet problem \eqref{ac:2.5}.  Moreover,
taking the limit in \eqref{ac:2.10}, we have
$$\|\nabla w\|_{L_q(\HS)} \leq C\|\bg\|_{L_q(\Omega_+)} 
+ CM_1\|\nabla w\|_{L_q(\HS)}.$$
Since $CM_1 \leq 1/2$, we have 
$\|\nabla w\|_{L_q(\HS)} \leq 2C\|\bg\|_{L_q(\Omega_+)}$. Thus, 
we have proved that the weak Dirichlet problem \eqref{ac:2.5}
admits at least one solution $w \in \hat H^1_{q,0}(\Omega_+)$
possessing the estimate \eqref{ac:2.4}. 
 The uniqueness follows from
the existence of solutions to the dual problem.  Thus, 
we have proved the unique existence of solutions of
Eq.\eqref{ac:2.5}. 

We now prove that $\nabla \bw \in H^1_q(\HS)^N$ provided that
$\dv \bg \in L_q(\HS)$. By Theorem \ref{thm.wd.2}, 
$\nabla w_1 \in H^1_q(\HS)^N$ and $w_1$ satisfies the estimate:
$$\|\nabla^2w_1\|_{L_q(\HS)} \leq C\|\dv\bg\|_{H^1_q(\Omega_+)}.
$$
Moreover, if we assume that $\nabla^2 w_j \in L_q(\HS)^N$
in addition, then applying Theorem \ref{thm.wd.2} to
\eqref{ac.2.9} and using \eqref{ad:7}, we see that
$\nabla^2w_{j+1} \in L_q(\HS)^{N^2}$ and 
\begin{equation}\label{ac:2.12}
\|\nabla^2w_{j+1}\|_{L_q(\HS)} \leq C\|\dv\bg\|_{L_q(\HS)}
+ CM_1\|\nabla^2 w_j\|_{L_q(\HS)} + C_A\|\nabla w_j\|_{L_q(\HS)}.
\end{equation}
And also, applying Theorem \ref{thm.wd.2} to the difference
$w_{j+1}-w_j$ and using \eqref{ad:7}, we have
\begin{align*}
&\|\nabla^2(w_{j+1}-w_j)\|_{L_q(\HS)}\\
&\quad \leq 
CM_1\|\nabla^2(w_j-w_{j-1})\|_{L_q(\HS)} 
+ C_A\|\nabla(w_j-w_{j-1})\|_{L_q(\HS)},
\end{align*}
which, combined with \eqref{ac:2.11}, leads to
\begin{align*}
&\|\nabla^2(w_{j+1}-w_j)\|_{L_q(\HS)} 
+ \|\nabla(w_j- w_{j-1})\|_{L_q(\HS)} \\
&\quad \leq 
CM_1\|\nabla^2(w_j-w_{j-1})\|_{L_q(\HS)} 
+ C(C_A+1)M_1\|\nabla(w_{j-1}-w_{j-2})\|_{L_q(\HS)}.
\end{align*}
Choosing $M_1 > 0$ so small that $CM_1 \leq 1/2$ and
$(C_A+1)M_1 \leq 1/2$, then we have 
$$\|\nabla^2(w_{j+1}-w_j)\|_{L_q(\HS)} 
+ \|\nabla(w_j- w_{j-1})\|_{L_q(\HS)}
\leq (1/2)^{j-1}L$$
with $L = \|\nabla^2(w_3-w_2)\|_{L_q(\HS)}
 + \|\nabla(w_2-w_1)\|_{L_q(\HS)}$.
 From this it follows that $\{\nabla^2 w_j\}_{j=1}^\infty$ is 
a Cauchy sequence in $L_q(\Omega)$, which yields that
$\nabla^2 w\in L_q(\HS)^{N^2}$.  Moreover, taking the limit
in \eqref{ac:2.12} and using \eqref{ac:2.4} gives that
$$\|\nabla^2 w\|_{L_q(\HS)} \leq C\|\dv\bg\|_{L_q(\HS)} 
+ (1/2)\|\nabla^2w\|_{L_q(\HS)} + C_A\|\bg\|_{L_q(\HS)}.
$$
which leads to \eqref{ac:2.3}. 

Applying what we have proved and using the estimate:
$$\|\dv \bh\|_{L_q(\HS)} + \|\bh\|_{L_q(\HS)} \leq C
\|h\|_{L_q(\HS)} \leq C(\|\dv\bff\|_{L_q(B^1_j)}
+ \|\nabla u\|_{L_q(B^1_j)},$$
which follows from \eqref{ad:7}, we have
$\nabla v = \nabla(\zeta^1_ju\circ\Phi_j) \in H^1_q(\HS)^N$
and
$$\|\nabla(\zeta^1_ju\circ\Phi_j)\|_{H^1_q(\HS)}
\leq C(\|\dv\bff\|_{L_q(B^1_j\cap\Omega)}
+ \|\nabla u\|_{L_q(B^1_j\cap\Omega)}).
$$
Since $\|u\|_{L_q(B^1_j\cap \Omega)} 
\leq c_1\|\nabla u\|_{L_q(B^1_j\cap \Omega)}$ as follows 
from Lemma \ref{lem:5.7}, we have
\begin{equation}\label{ad:9}
\|\zeta^1_j\nabla^2u\|_{L_q(\Omega)}
\leq C(\|\dv\bff\|_{L_q(B^1_j\cap \Omega)}
+ \|\nabla u\|_{L_q(B^1_j\cap\Omega)}).
\end{equation}
Combining \eqref{ad:3}, \eqref{ad:9} and 
\eqref{D.1*} gives
$$\|\nabla^2u\|_{L_q(\Omega)} 
\leq C(\|\dv\bff\|_{L_q(\Omega)} + \|\nabla u\|_{L_q(\Omega)})
\leq C(\|\dv\bff\|_{L_q(\Omega)} + \|\bff\|_{L_q(\Omega)}),
$$
which completes the proof of Theorem \ref{thm:ap-r}.


\subsection{Laplace-Beltrami operator}\label{subsec:p-4}

In this subsection, we
 introduce the  Laplace-Beltrami operatrs and 
some important formulas from differential geometry, 

Let $\Gamma$ be a hypersurface of class $C^3$ in $\BR^N$.
The $\Gamma$ is parametrized as 
$p=\phi(\theta)={}^\top(\phi_1(\theta), \ldots, \phi_N(\theta))$
 locally
at $p \in \Gamma$,  
where $\theta=(\theta_1, \ldots, \theta_{N-1})$ runs through a domain
$\Theta \subset \BR^{N-1}$. Let 
\begin{equation}\label{cg:0}
\tau_i = \tau_i(p) = \frac{\pd}{\pd \theta_i}\phi(\theta)
= \pd_i\phi\quad(i=1, \ldots, N-1).
\end{equation}
which forms a basis of the tangent space $T_p\Gamma$ of $\Gamma$ at
$p$. Let $\bn=\bn(p)$ denote the 
outer unit normal of $\Gamma$ at $p$. 
Notice that 
\begin{equation}\label{cg:1}
<\tau_i, \bn> = 0.
\end{equation}
 Here and in the following, 
$<\cdot, \cdot>$ denotes a standard inner product in $\BR^N$. 
To introduce the formula of $\bn$, we notice that
$$\det\left(\begin{matrix}
\frac{\pd \phi_1}{\pd\theta_1} & \cdots & \frac{\pd\phi_1}{\pd\theta_{N-1}} &
\frac{\pd\phi_1}{\pd\theta_k} \\
\vdots & \ddots & \vdots & \vdots \\
\frac{\pd \phi_N}{\pd\theta_1} & \cdots & \frac{\pd\phi_N}{\pd\theta_{N-1}} &
\frac{\pd\phi_N}{\pd\theta_k}
\end{matrix}\right)
=0 
$$
for any $k=1, \ldots, N-1$.  Thus, to satisfy \eqref{cg:1}, 
$\bn$ is defined by 
\begin{equation}\label{cg:3}
\bn = {}^\top\frac{(h_1, \ldots, h_N)}{H}
\end{equation}
with $H = \sqrt{\sum_{j=1}^N h_i^2}$ and 
$$H_i = \frac{\pd(\phi_1,\ldots, \hat\phi_i, \ldots, \phi_N)}
{\pd(\theta_1, \ldots,\theta_{N-1})}
= (-1)^{N+i}\det
\left(\begin{matrix}
\frac{\pd\phi_1}{\pd\theta_1} & \cdots & \frac{\pd\phi_1}{\pd\theta_{N-1}} \\
\vdots & \ddots & \vdots \\
\frac{\pd\phi_{i-1}}{\pd\theta_1} & \cdots & \frac{\pd\phi_{i-1}}{\pd\theta_{N-1}} \\
\frac{\pd\phi_{i+1}}{\pd\theta_1} & \cdots & \frac{\pd\phi_{i+1}}{\pd\theta_{N-1}}\\
\vdots & \ddots & \vdots \\
\frac{\pd\phi_{N}}{\pd\theta_1} & \cdots & \frac{\pd\phi_{N}}{\pd\theta_{N-1}}
\end{matrix}\right).
$$
For example, when $N=3$, 
\begin{align*}
\bn &= \frac{\frac{\pd\phi}{\pd\theta_1}\times \frac{\pd\phi}{\pd\theta_2}}
{\Bigl|\frac{\pd\phi}{\pd\theta_1}\times \frac{\pd\phi}{\pd\theta_2}
\Bigr|}
=H^{-1}{}^\top(h_1, h_2, h_3)
\quad
h_1 = 
\frac{\pd\phi_2}{\pd\theta_1}\frac{\pd\phi_3}{\pd\theta_2}-
\frac{\pd\phi_3}{\pd\theta_1}\frac{\pd\phi_2}{\pd\theta_2}, \\
h_2 &= 
\frac{\pd\phi_3}{\pd\theta_1}\frac{\pd\phi_1}{\pd\theta_2}-
\frac{\pd\phi_1}{\pd\theta_1}\frac{\pd\phi_3}{\pd\theta_2},
\enskip   
h_3 = 
\frac{\pd\phi_1}{\pd\theta_1}\frac{\pd\phi_2}{\pd\theta_2}-
\frac{\pd\phi_2}{\pd\theta_1}\frac{\pd\phi_1}{\pd\theta_2},
\enskip
H =\sqrt{ h_1^2 + h_2^2 + h_3^2}.
\end{align*}
Let
$$g_{ij} = g_{ij}(p) = <\tau_i, \tau_j>\quad(i, j=1, \ldots, N-1),
$$
and let $G$ be an $(N-1)\times(N-1)$ matrix whose $(i, j)^{\rm th}$
components are $g_{ij}$. The matrix $G$ is called the {\bf 
first fundamental 
form} of $\Gamma$. In the following, we employ Einstein's summation 
convention, which means that equal lower and upper indices are
to be summed. Since 
for 
any $\xi \in \BR^{N-1}$ with $\xi\not=0$, $<G\xi, \xi> = g_{ij}\xi^i\xi^j 
=<\xi^i\tau_i, \xi^j\tau_j> = |\xi^i\tau_i|^2 > 0$, 
$G$ is a positive symmetric
matrix, and therefore $G^{-1}$ exists. 
  Let $g^{ij}$ be the $(i, j)^{\rm th}$
component of $G^{-1}$ and let 
$$\tau^i = g^{ij}\tau_j.$$
Using 
$g^{ik}g_{kj} = g_{ik}g^{kj} = \delta^i_j$,  
where  $\delta^i_j$ 
are the Kronecker delta symbols defined by $\delta^i_i=1$ and 
$\delta^i_j=0$ for $i\not=j$,  we have
\begin{equation}\label{cg:dual}
<\tau_i, \tau^j> = <\tau^i, \tau_j> = \delta^i_j.
\end{equation} 
In fact, 
$<\tau_i, \tau^j> = <\tau_i, g^{jk}\tau_k> =
g^{jk}<\tau_i, \tau_k> = g^{jk}g_{ki} = \delta^j_i$. 
Thus, $\{\tau^j\}_{j=1}^{N-1}$ is a dual basis of $\{\tau_i\}_{i=1}^{N-1}$. 
In particular, we have 
\begin{equation}\label{cg:4}
\tau^i = g^{ij}\tau_j, \quad \tau_i = g_{ij}\tau^j,
\end{equation}

For any $\ba \in T_p\Gamma$, we write
$$\ba = a^i\tau_i = a_i\tau^j 
$$  
 By \eqref{cg:dual}, we have
$<\ba, \tau^i> = <a^j\tau_j, \tau^i> = a^j<\tau_j, \tau^i>=a^i$
and $<\ba, \tau_i> = a_j<\tau^j, \tau_i> = a_i$,  and so
\begin{equation}\label{cg:cont1}
\ba = <\ba, \tau^i>\tau_i = <\ba, \tau_i>\tau^i.
\end{equation}
In particular, by \eqref{cg:dual} and \eqref{cg:4}, we have
$$a_i = g_{ij}a^j, \quad a^i = g^{ij}a_j
$$
with $a_i=<\ba, \tau_i>$ and $a^i = <\ba, \tau^i>$. 
Notice that $\{\tau_1, \ldots, \tau_{N-1}, \bn\}$ forms a
basis of $\BR^N$. Namely, for any $N$ vector $\bb\in\BR^N$
we have
$$\bb = b^i\tau_i + <\bb, \bn>\bn = b_i\tau^i + <\bb, \bn>\bn$$
with $b^i=<\bb, \tau^i>$ and $b_i = <\bb, \tau_i>$. 

We next consider $\tau_{ij} = \pd_i\pd_j\phi = \pd_j\tau_i$.
Notice that $\tau_{ij} = \tau_{ji}$. 
Let
\begin{equation}\label{cg:5}
\Lambda^k_{ij} = <\tau_{ij}, \tau^k>, \quad
\ell_{ij} = <\tau_{ij}, \bn>.
\end{equation}
and then, we have
\begin{equation}\label{cg:6}
\tau_{ij} = \Lambda^k_{ij}\tau_k + \ell_{ij}\bn.
\end{equation}
Let $L$ be an $N-1\times N-1$ matrix whose $(i, j)^{\rm th}$ 
component is $\ell_{ij}$, which is 
called the {\bf second fundamental form} of $\Gamma$.
 Let
\begin{equation}\label{cg:7}
\CH(\Gamma) = \frac{1}{N-1}{\rm tr}\,(G^{-1}L) = 
\frac{1}{N-1}g^{ij}\ell_{ij}.
\end{equation}
The $\CH(\Gamma)$ is called the {\bf mean curvature} of 
$\CH(\Gamma)$. 

Let 
$g = \det G$.  One of the most important formulas
in this section is 
\begin{equation}\label{cg:8}
\pd_i(\sqrt{g}g^{ij}\tau_j) = \sqrt{g}g^{ij}\ell_{ij}\bn.
\end{equation}
This formula will be proved below. 

The $\Lambda^k_{ij}$ is called {\bf Christoffel symbols}. We know 
the formula:
\begin{equation}\label{cg:9}
\Lambda^r_{ij} = \frac12 g^{rk}(\pd_ig_{jk} + \pd_jg_{ki}
- \pd_kg_{ij}).
\end{equation}
In fact, 
\begin{align*}
\pd_ig_{jk} &= \pd_i<\tau_j, \tau_k>
= <\tau_{ij}, \tau_k> + <\tau_j, \tau_{ki}>, \\
\pd_jg_{ki} &= \pd_j<\tau_k, \tau_i>
= <\tau_{jk}, \tau_i> + <\tau_k, \tau_{ij}>, \\
\pd_kg_{ij} &= \pd_k<\tau_i, \tau_j>
= <\tau_{ki}, \tau_j> + <\tau_i, \tau_{jk}>,
\end{align*}
and so we have
$$\pd_ig_{jk} + \pd_jg_{ki}
- \pd_kg_{ij} = 2<\tau_{ij}, \tau_k>,
$$
which implies \eqref{cg:9}.  

We now prove \eqref{cg:8} by studying several steps. We first prove 
\begin{equation}\label{cg:8.1}
\pd_kg_{ij}= g_{jr}\Lambda^r_{ki} + g_{ir}\Lambda^r_{kj}.
\end{equation}
In fact, by \eqref{cg:6} and $<\tau_i, \bn>=0$, 
\begin{align*}
\pd_kg_{ij} &= <\tau_{ki}, \tau_j> + <\tau_i, \tau_{kj}>
= <\Lambda^r_{ki}\tau_r, \tau_j> + <\tau_i, \Lambda^r_{kj}\tau_r>
\\
&=g_{jr}\Lambda^r_{ki} + 
g_{ir}\Lambda^r_{kj}.
\end{align*}
We next prove
\begin{equation}\label{cg:8.2}
\pd_kg^{ij} = -g^{ir}\Lambda^j_{rk} - g^{jr}\Lambda^i_{rk},
\end{equation}
In fact, by $g^{ij}g_{j\ell} = \delta^i_\ell$,  we have
$$0 = \pd_k(g^{ij}g_{jk})
= (\pd_kg^{ij})g_{j\ell} + g^{ij}\pd_kg_{j\ell}.
$$
Thus, using \eqref{cg:8.1}, we have
\begin{align*}
(\pd_kg^{ij})g_{j\ell}&= -g^{ij}\pd_kg_{j\ell}
=-g^{ij}(g_{\ell r}\Lambda^r_{kj} + g_{jr}\Lambda^r_{k\ell})\\
&= -g^{ij}g_{\ell r}\Lambda^r_{kj} - \delta^i_r\Lambda^r_{k\ell}
= -g^{ij}g_{\ell r}\Lambda^r_{kj} - \Lambda^i_{k\ell}.
\end{align*}
Thus, 
\begin{align*}
\pd_kg^{im} & =  (\pd_kg^{ij})\delta^m_j
=(\pd_kg^{ij})g_{j\ell}g^{\ell m} 
= -g^{\ell m}(g^{ij}g_{\ell r}\Lambda^r_{kj} + 
\Lambda^i_{k\ell})\\
&= -\delta^m_rg^{ij}\Lambda^r_{kj}-g^{\ell m}\Lambda^i_{k\ell}
= -g^{ij}\Lambda^m_{kj} - g^{m\ell}\Lambda^i_{k\ell}
= -g^{i\ell}\Lambda^m_{k\ell} - g^{m\ell}\Lambda^i_{k\ell}.
\end{align*}
Thus, setting $m=j$ and $\ell=r$, we have \eqref{cg:8.2}, because 
$\Lambda^m_{k\ell}= \Lambda^m_{\ell k}$ and 
$\Lambda^i_{k\ell} = \Lambda^i_{\ell k}$.

We next prove
\begin{equation}\label{cg:8.3}
\pd_ig = 2\Lambda^j_{ij}g.
\end{equation}
Recall that $g = \det G$ and the $(i, j)^{\rm th}$ component of 
$G$ is $<\tau_i, \tau_j>$. From the definition of differentiation,
we have
\begin{align*}
&\frac{\det G(x+\be_i\Delta x_i)-\det G(x)}{\Delta x_i}\\
&= \frac{\det(G(x) + \pd_iG(x)\Delta x_i) - \det G(x)}{\Delta x_i}
+ O(\Delta x_i)\\
& = \frac{\det G(x)(\det(I+\pd_iG(x) G^{-1}(x)\Delta x_i)-1)}{\Delta
x_i} + O(\Delta x_i) \\
& = \det G(x)\,{\rm tr}\,(\pd_iG(x)G^{-1}(x)) + O(\Delta x_i).
\end{align*}
Thus, we have
$$\pd_i g=\det G\,{\rm tr}\,(\pd_iG G^{-1}).$$
Using \eqref{cg:6} and $<\tau_i, \bn>=0$, we have
\begin{align*}
{\rm tr}\,(\pd_iG G^{-1})
&= \pd_i(<\tau_j, \tau_k>)g^{kj}
= (<\tau_{ij}, \tau_k> + <\tau_j, \tau_{ik}>)g^{kj}\\
&= ((\Lambda^r_{ij}\tau_r, \tau_k> + 
<\Lambda^r_{ik}\tau_r, \tau_j>)g^{kj}
=(g_{rk}\Lambda^r_{ij} + g_{rj}\Lambda^r_{ik})
g^{kj} \\
&= 2\Lambda^j_{ij}
\end{align*}
Putting these two formulas gives \eqref{cg:8.3}.

We next prove
\begin{equation}\label{cg:8.4}
\pd_i(\sqrt{g}g^{ij}) = 
-\sqrt{g}g^{ik}\Lambda^j_{ik}.
\end{equation}
In fact, by \eqref{cg:8.2} and \eqref{cg:8.3}, 
\begin{align*}
\frac{1}{\sqrt{g}}\pd_i(\sqrt{g}g^{ij})
&= \frac{1}{2g}(\pd_ig)g^{ij} + \pd_ig^{ij}
=\frac{1}{2g}2\Lambda^\ell_{i\ell}gg^{ij}
-g^{ir}\Lambda^j_{ri}-g^{jr}\Lambda^i_{ri}\\
& = g^{jr}\Lambda^\ell_{r\ell}-g^{ir}\Lambda^j_{ri}
-g^{jr}\Lambda^\ell_{\ell r}
= -g^{ik}\Lambda^j_{ik}.
\end{align*}
Thus, we have \eqref{cg:8.4}.

We now prove \eqref{cg:8}. By \eqref{cg:6} and \eqref{cg:8.4},
\begin{align*}
\pd_i(\sqrt{g}g^{ij}\tau_j) &= \pd_i(\sqrt{g}g^{ij})\tau_j
+ \sqrt{g}g^{ij}\pd_i\tau_{ij} \\
& = -\sqrt{g}g^{ik}\Lambda^j_{ik}\tau_j 
+ \sqrt{g}g^{ij}(\Lambda^k_{ij}\tau_k + \ell_{ij}\bn)\\
& =-\sqrt{g}g^{ij}\Lambda^k_{ij}\tau_k + \sqrt{g}g^{ij}
\Lambda^k_{ij}\tau_k + \sqrt{g}g^{ij}\ell_{ij}\bn\\
& = \sqrt{g}g^{ij}\ell_{ij}\bn.
\end{align*}
Thus, we have\eqref{cg:8}. 

We now introduce the {\bf Laplace-Beltrami operator} 
$\Delta_\Gamma$  on
$\Gamma$, which is defined by
$$
\Delta_\Gamma f = \frac{1}{\sqrt{g}}\pd_i(\sqrt{g}g^{ij}
\pd_jf).
$$
By \eqref{cg:8.4}, we have
\begin{equation}\label{cg:10}
\Delta_\Gamma f = g^{ij}\pd_i\pd_jf -g^{ik}\Lambda^j_{ik}\pd_jf.
\end{equation}
By \eqref{cg:7} and \eqref{cg:8}, we have
$$
\Delta_\Gamma\phi = (N-1)\CH(\Gamma)\bn.
$$
Usually, we put $H(\Gamma) = (N-1)\CH(\Gamma)$, and so we have
\begin{equation}\label{cg:11}
\Delta_\Gamma x = H(\Gamma)\bn\quad \text{for $x \in \Gamma$.}
\end{equation}


One fundamental result for the Laplace-Beltrami operator is the following.
\begin{lem}\label{lem:3.2} Let $1 < q < \infty$.  Assume that 
$\Omega$ is a uniform $C^3$ domain and let $\Gamma$ be the boundary 
of $\Omega$.  Then, there exists a large number $m > 0$ such that 
for any $f \in W^{-1/q}_q(\Gamma)$, there exists a unique 
$v \in W^{2-1/q}_q(\Gamma)$ such that $v$ satisfies the equation
$$(m-\Delta_\Gamma)v = f\quad\text{on $\Gamma$} $$
and the estimate 
$$\|v\|_{W^{2-1/q}_q(\Gamma)} \leq C\|f\|_{W^{-1/q}_q(\Gamma)}
$$
for some constant $C > 0$. 
\end{lem}
\begin{rem}\thetag1
Lemma \ref{lem:3.2} was proved  by Amann-Hieber-Simonett
\cite[Theorem 10.3]{AHS} in the case where $\Gamma$ has finite covering, and 
by Shibata \cite[Theorem 2.2]{S3}  in the case where $\Gamma$ is the 
boundary of a uniform $C^2$ domain. 
\\
\thetag2 Let $(m-\Delta_\Gamma)^{-1}$ be defined by 
$(m-\Delta_\Gamma)^{-1}f = v$.
\end{rem}

\subsection{Parametrized surface}\label{sec:para}

Let $\phi_t : \Omega\to\BR^N$ be a injection map with suitable regularity for 
each time $t \geq 0$. Let $\Gamma$ be a $C^3$ hypersurface $\subset \Omega$ 
and set 
$$\Gamma_t = \{ x = \phi_t(y) \mid y \in \Gamma\}. $$
First we prove Theorem \ref{thm:rey}. Let $J(t)$ be the Jacobian of the 
map $x = \phi(y)$. Let $\bw(x, t) = ( \pd_t\phi_t)(\phi^{-1}(x))$.
We shall prove that
\begin{equation}\label{eq:ray}
\frac{\pd}{\pd t}J(t) = (\dv_x \bw(x, t))J(t).
\end{equation}
This formula is called a Reynols transport theorem. 
In fact, in view of the definition of differentiation, we consider
$$
\frac{J(t+\Delta t)-J(t)}{\Delta t}.
$$
Writing $\phi_t = {}^\top(x_1(y,t), \ldots, x_N(y,t))$, we have 
$J(t) = \det X(t)$ with 
$$X(t) = \left(\begin{matrix}
\frac{\pd x_1}{\pd y_1} & \ldots & \frac{\pd x_1}{\pd y_N} \\
\vdots & \ddots &\vdots \\
\frac{\pd x_N}{\pd y_1} & \ldots & \frac{\pd x_N}{\pd y_N}
\end{matrix}\right).
$$
By  mean value theorem, we write 
$X(t+\Delta t) = X(t) + X'(t)\Delta t + O((\Delta t)^2)$ with 
$$X'(t) = \left(\begin{matrix}
\frac{\pd^2 x_1}{\pd y_1 \pd t} & \ldots & \frac{\pd^2 x_1}{\pd y_N\pd t} \\
\vdots & \ddots &\vdots \\
\frac{\pd^2 x_N}{\pd y_1 \pd t} & \ldots & \frac{\pd^2 x_N}{\pd y_N\pd t}
\end{matrix}\right).
$$
Using this symbol,  we write
\begin{align*}
&\frac{J(t+\Delta t)-J(t)}{\Delta t}
= \frac{\det X(t+\Delta t) - \det X(t)}{\Delta t} \\ 
&\quad=\frac{\det (X(t) + \Delta t X'(t) + O((\Delta t)^2))-
\det X(t)}{\Delta t} \\
&\quad = \frac{\det X(t)(\det(\bI + \Delta t X'(t)X^{-1}(t)
+ O((\Delta t)^2)) - 1)}{\Delta t}\\
&\quad = (\det X(t)){\rm tr}(X'(t)X^{-1}(t)) + O(\Delta t).
\end{align*}
Since the $(i, j)^{\rm th}$ component of $X'(t)X(t)^{-1}$ is 
$$\sum_{k=1}^N\frac{\pd^2 x_i}{\pd y_k \pd t}\frac{\pd y_k}{\pd x_j}
= \frac{\pd w_i}{\pd x_j}, 
$$
where $\bw(x)= (\pd_t\phi)(\phi^{-1}(x)) = {}^\top(w_1, \ldots, w_N)$,
we have ${\rm tr}\,(X'(t)X^{-1})(t) = \dv_x\bw$, which shows \eqref{eq:ray}.

We also prove that
\begin{equation}\label{eq:energy}
\frac{\pd}{\pd t}|\Gamma_t| = -\int_{\Gamma_t} H(\Gamma_t)
<\bn_t, \dot\phi>\,d\sigma,
\end{equation}
where 
$|\Gamma_t|$ is the area of $\Gamma_t$,
$\dot\phi = \pd_t\phi$, 
and $d\sigma$ is the surface element of
$\Gamma_t$.  The formula $<\bn_t, \dot\phi>$ 
denotes the velocity of the evolution of 
$\Gamma_t$ with respect to $\bn_t$. 

To prove \eqref{eq:energy}, we have to parametrize
$\Gamma_t$.  
Let $\Gamma$ be parametrized as 
$y= y(\theta)$ for $\theta \in \Theta \subset \BR^{N-1}$,
and then the first fundamental form of 
$\Gamma_t$ is given by $G_t = (g_{ij}(t))$ with
$g_{ij}(t) = <\tau_i(t), \tau_j(t)>$, where 
$$\tau_i(t) = \frac{\pd \phi_t(y(\theta))}{\pd\theta_i}.$$
Let $g(t) = \det G_t$ and 
$G_t^{-1} = (g^{ij}(t))$.  Since the surface element of
$\Gamma_t$ is given by 
 $d\sigma = \sqrt{g}d\theta$, we have  
$$|\Gamma(t)| = \int_{\Gamma_t}\,d\sigma
= \int_\Theta \sqrt{g(t)}\,d\theta.$$
Thus, 
\begin{align*}
\frac{d}{dt}|\Gamma(t)| 
= \int_\Theta \frac{d}{dt}\sqrt{g(t)}
\,d\theta
= \int_\Theta\frac{\dot g}{2\sqrt{g}}\,d\theta
\end{align*}
To find $\dot g$, we calculate $(g(t+\Delta t) - g(t))/\Delta t$
as follows:
\begin{align*}
&\frac{g(t+\Delta t) - g(t)}{\Delta t} \\
&\quad = \frac{\det G(t+\Delta t) - \det G(t)}{\Delta t}\\
&\quad = \frac{\det (G(t) + \pd_t G(t)\Delta t + O((\Delta t)^2)-
\det G(t)}{\Delta t} \\
& \quad = \frac{\det(G(t)(\det(\bI +\pd_tG(t)G(t)^{-1}\Delta t
+ O((\Delta t)^2))-1)}{\Delta t}\\
&\quad = \det G(t)){\rm tr}(\pd_tG(t)G(t)^{-1}) + O(\Delta t)
\end{align*}
Thus, we have
\allowdisplaybreaks{
\begin{align*}
\dot g &= g\,{\rm tr}\, (\pd_tG(t)G(t)^{-1})
=g(<\dot \tau_i, \tau_j> + <\tau_i, \dot\tau_j>)g^{ji}\\
&= g(<\dot\tau_i, \tau_j>g^{ji} + < \dot\tau_i, \tau_j>g^{ij})
= 2gg^{ij}<\dot \tau_i, \tau_j>
\end{align*}
}
where we have used $g^{ij} = g^{ji}$. Putting these formulas
together gives 
\begin{align*}
\frac{d}{dt}|\Gamma(t)| &= \int_\Theta \sqrt{g}g^{ij}<\dot\tau_i, \tau_j>\,
d\theta\\
& = \int_{\Theta} \pd_i(\sqrt{g}g^{ij}<\dot \phi, \tau_j>)\,d\theta
- \int_\Theta<\pd_i(\sqrt{g} g^{ij}\tau_j), \dot\phi>\,d\theta \\
&= -\int_\Theta \frac{1}{\sqrt{g}}
<\pd_i(\sqrt{g} g^{ij}\pd_j\phi), \dot\phi>\,
\sqrt{g}\,d\theta
= -\int_{\Gamma_t} <\Delta_{\Gamma_t}\phi, \dot\phi>\,
d\sigma\\
& = -\int_{\Gamma_t}H(\Gamma_t)<\bn_t, \dot\phi>\,d\sigma,
\end{align*}
where we have used \eqref{cg:11}. 
This shows \eqref{eq:energy}

\subsection{Example of mean curvature}
Let $\Gamma$ be a closed hypersurface in $\BR^3$ defined by 
$|x| = r(\omega)$ for $\omega \in S_1 = \{\omega \in \BR^3 \mid
|\omega|=1\}$. We introduce the polar coordinates:
$$\omega_1 = \cos\varphi\sin\theta, \enskip
\omega_2 = \sin\varphi\sin\theta, \enskip \omega_3 = \cos\theta
$$
for $\varphi \in [0, 2\pi)$ and $\theta \in [0, \pi)$. And then, 
$\Gamma$ is represented by 
$$x_1=r(\varphi, \theta)\cos\varphi\sin\theta,
\enskip
x_2= r(\varphi, \theta)\sin\varphi\sin\theta, \enskip
x_3 = r(\varphi, \theta)\cos\theta.
$$
Let $H(\Gamma)$ be the doubled mean curvature of 
$\Gamma$. 
In the sequel, we will show the following well-known formula (cf. \cite{Sol1}):
\begin{equation}\label{ex:mean}\begin{aligned}
H(\Gamma) &= \frac{1}{r\sin\theta}
\Bigl\{\frac{\pd}{\pd\varphi}\Bigl(\frac{r_\varphi}
{\sin\theta\sqrt{r^2+|\nabla r|^2}}\Bigr)
+ \frac{\pd}{\pd\theta}\Bigl(\frac{\sin\theta r_\theta}{\sqrt{r^2
+ |\nabla r|^2}}\Bigr)\Bigr\}\\
&\quad-\frac{2}{\sqrt{r^2 + |\nabla r|^2}},
\end{aligned}\end{equation}
where $\nabla r= (r_\theta, r/\sin\theta)$, and so 
$|\nabla r|^2 = r_\theta^2 + (r_\varphi/\sin\theta)^2.$

From the definition, we have
$$H(\Gamma) = \sum_{ij=1}^2 g^{ij}\ell_{ij},
$$
where $g^{ij}$ is the $(i,j)^{\rm th}$ component of the inverse
matrix of the first fundamental form and $\ell_{ij}$ is
the $(i, j)^{\rm th}$ component of the second fundamental 
form. We have to find $g^{ij}$ and $\ell_{ij}$.
We first   find $g_{ij}$.  Since
\begin{equation}\label{ex:1.1}\begin{aligned}
\frac{\pd x}{\pd\varphi} & = r_\varphi\left(\begin{matrix}
\cos\varphi\sin\theta \\
\sin\varphi\sin\theta \\
\cos\theta\end{matrix}\right) + 
r\left(\begin{matrix}
-\sin\varphi\sin\theta \\
\cos\varphi\sin\theta \\
0\end{matrix}\right), \\
\frac{\pd x}{\pd\theta} & = r_\theta
\left(\begin{matrix} 
\cos\varphi\sin\theta \\
\sin\varphi\sin\theta \\
\cos\theta\end{matrix}\right)
+ r
\left(\begin{matrix} 
\cos\varphi\cos\theta \\
\sin\varphi\cos\theta \\
-\sin\theta\end{matrix}\right).
\end{aligned}\end{equation}
the first fundamental form $G$ is given by
$G=\left(\begin{matrix} g_{11} & g_{12} \\ g_{12} & g_{22}\end{matrix}\right)$
with 
\begin{align*}
g_{11} & = \frac{\pd x}{\pd\varphi}\cdot\frac{\pd x}{\pd\varphi}
= r_\varphi^2 + r^2\sin^2\theta, \\
g_{12} & = \frac{\pd x}{\pd \varphi}\cdot\frac{\pd x}{\pd \theta}
= r_\varphi r_\theta\\
g_{22} & =\frac{\pd x}{\pd\theta}\cdot\frac{\pd x}{\pd\theta}
= r_\theta^2 + r^2.
\end{align*}
And then ,
$$g = \det G = r^2\sin^2\theta(r^2 + r_\theta^2 + (r_\varphi/\sin\theta)^2)
= r^2\sin^2\theta(r^2+|\nabla r|^2),
$$
where $|\nabla r|^2 = r_\theta^2 + (r_\varphi/\sin\theta)^2$. Thus, 
$$G^{-1} = \frac{1}{g}\left(\begin{matrix} r^2_\theta + r^2 & -r_\varphi 
r_\theta \\
-r_\varphi r_\theta & r_\varphi^2 + r^2\sin^2\theta
\end{matrix}\right).
$$
and so 
\begin{equation}\label{ex:inverse}
g^{11} = \frac{r_\theta^2+ r^2}{g}, 
\enskip
g^{12}=g^{21} = -\frac{r_\varphi r_\theta}{g}, \enskip
g^{22} = \frac{r_\varphi^2+r^2\sin^2\theta}{g}.
\end{equation}
We next caluculate the second fundamental form of $\Gamma$.  For this purpose,
we first calculate the unit outer normal $\bn_{_\Gamma}$ to 
$\Gamma$, which is given by 
\begin{equation}\label{ex:normal}
\bn_{_\Gamma} = 
\frac{1}{\sqrt{g}}\left(
\begin{matrix} 
r_\varphi r\sin\varphi - r_\theta r\cos\varphi\sin\theta\cos\theta
+ r^2\cos\varphi\sin^2\theta \\
-r_\varphi r\cos\varphi -r_\theta r\sin\varphi\sin\theta\cos\theta
+ r^2\sin\varphi\sin^2\theta \\
r_\theta r\sin^2\theta + r^2\sin\theta\cos\theta
\end{matrix}\right).
\end{equation}
We next find the second fundamental form.  For this, first we 
calculate the second derivatives of $x$ as follows:
\begin{align}
\frac{\pd^2x}{\pd\varphi^2}
& = r_{\varphi\varphi}\left(\begin{matrix}
\cos\varphi\sin\theta \\
\sin\varphi\sin\theta \\
\cos\theta
\end{matrix}\right)
+ 2r_\varphi\left(
\begin{matrix}
-\sin\varphi\sin\theta \\
\cos\varphi\sin\theta \\
0
\end{matrix}\right)
+r\left(\begin{matrix}
-\cos\varphi\sin\theta \\
-\sin\varphi\sin\theta \\
0
\end{matrix}\right); \nonumber  \\
\frac{\pd^2 x}{\pd\varphi\pd\theta}
& = r_{\varphi \theta}
\left(\begin{matrix}
\cos\varphi\sin\theta \\
\sin\varphi \sin\theta \\
\cos\theta
\end{matrix}\right)
+ r_\varphi\left(\begin{matrix}
\cos\varphi\cos\theta\\
\sin\varphi\cos\theta \\
-\sin\theta
\end{matrix}\right)
 +r_\theta\left(
\begin{matrix}
-\sin\varphi\sin\theta \\
\cos\varphi\sin\theta \\ 0
\end{matrix}\right) \\
&\quad + r\left(\begin{matrix}
-\sin\varphi\cos\theta \\
\cos\varphi\cos\theta \\
0
\end{matrix}\right); \nonumber \\
\frac{\pd^2 x}{\pd\theta^2}
& = 
r_{\theta\theta}\left(\begin{matrix}
\cos\varphi\sin\theta \\
\sin\varphi\sin\theta \\
\cos\theta
\end{matrix}\right)
+ 2r_\theta\left(
\begin{matrix} \cos\varphi\cos\theta \\
\sin\varphi\cos\theta \\
-\sin\theta
\end{matrix}\right)
+ r\left(\begin{matrix}
-\cos\varphi\sin\theta \\
-\sin\varphi\sin\theta \\
-\cos\theta
\end{matrix}\right).\label{ex:2d}
\end{align}
Thus, we have
\begin{equation}\label{ex:second}
\begin{aligned}
\ell_{11} & = <\frac{\pd^2x}{\pd\varphi^2}, \bn_{_\Gamma}>\\
&= \frac{1}{\sqrt{g}}
(r^2r_{\varphi\varphi}\sin\theta-2rr_\varphi^2\sin\theta
+r^2r_\theta\sin^2\theta\cos\theta-r^3\sin^3\theta),\\
\ell_{12}& = <\frac{\pd^2x}{\pd\varphi\pd\theta}, \bn_{_\Gamma}>
\\
&= \frac{1}{\sqrt{g}}
(r^2r_{\varphi\theta}\sin\theta -2rr_\varphi r_\theta\sin\theta
-r^2r_\varphi\cos\theta), \\
\ell_{22} & = <\frac{\pd^2x}{\pd\theta^2}, \bn_{_\Gamma}> \\
&=\frac{1}{\sqrt{g}}(r^2r_{\theta\theta}\sin\theta
-2r_\theta^2r\sin\theta-r^3\sin\theta).
\end{aligned}
\end{equation}

We now calculate $H(\Gamma)$. Noting that 
$g^{12}=g^{21}$ and $\ell_{12} = \ell_{21}$, 
by \eqref{ex:inverse} and \eqref{ex:second}, 
we have
\begin{equation}\label{ex:carv1}\begin{aligned}
H(\Gamma) &= g^{11}\ell_{11} + 2g^{12}\ell_{12} + g^{22}\ell_{22} \\
&= \frac{1}{g^{3/2}}\{
((r^2_\theta+r^2)r_{\varphi\varphi} -2r_\varphi r_\theta r_{\varphi\theta}
+ r_\varphi^2 r_{\theta\theta} - 3rr_\varphi^2)r^2\sin\theta \\
&\quad\phantom{\times}
 + (r r_{\theta\theta} - 3r_\theta^2 - 2r^2)r^3\sin^3\theta
+ (r_\theta+r^2)r_\theta r^2\sin^2\theta\cos\theta \\
&\quad\phantom{\times}
+2r_\varphi^2r_\theta r^2\cos\theta\}
\end{aligned}\end{equation}
On the other hand, we have
\begin{align*}
&\frac{1}{r\sin\theta}\frac{\pd}{\pd\varphi}\frac{r_\varphi}{\sin\theta
\sqrt{r^2+|\nabla r|^2}}
= \frac{r^2\sin\theta((r^2+r^2_\theta)r_{\varphi\varphi}-rr^2_\varphi
-r_\varphi r_\theta r_{\theta\varphi})}{r^3\sin^3\theta
(r^2+|\nabla r|^2)^{3/2}}; \\
&\frac{1}{r\sin\theta}\frac{\pd}{\pd\theta}
\frac{\sin\theta r_\theta}{\sqrt{r^2+|\nabla r|^2}}
= \frac{1}{r^3\sin^3\theta (r^2+|\nabla r|^2)^{3/2}}\\
&\times\{r^2(r^2r_{\theta\theta} - r^3r_\theta^2)\sin^3\theta
+ r^2(r_\varphi^2r_{\theta\theta}-r_\varphi r_\theta r_{\theta\varphi})\sin\theta\\
&\qquad+r^2(r^2+r_\theta^2)r_\theta\sin^2\theta\cos\theta + 2r^2r_\theta
r_\varphi^2\cos\theta\}.
\end{align*}
Thus,  we have
\begin{align*}
\frac{1}{r\sin\theta}
&\Bigl\{\frac{\pd}{\pd\varphi}\Bigl(\frac{r_\varphi}
{\sin\theta\sqrt{r^2+|\nabla r|^2}}\Bigr)
+ \frac{\pd}{\pd\theta}\Bigl(\frac{\sin\theta r_\theta}{\sqrt{r^2
+ |\nabla r|^2}}\Bigr)\Bigr\}
-\frac{2}{\sqrt{r^2+|\nabla r|^2}} \\
&\quad= \frac{B}{r^3\sin^3\theta(r^2+|\nabla r|^2)^{3/2}}
\end{align*}
with
\begin{align*}
B &= r^2\sin\theta((r^2+r^2_\theta)r_{\varphi\varphi}-rr^2_\varphi
-r_\varphi r_\theta r_{\theta\varphi}) 
+ r^2(r^2r_{\theta\theta} - rr_\theta^2)\sin^3\theta \\
&+ r^2(r_\varphi^2r_{\theta\theta}-r_\varphi r_\theta r_{\theta\varphi})\sin\theta
+r^2(r^2+r_\theta^2)r_\theta\sin^2\theta\cos\theta + 2r^2r_\theta
r_\varphi^2\cos\theta\\
&-2r^3\sin^3\theta(r^2+r^2_\theta + r^2_\varphi\sin^{-2}\theta)\\
& = r^2((r^2+r_\theta^2)r_{\varphi\varphi} + r_\varphi^2r_{\theta\theta}
-3rr_\varphi^2-2r_\varphi r_\theta r_{\varphi\theta})\sin\theta \\
&+ r^3\sin^3\theta(rr_{\theta\theta}-3r_\theta^3-2r^2)
+r^2(r^2+r_\theta^2)r_\theta\sin^2\theta\cos\theta\\
&+ r^2r_\theta r_\varphi^2\cos\theta, 
\end{align*}
which, combined with \eqref{ex:carv1}, leads to  \eqref{ex:mean}.

\section{Free boundary problem with surface tension} \label{sec:2}

In this section, we consider the case where $\sigma$ is a positive constant
in \eqref{navier:1}, and we shall transform Eq. \eqref{navier:1} to
some problem formulated on a fixed domain by using the Hanzawa transform.
\subsection{Hanzawa transform}\label{subsec:2.1}
As $\Omega_t$ is unknown, we have to transform $\Omega_t$ to a fixed
domain $\Omega$.  In the following, we assume that 
the reference domain $\Omega$ is a uniform $C^3$ domain. Let
$\Gamma$ be the boundary of $\Omega$ and $\bn$ the unit outer normal
to $\Gamma$.  We may assume that $\bn$ is an $N$ vector of $C^3$ function
defined on $\BR^N$ satisfying the condition: $\|\bn\|_{H^3_\infty(\BR^N)}
< \infty$. We assume that $\Gamma_t$ is given by
\begin{equation}\label{boundary:1}
\Gamma_t = \{x = y + \rho(y, t)\bn + \xi(t) \mid 
y \in \Gamma\} \quad(t \in (0, T))
\end{equation}
with an unknown function $\rho(x, t)$, where $\xi(t)$ is a function 
depending solely on $t$.  Let $H_\rho(y, t)$ be a suitable extension of
$\rho(y,t)$  to $\BR^N$ such that for each $t \in (0, T)$ 
$\rho(y, t) = H_\rho(y, t)$ for $y \in \Gamma$ and 
\begin{equation}\label{sob:5.4}\begin{split}
C_1\|H_\rho(\cdot, t)\|_{H^k_q(\BR^N)} 
&\leq \|\rho(\cdot, t)\|_{W^{k-1/q}_q(\Gamma)}\leq 
C_2\|H_\rho(\cdot, t)\|_{H^k_q(\BR^N)},  \\
C_1\|\pd_tH_\rho(\cdot, t)\|_{H^\ell_q(\BR^N)} 
&\leq \|\pd_t\rho(\cdot, t)\|_{W^{\ell-1/q}_q(\Gamma)}\leq 
C_2\|\pd_tH_\rho(\cdot, t)\|_{H^\ell_q(\BR^N)} 
\end{split}\end{equation}
with some constants $C_1$ and $C_2$ for $k=1,2,3$ and 
$\ell=1,2$. 
We then define $\Omega_t$ by
$$
\Omega_t = \{x = y+\omega(y)H_\rho(y, t)\bn(y) + \xi(t) \mid y \in \Omega\}
\quad(t \in (0, T)),
$$
where $\omega(y)$ is a $C^\infty$ function which equals $1$ near $\Gamma$
and zero far from $\Gamma$.
For the notational simplicity, we set 
$\Psi_\rho(y, t) = \omega(y)H_\rho(y, t)\bn(y)$.  The transform:
\begin{equation}\label{map:1}
x = y + \Psi_\rho(y, t) + \xi(t)
\end{equation}
is called the Hanzawa transform, which was originally introduced by Hanzawa 
\cite{Hanzawa} to
treat classical solutions of the Stefan problem. 
Usually, $\xi(t)$ is also unknown functions and to prove the 
local well-posedness, we set $\xi(t) = 0$. 

In the following, we study how to change the equations and boundary conditions
under such transformation. 
Assume that 
\begin{equation}\label{trans:1}
\sup_{t \in (0, T)} \|\Psi_\rho(\cdot, t)\|_{H^1_\infty(\BR^N)} \leq \delta,
\end{equation}
where $\delta$ is a small positive number determined in such a way that
several conditions stated below will be satisfied. We first choose
$0 < \delta < 1$, and then the Hanzawa transform \eqref{map:1}
is injective for each $t \in (0, T)$.  In fact, let $x_i=
y_i + \Psi_\rho(y_i, t) + \xi(t)$ ($i=1,2$). We then have
\begin{align*}
&|x_1-x_2|  = |y_1 + \Psi_\rho(y_1, t) - (y_2 + 
\Psi_\rho(y_2, t)| \\
&\geq |y_1-y_2| 
- \|\nabla\Psi_\rho(\cdot, t)\|_{L_\infty(\BR^N)}
|y_1-y_2| \geq (1-\delta)|y_1- y_2|.
\end{align*}
Thus, the condition $0 < \delta <1$ implies if $y_1\not=y_2$ then
$x_1\not=x_2$, that is,  the Hanzawa transform is injective. 

We now set 
$$
\Omega_t = \{x = y+ \Psi_\rho(y, t) + \xi(t) \mid y \in \Omega\}
\quad(t \in (0, T).
$$
Notice that the Hanzawa transform maps $\Omega$ onto
$\Omega_t$ injectively. 

\subsection{Transformation of equations and the divergence 
free condition}\label{sub:1}
In this subsection, for the latter use we consider more general
transformation: $x = y + \Psi(y, t)$, where $\Psi(y, t)$ satisfies
the condition:
$$\sup_{t \in (0, T)}\|\Psi(\cdot, t)\|_{H^1_\infty(\BR^N)} 
\leq \delta$$
with some small $\delta > 0$. Moreover, we assume that
$\pd_t\Psi$ exists. 

Let $\bv$ and $\fp$ be solutions of Eq.\eqref{navier:1}, 
and let 
$$\bu(y, t) = \bv(y+\Psi(y,t), t), 
\quad \fq(y,t)
= \fp(y+\Psi(y, t), t).
$$
We now show 
that the first equation in \eqref{navier:1} is 
transformed to 
\begin{equation}\label{newchange:1} \pd_t\bu- 
\DV(\mu\bD(\bu) - \fq\bI)
= \bff(\bu, \Psi)
\quad\text{in $\Omega^T$},
\end{equation}
and the divergence free condition: $\dv \bv = 0$ 
in \eqref{navier:1} is transformed to
\begin{equation}\label{change:7}
\dv\bu = g(\bu, \Psi) = \dv\bg(\bu, \Psi) \quad \text{in 
$\Omega^T$}.  
\end{equation}
Where, $\bff(\bu, \Psi)$, $g(\bu,\Psi)$ and 
$\bg(\bu, \Psi)$ are suitable non-linear functions with
respect to $\bu$ and $\nabla\Psi$ given in \eqref{form:f} and 
\eqref{form:g}, below.

Let $\pd x/\pd y$ be the Jacobi matrix of the transformation
\eqref{map:1}, that is, 
$$\frac{\pd x}{\pd y} = \bI + \nabla\Psi(y, t) = (\delta_{ij}
+ \frac{\pd \Psi_{j}}{\pd y_i}),
\quad
\nabla\Psi = \left(\begin{matrix}
\pd_1\Psi_{1} & \pd_2\Psi_{1} & \ldots & \pd_N\Psi_{1} \\
\pd_1\Psi_{2} & \pd_2\Psi_{2} 
& \ldots & \pd_N\Psi_{2} \\
\vdots & \vdots &\ddots & \vdots \\
\pd_1\Psi_{N} & \pd_2\Psi_{N} & \ldots & \pd_N\Psi_{N}
\end{matrix}\right)
$$
where $\Psi(y, t) = {}^\top(\Psi_{1}(y, t), 
\ldots, \Psi_{N}(y, t))$, and 
$\pd_i\Psi_{j} = \dfrac{\pd \Psi_{j}}{\pd y_i}$. 
If $0 < \delta < 1$, then 
$$\Bigl(\frac{\pd x}{\pd y}\Bigr)^{-1} = \bI + \sum_{k=1}^\infty
(-\nabla\Psi(y, t))^k$$
exists, and therefore there exists an $N\times N$
 matrix $\bV_0(\bk)$ of $C^\infty$
functions defined on $|\bk| < \delta$ such that $\bV_0(0) = 0$ and 
\begin{equation}\label{trans:2}
\Bigl(\frac{\pd x}{\pd y}\Bigr)^{-1} = \bI + \bV_0(\nabla\Psi(y, t)).
\end{equation}
Here and in the following, $\bk = (k_{ij})$ and $k_{ij}$ are the 
variables corresponding to $\pd_i\Psi_{j}$.  

Let $V_{0ij}(\bk)$ be the $(i, j)^{\rm th}$ component of $\bV_0(\bk)$.
We then have 
\begin{equation}\label{change:1}
\nabla_x = (\bI + \bV_0(\bk))\nabla_y, \quad
\frac{\pd}{\pd x_i} = \sum_{j=1}^N(\delta_{ij} + V_{0ij}(\bk))
\frac{\pd}{\pd y_j},
\end{equation}
where $\nabla_z = {}^\top(\pd/\pd z_, \ldots, \pd/\pd z_N)$ 
for $z = x$ and $y$.
By \eqref{change:1}, we can write $\bD(\bv)$ as  
$\bD(\bv) = \bD(\bu) + \CD_\bD(\bk)\nabla\bu$ with 
\begin{equation} \label{div:1}\begin{split}
\bD(\bu)_{ij} &= \frac{\pd u_i}{\pd y_j} + \frac{\pd u_j}{\pd y_i},
\\
(\CD_\bD(\bk)\nabla\bu)_{ij} &= \sum_{k=1}^N
\Bigl(V_{0jk}(\bk)\frac{\pd u_i}{\pd y_k}
+ V_{0ik}(\bk)\frac{\pd u_j}{\pd y_k}\Bigr).
\end{split}\end{equation}
We next consider $\dv\bv$.  By \eqref{change:1}, we have
\begin{equation}\label{div:1-1}\dv_x\bv = \sum_{j=1}^N\frac{\pd v_j}{\pd x_j}
= \sum_{j,k=1}^N(\delta_{jk} + V_{0jk}(\bk))\frac{\pd u_j}{\pd y_k}
= \dv_y\bu + \bV_0(\bk):\nabla \bu.
\end{equation}
Let $J$ be the Jacobian of the transformation \eqref{map:1}.  
Choosing $\delta > 0$ small enough, we may assume 
that $J = J(\bk) = 1 + J_0(\bk)$, 
where $J_0(\bk)$ is a $C^\infty$ function defined for $|\bk| < \sigma$ such 
that $J_0(0) = 0$.  

To obtain another representation formula of $\dv_x\bv$,
we use the inner product $(\cdot, \cdot)_{\Omega_t}$.  For any 
test function $\varphi \in C^\infty_0(\Omega_t)$, 
we set $\psi(y) = \varphi(x)$. We then have 
\begin{align*}
&(\dv_x\bv, \varphi)_{\Omega_t} = -(\bv, \nabla\varphi)_{\Omega_t}
= -(J\bu, (\bI + \bV_0)\nabla_y\psi)_\Omega \\
&= (\dv((\bI + {}^\top\bV_0)J\bu), \psi)_\Omega
= (J^{-1}\dv((\bI + {}^\top\bV_0)J\bu), \varphi)_{\Omega_t},
\end{align*}
which, combined with \eqref{div:1-1}, leads to 
\begin{equation}\label{div:2}
\dv_x\bv = \dv_y\bu + \bV_0(\bk):\nabla\bu
= J^{-1}(\dv_y\bu + \dv_y(J{}^\top\bV_0(\bk)\bu)).
\end{equation}
Recalling that $J = 1 + J_0(\bk)$, we define $g(\bu, \Psi)$ and
$\bg(\bu, \Psi)$ by letting  
\begin{equation}\label{form:g}\begin{split}
g(\bu, \Psi) &= -(J_0(\bk)\dv\bu + (1+J_0(\bk))\bV_0(\bk):\nabla\bu), \\
\bg(\bu, \Psi) &= -(1+J_0(\bk)){}^\top\bV_0(\bk)\bu,
\end{split}\end{equation}
and then by \eqref{div:2} we see that 
the divergence free condition: $\dv\bv=0$ is transformed to
 Eq.\eqref{change:7}. 
In particular, it follows from \eqref{div:2} that 
\begin{equation}\label{div:3}
J_0\dv\bu + J\bV_0(\bk):\nabla\bu
= \dv(J{}^\top\bV_0(\bk)\bu).
\end{equation}

To derive  Eq. \eqref{newchange:1}, 
we first observe that
\begin{align}
&\sum_{j=1}^N\frac{\pd}{\pd x_j}(\mu\bD(\bv)_{ij} - \fp \delta_{ij})
\nonumber \\
&= \sum_{j,k=1}^N\mu(\delta_{jk} + V_{0jk})\frac{\pd}{\pd y_k}
(\bD(\bu)_{ij} + (\CD_\bD(\bk)\nabla\bu)_{ij})
-\sum_{j=1}^N(\delta_{ij} + V_{0ij})\frac{\pd \fq}{\pd y_j},
\label{change:5}
\end{align}
where we have used \eqref{div:1}.  Since 
$$\frac{\pd}{\pd t}[v_i(y + \Psi(y, t), t)]
= \frac{\pd v_i}{\pd t}(x, t) + 
\sum_{j=1}^N\frac{\pd\Psi_{j}}{\pd t}
\frac{\pd v_i}{\pd x_j}(x, t),
$$
we have 
$$\frac{\pd v_i}{\pd t} = \frac{\pd u_i}{\pd t}-\sum_{j,k=1}^N
\frac{\pd \Psi_{j}}{\pd t}
(\delta_{jk} + V_{0jk})\frac{\pd u_i}{\pd y_k},
$$
and therefore, 
\begin{equation}\label{change:3}
\frac{\pd v_i}{\pd t} + \sum_{j=1}^N v_j\frac{\pd v_i}{\pd x_j}
= \frac{\pd u_i}{\pd t}
 + \sum_{j,k=1}^N(u_j - \frac{\pd\Psi_{j}}{\pd t})
(\delta_{jk} + V_{0jk}(\bk))\frac{\pd u_i}{\pd y_k}.
\end{equation}
Putting \eqref{change:5} and \eqref{change:3} together gives 
\begin{align*}
0 & = \Bigl(\frac{\pd u_i}{\pd t} 
+ \sum_{j,k=1}^N(u_j - \frac{\pd\Psi_{j}}{\pd t})
(\delta_{jk} + V_{0jk}(\bk))\frac{\pd u_i}{\pd y_k}\Bigr) \\
&\quad - \mu\sum_{j,k=1}^N(\delta_{jk} + V_{0jk}(\bk))
\frac{\pd}{\pd y_k}(\bD(\bu)_{ij} + (\CD_\bD(\bk)\nabla\bu)_{ij}) \\
&\quad - \sum_{j=1}^N(\delta_{ij}+ V_{0ij}(\bk))\frac{\pd\fq}{\pd y_j}.
\end{align*}
Since $(\bI + \nabla\Psi)(\bI + \bV_0) = (\pd x/\pd y)(\pd y/\pd x) = \bI$, 
\begin{equation}\label{change:3*}
\sum_{i=1}^N(\delta_{mi} + \pd_m\Psi_{i})
(\delta_{ij} + V_{0ij}(\bk)) = \delta_{mj},
\end{equation}
and so we have
\begin{align*}
0=&\sum_{i=1}^N(\delta_{mi} + \pd_m\Psi_{i})
\Bigl(\frac{\pd u_i}{\pd t} 
+ \sum_{j,k=1}^N(u_j - \frac{\pd\Psi_{i}}{\pd t})
(\delta_{jk}+V_{0jk}(\bk))\frac{\pd u_i}{\pd y_k})\Bigr)
\\
&-\mu\sum_{i,j,k=1}^N(\delta_{mi} + \pd_m\Psi_{i})
(\delta_{jk}+V_{0jk}(\bk))\frac{\pd}{\pd y_k}
(\bD(\bu)_{ij} +(\CD_\bD(\bk)\nabla\bu)_{ij})\\
&-\frac{\pd \fq}{\pd y_m}.
\end{align*}
Thus, changing $i$ to $\ell$ and $m$ to $i$ in the formula above, 
we define an $N$-vector
of functions $\bff(\bu, \Psi)$ by letting 
\allowdisplaybreaks{
\begin{align}
&\bff(\bu, \Psi)|_i  = -\sum_{j,k=1}^N(u_j 
- \frac{\pd\Psi_{j}}{\pd t})
(\delta_{jk} + V_{0jk}(\bk))\frac{\pd u_i}{\pd y_k} \nonumber\\
&\quad-\sum_{\ell=1}^N\pd_i\Psi_{\ell}
\Bigl(\frac{\pd u_\ell}{\pd t} 
+ \sum_{j,k=1}^N(u_j - \frac{\pd\Psi_{j}}{\pd t})
(\delta_{jk} + V_{0jk}(\bk))\frac{\pd u_\ell}{\pd y_k}\Bigr) \nonumber \\
&\quad + \mu\Bigl(\sum_{j=1}^N \frac{\pd}{\pd y_j}(\CD_\bD(\bk)\nabla\bu)_{ij}
+ \sum_{j,k=1}^NV_{0jk}(\bk)\frac{\pd}{\pd y_k}(\bD(\bu)_{ij}
+ (\CD_\bD(\bk)\nabla\bu)_{ij})  \nonumber \\ 
&\quad + \sum_{j,k,\ell=1}^N\pd_i\Psi_{\ell}(\delta_{jk} + V_{0jk}(\bk))
\frac{\pd}{\pd y_k}(\bD(\bu)_{\ell j} + (\CD_\bD(\bk)\nabla\bu)_{\ell j})
\Bigr), \label{form:f}
\end{align}
}
where $\bff(\bu, \Psi)|_i$ denotes the $i^{\rm th}$ complonent of 
$\bff(\bu, \Psi)$.  We then see that Eq. \eqref{eq:momentum} is 
transformed to Eq.\eqref{newchange:1}.

\subsection{Transformation of the boundary conditions}
\label{subsec:2.3}

In this subsection, we consider the Hanzawa transform
given in Subsec.\ref{subsec:2.1}. 
To represent $\Gamma$ locally, we use the local coordinates 
near $x^1_\ell \in \Gamma$ such that 
\begin{equation}\label{loc:1}\begin{aligned}
\Omega \cap B^1_\ell & = \{y = \Phi_\ell(p) \mid
p \in \BR^N_+\} \cap B^1_\ell, \\
\Gamma \cap B^1_\ell & = 
\{y = \Phi_\ell(p', 0) \mid 
(p', 0) \in \BR^N_0\} \cap B^1_\ell. 
\end{aligned}\end{equation} 
Let $\{\zeta^1_\ell\}_{\ell\in \BN}$
be a partition of unity given in Proposition \ref{prop:lap}.
In the following we use the formula:
$$f = \sum_{\ell=1}^\infty \zeta^1_\ell f \quad\text{in $\Gamma$}$$
for any function, $f$, defined on $\Gamma$.

We write $\rho = \rho(y(p_1, \ldots, p_{N-1}, 0), t)$
in the following. 
By the chain rule, we have
\begin{equation}\label{chain:1}
\frac{\pd \rho}{\pd p_i} 
= \frac{\pd}{\pd p_i}\Psi_\rho(\Phi_{\ell}(p_1, \ldots, p_{N-1}, 0), t)
= \sum_{m=1}^N\frac{\pd \Psi_\rho}{\pd y_m}\frac{\pd \Phi_{\ell,m}}{\pd p_i}|_{p_N=0},
\end{equation}
where we have set $\Phi_\ell = {}^\top(\Phi_{\ell,1}, \ldots, \Phi_{\ell,N})$, 
and so, $\pd \rho/\pd p_i$ is defined in $B^1_\ell$ 
by letting
\begin{equation}\label{chain:2}
\frac{\pd \rho}{\pd p_i} 
= \sum_{m=1}^N\frac{\pd \Psi_\rho}{\pd y_m}\circ\Phi_\ell
\frac{\pd \Phi_{\ell, m}}{\pd p_i}.
\end{equation}


We first represent $\bn_t$. Recall that $\Gamma_t$ is given by
$x = y + \rho(y, t)\bn + \xi(t)$ for $y \in \Gamma$ (cf. \eqref{boundary:1}).
 Let 
$$\bn_t = a(\bn + \sum_{i=1}^{N-1}b_i\tau_i)
\quad\text{with $\tau_i = \frac{\pd}{\pd p_i}y 
=\frac{\pd}{\pd p_i}\Phi_\ell(p', 0)$}.
$$
These vectors $\tau_i$ ($i=1, \ldots, N-1$) 
form a basis of the tangent space of $\Gamma$ at $y=y(p_1,
\ldots, y_{N-1})$.
Since $|\bn_t|^2 =1$, we have
\begin{equation}\label{norm:1}
1  = a^2(1 + \sum_{i,j=1}^{N-1}g_{ij}b_ib_j)
\quad\text{with $g_{ij} = \tau_i\cdot\tau_j$}
\end{equation}
because $\tau_i\cdot\bn=0$. The vectors $\dfrac{\pd x}{\pd p_i}$
$(i=1, \ldots, N-1)$ form a basis of the tangent space of $\Gamma_t$,
and so  $\dfrac{\pd x}{\pd p_i} \cdot \bn_t=0$.  Thus, we have
\begin{equation}\label{normal:3.1.1}
0 = a(\bn + \sum_{j=1}^{N-1}b_j\tau_j)\cdot
(\frac{\pd y}{\pd p_i} + \frac{\pd \rho}{\pd p_i}\bn 
+ \rho\frac{\pd \bn}{\pd p_i}).
\end{equation}
Since $\bn\cdot\dfrac{\pd y}{\pd p_i} = \bn\cdot\tau_i= 0$,
$\dfrac{\pd \bn}{\pd p_i}\cdot\bn
= 0$ (because of $|\bn|^2=1$), and $\dfrac{\pd y}{\pd p_i}\cdot
\dfrac{\pd y}{\pd p_j} = \tau_i\cdot\tau_j=g_{ij}$,
by \eqref{normal:3.1.1} we have
\begin{equation}\label{norm:2}
\frac{\pd \rho}{\pd p_i} + \sum_{i=1}^{N-1}(g_{ij}
+ \rho\frac{\pd \bn}{\pd p_i}\cdot\tau_j)b_j = 0.
\end{equation} 
Let $H$ be an $(N-1)\times(N-1)$ matrix
whose $(i, j)^{\rm th}$ component is $\dfrac{\pd \bn}{\pd p_i}\cdot\tau_j$.
Since $\bn$ is defined in $\BR^N$ as an $N$ vector of $C^2$ functions with
$\|\bn\|_{H^2_\infty(\BR^N)} < \infty$,  we can write 
\begin{equation}
\frac{\pd \bn}{\pd p_i} \cdot \tau_j 
= \sum_{m=1}^N\frac{\pd \bn}{\pd y_m}\circ\Phi_\ell
\frac{\pd \Phi_{\ell, m}}{\pd p_i}\cdot
\frac{\pd \Phi_\ell}{\pd p_j}, \label{chain:3}
\end{equation}
 and then $H$ is defined in $\BR^N$ and  $\|H\|_{H^2_\infty(\BR^N)} 
\leq C$ with some constant $C$ independent of $\ell \in \BN$. 
Under the assumption \eqref{trans:1}, we may assume that the inverse
of $\bI + \rho HG^{-1}$ exists. Let 
$\nabla'_\Gamma \rho = (\dfrac{\pd \rho}{\pd p_1},\ldots, 
\dfrac{\pd \rho}{\pd p_{N-1}})$ with $\dfrac{\pd \rho}{\pd p_i}
= \dfrac{\pd}{\pd p_i}\rho(\Phi_\ell(p', 0))$, and then 
by \eqref{norm:2} we have 
\begin{equation}\label{norm:3}
b = -(G+\rho H)^{-1}\nabla'_\Gamma\rho
= -G^{-1}(\bI + \rho HG^{-1})^{-1}\nabla'_\Gamma\rho.
\end{equation}
Putting \eqref{norm:3} and \eqref{norm:1} together gives
\begin{align*}
a & = (1 + <Gb, b>)^{-1/2}\\
&= (1+<\bI + \rho HG^{-1})^{-1}\nabla'_\Gamma\rho, 
G^{-1}(\bI + \rho HG^{-1})^{-1}\nabla'_\Gamma\rho>)^{-1/2}.
\end{align*}
Thus, we have
\begin{equation}\label{repr:2.1} \bn_t = \bn
-\sum_{i,j=1}^{N-1}g^{ij}\tau_i\frac{\pd \rho}{\pd p_j} 
+\bV_\Gamma(\rho, \nabla'_\Gamma\rho)
\end{equation}
where we have set 
\begin{align*}
&\bV_\Gamma(\rho, \nabla'_\Gamma\rho) 
= -<G^{-1}((\bI+\rho H G^{-1})^{-1} -\bI)\nabla'_\Gamma\rho, \tau>
\\
&\quad + \{(1+ <(\bI+\rho HG^{-1})^{-1}\nabla'_\Gamma\rho, 
G^{-1}(\bI + \rho HG^{-1})^{-1}\nabla'_\Gamma\rho>)^{-1/2}-1\}\\
&\phantom{+ \{(1+ <(\bI+\rho HG^{-1})^{-1}\nabla'_\Gamma\rho, }\times
(\bn-<G^{-1}(\bI + \rho HG^{-1})^{-1}\nabla'_\Gamma\rho,\tau>).
\end{align*}
From \eqref{chain:2}, $\nabla'_\Gamma\rho$ is extended to $\BR^N$ 
by letting $\nabla'_\Gamma\rho = (\nabla \Phi_\ell)\nabla 
\Psi_\rho\circ\Phi_\ell:=\Xi_{\rho, \ell}$,
and so $\bV_\Gamma(\rho, \nabla'_\Gamma\rho)$ can be extended to 
$B^1_\ell$ by 
\begin{equation}\label{repr:2.1*}\begin{split}
&\bV_\Gamma(\rho, \nabla'_\Gamma\rho)
= -<G^{-1}((\bI+\Psi_\rho H G^{-1})^{-1} -\bI)
\Xi_{\rho, \ell}, \tau>
\\
& + \{(1+ <(\bI+\Psi_\rho HG^{-1})^{-1}\Xi_{\rho, \ell}, 
G^{-1}(\bI + \rho HG^{-1})^{-1}\Xi_{\rho, \ell}>)^{-1/2}-1\}\\
&\quad \quad \times
(\bn-<G^{-1}(\bI + \rho HG^{-1})^{-1}
\Xi_{\rho, \ell},\tau>), 
\end{split}\end{equation}
where $H$ is an $(N-1)\times(N-1)$ matrix whose $(i, j)^{\rm th}$
component, $H_{ij}$,  is given by 
$$H_{ij} = \sum_{m=1}^N\frac{\pd \bn}{\pd y_m}\circ\Phi_\ell
\frac{\pd \Phi_{\ell, m}}{\pd p_i}\cdot
\frac{\pd \Phi_\ell}{\pd p_j}.$$ 
Thus, we may write 
\begin{align*}
\bV_\Gamma(\rho, \nabla'_\Gamma\rho)
= \bV_{\Gamma, \ell}(\bar\bk)\bar\nabla \Psi_\rho\otimes\bar\nabla\Psi_\rho
\end{align*}
on $B^1_\ell$ with some function $\bV_{\Gamma, \ell}(\bar\bk)
= \bV_{\Gamma, \ell}(y, \bar\bk)$ defined on 
$B^1_\ell\times\{\bar\bk \mid |\bar\bk| \leq \delta\}$ 
with $\bV_{\Gamma, \ell}(0) = 0$ possessing
the estimate
$$\|(\bV_{\Gamma, \ell}(\bar\bk), \pd_{\bar\bk}
\bV_{\Gamma, \ell}(\cdot, \bk))
\|_{H^1_\infty(B^1_\ell)} \leq C
$$
with some constant $C$ independent of $\ell$.  Here and 
in the following $\bar\bk$ are the corresponding variables
to $\bar\nabla\Psi_\rho = (\Psi_\rho, \nabla\Psi_\rho)$. 
In view of \eqref{repr:2.1}, we have 
\begin{equation}\label{repr:2.1*}
\bn_t = \bn - \sum_{i,j=1}^{N-1}g^{ij}\tau_i\frac{\pd\rho}{\pd y_j}
+\bV_{\Gamma, \ell}(\bar\bk)\bar\nabla\Psi_\rho\otimes\bar\nabla\Psi_\rho
\quad\text{on $B^1_\ell \cap \Gamma$}.
\end{equation}
Let 
\begin{equation}\label{repr:2.1**}
\bV_\Gamma(\bar\bk) = \sum_{\ell=1}^\infty
\zeta^1_\ell \bV_{\Gamma, \ell}(\bar\bk),
\end{equation}
and then  we have
\begin{equation}\label{repr:2.1***}
\|(\bV_{\Gamma}(\bar\bk), \pd_{\bar\bk}\bV_\Gamma(\cdot, \bk))
\|_{H^1_\infty(\Omega)} \leq C
\end{equation}
for $|\bar\bk| \leq \delta$ with some constant $C>0$.

In view of \eqref{chain:2} and \eqref{repr:2.1}, 
the unit outer normal $\bn_t$ is also represented by 
$$\bn_t = \bn - \sum_{i, j=1}^{N-1}
\sum_{m=1}^Ng^{ij}\tau_i
\frac{\pd\Psi_\rho}{\pd y_m}\circ\Phi_\ell\frac{\pd \Phi_{\ell, m}}{\pd p_i}
+ \bV_{\Gamma, \ell}(\bar\bk)\bar\nabla\Psi_\rho\otimes\bar\nabla\Psi_\rho
$$
on $B^1_\ell$, and so we may write 
$$\bn_t = \bn + 
\tilde\bV_{\Gamma, \ell}(\bar\nabla\Psi_\rho)\bar\nabla\Psi_\rho
$$
for some functions $\tilde\bV_{\Gamma, \ell}(\bar\bk)
= \tilde\bV_{\Gamma, \ell}(y, \bar\bk)$ defined on 
$B^1_\ell\times\{ \bar\bk \mid |\bar\bk| \leq \delta\}$
possessing the estimate: 
\begin{equation}\label{normal:3.4}
\|(\tilde\bV_{\Gamma,\ell}, 
\pd_{\bar\bk}\tilde \bV_{\Gamma, \ell})(\cdot, \bar\kappa))
\|_{H^1_\infty(B^1_\ell)}
\leq C \quad\text{for $|\bar\kappa| \leq \delta$}
\end{equation}
with some constant $C$ independent of $\ell \in \BN$.  Thus, setting
\begin{equation}\label{normal:3.3}
\tilde \bV_\Gamma(\bar\bk) = 
\sum_{\ell=1}^\infty \zeta^1_\ell
\tilde\bV_{\Gamma, \ell}(\bar\bk),
\end{equation}
we have 
\begin{equation}\label{normal:3.1}
\bn_t = \bn + \tilde\bV_\Gamma(\bar\nabla\Psi_\rho)\bar\nabla\Psi_\rho
\end{equation}
and 
\begin{equation}\label{normal:3.2}
\|(\tilde\bV_\Gamma(\cdot, \bar\bk), \pd_{\bar\bk}\bV_\Gamma(\cdot, 
\bar\bk))\|_{H^1_\infty(\Omega)}
\leq C \quad\text{for $|\bar\kappa| \leq \delta$}.
\end{equation}

We now consider the kinematic equation: $V_{\Gamma_t} = \bv\cdot\bn_t$ in 
\eqref{navier:1}.  
Since $x = y + \rho(y, t)\bn + \xi(t)$, by \eqref{repr:2.1*} we have
\begin{align*}
V_{\Gamma_t} &=
\frac{\pd x}{\pd t}\cdot\bn_t\\
& = \sum_{\ell=1}^\infty\zeta^1_\ell
<(\pd_t\rho)\bn + \xi'(t), 
\bn - \sum_{i,j=1}^{N-1}g^{ij}\tau_i\frac{\pd\rho}{\pd p_j}
+ \bV_{\Gamma, \ell}(\bar\bk)\bar\nabla\Psi_\rho\otimes\bar\nabla\Psi_\rho> \\
& = \pd_t\rho + \xi'(t)\cdot\bn
+ \pd_t\rho<\bn, \bV_\Gamma(\bar\bk)
\bar\nabla\Psi_\rho\otimes\bar\nabla\Psi_\rho>\\
&- < \xi'(t) \mid \bar\nabla'_\Gamma \rho> 
+ <\xi'(t), \bV_\Gamma(\bar\bk)\bar\nabla\Psi_\rho\otimes
\bar\nabla\Psi_\rho>.
\end{align*}
Where, for any $N$-vector function $\bd$, we have set 
\begin{equation}\label{inner:1}
<\bd \mid \nabla'_\Gamma\rho> = \sum_{i,j=1}^{N-1}
g^{ij}<\tau_i, \bd>\frac{\pd \rho}{\pd p_j}.
\end{equation}
On the other hand, 
\begin{align*}
\bv\cdot\bn_t &= \sum_{\ell=1}^\infty\zeta^1_\ell \bu\cdot(\bn - 
\sum_{i,j=1}^{N-1}g^{ij}\tau_i\frac{\pd \rho}{\pd p_j}
+ \bV_{\Gamma, \ell}(\bar\bk)\bar\nabla\Psi_\rho\otimes\bar\nabla\Psi_\rho)\\
& = \bn\cdot\bu - <\bu \mid \nabla'_\Gamma\rho> + \bu\cdot\bV_\Gamma(\bar
\bk)\bar\nabla\Psi_\rho\otimes\bar\nabla\Psi_\rho. 
\end{align*}
From the consideration above, 
 the kinematic equation is transformed to the following equation:
\begin{equation}\label{kinematic:4}
\pd_t\rho +\xi'(t)\cdot\bn - \bn\cdot\bu +
<\bu \mid \nabla'_\Gamma\rho> 
= d(\bu, \rho) \quad\text{on $\Gamma\times(0, T)$}
\end{equation}
with
\begin{equation}\label{kinematic:4*}\begin{split}
d(\bu, \rho) &= \bu\cdot\bV_\Gamma(\bar\bk)\bar\nabla\Psi_\rho\otimes
\bar\nabla\Psi_\rho 
-\pd_t\rho<\bn, 
\bV_\Gamma(\bar\bk)\bar\nabla\Psi_\rho\otimes\bar\nabla\Psi_\rho>
\\
& + < \xi'(t) \mid \nabla'_\Gamma\rho> 
- <\xi'(t), \bV_\Gamma(\bar\bk)\bar\nabla\Psi_\rho\otimes\bar\nabla\Psi_\rho>.
\end{split}\end{equation}

We next consider the boundary condition:
\begin{equation}\label{52}
(\mu\bD(\bv) - \fp\bI)\bn_t  = \sigma H(\Gamma_t)\bn_t
\end{equation}
in Eq.\eqref{navier:1}. It is convenient  
to divide the formula in \eqref{52} into 
the tangential part and normal part on $\Gamma_t$
as follows:
\begin{gather}
\BPi_t\mu\bD(\bv)\bn_t = 0, \label{55} \\
<\mu\bD(\bv)\bn_t, \bn_t> - \fp 
= \sigma <H(\Gamma_t)\bn_t, \bn_t>.
\label{56}
\end{gather}
Where, $\BPi_t$ is defined by 
$\BPi_t\bd = \bd- < \bd, \bn_t>\bn_t$ for any 
$N$ vector $\bd$.

By \eqref{normal:3.1} we can write 
\begin{align}
\BPi_t\bd &= \BPi_0\bd 
+<\tilde \bV_{\Gamma}(\bar\nabla\Psi_\rho)\bar\nabla\Psi_\rho, \bd>\bn 
+ <\bn, \bd>\tilde \bV_{\Gamma}(\bar\nabla\Psi_\rho)
\bar\nabla\Psi_\rho \nonumber \\
&\,\,+ <\tilde \bV_{\Gamma}(\bar\nabla\Psi_\rho)\bar\nabla\Psi_\rho,\bd>
\tilde \bV_{\Gamma}(\bar\nabla\Psi_\rho)\bar\nabla\Psi_\rho,
\label{58**}
\end{align}
for any $\bd \in \BR^N$.
Thus, 
recalling \eqref{div:1}, we see that 
the boundary condition \eqref{56} is transformed to the following 
formula: 
\begin{equation}\label{58}
(\mu\bD(\bu)\bn)_\tau = \bh'(\bu, \Psi_\rho)
\quad\text{on $\Gamma\times(0, T)$}, 
\end{equation}
where we have set 
\begin{align}
\bh'(\bu, \Psi_\rho) 
& = - \mu
\{\BPi_0 \CD_\bD(\bk)\nabla \bu 
+ <\tilde \bV_{\Gamma}(\bar\bk)\bar\bk,\bD(\bu) 
+ \CD_\bD(\bk)\nabla\bu >\bn 
\nonumber \\
&\quad\,\,
+<\bn, \bD(\bu) + \CD_\bD(\bk)\nabla\bu>
\tilde\bV_{\Gamma}(\bar\bk)\bar\bk
\nonumber \\
&\quad\,\,+ <\tilde\bV_{\Gamma}(\bar\bk)\bar\bk, \bD(\bu) + 
\CD_\bD(\bk)\nabla\bu>\tilde\bV_{\Gamma,\ell}(\bar\bk)\bar\bk\}
\label{58*}
\end{align}
with $\bk = \nabla\Psi_\rho$ and $\bar\bk=(\Psi_\rho, \nabla\Psi_\rho)$. 
Here and in the following, for any $N$-vector $\ba$, we set
$$\ba_\tau = \ba - <\ba, \bn>\bn.$$

To consider the transformation of the 
 boundary condition \eqref{56}, 
we first consider the Laplace-Beltrami operator on
 $\Gamma_t$. From \eqref{cg:10}, we have
\begin{equation}\label{LB:4}
\Delta_{\Gamma_t}f = g^{ij}_t\pd_i\pd_jf
-g_t^{ij}\Lambda^j_{tik}\pd_jf,
\end{equation}
where $\Lambda^j_{tik}$ are Christoffel symbols of $\Gamma_t$ 
defined by \eqref{cg:5}.
Recall $x = y + \rho\bn+\xi(t)$.   The first fundamental form of $\Gamma_t
=(g_{tij})$ is 
given by
\begin{equation}\label{first:1}
g_{tij} =<\frac{\pd x}{\pd x_i}, \frac{\pd x}{\pd x_j}>
= g_{ij} + \tilde g_{ij}\rho + \frac{\pd\rho}{\pd p_i}\frac{\pd \rho}{\pd p_j}
\end{equation}
with 
$$\tilde g_{ij} = <\tau_i, \frac{\pd\bn}{\pd p_j}>
+ <\tau_j, \frac{\pd \bn}{\pd p_i}>
+ < \frac{\pd\bn}{\pd p_i}, \frac{\pd \bn}{\pd p_j}>.
$$
In view of \eqref{loc:1}, we can write $\tau_i$ and 
$\pd_j\bn$ as 
\begin{equation}\label{local:3.1}
\tau_i = \frac{\pd y}{\pd p_i} 
= \sum_{j=1}^{N-1}\frac{\pd \Phi_\ell(p)}{\pd p_j}, \quad
\frac{\pd \bn}{\pd p_j}
= \sum_{k=1}^N\frac{\pd \bn}{\pd y_k}\frac{\pd \Phi_{\ell k}(p)}{\pd p_j}.
\end{equation}
Let $H$ be an $(N-1)\times(N-1)$ matrix whose $(i, j)^{\rm th}$
components are $\tilde g_{ij}$, and then 
$$G_t = G+\rho H
+ \nabla_p\rho\otimes\nabla_p\rho=G(\bI + \rho G^{-1}H
+ G^{-1}\nabla_p\rho\otimes\nabla_p\rho).$$
Where, $\nabla_p\rho = (\frac{\pd \rho}{\pd p_1},
\ldots, \frac{\pd \rho}{\pd p_{N-1}})$. 
Choosing $\delta > 0$ small enough in \eqref{trans:1}, we know
that the inverse of $\bI + \rho G^{-1}H + \nabla_p\rho\otimes
\nabla_p\rho$ exists, and so
$$G_t^{-1} = (\bI+\rho G^{-1}H
+ G^{-1}\nabla_p\rho\otimes\nabla_p\rho))^{-1}G^{-1}
= G^{-1} - \rho G^{-1}HG^{-1} + O_2,
$$
that is,
$$g^{ij}_t = g^{ij} + \rho h^{ij} + O_2,$$
where $h^{ij}$ denotes the $(i, j)^{\rm th}$ component of 
$-G^{-1}HG^{-1}$. 
Here and in the following, 
$O_2$ denotes some nonlinear 
function with respect to $\rho$ and $\nabla_p\rho$ of the form:
\begin{equation}\label{residue:1}
O_2 = a_0\Psi_\rho^2 + \sum_{j=1}^Na_j\Psi_\rho\frac{\pd \Psi_\rho}
{\pd y_j} + \sum_{j,k=1}^Nb_{jk}
\frac{\pd \Psi_\rho}{\pd y_j}\frac{\pd \Psi_\rho}{\pd y_k}
\end{equation}
with suitable functions $a_j$, and $b_{jk}$ possessing the 
estimates
\begin{equation}\label{residue:2}
|(a_j, b_{jk})| \leq C, \quad 
|\nabla(a_j, b_{jk})| \leq C(|\nabla\Psi_\rho| + 
|\nabla^2\Psi_\rho|), \end{equation}
provided that \eqref{trans:1} holds. Moreover, 
the Christoffel symbols are given by 
$\Lambda^k_{tij} = g_t^{k\ell}<\tau_{tij}, \tau_{t\ell}>$,
where $\tau_{ti} = \frac{\pd x}{\pd i}$ and $\tau_{tij}
= \frac{\pd^2x}{\pd p_i\pd p_J}$.
Since
\begin{align*}
\tau_{tij} = \tau_{ij} + \pd_i\pd_j(\rho\bn)
=\tau_{ij} + \rho\pd_i\pd_j\bn + (\pd_i\rho)\pd_j\bn + 
(\pd_j\rho)\pd_i\bn + (\pd_i\pd_j\rho)\bn,
\end{align*}
we have
\begin{align*}
<\tau_{tij}, \tau_{t\ell}>
&=<\tau_{ij}, \tau_\ell> + \rho(<\pd_i\pd_j\bn, \tau_\ell>
+ <\tau_{ij}, \pd_\ell\bn>) \\ 
& +\pd_j\rho<\pd_i\bn, \tau_\ell> 
+ \pd_i\rho<\pd_j\bn, \tau_\ell>.
\end{align*}
Thus, 
\begin{align*}
\Lambda^k_{tij} &= <g^{k\ell}+h^{k\ell}\rho + O_2, 
<\tau_{ij}, \tau_\ell> + \rho(<\pd_i\pd_j\bn, \tau_\ell>
+ <\tau_{ij}, \pd_\ell\bn>) \\ 
&\quad +\pd_j\rho<\pd_i\bn, \tau_\ell> 
+ \pd_i\rho<\pd_j\bn, \tau_\ell>\\
& = \Lambda^k_{ij} + \rho g^{k\ell}
(<\pd_i\pd_j\bn, \tau_\ell>
+ <\tau_{ij}, \pd_\ell\bn>)\\
&\quad + \pd_j\rho\, g^{k\ell}<\pd_i\bn, \tau_\ell>
+ \pd_i\rho\, g^{k\ell}<\pd_j\bn, \tau_\ell>
+ O_2,
\end{align*}
and so we may write 
\begin{equation}\label{chris:1}
\Lambda^k_{tij} = \Lambda^k_{ij} + \rho A^k_{ij}
+ \pd_i\rho B^k_j + \pd_j\rho B^k_i + O_2
\end{equation}
with
\begin{align*}
A^k_{ij} & = g^{k\ell}
(<\pd_i\pd_j\bn, \tau_\ell>
+ <\tau_{ij}, \pd_\ell\bn>), \\
B^k_j & = g^{k\ell}<\pd_i\bn, \tau_\ell>,
\quad B^k_i = g^{k\ell}<\pd_j\bn, \tau_\ell>.
\end{align*}
Combining these formulas with \eqref{LB:4} gives
\begin{equation}\label{cg:2.4}
\Delta_{\Gamma_t} f = \Delta_{\Gamma}f 
+ \sum_{i,j=1}^{N-1}( h^{ij}\rho+ O_2)\pd_i\pd_j f 
+ \sum_{k=1}^{N-1}(h^k + O_2)\pd_kf
\end{equation}
with
\begin{align*}
h^k = -\sum_{i,j=1}^{N-1}
(g^{ij}A^k_{ij}\rho + g^{ij}B^k_j\pd_i\rho + 
g^{ij}B^k_i\pd_j\rho + \Lambda^k_{ij}h^{ij}\rho).
\end{align*}
Thus, we have
\begin{align*}
H(\Gamma_t)\bn_t &= 
\Delta_{\Gamma_t}(y + \rho\bn) \\
&= \Delta_\Gamma y + \rho\Delta_\Gamma\bn + \bn\Delta_\Gamma\rho
+ g^{ij}(\pd_i\rho\pd_j\bn + \pd_j\rho\pd_i\bn)\\
&\quad + \rho h^{ij}\pd_i\pd_jy + (\rho h^{ij}\pd_i\pd_j\rho)
\bn + h^k\pd_ky 
+ O_2\pd_i\pd_j\rho + O_2,
\end{align*}
which, combined with \eqref{repr:2.1} and $<\bn, \pd_j\bn>=0$, 
leads to
\begin{align*}
<H(\Gamma_t)\bn_t, \bn_t>
&= <\Delta_\Gamma y, \bn> + \rho<\Delta_\Gamma\bn, \bn>
+ \Delta_\Gamma \rho  +\rho h^{ij}<\pd_i\pd_jy,\bn>\\
&+ (h^{ij}\rho + O_2)\pd_i\pd_j\rho
-g^{ij}<\Delta_\Gamma y, \tau_i>\pd_j\rho  + O_2.
\end{align*} 
Noting that $\Delta_\Gamma y = H(\Gamma)\bn$, we have
\begin{align*}
<H(\Gamma_t)\bn_t, \bn_t>
&= \Delta_\Gamma \rho + H(\Gamma) + \rho<\Delta_\Gamma\bn,\bn>
\\
&+\rho h^{ij}<\pd_i\pd_jy,\bn> 
+ (h^{ij}\rho + O_2)\pd_i\pd_j\rho + O_2.
\end{align*}
Recalling that \eqref{residue:1} and \eqref{residue:2}, we see that 
there exists a function 
 $\bV'_\Gamma(\bar\bk)$  defined on $\overline\Omega \times
\{\bar\bk \mid |\bar\bk| \leq C\}$ satisfying the estimate:
$$\sup_{|\bar\bk| \leq\delta}
\|(\bV_\Gamma'(\cdot, \bar\bk),
\pd_{\bar\bk}\bV_\Gamma'(\cdot, \bar\bk))
\|_{H^1_\infty(\Omega)} \leq C $$
with some constant $C$, 
where $\pd_{\bar\bk}$ denotes the partial derivatives with respect to
variables $\bar\bk$,  
for which 
\begin{equation}\label{bdyfunc.1}
\sum_{\ell=1}^\infty\zeta^1_\ell(h^{ij}\rho + O_2)\pd_i\pd_j\rho + O_2)
= \bV'_\Gamma(\bar\bk)\bar\bk\otimes\bar{\bar\bk},
\end{equation}
where $\bar \bk = (\Psi_\rho, \nabla\Psi_\rho)$ and $\bar{\bar\bk}
= (\Psi_\rho, \nabla \Psi_\rho, \nabla^2\Psi_\rho)$. 
Therefore, we have
\begin{equation}\label{2.5.2}\begin{aligned}
<H(\Gamma_t)\bn_t, \bn_t> 
&= \Delta_\gamma\rho + H(\Gamma)
+ \BB\rho 
 + 
\bV'_\Gamma(\bar\nabla\Psi_\rho)\bar\nabla\Psi_\rho
\otimes\bar\nabla^2\Psi_\rho,
\end{aligned}\end{equation}
where $\bar\nabla^2\Psi_\rho = (\Psi_\rho, \nabla\Psi_\rho, 
\nabla^2\Psi_\rho)$ and we have set
\begin{equation}\label{bdyop:1}
\BB\rho = \{<\Delta_\Gamma\bn, \bn>
+ \sum_{\ell=1}^\infty\zeta^1_\ell(\sum_{i,j=1}^{N-1}
h^{ij}<\pd_i\pd_jy, \bn>)\}\rho.
\end{equation}

To turn to  Eq. \eqref{56}, in view of \eqref{div:1}, \eqref{normal:3.1},
and \eqref{normal:3.3},  we write
\begin{equation}\label{2.5.2*}\begin{split}
<\bn, \mu\bD(\bv)\bn_t> &= <\bn, \mu\bD(\bu)\bn> 
+ <\bn, \mu\bD(\bu)\bV_\Gamma(\bar\bk)\bar\bk> \\
&+ <\bn, \mu(\CD_\bD(\bk)\nabla\bu)(\bn+ \bV_\Gamma(\bar\bk)\bar\bk)>.
\end{split}\end{equation}
where $\bk = \nabla\Psi_\rho$ and $\bar\bk = (\Psi_\rho, 
\nabla\Psi_\rho)$. 
Thus, from  \eqref{2.5.2} and \eqref{2.5.2*}
it follows that the boundary condition \eqref{56}
is transformed to
\begin{equation}\label{2.5.3}
<\bn, \mu\bD(\bu)\bn> - (\fq + \sigma H(\Gamma))
-\sigma(\Delta_\Gamma\rho + \BB\rho) =
h_N(\bu, \Psi_\rho),
\end{equation}
where we have set
\begin{equation}\label{bdy:form}\begin{aligned}
&h_N(\bu, \Psi_\rho) = -<\bn, \mu\bD(\bu)\bV_\Gamma(\bar\bk)\bar\bk> 
 \\
&\quad-<\bn, \mu(\CD_\bD(\bk)\nabla\bu)(\bn+ \bV_\Gamma(\bar\bk)\bar\bk)>
+ \sigma \bV'_\Gamma(\bar\bk)
(\bar\bk, \bar{\bar\bk}),
\end{aligned}\end{equation}
where $\bk = \nabla\Psi_\rho$, $\bar\bk = \bar\nabla\Psi_\rho$,
and $\bar{\bar\bk} = \bar\nabla^2\Psi_\rho$. 
\subsection{Linearization principle}\label{subsec:2.4}

In this subsection, we give a linearization principle
\footnote{The linearization principle means  how to 
divide a nonlinear equation into a linear part and a non-linear
part. }
of Eq. \eqref{navier:1} for the local well-posedness
and the global well-posedness.
Putting \eqref{newchange:1}, \eqref{change:7}, \eqref{kinematic:4}, 
\eqref{58}, and \eqref{2.5.3} together, we see that Eq.\eqref{navier:1} 
is transformed to the equations:
\begin{equation}\label{linear:1}\left\{\begin{aligned}
\pd_t\bu - \DV(\mu\bD(\bu) - \fq\bI)  = \bff(\bu, \Psi_\rho)
&&\enskip&\text{in  $\Omega^T$},  \\
\dv\bu = g(\bu, \Psi_\rho) = \dv\bg(\bu, \Psi_\rho)
&&\enskip&\text{in  $\Omega^T$}, \\
\pd_t\rho +  \xi'(t)\cdot\bn- \bu\cdot\bn 
+ <\bu \mid \nabla'_\Gamma\rho> 
= d(\bu, \rho)&&\enskip&\text{on $\Gamma^T$},  \\
(\mu\bD(\bu)\bn)_\tau = \bh'(\bu, \Psi_\rho)
&&\enskip&\text{on $\Gamma^T$}, \\
<\mu\bD(\bu)\bn, \bn> -(\fq+\sigma H(\Gamma))
-\sigma(\Delta_\Gamma + \BB)\rho
= h_N(\bu, \Psi_\rho) &&\enskip&\text{on $\Gamma^T$}, \\
\bu|_{t=0}=\bu_0 \enskip\text{in $\Omega$},
\quad \rho|_{t=0} = \rho_0 \enskip
\text{on $\Gamma$}&.
\end{aligned}\right.\end{equation}
\vskip0.5pc\noindent
{\bf Local Well-posedness:} \vskip0.3pc
To prove the local well-posedness, the 
unique existence of local in time solutions, we take $\xi(t)=0$ and 
we assume that one of the following conditions holds:
\begin{equation}\label{linear:2}\text{
\thetag{a}} \quad \|H(\Gamma)\|_{W^{1-1/q}_q(\Gamma)} < \infty
\quad \text{or} 
\quad \text{\thetag{b}}\quad \|\nabla H(\Gamma)\|_{L_q(\Omega)} < \infty.
\end{equation} 
\begin{rem} For example, in the case where
$\Omega = B'_R\times\BR$ with $B'_R = \{x' \in 
\BR^{N-1} \mid |x'| < R\}$,  $H(\Gamma) = R^{-(N-2)}\times \chi_\BR(x_N)$,
where $\chi_{\BR}(x_N) = 1$ for any $x_N\in \BR$. And then, 
$\nabla H(\Gamma) = 0$. 
\end{rem}
We introduce the new pressure $\fp$ defined by letting
$\fp = \fq + \sigma H(\Gamma)$ in the \thetag{a} case. 
For the reference body
$\Omega$, we choose the boundary $\Gamma_0$ of 
the initial domain in such a way that
$$\Gamma_0 = \{x = y + h_0(y)\bn \mid y \in \Gamma\}$$
and $\|h_0\|_{B^{3-1/p-1/q}_{q,p}(\Gamma)}$ is small enough. 
But, we do not want to have any  restriction on the size of 
the initial velocity, $\bu_0$.  Thus, we approximate 
$\bu_0$ by $\bu_\kappa$ defined by letting
$$\bu_\kappa = \frac{1}{\kappa}\int^\kappa_0 T(s)\tilde\bu_0\,ds,$$
where $\tilde \bu_0$ is a suitable extension of $\bu_0$ to 
$\BR^N$ satisfying the condition:
\begin{equation}\label{linear:3}
\|\tilde\bu_0\|_{B^{2(1-1/p)}_{q,p}(\BR^N)} 
\leq C\|\bu_0\|_{B^{2(1-1/p)}_{q,p}(\Omega)},
\end{equation}
and $T(s)$ is some analytic semigroup defined on $\BR^N$ satisfying
the estimate:
\begin{gather}
\|T(\cdot)\tilde\bu_0\|_{L_\infty((0, \infty), 
B^{2(1-1/p)}_{q,p}(\BR^N))} \leq C\|\bu_0\|_{B^{2(1-1/p)}_{q,p}(\Omega)},
\nonumber \\
\|T(\cdot)\tilde\bu_0\|_{L_p((0, \infty), H^2_q(\BR^N))}
+ \|T(\cdot)\tilde\bu_0\|_{H^1_p((0, \infty), L_q(\BR^N))}
\leq C\|\bu_0\|_{B^{2(1-1/p)}_{q,p}(\Omega)}.
\label{linear:4}
\end{gather}
By \eqref{linear:4}
\begin{equation}\label{linear:5}\begin{split}
\|\bu_\kappa\|_{B^{2(1-1/p)}_{q,p}(\Omega)} 
&\leq C\|\bu_0\|_{B^{2(1-1/p)}_{q,p}(\Omega)}, \\
\|\bu_\kappa\|_{H^2_q(\Omega)} & \leq C\|\bu_0\|_{B^{2(1-1/p)}_{q,p}(\Omega)}
\kappa^{-1/p}.
\end{split}\end{equation}
From Eq. \eqref{linear:1},
the linearlization principle of Eq. \eqref{navier:1} is
the following equations: 
\begin{equation}\label{linear:6}\left\{\begin{aligned}
\pd_t\bu - \DV(\mu\bD(\bu) - \fq\bI)  = \bff(\bu, \Psi_\rho)
+ \kappa_b\nabla H(\Gamma)&&\enskip&\text{in $\Omega^T$}\, \\
\dv\bu = g(\bu, \Psi_\rho) = \dv\bg(\bu, \Psi_\rho)&
&\enskip&\text{in $\Omega^T$}\, \\
\pd_t\rho + <\bu_\kappa \mid \nabla'_\Gamma\rho> - \bu\cdot\bn 
= d(\bu, \Psi_\rho) + <\bu_\sigma - \bu \mid \nabla'_\Gamma\rho>&
&\enskip&\text{in $\Gamma^T$}\, \\
(\mu\bD(\bu)\bn)_\tau = \bh'(\bu, \Psi_\rho)&
&\enskip&\text{in $\Gamma^T$}\, \\
<\mu\bD(\bu)\bn, \bn> -\fq
-\sigma(\Delta_\Gamma + \BB)\rho
= h_N(\bu, \Psi_\rho)+ \kappa_a H(\Gamma)&
&\enskip&\text{in $\Gamma^T$}\, \\
\bu|_{t=0} = \bu_0 \enskip\text{in $\Omega$},
\quad  \rho|_{t=0} = \rho_0
\enskip\text{on $\Gamma$}&.
\end{aligned}\right.\end{equation}
where we have set 
\begin{align*}
&\kappa_a = \begin{cases} \sigma&\quad\text{in  the \thetag{a} case}, \\
0 &\quad\text{in  the \thetag{b} case},
\end{cases} \quad
\kappa_b = \begin{cases} 0&\quad\text{in  the \thetag{b} case}, \\
\sigma &\quad\text{in  the \thetag{a} case}.
\end{cases}
\end{align*}
{\bf Global well-posedness} \vskip0.3pc
In this lecture note, we only treat the bounded domain case.  The 
unbounded domain case was treated, for example,  by Saito and Shibata
\cite{Saito1, Saito2} in the $L_p$-$L_q$ framework.
 Assuming that the reference body
$\Omega_0$ is very closed to a ball and initial velocities
 are sufficiently small,
we shall prove the global wellposedness, the unique existence of 
global in time solutions. Here, the problem is formulated and 
the global well-posedness will be proved in Sect. \ref{sec:6}.

Let $B_R$ be the ball of radius $R$ centered at the origin and let
$S_R$ be the sphere of radius $R$ centered at the origin, 
that is $S_R$ is the boundary of $B_R$.  We assume
that
\begin{itemize}
\item[\thetag{A.1}]\qquad $|\Omega| = |B_R| = \dfrac{R^N\omega_N}{N}$, 
where $|D|$ denotes the Lebesgue measure of a Lebesgue measurable set 
$D$ in $\BR^N$ and $\omega_N$ is the area of $S_1$.
\item[\thetag{A.2}]\qquad$\displaystyle{\int_\Omega x\,dx = 0}$. 
\item[\thetag{A.3}]\qquad$\Gamma$ is a normal perturbation of $S_R$ given
by $\Gamma = \{x = y+\rho_0(y)\bn(y) \mid y \in S_R\}$ 
with a given small function $\rho_0(y)$ defined on $S_R$.
\end{itemize}
Here, $\bn = y/|y|$ is the unit outer normal to $S_R$. $\bn$ is extended to
$\BR^N$ by $\bn = R^{-1}y$. 
Let $\Gamma_t$ be given by 
\begin{align*}
\Gamma_t &= \{x = y + \rho(y, t)\bn + \xi(t) \mid y \in S_R\} \\
&= \{x = y + R^{-1}\rho(y, t)y + \xi(t) \mid y \in S_R\}
\end{align*}
where $\rho(y, t)$ is an unknown function with $\rho(y,0) = \rho_0(y)$
for $y \in S_R$. Let $\xi(t)$ be the barycenter point of the 
domain $\Omega_t$ defined by 
$$\xi(t) = \frac{1}{|\Omega|}\int_{\Omega_t} x\,dx.
$$
Here, $|\Omega_t| = |\Omega|$ follows eventually from 
the divergence free condition in 
Eq. \eqref{navier:1}. Notice that $\xi(t)$ is also unknown.  It follows from
the assumption \thetag{A.2} that $\xi(0) = 0$.  
Given a function $\rho$ defined on $S_R$, let $H_\rho$ be a suitable extension 
of $\rho$ such that $H_\rho = \rho$ on $S_R$ and the following
estimates hold:
\begin{equation}\label{sob:5.4.b}\begin{split}
C_1\|H_\rho(\cdot, t)\|_{H^k_q(B_R)} 
&\leq \|\rho(\cdot, t)\|_{W^{k-1/q}_q(S_R)}\leq 
C_2\|H_\rho(\cdot, t)\|_{H^k_q(B_R)},  \\
C_1\|\pd_tH_\rho(\cdot, t)\|_{H^\ell_q(B_R)} 
&\leq \|\pd_t\rho(\cdot, t)\|_{W^{\ell-1/q}_q(S_R)}\leq 
C_2\|\pd_tH_\rho(\cdot, t)\|_{H^\ell_q(B_R)} 
\end{split}\end{equation}
with some constants $C_1$ and $C_2$ for $k=1,2,3$ and 
$\ell=1,2$. 
We define the Hanzawa transform centered at $\xi(t)$ by
\begin{equation}\label{g.1.5}
x = \bh_z(y, t) = y + \Psi_\rho(y, t) + \xi(t)
\quad\text{for $y \in B_R$},
\end{equation}
where $\Psi_\rho(y, t) = R^{-1}H_\rho(y, t)y$, because  
$\bn = y/|y|$ for $y \in S_R$. 
And then, the linearization principle is  the following:
\begin{equation}\label{linear:8}\left\{\begin{aligned}
\pd_t\bu - \DV(\mu\bD(\bu) - \fq\bI) & = \bff(\bu, \Psi_\rho)
&\quad&\text{in $B_R^T$}, \\
\dv\bu = g(\bu, \Psi_\rho) &= \dv\bg(\bu, \Psi_\rho)
&\quad&\text{in $B_R^T$}, \\
\pd_t\rho - \bn\cdot P\bu & = d(\bu, \Psi_\rho)
&\quad&\text{on $S_R^T$}, \\
(\mu\bD(\bu)\bn)_\tau & = \bh'(\bu, \Psi_\rho)
&\quad&\text{on $S_R^T$}, \\
<\mu\bD(\bu)\bn, \bn> -\fq
-\sigma\CB\rho & = h_N(\bu, \Psi_\rho)
&\quad&\text{on $S_R^T$}, \\
(\bu, \rho)|_{t=0} & = (\bu_0, \rho_0)
&\quad&\text{in $B_R\times S_R$},
\end{aligned}\right.\end{equation}
where,  $B_R^T = B_R\times(0, T)$, 
$S_R^T = S_R\times(0, T)$, $\displaystyle{
P\bu = \bu - \frac{1}{|B_R|}\int_{B_R}\bu\,dy}$, $\CB = \Delta_{S_R}
+ \dfrac{N-1}{R^2}$, $\Delta_{S_R}$ is the Laplace-Beltrami
operator on $S_R$.  Moreover, 
$\bff(\bu, \Psi_\rho)$, $g(\bu, \Psi_\rho)$,
$\bg(\bu, \Psi_\rho)$, $\bh'(\bu, \Psi_h)$ are the same 
as in Eq. \eqref{linear:6}, 
but $h_N(\bv, \Psi_\rho)$ and 
$d(\bv, \Psi_\rho)$ are slightly different functions from these in
Eq. \eqref{linear:6} given latter. And, the last boundary condition in
\eqref{linear:8} will be given in Subsec. \ref{subsec.7.1}. 
When we prove the global wellosedness,
we assume that initial velocity, $\bu_0$,
is small enough, and so $\bu$ is small.  Thus, term: 
$\sum_{i,j=1}^{N-1}g^{ij}
<\tau_i, \bu>\frac{\pd \rho}{\pd p_j}$
in the kinematic equations can be put in  the right side. 

\section{Maximal $L_p$-$L_q$ regularity} \label{sec:3}
To prove the local and global well-posedness for Eq.\eqref{navier:1},
in this lecture note the maximal regularity theorem for the linear part 
(the left hand side of Eq. \eqref{linear:6} and \eqref{linear:8}) plays an 
essential role.  In this section, we study the maximal $L_p$-$L_q$
regularity theorem and the generation of analytic semigroup
associated with the Stokes equations with free boundary conditions.
\subsection{Statement of maximal regularity theorems}\label{subsec:3.1}
In this subsection, we shall state the maximal $L_p$-$L_q$ regularity
for the following three problems:
\begin{equation}\label{stokes:1}\left\{\begin{aligned}
\pd_t\bu - \DV(\mu\bD(\bu) - \fp\bI) = \bff,
\quad \dv\bu &= g = \dv \bg
&\quad&\text{in $\Omega^T$}, \\
(\mu\bD(\bu) - \fp\bI)\bn 
& = \bh &\quad&\text{on $\Gamma^T$}, \\
\bu& = \bu_0
&\quad&\text{in $\Omega$};
\end{aligned}\right.\end{equation}
\begin{equation}\label{stokes:2}\left\{\begin{aligned}
\pd_t\bv - \DV(\mu\bD(\bv) - \fq\bI) = 0,
\quad \dv\bv &= 0
&\quad&\text{in $\Omega^T$}, \\
\pd_t h - \bn \cdot\bv +\CF_1\bv & = d
&\quad&\text{on $\Gamma^T$}, \\
(\mu\bD(\bv) - \fq\bI)\bn - (\CF_2 h + \sigma\Delta_\Gamma h)\bn
& = 0 &\quad&\text{on $\Gamma^T$}, \\
(\bv, h)|_{t=0} & = (0, h_0)
&\quad&\text{in $\Omega\times\Gamma$};
\end{aligned}\right.\end{equation}
\begin{equation}\label{stokes:3}\left\{\begin{aligned}
\pd_t\bw - \DV(\mu\bD(\bw) - \fr\bI) = 0,
\quad \dv\bw &= 0
&\quad&\text{in $\Omega^T$}, \\
\pd_t \rho + A_\kappa\cdot\nabla'_\Gamma \rho 
-\bw\cdot\bn + \CF_1\bw& = d 
&\quad&\text{on $\Gamma^T$}, \\
(\mu\bD(\bw) - \fr\bI)\bn - (\CF_2\rho + \sigma\Delta_\Gamma\rho)\bn
& = 0 &\quad&\text{on $\Gamma^T$}, \\
(\bw, \rho)|_{t=0} & = (0, 0)
&\quad&\text{in $\Omega\times\Gamma$}.
\end{aligned}\right.\end{equation}
Here, $\CF_1$ and $\CF_2$ are linear operators such that 
\begin{equation}\label{cond-low}
\|\CF_1\bv\|_{W^{2-1/q}_q(\Gamma)} \leq M_0\|\bv\|_{H^1_q(\Omega)}, \quad
\|\CF_2 h\|_{H^1_q(\Omega)} \leq C_0\|h\|_{H^2_q(\Omega)}
\end{equation}
 with some constant $M_0$.   
If we consider the total problem:
\begin{equation}\label{stokes:4}\left\{\begin{aligned}
\pd_t U - \DV(\mu\bD(U) - P \bI) = \bff,
\quad \dv U &= g = \dv\bg 
&\quad&\text{in $\Omega^T$}, \\
\pd_t H + A_\kappa\cdot\nabla'_\Gamma H 
-U\cdot\bn + \CF_1U& = d 
&\quad&\text{on $\Gamma^T$}, \\
(\mu\bD(U) - P\bI)\bn - (\CF_2H + \sigma\Delta_\Gamma H)\bn
& = \bh &\quad&\text{on $\Gamma^T$}, \\
(U, H)|_{t=0} & = (\bu_0, h_0)
&\quad&\text{in $\Omega\times\Gamma$}.
\end{aligned}\right.\end{equation}
then, $U$, $P$ and $H$ are given by $U = \bu + \bv + \bw$
and $P = \fp + \fq + \fr$ and $H = h + \rho$, where $(\bu, \fp)$
is a solution of Eq.\eqref{stokes:1}, $(\bv, \fq, h)$
 a solution of Eq.\eqref{stokes:2} with $d=0$, 
and $(\bw, \fr, \rho)$ a solution of Eq.\eqref{stokes:3}
by replacing $d$ by $d+\bn\cdot\bu - \CF_1\bu - A_\kappa\cdot
\nabla'_\Gamma h$.


We now introduce a solenoidal space $J_q(\Omega)$ defined by 
\begin{equation}\label{sol:1}
J_q(\Omega) = \{\bff \in L_q(\Omega)^N \mid 
(\bff, \nabla\varphi)_\Omega = 0 \quad\text{for any 
$\varphi \in \hat H^1_{q',0}(\Omega)$}\}.
\end{equation}
%

Before stating the maximal $L_p$-$L_q$ regularity theorems for the
three equations given above, we introduce the assumptions on
$\mu$ and $A_\kappa$. \vskip0.5pc\noindent
{\bf Assumption on $\mu$, $\sigma$, and $A_\kappa$} \vskip0.3pc
There exist positive constants $m_0$, $m_1$, 
$m_2$, $m_3$, $a$ and $b$ for which, 
\begin{gather}
m_0 \leq \mu(x), \sigma(x) \leq m_1, \quad
|\nabla(\mu(x), \sigma(x))| \leq m_1
\quad\text{for any $x \in \overline{\Omega}$},\nonumber \\ 
|A_\kappa(x)| \leq m_2, 
\quad |A_\kappa(x)-A_\kappa(y)| \leq m_2|x-y|^a
\quad\text{for any $x$, $y \in \Gamma$},\nonumber\\ 
\|A_\kappa\|_{W^{2-1/r}_r(\Gamma)} \leq m_3\kappa^{-b}
\enskip(\kappa \in (0, 1)). \label{assump:3}
\end{gather}
Where, $r$ is an exponent in $(N, \infty)$. 

Moreover, in view of \eqref{assump:3}, 
for $=B^i_j=B_{r_0}(x^i_j)$ given in Proposition \ref{prop:lap}
 we assume that
\begin{align}
& |\mu(x) - \mu(x^0_j)| \leq M_1,
 \quad\text{for $x^0_j \in 
B^i_j$}, \label{assump:4.1} \\
& |\mu(x) - \mu(x^1_j)| \leq M_1,\enskip
 |\sigma(x) - \sigma(x^1_j)|
\leq M_1, \enskip |A_\kappa(x) - A_\kappa(x^1_j)| 
\leq M_1
\label{assump:4.2}
\end{align}
for any $x \in B^1_j$. 
Since $H^1_q(\Omega)$ is usually not dense in $\hat H^1_q(\Omega)$, 
it does not hold that $\dv \bu = \dv \bg$ implies $(\bu, \nabla\varphi)_\Omega
= (\bg, \nabla\varphi)_\Omega$ for all $\varphi \in \hat H^1_q(\Omega)$.
Of course, the opposite direction holds.  Thus, finally we introduce the 
following definition.
\begin{dfn}\label{dfn:1.4}  For $\bu$, $\bg\in L_q(\Omega)^N$, we say that 
$\dv \bu = \dv\bg$ in $\Omega$ if  $\bu - \bg 
\in J_q(\Omega)$.
\end{dfn}
To solve the divergence equation $\dv \bu = g$ in $\Omega$, it is necessary
to assume that $g$ is given by $g = \dv\bg$ for some $\bg$, and so we define 
the space $DI_q(G)$ by
$$DI_q(G) = \{(g, \bg) \mid g \in H^1_q(G), \enskip 
\bg \in L_q(G), \enskip g = \dv\bg \enskip\text{in $G$}\},
$$
where $G$ is any domain in $\BR^N$. 

Before stating the main results in this section, we give a definition
of the uniqueness. For Eq. \eqref{stokes:1}, the uniqueness is
defined as follows:
\begin{itemize}
\setlength{\itemsep}{0pt}
\item~ If $\bv$ and $\fp$ with 
\begin{align*}
\bv &\in L_p((0, T), H^2_q(\Omega)^N) \cap H^1_p((0, T), L_q(\Omega)^N), 
\\ 
\fp & \in L_p((0, T), H^1_q(\Omega) + \hat H^1_{q,0}(\Omega)),
\end{align*}
satisfy the homogeneous equations:
\begin{equation}\label{homoeq:1.1}
\left\{\begin{aligned}
\pd_t\bv - \DV(\mu\bD(\bv) - \fp\bI) = 0,
\quad \dv\bv &= 0
&\quad &\text{in $\Omega^T$}, \\
(\mu\bD(\bv) - \fp\bI)\bn& = 0 &\quad &\text{on $\Gamma^T$}, \\
\bv|_{t=0} & = 0 
&\quad &\text{on $\Omega$},
\end{aligned}\right.
\end{equation}
then, $\bv=0$ and $\fp=0$.
\end{itemize}
For Eq. \eqref{stokes:2} and \eqref{stokes:3},
the uniqueness is defined as follows: 
\begin{itemize}
\setlength{\itemsep}{0pt}
\item~ If $\bv$, $\fp$ and $\rho$ with 
\begin{align*}
\bv &\in L_p((0, T), H^2_q(\Omega)^N) \cap H^1_p((0, T), L_q(\Omega)^N), 
\\
\fp & \in L_p((0, T), H^1_q(\Omega) + \hat H^1_{q,0}(\Omega)), \\
\rho & \in L_p((0, T), W^{3-1/q}_q(\Gamma)) \cap H^1_p((0, T), W^{2-1/q}_q
(\Gamma))
\end{align*}
satisfy the homogeneous equations:
\begin{equation}\label{homoeq:1.1*}
\left\{\begin{aligned}
\pd_t\bv - \DV(\mu\bD(\bv) - \fp\bI) = 0,
\quad \dv\bv &= 0
&\quad &\text{in $\Omega^T$}, \\
\pd_t\rho + A_\kappa\cdot\nabla'_\Gamma\rho-\bv\cdot\bn
+ \CF_1\bv
& = 0 &\quad &\text{on $\Gamma^T$}, \\
(\mu\bD(\bv) - \fp\bI-((\CF_2 + \sigma\Delta_\Gamma)\rho)\bI)
\bn & = 0 &\quad &\text{on $\Gamma^T$}, \\
(\bv, \rho)|_{t=0} & = (0, 0) 
&\quad &\text{on $\Omega\times \Gamma$},
\end{aligned}\right.
\end{equation}
then, $\bv=0$, $\fp=0$, and $\rho=0$.
\end{itemize}
Notice that when $A_\kappa=0$, the uniquness is stated in the same
manner as in \eqref{homoeq:1.1*}.

We now state the maximal $L_p$-$L_q$ regularity theorems. 
\begin{thm}\label{thm:max.1}
Let $1 < p, q < \infty$ with $2/p + 1/q \not=1$.
 Assume that $\Omega$ is a uniform $C^2$
domain and that the weak Dirichlet problem is uniquely solvable
for index $q$. 
Then, there exists a $\gamma_0 > 0$ such that
the following assertion holds:  Let $\bu_0 \in B^{2(1-1/p)}_{q,p}(\Omega)^N$
be initial data for
Eq. \eqref{stokes:1} and let $\bff$, $g$, $\bg$, and $\bh$ be 
functions in the right side of Eq. \eqref{stokes:1} such that
$\bff \in L_p((0, T), L_q(\Omega)^N)$, and 
\begin{gather}
e^{-\gamma t}\bg
\in H^1_p(\BR, L_q(\Omega)^N), 
\quad e^{-\gamma t}g \in L_p(\BR, H^1_q(\Omega))
\cap H^{1/2}_p(\BR, L_q(\Omega)), \nonumber \\
e^{-\gamma t}\bh \in L_p(\BR, H^1_q(\Omega)^N)
\cap H^{1/2}_p(\BR, L_q(\Omega)^N)
 \label{increase:1}
\end{gather}
for any $\gamma \geq \gamma_0$. Assume that the compatibility
condition:
\begin{equation}
\bu_0 - \bg|_{t=0} \in J_q(\Omega)\quad
\text{and}  \quad \dv\bu_0 = g|_{t=0} \quad
 \text{in $\Omega$}, \label{compati:1}
\end{equation}
holds.  In addition, we assume that  the compatibility condition:
\begin{equation}
(\mu\bD(\bu_0)\bn)_\tau = (\bh|_{t=0})_\tau
\quad\text{on $\Gamma$}\label{compati:2}
\end{equation}
holds for $2/p + 1/q < 1$.  Then, problem \eqref{stokes:1} admits
 solutions $\bu$ and  $\fp$ with
\begin{align*}
\bu & \in L_p((0, T), H^2_q(\Omega)^N) \cap H^1_q((0, T), 
L_q(\Omega)^N),
\\
\fp &\in L_p((0, T), H^1_q(\Omega) + \hat H^1_{q,0}(\Omega))
\end{align*}
possessing the estimate:
\begin{align*}
&\|\bu\|_{L_p((0, T), H^2_q(\Omega))}
+ \|\pd_t\bu\|_{L_p((0, T), L_q(\Omega))} \\
&\quad
\leq Ce^{\gamma T}(\|\bu_0\|_{B^{2(1-1/p)}_{q,p}(\Omega)} 
 + \|\bff\|_{L_p((0, T), L_q(\Omega))} 
+\|e^{-\gamma t}\bg\|_{H^1_p(\BR, L_q(\Omega))}\\
&\quad+ \|e^{-\gamma  t}(g, \bh)\|_{L_p(\BR, H^1_q(\Omega))}
+ \|e^{-\gamma t}(g, \bh)\|_{H^{1/2}_p(\BR, L_q(\Omega))})
\end{align*}
for some constant $C > 0$.

Moreover,  if we assume that the weak Dirichlet problem is uniquely
solvable for $q' = q/(q-1)$ in addiion, then the uniqueness 
for Eq. \eqref{stokes:1} 
holds. 
\end{thm}
\begin{thm}\label{thm:max.2}
Let $1 < p, q < \infty$.
 Assume that $\Omega$ is a uniform $C^3$
domain and that the weak Dirichlet problem is uniquely solvable
for $q$.  
 Let 
$h_0 \in B^{3-1/p-1/q}_{q,p}(\Gamma)$ be initial data for
Eq. \eqref{stokes:2} and let $d$ be a function 
 in the right side of Eq. \eqref{stokes:2} such that
\begin{equation}
d \in L_p((0, T), W^{2-1/q}_q(\Gamma)). \label{increase:2}
\end{equation}
Then, problem \eqref{stokes:2} admits
 solutions $\bv$, $\fq$ and $h$ with
\begin{align*}
\bv & \in L_p((0, T), H^2_q(\Omega)^N) \cap H^1_p((0, T), 
L_q(\Omega)^N),
\\
\fq & \in L_p((0, T), H^1_q(\Omega) + \hat H^1_{q,0}(\Omega)), \\
h & \in L_p((0, T), W^{3-1/q}_q(\Gamma)) 
\cap H^1_p((0, T), W^{2-1/q}_q(\Gamma))
\end{align*}
possessing the estimate:
\begin{align*}
&\|\bv\|_{L_p((0, T), H^2_q(\Omega))}
+ \|\pd_t\bv\|_{L_p((0, T), L_q(\Omega))}\\ 
&\quad+ \|h\|_{L_p((0, T), W^{3-1/q}_q(\Gamma))} + 
\|\pd_th\|_{L_p((0, T), W^{2-1/q}_q(\Gamma))} \\
&\qquad\leq Ce^{\gamma T}(\|h_0\|_{B^{3-1/p-1/q}_{q,p}(\Gamma)} 
+ \|d\|_{L_p((0, T), W^{2-1/q}_q(\Gamma))})
\end{align*}
for some constants $C > 0$ and $\gamma$. 

Moreover, if we assme that the weak Dirichlet problem is uniquely
solvable for $q' = q/(q-1)$ in addition, then the uniqueness
for Eq. \eqref{stokes:2} 
holds. 
\end{thm}
\begin{thm}\label{thm:max.3}
Let $1 < p, q < \infty$. 
 Assume that $\Omega$ is a uniform $C^3$
domain and that the weak Dirichlet problem is uniquely solvable
for $q$.   Let $d$ be a function 
 in the right side of Eq. \eqref{stokes:3} such that
\begin{equation}
d \in L_p((0, T), W^{2-1/q}_q(\Gamma)). \label{increase:3}
\end{equation}
Then, problem \eqref{stokes:3} admits
unique solutions $\bw$, $\fr$ and $\rho$ with
\begin{align*}
\bw & \in L_p((0, T), H^2_q(\Omega)^N) \cap H^1_q((0, T), 
L_q(\Omega)^N), \\
\fr & \in L_p((0, T), H^1_q(\Omega) + \hat H^1_{q,0}(\Omega)), \\
\rho & \in L_p((0, T), W^{3-1/q}_q(\Gamma)) 
\cap H^1_p((0, T), W^{2-1/q}_q(\Gamma))
\end{align*}
possessing the estimate:
\begin{align*}
&\|\bw\|_{L_p((0, T), H^2_q(\Omega))}
+ \|\pd_t\bw\|_{L_p((0, T), L_q(\Omega))} \\
&\quad + \|\rho\|_{L_p((0, T), W^{3-1/q}_q(\Gamma))}
+ \|\pd_t\rho\|_{L_p((0, T), W^{2-1/q}_q(\Gamma))} \\
&\qquad
\leq Ce^{\gamma \kappa^{-b}T}
\|d\|_{L_p((0, T), W^{2-1/q}_q(\Gamma))}
\end{align*}
for some constants $C > 0$ and $\gamma$, where $\kappa \in (0, 1)$ and
$b$ is the constant appearing
in \eqref{assump:3}.

Moreover, if we assume that $\Omega$ is a uniform $C^3$
domain whose inside has a finite covering 
and that the weak Dirichlet problem is uniquely solvable
for $q'=q/(q-1)$ in addition, then    
the uniqueness for Eq. \eqref{stokes:3} holds. 
\end{thm}
\begin{rem} The uniqueness follows from the existence of solutions
to the dual problem in the case of Eq.\eqref{stokes:1} and Eq.\eqref{stokes:2}.
But, the uniqueness for the Eq.\eqref{stokes:3} follows from the 
{\it a priori} estimates, because we can not find a suitable
dual problem for Eq.\eqref{stokes:3}. Thus, in addition,
we need the assumption that the inside of $\Omega$ has a
finite covering. 
\end{rem}
Applying Theorems \ref{thm:max.1}, \ref{thm:max.2}, and 
\ref{thm:max.3}, we have the following
corollary.
\begin{cor}\label{cor:max.1}
Let $1<p, q < \infty$ with $2/p + 1/q \not=1$. 
Assume that $\Omega$ is a uniform 
$C^3$ domain and that the 
weak Dirichlet problem is uniquely solvable for $q$.
Let $\bu_0 \in B^{2(1-1/p)}_{q,p}(\Omega)^N$ and $h_0 
\in B^{3-1/p-1/q}_{q,p}(\Gamma)$ be initial data for Eq.\eqref{stokes:4} 
and let $\bff$, $g$, $\bg$, $d$, and $\bh$ be functions appearing in the
right side of Eq.\eqref{stokes:4} and satisfying the conditions:
\begin{align*}
\bff &\in L_p((0, T), L_q(\Omega)^N), \quad 
d \in L_p((0, T), W^{2-1/q}_q(\Gamma)), \\
e^{-\gamma t}g &\in L_p(\BR, H^1_q(\Omega)) \cap H^{1/2}_p(\BR, L_q(\Omega)), 
\quad
e^{-\gamma t}\bg \in H^1_p(\BR, L_q(\Omega)^N), 
\\
e^{-\gamma t}\bh
&\in L_p(\BR, H^1_q(\Omega)^N) \cap H^{1/2}_p(\BR, L_q(\Omega)^N).
\end{align*}
for any $\gamma \geq \gamma_0$ with some $\gamma_0$. 
Assume that the compatibility conditions \eqref{compati:1}
and \eqref{compati:2} are satisfied.  Then, problem \eqref{stokes:4}
admits  solutions $U$ $P$, and $H$ with 
\begin{align*}
U &\in L_p((0, T), H^2_q(\Omega)^N) \cap H^1_p((0, T), L_q(\Omega)^N), \\
P &\in L_p((0, T), H^1_q(\Omega) + H^1_{q,0}(\Omega)), \\
H &\in L_p((0, T), W^{3-1/q}_q(\Gamma)) \cap H^1_p((0, T), W^{2-1/q}_q(\Gamma))
\end{align*}
possessing the estimate:
\begin{align*}
&\|U\|_{L_p((0, T), H^2_q(\Omega))} + 
\|\pd_tU\|_{L_p((0, T), L_q(\Omega))}
+\|H\|_{L_p((0, T), W^{3-1/q}_q(\Gamma)} \\
&\quad + 
\|\pd_tH\|_{L_p((0, T), W^{2-1/q}_q(\Gamma))} \\
&\leq Ce^{2\gamma\kappa^{-b}T}\{\|\bu_0\|_{B^{2(1-1/p)}_{q,p}(\Omega)}
+ \kappa^{-b}\|h_0\|_{B^{3-1/p-1/q}_{q,p}(\Gamma)} \\
&\quad+ \|\bff\|_{L_p((0, T), L_q(\Omega))}
+ \|e^{-\gamma t}\pd_t\bg\|_{L_p(\BR, L_q(\Omega))}
 +\|e^{-\gamma t}(g, \bh)\|_{L_p(\BR, H^1_q(\Omega))} \\
&\quad + \|e^{-\gamma t}(g, \bh)\|_{H^{1/2}_p(\BR, L_q(\Omega))}
+ \|d\|_{L_p((0, T), W^{2-1/q}_q(\Gamma))}\}
\end{align*}
for any $\gamma \geq \gamma_0$ with some constant $C$ independent
of $\gamma$ and $\kappa$, where $\kappa \in (0, 1)$, and $b$ is the 
constant appearing in \eqref{assump:3}. 
\end{cor}
\pf As was mentioned after Eq.\eqref{stokes:4}, $U$, $P$ and $H$ are 
given by $U = \bu+\bv+\bw$, $P = \fp+\fq+\fr$ and 
$H = h + \rho$ Since $\bu$ and $\fp$ are solutions of Eq. \eqref{stokes:1},
by Theorem  \ref{thm:max.1} we have 
\begin{equation}\label{ccor:max.1}\begin{split}
&\|\bu\|_{L_p((0, T), H^2_q(\Omega))} + 
\|\pd_t\bu\|_{L_p((0, T), L_q(\Omega))} \\
& \leq Ce^{\gamma T}\{\|\bu_0\|_{B^{2(1-1/p)}_{q,p}(\Omega)}
+ \|\bff\|_{L_p((0, T), L_q(\Omega))}
+ \|e^{-\gamma t}\pd_t\bg\|_{L_p(\BR, L_q(\Omega))} \\
&\quad+\|e^{-\gamma t}(g, \bh)\|_{L_p(\BR, H^1_q(\Omega))}
+ \|e^{-\gamma t}(g, \bh)\|_{H^{1/2}_p(\BR, L_q(\Omega))}\}.
\end{split}\end{equation}
Since $\bv$, $\fq$ and $h$ are solutions of Eq. \eqref{stokes:2} with
$d=0$, appying Theorem \ref{thm:max.2} with $d=0$ 
yields that  
\begin{equation}\label{cor:max.2}\begin{split}
&\|\bv\|_{L_p((0, T), H^2_q(\Omega))} + 
\|\pd_t\bv\|_{L_p((0, T), L_q(\Omega))} 
+ \|h\|_{L_p((0, T), W^{3-1/q}_q(\Gamma))}\\
&\quad + \|\pd_th\|_{L_p((0, T), W^{2-1/q}_q(\Gamma))} 
\leq Ce^{\gamma T}\|h_0\|_{B^{3-1/p-1/q}_{q,p}(\Gamma)}.
\end{split}\end{equation}
Finally, recalling that $\bw$, $\fr$ and $\rho$ are solutions of 
Eq. \eqref{stokes:3} by replacing $d$ by $d + \bn\cdot\bu-\CF_1\bu-A_\kappa\cdot\nabla'_\Gamma h$, 
applying Theorem \ref{thm:max.3} and using the estimate
$$\|A_\kappa\cdot\nabla'_\Gamma h\|_{W^{2-1/q}_q(\Gamma)}
\leq C\|A_\kappa\|_{H^2_r(\Omega)}\|h\|_{W^{3-1/q}_q(\Gamma)},
$$
which follows from the Sobolev imbedding theorem 
and the assumption: $N < r < \infty$, we have
\begin{align*}
&\|\bw\|_{L_p((0, T), H^2_q(\Omega))} + 
\|\pd_t\bw\|_{L_p((0, T), L_q(\Omega))}\\
&\quad + \|\rho\|_{L_p((0, T), W^{3-1/q}_q(\Gamma))} 
+ \|\pd_t\rho\|_{L_p((0, T), W^{2-1/q}_q(\Gamma))} \\
&\leq Ce^{\gamma \kappa^{-b}T}
(\|h_0\|_{B^{3-1/p-1/q}_{q,p}(\Gamma)} \\
&\qquad + \|d+\bn\cdot\bu -\CF_1\bu - A_\kappa\cdot
\nabla'_\Gamma h\|_{L_p((0, T), W^{2-1/q}_q(\Gamma))}\}\\
&\leq Ce^{\gamma \kappa^{-b}T}\{
\kappa^{-b}\|h\|_{L_p((0, T), W^{3-1/q}_q(\Gamma)} + 
\|d\|_{L_p((0, T), W^{2-1/q}_q(\Gamma))}\\
&\qquad+ \|\bu\|_{L_p((0, T), H^2_q(\Omega))}\},
\end{align*}
which, combined with \eqref{ccor:max.1} and \eqref{cor:max.2},
leads to the required estimate.  This completes the proof of 
Corollary \ref{cor:max.1}. \qed

\subsection{$\CR$ bounded solution operators}\label{subsec:3.2*}
To prove Theorems \ref{thm:max.1}, \ref{thm:max.2}, 
and \ref{thm:max.3}, we use  $\CR$ bounded solution
operators associated with the following generalized resolvent 
problems:
\begin{equation}\label{res:1.1}\left\{\begin{aligned}
\lambda\bu - \DV(\mu\bD(\bu) - \fp\bI) = \bff, \quad
\dv \bu  = g =\dv\bg& &\quad&\text{in $\Omega$}, \\
(\mu\bD(\bu) - \fp\bI)\bn = \bh&
&\quad&\text{on $\Gamma$};
\end{aligned}\right.\end{equation}
\begin{equation}\label{res:1.2}\left\{\begin{aligned}
\lambda\bv - \DV(\mu\bD(\bv) - \fq\bI) = \bff, \quad
\dv \bv  = g=\dv\bg& &\quad&\text{in $\Omega$}, \\
\lambda \rho - \bv\cdot\bn + \CF_1\bv =d& &\quad&\text{on $\Gamma$}, \\
(\mu\bD(\bv) - \fq\bI)\bn - (\CF_2 \rho+\sigma\Delta_\Gamma \rho)\bn
 = \bh& &\quad&\text{on $\Gamma$}; 
\end{aligned}\right.\end{equation}
\begin{equation}\label{res:1.3}\left\{\begin{aligned}
\lambda\bv - \DV(\mu\bD(\bv) - \fq\bI) = \bff, \quad
\dv \bv =g=\dv\bg& &\quad&\text{in $\Omega$}, \\
\lambda \rho
+ A_\kappa\cdot\nabla'_\Gamma \rho - \bv\cdot\bn + \CF_1\bv =d
& &\quad&\text{on $\Gamma$}, \\
(\mu\bD(\bv) - \fq\bI)\bn - (\CF_2\rho+\sigma\Delta_\Gamma\rho)\bn
 = \bh& &\quad&\text{on $\Gamma$}.
\end{aligned}\right.\end{equation}
In the following, we consider Eq. \eqref{res:1.2} and 
\eqref{res:1.3} at the same time.  For this, we set 
$A_0 = 0$, and then Eq. \eqref{res:1.2} is represented
by Eq. \eqref{res:1.3} with $\kappa = 0$. 

We make a definition.
\begin{dfn}\label{def:2} Let $X$ and $Y$ be two Banach spaces.  A family
of operators $\CT \subset \CL(X, Y)$ is called $\CR$ bounded 
on $\CL(X, Y)$, if there exist constants $C > 0$ and 
$p \in [1, \infty)$ such that for each $n \in \BN$, 
$\{T_j\}_{j=1}^n \subset \CT$, and $\{f_j\}_{j=1}^n \subset X$, 
we have 
$$\|\sum_{k=1}^nr_kT_kf_k\|_{L_p((0, 1), Y)}
\leq C\|\sum_{k=1}^nr_kf_k\|_{L_p((0, 1), X)}.
$$
Here, the Rademacher functions $r_k$, $k \in \BN$, are given
by $r_k:[0,1] \to \{-1, 1\}$ $t \mapsto 
{\rm sign}\,(\sin 2^k\pi t)$.  The smallest such $C$ is called
$\CR$-bound of $\CT$ on $\CL(X, Y)$, which is denoted by  $\CR_{\CL(X, Y)}\CT$.
\end{dfn}
We introduce the definition of the uniqueness of solutions.
For Eq. \eqref{res:1.1}, the uniqueness is defined as follows:
\begin{itemize}
\item~ Let $\lambda \in U \subset \BC$.
If $\bu$ and $\fq$ with
$$\bu \in H^2_q(\Omega)^N, \quad \fp
\in H^1_q(\Omega) + \hat H^1_{q,0}(\Omega)$$
satisfy the homogeneous equations:
\begin{equation}\label{homoeq:1.2* }
\left\{\begin{aligned}
\lambda\bu - \DV(\mu\bD(\bu) - \fp\bI) = 0,
\quad \dv\bu &=0
&\quad &\text{in $\Omega$}, \\
(\mu\bD(\bu) - \fp\bI)\bn & = 0&\quad &\text{on $\Gamma$},
\end{aligned}\right.
\end{equation}
then $\bu=0$ and $\fp=0$.
\end{itemize}
For Eq. \eqref{res:1.2} and \eqref{res:1.3}, the uniqueness is defined as 
follows:
\begin{itemize}
\item~ Let $\lambda \in U \subset \BC$.
If $\bv$, $\fq$ and $\rho$ with
$$\bv \in H^2_q(\Omega)^N, \quad \fq
\in H^1_q(\Omega) + \hat H^1_{q,0}(\Omega), \quad
\rho \in W^{3-1/q}_q(\Gamma)$$
satisfy the homogeneous equations:
\begin{equation}\label{homoeq:1.2}
\left\{\begin{aligned}
\lambda\bv - \DV(\mu\bD(\bv) - \fq\bI) = 0,
\quad \dv\bv &=0
&\quad &\text{in $\Omega$}, \\
\lambda \rho + A_\kappa\cdot\nabla'_\Gamma \rho-\bv\cdot\bn + \CF_1\bv
& = 0 &\quad &\text{on $\Gamma$}, \\
(\mu\bD(\bv) - \fq\bI-((\CF_2 + \sigma\Delta_\Gamma)\rho)\bI)
\bn & = 0&\quad &\text{on $\Gamma$},
\end{aligned}\right.
\end{equation}
then $\bu=0$, $\fq=0$ and $\rho=0$.
\end{itemize}

We have the following theorems.
\begin{thm}\label{thm:rbdd:1.2} 
Let $1 < q < \infty$ and $0 < \epsilon_0 < \pi/2$.
Assume that $\Omega$ is a uniform $C^3$ domain and that the weak
Dirichlet problem is uniquely solvable for $q$.
Let $A_0=0$ and let $A_\kappa$ $(\kappa \in (0, 1))$ be
an $N-1$ vector of real valued functions satisfying 
\eqref{assump:3}.  
\\
\thetag1 $($Existence$)$~ 
Let   
\begin{align*}
X_q(\Omega) & = \{\bF = (\bff, d, \bh, g, \bg) \mid \bff
\in L_q(\Omega)^N, \enskip (g, \bg) \in DI_q(\Omega), \\
&\qquad \bh \in H^1_q(\Omega)^N, 
 \enskip d \in W^{2-1/q}_q(\Gamma)\}, \\
\CX_q(\Omega) & = \{F = (F_1, F_2, \ldots, F_7) \mid 
F_1, F_3, F_7 \in L_q(\Omega)^N, \enskip F_4 \in H^1_q(\Omega)^N, \\
&\qquad F_5 \in L_q(\Omega), \enskip F_6 \in H^1_q(\Omega),
\quad F_2 \in W^{2-1/q}_q(\Gamma)\}, \\
\Lambda_{\kappa, \lambda_0}
&= \begin{cases} \Sigma_{\epsilon_0, \lambda_0}
\quad&\text{for $\kappa = 0$}, \\
\BC_{+,\lambda_0} \quad&\text{for $\kappa \in (0, 1)$},
\end{cases}
\quad
\gamma_\kappa = \begin{cases}
1 \quad&\text{for $\kappa=0$}, \\
\kappa^{-b} \quad&\text{for $\kappa
\in (0, 1)$}.
\end{cases}
\end{align*}
Then, there exist a constant $\lambda_0 > 0$ 
 and operator families $\CA_1(\lambda)$, $\CP_1(\lambda)$ 
and $\CH_1(\lambda)$ with
\begin{align*}
\CA_1(\lambda) &\in \Hol(\Lambda_{\kappa, \lambda_0\gamma_\kappa},
\CL(\CX_q(\Omega), H^2_q(\Omega)^N)), \\
\CP_1(\lambda) &\in \Hol(\Lambda_{\kappa, \lambda_0\gamma_\kappa},
\CL(\CX_q(\Omega), H^1_q(\Omega)+\hat H^1_{q,0}(\Omega))), \\
\CH_1(\lambda) &\in \Hol(\Lambda_{\kappa, \lambda_0\gamma_\kappa},
\CL(\CX_q(\Omega), H^{3}_q(\Omega)))
\end{align*}
such that for every $\lambda \in \Lambda_{\kappa, \lambda_0\gamma_\kappa}$ and 
$(\bff,d, \bh, g, \bg)\in X_q(\Omega)$, 
$\bv=\CA_1 \bF_\lambda$ and $\fq 
= \CP_1(\lambda)\bF_\lambda$, and $\rho = \CH_1(\lambda)
\bF_\lambda$ are solutions of Eq. \eqref{res:1.3}, where
$$\bF_\lambda = (\bff,d, \lambda^{1/2}\bh, \bh,
 \lambda^{1/2}g, g, \lambda\bg).$$

Moreover, we have 
\begin{alignat*}2
\CR_{\CL(\CX_q(\Omega), H^{2-j}_q(\Omega)^N)}
(\{(\tau\pd_\tau)^\ell(\lambda^{j/2}\CA_1(\lambda))\mid
\lambda \in \Lambda_{\kappa, \lambda_0\gamma_\kappa}\}) &\leq r_b; \\
\CR_{\CL(\CX_q(\Omega), L_q(\Omega)^N)}
(\{(\tau\pd_\tau)^\ell(\nabla\CP_1(\lambda))\mid
\lambda \in \Lambda_{\kappa, \lambda_0\gamma_\kappa}\}) &\leq r_b; \\
\CR_{\CL(\CX_q(\Omega), H^{3-k}_q(\Omega))}
(\{(\tau\pd_\tau)^\ell(\lambda^k\CH_1(\lambda))\mid
\lambda \in \Lambda_{\kappa, \lambda_0\gamma_\kappa}\}) &\leq r_b
\end{alignat*}
for $\ell=0,1$, $j=0,1,2$, and 
$k=0,1$ with some constant $r_b > 0$. 
\vskip0.5pc\noindent
\thetag2 $($Uniqueness$)$~\thetag{i} 
When $\kappa=0$, if we assume that the weak Dirichlet problem
is uniquely solvable for $q'=q/(q-1)$ in addition,
then the uniqueness for Eq. \eqref{res:1.2} holds.
\\
\thetag{ii} When $\kappa \in (0, 1)$, if we assume that
 $\Omega$ is a uniformly $C^3$ domain whose
inside has a finite covering and that the weak Dirichlet
problem is uniquely solvable for $q'=q/(q-1)$ in addition,
 then the uniqueness for Eq. \eqref{res:1.3} holds. 
\end{thm} 
Here and in the following,
$\lambda = \gamma + i\tau \in \BC$. 
\begin{rem}\label{rem:3.3}
\thetag1 $F_1$, $F_2$, $F_3$, $F_4$, $F_5$, $F_6$,
and $F_7$  are 
variables corresponding to $\bff$,  $d$, 
$\lambda^{1/2}\bh$, $\bh$, $\lambda^{1/2}g$, $g$,
and, $\lambda\bg$,  respectively.\\
\thetag2 
We define the norms $\|\cdot\|_{X_q(\Omega)}$ and  
$\|\cdot\|_{\CX_q(\Omega)}$ by
\begin{align*}
\|(\bff,d, \bh, g, \bg)\|_{X_q(\Omega)}& = 
\|(\bff, \bg)\|_{L_q(\Omega)} + \|(g, \bh)\|_{H^1_q(\Omega)}
+ \|d\|_{W^{2-1/q}_q(\Gamma)} ; \\
\|(F_1, \ldots, F_7)\|_{\CX_q(\Omega)}
&= \|(F_1, F_3, F_7)\|_{L_q(\Omega)} + \|(F_4, F_6)\|_{H^1_q(\Omega)}
\\
&\phantom{\|(F_1, F_3, F_7)\|_{L_q(\Omega)}aaaaaaaaaa }
+ \|F_2\|_{W^{2-1/q}_q(\Gamma)}.
\end{align*}

\end{rem} 
\begin{thm}\label{thm:rbdd:1.1} 
Let $1 < q < \infty$ and $0 < \epsilon_0 < \pi/2$.\\
\thetag1 $($Existence$)$~ 
Assume that $\Omega$ is a uniform $C^2$ domain and that the weak
Dirichlet problem is uniquely solvable for $q$.  Let 
\begin{align*}
\tilde X_q(\Omega) & = \{\tilde\bF = (\bff, \bh, g, \bg) \mid \bff
\in L_q(\Omega)^N, \enskip (g, \bg) \in DI_q(\Omega), \enskip
\bh \in H^1_q(\Omega)^N\}, \\
\tilde\CX_q(\Omega) & = \{\tilde F = (F_1, F_3, \ldots, F_7) \mid 
F_1, F_3, F_7 \in L_q(\Omega)^N, \enskip F_3 \in L_q(\Omega), \\
&\qquad F_4 \in H^1_q(\Omega)^N, \enskip F_6 \in H^1_q(\Omega)\}.
\end{align*}
Then, there exist a constant $\lambda_0 > 0$ 
 and operator families $\CA_2(\lambda)$ and  $\CP_2(\lambda)$ with
\begin{align*}
\CA_2(\lambda) &\in \Hol(\Sigma_{\epsilon_0, \lambda_0},
\CL(\tilde \CX_q(\Omega), H^2_q(\Omega)^N)), \\
\CP_2(\lambda) &\in \Hol(\Sigma_{\epsilon_0, \lambda_0},
\CL(\tilde\CX_q(\Omega), H^1_q(\Omega)+\hat H^1_{q,0}(\Omega))),
\end{align*}
such that for every $\lambda \in \Sigma_{\epsilon_0, \lambda_0}$ and 
$(\bff, \bh, g, \bg) \in \tilde X_q(\Omega)$, 
$\bu=\CA_2(\lambda)\tilde\bF_\lambda$ and $\fp
= \CP_2(\lambda)\tilde\bF_\lambda$ are
 solutions of Eq. \eqref{res:1.1}, where
we have set 
$$\tilde\bF_\lambda = (\bff,\lambda^{1/2}\bh, \bh, g, \lambda^{1/2}g, 
\lambda\bg).$$ 

Moreover, we have 
\begin{alignat*}2
\CR_{\CL(\tilde\CX_q(\Omega), H^{2-j}_q(\Omega)^N)}
(\{(\tau\pd_\tau)^\ell(\lambda^{j/2}\CA_2(\lambda))\mid
\lambda \in \Sigma_{\epsilon_0, \lambda_0}\}) &\leq r_b; \\
\CR_{\CL(\tilde\CX_q(\Omega), L_q(\Omega)^N)}
(\{(\tau\pd_\tau)^\ell(\nabla\CP_2(\lambda))\mid
\lambda \in \Sigma_{\epsilon_0, \lambda_0}\}) &\leq r_b
\end{alignat*}
for $\ell=0,1$, $j=0,1,2$ with
some constant $r_b > 0$.  
\vskip0.5pc\noindent
$($Uniqueness$)$~ If we assume that the weak Dirichlet problem is 
uniquely solvable for $q' = q/(q-1)$ in addition, then the uniqueness
for Eq. \eqref{res:1.1} holds. 
\end{thm}
\subsection{Stokes operator and reduced Stokes operator}\label{subsec:3.3}

Since the pressure term has no time evolution, sometimes
it is convenient to eliminate the pressure terms,
for example to formulate the problem in the semigroup setting.  
For $\bu \in H^2_q(\Omega)^N$ and $h \in H^3_q(\Omega)$,
we introduce  functionals $K_0(\bu)$ and  $K(\bu, h)$. 
Let $K_0(\bu)  \in 
H^1_q(\Omega) + \hat H^1_{q,0}(\Omega)$ be a  unique 
solution of the weak Dirichlet problems:
\begin{equation}\label{wd:3.1.1} (\nabla K_0(\bu), \nabla\varphi)_\Omega
= (\DV(\mu\bD(\bu)) - \nabla\dv\bu, \nabla\varphi)_\Omega
\end{equation}
for any $\varphi \in \hat H^1_{q',0}(\Omega)$,
subject to 
$$K_0(\bu) = \mu<\bD(\bu)\bn, \bn> - \dv \bu\quad 
\text{on $\Gamma$}.$$ 
And, let $K(\bu, h)  \in 
H^1_q(\Omega) + \hat H^1_{q,0}(\Omega)$ be a unique solution of 
the weak Dirichlet problem:
\begin{equation}\label{wd:3.1.2} (\nabla K(\bu, h), \nabla\varphi)_\Omega
= (\DV(\mu\bD(\bu)) - \nabla\dv\bu, \nabla\varphi)_\Omega
\end{equation}
for any $\varphi \in \hat H^1_{q',0}(\Omega)$,
subject to 
$$K(\bu, h) = \mu<\bD(\bu)\bn, \bn>-(\CF_2h+\sigma \Delta_\Gamma h)
- \dv \bu\quad
\text{on $\Gamma$}.$$ 
By Remark \ref{rem:3.2}, we know the unique existence of 
$K_0(\bu)$ and $K(\bu, h)$ satisfying the estimates:
\begin{equation}\label{wd:4}\begin{aligned}
\|\nabla K_0(\bu)\|_{L_q(\Omega)}
&\leq C\|\nabla \bu\|_{H^1_q(\Omega)}, \\ 
\|\nabla K(\bu, h)\|_{L_q(\Omega)} &\leq C
(\|\nabla \bu\|_{H^1_q(\Omega)}
+ \|h\|_{W^{3-1/q}_q(\Gamma)})
\end{aligned}\end{equation}
for some constant depending on $C$.  We consider the reduced Stokes
equations:
\begin{equation}\label{rres:1.1} 
\left\{\begin{aligned}
\lambda \bu - \DV(\mu\bD(\bu)- K_0(\bu)\bI) & = \bff
&\quad&\text{in $\Omega$}, \\
(\mu\bD(\bu) - K_0(\bu)\bI)\bn
&=\bh&\quad&\text{on $\Gamma$}; 
\end{aligned}\right.\end{equation}
\begin{equation}\label{rres:1.2} 
\left\{\begin{aligned}
\lambda \bu - \DV(\mu\bD(\bu, h)-K(\bu, h)\bI) & = \bff
&\quad&\text{in $\Omega$}, \\
\lambda h + A_\kappa\cdot\nabla'_\Gamma h
- \bn\cdot\bu + \CF_1\bu& = d &\quad&\text{on $\Gamma$}, \\
(\mu\bD(\bu) - K(\bu, h)\bI)\bn - (\CF_2h+\sigma\Delta_\Gamma h)\bn
&=\bh&\quad&\text{on $\Gamma$}.
\end{aligned}\right.\end{equation}
Notice that both of 
the boundary conditions in Eq. \eqref{rres:1.1} and \eqref{rres:1.2}
are equivalent to 
\begin{equation}\label{rres:2}
(\mu\bD(\bu)\bn)_\tau = \bh_\tau\quad\text{and}\quad
\dv \bu = \bn\cdot\bh \quad\text{on $\Gamma$}.
\end{equation}
We now study the equivalence between Eq. \eqref{res:1.1}
and Eq. \eqref{rres:1.1}. The equivalence between Eq.\eqref{res:1.3}
and Eq.\eqref{rres:1.2} are similarly studied. 
We first assume that Eq. \eqref{res:1.1} is 
uniquely solvable. Given $\bff \in L_q(\Omega)^N$ in the right side of Eq.
\eqref{rres:1.1},  
let $g \in H^1_q(\Omega)$ be a unique solution of the variational equation:
\begin{equation}\label{wd:5} 
\lambda(g, \varphi)_\Omega + (\nabla g, \nabla\varphi)_\Omega
= (-\bff, \nabla\varphi)_\Omega
\quad\text{for any $\varphi \in H^1_{q',0}(\Omega)$}
\end{equation}
subject to $g = \bn\cdot\bh$ on $\Gamma$.  The unique existence of 
$g$ is guaranteed for $\lambda \in \Sigma_{\epsilon_0, \lambda_0}$ with 
some large $\lambda_0>0$. From \eqref{wd:5} it follows that
\begin{equation}\label{wd:5*}
(g, \varphi)_\Omega = 
(-\lambda^{-1}(\bff + \nabla g), \nabla\varphi)_\Omega,
\end{equation}
and so $(g, \bg) \in DI_q$ with $\bg = \lambda^{-1}(\bff + \nabla g)$.
Thus, from the assumption we know that Eq. \eqref{res:1.1} admits  
unique existence of solutions $\bu \in H^2_q(\Omega)^N$ and  
$\fq \in H^1_q(\Omega) + \hat H^1_{q,0}(\Omega)$. 
In view of the first equation in  Eq. \eqref{res:1.1} and Definition
\ref{dfn:1.4},  for any
$\varphi \in \hat H^1_{q',0}(\Omega)$ we have
\begin{align*}
(\bff, \nabla\varphi)_\Omega& = (\lambda\bu - \DV(\mu\bD(\bu) - \fp\bI), 
\nabla\varphi)_\Omega \\
& = (\lambda\bu, \nabla\varphi)_\Omega-(\nabla\dv\bu, \nabla\varphi)_\Omega
\\
&\qquad
-(\DV(\mu\bD(\bu))-\nabla\dv\bu, \nabla\varphi)_\Omega
+ (\nabla\fp, \nabla\varphi)_\Omega \\
& = \lambda(\lambda^{-1}(\bff + \nabla g), \nabla\varphi)_\Omega
-(\nabla g, \nabla\varphi)_\Omega 
+ (\nabla(\fp - K_0(\bu)), \nabla\varphi)_\Omega,
\end{align*}
and so, 
$$(\nabla(\fp-K_0(\bu)), \nabla\varphi)_\Omega = 0\quad
\text{for any $\varphi \in \hat H^1_{q',0}(\Omega)$}.
$$
Moreover, by the second equation in Eq. \eqref{res:1.1} 
and \eqref{wd:5}, we
$$\fp - K_0(\bu) = -\bh\cdot\bn +\dv\bu = -g + g =0 
\quad\text{on $\Gamma$}.
$$
Thus, the uniqueness implies that $\fp = K_0(\bu)$, which, combined
with Eq. \eqref{res:1.1},  shows that 
$\bu$ satisfy Eq. \eqref{rres:1.1}. 

Conversely, we assume that Eq. \eqref{rres:1.1} is uniquely
solvable.  Let $\bff \in L_q(\Omega)^N$ and $\bh \in H^1_q(\Omega)$.
Let $\theta \in H^1_q(\Omega) + \hat H^1_q(\Omega)$ be a unique solution
of the weak Dirichlet problem:
$$(\nabla \theta, \nabla\varphi)_\Omega = (\bff, \nabla\varphi)_\Omega
\quad\text{for any $\varphi \in \hat H^1_{q',0}(\Omega)$}, 
$$
subject to $\theta = \bn\cdot\bh$ on $\Gamma$. Setting $\fp = \theta
+ \fq$ in \eqref{res:1.1}, we then have 
\begin{alignat*}2
\lambda\bu - \DV(\mu\bD(\bu) - \fq\bI)
 = \bff-\nabla\theta, \quad \dv\bu =g = \dv \bg& &\quad&\text{in $\Omega$}, \\
(\mu\bD(\bu) - \fq\bI)\bn
= \bh-<\bh, \bn>\bn&
&\quad&\text{on $\Gamma$}.
\end{alignat*}
Let $\bff' = \bff-\nabla\theta$ and $\bh' = 
\bh-<\bh, \bn>\bn$.  We then have
\begin{equation}\label{rres:3}
\bh'\cdot\bn=0\quad\text{on $\Gamma$}\quad\text{and}\quad
(\bff', \nabla\varphi)_\Omega=0 \quad\text{for any
$\varphi \in \hat H^1_{q',0}(\Omega)$}. 
\end{equation}
Given $(g, \bg) \in DI_q(\Omega)^N$, let $K \in H^1_q(\Omega) + \hat H^1_{q,0}(\Omega)$
be a solution to the weak Dirichlet problem:
\begin{equation}\label{rres:4}
(\nabla K, \nabla\varphi)_\Omega = (\lambda\bg-\nabla g, \nabla\varphi)_\Omega
\quad\text{for any $\varphi \in \hat H^1_{q', 0}(\Omega)$},
\end{equation}
subject to $K = -g$ on $\Gamma$. 
Let $\bu$  be a solution of the equations:
\begin{equation}\label{rres:5}\left\{\begin{aligned}
\lambda\bu - \DV(\mu\bD(\bu) - K_0(\bu)\bI) & = \bff' + \nabla K&\quad
&\text{in $\Omega$}, \\
(\mu\bD(\bu) - K_0(\bu)\bI)\bn& =\bh' + g\bn
&\quad&\text{on $\Gamma$}.
\end{aligned}\right.\end{equation}
 By \eqref{wd:3.1.1}, \eqref{rres:3},
\eqref{rres:4} and \eqref{rres:5}, for any $\varphi
\in \hat H^1_{q,0}(\Omega)$ we have
\begin{equation}\label{rres:6}\begin{split} 
&(\nabla K, \nabla\varphi)_\Omega = (\bff' + \nabla K, \nabla\varphi)_\Omega\\
& = (\lambda\bu - \DV(\mu\bD(\bu) - K_0(\bu)\bI), \nabla\varphi)_\Omega \\
& = (\lambda \bu, \nabla\varphi)_\Omega-(\nabla\dv\bu, \nabla\varphi)_\Omega.
\end{split}\end{equation}
Since $H^1_{q',0}(\Omega)\subset \hat H^1_{q',0}(\Omega)$, 
by \eqref{rres:4} and \eqref{rres:6}
we have
$$\lambda(\dv \bu - g, \varphi)_\Omega + (\nabla(\dv\bu - g), 
\nabla\varphi)_\Omega = 0
\quad\text{for any $\varphi \in H^1_{q',0}(\Omega)$}.
$$
Here, we have used the fact that $\dv\bg = g \in H^1_q(\Omega)$.  
Recalling that $\bh'\cdot\bn=0$ and 
putting \eqref{rres:5} and the boundary condition of \eqref{wd:3.1.1}
 together gives 
$$g = (g\bn + \bh')\cdot\bn = <\mu\bD(\bu)\bn, \bn> - K_0(\bu)= \dv \bu$$
on $\Gamma$.  Thus, the uniqueness implies that $\dv\bu = g$ in
$H^1_q(\Omega)$, which, combined with \eqref{rres:4}
and \eqref{rres:6}, leads to 
$$(\bg, \nabla\varphi)_\Omega = (\bu, \nabla\varphi)_\Omega
\quad\text{for any $\varphi \in \hat H^1_{q', 0}(\Omega)$}
$$
because we may assume that $\lambda\not=0$.  
Since $K = - g$ on $\Gamma$, 
by \eqref{rres:6} we have
\begin{alignat*}2
\lambda\bu - \DV(\mu\bD(\bu) - (K_0(\bu)-K)\bI)
 = \bff', \quad \dv\bu = g = &\dv \bg &\quad&\text{in $\Omega$}, \\
(\mu\bD(\bu) - (K_0(\bu) - K)\bI)\bn
& = \bh'
&\quad&\text{on $\Gamma$}.
\end{alignat*}
 Thus, $\bu$ and $\fp = K_0(\bu) - K
-\theta$
 are required solutions of Eq. \eqref{res:1.1}.
\subsection{$\CR$-bounded solution operators for the reduced Stokes equations}
\label{subsec:2.2}
In the following theorem, we state
 the existence of $\CR$ bounded solution operators 
for the reduced Stokes equations \eqref{rres:1.2}. 
\begin{thm}\label{main:thm3}
Let $1 < q < \infty$ and $0 < \epsilon < \pi/2$.
Let $\Lambda_{\kappa, \lambda_0}$ be 
the set defined in Theorem \ref{thm:rbdd:1.2}.  Assume that 
$\Omega$ is a uniform $C^3$ domain and the weak Dirichlet problem
is uniquely solvable for $q$.  Let 
$A_0=0$ and $A_\kappa$ $(\kappa \in (0, 1))$ be an 
$N-1$ vector of real valued functions satisfying \eqref{assump:3}. 
Set  
\begin{equation}\label{rdata:0}\begin{split}
Y_q(\Omega) &= \{(\bff, d, \bh) \mid \bff\in L_q(\Omega)^N, 
 \quad  d \in W^{2-1/q}_q(\Gamma), \enskip 
\bh \in H^1_q(\Omega)^N\}, \\
\CY_q(\Omega)& = \{(F_1, \ldots, F_4) \mid
F_1, F_3 \in L_q(\Omega)^N, \enskip 
F_2 \in W^{2-1/q}_q(\Gamma)\\
&\phantom{= \{(F_1, \ldots, F_4) \mid
F_1, F_3 \in L_q(\Omega)^N,}\quad 
 F_4 \in H^1_q(\Omega)^N\}.
\end{split}\end{equation}
Then, there exist a constant $\lambda_* \geq 1$ 
and operator families: 
\begin{align*}
 \CA_r(\lambda) &\in 
{\rm Hol}\, (\Lambda_{\kappa, \lambda_*\gamma_\kappa}, 
\CL(\CY_q(\Omega), H^2_q(\Omega)^N)), \\
\CH_r(\lambda) &\in {\rm Hol}\, 
(\Lambda_{\kappa, \lambda_*\gamma_\kappa},
\CL(\CY_q(\Omega), H^3_q(\Omega)))
\end{align*}
such that for any $\lambda \in \Lambda_{\kappa, \lambda_0\gamma_\kappa}$ and 
$(\bff, d, \bh) \in Y_q(\Omega)$,  
$$\bu = \CA_r(\lambda)(\bff, d, \lambda^{1/2}\bh, \bh), \quad
h = \CH_r(\lambda)(\bff, d, \lambda^{1/2}\bh, \bh),$$ 
are  solutions of \eqref{rres:1.2},
and 
\begin{align*}
&\CR_{\CL(\CY_q(\Omega), H^{2-j}_q(\Omega)^N)}
(\{(\tau\pd_\tau)^\ell(\lambda^{j/2}\CA_r(\lambda)) \mid 
\lambda \in \Lambda_{\kappa, \lambda_*\gamma_\kappa}\}) \leq 
r_b, \\
&\CR_{\CL(\CY_q(\Omega), H^{3-k}_q(\Omega))}
(\{(\tau\pd_\tau)^\ell(\lambda^{k}\CH_r(\lambda)) \mid 
\lambda \in \Lambda_{\kappa, \lambda_*\gamma_\kappa}\}) \leq 
r_b, 
\end{align*}
for $\ell=0,1$, $j=0,1,2$ and $k=0,1$.  Here, $r_b$ is a
constant depending on $m_1$, $m_2$, $m_3$, $\lambda_0$,  $p$, $q$, and $N$, 
but independent of $\kappa \in (0, 1)$, and $\gamma_\kappa$ is the number
defined in Theorem \ref{thm:rbdd:1.2}.  

If we assume that $\Omega$ is a uniformly $C^3$ domain whose
inside has a finite covering and that the weak Dirichlet problem
is uniquely solvable for $q'=q/(q-1)$ in addition, then the uniqueness 
for Eq. \eqref{rres:1.2} holds. 

\end{thm}
\begin{remark}\label{rem:1} 
The norm of space $Y_q(\Omega)$ is defined by
$$\|(\bff, d, \bh)\|_{Y_q(\Omega)} = \|\bff\|_{L_q(\Omega)}
+ \|d\|_{W^{2-1/q}_q(\Gamma)} + \|\bh\|_{H^1_q(\Omega)};$$
and the norm of space $\CY_q(\Omega)$ is defined by
$$\|(F_1, F_2, F_3,  F_4)\|_{\CY_q(\Omega)}
= \|(F_1, F_3)\|_{L_q(\Omega)} +\|F_2\|_{W^{2-1/q}_q(\Gamma)}
+ \|F_4\|_{W^1_q(\Omega)}.
$$
\end{remark}
\begin{remark}\label{rem:2} 
By the equivalence of Eq.\eqref{res:1.3} and Eq. \eqref{rres:1.2}
that was pointed out in Subsec. \ref{subsec:3.3}, we see easily that 
Theorem \ref{thm:rbdd:1.2} follows immediately from 
Theorem \ref{main:thm3}.
\end{remark}
Concerning the existence of $\CR$ bounded solution operator for
Eq. \eqref{rres:1.1}, we have the following theorem. 
\begin{thm}\label{main:thm0}
Let $1 < q < \infty$ and $0 < \epsilon < \pi/2$.
  Assume that 
$\Omega$ is a uniform $C^2$ domain and the weak Dirichlet problem
is uniquely solvable for $q$.  
Set  
\begin{equation}\label{rdata:00}\begin{split}
\tilde\CY_q(\Omega) &= \{(\bff, \bh) \mid \bff\in L_q(\Omega)^N, \enskip 
\bh \in H^1_q(\Omega)^N\}, \\
\tilde \CY_q(\Omega)& = \{(F_1, F_3, F_4) \mid
F_1, F_3 \in L_q(\Omega)^N, \enskip 
 F_4 \in H^1_q(\Omega)^N\}.
\end{split}\end{equation}
Then, there exist a constant $\lambda_{**} \geq 1$ 
and operator family  $ \CA^0_r(\lambda)$ with
$$
 \CA^0_r(\lambda) \in 
{\rm Hol}\, (\Sigma_{\sigma, \lambda_{**}}, 
\CL(\tilde\CY_q(\Omega), H^2_q(\Omega)^N))$$
such that for any $\lambda \in \Sigma_{\sigma, \lambda_{**}}$ and 
$(\bff, \bh) \in \tilde Y_q(\Omega)$,  
$\bu = \CA_r^0(\lambda)(\bff, \lambda^{1/2}\bh, \bh)$ 
is a  solution of \eqref{rres:1.1},
and 
\begin{align*}
&\CR_{\CL(\tilde\CY_q(\Omega), H^{2-j}_q(\Omega)^N)}
(\{(\tau\pd_\tau)^\ell(\lambda^{j/2}\CA_r^0(\lambda)) \mid 
\lambda \in \Sigma_{\sigma, \lambda_{**}}\}) \leq 
r_b 
\end{align*}
for $\ell=0,1$ and $j=0,1,2$.  Here, $r_b$ is a
constant depending on $m_1$,  $\lambda_0$,  $p$, $q$, and $N$. 

If we assume that the weak Dirichlet problem is uniquely solvable
for $q'=q/(q-1)$ in addition, then the uniqueness for Eq. \eqref{rres:1.1}
holds. 
\end{thm}
\begin{remark} Theorem \ref{main:thm0} was proved by Shibata \cite{S1}
and can be proved by the same argument as in the proof of
Theorem \ref{main:thm3}.  Thus, we may omit its proof. 
\end{remark}

\subsection{Generation of $C^0$ analytic semigroup}\label{subsec:3.4}

In this subsection, we consider the following two initial-boundary value
problems for the reduced Stokes operator:
\begin{equation}\label{rstokes:1.1}\left\{
\begin{aligned}
\pd_t\bu - \DV(\mu\bD(\bu)-K_0(\bu)\bI) & = 0
&\quad&\text{in $\Omega^\infty$}, \\
(\mu\bD(\bu) - K_0(\bu)\bI)\bn & = 0
&\quad&\text{on $\Gamma^\infty$}, 
\\
\bu|_{t=0} & =\bu_0
&\quad&\text{in $\Omega$};
\end{aligned}\right.\end{equation}
\begin{equation}\label{rstokes:1.2}\left\{
\begin{aligned}
\pd_t\bv - \DV(\mu\bD(\bv)-K(\bv, h)\bI)  = 0&
&\enskip&\text{in $\Omega^\infty$}, \\
\pd_t h - \bn\cdot\bv + \CF_1\bv= 0& 
&\enskip&\text{on $\Gamma^\infty$}, \\
(\mu\bD(\bv) - K(\bv, h)\bI)\bn
-(\CF_2 h + \sigma\Delta_\Gamma h)\bn = 0&
&\enskip&\text{on $\Gamma^\infty$}, 
\\
\bv|_{t=0}=\bv_0\enskip\text{in $\Omega$}, \quad
h|_{t=0} = h_0\enskip \text{on $\Gamma$}&,
\end{aligned}\right.\end{equation}
where we have set $\Omega^\infty = \Omega\times(0, \infty)$ and 
$\Gamma^\infty=\Gamma\times(0, \infty)$. 

Let $\bu$ and $(\bv, h)$ be solutions of Eq.\eqref{rstokes:1.1}
and Eq.\eqref{rstokes:1.2}, respectively. Roughly speaking,
if $\bu \in J_q(\Omega)$
for any $t > 0$, then, $\bu$ and $\fp = K_0(\bu)$ are 
unique solutions of Eq. \eqref{stokes:1} with $g = \bg = \bh=0$.
And, if $\bv \in J_q(\Omega)$ for any $t > 0$, then 
$\bv$, $\fq = K(\bv, h)$, and $h$ are unique 
solutions of \eqref{stokes:2} with $d=0$ and $(\bv, h)|_{t=0}
= (\bv_0, h_0)$. 

Let us introduce  spaces and  operators to describe
\eqref{rstokes:1.1} and \eqref{rstokes:1.2}
in the semigroup setting.  Let 
\begin{equation}\label{rstokes:2}\begin{split}
\CD^1_q(\Omega) & = \{\bu \in J_q(\Omega) \cap H^2_q(\Omega)^N \mid
(\mu\bD(\bu)\bn)_\tau = 0 \enskip \text{on $\Gamma$}\}, \\
\CA^1_q\bu &= \DV(\mu(\bD(\bu) - K_0(\bu)\bI) \quad\text{for 
$\bu \in \CD^1_q(\Omega)$}; \\
\CH_q(\Omega) & = \{ (\bv, h) \mid \bv \in J_q(\Omega),
\enskip h \in W^{2-1/q}_q(\Gamma)\}, \\
\CD^2_q(\Omega) & = \{(\bv, h) \in \CH_q(\Omega) \mid 
\bv \in H^2_q(\Omega)^N, \enskip h \in W^{3-1/q}_q(\Gamma), \\
&\qquad (\mu\bD(\bv)\bn)_\tau = 0 \quad\text{on $\Gamma$}\},\\
\CA^2_q(\bv, h) & = 
(\DV(\mu\bD(\bv) - K(\bv, h)\bI), \bn\cdot\bv|_\Gamma)
\end{split}\end{equation}
for $(\bv, h) \in \CD^2_q(\Omega)$.
Since $\dv\bu=0$ for $\bu \in \CD^1_q(\Omega)$, by \eqref{rres:2}
$(\mu\bD(\bu)\bn)_\tau = 0$ is equivalent to
$(\mu\bD(\bu) - K_0(\bu)\bI)\bn=0$.
And also, for $(\bv, h) \in \CD^2_q(\Omega)$, 
$(\mu\bD(\bv)\bn)_\tau = 0$ is equivalent to 
$$(\mu\bD(\bv) - K(\bv, h)\bI)\bn - 
(\CF_2h+\sigma\Delta_\Gamma h) = 0.$$

Using the symbols defined  
in \eqref{rstokes:2}, we see that Eq.\eqref{rstokes:1.1} 
and Eq.\eqref{rstokes:1.2} are 
written as 
\begin{align}
\pd_t\bu - \CA^1_q\bu & = 0 \quad(t > 0), \quad \bu|_{t=0} = \bu_0;
\label{rstokes:3.1}\\
\pd_tU(t)-\CA^2_qU(t) &= 0 \quad(t > 0), \quad U(t)|_{t=0} = U_0, 
\label{rstokes:3.2}
\end{align}
where $\bu(t) \in \CD^1_q(\Omega)$ for $t > 0$ and 
$\bu_0 \in J_q(\Omega)$, and 
$U(t) = (\bv, h) \in \CD^2_q(\Omega)$ for $t > 0$ and 
$U_0 = (\bv_0, h_0) \in \CH_q(\Omega)$. The corresponding
resolvent problem to \eqref{rstokes:3.1} is that for any
$\bff \in J_q(\Omega)$ and $\lambda \in \Sigma_{\epsilon_0, \lambda_0}$
we find $\bu \in \CD^1_q(\Omega)$ uniquely solving the equation:
\begin{equation}\label{*rstokes:4.1} 
\lambda \bu - \CA^1_q\bu = \bff \quad\text{in $\Omega$} 
\end{equation}
possessing the estimate:
\begin{equation}\label{rstokes:5.1}
|\lambda|\|\bu\|_{L_q(\Omega)} + \|\bu\|_{H^2_q(\Omega)}
\leq C\|\bff\|_{L_q(\Omega)}.
\end{equation}
And also, the corresponding
resolvent problem to \eqref{rstokes:3.2} is that for any
$F \in \CH_q(\Omega)$ and $\lambda \in \Sigma_{\epsilon_0, \lambda_0}$
we find $U \in \CD^2_q(\Omega)$ uniquely solving the equation:
\begin{equation}\label{rstokes:4.2} 
\lambda U - \CA^2_qU = F \quad\text{in $\Omega\times \Gamma$} 
\end{equation}
possessing the estimate:
\begin{equation}\label{rstokes:5.2}
|\lambda|\|U\|_{\CH_q(\Omega)} + \|U\|_{\CD^2_q(\Omega)}
\leq C\|F\|_{\CH_q(\Omega)},
\end{equation}
where for $U = (\bv, h)$ we have set
$$\|U\|_{\CH_q(\Omega)} = \|\bv\|_{L_q(\Omega)} + \|h\|_{W^{2-1/q}_q(\Gamma)},
\enskip
\|U\|_{\CD^2_q(\Omega)} = \|\bv\|_{H^2_q(\Omega)} 
+ \|h\|_{W^{3-1/q}_q(\Gamma)}.
$$
Since $\CR$ boundedness implies boundedness as follows
from Definition \ref{def:2} with $n=1$, by Theorem \ref{main:thm0},
we know the 
unique existence of $\bu \in \CD^1_q(\Omega)$ satisfying \eqref{*rstokes:4.1} 
and \eqref{rstokes:5.1}.  And also, by Theorem \ref{main:thm3},
 we know the unique existence of $U\in 
\CD^2_q(\Omega)$ satisfying \eqref{rstokes:4.2}
and \eqref{rstokes:5.2}.  
Thus, using standard semigroup theory,
we have the following theorem.
\begin{thm}\label{thm:semi:1.1}
Let $ 1 < q < \infty$.  Assume that $\Omega$ is a uniform
$C^2$ domain in $\BR^N$ and that the weak Dirichlet problem
is uniquely solvable for $q$ and $q' = q/(q-1)$.  
Then, problem \eqref{rstokes:3.1}
generates a $C^0$ analytic semigroup $\{T_1(t)\}_{t\geq 0}$ on 
$J_q(\Omega)$ satisfying the estimates:
\begin{gather}
\|T_1(t)\bu_0\|_{L_q(\Omega)} + t(\|\pd_tT_1(t)\bu_0\|_{L_q(\Omega)}
+ \|T_1(t)\bu_0\|_{H^2_q(\Omega)})
 \leq Ce^{\gamma t}\|\bu_0\|_{L_q(\Omega)},
\label{semi:1.1}\\
\|\pd_tT_1(t)\bu_0\|_{L_q(\Omega)}
+ \|T_1(t)\bu_0\|_{H^2_q(\Omega)} 
\leq Ce^{\gamma t}\|\bu_0\|_{H^2_q(\Omega)}
\label{semi:1.2}
\end{gather}
for any $t > 0$ with some constants $C > 0$ and $\gamma > 0$.
\end{thm}
\begin{thm}\label{thm:semi:1.2}
Let $ 1 < q < \infty$.  Assume that $\Omega$ is a uniform
$C^3$ domain in $\BR^N$ and that the weak Dirichlet problem
is uniquely solvable for $q$ and $q' = q/(q-1)$.  
Then, problem \eqref{rstokes:3.2}
generates a $C^0$ analytic semigroup $\{T_2(t)\}_{t\geq 0}$ on 
$\CH_q(\Omega)$ satisfying the estimates:
\begin{gather}
\|T_2(t)U_0\|_{\CH_q(\Omega)} + t(\|\pd_tT_2(t)U_0\|_{\CH_q(\Omega)}
+ \|T_2(t)U_0\|_{\CD^2_q(\Omega)}) 
\leq Ce^{\gamma t}\|U_0\|_{\CH_q(\Omega)},
\label{semi:2.1}\\
\|\pd_tT_2(t)U_0\|_{\CH_q(\Omega)}
+ \|T_2(t)U_0\|_{\CD^2_q(\Omega)} 
\leq Ce^{\gamma t}\|U_0\|_{\CD^2_q(\Omega)}
\label{semi:2.2}
\end{gather}
for any $t > 0$ with some constants $C > 0$ and $\gamma > 0$.
\end{thm}
  
We now show the following maximal $L_p$-$L_q$ regularity theorem 
for Eq. \eqref{rstokes:1.1} and Eq.\eqref{rstokes:1.2}. 
\begin{thm}\label{thm:max:2.1} 
Let $1 < p, q < \infty$.  Assume that $\Omega$ is a uniform $C^2$ domain
in $\BR^N$ and that the weak Dirichlet problem is uniquely solvable
for $q$ and $q' = q/(q-1)$.  Let $\CD^1_{q,p}(\Omega)$ be a subspace of
$B^{2(1-1/p)}_{q, p}(\Omega)^N$ defined by 
$$\CD^1_{q,p}(\Omega) = (J_q(\Omega), \CD^1_q(\Omega))_{1-1/p, p},$$
where $(\cdot, \cdot)_{1-1/p,p}$ denotes a real interpolation functor. 
Then, there exists a $\gamma > 0$ such that 
for any initial data $\bu_0
\in \CD^1_{q,p}(\Omega)$, problem \eqref{rstokes:1.1} admits a unique
solution $\bu$  with
$$
e^{-\gamma t}\bu \in H^1_p((0, \infty), L_q(\Omega)^N) \cap L_p((0, \infty), 
H^2_q(\Omega)^N),
$$
possessing the estimate:
\begin{equation}\label{sei:3.1}
\|e^{-\gamma t}\pd_t\bu\|_{L_p((0, \infty), L_q(\Omega))}
+ \|e^{-\gamma t}\bu\|_{L_p((0, \infty), H^2_q(\Omega))}
\leq C\|\bu_0\|_{B^{2(1-1/p)}_{q,p}(\Omega)}.
\end{equation}
Here, for any Banach space $X$ with norm $\|\cdot\|_X$ we have set
$$\|e^{-\gamma t}f\|_{L_p((a, b), X)} = 
\Bigl(\int^b_a(e^{-\gamma t}\|f(t)\|_X)^p\,dt\Bigr)^{1/p}.
$$
\end{thm}
\begin{rem}\label{rem:3.1.1}
Since $\CD^1_{q,p}(\Omega) \subset B^{2(1-1/p)}_{q,p}(\Omega)$, 
in view of a boundary trace theorem, we see that 
for $\bu_0 \in \CD^1_{q,p}(\Omega)$, we have
$$\left\{\begin{aligned}
\bu_0 &\in J_q(\Omega), \quad(\mu\bD(\bu_0)\bn)_\tau 
= 0 \quad\text{on $\Gamma$} &\quad&\text{for
$\frac2p + \frac1q < 1$}, \\
\bu_0 &\in J_q(\Omega)
&\quad&\text{for
$\frac2p + \frac1q > 1$},
\end{aligned}\right.
$$
because $\bD(\bu_0) \in B^{1-2/p}_{q,p}(\Omega)$, and so
$\bD(\bu_0)|_\Gamma$ exists for $\frac2p + \frac1q < 1$,
but it does not exist for $\frac2p + \frac1q > 1$.
\end{rem}
\begin{thm}\label{thm:max:2.2} 
Let $1 < p, q < \infty$.  Assume that $\Omega$ is a uniform $C^3$ domain
in $\BR^N$ and that the weak Dirichlet problem is uniquely solvable
for $q$ and $q' = q/(q-1)$.  Let $\CD^2_{q,p}(\Omega)$ be a subspace of
$B^{2(1-1/p)}_{q, p}(\Omega)^N\times B^{3-1/p-1/q}_{q,p}(\Gamma)$ defined by 
$$\CD^2_{q,p}(\Omega) = (\CH_q(\Omega), \CD^2_q(\Omega))_{1-1/p, p}.$$
Then, there exists a $\gamma > 0$ such that 
for any initial data $(\bu_0, h_0)
\in \CD^2_{q,p}(\Omega)$, problem \eqref{rstokes:1.2} admits a unique
solution $U = (\bv, h)$ with
\begin{align*}
e^{-\gamma t}U &\in H^1_p((0, \infty), \CH_q(\Omega)) 
\cap L_p((0, \infty), \CD^2_q(\Omega))
\end{align*}
possessing the estimate:
\begin{equation}\label{sei:3.2}\begin{split}
&\|e^{-\gamma t}\pd_tU\|_{L_p((0, \infty), \CH_q(\Omega))}
+ \|e^{-\gamma t}U\|_{L_p((0, \infty), \CD^2_q(\Omega))} \\
&\quad
\leq C(\|\bu_0\|_{B^{2(1-1/p)}_{q,p}(\Omega)} 
+ \|h_0\|_{B^{3-1/p-1/q}_{q,p}(\Gamma)}).
\end{split}\end{equation}
\end{thm}
{\bf Proof of Theorem \ref{thm:max:2.2}}. \quad 
We only prove Theorem \ref{thm:max:2.2}, because
Theorem \ref{thm:max:2.1} can be proved by the same
argument. 
To prove Theorem \ref{thm:max:2.2}, we observe that 
\begin{align*}
&\Bigl(\int^\infty_0(e^{-\gamma t}\|\pd_tT_2(t)U_0\|_{\CH_q(\Omega)}
)^p\,dt\Bigr)^{1/p} \\
&\quad= \Bigl(\sum_{j=-\infty}^\infty \int^{2^{j+1}}_{2^j}
(e^{-\gamma t}\|\pd_tT_2(t)U_0\|_{\CH_q(\Omega)})^p\,dt\\
&\quad\leq \Bigl(\sum_{j=-\infty}^\infty(2^{j+1}-2^j)
(\sup_{t\in(2^j, 2^{j+1})}(e^{-\gamma t}\|\pd_tT_2(t)U_0
\|_{\CH_q(\Omega)})^p\Bigr)^{1/p}.
\end{align*}
We now introduce 
Banach spaces $\ell^s_p$ for $s \in \BR$ and $1\leq p \leq \infty$ 
which are sets of all sequences, $(a_j)_{j \in \BZ}$, such that
$\|(a_j)_{j \in \BZ}\|_{\ell^s_p} < \infty$, where we have set
$$\|(a_j)_{j\in \BZ}\|_{\ell^s_p}
= \begin{cases} \Bigl(\sum_{j=1}^\infty(2^{js}|a_j|)^p\Bigr)^{1/p}
\quad&\text{for $1 \leq p < \infty$}, \\
\sup \{2^{js}|a_j| \mid j \in \BZ\}\quad&\text{for $p=\infty$}.
\end{cases}
$$
Let 
$$a_j = \sup_{t \in (2^j, 2^{j+1})}e^{-\gamma t}\|\pd_tT_2(t)U_0
\|_{\CH_q(\Omega)},
$$
and then 
\begin{equation}\label{proof:3.1}\begin{split}
&\Bigl(\int^\infty_0(e^{-\gamma t}\|\pd_tT_2(t)\bu_0\|_{\CH_q(\Omega)}
)^p\,dt\Bigr)^{1/p}\\
&\quad
\leq \Bigl(\sum_{j=-\infty}^\infty(2^{j/p}a_j)^p\Bigr)^{1/p}
= \|(a_j)_{j \in \BZ}\|_{\ell^{1/p}_p}.
\end{split}\end{equation}
By real interpolation theory (cf. Bergh and L\"ofstr\"om 
\cite[5.6.Theorem]{BL}), we have $\ell^{1/p}_p
= (\ell^1_\infty, \ell^0_\infty)_{1-1/p, p}$. 
Moreover, by \eqref{semi:2.1} and  \eqref{semi:2.2}, we have 
$$\|(a_j)_{j\in \BZ}\|_{\ell^1_\infty} \leq C\|U_0\|_{\CH_q(\Omega)},
\quad 
\|(a_j)_{j \in \BZ}\|_{\ell^0_\infty} \leq C\|U_0\|_{\CD_q(\Omega)},
$$
and therefore, by real interpolation 
$$\|(a_j)_{j\in\BZ}\|_{\ell^{1/p}_p}=\|(a_j)_{j\in \BZ}
\|_{(\ell^1_\infty, \ell^0_\infty)_{1-1/p, p}}
\leq C\|U_0\|_{(\CH_q(\Omega), \CD_q(\Omega))_{1-1/p, p}}.
$$
Putting this and \eqref{proof:3.1} together gives  
$$\|e^{-\gamma t}\pd_tT_2(t)U_0\|_{L_p((0, \infty), \CH_q(\Omega))}
\leq C\|U_0\|_{\CD_{q,p}(\Omega)}.
$$
Analogously, we have
$$\|e^{-\gamma t}T_2(t)U_0\|_{L_p((0, \infty), \CD_q(\Omega))}
\leq C\|U_0\|_{\CD_{q,p}(\Omega)}.
$$
We now prove the uniqueness.  Let $U$ satisfy the homogeneous equation:
\begin{equation}\label{uniq:1}
\pd_tU - \CA^2_qU = 0 \quad(t > 0), \quad U|_{t=0} = 0,
\end{equation}
and the condition:
\begin{equation}\label{uniq:2}
e^{-\gamma t}U \in H^1_p((0, \infty), \CH_q(\Omega)) \cap L_p((0, \infty), 
\CD^2_q(\Omega)).
\end{equation}
Let $U_0$ be the zero extension of $U$ to $t < 0$.  In particular, 
from \eqref{uniq:1} it follows that
\begin{equation}\label{uniq:3}
\pd_tU_0 - \CA^2_qU_0 = 0 \quad(t \in \BR).
\end{equation}
For any $\lambda \in \BC$ with ${\rm Re}\,\lambda > \gamma$, we set
$$\hat U(\lambda) = \int^\infty_{-\infty} e^{-\lambda t}U_0(t)\,dt
= \int^\infty_0 e^{-\lambda t}U(t)\,dt.
$$
By \eqref{uniq:2} and H\"older's inequality,  we have
\begin{align*}
\|\hat U(\lambda)\|_{\CD_q(\Omega)} & 
\leq \Bigl(\int^\infty_0e^{-({\rm Re}\,\lambda - \gamma)tp'}\,dt
\Bigr)^{1/{p'}}\|e^{-\gamma t}U\|_{L_p((0, \infty), \CD^2_q(\Omega))}\\
&\leq (({\rm Re}\,\lambda-\gamma)p')^{-1/{p'}}
\|e^{-\gamma t}U\|_{L_p((0, \infty), \CD^2_q(\Omega))}.
\end{align*}
Since $\lambda\hat U(\lambda) = \int^\infty_0e^{-\lambda t}\pd_tU(t)\,dt$, 
we also have
\begin{align*}
\|\lambda \hat U(\lambda)\|_{\CH_q(\Omega)} 
\leq (({\rm Re}\,\lambda-\gamma)p')^{-1/{p'}}
\|e^{-\gamma t}\pd_t U\|_{L_p((0, \infty), \CH_q(\Omega))}.
\end{align*}
Thus, by \eqref{uniq:3}, we have $\hat U(\lambda) \in \CD^2_q(\Omega)$
satisfies the homogeneous equation:
$$\lambda \hat U(\lambda) - \CA^2_q\hat U(\lambda) = 0 \quad
\text{in $\Omega\times\Gamma$}.
$$
Since the uniqueness of the resolvent problem holds for
$\lambda \in \Sigma_{\epsilon_0, \lambda_0}$, we have 
$\hat U(\lambda) = 0$ for any $\lambda \in \BC$
with ${\rm Re}\, \lambda > \max(\lambda_0, \gamma)$.  By the Laplace
inverse transform, we have $U_0(t) = 0$ for $t \in \BR$, 
that is $U(t) = 0$ for $t >0$, which shows the uniqueness.
This completes the proof of Theorem \ref{thm:max:2.2}.
\qed

\subsection{Proof of maximal regularity theorem}\label{subsec:3.5}

We first prove Theorem \ref{thm:max.1}
with the help of Theorem \ref{thm:rbdd:1.1}.  \vskip0.5pc
{\bf Proof of Theorem \ref{thm:max.1}}. \quad
The key tool in the proof
of Theorem \ref{thm:max.1} is the Weis operator valued Fourier multiplier
theorem.  To state it we need to make a few definitions. For a Banach 
space, $X$,  $\CD(\BR, X)$ denotes the space of $X$-valued
$C^\infty(\BR)$ functions with compact support  
$\CD'(\BR, X) = \CL(\CD(\BR), X)$  the space of 
$X$-valued distributions.  And also, 
$\CS(\BR, X)$ denotes the space of $X$-valued rapidly decreasing 
functions and $\CS'(\BR, X) = \CL(\CS(\BR), X)$ 
the space of $X$-valued tempered distributions. 
Let $Y$ be another Banach space.  Then, given $m
\in L_{1, {\rm loc}}(\BR, \CL(X, Y))$, we define an operator
$T_m: \CF^{-1}\CD(\BR, X) \to \CS'(\BR, Y)$ by letting
\begin{equation}\label{four:1} T_m\phi
= \CF^{-1}[m\CF[\phi]]
\quad\text{for all $\CF\phi \in \CD(\BR, X)$},
\end{equation}
where $\CF$ and $\CF^{-1}$ denote the Fourier transform and its
inversion formula, respectively. 
\begin{dfn}\label{def:3} A Banach space $X$ is said to be a UMD Banach
space, if the Hilbert transform is bounded on 
$L_p(\BR, X)$ for some (and then all) $p \in (1, \infty)$.  Here,
the Hilbert transform $H$ operating on $f \in \CS(\BR, X)$
is defined by 
$$[Hf](t) = \frac1\pi\lim_{\epsilon\to0+}
\int_{|t-s| > \epsilon} \frac{f(s)}{t-s}\,ds
\quad(t \in \BR).$$
\end{dfn}
\begin{thm}[Weis \cite{Weis}]\label{thm:weis} Let $X$ and $Y$ be
two UMD Banach spaces and $1 < p < \infty$.  Let $m$ be a function in 
$C^1(\BR\setminus\{0\}, \CL(X, Y))$ such that
\begin{align*}
\CR_{\CL(X, Y)}(\{m(\tau) \mid \tau \in \BR\setminus\{0\}\}) &= 
\kappa_0 < \infty, \\
\CR_{\CL(X, Y)}(\{\tau m'(\tau) \mid \tau \in \BR\setminus\{0\}\}) &= 
\kappa_1 < \infty.
\end{align*}
Then, the operator $T_m$ defined in \eqref{four:1} is 
extended to a bounded linear operator from $L_p(\BR, X)$ into
$L_p(\BR, Y)$.  Moreover, denoting this extension by $T_m$, we have
$$\|T_mf\|_{L_p(\BR, Y)} \leq C(\kappa_0+\kappa_1)\|f\|_{L_p(\BR, X)}
\quad\text{for all $f \in L_p(\BR, X)$}$$
with some positive constant $C$ depending on $p$.
\end{thm}

We now construct a solution of Eq. \eqref{stokes:1}.
Let $\bff_0$ be the zero extension of $\bff$ outside of $(0, T)$, that is
$\bff_0(t) = \bff(t)$ for $t \in (0, T)$ and $\bff_0(t) = 0$ for 
$t \not\in (0, T)$.  
Notice that $\bff_0$, $g$, $\bg$,  and $\bh$ 
 are defined on the whole line $\BR$.  Thus, 
we first consider the equations:
\begin{equation}\label{max:1}\left\{\begin{aligned}
\pd_t\bu_1 - \DV(\mu\bD(\bu_1) - \fq_1\bI) &= \bff_0 
&\quad&\text{in $\Omega\times\BR$}, \\
\dv\bu_1 &= g =\dv \bg
&\quad&\text{in $\Omega\times\BR$}, \\
(\mu\bD(\bu_1) - \fq_1\bI)\bn 
& = \bh
&\quad&\text{in $\Gamma\times\BR$}. 
\end{aligned}\right.\end{equation}
Let $\CF_L$ be the Laplace transform with respect to a time variable $t$
defined by
$$\hat f(\lambda) = \CF_L[f](\lambda)
= \int_{\BR}e^{-\lambda t}f(t)\,dt$$
for $\lambda = \gamma + i\tau \in \BC$.   Obviously, 
$$\CF_L[f](\lambda ) = \int_{\BR}e^{-i\tau t}e^{-\gamma t}f(t)\,dt
= \CF[e^{-\gamma t}f](\tau). 
$$
Applying the Laplace transform to Eq. \eqref{max:1} gives
\begin{equation}\label{max:2}\left\{\begin{aligned}
\lambda\hat\bu_1 - \DV(\mu\bD(\hat \bu_1) - \hat \fq_1\bI) &= \hat \bff_0 
&\quad&\text{in $\Omega$}, \\
\dv\hat\bu_1 &= \hat g = \dv \hat\bg
&\quad&\text{in $\Omega$}, \\
(\mu\bD(\hat \bu_1) - \hat \fq_1\bI)\bn 
& = \hat \bh
&\quad&\text{in $\Gamma$}. 
\end{aligned}\right.\end{equation}
Applying Theorem \ref{thm:rbdd:1.1}, we have 
$\hat \bu_1 = \CA_2(\lambda)\bF'_\lambda$ and 
$\hat \fq_1 = \CP_2(\lambda)\bF'_\lambda$ 
for $\lambda \in \Sigma_{\epsilon_0, \lambda_0}$,
where 
$$\bF'_\lambda = (\hat \bff_0(\lambda), 
\lambda^{1/2}\hat \bh(\lambda), \hat \bh(\lambda)
\lambda^{1/2}\hat g(\lambda),  \hat g(\lambda), 
\lambda \hat \bg(\lambda)).$$ 
Let $\CF_L^{-1}$ be the inverse Laplace transform defined by
$$\CF_L^{-1}[g](t) = \frac{1}{2\pi}\int_\BR e^{\lambda t}g(\tau)\,d\tau 
= e^{\gamma t}\frac{1}{2\pi}\int_\BR e^{i\tau t}g(\tau)\,d\tau
$$
for $\lambda = \gamma + i\tau \in \BC$. Obviously, 
$$\CF_L^{-1}[g](t) = e^{\gamma t}\CF^{-1}[g](t),
\quad \CF_L\CF_L^{-1} = \CF_L^{-1}\CF_L = \bI.$$
Setting 
$$\Lambda_\gamma^{1/2}f = \CF_L^{-1}[\lambda^{1/2}\CF_L[f]]
= e^{\gamma t}\CF^{-1}[\lambda^{1/2}\CF[e^{-\gamma t}f]],
$$
and
using the facts that
$$\lambda \hat \bg(\lambda) = \CF_L[\pd_t\bg](\lambda),
\quad \lambda^{1/2} \hat f(\lambda) = \CF_L[\Lambda^{1/2}_\gamma f] 
= \CF[e^{-\gamma t}\Lambda_\gamma^{1/2}f]$$
for $f \in \{g, \bh\}$ 
we define $\bu_1$ and $\fq_1$  by 
\begin{align*}
\bu_1(\cdot, t) & = \CF_L[\CA_2(\lambda)\bF'_\lambda]
= e^{\gamma t}\CF^{-1}[\CA_2(\lambda)\CF[e^{-\gamma t}F'(t)](\tau)], \\
\fq_1(\cdot, t) & = \CF_L[\CP_2(\lambda)\bF'_\lambda]
= e^{\gamma t}\CF^{-1}[\CP_2(\lambda)\CF[e^{-\gamma t}F'(t)](\tau)], 
\end{align*}
with $F'(t) = (\bff_0, \Lambda_\gamma^{1/2}\bh, \bh,
\Lambda_\gamma^{1/2}g, g, \pd_t\bg)$, 
where $\gamma$ is chosen as $\gamma > \lambda_0$,
and so $\gamma + i\tau \in \Sigma_{\epsilon_0, \lambda_0}$ for any
$\tau \in \BR$.  By Cauchy's theorem in the theory of one complex
variable, $\bu_1$ and $\fq_1$  are independent of 
choice of $\gamma$ whenever $\gamma > \lambda_0$
and the condition \eqref{increase:1} is satisfied 
for $\gamma > \lambda_0$.   Noting that
\begin{align*}
\pd_t\bu_1 &= \CF_L^{-1}[\lambda\CA_1(\lambda)\bF'_\lambda]
=e^{\gamma t}\CF^{-1}[\lambda \CA_1(\lambda)\CF[e^{-\gamma t}F'(t)](\tau)], 
\end{align*}
and applying Theorem \ref{thm:weis}, we have 
\begin{equation}\label{max:3}\begin{split}
&\|e^{-\gamma t}\pd_t\bu_1\|_{L_p(\BR, L_q(\Omega))}
+ \|e^{-\gamma t}\bu_1\|_{L_p(\BR, H^2_q(\Omega))}
+ \|e^{-\gamma t}\nabla\fq_1\|_{L_p((\BR, L_q(\Omega))}\\
& \leq Cr_b\|e^{-\gamma t}F'\|_{L_p(\BR, \CH_q(\Omega))}\\
&\leq Cr_b\{\|e^{-\gamma t}\bff\|_{_p((0, T), L_q(\Omega))}
+ 
\|e^{-\gamma t}g\|_{_p(\BR, H^1_q(\Omega))}
+\|e^{-\gamma t}\Lambda_\gamma^{1/2}g\|_{_p(\BR, L_q(\Omega))} \\
&
+\|e^{-\gamma t}\pd_t\bg\|_{_p(\BR, L_q(\Omega))}
+\|e^{-\gamma t}\bh\|_{_p(\BR, H^1_q(\Omega))}
+ 
\|e^{-\gamma t}\Lambda_\gamma^{1/2}\bh\|_{_p(\BR, L_q(\Omega))}\}.
\end{split}\end{equation}

We now write solutions $\bu$ and $\fq$  of Eq. \eqref{stokes:1} by
$\bu= \bu_1 + \bu_2$ and $\fq = \fq_1 + \fq_2$, 
where $\bu_2$ and $q_2$  are solutions of the 
following equations: 
\begin{equation}\label{max:4}\left\{\begin{aligned}
\pd_t\bu_2 - \DV(\mu\bD(\bu_2) \fq_2\bI) = 0,
\quad \dv \bu_2=0&
&\quad&\text{in $\Omega\times(0, \infty)$}, \\
(\mu\bD(\bu_2) - \fq_2\bI)\bn = 0&
&\quad&\text{on $\Gamma\times(0, \infty)$}, \\
\bu_2 = \bu_0-\bu_1|_{t=0}&
&\quad&\text{in $\Omega$}.
\end{aligned}\right.\end{equation}
Notice that $\dv\bu_2 = 0$ in $\Omega\times(0, \infty)$ means that
$\bu_2 \in J_q(\Omega)$ for any $t > 0$.  By real interpolation theory,
we know that
\begin{equation}\label{max:5}\begin{split}
&\sup_{t\in (0, \infty)} e^{-\gamma t}\|\bu_1(t)\|_{B^{2(1-1/p)}_{q,p}(\Omega)}
\\
&\quad \leq C(\|e^{-\gamma t}\bu_1\|_{L_p((0, \infty), H^2_q(\Omega))}
+ \|e^{-\gamma t}\pd_t\bu_1\|_{L_p((0, \infty), L_q(\Omega))}).
\end{split}\end{equation}
In fact, this inequality follows from the following theory
 (cf. Tanabe[p.1]\cite{tanabe}):  Let $X_1$ and $X_2$
be two Banach spaces such that $X_2$ is a dense subset of
$X_1$, and then 
\begin{equation}\label{tanabe:1}
L_p((0, \infty), X_2) \cap H^1_p((0, \infty), X_1)
\subset C([0, \infty), (X_1, X_2)_{1-1/p, p}),
\end{equation}
and
\begin{equation}\label{tanabe:2}
\sup_{t\in (0, \infty)} \|u(t)\|_{(X_1, X_2)_{1-1/p,p}}
\leq C(\|u\|_{L_p((0, \infty), X_2)}
+ \|\pd_tu\|_{L_p((0, \infty), X_1)}).
\end{equation}
Since $B^{2(1-1/p)}_{q,p}(\Omega) = (L_p(\Omega), H^2_q(\Omega))_{1-1/p,p}$, 
we have \eqref{max:5}.  
Thus, 
$$\bu_0 - \bu_1|_{t=0}
\in B^{2(1-1/p)}_{q,p}(\Omega).
$$
By the compatibility condition \eqref{compati:1} and 
\eqref{max:1}, we have
$$(\bu_0-\bu_1|_{t=0}, \nabla\varphi)
= (\bu_0 - \bg|_{t=0}, \nabla\varphi) = 0
\quad \text{for any $\varphi \in \hat H^1_{q',0}(\Omega)$}.
$$
Moreover, if $2/p + 1/q < 1$, then by the compatibility condition
\eqref{compati:2} and \eqref{max:1}, we have
$$(\mu\bD(\bu_0-\bu_1|_{t=0})\bn)_\tau 
= (\mu\bD(\bu_0)\bn)_\tau - (\bh|_{t=0})_\tau = 0 \quad
\text{on $\Gamma$}.
$$
Thus, if $2/p + /q \not=1$, then $\bu_0- \bu_1|_{t=0}
 \in \CD^1_{q,p}(\Omega)$. Applying Theorem
\ref{thm:max:2.1}, we see that there exists a $\gamma'> 0$ such 
that Eq. \eqref{max:4} admits unique solutions $\bu_2$  with
$\fq_2 = K_0(\bu_2)$ and 
\begin{equation}\label{max:6}\begin{split}
e^{-\gamma' t}\bu_2 &\in H^1_p((0, \infty), L_q(\Omega)^N) \cap 
L_p((0, \infty), H^2_q(\Omega)^N)
\end{split}\end{equation}
possessing the estimate:
\begin{equation}\label{max:7}\begin{split}
&\|e^{-\gamma' t}\pd_t\bu_2\|_{L_p((0, \infty), L_q(\Omega))}
+ \|e^{-\gamma't}\bu_2\|_{L_p((0, \infty), H^2_q(\Omega))} \\
&\quad \leq C\|\bu_0-\bu_1|_{t=0}\|_{B^{2(1-1/p)}_{q,p}(\Omega)}.
\end{split}\end{equation}
Thus, setting $\bu= \bu_1 + \bu_2$ and $\fq = \fq_1 + K_0(\bu_2)$ 
 and choosing $\gamma_0$ in such a way that
$\gamma_0 > \max(\lambda_0, \gamma')$, by \eqref{max:1}, \eqref{max:3},
\eqref{max:4}, \eqref{max:5}, \eqref{max:6} and \eqref{max:7}, we see
that $\bu$ and $\fq$  are required solutions of Eq. \eqref{stokes:1}.
Employing the same argument as in the proof of the uniqueness of
Theorem \ref{thm:max:2.2}, we can show the uniqueness. 
This completes the proof of  Theorem \ref{thm:max.1}. \qed
\vskip0.5pc
Employing the same argument as above, we can show Theorem \ref{thm:max.2},
and so we may omit the proof of Theorem \ref{thm:max.2}. Thus, we finally
give a \vskip0.5pc
{\bf Proof of Theorem \ref{thm:max.3}}\quad Let $d_0$ be the zero extension
of $d$ outside of $(0, T)$, that is $d_0(t) = d(t)$ for $t \in (0, T)$ and 
$d_0(t) = 0$.  Employing the same argument as in the proof of
Theorem \ref{thm:max.3} above, we can show the existence of solutions, 
$\bw$, $\fr$ and $\rho$,  of the equations:
\begin{equation}\label{stokes:3*}\left\{\begin{aligned}
\pd_t\bw - \DV(\mu\bD(\bw) - \fr\bI) = 0,
\quad \dv\bu &= 0
&\quad&\text{in $\Omega\times\BR$}, \\
\pd_t \rho + A_\kappa\cdot\nabla'_\Gamma \rho 
-\bw\cdot\bn + \CF_1\bw& = d_0
&\quad&\text{on $\Gamma\times\BR$}, \\
(\mu\bD(\bw) - \fr\bI)\bn - (\CF_2\rho + \sigma\Delta_\Gamma\rho)\bn
& = 0 &\quad&\text{on $\Gamma\times\BR$},
\end{aligned}\right.\end{equation}
possessing the estimate:
\begin{equation}\label{proof:max:1}\begin{split}
&\|e^{-\gamma t}\bw\|_{L_p(\BR, H^2_q(\Omega))}
+ \|e^{-\gamma t}\pd_t\bw\|_{L_p(\BR, L_q(\Omega))} \\
&\quad+ \|e^{-\gamma t}\rho\|_{L_p(\BR, W^{3-1/q}_q(\Gamma))}
+ \|e^{-\gamma t}\pd_t\rho\|_{L_p(\BR, W^{2-1/q}_q(\Gamma))}\\
&\qquad\leq C\|e^{-\gamma t}d_0\|_{L_p(\BR, W^{2-1/q}_q(\Gamma))}
\leq C\|d\|_{L_q((0, T), W^{2-1/q}_q(\Gamma)}
\end{split}\end{equation}
for any $\gamma \geq \lambda_0\kappa^{-b}$,
where we have used \eqref{tanabe:1} and \eqref{tanabe:2}. 
In particular, for any $\epsilon > 0$ and 
$\gamma \geq \lambda_0\kappa^{-b}$,  we have
\begin{align*}
&\|e^{-\gamma t}\bw\|_{L_p((-\infty, -\epsilon), H^2_q(\Omega))}
+ \|e^{-\gamma t}\rho\|_{L_p((-\infty, -\epsilon), 
W^{3-1/q}_q(\Gamma))} \\
&\leq 
\|e^{-\gamma t}\bw\|_{L_p(\BR, H^2_q(\Omega))}
+ \|e^{-\gamma t}\rho\|_{L_p(\BR, W^{3-1/q}_q(\Gamma))}
 \leq C\|d\|_{L_p((0, T), W^{2-1/q}_q(\Gamma))}.
\end{align*}
By the monotonicity of $e^{-\gamma t}$, we have 
\begin{align*}
&\|e^{-\gamma t}\bw\|_{L_p((-\infty, -\epsilon), H^2_q(\Omega))}
+ \|e^{-\gamma t}\rho\|_{L_p((-\infty, -\epsilon), 
W^{3-1/q}_q(\Gamma))} \\
&\quad
\geq e^{\gamma \epsilon}
(\|\bw\|_{L_p((-\infty, -\epsilon), H^2_q(\Omega))}
+ \|\rho\|_{L_p((-\infty, -\epsilon), 
W^{3-1/q}_q(\Gamma))})
\end{align*}
Putting these inequalities together gives
\begin{align*}
&\|\bw\|_{L_p((-\infty, -\epsilon), H^2_q(\Omega))}
+ \|\rho\|_{L_p((-\infty, -\epsilon), 
W^{3-1/q}_q(\Gamma))} \\
&\quad
\leq Ce^{-\gamma \epsilon}\|d\|_{L_p((0, T), W^{2-1/q}_q(\Gamma))}.
\end{align*}
for any $\gamma \geq \lambda_0\kappa^{-b}$.
Thus, letting $\gamma \to \infty$, we have
$$\|\bw\|_{L_p((-\infty, -\epsilon), H^2_q(\Omega))}
+ \|\rho\|_{L_p((-\infty, -\epsilon), 
W^{3-1/q}_q(\Gamma))} =0.$$
Since $\epsilon>0$ is chosen arbitrarily, we have
$$\|\bw\|_{L_p((-\infty, 0), H^2_q(\Omega))}
+ \|\rho\|_{L_p((-\infty, 0), 
W^{3-1/q}_q(\Gamma))} =0,$$
which shows that $\bw = 0$ and $\rho=0$ for $t < 0$,
because 
$$\bw \in C(\BR, B^{2(1-1/p)}_{q,p}(\Omega)), \quad 
\rho \in C(\BR, W^{3-1/p-1/q}_{q,p}(\Gamma)).
$$ 

By the monotonicity of $e^{-\gamma t}$, we have 
\begin{align*}
&\|e^{-\gamma t}\bw\|_{L_p(\BR, H^2_q(\Omega))}
+ \|e^{-\gamma t}\pd_t\bw\|_{L_p(\BR, L_q(\Omega))} \\
&\geq 
\|e^{-\gamma t}\bw\|_{L_p((0, T), H^2_q(\Omega))}
+ \|e^{-\gamma t}\pd_t\bw\|_{L_p((0, T), L_q(\Omega))}\\
&\geq e^{-\gamma T}(
\|\bw\|_{L_p((0, T), H^2_q(\Omega))}
+ \|\pd_t\bw\|_{L_p((0, T), L_q(\Omega))}).
\end{align*}
Similarly, we have
\begin{align*}
&\|e^{-\gamma t}\rho\|_{L_p(\BR, W^{3-1/q}_q(\Gamma))}
+ \|e^{-\gamma t}\pd_t\rho\|_{L_p(\BR, W^{2-1/q}_q(\Gamma))} \\
&\quad \geq e^{-\gamma T}(\|\rho\|_{L_p((0, T), W^{3-1/q}_q(\Gamma))}
+ \|\pd_t\rho\|_{L_p((0, T), W^{2-1/q}_q(\Gamma))}).
\end{align*}
Thus, by \eqref{proof:max:1}, we have
\begin{align*}
&\|\bw\|_{L_p((0, T), H^2_q(\Omega))}
+ \|\pd_t\bw\|_{L_p((0, T), L_q(\Omega))} \\
&\quad + 
\|\rho\|_{L_p((0, T), W^{3-1/q}_q(\Gamma))}
+ \|\pd_t\rho\|_{L_p((0, T), W^{2-1/q}_q(\Gamma))}\\
&\qquad
\leq Ce^{\gamma T}
\|d\|_{L_p((0, T), W^{2-1/q}_q(\Gamma))}
\end{align*}
for any $\gamma \geq \lambda_0\kappa^{-b}$. 
This completes the proof of the existence part of 
Theorem \ref{thm:max.3}. 
Employing the same argument as in the proof of the uniqueness of
Theorem \ref{thm:max:2.2}, we can show the uniqueness. 
This completes the proof of  Theorem \ref{thm:max.3}. \qed 

\section{$\CR$ bounded solution operators}\label{sec:5.0}
In this section, we mainly prove Theorem \ref{main:thm3}.  
The operators  $\CF_1$ and $\CF_2$ can be treated by 
perturbation method, and so the subsections 
\ref{subsec:5.1} -- \ref{newsubsec:6.4} below devote to proving 
the existence part for Eq. \eqref{rres:1.2} for the following
equations:
\begin{equation}\label{eq:7.0} 
\left\{\begin{aligned}
\lambda \bu - \DV(\mu\bD(\bu, h)-K(\bu, h)\bI) & = \bff
&\quad&\text{in $\Omega$}, \\
\lambda h + A_\kappa\cdot\nabla'_\Gamma h
- \bn\cdot\bu& = d &\quad&\text{on $\Gamma$}, \\
(\mu\bD(\bu) - K(\bu, h)\bI)\bn - (\sigma\Delta_\Gamma h)\bn
&=\bh&\quad&\text{on $\Gamma$}.
\end{aligned}\right.\end{equation}
And then, in Subsec. \ref{subsec:new5}, we prove the 
exsitence part of Theorem \ref{main:thm3} for
Eq. \eqref{rres:1.2} by using a perturbation method.

The existence part of Theorem \ref{main:thm0}
can be proved in the 
same manner as in the proof of the existence part of
Theorem \ref{main:thm3} and also 
has been proved by Shibata \cite{S1}, and so we omit its proof.

Concerning the uniqueness part, 
we first prove Theorem \ref{main:thm0} in Subsec. \ref{subsec:5.9}.
Finally, the uniqueness part of Theorem \ref{main:thm3} will be 
proved in Subsec. \ref{sec:5.10} by
showing {\it a apriori} estimates for Eq. \eqref{rres:2} 
under the assumption that $\Omega$ is a uniform $C^3$ domain
whose inside has a finite covering.

\subsection{Model Problem in $\BR^N$; 
Constant $\mu$ case}\label{subsec:5.1}
In this subsection, 
we assume that $\mu$ is a constant satisfying the assumption
\eqref{assump:3}, that is  $m_0 \leq \mu \leq m_1$. 
Given $\bu \in H^2_q(\BR^N)^N$, let 
$u = K_0(\bu)$ be a unique solution of the weak Laplace problem:
\begin{equation}\label{3.1.1}
(\nabla u, \nabla\varphi)_{\BR^N} = (\DV(\mu\bD(\bu)) - \nabla\dv\bu, 
\nabla\varphi)_{\BR^N}.
\end{equation}
for any $\varphi \in \hat H^1_{q'}(\BR^N)$.
In this subsection, we consider the resolvent problem:
\begin{equation}\label{3.1.2} 
\lambda\bu - \DV(\mu\bD(\bu) - K_0(\bu)\bI) = \bff\quad\text{in $\BR^N$},
\end{equation}
and prove the following theorem.
\begin{thm}\label{thm:4.1} Let $1 < q < \infty$, $0 < \epsilon < \pi/2$,
and $\lambda_0 > 0$.  Then, there exists an operator family
$\CA_0(\lambda) \in {\rm Hol}\,(\Sigma_{\epsilon, \lambda_0}, 
\CL(L_q(\BR^N)^N, H^2_q(\BR^N)^N))$ such that for any 
$\lambda=\gamma + i\tau \in \Sigma_{\epsilon, \lambda_0}$ and 
$\bff \in L_q(\BR^N)^N$, $\bu = \CA_0(\lambda)\bff$ is a unique solution of
Eq. \eqref{3.1.2} and 
\begin{equation}\label{3.1.7}
\CR_{\CL(L_q(\BR^N)^N, H^{2-j}_q(\BR^N)^N)}
(\{(\tau\pd_\tau)^\ell(\lambda^{j/2}\CA_0(\lambda)) \mid
\lambda \in \Sigma_{\epsilon, \lambda_0}\}) \leq r_b(\lambda_0)
\end{equation}for $\ell= 0,1$ and $j=0,1,2$, where $r_b(\lambda_0)$ is a 
constant depending on $\epsilon$, $\lambda_0$, $m_0$, $m_1$, 
$q$ and $N$, but independent of $\mu \in [m_0, m_1]$. 
\end{thm}
\pf We first consider the Stokes equations:
\begin{equation}\label{3.1.3} \lambda\bu - \DV(\mu(\bD(\bu) - \fq\bI)
= \bff, \quad \dv\bu = g = \dv\bg
\quad\text{in $\BR^N$}.
\end{equation}
Since $\DV(\mu\bD(\bu) - \fq\bI) = \mu\Delta\bu 
+ \mu\nabla\dv\bu  -\nabla\fq$, applying $\dv$ to \eqref{3.1.3},
we have
$$\lambda\dv \bg - 2\mu\Delta g + \Delta\fq = \dv\bff,$$
and so, 
$$\fq = 2\mu g + \Delta^{-1}(\dv\bff - \lambda\dv\bg).$$
Combining this with \eqref{3.1.3} gives 
\begin{equation}\label{3.1.4} \lambda\bu - \mu\Delta\bu =
\bff - \nabla\Delta^{-1}\dv\bff - \mu\nabla g + \lambda\nabla\Delta^{-1}\dv\bg.
\end{equation}
We now look for a solution formula for Eq. \eqref{3.1.2}.  Let $g$ 
be a solution of the variational problem:
$$(\lambda g, \varphi)_{\BR^N} + (\nabla g, \nabla\varphi)_{\BR^N}
= (-\bff, \nabla\varphi)_{\BR^N} \quad\text{for any $\varphi 
\in \hat H^1_{q'}(\BR^N)$},
$$
and then this $g$ is given by $g = (\lambda-\Delta)^{-1}\dv\bff$. 
According to \eqref{wd:5*}, we set $\bg = \lambda^{-1}(\bff + \nabla g)$.  
Inserting these
formulas into \eqref{3.1.4} gives 
$$\lambda\bu - \mu\Delta\bu = \bff-(\mu-1)\nabla g
= \bff-(\mu-1)(\lambda-\Delta)^{-1}\nabla\dv\bff.
$$
Thus, we have
$$\bu = \CF^{-1}_\xi\Bigl[\frac{\CF[\bff](\xi)}{\lambda + \mu|\xi|^2}\Bigr]
+(\mu-1)\CF^{-1}_\xi\Bigl[\frac{\xi\xi\cdot\CF[\bff](\xi)}{(\lambda+\mu|\xi|^2)
(\lambda+|\xi|^2)}\Bigr],
$$
where $\CF$ and $\CF_\xi^{-1}$ denote the Fourier transform 
and its inversion formula defined by 
$$\CF[f](\xi) = \int_{\BR^N}e^{-ix\cdot\xi}f(x)\,dx,
\quad
\CF_\xi^{-1}[g(\xi)](x) = \frac{1}{(2\pi)^N}
\int_{\BR^N}e^{ix\cdot\xi}g(\xi)\,d\xi.
$$
Thus, we define an operator family $\CA_0(\lambda)$ acting
on $\bff \in L_q(\BR^N)^N$ by 
$$\CA_0(\lambda) \bff = \CF^{-1}_\xi\Bigl[
\frac{\CF[\bff](\xi)}{\lambda + \mu|\xi|^2}\Bigr]
+(\mu-1)\CF^{-1}_\xi\Bigl[\frac{\xi\xi\cdot\CF[\bff](\xi)}{(\lambda+\mu|\xi|^2)
(\lambda+|\xi|^2)}\Bigr].
$$
To prove the $\CR$-boundedness of $\CA_0(\lambda)$, we use the 
following lemma.
\begin{lem}\label{lem:fund1} Let $0 < \epsilon < \pi/2$.  Then, 
for any $\lambda \in \Sigma_\epsilon$ and $x \in [0, \infty)$, we have
\begin{equation}\label{3.1.5} |\lambda + x| \geq(\sin\frac{\epsilon}{2})
(|\lambda|+ x).
\end{equation}\end{lem}
\pf Representing $\lambda = |\lambda| e^{i\theta}$ and using 
$\cos\theta \geq \cos(\pi-\epsilon) = -\cos\epsilon$ for
$\lambda \in \Sigma_\epsilon$, we have
\eqref{3.1.5}. \qed
\begin{lem}\label{lem:fund2} Let $1 < q < \infty$ and let $U$ be a subset of 
$\BC$.  Let $m= m(\lambda, \xi)$ be a function defined on 
$U\times(\BR^N\setminus\{0\})$ which is infinitely 
differentiable with respect to $\xi \in \BR^N\setminus\{0\}$ for 
each $\lambda \in U$.  Assume that for any multi-index $\alpha 
\in \BN_0^N$ there exists a constant $C_\alpha$ depending on $\alpha$ such that
\begin{equation}\label{3.1.6}
|\pd_\xi^\alpha m(\lambda, \xi)| \leq C_\alpha|\xi|^{-|\alpha|}
\end{equation}
for any $(\lambda, \xi) \in U\times(\BR^N\setminus\{0\})$.
Set 
$$\bb(m) = \max_{|\alpha| \leq N+1} C_\alpha.$$
Let $K_\lambda$ be an operator defined by 
$$K_\lambda f = \CF_\xi^{-1}[m(\lambda, \xi)\CF[f](\xi)].
$$
Then, the operator family $\{K_\lambda\mid \lambda \in U\}$
is $\CR$-bounded on $\CL(L_q(\BR^N))$ and 
$$\CR_{\CL(L_q(\BR^N))}(\{K_\lambda \mid \lambda \in U\})
\leq C_{N,q}\bb(m)$$
for some constant $C_{q,N}$ depending solely on $q$ and $N$.
\end{lem}
\pf 
Lemma \ref{lem:fund2} was proved by Enomoto and Shibata \cite[Theorem 3.3]{ES1}
and Denk and Schnaubelt \cite[Lemma 2.1]{DS}. \qed
\vskip0.5pc
By Lemma \ref{lem:fund1}, we have 
\begin{align*}
\Bigl|\pd_\xi^\alpha\frac{\lambda^{j/2}\xi^\beta}{\lambda + \mu|\xi|^2}
\Bigr|& \leq C_\alpha|\xi|^{-|\alpha|}\lambda_0^{k/2}, \\
\Bigl|\pd_\xi^\alpha\frac{\xi_\ell\xi_m\lambda^{j/2}\xi^\beta}
{(\lambda + \mu|\xi|^2)(\lambda + |\xi|^2)}\Bigr|
&\leq C_\alpha|\xi|^{-|\alpha|}\lambda_0^{-k/2}
\quad(\ell,  m=1, \ldots, N)
\end{align*}
for any $j \in \BN_0$, $k \in \BN_0$ and $\beta \in \BN_0^N$ such that
$j+k+|\beta|=2$ and for any $\alpha \in \BN_0^N$ and 
$(\lambda, \xi) \in \Sigma_{\epsilon, \lambda_0}\times
(\BR^N\setminus\{0\})$. Thus, by Lemma \ref{lem:fund2}, we have
\eqref{3.1.7}, which completes the proof of Theorem \ref{thm:4.1}.
\qed
\vskip0.5pc
We conclude  this subsection by introducing  some fundamental properties of
$\CR$-bounded operators and Bourgain's results concerning
 Fourier multiplier theorems with scalar multiplieres.
\begin{prop}\label{prop:4.1}
\thetag{a} Let $X$ and $Y$ be Banach spaces,
and let $\CT$ and $\CS$ be
$\CR$-bounded families in $\CL(X,Y)$.
Then, $\CT+ \CS = \{T + S \mid
T \in \CT, \enskip S \in \CS\}$ is also an $\CR$-bounded
family in
$\CL(X,Y)$ and
$$\CR_{\CL(X,Y)}(\CT + \CS) \leq
\CR_{\CL(X,Y)}(\CT) + \CR_{\CL(X,Y)}(\CS).$$

\thetag{b} Let $X$, $Y$ and $Z$ be Banach spaces,  and let
$\CT$ and $\CS$ be $\CR$-bounded families in $\CL(X, Y)$ and
$\CL(Y, Z)$, respectively.  Then, $\CS\CT = \{ST \mid
T \in \CT, \enskip S \in \CS\}$ also an $\CR$-bounded
family in $\CL(X, Z)$ and
$$\CR_{\CL(X, Z)}(\CS\CT) \leq \CR_{\CL(X,Y)}(\CT)\CR_{\CL(Y, Z)}(\CS).$$

\thetag{c} Let $1 < p, \, q < \infty$ and let $D$ be a domain in $\BR^N$.
Let $m=m(\lambda)$ be a bounded function defined on a subset
$U$ of $\BC$ and let $M_m(\lambda)$ be a map defined by
$M_m(\lambda)f = m(\lambda)f$ for any $f \in L_q(D)$.  Then,
$\CR_{\CL(L_q(D))}(\{M_m(\lambda) \mid \lambda \in U\})
\leq C_{N,q,D}\|m\|_{L_\infty(U)}$.

\thetag{d} Let $n=n(\tau)$ be a $C^1$-function defined on $\BR\setminus\{0\}$
that satisfies the conditions $|n(\tau)| \leq \gamma$
and $|\tau n'(\tau)| \leq \gamma$ with some constant $c > 0$ for any
$\tau \in \BR\setminus\{0\}$.  Let $T_n$ be an
 operator-valued Fourier multiplier defined by $T_n f = \CF^{-1}(n \CF[f])$
for any $f$ with $\CF[f] \in \CD(\BR, L_q(D))$.  Then,
$T_n$ is extended to a bounded linear operator
from $L_p(\BR, L_q(D))$ into itself.  Moreover,
denoting this extension also by $T_n$, we have
$$\|T_n\|_{\CL(L_p(\BR, L_q(D)))} \leq C_{p,q,D}\gamma.
$$
\end{prop}
\pf 
The assertions a) and b) follow from
\cite[p.28, Proposition 3.4]{DHP}, and
the assertions c) and d) follow from
 \cite[p.27, Remarks 3.2]{DHP}
(see also Bourgain \cite{Bourgain}).
\qed

\subsection{Perturbed problem in $\BR^N$}\label{subsec:3.2}

In this subsection, we consider the case where  $\mu(x)$ 
is a real valued function satisfying \eqref{assump:3}. 
Let $x_0$ be any point in $\Omega$ and 
let $d_0$ be a positive number such that 
$B_{d_0}(x_0) \subset \Omega$. In view of \eqref{assump:3}, we assume 
that 
\begin{equation}\label{p.1.1} |\mu(x)-\mu(x_0)| \leq m_1M_1
\quad\text{for $x \in B_{d_0}(x_0)$},
\end{equation}
where we have set $M_1 = d_0$.  We assume that $M_1 \in (0, 1)$ below.
Let $\varphi$ be a function in $C^\infty_0(\BR^N)$ which equals $1$ for 
$x \in B_{d_0/2}(x_0)$ and $0$ outside of $B_{d_0}(x_0)$.  Let
\begin{equation}\label{p.1.2} \tilde\mu(x) = \varphi(x)\mu(x)
+(1-\varphi(x))\mu(x_0).
\end{equation}
Let $\tilde K_0(\bu) \in \hat H^1_q(\BR^N)$ be a unique solution of
the weak Laplace problem:
\begin{equation}\label{p.1.3}
(\nabla u, \nabla\varphi)_{\BR^N} = (\DV(\tilde\mu \bD(\bu))-\nabla\dv\bu, 
\nabla\varphi)_{\BR^N}
\quad\text{for any $\varphi \in \hat H^1_{q'}(\BR^N)$}.
\end{equation}
We consider the resolvent problem:
\begin{equation}\label{p.1.4}
\lambda\bu - \DV(\tilde\mu\bD(\bu) - \tilde K_0(\bu)\bI) = \bff
\quad\text{in $\BR^N$}.
\end{equation}
We shall prove the following theorem.
\begin{thm}\label{thm:p.1} Let $1 < q <\infty$ and $0 < \epsilon < \pi/2$.
Then, there exist $M_1 \in (0, 1)$, $\lambda_0 \geq 1$ and 
an operator family $\tilde \CA_0(\lambda)$ with
$$\tilde \CA_0(\lambda) \in {\rm Hol}\,(\Sigma_{\epsilon, \lambda_0}, 
\CL(L_q(\BR^N)^N, H^2_q(\BR^N)^N))
$$
such that for any $\lambda \in \Sigma_{\epsilon, \lambda_0}$ and $\bff
\in L_q(\BR^N)^N$, $\bu = \tilde \CA(\lambda)\bff$ is  a unique solution
of Eq. \eqref{p.1.4}, and 
$$\CR_{\CL(L_q(\BR^N)^N, H^{2-j}_q(\BR^N)^N)}
(\{(\tau\pd_\tau)^\ell(\lambda^{j/2}\tilde \CA_0(\lambda)) \mid
\lambda \in \Sigma_{\epsilon, \lambda_0}\}) \leq \tilde r_b
$$
for $\ell=0,1$ and $j=0,1,2$.  Where, $\tilde r_b$ 
is a constant independent of $M_1$ and $\lambda_0$.
\end{thm}
\pf
Let $u=K_{x_0}(\bu) \in \hat H^1_q(\BR^N)$ 
be a unique solution of the 
weak Laplace equation:
\begin{equation}\label{p.1.5} (\nabla u, \nabla\varphi)_{\BR^N}
=(\DV(\mu(x_0)\bD(\bu) - \nabla\dv\bu, \nabla\varphi)_{\BR^N}
\end{equation}
for any $\varphi \in \hat H^1_{q'}(\BR^N)$. 
We consider the resolvent problem:
\begin{equation}\label{p.1.6} \lambda\bu - \DV(\mu(x_0)\bD(\bu) - 
K_{x_0}(\bu)\bI) = \bff \quad\text{in $\BR^N$}.
\end{equation}
Let $\CB_{x_0}(\lambda) \in {\rm Hol}\,(\Sigma_{\epsilon, 1}, \CL(L_p(\BR^N)^N,
H^2_q(\BR^N)^N))$ be a solution operator of Eq. \eqref{p.1.6} such that 
for any $\lambda \in \Sigma_{\epsilon, 1}$ and $\bff \in L_q(\BR^N)^N$, $\bu = \CB_{x_0}(\lambda)\bff$ is a unique solution of Eq.\eqref{p.1.6} and 
\begin{equation}\label{p.1.7}
\CR_{\CL(L_q(\BR^N)^N, H^{2-j}_q(\BR^N)^N)}
(\{(\tau\pd_\tau)^\ell(\lambda^{j/2}\CB_{x_0}(\lambda)) \mid 
\lambda \in \Sigma_{\epsilon, 1}\}) \leq \gamma_0
\end{equation}
for $\ell=0,1$ and $j=0,1,2$, where $\gamma_0$ is a constant independent of 
$M_1$ and $\nabla\varphi$. Such an operator is given in Theorem \ref{thm:4.1}
with $\mu = \mu(x_0)$ and $\lambda_0 = 1$.  Inserting the formula:
$\bu = \CB_{x_0}(\lambda)\bff$ into \eqref{p.1.4} gives
\begin{equation}\label{p.1.8} \lambda\bu - \DV(\tilde\mu(x)\bD(\bu) - 
\tilde K_0(\bu)\bI) = \bff - \CR(\lambda)\bff
\quad\text{in $\BR^N$},
\end{equation}
where we have set
\begin{equation}\label{p.1.9}\begin{split}
\CR(\lambda)\bff &= \DV(\tilde\mu(x)\bD(\CB_{x_0}(\lambda)\bff) - \mu(x_0)
\bD(\CB_{x_0}(\lambda)\bff)) \\
&- \nabla(\tilde K_0(\CB_{x_0}(\lambda)\bff)
- K_{x_0}(\CB_{x_0}(\lambda)\bff)).
\end{split}\end{equation}
We shall estimate $\CR(\lambda)\bff$. 
For any $\varphi \in \hat H^1_{q'}(\BR^N)$, by \eqref{p.1.3} and
\eqref{p.1.5}, we have
\begin{align*}
&(\nabla(\tilde K_0(\CB_{x_0}(\lambda)\bff) - K_{x_0}(\CB_{x_0}(\lambda)\bff)),
\nabla\varphi)_{\BR^N} \\
&= 
((\DV((\tilde\mu(x)-\mu(x_0))\bD(\CB_{x_0}(\lambda)\bff)),
\nabla\varphi)_{\BR^N}.
\end{align*}
Since $\tilde\mu(x) - \mu(x_0) = \varphi(x)(\mu(x)-\mu(x_0))$, by \eqref{p.1.1}
and \eqref{assump:3}, we have
\begin{align*}
&\|\DV((\tilde\mu(x)-\mu(x_0))\bD(\CB_{x_0}(\lambda)\bff)\|_{L_q(\BR^N)}
\\
&\quad
\leq M_1\|\nabla^2\CB_{x_0}(\lambda)\bff\|_{L_q(\BR^N)}
+ C_{m_1, \nabla\varphi}\|\nabla \CB_{x_0}(\lambda)\bff\|_{L_q(\BR^N)}.
\end{align*}
Here and in the following, $C_{m_1, \nabla\varphi}$ denotes a generic constant
depending on $m_1$ and $\|\nabla\varphi\|_{L_\infty(\BR^N)}$.  Thus, 
we have
\begin{equation}\label{p.1.10}
\|\CR(\lambda)\bff\|_{L_q(\BR^N)} \leq 
CM_1\|\nabla^2\CB_{x_0}(\lambda)\bff\|_{L_q(\BR^N)}
+ C_{m_1, \nabla\varphi}\|\nabla
\CB_{x_0}(\lambda)\bff\|_{L_q(\BR^N)}.
\end{equation}
Here and in the following, $C$ denotes a generic constants independent
of $M_1$, $m_1$, and $\|\nabla\varphi\|_{L_\infty(\BR^N)}$.  Let 
$\lambda_0$ be any number $\geq 1$ and let $n \in \BN$,
$\{\lambda_k\}_{k=1}^n \subset (\Sigma_{\epsilon, \lambda_0})^n$, 
and $\{F_k\}_{k=1}^n \subset (L_q(\BR^N)^N)^n$. By \eqref{p.1.10}, 
\eqref{p.1.7} and Proposition \ref{prop:4.1}, we have
\allowdisplaybreaks{
\begin{align*}
&\int^1_0\|\sum_{k=1}^n r_k(u)\CR(\lambda_k)\bff_k\|_{L_q(\BR^N)}^q\,du \\
& \leq 2^{q-1}M_1^q\int^1_0\|\sum_{k=1}^nr_k(u)\nabla^2
\CB_{x_0}(\lambda_k)\bff_k\|_{L_q(\BR^N)}^q\,du \\
& + 2^{q-1}C_{m_1, \nabla\varphi}^q
\int^1_0\|\sum_{k=1}^nr_k(u)\nabla
\CB_{x_0}(\lambda_k)\bff_k\|_{L_q(\BR^N)}^q\,du \\
& \leq 2^{q-1}M_1^q\int^1_0\|\sum_{k=1}^nr_k(u)\nabla^2
\CB_{x_0}(\lambda_k)\bff_k\|_{L_q(\BR^N)}^q\,du \\
& + 2^{q-1}C_{m_1, \nabla\varphi}^q \lambda_0^{-q/2}
\int^1_0\|\sum_{k=1}^nr_k(u)\lambda_k^{1/2}\nabla
\CB_{x_0}(\lambda_k)\bff_k\|_{L_q(\BR^N)}^q\,du \\
& \leq 2^{q-1}(M_1^q + C_{m_1, \nabla^\varphi}^q\lambda_0^{-q/2})
\gamma_0^q\int^1_0\|\sum_{k=1}^n r_k(u)\bff_k\|_{L_q(\BR^N)}^q\,du.
\end{align*}
}
Choosing $M_1$ so small that $2^{q-1}M_1^q\gamma_0^q \leq (1/2)(1/q)^q$ and 
$\lambda_0 \geq 1$ so large that $2^{q-1}C_{m_1, \nabla\varphi}^q\gamma_0^q 
\lambda_0^{-q/2} \leq (1/2)(1/2)^q$, we have
$$\CR_{\CL(L_q(\BR^N)}(\{\CR(\lambda) \mid \lambda \in \Sigma_{\epsilon, 
\lambda_0}\}) \leq 1/2.
$$
Analogously, we have 
$$\CR_{\CL(L_q(\BR^N)}(\{\tau\pd_\tau\CR(\lambda) 
\mid \lambda \in \Sigma_{\epsilon, 
\lambda_0}\}) \leq 1/2.
$$
Thus, $(\bI - \CR(\lambda))^{-1} = \bI + \sum_{j=1}^\infty \CR(\lambda)^j$
exists and
\begin{equation}\label{p.1.11*}
\CR_{\CL(L_q(\BR^N)}(\{(\tau\pd_\tau)^\ell(\bI - \CR(\lambda))^{-1} 
\mid \lambda \in \Sigma_{\epsilon, \lambda_0}\}) \leq 4
\quad\text{for $\ell=0,1$}. 
\end{equation}
Setting $\tilde \CA_0(\lambda) = \CB_{x_0}(\lambda)(\bI - \CR(\lambda))^{-1}$,
by \eqref{p.1.7}, \eqref{p.1.11*} and Propsoition \ref{prop:4.1},
we see that $\tilde \CA_0(\lambda)$ is a solution operator
satisfying the required properties with $\tilde r_b = 4\gamma_0$. 

To prove the uniqueness of solutions of Eq. \eqref{p.1.4}, let 
$\bu \in H^2_q(\BR^N)^N$ be a solution of the homogeneous equatuion:
$$\lambda\bu - \DV(\tilde\mu\bD(\bu) - \tilde K_0(\bu)\bI) = 0
\quad\text{in $\BR^N$}.
$$
And then, $\bu$ satisfies the non-homogeneous equation:
\begin{equation}\label{p.1.11}
\lambda \bu - \DV(\mu(x_0)\bD(\bu) - K_{x_0}(\bu)\bI) = R\bu
\quad\text{in $\BR^N$},
\end{equation}
where we have set 
$$R\bu = -\DV((\tilde \mu(x) -\mu(x_0))\bD(\bu))
+ \nabla(\tilde K_0(\bu) - K_{x_0}(\bu)). 
$$
Analogously to the proof of \eqref{p.1.1}, we have
\begin{equation}\label{p.1.12}
\|R\bu\|_{L_q(\BR^N)} \leq CM_1\|\nabla^2\bu\|_{L_q(\BR^N)}
+ C_{m_1, \nabla\varphi}\|\nabla\bu\|_{L_q(\BR^N)}.
\end{equation}
On the other hand, applying Theorem \ref{thm:4.1} to \eqref{p.1.11}
for $\lambda \in \Sigma_{\epsilon, 1}$, we have
\begin{equation}\label{p.1.13}
|\lambda|\|\bu\|_{L_q(\BR^N)} + |\lambda|^{1/2}\|\bu\|_{H^1_q(\BR^N)}
+ \|\bu\|_{H^2_q(\BR^N)} \leq C\|R\bu\|_{L_q(\BR^N)}.
\end{equation}
Combining \eqref{p.1.12} and \eqref{p.1.13} gives 
$$(\lambda_0^{1/2}-CC_{m_1, \nabla\varphi})\|\bu\|_{H^1_q(\BR^N)}
+ (1-CM_1)\|\bu\|_{H^2_q(\BR^N)} \leq 0.
$$
Choosing $M_1 \in (0, 1)$ so small that $1 - CM_1 > 0$ and 
$\lambda_0 \geq 1$ so large that 
$\lambda_0^{1/2}-CC_{m_1, \nabla\varphi} > 0$, we have $\bu=0$. 
This proves the uniqueness, and therefore we have proved
Theorem \ref{thm:p.1}
\qed
\subsection{Model Problem in $\BR^N_+$} \label{subsec:4}
In this section, we assume that $\mu$, $\delta$, and $A_\kappa$ ($\kappa
\in [0, 1)$) are constants and 
an $N-1$ constant vector satisfying the conditions:
\begin{equation}\label{4.1}
m_0 \leq \mu, \sigma \leq m_1, \quad A_0 = 0,
\quad |A_\kappa| \leq m_2 \enskip (\kappa \in (0, 1)).
\end{equation}
Let 
\begin{align*}
\BR^N_+ &= \{(x_1, \ldots, x_N) \in \BR^N \mid x_N > 0\}, \\
\BR^N_0 &= \{(x_1, \ldots, x_N) \in \BR^N \mid x_N=0\},\\
\bn_0 &= {}^\top(0, \ldots, 0, -1).
\end{align*}
Given $\bu \in H^2_q(\BR^N_+)^N$ and $h \in W^{3-1/q}_q(\BR^N_0)$, let
$K(\bu, h) \in H^1_q(\BR^N_+) + \hat H^1_{q,0}(\BR^N_+)$ 
be a unique solution of the weak Dirichlet problem:
\begin{equation}\label{4.2}
(\nabla K(\bu, h), \nabla\varphi)_{\BR^N_+}
= (\DV(\mu\bD(\bu))-\nabla\dv\bu, \nabla\varphi)_{\BR^N_+}
\end{equation}
for any $\varphi \in \hat H^1_{q',0}(\BR^N)$,
subject to $K(\bu, h) = <\mu\bD(\bu)\bn_0, \bn_0> - \sigma\Delta' h -\dv\bu$
on $\BR^N_0$, where $\Delta'h = \sum_{j=1}^{N-1}\pd^2h/\pd x_j^2$.
In this section, we consider the half space problem:
\begin{equation}\label{4.3}\left\{\begin{aligned}
\lambda \bu - \DV(\mu\bD(\bu) - K(\bu, h)\bI) &=\bff
&\quad&\text{in $\BR^N_+$}, \\
\lambda h + A_\kappa\cdot\nabla'h - \bu\cdot\bn_0 &=d
&\quad&\text{on $\BR^N_0$}, \\
(\mu\bD(\bu)-K(\bu, h)\bI)\bn_0 - \sigma(\Delta'h)\bn_0
& = \bh &\quad&\text{on $\BR^N_0$}, 
\end{aligned}\right.\end{equation}
where $\nabla' = (\pd_1, \ldots, \pd_{N-1})$. The last equations 
in \eqref{4.3} are equivalent to 
$$(\mu\bD(\bu)\bn_0)_\tau = \bh_\tau\quad\text{and}\quad \dv \bu = 
\bh\cdot\bn_0
\quad \text{on $\BR^N_0$}.$$
Where, we have set  $\bh_\tau=\bh-<\bh, \bn_0>\bn_0$. . 
We shall show the following theorem
\begin{thm}\label{main:half1} Let $1 < q < \infty$, let $\mu$, $\sigma$,
and  $A_\kappa$ are constants and an $N-1$ constant vector satisfying 
the conditions in \eqref{4.1}. 
Let $\Lambda_{\kappa, \lambda_0}$  the set 
defined in Theorem \ref{thm:rbdd:1.2}. 
 Let $Y_q(\BR^N_+)$ and $\CY_q(\BR^N_+)$ be
spaces defined by replacing $\Omega$ and $\Gamma$ by $\BR^N_+$ and 
$\BR^N_0$ in Theorem \ref{main:thm3}. 
Then,  there exist a constant $\lambda_0 \geq 1$ 
and operator families: 
\begin{equation}\label{r-est:0}\begin{aligned}
 \CA_0(\lambda) &\in 
{\rm Hol}\, (\Lambda_{\kappa, \lambda_0},  
\CL(\CY_q(\BR^N_+), H^2_q(\BR^N_+)^N)), \\ 
\CH_0(\lambda) &\in {\rm Hol}\, (\Lambda_{\kappa, \lambda_0},
\CL(\CY_q(\BR^N_+), H^3_q(\BR^N_+)))
\end{aligned}\end{equation}
such that for any $\lambda=\gamma+i\tau \in \Lambda_{\kappa, \lambda_0}$ and 
$(\bff, d, \bh) \in Y_q(\BR^N_+)$,  
$$\bu = \CA_0(\lambda)(\bff, d, \lambda^{1/2}\bh, \bh), \quad
h = \CH_0(\lambda)(\bff, d, \lambda^{1/2}\bh, \bh),$$ 
are unique solutions of \eqref{4.3},
and 
\begin{equation}\label{r-est:1}\begin{split}
&\CR_{\CL(\CY_q(\BR^N_+), H^{2-j}_q(\BR^N_+)^N)}
(\{(\tau\pd_\tau)^\ell(\lambda^{j/2}\CA_0(\lambda)) \mid 
\lambda \in \Lambda_{\kappa, \lambda_0}\}) \leq 
r_b, \\
&\CR_{\CL(\CY_q(\BR^N_+), H^{3-k}_q(\BR^N_+))}
(\{(\tau\pd_\tau)^\ell(\lambda^{k}\CH_0(\lambda)) \mid 
\lambda \in \Lambda_{\kappa, \lambda_0}\}) \leq 
r_b, 
\end{split}\end{equation}
for $\ell=0,1$, $j=0,1,2$ and $k=0,1$.  Here, $r_b$ is a
constant depending on $m_0$, $m_1$, $m_2$, 
$\lambda_0$, $q$, and $N$.
\end{thm}
\begin{remark}
In this section, what the constant $c$ depends on $m_0$, $m_1$, $m_2$
means that the constant $c$ depends on $m_0$, $m_1$, $m_2$, 
but is independent of $\mu$, $\sigma$ and $A_\kappa$ whenever
$\mu \in [m_0, m_1]$,  $\sigma \in [m_0, m_1]$, and $|A_\kappa| \leq m_2$
for $\kappa \in [0, 1)$. 
\end{remark}
To prove Theorem \ref{main:half1}, 
as an auxiliary problem, we first consider the following equations:
\begin{equation}\label{4.3.1}\left\{\begin{aligned}
\lambda \bv - \DV(\mu\bD(\bv) - \theta\bI) =0,
\quad \dv \bv & = 0
&\quad&\text{in $\BR^N_+$}, \\
(\mu\bD(\bv) - \theta\bI)\bn_0 &= \bh
&\quad&\text{on $\BR^N_0$}, 
\end{aligned}\right.\end{equation}
and we shall prove the following theorem,  
which was essentially proved by Shibata and Shimizu \cite{SS4}.
\begin{thm}\label{lem:sol1} Let $1 < q < \infty$, $\epsilon \in (0, \pi/2)$,
and $\lambda_0 > 0$. Let 
\begin{align*}
\CY'_q(\BR^N_+) &= \{(G_1, G_2) \mid G_1 \in L_q(\BR^N_+)^N,
\quad G_2 \in H^1_q(\BR^N)^N\}, \\
\hat H^1_q(\BR^N_+) &= \{\theta \in L_{q, {\rm loc}}(\BR^N_+) \mid 
\nabla\theta \in L_q(\BR^N_+)\}.
\end{align*}
 Then, there exists a solution operator $\CV(\lambda)$ with 
$$\CV(\lambda) \in \Hol(\Sigma_{\epsilon,\lambda_0},  
\CL(\CY'(\BR^N_+), H^2_q(\BR^N_+)^N))$$
such that
for any $\lambda=\gamma + i\tau \in \Sigma_{\epsilon, \lambda_0}$ and 
$\bh \in H^1_q(\BR^N_+)^N$, 
$\bv = \CV(\lambda)(\lambda^{1/2}\bh, \bh)$ 
are  unique solutions of
Eq. \eqref{4.3.1} with some $\theta \in \hat H^1_q(\HS)$ and 
\begin{align*}
\CR_{\CL(\CY'_q(\BR^N_+), H^{2-j}_q(\BR^N_+)^N)}
(\{(\tau\pd_\tau)^\ell(\lambda^{j/2}\CV(\lambda)) \mid 
\lambda \in \Sigma_{\epsilon, \lambda_0}\}) &\leq r_b(\lambda_0)
\end{align*}
for $\ell=0,1$, and $j=0,1,2$. Here, $r_b(\lambda_0)$ 
is a constant depending on $m_0$, $m_1$, $m_2$, $\epsilon$,
$\lambda_0$, $N$, and $q$.
\end{thm}
\pf To prove Theorem \ref{lem:sol1}, we start with the 
solution formulas of Eq. \eqref{4.3.1}, which were obtained in
Shibata and Shimizu \cite{SS4} essentially, but for the sake of
the completeness of the paper as much as possible and also for the later use,
we will derive them in the following.  
Applying the partial Fourier transform with respect to
$x' = (x_1, \ldots, x_{N-1})$ to Eq. \eqref{4.3}, we have
\begin{equation}\label{4.4}\left\{\begin{aligned}
\lambda \hat v_j + \mu|\xi'|^2 -\pd_N^2\hat v_j 
+ i\xi_j\hat \theta &=0, 
\\ 
\lambda \hat v_N + \mu|\xi'|^2 -\pd_N^2\hat v_N 
+ \pd_N\hat \theta &=0
&\quad&(x_N > 0), \\
\sum_{j=1}^{N-1} i\xi_j\hat v_j + \pd_N\hat v_N & = 0 
&\quad&(x_N > 0), \\
\mu(\pd_N\hat v_j + i\xi_j\hat v_N) = g_j,\quad 
2\mu \pd_N\hat v_N - \hat \theta &=g_N&\quad&\text{for $x_N=0$}.
\end{aligned}\right.\end{equation}
Here, for $f=f(x', x_N)$, $x'=(x_1, \ldots, x_{N-1}) \in \BR^{N-1}$, 
$x_N \in (a, b)$, $\hat f$ denotes 
the partial Fourier transform  of $f$ with respect
to $x'$  defined by 
\begin{align*}
\hat f(\xi', x_N)= \CF'[f(\cdot, x_N)](\xi') 
= \int_{\BR^{N-1}} e^{-ix'\cdot\xi'}f(x', x_N)\,dx'
\end{align*}
with $\xi' = (\xi_1, \ldots, \xi_{N-1}) \in \BR^{N-1}$ and
 $x'\cdot\xi' = \sum_{j=1}^{N-1}x_j\xi_j$, and we have
set $g_j = -\hat h_j(\xi', 0)$.  
To obtain solution formula, we set 
$$\hat v_j = \alpha_je^{-Ax_N} + \beta_je^{-Bx_N},
\quad \hat\theta = \omega e^{-Ax_N}$$
with $A = |\xi'|$ and $B = \sqrt{\lambda\mu^{-1} + |\xi'|^2}$,  
and then from \eqref{4.4} we have
\begin{align}
&\mu\alpha_j(B^2-A^2) + i\xi_i\omega = 0, \quad
\mu\alpha_N(B^2-A^2) - A\omega = 0, \label{1}\\
&\sum_{k=1}^{N-1}i\xi_k\alpha_k - A\alpha_N=0, \quad
\sum_{k=1}^{N-1}i\xi_k\beta_k - B\beta_N = 0, 
\label{2}\\
&\mu\{(A\alpha_j + B\beta_j) - i\xi_j(\alpha_N + \beta_N)\}
= g_j, \label{3}\\
&2\mu(A\alpha_N + B\beta_N) + \omega = g_N.
\label{4}
\end{align}
The solution formula of Eq. \eqref{4.3} was given in Shibata and Shimizu
\cite{SS4}, but there is an error in the formula in \cite[\thetag{4.17}]{SS4}
such as  
$$\mu\{(A\alpha_j + B\beta_j) + i\xi_j(\alpha_N + \beta_N)\}
= \hat h_j(\xi', 0),$$
which should read 
$$\mu\{(A\alpha_j + B\beta_j) - i\xi_j(\alpha_N + \beta_N)\}
= -\hat h_j(\xi', 0)$$
as \eqref{3} above. The formulas obtained in \cite{SS4} are
 correct, but we repeat here how to obtain $\alpha_j$, $\beta_j$ and $\omega$,
 because this error confuses readers. 

We first drive $2\times2$ system of equations with respect to
$\alpha_N$ and $\beta_N$. Multiplying \eqref{3} with $i\xi_j$, 
summing up the resultant formulas from $j=1$
through $N-1$ and writing $i\xi'\cdot m' = \sum_{j=1}^{N-1}i\xi_jm_j$
for $m_j \in \{\alpha_j, \beta_j, g_j\}$ give 
\begin{equation}\label{8}
 Ai\xi'\cdot\alpha' + Bi\xi'\cdot\beta'
+A^2(\alpha_N+\beta_N) = \mu^{-1}i\xi'\cdot g'.
\end{equation} 
By \eqref{2}, 
\begin{equation}\label{7}
i\xi'\cdot\alpha' = A\alpha_N,
\quad i\xi'\cdot\beta' = B\beta_N,
\end{equation}
 which, combined with
\eqref{8},  leads to 
\begin{equation}\label{9}
2A^2 \alpha_N +(A^2 + B^2)\beta_N = \mu^{-1}i\xi'\cdot g'.
\end{equation}
By \eqref{1}, 
\begin{equation}\label{10}
\omega = \frac{\mu (B^2-A^2)}{A}\alpha_N,
\end{equation}
which, combined with \eqref{4}, leads to
\begin{equation}\label{11}
(A^2 + B^2)\alpha_N + 2AB\beta_N = \mu^{-1}Ag_N.
\end{equation}
Thus, setting
$$\CL = \left(\begin{matrix}
A^2 + B^2 & 2A^2 \\ 
2AB & A^2 + B^2 
\end{matrix}\right)
\quad(\text{Lopatinski matrix}),
$$
we have
$$\CL\left(\begin{matrix} \beta_N \\ \alpha_N
\end{matrix}\right) = 
\left(\begin{matrix} \mu^{-1}i\xi'\cdot g' \\
\mu^{-1}Ag_N \end{matrix}\right).
$$
Since 
$$\det \CL = (A^2+B^2)^2-4A^3B = A^4-4A^3B + 2A^2B^2 + B^4
=(B-A)D(A,B)$$
with 
$$D(A,B) = B^3 + AB^2 + \textcolor{red}{3A^2B^2
\to 3A^2B} - A^3,$$
we have 
$$\CL^{-1} = \frac{1}{(B-A)D(A,B)}\left(\begin{matrix}
A^2 + B^2 & -2A^2 \\
-2AB & A^2 + B^2
\end{matrix}\right).
$$
Thus, we have 
\begin{equation}\label{16}\begin{split}
\beta_N & = \frac{1}{\mu(B-A)D(A,B)}((A^2 + B^2)i\xi'\cdot g'
-2A^3g_N), \\
\alpha_N & = \frac{-1}{\mu(B-A)D(A,B)}
(2ABi\xi'\cdot g' - (A^2+B^2)Ag_N)).
\end{split}\end{equation}
In particular, 
$$\hat v_N = \alpha_Ne^{-Ax_N} + \beta_Ne^{-Bx_N}
= \alpha_N(e^{-Ax_N} - e^{-Bx_N})
+(\alpha_N+\beta_N)e^{-Bx_N}.
$$
We have
\begin{equation}\label{12}\begin{split}
&\alpha_N+\beta_N \\
& = \frac{1}{(B-A)D(A,B)}
((A^2+B^2-2AB)i\xi'\cdot g' +((A^2+B^2)A-2A^3)g_N)\\
&= \frac{1}{\mu(B-A)D(A,B)}
((B-A)^2i\xi'\cdot g' +A(B^2-A^2)g_N)\\
&=\frac{1}{\mu(B-A)D(A,B)}
((B-A)^2i\xi'\cdot g' +A(B-A)(A+B)g_N)\\
& = \frac{1}{\mu D(A,B)}
((B-A)i\xi'\cdot g'+A(A+B)g_N).
\end{split}\end{equation}
Setting
$$\CM(x_N) = \frac{e^{-Bx_N} - e^{-Ax_N}}{B-A},$$
we have
\begin{equation}\label{12*}\begin{aligned}
\hat v_N&= \frac{ A}{\mu D(A,B)}\CM(x_N)
(2Bi\xi'\cdot g' - (A^2+B^2)g_N)\\
&\quad+ \frac{e^{-Bx_N}}{\mu D(A,B)}
((B-A)i\xi'\cdot g' + A(A+B)g_N).
\end{aligned}\end{equation}
By \eqref{10} and \eqref{16}, 
\begin{align*}
&\omega = \frac{\mu(B^2-A^2)}{A}\alpha_N\\
&= \frac{\mu(B^2-A^2)}{A}
\frac{-1}{\mu(B-A)D(A,B)}
(2ABi\xi'\cdot g' - (A^2+B^2)Ag_N))
\\
& = -\frac{(A+B)}{ D(A,B)}
(2Bi\xi'\cdot g' - (A^2+B^2)g_N))
\end{align*}
and so
\begin{equation}\label{16*}
\hat \theta = 
-\frac{(A+B)e^{-Ax_N}}{ D(A,B)}
(2Bi\xi'\cdot g' - (A^2+B^2)g_N)).
\end{equation}
By \eqref{1}, 
\begin{equation}\label{13}\begin{aligned}
\alpha_j &= -\frac{i\xi_j}{\mu(B^2 - A^2)}\omega\\
&= \frac{i\xi_j}{\mu(B^2-A^2)}\frac{A+B}{D(A,B)}
(2Bi\xi'\cdot g' - (A^2+B^2)g_N))\\
&=\frac{i\xi_j}{\mu(B-A)D(A,B)}
(2Bi\xi'\cdot g' - (A^2+B^2)g_N).
\end{aligned}\end{equation}
By \eqref{3}
$$\beta_j = \frac{1}{\mu B}g_j + \frac{1}{B}
(i\xi_j(\alpha_N + \beta_N)-A\alpha_j).$$
By \eqref{12} and \eqref{13}
\begin{align*}
&i\xi_j(\alpha_N + \beta_N)-A\alpha_j \\
&= 
\frac{i\xi_j}{\mu(B-A)D(A,B)}\{(B-A)^2i\xi'\cdot g'
+A(B-A)(A+B)g_N\\
&\phantom{\frac{i\xi_j}{\mu(B-A)D(A,B)}\{}\quad
 -A
(2Bi\xi'\cdot g' - (A^2+B^2)g_N)\}\\
&=
\frac{i\xi_j}{\mu(B-A)D(A,B)}\{(A^2-4AB+B^2)i\xi'\cdot g'
+2AB^2g_N)\},
\end{align*}
and therefore
\begin{equation}\label{14}
\beta_j = \frac{1}{\mu B}g_j
+\frac{i\xi_j}{\mu(B-A)D(A,B)B}\{(A^2-4AB+B^2)i\xi'\cdot g'
+2AB^2g_N)\}.
\end{equation}
Combining \eqref{13} and \eqref{14} gives
\begin{align*}
\hat v_j &= \frac{e^{-Bx_N}}{\mu B}g_j
+ \frac{i\xi_j e^{-Ax_N}}
{\mu(B-A)D(A,B)}\{2Bi\xi'\cdot g' - (A^2+B^2)g_N\}\\
&+ \frac{i\xi_j e^{-Bx_N}}
{\mu(B-A)D(A,B)B}\{(A^2-4AB+B^2)i\xi'\cdot g'
+2AB^2g_N)\}\\
& = \frac{1}{\mu B}g_j
+ Ii\xi'\cdot g' + II g_N,
\end{align*}
with
\begin{align*}
I & = \frac{i\xi_je^{-Ax_N}}{\mu(B-A)D(A,B)}
2B + \frac{i\xi_je^{-Bx_N}}{\mu(B-A)D(A,B)B}(A^2-4AB + B^2),\\
II & = -\frac{i\xi_j e^{-Ax_N}}{\mu(B-A)D(A,B)}(A^2+B^2)
+ \frac{i\xi_j e^{-Bx_N}}
{\mu(B-A)D(A,B)}2AB
\end{align*}
We proceed as follows:
\begin{align*}
I & = \frac{i\xi_j(e^{-Ax_N}-e^{-Bx_N})}{\mu(B-A)D(A,B)}2B
+ \frac{i\xi_je^{-Bx_N}}{\mu(B-A)D(A,B)B}(A^2-4AB + 3B^2)\\
& = -\frac{2i\xi_jB\CM(x_N)}{\mu D(A,B)}
+ \frac{i\xi_j(3B-A)e^{-Bx_N}}{\mu D(A,B)B}; \\
II & = -\frac{i\xi_j (e^{-Ax_N}-e^{-Bx_N})}{\mu(B-A)D(A,B)}(A^2+B^2)
- \frac{i\xi_j e^{-Bx_N}(A^2-2AB+B^2)}
{\mu(B-A)D(A,B)}
\\
& = \frac{i\xi_j (A^2+B^2)\CM(x_N)}{\mu D(A,B)}
 -\frac{i\xi_je^{-Bx_N}(B-A)}{\mu D(A,B)}.
\end{align*}
Therefore, we have
\begin{equation}\label{15}\begin{aligned}
\hat v_j  &= \frac{e^{-Bx_N}}{\mu B}g_j
-\frac{i\xi_j \CM(x_N)}{\mu D(A,B)}(2Bi\xi'\cdot g'-(A^2+B^2)g_N) 
\\
&\quad
 + \frac{i\xi_je^{-Bx_N}}{\mu D(A,B)B}((3B-A)i\xi'\cdot g' - B(B-A)g_N).
\end{aligned}\end{equation}

To define solution operators for Eq. \eqref{4.3}, 
we make preparations. 
\begin{lem}\label{lem:4.3}
Let $s \in \BR$ and $0 < \epsilon < \pi/2$. 
Then, there exists a positive constant $c$ depending
on $\epsilon$, $m_1$ and $m_2$ such that
\begin{gather}
c(|\lambda|^{1/2} + A) \leq {\rm Re}\,B \leq |B| \leq (\mu^{-1}|\lambda|)^{1/2}
+ A, \label{ineq:4.1.1}\\
c(|\lambda|^{1/2} + A)^3 \leq |D(A,B)| \leq 6((\mu^{-1}|\lambda|)^{1/2}
+ A)^3. \label{ineq:4.1.2}
\end{gather}
for any $\lambda \in \Sigma_\epsilon$ and $\mu \in [m_1, m_2]$. 
\end{lem}
\pf  The inequality in the left side of \eqref{ineq:4.1.1}
 follows immediately from 
Lemma \ref{lem:fund1}.  Notice that
\begin{align*}
D(A,B) &= B^3 +3A^2B + AB^2 - A^3 
= B(B^2 + 2A^2) + A(A^2+\mu^{-1}\lambda)-A^3 \\
&= B(\mu^{-1}\lambda + 4A^2) + \mu^{-1}A\lambda.
\end{align*}
If we consider the angle of $B(\mu^{-1}\lambda + 4A^2)$ and 
$-\mu^{-1}A\lambda$, then we see easily that $D(A,B) \not=0$.  
Thus, studying the following 
three cases: $R_1|\lambda|^{1/2} \leq A$, $R_1A \leq |\lambda|^{1/2}$ and 
$R_1^{-1}A\leq |\lambda|^{1/2} \leq R_1A$ for sufficient large $R_1 > 0$,
we can prove the inequality in the left side of \eqref{ineq:4.1.2}. 
The detailed proof was given in Shibata and Shimizu \cite{SS1}. 
The independence of the constant c of $\lambda \in \Sigma_\epsilon$
and $\mu \in [m_0, m_1]$ 
follows  from the 
homogeneity: $\sqrt{\mu^{-1}(m^2\lambda) + (mA)^2} = m\sqrt{\mu^{-1}\lambda
+ A^2}$ and $D(mA, mB) = m^3D(A,B)$ for any $m>0$ and the compactness of
the interval $[m_0, m_1]$. 
\qed\vskip0.5pc

To introduce the key tool of proving the $\CR$ boundedness in
the half space, we make a definition.
\begin{dfn}\label{def:4.1} Let $V$ be a domain in $\BC$, 
let $\Xi = V\times (\BR^{N-1}\setminus\{0\})$,  and let $m: 
\Xi \to \BC$; $(\lambda, \xi') \mapsto m(\lambda, \xi')$ be 
$C^1$ with respect to $\tau$, where $\lambda = \gamma + i\tau \in V$, 
and $C^\infty$ with respect to $\xi' \in \BR^{N-1}\setminus\{0\}$. 
\begin{itemize}
\item[\thetag1]~$m(\lambda, \xi')$ is called a multiplier of order
$s$ with type $1$ on $\Xi$, if the estimates:
\begin{align*}
|\pd_{\xi'}^{\kappa'}m(\lambda, \xi')| &\leq C_{\kappa'}(|\lambda|^{1/2}
+ |\xi'|)^{s-|\kappa'|}, \\
|\pd_{\xi'}^{\kappa'}(\tau\pd_\tau m(\lambda, \xi'))| 
&\leq C_{\kappa'}(|\lambda|^{1/2} + |\xi'|)^{s-|\kappa'|}
\end{align*}
hold for any multi-index $\kappa \in \BN_0^{N}$ and 
$(\lambda, \xi') \in \Xi$ with some constant $C_{\kappa'}$ depending
solely on $\kappa'$ and $V$. 
\item[\thetag2]~$m(\lambda, \xi')$ is called a multiplier of order
$s$ with type $2$ on $\Xi$, if the estimates:
\begin{align*}
|\pd_{\xi'}^{\kappa'}m(\lambda, \xi')| &\leq C_{\kappa'}(|\lambda|^{1/2}
+ |\xi'|)^s|\xi'|^{-|\kappa'|}, \\
|\pd_{\xi'}^{\kappa'}(\tau\pd_\tau m(\lambda, \xi'))| 
&\leq C_{\kappa'}(|\lambda|^{1/2} + |\xi'|)^s|\xi'|^{-|\kappa'|}
\end{align*}
hold for any multi-index $\kappa \in \BN_0^{N}$ and 
$(\lambda, \xi') \in \Xi$ with some constant $C_{\kappa'}$ depending
solely on $\kappa'$ and $V$. 
\end{itemize}
Let $\bM_{s,i}(V)$ be the set of all multipliers of order $s$ with type $i$
on $\Xi$ for $i=1,2$. For $m \in \bM_{s,i}(V)$, we set  
$M(m, V) = \max_{|\kappa'| \leq N} C_{\kappa'}$. 
\end{dfn}
Let  $\CF^{-1}_{\xi'}$ be the inverse partial 
Fourier transform defined by
$$\CF^{-1}_{\xi'}[f(\xi', x_N)](x') = \frac{1}{(2\pi)^{N-1}}\int_{\BR^{N-1}} 
e^{ix'\cdot\xi'}f(\xi', x_N)\,d\xi'.
$$
Then, we have the following two lemmata which were proved essentially by  
Shibata and Shimizu \cite[Lemma 5.4 and Lemma 5.6]{SS4}.
\begin{lem} \label{lem:4.1} Let $0 < \epsilon < \pi/2$, $1 < q < \infty$,
and $\lambda_0 > 0$.  Given  
$m \in \bM_{-2,1}(\Lambda_{\kappa, \lambda_0})$, we define an 
operator $L(\lambda)$  by 
$$[L(\lambda)g](x)  = \int^\infty_0\CF^{-1}_{\xi'}
[m(\lambda, \xi')\lambda^{1/2}e^{-B(x_N+y_N)}\hat g(\xi', y_N)](x')\,dy_N.
$$
Then, we have
$$\CR_{\CL(L_q(\BR^N_+), H^{2-j}_q(\HS)^N)}(\{(\tau\pd_\tau)^\ell
(\lambda^{j/2}\pd_x^\alpha L_i(\lambda))
\mid \lambda \in \Lambda_{\kappa, \lambda_0}\})
\leq r_b(\lambda_0) 
$$
for any $\ell=0,1$ and $j=0,1,2$.   Where $\tau$ denotes the imaginary part
of $\lambda$, and  $r_b(\lambda_0)$ is a constant depending on 
$M(m, \Lambda_{\kappa, \lambda_0})$, $\epsilon$, $\lambda_0$, $N$,
and $q$.    
\end{lem}
\begin{lem} \label{lem:4.2} Let $0 < \epsilon < \pi/2$, $1 < q < \infty$,
and $\lambda_0 > 0$.  Given  
$m \in \bM_{-2,2}(\Lambda_{\kappa, \lambda_0})$, we define
operators $L_i(\lambda)$ $(i=1, \ldots, 4)$ by 
\begin{align*}
[L_1(\lambda)g](x) & = \int^\infty_0\CF^{-1}_{\xi'}
[m(\lambda, \xi')Ae^{-B(x_N+y_N)}\hat g(\xi', y_N)](x')\,dy_N, \\
[L_2(\lambda)g](x) & = \int^\infty_0\CF^{-1}_{\xi'}
[m(\lambda, \xi')Ae^{-A(x_N+y_N)}\hat g(\xi', y_N)](x')\,dy_N, \\
[L_3(\lambda)g](x) & = \int^\infty_0\CF^{-1}_{\xi'}
[m(\lambda, \xi')A^2\CM(x_N+y_N)\hat g(\xi', y_N)](x')\,dy_N, \\
[L_4(\lambda)g](x) & = \int^\infty_0\CF^{-1}_{\xi'}
[m(\lambda, \xi')\lambda^{1/2}A\CM(x_N+y_N)
\hat g(\xi', y_N)](x')\,dy_N.
\end{align*}
Then, we have
$$\CR_{\CL(L_q(\BR^N_+), H^{2-j}_q(\HS)^N)}(\{(\tau\pd_\tau)^\ell
(\lambda^{j/2}\pd_x^\alpha L_i(\lambda))
\mid \lambda \in \Lambda_{\kappa, \lambda_0}\})
\leq r_b(\lambda_0)
$$
for $\ell=0,1$ and $j=0,1,2$.
Where $\tau$ denotes the imaginary part
of $\lambda$, and  $r_b(\lambda_0)$ is a constant depending on 
$M(m, \Lambda_{\kappa, \lambda_0})$, $\epsilon$, $\lambda_0$, $N$,
and $q$.   
\end{lem}
To construct solution operators, we use the following
lemma.
\begin{lem}\label{lem:solop} Let $0 < \epsilon <\pi/2$, 
$1 < q < \infty$ and $\lambda_0 > 0$. Given multipliers,
$n_1 \in \bM_{-2,1}(\Lambda_{\kappa, \lambda_0})$, 
$n_2 \in \bM_{-2,2}(\Lambda_{\kappa, \lambda_0})$, 
and $n_3 \in \bM_{-1,2}(\Lambda_{\kappa, \lambda_0})$, 
we define operators $T_i(\lambda)$ $(i=1, 2, 3)$ by
\begin{align*}
T_1(\lambda)h & = \CF^{-1}_{\xi'}[\lambda^{1/2}e^{-Bx_N}n_1(\lambda, \xi')
\hat h(\xi', 0)](x'), \\
T_2(\lambda)h & = \CF^{-1}_{\xi'}[Ae^{-Bx_N}n_2(\lambda, \xi')
\hat h(\xi', 0)](x'), \\
T_3(\lambda)h & = \CF^{-1}_{\xi'}[A\CM(x_N)n_3(\lambda, \xi')
\hat h(\xi', 0)](x').
\end{align*}
Let 
$$\CZ_q(\BR^N_+) = \{(G_3, G_4) \mid G_3 \in L_q(\HS), \enskip
G_4 \in H^1_q(\HS)\}. $$
Then, there exist operator families $\CT_i(\lambda)
\in {\rm Hol}\,(\Lambda_{\kappa, \lambda_0}, \CL(\CY'_q(\HS), 
H^2_q(\HS)))$ such that for any $\lambda = \gamma + i\tau
\in \Lambda_{\kappa, \lambda_0}$
and $h \in H^1_q(\HS)$, $T_i(\lambda)h = \CT_i(\lambda)(\lambda^{1/2}h, h)$ 
and 
\begin{equation}\label{lem:eq.1}
\CR_{\CL(\CY'_q(\HS), H^{2-j}_q(\HS))}
(\{(\tau\pd_\tau)^\ell(\lambda^{j/2}\CT_i(\lambda)) \mid
\lambda \in \Lambda_{\kappa, \lambda_0}\}) 
\leq r_b(\lambda_0)
\end{equation}
for $\ell=0,1$, $j=0,1,2$.  Where $r_b(\lambda_0)$ is a constant depending on 
$M(n_i, \Lambda_{\kappa, \lambda_0})$ $(i=1,2,3)$, 
$\epsilon$, $\lambda_0$, $N$,
and $q$.
\end{lem}
\pf By Volevich's trick we write 
\begin{align*}
&T_1(\lambda)h \\
& = -\int^\infty_0\CF_{\xi'}^{-1}[\frac{\pd}{\pd y_N}
(\lambda^{1/2}e^{-B(x_N+y_N)}n_1(\lambda, \xi')\hat h(\xi', y_N))](x')\,dy_N
\\
& = -\int^\infty_0\CF^{-1}_{\xi'}[\lambda^{1/2}e^{-B(x_N+y_N)}
n_1(\lambda, \xi')\pd_N\hat h(\xi', y_N)](x')\,dy_N \\
& \quad + \int^\infty_0\CF^{-1}_{\xi'}[\lambda^{1/2}e^{-B(x_N+y_N)}
\frac{\lambda^{1/2}}{\mu B}n_1(\lambda, \xi')\lambda^{1/2}\hat h(\xi', y_N)
](x')\,dy_N\\
& \quad- \sum_{j=1}^{N-1}\int^\infty_0\CF_{\xi'}^{-1}[Ae^{-B(x_N+y_N)}
\frac{\lambda^{1/2}}{B}\frac{i\xi_j}{A}n_1(\lambda, \xi')
\CF[\pd_j h(\cdot, y_N)]](x')\,dy_N,
\end{align*}
where we have used the formula:
$$B = \frac{\mu^{-1}\lambda + A^2}{\mu B}= \frac{\lambda}{\mu B}
-\sum_{j=1}^{N-1}\frac{A}{B}\frac{i\xi_j}{A}i\xi_j.
$$
Let
\begin{align*}
&\CT_1(\lambda)(G_3, G_4)\\ 
& = -\int^\infty_0\CF^{-1}_{\xi'}[\lambda^{1/2}e^{-B(x_N+y_N)}
n_1(\lambda, \xi')\CF[\pd_N G_4(\cdot, y_N)]](x')\,dy_N \\
& \quad + \int^\infty_0\CF^{-1}_{\xi'}[\lambda^{1/2}e^{-B(x_N+y_N)}
\frac{\lambda^{1/2}}{\mu B}n_1(\lambda, \xi')
\CF[G_3(\cdot, y_N)]](x')\,dy_N\\
& \quad- \sum_{j=1}^{N-1}\int^\infty_0\CF_{\xi'}^{-1}[Ae^{-B(x_N+y_N)}
\frac{\lambda^{1/2}}{B}\frac{i\xi_j}{A}n_1(\lambda, \xi')
\CF[\pd_j G_4(\cdot, y_N)]](x')\,dy_N, 
\end{align*}
and then, $T_1(\lambda)h = \CT_1(\lambda)(\lambda^{1/2}h, h)$. 
Moreover, Lemma \ref{lem:4.1} and Lemma \ref{lem:4.2} yield
\eqref{lem:eq.1} with $j=1$, because
\begin{align*}
n_1(\lambda, \xi') &\in \bM_{-2,1}(\Lambda_{\kappa,\lambda_0}), \quad
\frac{\lambda^{1/2}}{\mu B}n_1(\lambda, \xi') 
\in \bM_{-2,1}(\Lambda_{\kappa,\lambda_0}), \\
\frac{\lambda^{1/2}}{B}\frac{i\xi_j}{A}n_1(\lambda, \xi')
 &\in \bM_{-2,2}(\Lambda_{\kappa,\lambda_0}).
\end{align*}
Analogously, we can prove the existence of $\CT_2(\lambda)$.

To construct $\CT_3(\lambda)$, we use the formula:
$$\frac{\pd}{\pd x_N}\CM(x_N) = -e^{-Bx_N}-A\CM(x_N), $$
and then, by  Volevich's trick we have
\begin{align*}
&T_3(\lambda)h \\
&= -\int^\infty_0\CF_{\xi'}^{-1}[\frac{\pd}{\pd y_N}
(A\CM(x_N+y_N)n_3(\lambda, \xi')\hat h(\xi', y_N))](x')\,dy_N
= -I + II
\end{align*}
with 
\begin{align*}
I & = \int^\infty_0\CF^{-1}_{\xi'}[A\CM(x_N+y_N)n_3(\lambda, \xi')
\pd_N\hat h(\xi', y_N)](x')\,dy_N; \\
II & = \int^\infty_0\CF^{-1}_{\xi'}[
(Ae^{-B(x_N+y_N)} + A^2\CM(x_N+y_N))
n_3(\lambda, \xi')\hat h(\xi', y_N)](x')\,dy_N.
\end{align*}
Using the formula:
$$1 = \frac{B^2}{B^2} = \frac{\lambda^{1/2}}{\mu B^2}\lambda^{1/2}
+\frac{A}{B^2}A
= \frac{\lambda^{1/2}}{\mu B^2}\lambda^{1/2}
- \sum_{j=1}^{N-1}\frac{i\xi_j}{B^2}i\xi_j,
$$
we have 
\allowdisplaybreaks{
\begin{align*}
I &= \int^\infty_0\CF^{-1}_{\xi'}[\lambda^{1/2}A\CM(x_N+y_N)
\frac{\lambda^{1/2}}{\mu B^2}n_3(\lambda, \xi') \pd_N\hat h(\xi', y_N)]
(x')\,dy_N \\
&+ \int^\infty_0\CF^{-1}_{\xi'}[A^2\CM(x_N+y_N)
\frac{A}{ B^2}n_3(\lambda, \xi') \pd_N\hat h(\xi', y_N)]
(x')\,dy_N; \\
II & = \int^\infty_0\CF^{-1}_{\xi'}[
(Ae^{-B(x_N+y_N)} + A^2\CM(x_N+y_N))
\\
&\phantom{aaaaaaaaaaaaaaaaaaaa}\times
\frac{\lambda^{1/2}}{\mu B^2}n_3(\lambda, \xi')        
\lambda^{1/2}\hat h(\xi', y_N)](x')\,dy_N \\
&-\sum_{j=1}^{N-1}
\int^\infty_0\CF^{-1}_{\xi'}[
(Ae^{-B(x_N+y_N)} + A^2\CM(x_N+y_N))\\
&\phantom{aaaaaaaaaaaaaaaaaaaa}\times
\frac{i\xi_j}{B^2}n_3(\lambda, \xi')
\CF[\pd_jh(\cdot, y_N)]](x')\,dy_N.
\end{align*}
}
Let 
\allowdisplaybreaks{
\begin{align*}
\CT_3&
(\lambda)(G_3, G_4) \\
&= -\int^\infty_0\CF^{-1}_{\xi'}[\lambda^{1/2}A\CM(x_N+y_N)
\frac{\lambda^{1/2}}{\mu B^2}n_3(\lambda, \xi') 
\CF[\pd_N G_4(\cdot, y_N)]]
(x')\,dy_N \\
&- \int^\infty_0\CF^{-1}_{\xi'}[A^2\CM(x_N+y_N)
\frac{A}{ B^2}n_3(\lambda, \xi')\CF[\pd_N G_4(\cdot, y_N)]]
(x')\,dy_N \\
& + \int^\infty_0\CF^{-1}_{\xi'}[
(Ae^{-B(x_N+y_N)} + A^2\CM(x_N+y_N))\\
&\phantom{aaaaaaaaaaaaaaa}\times
\frac{\lambda^{1/2}}{\mu B^2}n_3(\lambda, \xi')
\CF[G_3(\cdot, y_N)]](x')\,dy_N \\
&-\sum_{j=1}^{N-1}
\int^\infty_0\CF^{-1}_{\xi'}[
(Ae^{-B(x_N+y_N)} + A^2\CM(x_N+y_N))\\
&\phantom{aaaaaaaaaaaaaaa}\times
\frac{i\xi_j}{B^2}n_3(\lambda, \xi')
\CF[\pd_jG_4(\cdot, y_N)]](x')\,dy_N,
\end{align*}
}
and then $T_3(\lambda)h = \CT_3(\lambda)(\lambda^{1/2}h, h)$.
Moreover, Lemma \ref{lem:4.2} yields \eqref{lem:eq.1}
for $j=3$, because 
\begin{align*}
&n_3(\lambda)\in 
\bM_{-2,2}(\Lambda_{\kappa, \lambda_0}),
\quad
 \frac{\lambda^{1/2}}{\mu B^2}n_3(\lambda, \xi')
\in 
\bM_{-2,2}(\Lambda_{\kappa, \lambda_0}), \\ 
&\frac{i\xi_j}{B^2}n_3(\lambda, \xi') \in 
\bM_{-2,2}(\Lambda_{\kappa, \lambda_0}).
\end{align*}
This completes the proof of Lemma \ref{lem:solop}. 
\qed
\vskip0.5pc


\noindent
{\bf Continuation of Proof of Theorem \ref{lem:sol1}}.~ 
Let 
$$v_j(x) = \CF^{-1}_{\xi'}[\hat v_j(\xi', x_N)](x'),
$$
and then by \eqref{12*} and \eqref{15} we have
\begin{align*}
v_N & = \CF^{-1}_{\xi'}\Bigl[\frac{ A}{\mu D(A,B)}\CM(x_N)
((A^2+B^2)\hat h_N(\xi', 0)
-2B\sum_{\ell=1}^{N-1}i\xi_\ell \hat h_\ell(\xi', 0))\Bigr](x')\\
&- \CF^{-1}_{\xi'}\Bigl[\frac{Ae^{-Bx_N}}{\mu D(A,B)}
((B-A)\sum_{\ell=1}^{N-1}\frac{i\xi_\ell}{A}
\hat h_\ell(\xi', 0) + (A+B)\hat h_N(\xi', 0))\Bigr](x');\\
v_k & = -\CF^{-1}_{\xi'}\Bigl[\frac{\lambda^{1/2}}{\mu^2 B^3}
\lambda^{1/2}e^{-Bx_N}\hat h_k(\xi', 0)\Bigr](x') 
- \CF^{-1}_{\xi'}\Bigl[\frac{A}{\mu B^3}Ae^{-Bx_N}\hat h_k(\xi', 0)\Bigr](x')
\\
&+\CF^{-1}_{\xi'}\Bigl[A \CM(x_N)\frac{i\xi_k}{A}
\frac{1}{\mu D(A,B)}(2B\sum_{\ell=1}^{N-1}i\xi_\ell
\hat h_\ell(\xi', 0)\\
&\phantom{aaaaaaaaaaaaaaaaaaaaaaaaaaaaaaaaaaa}
-(A^2+B^2)\hat h_N(\xi', 0))\Bigr](x')\\
&- \CF^{-1}_{\xi'}\Bigl[Ae^{-Bx_N}
\frac{i\xi_k}{A}\frac{1}{\mu D(A,B)B}((3B-A)
\sum_{\ell=1}^{N-1}i\xi_\ell\hat h_\ell(\xi', 0)\\
&\phantom{aaaaaaaaaaaaaaaaaaaaaaaaaaaaaaaaaaa}
 - B(B-A)\hat h_N(\xi', 0))
\Bigr](x'),
\end{align*}
for $k=1, \ldots, N-1$, 
where we have used the formula $\dfrac{1}{\mu B} = \dfrac{\lambda }{ \mu^2 B^3}
+ \dfrac{A^2}{\mu B^3}$ to treat the first term of $\hat v_j$ in \eqref{15}.
Since
\begin{align*}
\frac{Bi\xi_\ell}{\mu D(A,B)},\enskip \frac{A^2+B^2}{\mu D(A,B)},\enskip
\frac{i\xi_k}{A}\frac{Bi\xi_\ell}{\mu D(A, B)}, \enskip
\frac{i\xi_k}{A}\frac{A^2+B^2}{\mu D(A,B)} &\in 
\bM_{-1,2}(\Sigma_{\epsilon, \lambda_0}), \\
\frac{B-A}{\mu D(A,B)}\frac{i\xi_\ell}{A}, \enskip
\frac{A+B}{\mu D(A,B)}, \enskip \frac{A}{\mu B^3}
&\in \bM_{-2,2}(\Sigma_{\epsilon, \lambda_0})\\
\frac{i\xi_k}{A}\frac{(3B-A)i\xi_\ell}{\mu D(A,B)B},
\enskip \frac{i\xi_k}{A}\frac{B(B-A)}{\mu D(A,B)B}
&\in \bM_{-2,2}(\Sigma_{\epsilon, \lambda_0}), 
\end{align*}
and $\dfrac{\lambda^{1/2}}{\mu^2 B^3} \in \bM_{-2,1}
(\Sigma_{\epsilon, \lambda_0})$, by Lemma \ref{lem:solop} we have
Theorem \ref{lem:sol1}. \qed \vskip0.5pc 
We next consider the equations:
\begin{equation}\label{4.20}\left\{\begin{aligned}
\lambda \bw - \DV(\mu\bD(\bw) - \fq\bI) =0,
\quad \dv \bw & = 0
&\quad&\text{in $\BR^N_+$}, \\
\lambda h + A_\kappa\cdot \nabla'h - \bw\cdot\bn_0 &= d
&\quad&\text{on $\BR^N_0$}, \\
(\mu\bD(\bw) - \fq\bI)\bn_0 -\sigma (\Delta'h)\bn_0  &= 0
&\quad&\text{on $\BR^N_0$}.
\end{aligned}\right.\end{equation}
We shall prove the following theorem.
\begin{thm}\label{lem:sol2}
Let $1 < q < \infty$ and $\epsilon \in (0, \pi/2)$. 
Then, there exist a $\lambda_1 > 0$ and solution operators
 $\CW(\lambda)$
and $\CH_\kappa(\lambda)$ with
\begin{align*}
\CW(\lambda) &\in \Hol(\Lambda_{\kappa, \lambda_1},
\CL(H^2_q(\BR^N_+), H^2_q(\BR^N_+)^N)),\\
\CH_\kappa(\lambda)&
\in \Hol(\Lambda_{\kappa, \lambda_1}, \CL(H^2_q(\BR^N_+), H^3_q(\BR^N_+))), 
\end{align*}
such that for any $\lambda = \gamma + i\tau
\in \Lambda_{\kappa, \lambda_1}$ and 
$d \in H^2_q(\BR^N_+)$, 
$\bw = \CW(\lambda)d$ and $h = \CH_\kappa(\lambda)d$
 are unique solutions of Eq. \eqref{4.20}
with some $\fq \in \hat H^1_q(\Omega)$, and 
\begin{align*}
\CR_{\CL(H^2_q(\BR^N_+), H^{2-k}_q(\BR^N_+)^N)}
(\{(\tau\pd_\tau)^\ell(\lambda^{k/2}\CW(\lambda))\mid
\lambda \in \Lambda_{\kappa, \lambda_1}\}) &\leq r_b(\lambda_1), \\
\CR_{\CL(H^2_q(\BR^N_+), H^{3-m}_q(\BR^N_+))}
(\{(\tau\pd_\tau)^\ell(\lambda^m\CH_\kappa(\lambda))\mid
\lambda \in \Lambda_{\kappa, \lambda_1}\}) &\leq r_b(\lambda_1)
\end{align*}
for $\ell=0,1$, $k=0,1,2$, and $m=0,1$, where $r_b(\lambda_1)$ is a
constant depending on $m_0$, $m_1$, $m_2$,
$\epsilon$, $\lambda_1$, $N$, and $q$.
\end{thm}
\pf We start with solution formulas. 
Applying the partial Fourier transform to Eq. \eqref{4.20},
 we have the following generalized resolvent problem: 
\allowdisplaybreaks{
\begin{alignat}2
\lambda\hat w_j + \mu |\xi'|^2\hat w_j - \mu\pd_N^2\hat w_j + i\xi_j \hat \fq
 = 0& &\quad&(x_N>0),\nonumber \\
\lambda\hat w_N + \mu |\xi'|^2\hat w_N - \mu\pd_N^2\hat w_N + \pd_N\hat \fq
 = 0& &\quad&(x_N>0),\nonumber \\
\sum_{j=1}^{N-1}i\xi_j\hat w_j + \pd_N \hat w_N  = 0& &\quad&(x_N > 0),
\nonumber \\
\mu(\pd_N\hat w_j(0) + i\xi_j\hat w_N(0))  = 0,\quad
2\mu\pd_N\hat w_N - \fq  = -\sigma A^2\hat h& &\quad&\text{for $x_N=0$}, 
\nonumber  \\
\lambda \hat h +\sum_{j=1}^{N-1}i\xi_jA_{\kappa j}\hat h + \hat w_N = \hat d&
&\quad&\text{for $x_N=0$}.
\label{4.21}
\end{alignat}
}
Where, we have set $A_\kappa = (A_{\kappa 1}, \ldots, A_{\kappa N-1})$. 
Using the solution formulas given in \eqref{12*} and \eqref{15} with 
$g_j = 0$ ($j=1,\ldots, N-1$), and $g_N= \sigma A^2 \hat h$, 
we have
\begin{equation}\label{4.22}\begin{split}
\hat w_j & = 
 \frac{i\xi_j\CM(x_N)}{\mu D(A,B)}
\sigma A^2(A^2+B^2)\hat h
 - \frac{i\xi_j e^{-Bx_N}}{\mu D(A,B)} 
\sigma A^2(B-A)\hat h, \\
\hat w_N & = -\frac{A \CM(x_N)}{\mu D(A,B)}\sigma A^2(A^2+B^2)\hat h 
 + \frac{e^{-Bx_N}}{\mu D(A, B)}\sigma A^3(A+B)\hat h, \\
\hat \fq & = - \frac{(A+B)A^2(A^2+B^2)e^{-Ax_N}}{D(A,B)}\hat h
\end{split}\end{equation}
Inserting the formula of $\hat w_N|_{x_N=0}$
 into the last equation in \eqref{4.21},
we have
$$(\lambda + i\xi'\cdot A_\kappa)\hat h 
+ \frac{\sigma A^3(A+B)}{\mu D(A,B)}\hat h = \hat d,
$$
where we have set $i\xi'\cdot A_\kappa = \sum_{j=1}^{N-1}i\xi_jA_{\kappa j}$, 
which implies that
\begin{equation}\label{4.23}
\hat h = \frac{\mu D(A,B)}{E_\kappa} \hat d
\end{equation}
with $E_\kappa=\mu(\lambda + i\xi'\cdot A_\kappa)D(A, B) + \sigma A^3(A+B)$.
Thus, we have the following solution formulas:
\begin{equation}\label{4.24}\begin{split}
\hat w_j & = i\xi_j\CM(x_N)\frac{\sigma A^2(A^2+B^2)}{E_\kappa}\hat d
-i\xi_je^{-Bx_N}\frac{\sigma A^2(B-A)}{E_\kappa}\hat d, \\
\hat w_N & = -A\CM(x_N)\frac{\sigma A^2(A^2+B^2)}{E_\kappa}\hat d
+ e^{-Bx_N}\frac{\sigma A^3(A+B)}{E_\kappa}\hat d, \\
\hat q & = -\frac{\mu(A+B)A^2(A^2+B^2)e^{-Ax_N}}{E_\kappa}\hat d.
\end{split}\end{equation}
\vskip0.5pc\noindent
Concerning the estimation for $E_\kappa$, we have the following lemma.
\begin{lem}\label{lem:4.3.1} \thetag1 ~ Let $0 < \epsilon < \pi/2$ and 
let $E_0$ be the function defined 
in \eqref{4.23} with $A_0=0$.  
Then, there exists a $\lambda_1 > 0$ and $c>0$ 
such that the estimate: 
\begin{equation}\label{eq:4.3.6}
|E_0| \geq c(|\lambda|+A)(|\lambda|^{1/2} + A)^3
\end{equation}
holds for $(\lambda, \xi') \in \Sigma_{\epsilon, \lambda_1}
\times (\BR^{N-1}\setminus\{0\})$.
\vskip0.5pc\noindent
\thetag2 Let $\kappa \in (0, 1)$ and let $E_\kappa$ be the function defined
in \eqref{4.23}.  Then, there exists a  
$\lambda_1 > 0$ and $c>0$ such that 
\begin{equation}\label{eq:4.3.7}
|E_\kappa| \geq c(|\lambda|+A)(|\lambda|^{1/2} + A)^3
\end{equation}
holds for $(\lambda, \xi') \in \BC_{+, \lambda_1}
\times (\BR^{N-1}\setminus\{0\})$.

Where, the constant $c$ in \thetag1 and \thetag2 depends on $\lambda_1$,
$m_0$, $m_1$, and $m_2$.
\end{lem}
\pf
We first study the case where $|\lambda| \geq R_1A$ for large $R_1 > 0$.
Note that  $|\lambda| \geq \lambda_1$. 
Since $|B| \leq A + \mu^{-1/2}|\lambda|^{1/2}$
and since $\Lambda_{\kappa, \lambda_1} \subset \Sigma_\epsilon$, 
by Lemma \ref{lem:4.3} we have 
\begin{align*}
|E_\kappa| 
&\geq \mu|\lambda||D(A,B)| - \mu|A_\kappa||A||D(A,B)|
- \sigma A^3(A+\mu^{-1/2}|\lambda|^{1/2}) \\
& \geq c\mu|\lambda|(|\lambda|^{1/2}+A)^3 - 
\mu m_2 CR_1^{-1}|\lambda|(|\lambda|^{1/2} + A)^3\\
&\quad -\sigma R^{-1}_1|\lambda|(|\lambda|^{1/2}+ A)^3
-\mu^{-1/2}\sigma|\lambda|^{1/2}(|\lambda|^{1/2}+ A)^3\\
&\geq(c\mu/2)|\lambda|(|\lambda|^{1/2}+A)^3
+ ((c\mu/2) - \mu m_2C R_1^{-1} \\
&\quad - \sigma R_1^{-1} 
- \sigma/(\mu|\lambda|)^{1/2})|\lambda|(|\lambda|^{1/2}+A)^3.
\end{align*}
Thus, choosing $R_1 > 0$ and $\lambda_1 > 0$ so large that
$(c\mu/4) - \mu m_2 CR_1^{-1}-\sigma R_1^{-1} \geq 0$ and 
$(c\mu/4) - \sigma/(\mu\lambda_1)^{1/2}
\geq 0$, we have
\begin{equation}\label{4.25}
|\tilde E_\kappa| \geq (c\mu/2)|\lambda|(|\lambda|^{1/2} + A)^3
\geq (c\mu/4)(|\lambda| + R_1A)(|\lambda|^{1/2} + A)^3
\end{equation}
provided that $|\lambda| \geq R_1A$ and 
$\lambda \in \Lambda_{\kappa, \lambda_1}$. When $\kappa=0$,
we may assume that 
$m_2=0$ above.

We now consider the case where $|\lambda| \leq R_1A$.
We assume that
$\lambda 
\in \Sigma_{\epsilon, \lambda_1}$.  In this case, 
we have $A \geq R_1^{-1}|\lambda|^{1/2}\lambda_1^{1/2}$, and so,
 choosing $\lambda_1$ large enough, 
we have $B = A(1+O(\lambda_1)^{-1/2}))$.  In particular, 
$D(A, B) = 4A^3(1+O(\lambda_1^{-1/2}))$. Thus, we have 
$$E_\kappa = 4\mu(\lambda + i\xi'\cdot A_\kappa)
A^3(1 + O(\lambda_1^{-1/2})) + 2\sigma A^4(1+O(\lambda_1^{-1/2})).
$$
We first consider the case where $\kappa=0$.  
Using Lemma \ref{lem:fund1}, we have
\begin{align*}
|E_0| &\geq |4\mu \lambda A^3 + 2\sigma A^4| - 4\mu|\lambda|A^3O(\lambda_1^{-1/2})
- 2\sigma A^4O(\lambda_1^{-1/2}) \\
& \geq (\sin \epsilon)(4\mu|\lambda|A^3 + 2\sigma A^4)
- O(\lambda_1^{-1/2})(4\mu|\lambda|A^3 + 2\sigma A^4).
\end{align*}
Thus, choosing $\lambda_1 > 0$ so large that $(\sin\epsilon/2) 
- O(\lambda_1^{-1/2}) \geq 0$,
we have
\begin{align*}
|E_0| &\geq (\sin\epsilon/2)(4\mu|\lambda|A^3 + 2\sigma A^4)
\geq c(|\lambda| + A)A^3 \\
&\geq c/2^3(|\lambda|+A)(A + R_1^{-1}\lambda_1^{1/2}|\lambda|^{1/2})^3.
\end{align*}
This completes the proof of \thetag1.  

We next consider the case of $\kappa \in (0, 1)$. 
 Taking the real part gives
\begin{align*}
{\rm Re}\, E_\kappa = 4\mu({\rm Re}\,\lambda) A^3(1+O(\lambda_1^{-1/2}))
+ &O(\lambda_1^{-1/2})({\rm Im}\,\lambda + A_\kappa\cdot\xi')A^3\\
+ 2\sigma A^4(1+&O(\lambda_1^{-1/2})).
\end{align*}
Since ${\rm Re}\,\lambda \geq \lambda_1 > 0$ and 
$|\lambda| \leq R_1A$, we have
$$
{\rm Re}\, E_\kappa 
\geq 2\sigma A^4 - (4\mu(m_2+R_1)+2\sigma)O(\lambda_1^{-1/2})A^4,
$$
and so, 
choosing $\lambda_1 > 0$ so large that 
$\sigma - (4\mu(m_2+R_1)+2\sigma)O(\lambda_1^{-1/2}) \geq 0$,
we have
$$|E_\kappa| \geq {\rm Re}\, E_\kappa  
\geq \sigma A^4 \geq (\sigma/2^4)(A+R_1^{-1}|\lambda|)
(A + R_1^{-1}\lambda_1^{1/2}|\lambda|^{1/2})^3.
$$
This completes the proof of Lemma \ref{lem:4.3.1}. 
\qed
\vskip0.5pc\noindent
{\bf Continuation of Proof of Theorem \ref{lem:sol2}}.~ 
Let $w_j = \CF^{-1}_{\xi'}[\hat w_j]$, 
$\fq = \CF^{-1}_{\xi'}[\hat \fq]$ and 
$\eta = \varphi(x_N)\CF^{-1}_{\xi'}[e^{-Ax_N} \hat h]$,
where $\varphi(x_N) \in C^\infty_0(\BR)$ equals to
$1$ for $x_N \in (-1, 1)$ and $0$ for $x_N \not\in
[-2, 2]$.   Notice that
$\eta|_{x_N=0} = h$. 

Let $w_j(x) = \CF^{-1}_{\xi'}[\hat w_j(\xi', x_N)](x')$.
In view of \eqref{4.24} and Volevich's trick, we define $\CW_j(\lambda)$
by 
\begin{align*}
\CW_j(\lambda)d & = \int^\infty_0\CF^{-1}_{\xi'}\Bigl[-(Ae^{-B(x_N+y_N)}
+ A^2\CM(x_N+y_N))\\
&\phantom{= \int^\infty_0\CF^{-1}_{\xi'}\Bigl[-(Ae^{-B(x_N+y_N)}}
\frac{i\xi_j}{A}\frac{\sigma(A^2+B^2)}{E_\kappa}\CF[\Delta'd](\xi', y_N)\\
&\quad 
+ Ae^{-B(x_N+y_N)}\frac{i\xi_j}{A}\frac{\sigma B(B-A)}{E_\kappa}
\CF'[\Delta'd](\xi', y_N)\Bigr](x')\,dy_N\\
&+\int^\infty_0\CF^{-1}_{\xi'}\Bigl[-A^2\CM(x_N+y_N)
\frac{\sigma(A^2+B^2)}{E_\kappa}\CF'[\pd_j\pd_Nd](\xi', y_N)\\
&\quad+Ae^{-B(x_N+y_N)}\frac{\sigma A(B-A)}{E_\kappa}\CF'[\pd_j\pd_Nd]
(\xi', y_N)\Bigr](x')\,dy_N,
\end{align*}
where we have used $\CF'[\Delta'd](\xi', y_N) = -A^2\hat d(\xi', y_N)$.
We have $\CW_j(\lambda)d = w_j$. 
By  Lemma \ref{lem:4.3.1}, we see that 
$$\frac{A^2+B^2}{E_\kappa}, 
\quad \frac{A^2+B^2}{E_\kappa}\frac{\xi_j}{A}, \quad
\frac{A(B-A)}{E_\kappa}, \quad \frac{B(B-A)}{E_\kappa}\frac{\xi_j}{A}
$$
belong to $\bM_{-2,2}(\Lambda_{\kappa, \lambda_1})$, and so by 
Lemma \ref{lem:4.2}, we have 
$$\CR_{\CL(H^2_q(\BR^N_+), H^{2-k}_q(\BR^N_+))}
(\{(\tau\pd_\tau)^\ell(\lambda^{k/2}\CW_j(\lambda))
\mid\lambda \in \Lambda_{\kappa, \lambda_1}\}) \leq r_b(\lambda_1)
$$
for $\ell=0,1$ and $k=0,1,2$, where $r_b(\lambda_1)$
is a constant depending on $m_0$,
$m_1$, $m_2$ and $\lambda_1$.

Analogously, we have 
$$\CR_{\CL(H^2_q(\BR^N_+), H^{2-k}_q(\BR^N_+))}
(\{(\tau\pd_\tau)^\ell(\lambda^{k/2}\CW_N(\lambda))
\mid\lambda \in \Lambda_{\kappa, \lambda_1}\}) \leq r_b(\lambda_1)
$$
for $\ell=0,1$ and $k=0,1,2$. 
Thus, our final task is to construct $\CH_\kappa(\lambda)$. 
In view of \eqref{4.23}, we define 
$\CH_\kappa(\lambda)$ acting on $d \in H^2_q(\BR^N_+)$ by 
$$\CH_\kappa(\lambda)d 
= \varphi(x_N)\CF^{-1}_{\xi'}\Bigl[e^{-Ax_N}\frac{\mu D(A,B)}{E_\kappa}
\hat d(\xi',0)\Bigr](x').$$
Since $\varphi(x_N)$ equals one for $x_N \in (-1, 1)$, we have
$\CH_\kappa(\lambda)d|_{x_N=0} = h$.  Recalling the definition
of $\hat h$ given in \eqref{4.23} and using Volevich's trick,
we have $\CH_\kappa(\lambda)d = 
\varphi(x_N)\{\Omega_\kappa(\lambda)d + \CH_\kappa^2(\lambda)d\}$
with 
\begin{align*}
\Omega_\kappa(\lambda)d 
& = \int^\infty_0\CF^{-1}_{\xi'}\Bigl[Ae^{-A(x_N+y_N)}
\frac{\mu D(A,B)}{ E_\kappa}\varphi(y_N)
\hat d(\xi', y_N)\Bigr](x')\, dy_N, \\
\CH_\kappa^2(\lambda)d 
& = -\int^\infty_0\CF^{-1}_{\xi'}\Bigl[e^{-A(x_N+y_N)}
\frac{\mu D(A,B)}{ E_\kappa}
\pd_N(\varphi(y_N)\hat d(\xi', y_N))\Bigr](x')\, dy_N.
\end{align*}
We use the following lemma. 
\begin{lem}\label{lem:4.3.2} Let $\Lambda$ be a domain in
$\BC$ and let $1 < q < \infty$.  Let
$\varphi$ and $\psi$ be two $C^\infty_0((-2, 2))$ functions. 
Given $m \in \bM_{0, 2}(\Lambda)$, we define  operators
$L_6(\lambda)$ and $L_7(\lambda)$ acting on $g \in L_q(\BR^N_+)$ by 
\begin{align*}
[L_6(\lambda)g](x) &=\varphi(x_N)\int^\infty_0\CF^{-1}_{\xi'}
\Bigl[e^{-A(x_N+y_N)}m(\lambda, \xi')\hat g(\xi', y_N)\psi(y_N)\Bigr]\,
dy_N, \\
[L_7(\lambda)g](x) &=\varphi(x_N)\int^\infty_0\CF^{-1}_{\xi'}
\Bigl[Ae^{-A(x_N+y_N)}m(\lambda, \xi')\hat g(\xi', y_N)\psi(y_N)\Bigr]\,
dy_N.
\end{align*}
Then, 
\begin{equation}\label{4.27}
\CR_{\CL(L_q(\BR^N_+)}(\{(\tau\pd_\tau)^\ell L_k(\lambda) \mid
\lambda \in \Lambda\}) \leq r_b
\end{equation}
for $\ell=0,1$ and $k=6,7$, where 
 $r_b$ is a constant depending on 
$M(m, \Lambda)$.  Here, $M(m, \Lambda)$ is the number
defined in Definition \ref{def:4.1}.   
\end{lem}
\pf Using the assertion for $L_2(\lambda)$ in 
 Lemma \ref{lem:4.2}, we can show 
\eqref{4.27} immediately for $k=7$, and so we show \eqref{4.27}
only in the case that $k=6$ below. 
In view of Definition \ref{def:2}, for any $n \in \BN$, we take
$\{\lambda_j\}_{j=1}^n \subset \Lambda$, $\{g_j\}_{j=1}^n
\subset L_q(\BR^N_+)$, and $r_j(u)$ ($j=1, \ldots, n$) are  
Rademacher functions. For the notational simplicity,  we set
\begin{align*}
|||L_6(\lambda)g||| &= \|\sum_{j=1}^nr_j(u)L_6(\lambda_j)g_j\|_{L_q((0,1),
L_q(\BR^N))}\\ 
&= \Bigl(\int^1_0\|\sum_{j=1}^nr_j(u)L_6(\lambda_j)g_j\|_{L_q(\BR^N_+)}^q\,du
\Bigr)^{1/q}.
\end{align*}
By the Fubini-Tonelli theorem, we have
\begin{align*}
|||L_6(\lambda)g|||^q
&= \int^1_0\int^\infty_0\int_{\BR^{N-1}}
|\sum_{j=1}^nr_j(u)L_6(\lambda_j)g_j|^q\,dy'dx_N\,du\\
&= \int^\infty_0\Bigl(\int^1_0\|\sum_{j=1}^nr_j(u)L_6(\lambda_j)
g_j\|_{L_q(\BR^{N-1})}^q\,du\Bigr)\,dx_N.
\end{align*}
Since
$$|\pd_{\xi'}^{\alpha'}(e^{-A(x_N+y_N)}m_0(\lambda, \xi'))|
\leq C_{\alpha'}|\xi'|^{-|\alpha'|}$$
for any $x_N \geq0$, $y_N \geq 0$, $(\lambda, \xi')
\in \Lambda\times(\BR^{N-1}\setminus\{0\})$, and 
$\alpha' \in \BN^{N-1}$, by Theorem \ref{thm:4.1} 
we have\allowdisplaybreaks{
\begin{align}
\int^1_0\|\sum_{j=1}^nr_j(u)&\CF_{\xi'}^{-1}
\Bigl[e^{-A(x_N+y_N)}m(\lambda_j, \xi')\hat g_j(\xi', y_N)
\Bigr](y')\|_{L_q(\BR^{N-1})}^q\,du \nonumber \\
&\leq C_{N,q}M(m, \Lambda)\int^1_0\|\sum_{j=1}^nr_j(u)g_j(\cdot, y_N)
\|_{L_q(\BR^{N-1})}^q\,du.\label{4.28}
\end{align}
For any $x_N \geq 0$, by Minkowski's integral inequality, 
Lemma \ref{lem:fund2}, and H\"older's inequality, we have
\allowdisplaybreaks{
\begin{align*}
&\Bigl(\int^1_0\|\sum_{j=1}^nr_j(u)L_6(\lambda_j)g_j
\|_{L_q(\BR^{N-1})}^q\,du\Bigr)^{1/q}\\
& =|\varphi(x_N)|\Bigl(\int^1_0\|
\int^\infty_0\CF^{-1}_{\xi'}[\sum_{j=1}^nr_j(u)e^{-A(x_N+y_N)}
m(\lambda_j, \xi')\hat g_j(\xi', y_N)](y')\\
&\phantom{=|\varphi(x_N)|\Bigl(\int^1_0\|
\int^\infty_0\CF^{-1}_{\xi'}}
\psi(y_N)\,dy_N\|_{L_q(\BR^{N-1})}^q\,du
\Bigr)^{1/q}\\
&\leq
|\varphi(x_N)|\Bigl(\int^1_0
\int^\infty_0\|\CF^{-1}_{\xi'}[\sum_{j=1}^nr_j(u)e^{-A(x_N+y_N)}
m(\lambda_j, \xi')\hat g_j(\xi', y_N)](y')\\
&\phantom{=|\varphi(x_N)|\Bigl(\int^1_0\|
\int^\infty_0\CF^{-1}_{\xi'}}
\psi(y_N)\|_{L_q(\BR^{N-1})}\,dy_N\Bigr)^q\,du
\Bigr)^{1/q}\\
&\leq 
|\varphi(x_N)|\int^\infty_0\Bigl(
\int^1_0\|\CF^{-1}_{\xi'}[\sum_{j=1}^nr_j(u)e^{-A(x_N+y_N)}
m(\lambda_j, \xi')\\
&\phantom{=|\varphi(x_N)|\Bigl(\int^1_0\|
\int^\infty_0\CF^{-1}_{\xi'}}
\hat g_j(\xi', y_N)](y')\|_{L_q(\BR^{N-1})}^q\,du
\Bigr)^{1/q}|\psi(y_N)|\,dy_N
\\
&\leq C_{N,q}M(m, \Lambda)|\varphi(x_N)|
\int^\infty_0\Bigl(\int^1_0\|\sum_{j=1}^nr_j(u)g_j(\cdot, y_N)
\|_{L_q(\BR^{N-1})}^q\,du\Bigr)^{1/q}\\
&\phantom{aaaaaaaaaaaaaaaaaaaaaaaaaaaaaaaaaaa}
\times|\psi(y_N)|\,dy_N\\
&\leq C_{N,q}M(m, \Lambda)|\varphi(x_N)|\Bigl(\int^\infty_0\int^1_0
\|\sum_{j=1}^nr_j(u)g_j(\cdot, y_N)\|_{L_q(\BR^{N-1})}^q\,du\,dy_N
\Bigr)^{1/q}\\
&\phantom{aaaaaaaaaaaaaaaaaaaaaaaaaaaaaaaaaaa}
\times
\Bigl(\int^\infty_0|\psi(y_N)|^{q'}\,dy_N\Bigr)^{1/{q'}}\\
& = C_{n,q}M(m, \Lambda)|\varphi(x_N)|\Bigl(\int^1_0
\|\sum_{j=1}^nr_j(u)g_j(\cdot, y_N)\|_{L_q(\BR^N_+)}^q\,du\Bigr)^{1/q}
\\
&\phantom{aaaaaaaaaaaaaaaaaaaaaaaaaaaaaaaaaaa}
\times
\Bigl(\int^\infty_0|\psi(y_N)|^{q'}\,dy_N\Bigr)^{1/{q'}}
\end{align*}
Putting these inequalities together and using H\"older's inequality
gives
\begin{align*}
&\int^1_0\|\sum_{j=1}^nr_j(u)L_6(\lambda_j)g_j\|_{L_q(\BR^N_+)}^q\,du\\
&\quad\leq(C_{n,q}M(m, \Lambda))^q \int^\infty_0|\varphi(x_N)|^qdx_N
\int^1_0\|\sum_{j=1}^n
r_j(u)g_j\|_{L_q(\BR^N_+)}^q\,du\,\\
&\phantom{\quad\leq(C_{n,q}M(m, \Lambda))^q \int^\infty_0|\varphi(x_N)|^qdx_N
\int^1_0\|\sum_{j=1}^n}
\times\Bigl(
\int^\infty_0|\psi(y_N)|^{q'}\,dy_N\Bigr)^{q/{q'}},
\end{align*}
}
and so, we have 
\begin{align*}
&\|\sum_{j=1}^nr_jL_6(\lambda_j)g_j\|_{L_q((0, 1), L_q(\BR^N_+))}\\
&\quad
\leq C_{n,q}M(m, \Lambda)\|\varphi\|_{L_q(\BR)}\|\psi\|_{L_{q'}(\BR)}
\|\sum_{j=1}^nr_jg_j\|_{L_q((0, 1), L_q(\BR^N_+))}.
\end{align*}
This shows Lemma \ref{lem:4.3.2}. 
\qed
\vskip0.5pc\noindent
{\bf Continuation of Proof of Theorem \ref{lem:sol2}}.~For $(j, \alpha', k)
\in \BN_0\times\BN_0^{N-1}\times \BN_0$ with $j+|\alpha'| + k \leq 3$
and $j=0,1$, we write
$$\lambda^j\pd_{x'}^{\alpha'}\pd_N^k\CH_\kappa(\lambda)d
= \sum_{n=0}^k{}_kC_n(\pd_N^{k-n}\varphi(x_N))
[\lambda^j\pd_{x'}^{\alpha'}\pd_N^n\Omega_\kappa(\lambda)d
+ \lambda^j\pd_{x'}^{\alpha'}\pd_N^n\CH^2_\kappa(\lambda)d],
$$
and then 
\begin{align*}
&\lambda^j\pd_{x'}^{\alpha'}\pd_N^n\Omega_\kappa(\lambda)d  \\
&\quad=\int^\infty_0\CF^{-1}_{\xi'}\Bigl[Ae^{-A(x_N+y_N)}
\frac{\mu\lambda^j(i\xi')^{\alpha'}(-A)^nD(A,B)}
{(1+A^2)E_\kappa}\varphi(y_N)\\
&\phantom{\quad=\int^\infty_0\CF^{-1}_{\xi'}\Bigl[Ae^{A(x_N+y_N)}}
\phantom{aaaaaaaaaa}
\CF'[(1-\Delta')d](\xi', y_N)\Bigr](x')\,dy_N; \\
&\lambda^j\CH^2_\kappa(\lambda)d 
= \int^\infty_0\CF^{-1}_{\xi'}\Bigl[
e^{-A(x_N+y_N)}\frac{\mu\lambda^jD(A,B)}{E_\kappa}\\
&\phantom{\quad=\int^\infty_0\CF^{-1}_{\xi'}\Bigl[Ae^{-A(x_N+y_N)}}
\phantom{aaaaaaaaaa}
\pd_N(\varphi(y_N)\hat d(\xi', y_N))\Bigr](x')\,dy_N;
\\
&\lambda^j\pd_{x'}^{\alpha'}\pd_N^n\CH^2_\kappa(\lambda)d  \\
&\quad=\int^\infty_0\CF^{-1}_{\xi'}\Bigl[e^{-A(x_N+y_N)}
\frac{\mu\lambda^j(i\xi')^{\alpha'}(-A)^nD(A,B)}
{(1+A^2)E_\kappa}\\
&\phantom{\quad=\int^\infty_0\CF^{-1}_{\xi'}\Bigl[Ae^{-A(x_N+y_N)}}
\phantom{aaaaaaaaaa}
\pd_N(\varphi(y_N)\hat d(\xi', y_N))](x')\,dy_N \\
&\quad
-\sum_{k=1}^{N-1}\int^\infty_0\CF^{-1}_{\xi'}\Bigl[Ae^{-A(x_N+y_N)}
\frac{\mu\lambda^j(i\xi')^{\alpha'}(-A)^nD(A,B)}
{(1+A^2)E_\kappa}\frac{i\xi_j}{A}\\
&\phantom{\quad=\int^\infty_0\CF^{-1}_{\xi'}\Bigl[Ae^{-
A(x_N+y_N)}}
\phantom{aaaaaaaaaa}
\pd_N(\varphi(y_N)
\CF[\pd_jd(\cdot, y_N)](\xi')\}\Bigr](x')\,dy_N
\end{align*}
for $|\alpha'| + n \geq 1$. Where, we have used the formula:
$$1 = \frac{1+A^2}{1+A^2} = \frac{1}{1+A^2} - \sum_{j=1}^{N-1}\frac{A}{1+A^2}
\frac{i\xi_j}{A}
i\xi_j
$$
in the third equality. 
 By Lemma \ref{lem:4.3} and Lemma \ref{lem:4.3.1}, 
we see that multipliers:
\begin{align*}
&\frac{\lambda^j(i\xi')^{\alpha'}A^nD(A,B)}{(1+A^2)E_\kappa}, \quad
\frac{\lambda^j D(A,B)}{E_\kappa}, \\ 
&\frac{\lambda^j(i\xi')^{\alpha'}A^nD(A,B)}{(1+A^2)E_\kappa}, \enskip
\frac{\lambda^j(i\xi')^{\alpha'}A^nD(A,B)}{(1+A^2)E_\kappa}\frac{\xi_j}{A}
\end{align*}
belong to $\bM_{0,2}(\Lambda_{\kappa, \lambda_1})$, 
because $j+|\alpha'| + n \leq 3$ and 
$j=0,1$.  Thus, using Lemma \ref{lem:4.3.2}, we see that 
for any $n \in \BN$,
$\{\lambda_j\}_{j=1}^n \subset \Lambda_{\kappa, \lambda_1}$, 
and $\{d_j\}_{j=1}^n
\subset H^2_q(\BR^N_+)$, the inequality:
\begin{align*}
&\|\sum_{\ell=1}^nr_\ell(\cdot)(\pd_N^{k-n}\varphi)
(\lambda_\ell)^j\pd_{x'}^{\alpha'}\pd_N^n
(\Omega_\kappa(\lambda_\ell), \CH^2_\kappa(\lambda_\ell))d_\ell
\|_{L_q((0, 1), L_q(\BR^N))}\\
&\quad
\leq C\|\sum_{\ell=1}^nr_\ell(\cdot)d_\ell\|_{L_q((0, 1), H^2_q(\BR^N_+))}
\end{align*}
holds, which leads to
\begin{align*}
&\|\sum_{\ell=1}^nr_\ell(\cdot)(\lambda_\ell)^j
\pd_{x'}^{\alpha'}\pd_N^k\CH_\kappa(\lambda_\ell)d_\ell
\|_{L_q((0,1), L_q(\BR^N_+))}\\
&\quad
\leq C\|\sum_{\ell=1}^nr_\ell(u)d_\ell\|_{L_q((0, 1), H^2_q(\BR^N_+))}.
\end{align*}
Here, $C$ is a constant depending on $N$, $q$, $m_0$, $m_1$, and 
$m_2$. This shows that
$$\CR_{\CL(H^2_q(\BR^N_+), H^{3-k}_q(\BR^N_+))}
(\{\lambda^k\CH_\kappa(\lambda) \mid
\lambda \in \Lambda_{\kappa, \lambda_1}\}) 
\leq r_b$$
for $k=0,1$.  Here, $r_b(\lambda_1)$ is a constant depending on 
$N$, $q$, $m_0$, $m_1$, and $m_2$, but independent of
$\mu, \sigma \in [m_0, m_1]$ and $|A_\kappa| \leq m_2$
for $\kappa \in [0, 1)$.
 Analogously,
we have 
$$\CR_{\CL(H^2_q(\BR^N_+), H^{3-k}_q(\BR^N_+))}
(\{\tau\pd_\tau(\lambda^k\CH_\kappa(\lambda)) \mid
\lambda \in \Lambda_{\kappa, \lambda_1}\}) 
\leq r_b(\lambda_1)$$
for $k=0,1$. This completes the proof of Theorem \ref{lem:sol2}.
\qed 
\vskip0.5pc\noindent
{\bf Proof of Theorem \ref{main:half1}}.~To prove Theorem \ref{main:half1},
in view of the consideration in Subsec. \ref{subsec:3.3},
 we first  consider the equation:
\begin{equation}\label{4.31}
\dv\bv = g \quad\text{in $\HS$}.
\end{equation}
Where, $g$ is a solution of 
the variational equation:
\begin{equation}\label{4.30} \lambda(g, \varphi)_\HS
+ (\nabla g, \nabla\varphi)_\HS = (-\bff, \nabla\varphi)_\HS
\quad\text{for any $\varphi \in H^1_{q',0}(\HS)$},
\end{equation}
subject to $g = \rho$ on $\Gamma$.
We have the following theorem. 
\begin{lem}\label{lem:p.4.1} Let $1 < q < \infty$, $0 < \epsilon
< \pi/2$, and $\lambda_0 > 0$.  Let 
\begin{align*}
Y''_q(\HS) & =\{(\bff, \rho) \mid \bff \in L_q(\HS)^N,
\enskip \rho \in H^1_q(\HS)\}, \\
\CY''_q(\HS) & = \{(F_1, G_3, G_4) \mid F_1 \in L_q(\HS)^N,
\enskip G_3 \in L_q(\HS), \enskip G_4 \in H^1_q(\HS)\}.
\end{align*}
Let $g$ be a solution of the variational problem 
\eqref{4.30}.  Then, there exists an operator family $\CB_0(\lambda)
\in \Hol(\Sigma_{\epsilon, \lambda_0}, \CL(\CY''_q(\HS),
H^2_q(\HS)^N))$ such that for any $\lambda \in \Sigma$ and 
$(\bff, \rho) \in Y''_q(\HS)$, problem \eqref{4.31} admits a
solution $\bv = \CB_0(\lambda)(\bff, \lambda^{1/2}\rho, \rho)$, 
and
$$\CR_{\CL(\CY''_q(\HS), H^{2-j}_q(\HS)^N)}
(\{(\tau\pd_\tau)^\ell(\lambda^{j/2}\CB_0(\lambda)) \mid 
\lambda \in \Sigma_{\epsilon, \lambda_0}\}) \leq r_b(\lambda_0)
$$
for $\ell=0,1$ and $j=0,1,2$, where $r_b(\lambda_0)$ is a constant depending
on $\epsilon$, $\lambda_0$, $N$, and $q$.
\end{lem}
\pf This lemma was proved in Shibata
\cite[Lemma 9.3.10]{S4}, but for the sake of completeness of 
this lecture note as much as possible, 
we give a proof. Let $g_1$ be a solution 
of the equation:
$$(\lambda-\Delta)g_1 = \dv \bff \quad\text{in $\HS$}, \quad g_1|_{x_N=0}=0,
$$
and let $g_2$ be a solution of the equation:
$$(\lambda-\Delta)g_2 = 0 \quad\text{in $\HS$}, \quad g_2|_{x_N=0}=\rho.
$$
And then, $g = g_1 + g_2$ is a solution of Eq. \eqref{4.30}. To construct
$g_1$ and $g_2$, we use  the even extension, $f^e$, 
and odd extension, $f^o$,  of a function, $f$ introduced in
\eqref{4.ext.1}. Let $\bff = {}^\top(f_1, \ldots, f_N)$.  Notice that
$(\dv\bff)^o = \sum_{j=1}^{N-1}\pd_jf_j^o + \pd_Nf_N^e$.  We define 
$g_1$ by letting
$$g_1 = \CF^{-1}_\xi\Bigl[\frac{\CF[(\dv\bff)^o](\xi)}{\lambda+|\xi|^2}
\Bigr] = \CF^{-1}_\xi\Bigl[\frac{\sum_{k=1}^{N-1}i\xi_k\CF[f_k^o](\xi)
+i\xi_N\CF[f^e_N](\xi)}{\lambda+|\xi|^2}
\Bigr].
$$
And also, the $g_2$ is defined by
\begin{equation}\label{add:5.1}
g_2(x) = \CF^{-1}_{\xi'}[e^{-B_0x_N}\hat \rho(\xi', 0)](x') = 
\frac{\pd h}{\pd x_N},
\end{equation}
where we have set $B_0 = \sqrt{\lambda + |\xi'|^2}$ and 
$h(x) = -\CF^{-1}_{\xi'}[B_0^{-1}e^{-B_0x_N}\hat \rho(\xi', 0)](x')$.
By Volevich's trick, we have
\begin{align*}
h(x)& = \int^\infty_0\CF^{-1}_{\xi'}[B_0^{-1}e^{-B_0(x_N+y_N)}
\widehat{(\pd_N\rho)}(\xi', y_N)](x')\,dy_N \\
&-\int^\infty_0\CF^{-1}_{\xi'}[e^{-B_0(x_N+y_N)}\hat \rho(\xi', y_N)]
(x')\,dy_N \\
& = \int^\infty_0\CF^{-1}_{\xi'}[\frac{\lambda^{1/2}}{B_0^3}
\lambda^{1/2}e^{-B_0(x_N+y_N)}
\widehat {(\pd_N\rho)}(\xi', y_N)](x')\,dy_N \\
&+ \int^\infty_0\CF^{-1}_{\xi'}[\frac{A}{B_0^3}
Ae^{-B_0(x_N+y_N)}
\widehat {(\pd_N\rho)}(\xi', y_N)](x')\,dy_N \\
&-\int^\infty_0\CF^{-1}_{\xi'}[\frac{1}{B_0^2}
\lambda^{1/2}e^{-B_0(x_N+y_N)}\widehat{(\lambda^{1/2}\rho)}(\xi', y_N)]
(x')\,dy_N \\
&+\sum_{j=1}^{N-1}
\int^\infty_0\CF^{-1}_{\xi'}[\frac{1}{B_0^2}\frac{i\xi_j}{A}
Ae^{-B_0(x_N+y_N)}\widehat{(\pd_j\rho)}(\xi', y_N)]
(x')\,dy_N. 
\end{align*}
Let $\CZ_q(\HS)$ be the same space as in Lemma \ref{lem:solop}. 
We then define an operator $H(\lambda)$ 
acting on $(G_3, G_4) \in \CZ_q(\HS)$ by setting
\begin{align*}
H(\lambda)(G_3,G_4) &  
= \int^\infty_0\CF^{-1}_{\xi'}[\frac{\lambda^{1/2}}{B_0^3}
\lambda^{1/2}e^{-B_0(x_N+y_N)}
\widehat {\pd_NG_4}(\xi', y_N)](x')\,dy_N \\
&+ \int^\infty_0\CF^{-1}_{\xi'}[\frac{A}{B_0^3}
Ae^{-B_0(x_N+y_N)}
\widehat {\pd_NG_4}(\xi', y_N)](x')\,dy_N \\
&-\int^\infty_0\CF^{-1}_{\xi'}[\frac{1}{B_0^2}
\lambda^{1/2}e^{-B_0(x_N+y_N)}\widehat{G_3}(\xi', y_N)]
(x')\,dy_N \\
&+\sum_{j=1}^{N-1}
\int^\infty_0\CF^{-1}_{\xi'}[\frac{1}{B_0^2}\frac{i\xi_j}{A}
Ae^{-B_0(x_N+y_N)}\widehat{\pd_jG_4}(\xi', y_N)]
(x')\,dy_N. 
\end{align*}
By Lemma \ref{lem:4.1} and Lemma \ref{lem:4.2}, we see
that
\begin{equation}\label{add:5.2}\begin{aligned}
&H(\lambda) \in {\rm Hol}\,(\Sigma_\epsilon,
\CL(\CY'_q(\HS), H^2_q(\HS))), 
\quad  
h = H(\lambda)(\lambda^{1/2}\rho, \rho), \\
&
\CR_{\CL(\CY'_q(\HS), H^{2-j}_q(\HS))}(\{(\tau\pd_\tau)^\ell 
(\lambda^{j/2}H(\lambda)) \mid \lambda \in \Sigma_{\epsilon, \lambda_0} \})
\leq r_b(\lambda)
\end{aligned}\end{equation}
for $\ell=0,1$ and $j=0,1,2$. Moreover, we have
\begin{equation}\label{add:5.3}
(\lambda-\Delta)H(\lambda)(G_3, G_4) = 0 \quad\text{in $\HS$}.
\end{equation}

Let $\bv_1$ be an $N$ vector of functions defined by
$$\bv_1 = -\CF^{-1}_\xi\Bigl[\frac{i\xi\CF[g_1](\xi)}{|\xi|^2}\Bigr]
= -\CF^{-1}_\xi\Bigl[\frac{\xi(\sum_{k=1}^{N-1}\xi_k\CF[f^o_k](\xi)
+ \xi_N\CF[f^e_N](\xi))}{(\lambda+|\xi|^2)|\xi|^2}\Bigr].$$
We see that $\dv \bv_1 = g_1$ in $\BR^N_+$.  Moreover, 
by Lemma \ref{lem:fund1} and Lemma \ref{lem:fund2}, there exists an 
operator family $\CB^1_0(\lambda) \in {\rm Hol}\,(\Sigma_\epsilon, 
\CL(L_q(\BR^N_+), H^2_q(\BR^N_+)^N))$ such that $\bv_1 = \CB^1_0(\lambda)\bff$
and 
$$\CR_{\CL(L_q(\HS)^N, H^{2-j}_q(\HS)^N)}
(\{(\tau\pd_\tau)^\ell(\lambda^{j/2}\CB^1_0(\lambda)) \mid 
\lambda \in \Sigma_{\epsilon, \lambda_0}\}) \leq r_b(\lambda_0)
$$
for $\ell = 0, 1$, $j=0,1,2$, and $\lambda_0 > 0$, where 
$r_b(\lambda_0)$ is a constant depending on $\epsilon$, $\lambda_0$, 
$N$ and $q$. 

Let 
$$v_{2j} = \CF^{-1}_\xi\Bigl[\frac{\xi_j\xi_N\CF[h^e](\xi)}{|\xi|^2}\Bigl]
= -\CF^{-1}_\xi\Bigl[\frac{i\xi_j\CF[g_2^o](\xi)}{|\xi|^2}\Bigr]
\quad(j=1, \ldots, N)
$$
and let $\bv_2 = {}^\top(v_{21}, \ldots, v_{2N})$, and then by \eqref{add:5.1}
 we have 
$\dv \bv_2= g_2$ in $\HS$.  Thus, we define an operator 
$\CB^2_0(\lambda) = (\CB^2_{01}(\lambda), \ldots, \CB^2_{0N}(\lambda))$ 
acting on $(G_3, G_4) \in \CZ_q(\HS)$ by 
$$\CB^2_{0j}(\lambda)(G_3, G_4)
 = \CF^{-1}_\xi\Bigl[\frac{\xi_j\xi_N\CF[
H(\lambda)(G_3, G_4)^e](\xi)}{|\xi|^2}\Bigl].
$$ 
By \eqref{add:5.3}, we have $\bv_2= \CB^2_0(\lambda)(\lambda^{1/2}\rho, \rho)$.
Noting that $\pd_Nf^e = (\pd_Nf)^o$ and $\pd_k\pd_Nf^e
= (\pd_k\pd_Nf)^0$ ($k-1, \ldots, N-1$), we have 
\begin{align*}
\lambda\CB^2_{0j}(\lambda)(G_3, G_4)
&=\CF_\xi^{-1}\Bigl[\frac{\xi_j\xi_N\CF[\lambda H(\lambda)(G_3, G_4)^e](\xi)}
{|\xi|^2}\Bigr], \\
\lambda^{1/2}\nabla\CB^2_{0j}(\lambda)(G_3, G_4)
&=\CF_\xi^{-1}\Bigl[\frac{\xi_j\xi\CF[\lambda^{1/2}
(\pd_NH(\lambda)(G_3, G_4))^o](\xi)}
{|\xi|^2}\Bigr], \\
\pd_k\nabla\CB^2_{0j}(\lambda)(G_3, G_4)
&=\CF_\xi^{-1}\Bigl[\frac{\xi_j\xi\CF[\lambda^{1/2}
(\pd_k\pd_NH(\lambda)(G_3, G_4))^o](\xi)}
{|\xi|^2}\Bigr],
\end{align*} 
Moreover, since  $\pd^2_NH(\lambda)(\lambda)(G_3, G_4)
 = -(\lambda-\sum_{j=1}^{N-1}\pd_j^2)H(\lambda)(\lambda)(G_3, G_4)$
as follows from \eqref{add:5.3}, we have
$$
\pd_N^2\CB^2_0(\lambda)(G_3, G_4)
  = \CF^{-1}_\xi\Bigl[\frac{\xi_N^2}{|\xi|^2}
\CF[((\lambda -\Delta') H(\lambda)(G_3, G_4))^e](\xi)\Bigr]. 
$$
Thus, by \eqref{add:5.2} and the Fourier multiplier
theorem, we see that
$$\CR_{\CL(\CY'_q(\HS), H^{2-j}_q(\HS)^N)}
(\{(\tau\pd_\tau)^\ell(\lambda^{j/2}\CB^2_0(\lambda)) \mid
\lambda \in \Sigma_{\epsilon, \lambda_0}\}) \leq r_b(\lambda_0)
$$
for $\ell=0,1$, $j=0,1,2$, and $\lambda_0 > 0$, where $r_b(\lambda_0)$ is a 
constant depending on $\epsilon$, $\lambda_0$, $N$, and $q$. 
Since $\bv = \bv_1 + \bv_2$ is a solution of  Eq. \eqref{4.31}, 
setting $\CB_0(\lambda)(F_1, G_3, G_4) = \CB^1_0(\lambda)F_1
+ \CB^2_0(\lambda)(G_3, G_4)$, 
we see that $\CB_0(\lambda)$ is the 
required operator, which completes the proof of Lemma \ref{lem:p.4.1}.
\qed
\vskip0.5pc
We now prove Theorem \ref{main:half1}. 
Let $(\bff, d, \bh)
\in Y_q(\HS)$.  Let $g$ be a solution of 
Eq. \eqref{4.30} with $\rho=\bn_0\cdot\bh$, and let 
$\bu$, $\fq$ and $h$ be solutions of the equations:
\begin{equation}\label{eq:4.33}
\left\{\begin{aligned}
\lambda\bu - \DV(\mu\bD(\bu) - \fq\bI) = \bff,
\quad \dv\bu &= g = \dv\bg
&\quad &\text{in $\HS$}, \\
\lambda h + A_\kappa\cdot\nabla'_\Gamma h-\bu\cdot\bn_0
& = d &\quad &\text{on $\BR^N_0$}, \\
(\mu\bD(\bu) - \fq\bI-\sigma(\Delta'h)\bI)
\bn_0 & = \bh &\quad &\text{on $\BR^N_0$}.
\end{aligned}\right.
\end{equation}
Then, according to what pointed out
in Subsec. \ref{subsec:3.3}, $\bu$ and $h$ are solutions of Eq. \eqref{4.3}.  
Thus, we shall look for $\bu$, $\fq$ and $h$ below.
In view of \eqref{wd:5} and \eqref{wd:5*}, 
applying Lemma \ref{lem:p.4.1} with $\rho = \bn_0\cdot\bh$, 
we define $\bu_0$ by $\bu_0= \CB_0(\bff, \lambda^{1/2}\bn_0\cdot\bh, 
\bn_0\cdot\bh)$. Notice that 
$\dv\bu_0 = g = \dv \bg_0$ with 
$\bg = \lambda^{-1}(\bff + \nabla g)$.
We then look for $\bw_0$, $\fq$, and $h$ satisfying the equations:
\begin{equation}\label{4.33}\left\{\begin{aligned}
\lambda \bw_0 - \DV(\mu\bD(\bw_0) - \fq\bI) = \bff - \bff_0, 
\quad \dv \bw_0 = 0& &\quad&\text{in $\HS$}, \\
\lambda h + A_\kappa\cdot\nabla'h - \bw_0\cdot\bn_0
= d+d_0& &\quad&\text{on $\BR^N_0$}, \\
(\mu\bD(\bw_0) - \fq\bI)\bn_0-\sigma(\Delta'h)\bn_0=
\bh-\bh_0& &\quad&\text{on $\BR^N_0$},
\end{aligned}\right. \end{equation}
where we have set
$$\bff_0 = \lambda \bu_0 - \DV(\mu\bD(\bu_0)), \quad
d_0 = \bu_0\cdot\bn_0, 
\quad \bh_0 = \mu\bD(\bu_0)\bn_0.
$$
To solve Eq. \eqref{4.33}, we  first consider the equations:
\begin{equation}\label{4.35}\left\{\begin{aligned}
\lambda \bU_1 - \DV(\mu\bD(\bU_1)-P_1\bI) = \bF, \quad
\dv\bU_1 = 0& &\quad&\text{in $\HS$}, \\
\pd_N (\bU_1\cdot\bn_0) = 0, \quad P_1=0& &\quad&\text{on $\BR^N_0$}.
\end{aligned}\right.\end{equation}
 For $\bF = {}^\top(F_1, \ldots, F_N) \in L_q(\HS)^N$, let 
$\tilde \bF = {}^\top(F_1^e, \ldots, F_{N-1}^e, F_N^o)$. Let $\CB_1(\lambda)$
and $\CP_1(\lambda)$ be operators acting on $\bF \in L_q(\BR^N_+)^N$
defined by 
\begin{align*}
\CB_1(\lambda)\bF &= \CF^{-1}_{\xi}\Bigl[\frac{\CF[\tilde\bF](\xi) 
- \xi\xi\cdot\CF[\tilde\bF](\xi)|\xi|^{-2}}{\lambda+\mu|\xi|^2}\Bigr],
\\
\quad \CP_1(\lambda)\bF &= 
\CF^{-1}_{\xi}\Bigl[\frac{\xi\cdot\CF[\tilde\bF](\xi)}{|\xi|^2}
\Bigr].
\end{align*}
As was seen in Shibata and Shimizu \cite[p.587]{SS4} or \cite[Proof of
Theorem 4.3]{SS1}, $\bU_1 = \CB_1(\lambda)\bF$ and $P_1 = \CP_1(\lambda)\bF$
satisfy Eq. \eqref{4.35}.  Moreover, employing the same argument as in 
Sect. \ref{subsec:5.1}, by Lemma \ref{lem:fund1} and Lemma \ref{lem:fund2}, 
we see that 
\begin{align*}
\CB_1(\lambda) &\in \Hol(\Sigma_{\epsilon, \lambda_0}, 
\CL(L_q(\BR^N_+)^N, H^2_q(\HS)^N)), \\ 
\CP_1(\lambda) &\in \Hol(\Sigma_{\epsilon, \lambda_0}, \CL(L_q(\HS)^N, 
\hat H^1_q(\HS)))
\end{align*}
for any $\epsilon \in (0, \pi/2)$ and $\lambda_0 > 0$, and
moreover
\begin{align*}
\CR_{\CL(L_q(\HS)^N, H^{2-j}_q(\HS)^N)}
(\{(\tau\pd_\tau)^\ell(\lambda^{j/2}\CB_1(\lambda)) \mid 
\lambda \in \Sigma_{\epsilon, \lambda_0}\}) &\leq r_b(\lambda_0)
\end{align*}
for $\ell=0,1$ and $j=0,1,2$, where $r_b(\lambda_0)$ is a constant depending on
$\epsilon$, $\lambda_0$, $m_0$, and $m_1$. In particular, 
we set 
\begin{equation}\label{4.37} \bu_1 = \CB_1(\lambda)(\bff-\bff_0),
\quad \fq_1 = \CP_1(\lambda)(\bff-\bff_0).
\end{equation}
We now let $\bu = \bu_0+\bu_1+ \bU_2$ and $\fq = \fq_1 + P_2$, and 
then 
\begin{equation}\label{4.38.1}\left\{\begin{aligned}
\lambda \bU_2 - \DV(\mu\bD(\bU_2) - P_2\bI) = 0, 
\quad \dv \bU_2 = 0& &\quad&\text{in $\HS$}, \\
\lambda h + A_\kappa\cdot\nabla'h - \bU_2\cdot\bn_0
= d+d_2& &\quad&\text{on $\BR^N_0$}, \\
(\mu\bD(\bU_2) - P_2\bI)\bn_0-\sigma (\Delta'h)\bn_0=
\bh-\bh_2& &\quad&\text{on $\BR^N_0$},
\end{aligned}\right. \end{equation}
where we have set 
$$d_2 = \bn_0\cdot(\bu_0+\bu_1), \quad \bh_2= \mu\bD(\bu_0+\bu_1).
$$
Thus, for $\bH \in H^1_q(\HS)^N$ we consider the equations: 
\begin{equation}\label{4.38.2}\left\{\begin{aligned}
\lambda \bU_2 - \DV(\mu\bD(\bU_2) - P_2\bI) = 0, 
\quad \dv \bU_2 = 0& &\quad&\text{in $\HS$},  \\
(\mu\bD(\bU_2) - P_2\bI)\bn_0=
\bH& &\quad&\text{on $\BR^N_0$},
\end{aligned}\right. \end{equation}
and then by Theorem \ref{lem:sol1}, we see 
that 
$\bU_2 = \CV(\lambda)(\lambda^{1/2}\bH, \bH)$ 
is a unique solutions of Eq. \eqref{4.38.2} with
some $P_2 \in \hat H^1_q(\HS)$. In particular, we set
$\bu_2 = \CV(\lambda)(\lambda^{1/2}(\bh-\bh_2), (\bh-\bh_2))$.

We finally let  $\bu=\bu_0+\bu_1 + \bu_2 + \bu_3$ and  
$\fq = \fq_1 + \fq_2 + \fq_3$, and then $\bu_3$, $\fq_3$ and $h$ 
are solutions of the equations:
\begin{equation}\label{4.38.3}\left\{\begin{aligned}
\lambda \bu_3 - \DV(\mu\bD(\bu_3) - \fq_3\bI) = 0, 
\quad \dv \bu_3 = 0& &\quad&\text{in $\HS$}, \\
\lambda h + A_\kappa\cdot\nabla'h - \bu_3\cdot\bn_0
= d+d_3& &\quad&\text{on $\BR^N_0$}, \\
(\mu\bD(\bu_3) - \fq_3\bI)\bn_0-\sigma (\Delta'h)\bn_0=0& 
&\quad&\text{on $\BR^N_0$},
\end{aligned}\right. \end{equation}
where $d_3= \bn_0\cdot(\bu_0+\bu_1+\bu_2)$. 
Setting 
$\CW(\lambda) = {}^\top(\CW_1(\lambda), \ldots, \CW_N(\lambda))$, 
by Theorem \ref{lem:sol2}, 
we see that $\bu_3 = \CW(\lambda)(d+d_3)$ 
and $h = \CH_\kappa(\lambda)(d+d_3)$ are unique solutions of
Eq. \eqref{4.38.3} with some $\fq_3 \in \hat H^1_q(\HS)$.  
Since the composition of two 
$\CR$-bounded operators   is also $\CR$ bounded as follows
from Proposition \ref{prop:4.1}, we see easily that 
given $\epsilon \in (0, \pi/2)$, there exist $\lambda_1 > 0$
and operator families $\CA_0(\lambda)$ and $\CH_0(\lambda)$ satisfying
\eqref{r-est:0} 
such that $\bu = \CA_0(\lambda)(\bff, d, \lambda^{1/2}\bh, \bh)$
and $h=\CH_0(\lambda)(\bff, d, \lambda^{1/2}\bh, \bh)$ are unique 
solutions of Eq. \eqref{4.3}, and moreover the estimate
\eqref{r-est:1} holds. This completes the proof of Theorem \ref{main:half1}.

\subsection{Problem in a bent half space}\label{sec:5} 
Let $\Phi:\BR^N\to\BR^N$ : $x \to y=\Phi(x)$ be a bijection of
$C^1$ class and let $\Phi^{-1}$ be its inverse map.  We assume that 
$\nabla \Phi$ and $\nabla\Phi^{-1}$ have the forms: 
$\nabla\Phi = \CA+B(x)$ and $\nabla\Phi^{-1} = \CA_{-1} + B_{-1}(y)$,
where $\CA$ and $\CA_{-1}$ are $N\times N$ orthogonal matrices with 
constant coefficients and $B(x)$ and $B_{-1}(y)$ are matrices of 
functions in $C^2(\BR^N)$ such that 
\begin{equation}\label{assump:4}\begin{aligned}
&\|(B, B_{-1})\|_{L_\infty(\BR^N)} \leq M_1, \quad 
\|\nabla(B, B_{-1})\|_{L_\infty(\BR^N)} \leq C_A, \\
&\|\nabla^2(B, B_{-1})\|_{L_\infty(\BR^N)} \leq M_2.
\end{aligned}\end{equation}
Here, $C_A$ is a constant depending on constants 
$A$, $a_1$, $a_2$ appearing in Definition \ref{dfn:1}. 
We choose $M_1 > 0$ small enough and  $M_2$ large enough
 eventually, and so 
we may assume that $0 < M_1 \leq 1 \leq C_A \leq M_2$. Let
$\Omega_+ = \Phi(\BR^N_+)$ and $\Gamma_+ = \Phi(\BR^N_0)$.  
In the sequel,  $a_{ij}$ and $b_{ij}(x)$ denote the $(i, j)^{\rm th}$
element of $\CA_{-1}$ and $(B_{-1}\circ\Phi)(x)$.

Let $\bn_+$ be the unit outer normal to $\Gamma_+$.
Setting $\Phi^{-1}={}^\top(\Phi_{-1,1}, \ldots, \Phi_{-1.N})$, 
we see that $\Gamma_+$ is represented by $\Phi_{-,N}(y)=0$,
which yields that
\begin{equation}\label{normal:1*}
\bn_+(x) = -\frac{(\nabla \Phi_{-1,N})\circ\Phi(x)}
{|\nabla \Phi_{-1,N})\circ\Phi(x)|}
=-\frac{{}^\top(a_{N1}+b_{N1}(x), \ldots, a_{NN} 
+b_{NN}(x))} 
{(\sum_{j=1}^N(a_{Nj} + b_{Nj}(x))^2)^{1/2}}.
\end{equation}
Obviously, $\bn_+$ is defined on $\BR^N$ and 
$\bn_+$ denotes the unit outer normal to $\Gamma_+$ for $y=\Phi(x',0)
\in \Gamma_+$. 
By \eqref{assump:4}, writing  
\begin{equation}\label{res:0}\bn_+
=-{}^\top(a_{N1}, \ldots, a_{NN}) +\bb_+(x),
\end{equation}
we see that $\bb_+$ is an $N$-vector defined on $\BR^N$,
which satisfies the estimates:
\begin{equation}\label{res:1}
\|\bb_+\|_{L_\infty(\BR^N)} \leq C_NM_1, \enskip
\|\nabla\bb_+\|_{L_\infty(\BR^N)} \leq C_NC_A, \enskip 
\|\nabla^2\bb_+\|_{L_\infty(\BR^N)} \leq C_{M_2}.
\end{equation}
We next give the Laplace-Beltrami operator on $\Gamma_+$.  Let 
$$g_{+ij}(x) = \frac{\pd\Phi}{\pd x_i}(x)\cdot\frac{\pd\Phi}{\pd x_j}(x)=
\sum_{k=1}^N(a_{ik}+b_{ik}(x))(a_{jk} + b_{jk}(x))
= \delta_{ij} + \tilde g_{ij}(x)$$
with $\tilde g_{ij} = \sum_{k=1}^N(a_{ik}b_{jk}(x) + a_{jk}b_{ik}(x)
+b_{ik}(x)b_{jk}(x))$.  Since $\Gamma_+$ is given by $y_N = \Phi(x', 0)$, 
letting $G(x)$ be an $N\times N$ matrix whose $(i, j)^{\rm th}$ element are 
$g_{ij}(x)$, we see that $G(x',0)$ is  the 1st fundamental matrix of $\Gamma_+$. 
Let 
$g_+ : =\sqrt{\det G}$ and let $g^{ij}_+(x)$ 
denote the $(i, j)^{\rm th}$ component of the inverse matrix,
$G^{-1}$, of $G$.  By \eqref{assump:4}, 
we can write
$$g_+ = 1 + \tilde g_+, \quad g^{ij}_+(x) = \delta_{ij} 
+ \tilde g^{ij}_+(x)
$$
with 
\begin{equation}\label{res:3}\begin{aligned}
&\|(\tilde g_+, \tilde g_+^{ij})\|_{L_\infty(\BR^N)} \leq C_NM_1,\quad
\|\nabla(\tilde g_+, \tilde g^{ij}_+)\|_{L_\infty(\BR^N)}\leq C_NC_A,
\\
&\|\nabla^2(\tilde g_+, \tilde g^{ij}_+)\|_{L_\infty(\BR^N)}
\leq C_{M_2}.
\end{aligned}\end{equation}
The Laplace-Beltrami operator $\Delta_{\Gamma_+}$ is given by
\begin{equation}\label{res:4}\begin{aligned}
(\Delta_{\Gamma_+}f)(y)& = \sum_{i,j=1}^{N-1}\frac{1}{g_+(x',0)}
\frac{\pd}{\pd x_i}\{g_+(x',0)g^{ij}_+(x', 0)
\frac{\pd}{\pd x_j}f(\Phi(x',0))\}\\
&=\Delta'f(\Phi(x',0)) + \CD_+f
\end{aligned}\end{equation}
for $y=\Phi(x',0) \in \Gamma_+$. Where,  
$$(\CD_+ f)(y) = \sum_{i,j=1}^{N-1}\tilde g^{ij}(x)
\frac{\pd^2f\circ\Phi}{\pd x_i\pd x_j}(x)
+ \sum_{j=1}^{N-1}g^j(x)\frac{\pd f\circ\Phi}{\pd x_j}(x)
\quad\text{for $y = \Phi(x)$}$$
with
$$g^j(x) = \frac{1}{g_+(x)}
\sum_{i=1}^{N-1}\frac{\pd}{\pd x_i}(g_+(x)g^{ij}(x)).
$$
By \eqref{res:3}
\begin{equation}\label{res:2}
\|\CD_+f\|_{H^1_q(\HS)} \leq C_NM_1\|\nabla^3 f\|_{L_q(\HS)}
+ C_{M_2}\|f\|_{H^2_q(\HS)}.
\end{equation}

We now formulate  problem treated in this section. 
Let $y_0$ be any point of $\Gamma_+$ and 
let $d_0$ be a positive number such that
\begin{equation}\label{s5.1}\begin{split}
& |\mu(y)- \mu(x_0)|, 
|\sigma(y)- \sigma(y_0)| \leq m_1M_1,
\quad\text{for any $y \in \overline{\Omega_+} \cap B_{d_0}(y_0)$};
\\
& |A_\kappa(y)-A_\kappa(y_0)| \leq m_2M_1 \quad\text{for any $y \in 
\Gamma_+ \cap B_{d_0}(y_0)$}.
\end{split}\end{equation}
In addition, $\mu$, $\sigma$, and $A_\kappa$ satisfy the following 
conditions: 
\begin{equation}\label{s5.2}\begin{split}
&m_0\leq \mu(y), \sigma(y) \leq m_1,
\quad|\nabla\mu(y)|, |\nabla \sigma(y)| \leq m_1 
\quad\text{for any $y \in \overline{\Omega_+}$}, \\
&|A_\kappa(y)| \leq m_2 \quad\text{for any $y \in \Gamma_+$},
\quad 
\|A_\kappa\|_{W^{2-1/q}_r(\Omega_+)} \leq m_3\kappa^{-b}
\end{split}\end{equation}
for any $\kappa \in (0, 1)$. 
In view of \eqref{assump:3} 
 and \eqref{s5.2}, to have \eqref{s5.1} for given $M_1 \in (0, 1)$
 it suffices to choose $d_0 > 0$ in such a way that 
$d_0 \leq M_1$ and $d_0^a \leq M_1$. 
We assume that $N < r < \infty$ according to \eqref{assump:3} and let 
 $A_0 = 0$.  Let $\varphi(y)$ be a function in $C^\infty_0(\BR^N)$
which equals $1$ for $y \in B_{d_0/2}(y_0)$ and $0$ in the outside of 
$B_{d_0}(y_0)$.  We assume that $\|\nabla\varphi\|_{H^1_\infty(\BR^N)}
\leq M_2$. Let 
\begin{align*}
\mu_{y_0}(y) &= \varphi(y)\mu(y) + (1-\varphi(y))\mu(y_0), \\
\sigma_{y_0}(y) &= \varphi(y)\sigma(y) + (1-\varphi(y))\sigma(y_0),
\\
A_{\kappa, y_0}(y) &= \varphi(y)A_{\kappa}(y) + (1-\varphi(y))A_{\kappa}(y_0).
\end{align*}
In the following, $C$ denotes a generic constant depending on $m_0$, $m_1$,
$m_2$, $m_3$, $N$, $\epsilon$ and $q$, and $C_{M_2}$ denotes a 
generic constant depending on $M_2$, $m_0$, $m_1$,
$m_2$, $m_3$, $N$, $\epsilon$ and $q$.

Given 
$\bv \in H^2_q(\Omega_+)^N$ and $h \in H^3_q(\Omega_+)$, 
let $K_b(\bv, h)$ be a unique solution of the weak Dirichlet
problem:
\begin{equation}\label{weakD-bent}
(\nabla K_b(\bv, h), \nabla\varphi)_{\Omega_+}
= (\DV(\mu_{y_0}\bD(\bv))-\nabla\dv\bv, \nabla\varphi)_{\Omega_+} 
\end{equation}
for any $\varphi \in \hat H^1_{q',0}(\Omega_+)$ 
subject to 
$$K_b(\bv, h) = <\mu_{y_0}\bD(\bv)\bn_+, \bn_+>
 -\sigma_{y_0} \Delta_{\Gamma_+}h
-\dv\bv\quad\text{ on $\Gamma_+$}.
$$  
We then consider the following equations:
\begin{equation}\label{5.1} \left\{\begin{aligned}
\lambda\bv - \DV(\mu_{y_0}\bD(\bv) - K_b(\bv, h)\bI) &= \bg 
&\quad&\text{in $\Omega_+$}, \\
\lambda h + A_{\kappa, y_0}\cdot\nabla_{\Gamma_+}h - \bv\cdot\bn_+
&= g_d &\quad&\text{on $\Gamma_+$}, \\
(\mu_{y_0}\bD(\bv)- K_b(\bv, h)\bI)\bn_+  - \sigma_{y_0}
(\Delta_{\Gamma_+}h)\bn_+ &=
\bg_b &\quad&\text{on $\Gamma_+$}.
\end{aligned}\right.\end{equation}
The following theorem is a main result in this section. 
\begin{thm}\label{thm:b} Let $1 < q < \infty$ and $0 < \epsilon
< \pi/2$. Let $\gamma_\kappa$ be the number 
defined in Theorem \ref{thm:rbdd:1.2}. 
Then, there exist $M_1 \in (0,1)$, $\tilde\lambda_0 \geq 1$ and operator
families $\CA_b(\lambda)$ and $\CH_b(\lambda)$ with 
\begin{align*}
\CA_b(\lambda) &\in \Hol(\Lambda_{\kappa, \tilde\lambda_0\gamma_\kappa}, 
\CL(\CY_q(\Omega_+), H^2_q(\Omega_+)^N)), \\
\CH_b(\lambda) & \in \Hol(\Lambda_{\kappa, \tilde\lambda_0\gamma_\kappa}, 
\CL(\CY_q(\Omega_+), H^3_q(\Omega_+)))
\end{align*}
such that for any $\lambda=\gamma + i\tau 
\in \Lambda_{\kappa, \tilde\lambda_0\gamma_\kappa}$ and 
$(\bg, g_d, \bg_b) \in Y_q(\Omega_+)$, 
$$\bu = \CA_b(\lambda)(\bg, g_d, \lambda^{1/2}\bg_b, \bg_b),\quad
h = \CH_b(\lambda)(\bg, g_d, \lambda^{1/2}\bg_b, \bg_b)
$$
are unique solutions of Eq. \eqref{5.1}, and 
\begin{equation}\label{ad.5.9.2}\begin{split}
\CR_{\CL(\CY_q(\Omega_+), H^{2-j}_q(\Omega_+)^N)}
(\{(\tau\pd_\tau)^\ell (\lambda^{j/2}\CA_b(\lambda)) \mid
\lambda \in \Lambda_{\kappa, \tilde\lambda_0\gamma_\kappa}\}) \leq r_b, \\
\CR_{\CL(\CY_q(\Omega_+), H^{3-k}_q(\Omega_+))}
(\{(\tau\pd_\tau)^\ell (\lambda^k\CH_b(\lambda)) \mid
\lambda \in \Lambda_{\kappa, \tilde\lambda_0\gamma_\kappa}\}) \leq r_b, 
\end{split}\end{equation}
for $\ell=0,1$, $j=0,1,2$, and $k=0,1$. Where, $r_b$ is a constant
depeding on $m_0$, $m_1$, $m_2$, $N$, $q$, and $\epsilon$, but
independent of $M_1$ and $M_2$, and moreover, $\tilde\lambda_0$
is a constant depending on $M_2$.
\end{thm}
Below, we shall prove Theorem \ref{thm:b}.  By the change of variables
$y=\Phi(x)$, we transform Eq. \eqref{5.1} 
to a problem in the half-space. Let 
\begin{align*}
y_0 &= \Phi(x_0), \quad
\tilde\mu(x) = \varphi(\Phi(x))\mu(\Phi(x)), \quad
\tilde\sigma(x) = \varphi(\Phi(x))\sigma(\Phi(x)), \\
\tilde A_\kappa(x) &= \varphi(\Phi(x))A_\kappa(\Phi(x)).
\end{align*}
Notice that
\begin{align*}
\mu_{y_0}(\Phi(x))&= \mu(y_0) + \tilde\mu(x)-\tilde\mu(x_0), \quad
\sigma_{y_0}(\Phi(x))= \sigma(y_0) + \tilde\sigma(x)-\tilde\sigma(x_0), \\
A_\kappa(\Phi(x',0)) &
= A_\kappa(y_0) + \tilde A_\kappa(x)-\tilde A_\kappa(x_0).
\end{align*}
We may assume that $m_1$, $m_2$, $m_3 \leq M_2$.  Recalling
that $\|\nabla\varphi\|_{H^1_\infty(\BR^N)} \leq M_2$, 
by \eqref{s5.1} and \eqref{s5.2} we have
\begin{equation}\label{st.2*}\begin{aligned}
&|\tilde\mu(x)-\tilde\mu(x_0)| \leq m_1M_1, \quad 
|\tilde\sigma(x) - \tilde\sigma(x_0)| \leq m_1M_1, \\
&\|(\tilde\mu, \tilde\sigma)\|_{L_\infty(\BR^N)} \leq m_1, \quad 
\|\nabla(\tilde\mu, \tilde\sigma)\|_{L_\infty(\BR^N)}
\leq C_{M_2}, \\
&|\tilde A_\kappa(x) - \tilde A_\kappa(x_0)| \leq m_2M_1, \quad
\|\tilde A_\kappa\|_{L_\infty(\BR^N_0)} \leq m_2, \\
&\|\nabla \tilde A_\kappa\|_{W^{1-1/q}_q(\BR^N_0)} \leq C_{M_2}\kappa^{-b}
\quad\text{for $\kappa \in (0, 1)$.}
\end{aligned}\end{equation}
Since $x = \Phi^{-1}(y)$, we have
\begin{equation}\label{change:1*}
\frac{\pd}{\pd y_j} = \sum_{k=1}^N(a_{kj}+ b_{kj}(x))\frac{\pd}{\pd x_k}
\end{equation}
where $(\nabla \Phi^{-1})(\Phi(x)) = (a_{ij} + b_{ij}(x))$.  Let 
$$\fg := \det \nabla \Phi, \quad \tilde \fg = \fg -1.$$
By \eqref{assump:4}, 
\begin{equation}\label{res:5}
\|\tilde\fg\|_{L_\infty(\BR^N)} \leq C_NM_1, \quad
\|\nabla\tilde\fg\|_{L_\infty(\BR^N)} \leq C_NC_A, \quad 
\|\nabla^2\tilde\fg\|_{L_\infty(\BR^N)} \leq C_{M_2}.
\end{equation}
By the change of variables: $y = \Phi(x)$, the weak Dirichlet problem:
$$(\nabla u, \nabla\varphi)_{\Omega_+}
= (\bk, \nabla\varphi)_{\Omega_+}
\quad\text{for any $\varphi \in \hat H^1_{q',0}(\Omega_+)$},
$$
subject to $u=k$ on $\Gamma_+$, is transformed to the following
variational problem:
\begin{equation}\label{var:1}
(\nabla v, \nabla\psi)_\HS + (\CB^0\nabla v, \nabla\psi)_\HS
= (\bh, \nabla\psi)_\HS
\quad\text{for any $\psi \in \hat H^1_{q',0}(\HS)$},
\end{equation}
subject to $v = h$, where $\bh = \fg(\CA_{-1}
+B_{-1}\circ\Phi)\bk\circ\Phi$ and $h=k\circ\Phi$.  Moreover, 
$\CB^0$ is an $N\times N$ matrix whose $(\ell, m)^{\rm th}$ component, 
$\CB_{\ell m}^0$, is given by
$$\CB^0_{\ell m} = \tilde \fg \delta_{\ell m} + 
\fg\sum_{j=1}^N(a_{\ell j}b_{mj}(x) + a_{mj}b_{\ell j}(x)
+ b_{\ell j}b_{mj}(x)).
$$
By \eqref{assump:4}, we have
\begin{equation}\label{res:6}\begin{aligned}
&\|\CB^0_{\ell m}\|_{L_\infty(\BR^N)} \leq C_NM_1, \quad
\|\nabla\CB_{\ell m}^0\|_{L_\infty(\BR^N)} 
\leq C_NC_A, \\ 
&\|\nabla^2\CB_{\ell m}^0\|_{L_\infty(\BR^N)} 
\leq C_{M_2}.
\end{aligned}\end{equation}
\begin{lem}\label{lem:1.4.2} Let $1 < q < \infty$.  Then, 
there exist an $M_1 \in (0, 1)$ and an operator $\CK_1$ with 
$$\CK_1 \in \CL(L_q(\HS)^N, H^1_q(\HS) + \hat H^1_{q,0}(\HS))
$$
such that for any $\bff \in L_q(\HS)^N$and 
$f \in H^1_q(\HS)$, $v = \CK_1(\bff, f)$ is a unique solution
of the variational problem:
\begin{equation}\label{1.97} (\nabla v, \nabla\psi)_\HS
+ (\CB^0\nabla v, \nabla\psi)_\HS =(\bff, \nabla\psi)_\HS
\quad\text{for any $\psi \in \hat H^1_{q',0}(\HS)$},
\end{equation}
subject to $v= f$ on $\BR^N_0$, which possesses the estimate:
\begin{equation}\label{1.98}
\|\nabla v\|_{L_q(\HS)} \leq C_{M_2}(\|\bff\|_{L_q(\HS)}
+ \|f\|_{H^1_q(\HS)}).
\end{equation}
\end{lem}
\pf We know the unique existence theorem of the
variational problem: 
$$(\nabla v, \nabla \psi)_\HS = (\bff, \nabla\psi)_\HS
\quad \text{for any $\psi \in \hat H^1_{q',0}(\HS)$},$$
subject to $v = f$ on $\HS$. Thus, choosing $M_1 > 0$ small 
enough in \eqref{res:6} and using the Banach fixed point theorem, we 
can easily prove the lemma. \qed
\vskip0.5pc
Using the change of the unknown functions: $\bu = \CA_{-1}\bv\circ\Phi$ as 
well as the change of variable: $y=\Phi(x)$, we will derive the problem
in $\HS$ from \eqref{5.1}. Noting that $\CA = {}^\top\CA_{-1}$, 
by \eqref{change:1*} we have
\begin{equation}\label{stress:1}
D_{ij}(\bv)
= \sum_{k,\ell=1}^Na_{ki}a_{\ell j}D_{k\ell}(\bu) + b^d_{ij}:\nabla\bu
\end{equation}
with
$b^d_{ij}:\nabla\bu = \sum_{k,\ell=1}^N
a_{kj}b_{\ell i}D_{k\ell}(\bu)$. 
Setting $\bb_+(x) = {}^\top(b_{+1}, \ldots, b_{+N})$ in \eqref{res:0}, 
by \eqref{res:0} we have
\begin{equation}\label{stress:2}
<\bD(\bv) \bn_+, \bn_+> = <\bD(\bu)\bn_0, \bn_0> + \CB^1:\nabla\bu
\end{equation}
where we have set 
\begin{align*}
\CB^1:\nabla\bu &= -2\sum_{i,j=1}^Na_{ji}b_{+i}D_{jN}(\bu) 
+ \sum_{i,j,k,\ell=1}^Na_{ki}a_{\ell j}b_{+i}b_{+j}D_{k\ell}(\bu)
 \\
& + \sum_{i,j=1}^N(b^d_{ij}:\nabla \bu)(a_{Ni} + b_{+i})(a_{Nj}+b_{+j}).
\end{align*}
By \eqref{assump:4}, we have
\begin{equation}\label{res:7}\begin{split} 
\|\CB^1:\nabla\bu\|_{L_q(\HS)} &\leq C_NM_1\|\nabla\bu\|_{L_q(\HS)},\\
\|\CB^1:\nabla \bu\|_{H^1_q(\HS)} &\leq C_N\{M_1\|\nabla^2 \bu\|_{L_q(\HS)}
+ C_A\|\bu\|_{H^1_q(\HS)}\}.
\end{split}\end{equation}
And also, 
\begin{equation}\label{stress:3}
\dv \bv = \dv\bu + \CB^2:\nabla\bu\quad\text{
with }\quad
\CB^2:\nabla\bu = \sum_{\ell, k=1}^M (\sum_{j=1}^N b_{kj}a_{\ell j})
\frac{\pd u_\ell}{\pd x_k}.
\end{equation}
By \eqref{assump:4}, we have
\begin{equation}\label{res:8} \begin{split}
\|\CB^2:\nabla\bu\|_{L_q(\HS)} &\leq C_NM_1\|\nabla\bu\|_{L_q(\HS)},\\
\|\CB^2:\nabla \bu\|_{H^1_q(\HS)} &\leq C_N\{M_1\|\nabla^2 \bu\|_{L_q(\HS)}
+ C_A\|\bu\|_{H^1_q(\HS)}\}.
\end{split}\end{equation}
By \eqref{stress:1}, we have 
\begin{equation}\label{stress:4}
\CA_{-1}\DV(\mu_{y_0}\bD(\bv)) = \DV(\mu(y_0)\bD(\bu)) + \CR^1:\bu
\end{equation}
with $\CR^1:\bu = {}^\top(\CR^1:\bu|_1, \ldots, \CR^1:\bu|_N)$, and 
\begin{align*}
&\CR^1:\bu|_s \\
&= \sum_{k=1}^N\frac{\pd}{\pd x_k}
\{(\tilde \mu(x)-\tilde \mu(x_0))
D_{sk}(\bu)\} + \sum_{i,j,k=1}^Na_{si}a_{kj}\frac{\pd}{\pd x_k}(\tilde\mu(x)b^d_{ij}:\nabla\bu) \\
&+ \sum_{j,k,\ell,m=1}^Na_{\ell j}b_{kj}\frac{\pd}{\pd x_k}
(\tilde\mu D_{s\ell}(\bu))
+ \sum_{i,j,k=1}^Na_{si}b_{kj}\frac{\pd}{\pd x_k}
(\tilde \mu b^d_{ij}:\nabla\bu).
\end{align*}
By \eqref{assump:4} and \eqref{st.2*}, 
\begin{equation}\label{res:9.1}
\|\CR^1:\bu\|_{L_q(\HS)} \leq C_Nm_1M_1\|\nabla^2\bu\|_{L_q(\HS)} 
+C_{M_2}\|\bu\|_{H^1_q(\HS)}. 
\end{equation}
And also, by \eqref{stress:1}
$$
(\CA_{-1}+B_{-1}\circ\Phi^{-1})
\DV(\mu_{y_0}\bD(\bv)) = \DV(\mu(y_0)\bD(\bu)) + \CR^2:\bu
$$
with $\CR^2:\bu = (\CR^2:\bu|_1, \ldots, \CR^2:\bu|_N)$ and 
\begin{align*}
\CR^2:\bu|_s &= \CR^1:\bu|_s \\
&+ \sum_{i,j,k=1}^Nb_{si}(a_{kj} + b_{kj})\frac{\pd}{\pd x_k}
[\tilde \mu(x)\{\sum_{\ell, m=1}^Na_{\ell i}a_{mj}
D_{\ell m}(\bu) + b^d_{ij}:\nabla\bu\}].
\end{align*}
By \eqref{assump:4}
\begin{equation}\label{res:9.2}
\|\CR^2:\bu\|_{L_q(\HS)} \leq C_Nm_1M_1\|\nabla^2\bu\|_{L_q(\HS)} 
+ C_{M_2}\|\bu\|_{H^1_q(\HS)}.
\end{equation}
And also, we have
$$
(\CA_{-1} + B_{-1}\circ\Phi^{-1})(\nabla\dv\bv)\circ\Phi 
= \nabla \dv\bu + \CR^3:\bu
$$
with $\CR^3:\bu = (\CR^3:\bu|_1, \ldots, \CR^3:\bu|_N)$ and 
\begin{align*}
\CR^3:\bu|_s &= \frac{\pd}{\pd x_s}(\CB^2:\nabla\bu)\\
& + 
\sum_{k=1}^N\{\sum_{i=1}^N(
a_{si}b_{ki} +b_{si}(a_{ki} + b_{ki}))\}\frac{\pd}{\pd x_k}
(\dv\bu + \CB^2:\nabla\bu).
\end{align*}
By \eqref{assump:4}
\begin{equation}\label{res:10}
\|\CR^3:\bu\|_{L_q(\HS)} \leq C_Nm_1(M_1\|\nabla^2\bu\|_{L_q(\HS)} 
+ C_A\|\bu\|_{H^1_q(\HS)}).
\end{equation}
Let 
$$\bff(\bu) : = \fg(\CA_{-1}+B_{-1}\circ\Phi)(\DV(\mu_{y_0}\bD(\bv))
 - \nabla\dv\bv)\circ\Phi, $$
and then 
$$\bff(\bu) = \DV(\mu(y_0)\bD(\bu)) - \nabla\dv\bu + \CR^4:\bu
$$
with $\CR^4:\bu=(\CR^4:\bu|_1, \ldots \CR^4:\bu|_N)$ and 
$$\CR^4:\bu|_s = \tilde\fg(\DV(\mu(y_0)\bD(\bu)) - \nabla\dv\bu)
+ \fg\CR^2:\bu-\fg\CR^3:\bu.$$
By \eqref{assump:4}, \eqref{st.2*}, \eqref{res:5}, \eqref{res:9.1}, 
 \eqref{res:9.2}, 
 and \eqref{res:10},
\begin{equation}\label{res:11}
\|\CR^4:\bu\|_{L_q(\HS)} \leq C_N(m_1+1)M_1\|\nabla^2\bu\|_{L_q(\HS)} 
+ C_{M_2}\|\bu\|_{H^1_q(\HS)}.
\end{equation}
In view of \eqref{res:4}, \eqref{stress:2} and \eqref{stress:3}, 
setting 
\begin{align*}
\rho &= h\circ\Phi, \\
f(\bu, \rho) &= <(\tilde\mu(x)-\tilde\mu(x_0))
\bD(\bu)\bn_0, \bn_0> + \tilde\mu\CB^1:\nabla\bu \\
&- 
(\tilde\sigma(x)-\tilde\sigma(x_0))\Delta'\rho
-\tilde\sigma(x) \CD_+\rho - \CB^2:\nabla\bu,
\end{align*}
we have
\begin{align*}
&<\mu_{y_0}\bD(\bv)\bn_+, \bn_+> - \sigma_{y_0}
 \Delta_{\Gamma_+}h - \dv \bv \\
&
\quad= <\mu(y_0)\bD(\bu)\bn_0, \bn_0> - 
\sigma(y_0) \Delta'\rho - \dv \bu + f(\bu, \rho).
\end{align*}
Thus, $K_1(\bu, \rho) = K_b(\bv, h)\circ\Phi$ satisfies the variational 
equation:
\begin{align*}
&(\nabla K_1(\bu, \rho), \nabla\psi)_\HS
+(\CB^0\nabla K_1(\bu, \rho), \nabla\psi)_\HS\\
&\quad
= (\DV(\mu(y_0)\bD(\bu))-\nabla \dv\bu +\CR^4:\bu, \nabla\psi)_\HS
\end{align*}
for any $\psi\in \hat H^1_{q',0}(\HS)$,
subject to 
$$K_1(\bu, \rho) = <\mu(y_0)\bD(\bu)\bn_0, \bn_0> - 
\sigma(y_0) \Delta'\rho - \dv \bu + f(\bu, \rho)
\quad\text{on $\BR^N_0$}.$$

Let $\tilde K(\bu, \rho) \in H^1_q(\HS) + \hat H^1_{q,0}(\HS)$ be a unique
solution of the weak Dirichlet problem:
$$(\nabla \tilde K(\bu, \rho), \nabla\psi)_\HS = (\DV(\mu(y_0)\bD(\bu))-
\nabla\dv\bu, \nabla\psi)_\HS  
$$
for any $\psi \in \hat H^1_{q',0}(\HS)$
subject to 
$$\tilde K(\bu, \rho) = <\mu(y_0)\bD(\bu)\bn_0, \bn_0> 
- \sigma(y_0)\Delta'\rho-\dv \bu
\quad\text{ on $\BR^N_0$}.
$$ 
Setting $K_1(\bu, \rho) = \tilde K(\bu, \rho) + K_2(\bu, \rho)$, 
we then see that $K_2(\bu, \rho)$ satisfies the variational equation:
$$(\nabla K_2(\bu, \rho), \nabla\psi)_\HS +
(\CB^0\nabla K_2(\bu, \rho), \nabla\psi)_\HS 
= (\CR^4:\bu - \CB^0\nabla \tilde K(\bu, \rho), \nabla\psi)_\HS 
$$
for any $\varphi \in \hat H^1_{q',0}(\HS)$,
subject to $K_2(\bu, \rho) = f(\bu, \rho)$ on $\BR^N_0$. In view of 
Lemma \ref{lem:1.4.2}, we have
$$K_2(\bu, \rho) = \CK_1(\CR^4:\bu - 
\CB^0\nabla \tilde K(\bu, \rho), f(\bu, \rho)).
$$
By Lemma \ref{lem:1.4.2},  \eqref{res:6}, 
\eqref{res:7}, \eqref{res:8}, 
\eqref{res:2}, and  \eqref{res:11}, we have
\begin{equation}\label{res:12}\begin{split}
&\|\nabla K_2(\bu, \rho)\|_{L_q(\HS)} \\
&\quad \leq C_N(1+m_1)
M_1(\|\nabla^2\bu\|_{L_q(\HS)} + \|\nabla^3\rho\|_{L_q(\HS)})
\\
&\quad\quad
+ C_{M_2}(\|\bu\|_{H^1_q(\HS)} + \|\rho\|_{H^2_q(\HS)}).
\end{split}\end{equation}
Since 
\begin{align*}
\CA_{-1}\nabla K_b(\bv, h)|_s &= \sum_{i,k=1}^N a_{si}(a_{ki}+b_{ki})
\frac{\pd}{\pd x_k}K_1(\bu, \rho) \\
& = \frac{\pd}{\pd x_s} \tilde K(\bu, \rho) 
+ \sum_{k=1}^N(\sum_{i=1}^Na_{si}b_{ki})\frac{\pd}{\pd x_k}
\tilde K(\bu, \rho) \\
&+ \sum_{k=1}^N(\delta_{ks} + \sum_{i=1}^N a_{si}b_{ki})
\frac{\pd}{\pd x_k}K_2(\bu, \rho),
\end{align*} 
by \eqref{stress:4} we see that the first equation of Eq.\eqref{5.1} 
is transformed to
$$\lambda\bu -\DV(\mu(y_0)\bD(\bu)-\tilde K(\bu, \rho)\bI) + \CR^5(\bu, \rho)
= \bh\quad\text{in $\HS$}, \\
$$
where $\bh = \CA_{-1}\bg\circ\Phi$, 
$R^5(\bu, \rho) = (\CR^5(\bu, \rho)|_1, \ldots, \CR^5(\bu, \rho)|_N)$, 
and
\begin{align*}
\CR^5(\bu, \rho)|_s &= -\CR^1:\bu|_s 
+ \sum_{k=1}^N(\sum_{i=1}^Na_{si}b_{ki})\frac{\pd}{\pd x_k}
\tilde K(\bu, \rho) \\
&+ \sum_{k=1}^N(\delta_{ks} + \sum_{i=1}^N a_{si}b_{ki})
\frac{\pd}{\pd x_k}K_2(\bu, \rho).
\end{align*}
By \eqref{res:0}, we have
\begin{align*}
\bv\cdot\bn_+ &= -({}^\top\CA_{-1}\bu)\cdot{}^\top(a_{N1}, \ldots, a_{NN})
 + ({}^\top\CA_{-1}\bu)\cdot\bb_+ \\
&= \bu\cdot\bn_0 + \bu\cdot(\CA_{-1}\bb_+),
\end{align*}
and so the second equation of Eq.\eqref{5.1} is transformed
to
$$\lambda\rho + A_\kappa(y_0)\cdot\nabla'\rho - \bu\cdot\bn_0
+ \CR^6_\kappa(\bu, \rho) = h_d$$
with $h_d = g_d\circ\Phi$ and 
\begin{align*}
\CR_0^6(\bu, \rho) &= - \bu\cdot(\CA_{-1}\bb_+)
\quad\text{for $\kappa = 0$},  \\
\CR^6_\kappa(\bu, \rho) &= 
 (\tilde A_\kappa(x)-\tilde A_\kappa(x_0))\nabla'\rho
- \bu\cdot(\CA_{-1}\bb_+) \quad\text{for $\kappa \in (0, 1)$}. 
\end{align*}
 By \eqref{res:0} and \eqref{stress:1},
 we have
$\CA_{-1}\mu_{y_0}\bD(\bv)\bn_+= \mu(y_0)\bD(\bu)\bn_0 + \CR^7_1(\bu)$, where 
 $\CR^7_1(\bu)$ is an $N$- vector of functions whose $s^{\rm th}$ component,
$\CR^7_1(\bu)|_s$, is defined by
\begin{align*}
&\CR^7_1(\bu)|_s = -(\tilde\mu(x)-\tilde\mu(x_0))D_{sN}(\bu)\\
&+(\mu(y_0)+\tilde\mu(x)-\tilde\mu(x_0))
\sum_{i,j=1}^N(a_{ij}b_{+j}D_{si}(\bu)
+a_{si}b_{ij}^d:\nabla\bu(-a_{Nj} + b_{+j})).
\end{align*}
By \eqref{res:0}, 
$$\CA_{-1}K_b(\bv, h)\bn_+ = \tilde K(\bu, \rho)\bn_0 + 
\tilde K(\bu, \rho)\CA_{-1}\bb_+
+ K_2(\bu,\rho)(\bn_0+\CA_{-1}\bb_+). $$
By \eqref{res:4}, 
\begin{align*}
\CA_{-1}\sigma_{y_0}(\Delta_{\Gamma_+}\bh)\bn_+ &= 
\sigma(y_0)(\Delta'\rho)\bn_0 
+ (\tilde\sigma(x)-\tilde\sigma(x_0))(\Delta'\rho)\bn_0 \\
&+ \tilde\sigma(x)\{
(\Delta'\rho)(\CA_{-1}\bb_+) + (\CD_+\rho)(\bn_0 + \CA_{-1}\bb_+)\}.
\end{align*}
Putting formulas above together yields that the third equation of
Eq.\eqref{5.1} is transformed to the equation: 
$$(\mu(y_0)\bD(\bu) - \tilde K(\bu, \rho)\bI)\bn_0 - \sigma(y_0)
(\Delta'\rho)\bn_0
+ \CR^7(\bu, \rho) = \bh_b
\quad\text{on $\BR^N_0$},
$$
where $\bh_b = \CA_{-1}\bg_b\circ\Phi$, and 
\begin{align*}
&\CR^7(\bu, \rho) = \CR^7_1(\bu, \rho) - \tilde K(\bu, \rho)(\CA_{-1}\bb_+)
- K_2(\bu, \rho)(\bn_0 + \CA_{-1}\bb_+)\\
&- (\tilde\sigma(x)-\tilde\sigma(x_0))(\Delta'\rho)\bn_0 
- \tilde\sigma(x)\{
(\Delta'\rho)(\CA_{-1}\bb_+) + (\CD_+\rho)(\bn_0 + \CA_{-1}\bb_+)\}.
\end{align*}
Summing up, we have seen that Eq.\eqref{5.1} is transformed to
the following equations:
\begin{equation}\label{5.2}\left\{\begin{aligned}
\lambda\bu - \DV(\mu(y_0)\bD(\bu) - \tilde K(\bu, \rho)\bI) 
+ \CR^5(\bu, \rho)
&= \bh &\quad&\text{in $\HS$}, \\
\lambda \rho + A_\kappa(y_0)\cdot\nabla'\rho - \bu\cdot\bn_0 
+R^6_\kappa(\bu,\rho) &= h_d &\quad&\text{on $\BR^N_0$}, \\
(\mu(y_0)\bD(\bu) - \tilde K(\bu, \rho)\bI)\bn_0 
- \sigma(y_0)(\Delta'\rho)\bn_0
+ \CR^7(\bu, \rho) &= \bh_b &\quad&\text{on $\BR^N_0$},
\end{aligned}\right.\end{equation}
where $\bh= \CA_{-1}\bg\circ\Phi$, $h_d = g_d\circ\Phi$, 
$\bh_d = \CA_{-1}\bg_d\circ\Phi$, and 
$\CR^5(\bu, \rho)$, $\CR^6_\kappa(\bu, \rho)$
 and $\CR^7(\bu, \rho)$ are linear in $\bu$ and $\rho$ and satisfy
the estimates:
\begin{equation}\label{5.3}\begin{split}
\|\CR^5(\bu, \rho)\|_{L_q(\HS)} & \leq CM_1(\|\nabla^2\bu\|_{L_q(\HS)}
+ \|\nabla^3\rho\|_{L_q(\HS)})\\
&\phantom{aaaaa}
+ C_{M_2}(\|\bu\|_{H^1_q(\HS)} + \|\rho\|_{H^2_q(\HS)}), \\
\|\CR^6_0(\bu, \rho)\|_{W^{2-1/q}_q(\BR^N_0)}
& \leq CM_1\|\nabla^2\bu\|_{L_q(\HS)} + C_{M_2}\|\bu\|_{H^1_q(\HS)}, \\
\|\CR^6_\kappa(\bu, \rho)\|_{W^{2-1/q}_q(\BR^N_0)}
& \leq CM_1(\|\nabla^2\bu\|_{L_q(\HS)} + \|\nabla^3\rho\|_{L_q(\HS)})\\
&\phantom{aaaaa}
+ C_{M_2}(\|\bu\|_{H^1_q(\HS)} + \kappa^{-b}\|\rho\|_{H^2_q(\HS)}),\\
\|\CR^7(\bu, \rho)\|_{L_q(\HS)} & \leq CM_1(\|\nabla\bu\|_{L_q(\HS)}
+ \|\nabla^2\rho\|_{L_q(\HS)})\\
&\phantom{aaaaa}
+ C_{M_2}(\|\bu\|_{L_q(\HS)} + \|\rho\|_{H^1_q(\HS)}), \\
\|\CR^7(\bu, \rho)\|_{H^1_q(\HS)} & \leq CM_1(\|\nabla^2\bu\|_{L_q(\HS)}
+ \|\nabla^3\rho\|_{L_q(\HS)})\\
&\phantom{aaaaa}
+ C_{M_2}(\|\bu\|_{H^1_q(\HS)} + \|\rho\|_{H^2_q(\HS)}).
\end{split}\end{equation}
Here and in the following, $C$ denotes a generic constant depending on
$N$, $q$, $m_1$, and $m_2$ and $C_{M_2}$ a generic constant depending on
$N$, $q$, $m_1$, $m_2$, $m_3$ and $M_2$. 
By Theorem \ref{main:half1}, there exists a large numbger $\lambda_0$ 
and operator families $\CA_0(\lambda)$ and $\CH_0(\lambda)$ with
\begin{align*}
\CA_0(\lambda) & \in {\rm Hol}\,(\Lambda_{\kappa, \lambda_0}, 
\CL(\CY(\HS), H^2_q(\HS)^N)), \\
\CH_0(\lambda) & \in {\rm Hol}\,(\Lambda_{\kappa, \lambda_0},
\CL(\CY(\HS), H^3_q(\HS)))
\end{align*}
such that for any $\lambda \in \Lambda_{\kappa, \lambda_0}$
and $(\bff, d, \bh) \in Y_q(\HS)$, $\bu$ and $\rho$ with
$$\bu = 
\CA_0(\lambda)F_\lambda(\bff, d, \bh), \quad 
\rho = \CH_0(\lambda)F_\lambda(\bff, d, \bh),$$
where $F_\lambda(\bff, d, \bh) = 
(\bff, d, \lambda^{1/2}\bh, \bh)$, are unique solutions of the equations:
$$\left\{\begin{aligned}
\lambda \bu - \DV(\mu(y_0)\bD(\bu) - \tilde K(\bu, \rho)\bI) &=\bff
&\quad&\text{in $\BR^N_+$}, \\
\lambda \rho + A_\kappa(y_0)\cdot\nabla'\rho - \bu\cdot\bn_0 &=d
&\quad&\text{on $\BR^N_0$}, \\
(\mu(y_0)\bD(\bu) - \tilde K(\bu, \rho)\bI)\bn_0 
-\sigma(y_0)(\Delta'\rho)\bn_0
& = \bh,
&\quad&\text{on $\BR^N_0$},
\end{aligned}\right.$$
and
\begin{align*}
\CR_{\CL(\CY(\HS), H^{2-j}_q(\HS)^N)}
(\{(\tau\pd_\tau)^s(\lambda^{j/2}\CA_0(\lambda)) \mid 
\lambda \in \Lambda_{\kappa, \lambda_0}\}) &\leq r_b, \\
\CR_{\CL(\CY(\HS), H^{3-k}_q(\HS))}
(\{(\tau\pd_\tau)^s(\lambda^k\CH_0(\lambda)) \mid 
\lambda \in \Lambda_{\kappa, \lambda_0}\}) &\leq r_b
\end{align*}
for $s=0,1$, $j=0,1,2$, and $k=0,1$.  Where, $r_b$ is a constant
depending on $\epsilon$, $N$, $m_1$, and $m_2$.

 Let $\bu
= \CA_0(\lambda)F_\lambda(\bh, h_d, \bh_b)$ and 
$\rho = \CH_0(\lambda)F_\lambda(\bh, h_d, \bh_b)$ in \eqref{5.2}. 
Then, Eq.\eqref{5.2} is rewritten as 
\begin{equation}\label{5.5}\left\{\begin{aligned}
&\lambda\bu - \DV(\mu(y_0)\bD(\bu) - \tilde K(\bu, \rho)\bI) + 
\CR^5(\bu, \rho) \\
&\quad = \bh
+ \CR^8(\lambda)F_\lambda(\bh, h_d, \bh_b)&&\text{in $\HS$}, \\
&\lambda\rho + A_\kappa\cdot\nabla'\rho-\bu\cdot\bn_0
+ R^6_\kappa(\bu, \rho) \\
&\quad =h_d
+ \CR^8_d(\lambda)F_\lambda(\bh, h_d, \bh_b)&&\text{in $\HS$}, \\
&(\mu(y_0)\bD(\bu) - \tilde K(\bu, \rho)\bI)\bn_0
- \sigma(y_0)(\Delta'\rho)\bn_0
+ \CR^7(\bu, \rho)\\
&\quad = \bh_b + \CR^8_b(\lambda)F_\lambda(\bh, h_d, \bh_b)
&&\text{on $\BR^N_0$},
\end{aligned}\right.\end{equation}
where we have set 
\begin{align*}
&\CR^8(\lambda)(F_1, F_2, F_3, F_4) \\ 
&\quad
= \CR^5(\CA_0(\lambda)(F_1, F_2, F_3, F_4),
\CH_0(\lambda)(F_1, F_2, F_3, F_4)), \\
&\CR^8_d(\lambda)(F_1, F_2, F_3, F_4) \\ 
&\quad
= \CR^6_\kappa(\CA_0(\lambda)(F_1, F_2, F_3, F_4), 
\CH_0(\lambda)(F_1, F_2, F_3, F_4)), \\
&\CR^8_b(\lambda)(F_1, F_2, F_3, F_4)  \\
&\quad= \CR^7(\CA_0(\lambda)(F_1, F_2, F_3, F_4), 
\CH_0(\lambda)(F_1, F_2, F_3, F_4)).
\end{align*}
Let 
$$\CR^9(\lambda)F = (\CR^8(\lambda)F, 
 \CR^8_d(\lambda)F,  \CR^8_b(\lambda)F)
$$
for $F = (F_1, F_2, F_3, F_4) \in \CY_q(\HS)$.  
Notice that 
$$\CR^9(\lambda)F = (\CR^8(\lambda)F,  \CR^8_d(\lambda)F,  
\lambda^{1/2}\CR^8_b(\lambda)F, \CR^8_b(\lambda)F) \in \CY_q(\HS), 
$$
for $F = (F_1, F_2, F_3, F_4) \in \CY_q(\HS)$
 and that the right side of Eq.\eqref{5.5} is written as 
$(\bh, h_d, \bh_b) + F^9(\lambda)F_\lambda(\bh, h_d, \bh_b)$.
By \eqref{5.3}, \eqref{res:1}, Proposition \ref{prop:4.1}, 
and Theorem \ref{main:half1}, 
we have
\begin{equation}\label{5.4*} \CR_{\CL(\CY_q(\HS))}
(\{(\tau\pd_\tau)^\ell(F_\lambda\CR^9(\lambda)) \mid
\lambda \in \Lambda_{\kappa, \lambda_1}\}) 
\leq CM_1 + C_{M_2}(\lambda_1^{-1/2} + \lambda_1^{-1}\gamma_\kappa)
\end{equation}
for any $\lambda_1 \geq \lambda_0$. Here and in the following,
 $C$ denotes a generic constant depending on $N$, $\epsilon$, $m_1$, $m_2$,
and $C_A$, 
and $C_{M_2}$ denotes a generic constant 
depending on $N$, $\epsilon$, $m_1$, $m_2$, $m_3$, $C_A$, and 
$M_2$.   Choosing $M_1$ so small that
$CM_1 \leq 1/4$ and choosing $\lambda_1 > 0$ so large that
$C_{M_2}\lambda_1^{-1/2} \leq 1/8$ and 
$C_{M_2}\lambda_1^{-1}\gamma_\kappa \leq 1/8$, by \eqref{5.4*} we have
\begin{equation}\label{5.5.1}
\CR_{\CL(\CY_q(\HS))}(\{(\tau\pd_\tau)^\ell(F_\lambda \CR^9(\lambda))
\mid \lambda \in \Lambda_{\kappa, \lambda_1}\}) 
\leq 1/2
\end{equation}
for $\ell=0,1$. Since $\gamma_\kappa \geq 1$ and we may assume that
$C_{M_2} \geq 1$, if $\lambda_1 \geq 64C^2_{M_2}\gamma_\kappa$,
then $C_{M_2}\lambda_1^{-1/2} \leq 1/8$ and 
$C_{M_2}\lambda_1^{-1}\gamma_\kappa \leq 1/8$.

Recall that for $F = (F_1, F_2, F_3, F_4) \in \CY_q(\HS)$ 
and $(\bh, h_d, \bh_b)
\in Y_q(\HS)$, 
\begin{equation}\label{5.7}\begin{split}
\|(F_1, F_2, F_3, F_4)\|_{\CY_q(\HS)} &= \|(F_1, F_3)\|_{L_q(\HS)} 
+ \|F_2\|_{W^{2-1/q}_q(\BR^N_0)} + \|F_4\|_{H^1_q(\HS)},
\\
\|(\bh, h_d, \bh_b)\|_{X_q(\HS)} 
&= \|\bh\|_{L_q(\HS)} + \|h_d\|_{W^{2-1/q}_q(\BR^N_0)}
+ \|\bh_b\|_{H^1_q(\HS)}
\end{split}\end{equation}
(cf. Remark \ref{rem:1}, where $\Omega$ should be replaced by $\HS$).
By \eqref{5.5.1} we have
\begin{equation}\label{5.6}
\|F_\lambda(\CR^9(\lambda)F_\lambda(\bh, h_d, \bh_b))\|_{\CY_q(\HS)}
\leq (1/2)\|F_\lambda(\bh, h_d, \bh_b)\|_{\CY_q(\HS)}.
\end{equation}
In view of \eqref{5.7}, when $\lambda\not=0$, $\|\CF_\lambda(\bh, h_d, \bh_b)
\|_{\CY_q(\HS)}$ is an equivalent norm to $\|(\bh, h_d, \bh_b)\|_{X_q(\HS)}$.  
Thus, by \eqref{5.6} $(\bI + \CR^9(\lambda)F_\lambda)^{-1}
= \sum_{j=1}^\infty(-\CR^9(\lambda)F_\lambda)^j$ exists in $\CL(X_q(\HS))$.
Setting 
\begin{equation}\label{5.8}\begin{aligned}
\bu &= \CA_0(\lambda)F_\lambda(\bI + 
\CR^9(\lambda)F_\lambda)^{-1}(\bh, h_d, \bh_b), \\
\rho &= \CH_0(\lambda)F_\lambda(\bI + 
\CR^9(\lambda)F_\lambda)^{-1}(\bh, h_d, \bh_b), 
\end{aligned}\end{equation}
by \eqref{5.5} we see that $\bu$ and $\rho$ are solutions of
Eq.\eqref{5.2}.  In view of  \eqref{5.5}, 
$(\bI + F_\lambda\CR^9(\lambda))^{-1} = \sum_{j=0}^\infty
(-F_\lambda\CR^9(\lambda))^j$ exists in $\CL(\CY_q(\HS))$,
and  
\begin{equation}\label{5.9}
\CR_{\CL(\CY_q(\HS))}(\{(\tau\pd_\tau)^\ell
(\bI + F_\lambda\CR^9(\lambda))^{-1} \mid \lambda 
\in \Lambda_{\kappa, \lambda_1}\}) \leq 4
\end{equation}
for $\ell=0,1$. 
Since
\begin{align*}
F_\lambda(\bI + \CR^9(\lambda)F_\lambda)^{-1}
&= F_\lambda\sum_{j=0}^\infty(-\CR^9(\lambda)F_\lambda)^j
= (\sum_{j=0}^\infty(-F_\lambda\CR^9(\lambda))^j)F_\lambda\\ 
&= (\bI + F_\lambda\CR^9(\lambda))^{-1}F_\lambda,
\end{align*}
defining operators $\CA_1(\lambda)$ and $\CH_1(\lambda)$ acting on
$F = (F_1, F_2, F_3, F_4) \in \CY_q(\HS)$ by 
$$\CA_1(\lambda)F = \CA_0(\lambda)(\bI + F_\lambda\CR^9(\lambda))^{-1}F,
\quad 
\CH_1(\lambda)F_1 = \CH_0(\lambda)(\bI + F_\lambda\CR^9(\lambda))^{-1}F,
$$
by \eqref{5.8} $\bu = \CA_1(\lambda)F_\lambda(\bh, h_d, \bh_b)$ and 
$\rho = \CH_1(\lambda)F_\lambda(\bh, h_d, \bh_b)$ are  solutions of
Eq.\eqref{5.2}.  Moreover, by \eqref{5.9} and Theorem \ref{main:half1}
\begin{equation}\label{ad.5.9.1}\begin{split}
\CR_{\CL(\CY_q(\HS), H^{2-j}_q(\HS)^N)}
(\{(\tau\pd_\tau)^\ell(\lambda^{j/2}\CA_1(\lambda)) \mid 
\lambda \in \Lambda_{\kappa, \lambda_1\gamma_\kappa}\}) & \leq 4r_b, \\
\CR_{\CL(\CY_q(\HS), H^{3-k}_q(\HS))}
(\{(\tau\pd_\tau)^\ell(\lambda^k\CH_1(\lambda)) \mid 
\lambda \in \Lambda_{\kappa, \lambda_1\gamma_\kappa}\}) & \leq 4r_b,
\end{split}\end{equation}
for $\ell=0,1$, $j=0,1,2$ and $k=0,1$. Recalling that
\begin{align*}
\bv &= ({}^\top\CA_{-1}\bu)\circ\Phi^{-1}, \quad
h = \rho\circ\Phi^{-1}, \\ 
\bh & =\CA_{-1}\bg\circ\Phi, \quad h_d = g_d\circ\Phi, \quad  
\bh_d = \CA_{-1}\bg_d\circ\Phi,
\end{align*}
 we define operators
$\CA_b(\lambda)$ and $\CH_b(\lambda)$ acting on
$F = (F_1, F_2, F_3, F_4) \in \CY_q(\Omega_+)$ by 
\begin{align*}
&\CA_b(F_1, F_2, F_3, F_4)  \\
&\phantom{\CH_b(F_1, F_2)}
= 
{}^\top\CA_{-1}[\CA_1(\lambda)(\CA_{-1}F_1\circ\Phi, F_2\circ\Phi, 
\CA_{-1}F_3\circ\Phi, 
F_4\circ\Phi)]\circ\Phi^{-1}, \\
&\CH_b(F_1, F_2, F_3, F_4)  = 
[\CH_1(\lambda)(\CA_{-1}F_1\circ\Phi, F_2\circ\Phi, 
\CA_{-1}F_3\circ\Phi, 
F_4\circ\Phi)]\circ\Phi^{-1}.
\end{align*}
Obviously, given any $(\bg, g_d, \bg_b) \in Y_q(\Omega_+)$, 
$\bu = \CA_b(\lambda)F_\lambda(\bg, g_d, \bg_b)$ and 
$h = \CH_b(\lambda)F_\lambda(\bg, g_d, \bg_b)$ are solutions of
Eq.\eqref{5.1}.  From \eqref{assump:4} 
we have 
\begin{align*}
\|g\circ\Phi^{-1}\|_{H^\ell_q(\Omega_+)} 
& \leq C_A\|g\|_{H^\ell_q(\HS)} 
\quad\text{for $\ell = 0, 1, 2$}, \\
\|\nabla^3(g\circ\Phi^{-1})\|_{L_q(\Omega_+)} 
& \leq C_A\|\nabla^2g\|_{H^1_q(\HS)}
+ C_{M_2}\|\nabla g\|_{L_q(\HS)}, \\
\|h\circ\Phi\|_{H^\ell_q(\HS)} & \leq C_A\|h\|_{H^\ell_q(\Omega_+)}
\quad\text{for $\ell=0,1,2$}.
\end{align*}
and so, in view of \eqref{ad.5.9.1} 
we can choose $\tilde\lambda_0 \geq \lambda_1$ suitably
large such that $\CA_b(\lambda)$ and $\CH_b(\lambda)$ satisfy
the estimates:
\begin{align*}
\CR_{\CL(\CY_q(\Omega_+), H^{2-j}_q(\Omega_+)^N)}
(\{(\tau\pd_\tau)^\ell(\lambda^{j/2}\CA_b(\lambda)) \mid 
\lambda \in \Lambda_{\kappa, \lambda_1\gamma_\kappa}\}) & \leq Cr_b, \\
\CR_{\CL(\CY_q(\HS), H^{3-k}_q(\HS))}
(\{(\tau\pd_\tau)^\ell(\lambda^k\CH_b(\lambda)) \mid 
\lambda \in \Lambda_{\kappa, \lambda_1\gamma_\kappa}\}) & \leq Cr_b,
\end{align*}
for $\ell=0,1$, $j=0,1,2$ and $k=0,1$, 
where $C$ and $r_b$ are constants  independent of $M_2$. 
 This completes the existence part of Theorem \ref{thm:b}.

The uniqueness can be proved by showing {\it a priori} estimates of solutions
of Eq. \eqref{5.2} with $(\bh, h_d, \bh_b) = (0, 0, 0)$
in the same manner as in the proof of Theorem \ref{thm:p.1}.
This completes  the proof of Theorem \ref{thm:b} without the operators 
$\CL_1$ and $\CL_2$.

\subsection{Some preparation for the proof of Theorem \ref{main:thm3}}

In the following, we use the symbols given in Proposition
\ref{prop:lap} in Subsec. \ref{subsec:p-1} and we write  
$\Omega_j=\Phi_j(\HS)$, and  $\Gamma_j = 
\Phi_j(\BR^N_0)$ for the sake of simplicity. 
Recall that $B^i_j = B_{r_0}(x^i_j)$. 
In view of the assumptions \eqref{assump:3}, \eqref{assump:4.1},
and \eqref{assump:4.2}, 
we may assume that 
\begin{align}
& |\mu(x)- \mu(x^i_j)| \leq M_1 \quad
\text{for any $x \in B^i_j$}; \nonumber \\
&|\sigma(x)- \sigma(x^1_j)| \leq M_1
\quad\text{for any $x \in \Gamma_j \cap B^i_j$}; 
\nonumber \\
&|A_\kappa(x)-A_\kappa(x^1_j)| \leq M_1 
\quad\text{for any $x \in \Gamma_j \cap B^i_j$}; 
\label{R.5.1*}\\
&m_0\leq \mu(x), \sigma(x) \leq m_1,
\quad|\nabla\mu(x)|, |\nabla\sigma(x)| \leq m_1 
\quad\text{for any $x \in \overline{\Omega}$},  \nonumber \\
&|A_\kappa(x)| \leq m_2 \quad\text{for any $x \in \Gamma$},
\quad 
\|A_\kappa\|_{W^{2-1/q}_r(\Gamma)} \leq m_3\kappa^{-b}
\label{R.5.1}
\end{align}
for any $\kappa \in (0, 1)$.
Here, $m_0$, $m_1$, $m_2$, $m_3$, $b$ and $r$ are constants given
in \eqref{assump:3}. 

We next prepare some propositions used to construct a parametrix. 
\begin{prop}\label{prop:5.1}
Let $X$ be a Banach space and $X^*$ its dual space,
	while $\|\cdot\|_X$, $\|\cdot\|_{X^*}$, and $<\cdot,\cdot>$
denote the norm of $X$, the norm of $X^*$,
	and the duality pairing between of $X$ and $X^*$, respectively.
	Let $n\in\BN$, $l=1,\dots,n$, and $\{a_l\}_{l=1}^n\subset\BC$,
	and let $\{f_j^l\}_{j=1}^\infty$ be sequences 
in $X^*$ and $\{g_j^l\}_{j=1}^\infty$, $\{h_j\}_{j=1}^\infty$ 
be sequences of positive numbers.
	Assume that there exist maps $\CN_j:X\to[0,\infty)$ 
such that
	\begin{equation*}
		|<f_j^l,\varphi>|\leq M_3 g_j^l\CN_j(\varphi) 
\quad (l=1,\dots,n), \quad
		\Big|\Big<\sum_{l=1}^n a_l f_j^l,\varphi\Big>\Big|
 \leq M_3h_j\CN_j(\varphi)
	\end{equation*}
	for any $\varphi\in X$ with 
some positive constant $M_3$ independent of $j\in\BN$ and $l=1,\dots,n$.
	If
	\begin{equation*}
		\sum_{j=1}^\infty \left(g_j^l\right)^q < \infty, \quad 
		\sum_{j=1}^\infty \left(h_j\right)^q < \infty, \quad 
		\sum_{j=1}^\infty\left(\CN_j(\varphi)\right)^{q'}
 \leq \left(M_4\|\varphi\|_{X}\right)^{q'}
	\end{equation*}
	with $1<q<\infty$ and $q'=q/(q-1)$ for some positive constant $M_4$,
	then the infinite sum $f^l = \sum_{j=1}^\infty f_j^l$ 
exists in the strong topology of $X^*$ and
	\begin{equation}\label{160201_1}
		\|f^l\|_{X^*} \leq M_3M_4
\Big(\sum_{j=1}^\infty\big(g_j^l\big)^q \Big)^{1/q}, \quad 
		\Big\|\sum_{l=1}^n a_l f^l \Big\|_{X^*}
 \leq M_3 M_4 \Big(\sum_{j=1}^\infty\big(h_j\big)^q\Big)^{1/q}.
	\end{equation}
\end{prop}
\pf For a proof, see Proposition 9.5.2 in Shibata \cite{S4}. \qed
\vskip0.5pc
The following propositions are used to define the infinite sum of
 $\CR$-bounded operator families defined on $\BR^N$ and $\Omega_j$.  
\begin{prop}
\label{prop:5.2} Let $1 < q < \infty$, $i=0,1$, and 
$n \in \BN_0$. Set $\CH^0_j=\BR^N$ and 
$\CH^1_j = \Omega_j$.  Let $\eta^i_j$ be a function in $C^\infty_0(B^i_j)$
such that $\|\eta^i_j\|_{H^n_\infty(\BR^N)} \leq c_1$ for any 
$j \in \BN$ with some constant $c_1$ independent of 
$j \in \BN$.   
 Let $f_j$ $(j \in \BN)$ be  elements in $H^n_q(\CH^i_j)$
such that $\sum_{j=1}^\infty \|f_j\|_{H^n_q(\CH^i_j)}^q < \infty$.
Then, $\sum_{j=1}^\infty\eta^i_jf_j$ converges some $f \in 
H^n_q(\Omega)$ strongly in $H^n_q(\Omega)$, and 
$$\|f\|_{H^n_q(\Omega)} \leq C_q\{\sum_{j=1}^\infty 
\|f_j\|_{H^n_q(\CH^i_j)}^q\}^{1/q}.$$ 
\end{prop}
\pf For a proof, see Proposition 9.5.3 in Shibata \cite{S4}. \qed
\begin{prop}\label{prop:5.3} 
Let $1 < q < \infty$ and $n=2,3$.  Then we have the following 
assertions. \\
\thetag1~ There exist extension maps $\bT^n_j: W^{n-1/q}_q(\Gamma_j)
\to H^n_q(\Omega_j)$ such that for any $h \in W^{n-1/q}_q(\Gamma_j)$,
$\bT^n_jh = h$ on $\Gamma_j$ and $\|\bT^n_jh\|_{H^n_q(\Omega_j)}
\leq C\|h\|_{W^{n-1/q}_q(\Gamma_j)}$ with some constant
$C>0$ independent of $j \in \BN$. \\
\thetag2~ There exists an extension map $\bT^n_\Gamma : W^{n-1/q}_q(\Gamma)
\to H^n_q(\Omega)$ such that for $h \in W^{n-1/q}_q(\Gamma)$,
$\bT^n_\Gamma h = h$ on $\Gamma$ and $\|\bT^n_\Gamma h\|_{H^n_q(\Omega)}
\leq C\|h\|_{W^{n-1/q}_q(\Gamma)}$ with some constant $C > 0$. 
\end{prop}
\pf For a proof, see Proposition 9.5.4 in Shibata \cite{S4}. \qed
\begin{prop}\label{prop:5.4}
Let $1 < q < \infty$ and $n=2,3$ and let $\eta_j \in 
C^\infty_0(B^1_j)$ $(j \in \BN)$ with  
$\|\eta_j\|_{H^n_\infty(\BR^N)}\leq c_2$ for some 
constant $c_2$ independent of $j \in \BN$. Then, 
we have the following two assertions: \\
\thetag1~Let $f_j$ $(j \in \BN)$ be functions in $W^{n-1/q}_q(\Gamma_j)$  
satisfying the condition:
$\sum_{j=1}^\infty \|f_j\|_{W^{n-1/q}_q(\Gamma_j)}^q < \infty$, 
and 
then the infinite sum $\sum_{j=1}^\infty \eta_jf_j$ converges to some $f
\in W^{n-1/q}_q(\Gamma)$ strongly in $W^{n-1/q}_q(\Gamma)$ and
$$\|f\|_{W^{n-1/q}_q(\Gamma)} \leq C_q\{\sum_{j=1}^\infty 
\|f_j\|_{W^{n-1/q}_q(\Gamma_j)}^q\}^{1/q}.$$
\thetag2~For any $h \in W^{n-1/q}_q(\Gamma)$, 
$$\sum_{j=1}^\infty \|\eta_jh\|_{W^{n-1/q}_q(\Gamma_j)}^q
\leq C\|h\|_{W^{n-1/q}_q(\Gamma)}^q.
$$
\end{prop}
\pf For a proof, see Proposition 9.5.5 in Shibata \cite{S4}. \qed

\subsection{Parametrix of solutions of Eq. \eqref{rres:1.1}}
\label{subsec:5.6}

In this subsection, we construct a parametrix for Eq. \eqref{eq:7.0}.  
Let $\{\zeta^i_j\}_{j \in \BN}$ and 
$\{\tilde\zeta^i_j\}_{j\in\BN}$ $(i = 0,1)$ be sequences of $C^\infty_0$ functions given in Proposition
\ref{prop:lap}, 
and let $(\bff, d, \bh) \in Y_q(\Omega)$  (cf. \eqref{rdata:0}). 
Recall that    
$\Omega_j=\Phi_j(\HS)$ and  $\Gamma_j = 
\Phi_j(\BR^N_0)$.  Let 
\begin{align*}
\mu^i_j(x) &= \tilde\zeta^i_j(x)\mu(x)+(1-\tilde\zeta^i_j(x))\mu(x^i_j),
\\ 
\sigma_j(x) &= \tilde\zeta^1_j(x)\sigma(x)+(1-\tilde\zeta^1_j(x))
\sigma(x^1_j), \\
A_{\kappa, j}(x) &= \tilde\zeta^1_j(x)A_\kappa(x)
+(1-\tilde\zeta^1_j(x))A_\kappa(x^1_j).
\end{align*}
Notice that
$$\zeta^i_j\mu = \zeta^i_j\mu^i_j, \quad \zeta^1_j\sigma = \zeta^1_j\sigma_j,
\quad \zeta^1_jA_\kappa = \zeta^1_jA_{\kappa, j},$$
because $\tilde\zeta^1_j = 1$ on ${\rm supp}\, \zeta^1_j$. 
We consider the equations:
\begin{align}
\,\,\lambda \bu^0_j - \DV(\mu^0_j\bD(\bu^0_j) - K_{0j}(\bu^0_j)\bI) 
= \tilde\zeta^0_j\bff \quad\,\quad\text{in $\BR^N$};\,& \label{126}\\
\left\{\begin{aligned}
\lambda \bu^1_j - \DV(\mu^1_j\bD(\bu^1_j)-K_{1j}(\bu^1_j, h_j)\bI)
&= \tilde\zeta^1_j\bff &\quad&\text{in $\Omega_j$}, \\
\lambda h_j + A_{\kappa, j}\cdot\nabla_{\Gamma_j}h_j-\bn_j\cdot\bu_j
&=\tilde\zeta^1_jd&\quad&\text{on $\Gamma_j$}, \\
(\mu^1_j\bD(\bu^1_j)-K_{1j}(\bu^1_j, h_j)\bI)\bn_j - 
\sigma_j(\Delta_{\Gamma_j}h_j)\bn_j & =\tilde\zeta^1_j\bh
&\quad&\text{on $\Gamma_j$}.
\end{aligned}\right.& \label{127}
\end{align}
Here, for $\bu \in H^2_q(\BR^N)^N$, 
$K_{0j}(\bu) \in \hat H^1_q(\BR^N)$ 
denotes a unique solution of the weak Laplace equation:
\begin{equation}\label{131} (\nabla K_{0j}(\bu), \nabla\varphi)_{\BR^N}
= (\DV(\mu^0_j\bD(\bu))-\nabla\dv\bu,
\nabla\varphi)_{\BR^N}
\end{equation}
for any $\varphi \in \hat H^1_{q'}(\BR^N)$.
And, for $\bu \in H^2_q(\Omega_j)$ and 
$h \in H^3_q(\Omega_j)$, $K_{1j}(\bu, h_j)
\in H^1_q(\Omega_j) + \hat H^1_{q,0}(\Omega_j)$ denotes
 a unique solution of the weak Dirichlet
problem:
\begin{equation}\label{132} (\nabla K_{1j}(\bu, h), 
\nabla\varphi)_{\Omega_j}
= (\DV(\mu^1_j\bD(\bu))-\nabla\dv\bu,
\nabla\varphi)_{\Omega_j}
\end{equation}
for any $\varphi \in \hat H^1_{q',0}(\Omega_j)$, 
subject to 
$$K_{1j}(\bu, h) = <\mu^1_j\bD(\bu)\bn_j, \bn_j>
-\dv\bu - \sigma_j\Delta_{\Gamma_j}h
\quad \text{on $\Gamma_j$}.
$$ 
Moreover, we denote  the unit outer normal to $\Gamma_j$ 
by $\bn_j$, which is defined
on $\BR^N$ and satisfies the estimate:
$$\|\bn_j\|_{L_\infty(\BR^N)} \leq C, \quad 
\|\nabla \bn_j\|_{L_\infty(\BR^N)} \leq C_A, \quad 
\|\nabla^2\bn_j\|_{L_\infty(\BR^N)} \leq C_{M_2}.
$$
Let 
$\nabla_{\Gamma_j} 
= (\pd_1, \ldots, \pd_{N-1})$ with $\pd_j = \pd/\pd x_j$ for
$y = \Phi_j(x', 0) \in \Gamma_j$ and let   
$\Delta_{\Gamma_j}$ be the Laplace-Beltrami operator on $\Gamma_j$,
which have the form:
$$\Delta_{\Gamma_j}f = \Delta'f + \CD_{\Gamma_j}f 
\quad\text{on $\Phi^{-1}_j(\Gamma_j)$},$$
where $\Delta'f = \sum_{j=1}^{N-1}\pd_j^2f$ and $\CD_{\Gamma_j}f
= \sum_{k.\ell=1}^{N-1}a^j_{k\ell}\pd_k\pd_\ell f
+\sum_{k=1}^{N-1}a^j_k\pd_kf$,
and $a^j_{k\ell}$ and $a^j_k$ satisfy the following estimates:
\begin{align*}
&\|a^j_{k\ell}\|_{L_\infty(\BR^N)} \leq CM_1, 
\quad\|(\pd_1a^j_{k\ell}, \ldots, \pd_{N-1}a^j_{k\ell}, 
 a^j_k)\|_{L_\infty(\BR^N)} \leq C_A, \\
&\|(\pd_1a^j_{k\ell}, \ldots, \pd_{N-1}a^j_{k\ell}, 
 a^j_k)\|_{H^1_\infty(\BR^N)} \leq C_{M_2}.
\end{align*}
Notice that $\bn_j = \bn$ and $\Delta_{\Gamma_j} = \Delta_{\Gamma}$
on $\Gamma_j \cap B^1_j = \Gamma \cap B^1_j$. 
We know the existence of $K_{0j}(\bu^0_j) \in \hat H^1_q(\BR^N)
$ possessing the estimate:
\begin{equation}\label{pres:1}\|\nabla K_{0j}(\bu^0_j)\|_{L_q(\BR^N)} 
\leq C\|\nabla \bu^0_j\|_{H^1_q(\BR^N)}.
\end{equation}
Let $\rho$ be a function in $C^\infty_0(B_{r_0})$ such that 
$\int_{\BR^N}\rho\,dx=1$.  Below, this $\rho$ is fixed. 
Since $K_{0j}(\bu^0_j) + c$ also satisfy the variational equation \eqref{131}
for any constant $c$, we may assume that 
\begin{equation}\label{pres:2}
\int_{B^0_j} K_{0j}(\bu^0_j)\rho(x-x^0_j)\,dx = 0.
\end{equation}
Moreover, 
choosing $M_1 \in (0, 1)$ suitably small, 
we have the unique existence of solutions $K_{1j}(\bu^1_j, h_j)
\in H^1_q(\Omega_j) + \hat H^1_{q,0}(\Omega_j)$ 
of Eq.\eqref{132} possessing the estimates:
\begin{equation}\label{press:3}
\|\nabla K_{1j}(\bu^1_j, h_j)\|_{L_q(\Omega_j)} 
\leq C(\|\nabla \bu^1_j\|_{H^1_q(\Omega_j)}
+ \|h_j\|_{W^{3-1/q}_q(\Gamma_j)}). 
\end{equation}
Let $Y_q(\Omega_j)$ and $\CY_q(\Omega_j)$ be the spaces defined
in \eqref{rdata:0} replacing $\Omega$ by $\Omega_j$. 
By Theorem \ref{thm:4.1} and Theorem \ref{thm:b}, there exist constants
$M_1 \in (0, 1)$ and 
$\lambda_0 \geq 1$, which are independent of $j \in \BN$, 
 and operator families
\begin{align*}
\CS_{0j}(\lambda) &\in \Hol(\Sigma_{\epsilon, \lambda_0}, \CL(L_q(\BR^N)^N, 
H^2_q(\BR^N)^N)),\\
\CS_{1j}(\lambda) &\in \Hol(\Lambda_{\kappa, \lambda_0\gamma_\kappa}, 
\CL(\CY_q(\Omega_j), H^2_q(\Omega_j)^N)), \\
\CH_j(\lambda) &\in \Hol(\Lambda_{\kappa,\lambda_0\gamma_\kappa},
\CL(\CY_{q}(\Omega_j), H^3_q(\Omega_j)))
\end{align*}
such that for each $j \in \BN$,
Eq.\eqref{126} admits  a unique solution
$\bu^0_j = \CS_{0j}(\lambda)\tilde\zeta^0_j\bff$ and 
Eq.\eqref{127} admits unique solutions 
$\bu^1_j = \CS_{1j}(\lambda)\tilde\zeta^1_j\bF_\lambda(\bff, d, \bh)$
and $h_j = \CH_j(\lambda)\tilde\zeta^1_j\bF_\lambda(\bff, d, \bh)$, where 
$\bF_\lambda(\bff, d, \bh) = (\bff, d, \lambda^{1/2}\bh, \bh)$,
and $\tilde\zeta^1_jF_\lambda(\bff, d, \bh) 
= (\tilde\zeta^1_j\bff,  \tilde\zeta^1_jd,
\lambda^{1/2}\tilde\zeta^1_j\bh, \tilde\zeta^1_j\bh)$. 
Moreover, there exists a number $r_b > 0$ independent of 
$M_1$, $M_2$, and 
$j \in \BN$ such that 
\begin{equation}\label{136}\begin{split}
\CR_{\CL(L_q(\BR^N)^N, H^{2-k}_q(\BR^N)^N)}
(\{(\tau\pd_\tau)^\ell(\lambda^{k/2}\CS_{0j}(\lambda)) \mid
\lambda \in \Sigma_{\epsilon, \lambda_0}\}) & \leq r_b, \\
\CR_{\CL(\CY_q(\Omega_j), H^{2-k}_q(\Omega_j)^N)}
(\{(\tau\pd_\tau)^\ell(\lambda^{k/2}\CS_{1j}(\lambda)) \mid
\lambda \in \Lambda_{\epsilon, \lambda_0\gamma_\kappa}\}) & \leq r_b, \\
\CR_{\CL(\CY_q(\Omega_j), H^{3-n}_q(\Omega_j))}
(\{(\tau\pd_\tau)^\ell(\lambda^{n}\CH_j(\lambda)) \mid
\lambda \in \Lambda_{\kappa, \lambda_0\gamma_\kappa}\}) & \leq r_b,
\end{split}\end{equation}
for $\ell=0,1$, $j \in \BN$, $k=0,1,2$, and $n=0,1$. Notice that
$\lambda_0\gamma_\kappa \geq \lambda_0$. 
 
By \eqref{136}, we have
\begin{equation}\label{137}\begin{split}
&|\lambda|\|\bu^0_j\|_{L_q(\BR^N)}
+ |\lambda|^{1/2}\|\bu^0_j\|_{H^1_q(\BR^N)}
+ \|\bu^0_j\|_{H^2_q(\BR^N)} \leq r_b\|\zeta^0_j\bff\|_{L_q(\BR^N)},
\\
&|\lambda|\|\bu^1_j\|_{L_q(\Omega_j)}
+ |\lambda|^{1/2}\|\bu^1_j\|_{H^1_q(\Omega_j)}
+ \|\bu^1_j\|_{H^2_q(\Omega_j)}\\
&\quad
+ |\lambda|\|h_j\|_{H^2_q(\Omega_j)}
+\|h_j\|_{H^3_q(\Omega_j)}\\
& \leq r_b(\|\tilde\zeta^1_j\bff\|_{L_q(\Omega_j)}
+\|\tilde\zeta^1_jd\|_{W^{2-1/q}_q(\Gamma_j)}
+ |\lambda|^{1/2}\|\bh\|_{L_q(\Omega_j)}
+ \|\bh\|_{H^1_q(\Omega_j)})
\end{split}\end{equation}
for $\lambda \in \Sigma_{\kappa, \lambda_0\gamma_\kappa}$. 
Let 
\begin{equation}\label{138}
\bu = \sum_{i=0}^1\sum_{j=1}^\infty \zeta^i_j\bu^i_j,
\quad h = \sum_{j=1}^\infty \zeta^1_jh_j.
\end{equation}
Then, by \eqref{126}, \eqref{127}, \eqref{137}, Proposition \ref{prop:5.2}, and
Proposition \ref{prop:5.4}, we have 
$\bu \in H^2_q(\Omega)^N$, $h \in H^3_q(\Omega)$, and 
\begin{align*}
&|\lambda|\|\bu\|_{L_q(\Omega)} + |\lambda|^{1/2}\|\bu\|_{H^1_q(\Omega)}
+\|\bu\|_{H^2_q(\Omega)} + |\lambda|\|h\|_{H^2_q(\Omega)}
+ \|h\|_{H^3_q(\Omega)} \\
&\quad\leq C_qr_b(\|\bff\|_{L_q(\Omega)} + \|d\|_{W^{2-1/q}_q(\Gamma)}
+ |\lambda|^{1/2}\|\bh\|_{L_q(\Omega)} + \|\bh\|_{H^1_q(\Omega)}\}
\end{align*}
for $\lambda \in \Lambda_{\kappa, \lambda_0\gamma_\kappa}$.

Moreover, we have
\begin{equation}\label{140}\left\{\begin{aligned}
\lambda\bu - \DV(\mu\bD(\bu) - K(\bu, h)\bI) 
&= \bff - V^1(\lambda)(\bff, d, \bh)&\quad&
\text{in $\Omega$}, \\
\lambda h + A_\kappa\cdot\nabla'_\Gamma h - \bu\cdot\bn 
& = d - V^2_\kappa(\lambda)(\bff, d, \bh)
&\quad&\text{on $\Gamma$}, \\
(\mu\bD(\bu) - K(\bu, h)\bI -(\sigma\Delta_\Gamma h)\bI)\bn
&= \bh - V^3(\lambda)(\bff, d, \bh)
&\quad&\text{on $\Gamma$},
\end{aligned}\right.\end{equation}
where we have set
\begin{align*}
V^1(\lambda)(\bff, d, \bh)& = V^1_1(\lambda)(\bff, d, \bh)
+ V^1_2(\lambda)(\bff, d, \bh),\\
V^1_1(\lambda)(\bff, d, \bh)& = \sum_{i=0}^1\sum_{j=1}^\infty[
\DV(\mu(\bD(\zeta^i_j\bu^i_j)-\zeta^i_j\bD(\bu^i_j)))\\
&\quad
+ \DV(\zeta^i_j\mu^i_j\bD(\bu^i_j))-\zeta^i_j\DV(\mu^i_j\bD(\bu^i_j))],\\
V^1_2(\lambda)(\bff, d, \bh)& = \nabla K(\bu, h) - 
\sum_{j=1}^\infty \zeta^0_j\nabla K_{0j}(\bu^0_j)
-\sum_{j=1}^\infty \zeta^1_j\nabla K_{1j}(\bu^1_j, h_j),\\
V^2_\kappa(\lambda)(\bff, d, \bh)&
= \sum_{j=1}^\infty A_\kappa(x)\cdot((\nabla'_\Gamma\zeta^1_j)h_j), \\
V^3(\lambda)(\bff, d, \bh)& = V^3_1(\lambda)(\bff, d, \bh) - 
V^3_2(\lambda)(\bff, d, \bh) - V^3_3(\lambda)(\bff, d, \bh), \\
V^3_1(\lambda)(\bff, d, \bh) & =
\sum_{j=1}^\infty \mu(\bD(\zeta^1_j\bu^1_j)-\zeta^1_j\bD(\bu^1_j))\bn, \\
V^3_2(\lambda)(\bff, d, \bh) 
&= \{\sum_{j=1}^\infty\zeta^1_jK_{1j}(\bu^1_j, h_j)
-K(\bu, h)\}\bn, \\
V^3_3(\lambda)(\bff, d, \bh) & = 
\sum_{j=1}^\infty \sigma(\Delta_{\Gamma_j}(\zeta^1_jh^1_j)-
\zeta^1_j\Delta_{\Gamma_j}h^1_j).
\end{align*}

For $F = (F_1, F_2, F_3, F_4) \in \CY_q(\Omega)$, we define operators
$\CA_p(\lambda)$ and $\CB_p(\lambda)$ acting on $F$ by 
\begin{equation}\label{144*}\begin{aligned}
\CA_p(\lambda)F &= \sum_{j=1}^\infty
\zeta^0_j\CS^0_j(\lambda)\tilde\zeta^0_jF_1
+ \sum_{j=1}^\infty\zeta^1_j\CS_{1j}(\lambda)\tilde\zeta^1_jF, \\
\CB_p(\lambda)F  &= \sum_{j=1}^\infty \zeta^1_j\CH_j(\lambda)\tilde\zeta^1_jF.
\end{aligned}\end{equation}
Then, by Proposition \ref{prop:5.2} and \eqref{136}, we have 
$\bu = \CA_p(\lambda)\bF_\lambda(\bff, d, \bh)$, $h = \CH_p(\lambda)
\bF_\lambda(\bff, d, \bh)$, and 
\begin{equation}\label{144}\begin{split}
&\CA_p(\lambda)\in \Hol(\Lambda_{\kappa, \lambda_1}, 
\CL(\CY_q(\Omega), H^2_q(\Omega)^N)), \\ 
&\CB_p(\lambda) \in \Hol(\Lambda_{\kappa, \lambda_1}, 
\CL(\CY_q(\Omega), H^3_q(\Omega))), \\
&\CR_{\CL(\CY_q(\Omega), H^{2-j}_q(\Omega)^N)}
(\{(\tau\pd_\tau)^\ell(\lambda^{j/2}\CA_p(\lambda))
\mid \lambda \in \Lambda_{\kappa, \lambda_1}\}) \\
&\phantom{\CR_{\CL(\CY_q(\Omega), H^{3-k}_q(\Omega))}
(\{(\tau\pd_\tau)^\ell(\lambda^k\CB_p(\lambda))}
\leq 
(C + C_{M_2}\lambda_1^{-1/2})r_b, \\
&\CR_{\CL(\CY_q(\Omega), H^{3-k}_q(\Omega))}
(\{(\tau\pd_\tau)^\ell(\lambda^k\CB_p(\lambda))
\mid \lambda \in \Lambda_{\kappa, \lambda_1}\})\\
&\phantom{\CR_{\CL(\CY_q(\Omega), H^{3-k}_q(\Omega))}
(\{(\tau\pd_\tau)^\ell(\lambda^k\CB_p(\lambda))}
 \leq 
(C + C_{M_2}\lambda_1^{-1})r_b
\end{split}\end{equation}
for $\ell=0,1$, $j=0,1,2$, and $k=0,1$ for any 
$\lambda_1 \geq \lambda_0\gamma_\kappa$.

\subsection{Estimates of the remainder terms}\label{subsec:n5.3}
For $F = (F_1, F_2, F_3, F_4) \in \CY_q(\Omega)$, let
\begin{align*}
\CV^1(\lambda)F& = \CV^1_1(\lambda)F
+ \CV^1_2(\lambda)F,\\
\CV^1_1(\lambda)F& = \sum_{j=1}^\infty[
\DV(\mu(\bD(\zeta^0_j\CS_{0j}(\lambda)\tilde\zeta^0_jF_1)
-\zeta^0_j\bD(\CS_{0j}(\lambda)\tilde\zeta^0_jF_1))) \\
&+ \DV(\zeta^0_j\mu^0_j\bD(\CS_{0j}(\lambda)\tilde\zeta^0_jF_1))
-\zeta^0_j\DV(\mu^0_j\bD(\CS_{0j}(\lambda)\tilde\zeta^0_jF_1))]\\
& +\sum_{j=1}^\infty[
\DV(\mu(\bD(\zeta^1_j\CS_{1j}(\lambda)\tilde\zeta^1_jF)
-\zeta^1_j\bD(\CS_{1j}(\lambda)\tilde\zeta^1_jF))) \\
&+ \DV(\zeta^1_j\mu^1_j\bD(\CS_{1j}(\lambda)\tilde\zeta^1_jF))
-\zeta^1_j\DV(\mu^1_j\bD(\CS_{1j}(\lambda)\tilde\zeta^1_jF))],
\\
\CV^1_2(\lambda)F& = \nabla K(\CA_p(\lambda)F, \CB_p(\lambda)F) - 
\sum_{j=1}^\infty \zeta^0_j\nabla K_{0j}(\CS_{0j}(\lambda)\tilde\zeta^0_jF_1)
\\
&-\sum_{j=1}^\infty \zeta^1_j\nabla K_{1j}(\CS_{1j}(\lambda)\tilde\zeta^1_jF, 
\CH_j(\lambda)\tilde\zeta^1_jF),\\
\CV^2_\kappa(\lambda)F&
= \sum_{j=1}^\infty A_\kappa(x)\cdot((\nabla'_\Gamma\zeta^1_j)
\CH_j(\lambda)\tilde\zeta^1_jF), \\
\CV^3(\lambda)F& = \CV^3_1(\lambda)F + 
\CV^3_2(\lambda)F + \CV^3_3(\lambda)F, \\
\CV^3_1(\lambda)F & =
\sum_{j=1}^\infty \mu(\bD(\zeta^1_j\CS_{1j}(\lambda)\tilde\zeta^1_jF)
-\zeta^1_j\bD(\zeta^1_j\CS_{1j}(\lambda)\tilde\zeta^1_jF))\bn \\
\CV^3_2(\lambda)F& 
= \{\sum_{j=1}^\infty\zeta^1_jK_{1j}(\CS_{1j}(\lambda)\tilde\zeta^1_jF,
\CH_j(\lambda)\tilde\zeta^1_jF)
-K(\CA_p(\lambda)F,  \CB_p(\lambda)F)\}\bn, \\
\CV^3_3(\lambda)F & = 
\sum_{j=1}^\infty\sigma(\Delta_{\Gamma_j}(\zeta^1_j
\CH_j(\lambda)\tilde\zeta^1_jF)-
\zeta^1_j\Delta_{\Gamma_j}(\CH_j(\lambda)\tilde\zeta^1_jF))
\end{align*}
Notice that $\CV^2_\kappa(\lambda)F = 0$ for $\kappa = 0$.

Let 
\begin{align*}
V(\lambda)(\bff, d, \bh) &= (V^1(\lambda)(\bff, d, \bh), 
V^2(\lambda)(\bff, d, \bh), V^3(\lambda)(\bff, d, \bh)), \\
\CV(\lambda)F &= (\CV^1(\lambda)F, \CV^2_\kappa(\lambda)F, \CV^3(\lambda)F).
\end{align*}
Since $\bu^0_j = \CS_{0j}(\lambda)\tilde\zeta^0_j\bff$,
$\bu^1_j = \CS_{1j}(\lambda)\tilde\zeta^1_j\bF_\lambda(\bff, d, \bh)$, 
and $h_j = \CH_j(\lambda)\zeta^1_j\bF_\lambda\tilde(\bff, d, \bh)$,
we have 
\begin{equation}\label{145}
V(\lambda)(\bff, d, \bh) = \CV(\lambda)
\bF_\lambda(\bff, d, \bh).
\end{equation}
In what follows, we shall prove that
\begin{equation}\label{146}\begin{aligned}
&\CR_{\CL(\CY_q(\Omega))}(\{(\tau\pd_\tau)^\ell
(\bF_\lambda\CV(\lambda)) \mid 
\lambda \in \Lambda_{\kappa, \tilde\lambda_0} \}) \\
&\quad\leq C_{q}r_b(\epsilon + C_{M_2, \epsilon}
(\tilde\lambda_0^{-1}\gamma_\kappa + \tilde\lambda_0^{-1/2}))
\end{aligned}
\end{equation}
for $\ell=0,1$ and $\tilde\lambda_0 \geq \lambda_0\gamma_\kappa$, 
where $\gamma_\kappa$ is the number given in Theorem \ref{thm:rbdd:1.2}. 

To prove \eqref{146}, we use Proposition \ref{prop:lap}, 
Proposition \ref{prop:4.1}, Propositions \ref{prop:5.1}--\ref{prop:5.4}, 
\eqref{res:4}, \eqref{res:2}, \eqref{R.5.1*}, \eqref{R.5.1}, \eqref{D.1*}
and \eqref{136}. 
 In the following, $\tilde\lambda_0$ is any number 
such that $\tilde\lambda_0 \geq \lambda_0\gamma_\kappa$.  
We start with the following
estimate of $\CV^1_1(\lambda)$: 
\begin{equation}\label{est:6.1}
\CR_{\CL(\CY_q(\Omega), L_q(\Omega)^N)}(\{(\tau\pd_\tau)^\ell
\CV^1_1(\lambda) \mid \lambda \in \Sigma_{\kappa, \tilde\lambda_0}\})
\leq C_{M_2}r_b\tilde\lambda_0^{-1/2}
\end{equation}
for $\ell=0,1$.
In fact, since $D_{\ell, m}(\zeta^i_j\bu) - \zeta^i_jD_{\ell m}(\bu)
= (\pd_\ell\zeta^i_j)u_m + (\pd_m\zeta^i_j)u_\ell$, 
and $\dv(\zeta^i_j \bu) - \zeta^i_j\dv\bu
= \sum_{k=1}^N(\pd_k\zeta^i_j)u_k$, 
for any $n \in \BN$, $\{\lambda_\ell\}_{\ell=1}^n
\subset \Lambda_{\kappa, \tilde\lambda_0}^n$, 
and $\{F_\ell = (F_{1\ell}, F_{2\ell}, F_{3\ell},
F_{4\ell})\}_{\ell=1}^n \subset \CY_q(\Omega)^n$, we have
\begin{align*}
&\int^1_0\|\sum_{\ell=1}^nr_\ell(u)\CV^1_1(\lambda_\ell)F_\ell\|_{L_q(\Omega)}^q\,du
\\
&\leq 
C_q^qM_2^q\sum_{j=1}^\infty\{
\int^1_0\|\sum_{\ell=1}^nr_\ell(u)\CS_{0j}(\lambda_\ell)
\tilde\zeta^0_jF_{1\ell}\|_{H^1_q(\BR^N)}^q\,du \\
&\phantom{\leq 
C_q^qM_2^q\tilde\lambda_0^{-q/2}\sum_{j=1}^\infty\{
aaaaa}
+ \int^1_0\|\sum_{\ell=1}^nr_\ell(u)S_{1j}(\lambda_\ell)\tilde\zeta^1_jF_\ell
\|_{H^1_q(\Omega_j)}^q\,du \}\\
&\leq 
C_q^qM_2^q\tilde\lambda_0^{-q/2}\sum_{j=1}^\infty\{
\int^1_0\|\sum_{\ell=1}^nr_\ell(u)\lambda_\ell^{1/2}\CS_{0j}(\lambda_\ell)
\tilde\zeta^0_jF_{1\ell}\|_{H^1_q(\BR^N)}^q\,du \\
&\phantom{\leq 
C_q^qM_2^q\tilde\lambda_0^{-q/2}\sum_{j=1}^\infty\{
aaaaa}
+ \int^1_0\|\sum_{\ell=1}^nr_\ell(u)\lambda_\ell^{1/2}
S_{1j}(\lambda_\ell)\tilde\zeta^1_jF_\ell
\|_{H^1_q(\Omega_j)}^q\,du \}\\
&\leq C_q^qM_2^q\tilde\lambda_0^{-q/2}r_b^q\sum_{j=1}^\infty\{
\int^1_0\|\sum_{\ell=1}^nr_\ell(u)\tilde\zeta^0_jF_{1\ell}
\|_{L_q(\BR^N)}^q\,du\\
&\phantom{\leq 
C_q^qM_2^q\tilde\lambda_0^{-q/2}\sum_{j=1}^\infty\{
aaaaa}
+
\int^1_0\|\sum_{\ell=1}^nr_\ell(u)\tilde\zeta^1_jF_{\ell}
\|_{\CY_q(\Omega_j)}^q\,du\} \\
&\leq C_q^{2q}M_2^q\tilde\lambda_0^{-q/2}r_b^q\int^1_0
\|\sum_{\ell=1}^n r_\ell(u)F_\ell\|_{L_q(\Omega)}^q\,du.
\end{align*}
This shows that
$$\CR_{\CL(\CY_q(\Omega), L_q(\Omega)^N)}(\{
\CV^1_1(\lambda) \mid \lambda \in \Sigma_{\kappa, \tilde\lambda_0}\})
\leq C_{M_2}r_b\tilde\lambda_0^{-1/2}.
$$
Analogously, we can show that
$$\CR_{\CL(\CY_q(\Omega), L_q(\Omega)^N)}(\{
\tau\pd_\tau\CV^1_1(\lambda) \mid \lambda 
\in \Sigma_{\kappa, \tilde\lambda_0}\})
\leq C_{M_2}r_b\tilde\lambda_0^{-1/2},
$$
and therefore we have \eqref{est:6.1}. 

For $r \in (N, \infty)$ and $q \in (1, \infty)$, by the extension of 
functions defined on $\Gamma_j$ to $\Omega_j$ and Sobolev's imbedding
theorem, we have
$$\|ab\|_{W^{2-1/q}_q(\Gamma_j)} \leq C_{q, r, K}
\|a\|_{H^2_q(\Omega_j)}\|b\|_{W^{2-1/q}_q(\Gamma_j)}
$$
for any $a \in H^2_r(\Omega)$ and $b \in W^{2-1/q}_q(\Gamma_j)$.
 Applying this inequality, we have
\begin{align*}
\|A_\kappa\cdot((\nabla'_\Gamma\zeta^1_j)\CH_j(\lambda)\tilde\zeta^1_jF)
\|_{W^{2-1/q}_q(\Gamma_j)}
\leq C{q,r}M_2m_3\kappa^{-b}\|\CH_j(\lambda)\tilde\zeta^1_jF
\|_{H^2_q(\Omega_j)}
\end{align*}
for $\kappa \in (0, 1)$. Thus, employing the same argument as in the proof of
\eqref{est:6.1}, we have
\begin{align*}
\CR_{\CL(\CY_q(\Omega), W^{2-1/q}_q(\Gamma))}
(\{(\tau\pd_\tau)^\ell\CV^2_\kappa(\lambda) \mid
\lambda \in \Lambda_{\kappa, \tilde\lambda_0}\})
&\leq C_{M_2}r_b\tilde\lambda_0^{-1}\kappa^{-b}
\end{align*}
for $\ell=0,1$ and $\kappa \in (0, 1)$.

Employing the same argument as in the proof of \eqref{est:6.1}, 
we also have 
\begin{align*}
\CR_{\CL(\CY_q(\Omega), L_q(\Gamma)^N)}
(\{(\tau\pd_\tau)^\ell(\lambda^{1/2}\CV^3_m(\lambda)) \mid
\lambda \in \Lambda_{\kappa, \tilde\lambda_0}\})
&\leq C_{M_2}r_b\tilde\lambda_0^{-1/2}, \\
\CR_{\CL(\CY_q(\Omega), H^1_q(\Gamma)^N)}
(\{(\tau\pd_\tau)^\ell\CV^3_m(\lambda) \mid
\lambda \in \Lambda_{\kappa, \tilde\lambda_0}\})
&\leq C_{M_2}r_b\tilde\lambda_0^{-1/2},
\end{align*}
for $\ell=0,1$,  $m=1$ and $3$. Noting that $\mu\zeta^1_j = \mu^1_j\zeta^1_j$
and $\sigma\zeta^1_j = \sigma^1_j\zeta^1_j$, we have  
\begin{equation}\label{diff:1}\begin{split}
&\sum_{j=1}^\infty\zeta^1_jK_{1j}(\CS_{1j}(\lambda)\tilde\zeta^1_jF, 
\CH_j(\lambda)\tilde\zeta^1_jF)
- K(\CA_p(\lambda)F, \CB_p(\lambda)F) \\
&\quad 
= \sum_{j=1}^\infty\mu
<\zeta^1_j\bD(\CS_{1j}(\lambda)\tilde\zeta^1_jF)
-\bD(\zeta^1_j\CS_{1j}(\lambda)\tilde\zeta^1_jF), \bn> \\
&\quad
-\sum_{j=1}^\infty\{\zeta^1_j\dv\CS_{1j}(\lambda)\tilde\zeta^1_jF - 
\dv(\zeta^1_j\CS_{1j}(\lambda)\tilde\zeta^1_jF)\}\\
&\quad - \sum_{j=1}^\infty\sigma
\{\zeta^1_j\Delta_{\Gamma_j}
(\CH_j(\lambda)\tilde\zeta^1_jF)-
\Delta_{\Gamma_j}(\zeta^1_j\CH_j(\lambda)\tilde\zeta^1_jF)\}
\end{split}\end{equation}
on $\Gamma$, where we have used $\Delta_\Gamma=\Delta_{\Gamma_j}$ and 
$\bn=\bn_j$ on $\Gamma_j\cap B^1_j$.
Employing the same argument as 
in the proof of \eqref{est:6.1}, we have
\begin{align*}
\CR_{\CL(\CY_q(\Omega), L_q(\Gamma)^N)}
(\{(\tau\pd_\tau)^\ell(\lambda^{1/2}\CV^3_2(\lambda)) \mid
\lambda \in \Lambda_{\kappa, \tilde\lambda_0}\})
&\leq C_{M_2}r_b\tilde\lambda_0^{-1/2}; \\
\CR_{\CL(\CY_q(\Omega), H^1_q(\Gamma)^N)}
(\{(\tau\pd_\tau)^\ell\CV^3_2(\lambda) \mid
\lambda \in \Lambda_{\kappa, \tilde\lambda_0}\})
&\leq C_{M_2}r_b\tilde\lambda_0^{-1/2}
\end{align*}
for $\ell=0, 1$.

The final task is to prove that
\begin{equation}\label{est:6.2} 
\CR_{\CL(\CY_q(\Omega), L_q(\Gamma)^N)}
(\{(\tau\pd_\tau)^\ell \CV^1_2(\lambda) \mid
\lambda \in \Lambda_{\kappa, \tilde\lambda_0}\})
\leq C_{q,r}(\epsilon + C_{M_2, \epsilon}\tilde\lambda_0^{-1/2})r_b
\end{equation}
for $\ell=0,1$. For this purpose, we use Lemma \ref{lem:5.7} and
the following lemma. 
\begin{lem}\label{lem:63.3.6} Let $1 < q < \infty$. For 
$\bu \in H^2_q(\BR^N)$, let $K_{0j}(\bu)$ be a unique solution of the weak
Laplace equation \eqref{131} satisfying \eqref{pres:2}.  Then, we have
\begin{equation}\label{63.43}
\|K_{0j}(\bu)\|_{L_q(B^0_j)} \leq C\|\nabla\bu\|_{L_q(\BR^N)}.
\end{equation}
\end{lem}
\pf Let $\rho$ be the same function in \eqref{pres:2}. Let 
$\psi$ be any function in $C^\infty_0(B^0_j)$ and we set 
$$\tilde\psi(x) = \psi(x) - \rho(x-x^0_j)\int_{\BR^N}\psi(y)\,dy.
$$
Then, 
\begin{equation}\label{63.46}
\tilde\psi \in C^\infty_0(B^0_j), \quad \int_{\BR^N}\tilde\psi\,dx = 0, \quad
\|\tilde\psi\|_{L_{q'}(B^0_j)} \leq C_{q'}\|\psi\|_{_{q'}(B^0_j)}.
\end{equation}
Moreover, 
\begin{equation}\label{63.47}
\tilde\psi \in \hat H^1_q(\BR^N)^* = \hat H^{-1}_{q'}(\BR^N),
\quad
\|\tilde\psi\|_{\hat H^{-1}_{q'}(\BR^N)} \leq C_{q'}\|\psi\|_{L_{q'}(\BR^N)}.
\end{equation}
In fact,  by Lemma \ref{lem:5.7}, for any $\varphi \in \hat H^1_{q'}(\BR^N)$,
there exists a constant $e_j$ for which 
$$\|\varphi-e_j\|_{L_q(B^0_j)} \leq c_{q}\|\nabla\varphi\|_{L_q(B^0_j)}.
$$
Thus, by \eqref{63.46}, we have
\begin{align*}
|(\tilde\psi, \varphi)_{\BR^N}| = |(\tilde\psi, \varphi-e_j)_{\BR^N}|
&\leq \|\tilde\psi\|_{L_{q'}(B^0_j)}\|\varphi-e_j\|_{L_q(B^0_j)}\\
&\leq C_q\|\tilde\psi\|_{_{q'}(B^0_j)}\|\nabla\varphi\|_{L_q(B^0_j)},
\end{align*}
which yields \eqref{63.47}. Let $\Psi$ be a function in 
$\hat H^1_{q'}(\BR^N)$ such that $\nabla\Psi \in H^1_{q'}(\BR^N)^N$, 
\begin{equation}\label{63.48*}\begin{split}
(\nabla\Psi, \nabla\theta)_{\BR^N} &= (\tilde\psi, \theta)_{\BR^N}
\quad\text{for any $\theta \in \hat H^1_q(\BR^N)$}, \\
\|\nabla\Psi\|_{H^1_{q'}(\BR^N)} &\leq C(\|\tilde\psi\|_{L_{q'}(\BR^N)}
+ \|\tilde\psi\|_{\hat H^{-1}_{q'}(\BR^N)}).
\end{split}\end{equation}
By \eqref{63.46} and \eqref{63.47}, we have
\begin{equation}\label{63.48}
\|\nabla\Psi\|_{H^1_{q'}(\BR^N)} \leq C_{q'}\|\psi\|_{L_{q'}(\BR^N)}.
\end{equation}
By \eqref{pres:2}, \eqref{63.48*}, and the divergence theorem of 
Gau\ss, we have
\begin{align*}
&(K_{0j}(\bu), \psi)_{\BR^N} = (K_{0j}(\bu), \tilde\psi)_{\BR^N}
= (\nabla K_{0j}(\bu), \nabla\Psi)_{\BR^N} \\
&=
(\DV(\mu^0_j\bD(\bu))-\nabla\dv\bu, \nabla\Psi)_{\BR^N}
 = -(\mu^0_j\bD(\bu), \nabla^2\Psi)_{\BR^N} +
(\dv\bu, \Delta\Psi)_{\BR^N},
\end{align*}
and therefore by \eqref{63.48} 
$$|(K_{0j}(\bu), \psi)_{\BR^N}| \leq C
\|\nabla \bu\|_{L_q(\BR^N)}\|\psi\|_{L_{q'}(\BR^N)},
$$
which proves \eqref{63.43}.  This completes the proof of
Lemma \ref{lem:63.3.6}. \qed
\begin{lem}\label{lem:63.3.4} Let $1 < q < \infty$.  
For $\bu \in H^2_q(\Omega_j)$
and $h \in H^3_q(\Omega_j)$, let $K_{1j}(\bu, h)
\in H^1_q(\Omega_q) + \hat H^1_{q,0}(\Omega_j)$ be a unique solution of the
weak Dirichlet problem \eqref{132}. Then, we have
\begin{align*}
\|K_{1j}(\bu, h)\|_{L_q(\Omega_j\cap B^1_j)}  
&\leq C(\|\nabla\bu\|_{L_q(\Omega_j)}
+ \|h\|_{H^2_q(\Omega_j)} \\
&\quad
+ \|\nabla^2\bu\|_{L_q(\Omega_j)}^{1/q}\|\nabla\bu\|_{L_q(\Omega_j)}^{1-1/q}
+ \|h\|_{H^3_q(\Omega_j)}^{1/q}\|h\|_{H^2(\Omega_j)}^{1-1/q}).
\end{align*}
Here, the constant $C$ depends on $q$ and $C_A$. 
\end{lem}
\begin{remark} By Young's inequality, we have
\begin{equation}\label{young:1}\begin{aligned}
\|K_{1j}(\bu, h)\|_{L_q(\Omega_j\cap B^1_j)} 
&\leq \epsilon(\|\nabla^2\bu\|_{L_q(\Omega_j)}+
\|h\|_{H^3_q(\Omega_j)}) \\
&\quad + C_{\epsilon}(
\|\nabla\bu\|_{L_q(\Omega_j)}
+ \|h\|_{H^2_q(\Omega_j)})
\end{aligned}\end{equation}
for any $\epsilon \in (0, 1)$ with some constant $C_{\epsilon, q}$
depending on $\epsilon$ and $q$.
\end{remark}
\pf For a proof, see Lemma 3.4 in Shibata \cite{S0}. \qed
\vskip0.5pc
To prove \eqref{est:6.2}, we divide $\CV^1_2(\lambda)$ into two parts
as $\CV^1_2(\lambda) = \nabla \CV^1_{21}(\lambda) + \CV^1_{22}(\lambda)$, where
\begin{align*}
\CV^1_{21}(\lambda)F & = K(\CA_p(\lambda)F, \CB_p(\lambda)F)
-\sum_{j=1}^\infty\zeta^0_j K_{0j}(\CS_{0j}(\lambda)\tilde\zeta^0_jF_1)
\\
&\quad
- \sum_{j=1}^\infty\zeta^1_j K_{1j}(\CS_{1j}(\lambda)\tilde\zeta^0_jF,
\CH_j(\lambda)\tilde\zeta^1_jF), \\
\CV^1_{22}(\lambda)F & =
\sum_{j=1}^\infty(\nabla\zeta^0_j) K_{0j}(\CS_{0j}(\lambda)\tilde\zeta^0_jF_1)
\\
&\quad
+ \sum_{j=1}^\infty\nabla(\zeta^1_j)K_{1j}(\CS_{1j}(\lambda)\tilde\zeta^0_jF,
\CH_j(\lambda)\tilde\zeta^1_jF).
\end{align*}
By  \eqref{wd:3.1.2}, \eqref{131}, and \eqref{132}, for any 
$\varphi \in \hat H^1_{q'0}(\Omega)$ we have
$(\nabla \CV^1_{21}F, \nabla\varphi)_\Omega = I - II$, where 
\allowdisplaybreaks{
\begin{align*}
I & = (\DV(\mu\bD(\CA_p(\lambda)F))
- \nabla\dv(\CA_p(\lambda)F), \nabla\varphi)_\Omega, \\
II & = \sum_{j=1}^\infty((\nabla\zeta^0_j)K_{0j}(
\CS_{0j}(\lambda)\tilde\zeta^0_j F_1), \nabla(\varphi - e_j))_\Omega \\
&+ \sum_{j=1}^\infty((\nabla\zeta^1_j)K_{1j}
(\CS_{1j}(\lambda)\tilde\zeta^1_jF, 
\CH_j(\lambda)\tilde\zeta^1_jF), \nabla\varphi)_\Omega\\
& + \sum_{j=1}^\infty(\nabla K_{0j}(\CS_{0j}(\lambda)\tilde\zeta^0_jF_1), 
\nabla(\zeta^0_j(\varphi-e_j)))_\Omega\\
&+ \sum_{j=1}^\infty(\nabla K_{1j}(\CS_{1j}(\lambda)\tilde\zeta^1_jF,
 \CH_j(\lambda)\tilde\zeta^1_jF), \nabla(\zeta^1_j\varphi))_\Omega\\
&-\sum_{j=1}^\infty((\nabla\zeta^0_j)\nabla K_{0j}(\CS_{0j}(\lambda)
\tilde\zeta^0_j F_1), \varphi - e_j)_\Omega \\
&- \sum_{j=1}^\infty((\nabla\zeta^1_j)\nabla K_{1j}
(\CS_{1j}(\lambda)\tilde\zeta^1_jF, \CH_{j}(\lambda)\tilde\zeta^1_jF),
\varphi)_\Omega.
\end{align*}
}
Here and in the following, $e_j=c^0_j(\varphi)$ 
are constants given in Lemma \ref{lem:5.7}.
By the definition \eqref{144*}, we have
\allowdisplaybreaks{
\begin{align*}
I &= \sum_{j=1}^\infty(\DV(\mu\bD(\zeta^0_j\CS_{0j}(\lambda)
\tilde\zeta^0_j F_1))- \nabla\dv(\zeta^0_j\CS_{0j}(\lambda)
\tzeta^0_jF_1), \nabla\varphi)_{\BR^N}\\
&+ \sum_{j=1}^\infty(\DV(\mu\bD(\zeta^1_j\CS_{1j}(\lambda)\tzeta^1_jF))
-\nabla\dv(\zeta^1_j\CS_{1j}(\lambda)\tzeta^1_jF), \nabla\varphi)_\Omega\\
& = \sum_{j=1}^\infty(\zeta^0_j\DV(\mu^0_j\bD(\CS_{0j}(\lambda)
\tzeta^0_jF_1)) - \zeta^0_j\nabla\dv(\CS_{0j}(\lambda)\tzeta^0_jF_1), 
\nabla\varphi)_{\BR^N} \\
& + \sum_{j=1}^\infty(\zeta^1_j\DV(\mu^1_j\bD(\CS_{1j}(\lambda)
\tzeta^1_jF)) - \zeta^1_j\nabla\dv(\CS_{1j}(\lambda)\tzeta^1_jF), 
\nabla\varphi)_{\Omega} + III,
\end{align*}
}
where 
\allowdisplaybreaks{
\begin{align}
III &
= \sum_{j=1}^\infty(\DV(\mu\bD(\zeta^0_j\CS_{0j}(\lambda)
\tzeta^0_jF_1))-\zeta^0_j\DV(\mu\bD(\CS_{0j}(\lambda)\tzeta^0_jF_1)),
\nabla\varphi)_{\BR^N} \nonumber\\
& -\sum_{j=1}^\infty(\nabla\dv(\zeta^0_j\CS_{0j}(\lambda)\tzeta^0_jF_1)
-\zeta^0_j\nabla\dv(\CS_{0j}(\lambda)\tzeta^0_jF_1), 
\nabla\varphi)_\Omega \nonumber\\
&+ \sum_{j=1}^\infty(\DV(\mu\bD(\zeta^1_j\CS_{1j}(\lambda)
\tzeta^1_jF))-\zeta^1_j\DV(\mu\bD(\CS_{1j}(\lambda)\tzeta^1_jF)),
\nabla\varphi)_{\BR^N} \nonumber\\
& -\sum_{j=1}^\infty(\nabla\dv(\zeta^1_j\CS_{1j}(\lambda)\tzeta^1_jF)
-\zeta^1_j\nabla\dv(\CS_{1j}(\lambda)\tzeta^1_jF), 
\nabla\varphi)_\Omega.
\label{weak:0}
\end{align}}
Since $\zeta^0_j(\varphi-e_j) \in \hat H^1_{q',0}(\BR^N)$,
and $\zeta^1_j\varphi \in \hat H^1_{q',0}(\Omega_j)$, 
by \eqref{131} and \eqref{132}, we have
\allowdisplaybreaks{
\begin{align*}
II & = \sum_{j=1}^\infty(\zeta^0_j\DV(\mu^0_j\bD(\CS_{0j}(\lambda)
\tzeta^0_jF_1)) - \zeta^0_j\nabla\dv(\CS_{0j}(\lambda)\tzeta^0_jF_1), 
\nabla\varphi)_{\BR^N} \\
& + \sum_{j=1}^\infty(\zeta^1_j\DV(\mu^1_j\bD(\CS_{1j}(\lambda)
\tzeta^1_jF)) - \zeta^1_j\nabla\dv(\CS_{1j}(\lambda)\tzeta^1_jF), 
\nabla\varphi)_{\Omega} + IV,
\end{align*}}
where we have set 
\allowdisplaybreaks{
\begin{align*}
&IV=\\
&\sum_{j=1}^\infty\Bigl\{2(K_{0j}(\CS_{0j}(\lambda)
\tilde\zeta^0_jF_1)(\nabla\zeta^0_j), \nabla\varphi)_\Omega 
+(K_{0j}(\CS_{0j}(\lambda)\tilde\zeta^0_jF_1)(\Delta\zeta^0_j), 
\varphi-e_j)_\Omega \\
&- (\mu^0_j\bD(\CS_{0j}(\lambda)\tilde\zeta^0_jF_1):(\nabla^2\zeta^0_j), 
\varphi-e_j)_\Omega 
-(\mu^0_j\bD(\CS_{0j}(\lambda)\tilde\zeta^0_jF_1)(\nabla\zeta^0_j), 
\nabla\varphi)_\Omega \\
&+(\dv (\CS_{0j}(\lambda)\tilde\zeta^0_jF_1)(\Delta\zeta^0_j), 
\varphi-e_j)_\Omega 
+ (\dv (\CS_{0j}(\lambda)\tilde\zeta^0_jF_1)(\nabla\zeta^0_j),
\nabla\varphi)_\Omega \\
&+2(K_{1j}(\CS_{1j}(\lambda)\tilde\zeta^1_jF, \CH_j(\lambda)\tilde\zeta^1_jF)
(\nabla\zeta^1_j), 
\nabla\varphi)_\Omega \\
&+ (K_{1j}(\CS_{1j}(\lambda)\tilde\zeta^1_jF, \CH_j(\lambda)\tilde\zeta^1_jF)
(\Delta\zeta^1_j), \varphi)_\Omega \\
& - (\mu^1_j\bD(\CS_{1j}(\lambda)\tilde\zeta^1_jF_1):(\nabla^2\zeta^1_j), 
\varphi)_\Omega 
-(\mu^1_j\bD(\CS_{1j}(\lambda)\tilde\zeta^1_jF_1)(\nabla\zeta^1_j), 
\nabla\varphi)_\Omega \\
&+(\dv (\CS_{1j}(\lambda)\tilde\zeta^1_jF_1)(\Delta\zeta^1_j), 
\varphi)_\Omega 
+ (\dv (\CS_{1j}(\lambda)\tilde\zeta^1_jF_1)(\nabla\zeta^1_j),
\nabla\varphi)_\Omega\Bigr\}.
\end{align*}}
Thus, we have
\begin{equation}\label{weak:1}
(\nabla\CV^1_{21}(\lambda)F, \nabla\varphi)_\Omega = III + IV.
\end{equation}
We let define operators $\CL(\lambda)$ and $\CM(\lambda)$ acting on
$F \in \CY_q(\Omega)$ by the following formulas:
\allowdisplaybreaks{
\begin{align*}
&\CL(\lambda)F  =  \sum_{j=1}^\infty(\DV(\mu\bD(\zeta^0_j\CS_{0j}(\lambda)
\tzeta^0_jF_1))-\zeta^0_j\DV(\mu\bD(\CS_{0j}(\lambda)\tzeta^0_jF_1))\\
&\quad-\sum_{j=1}^\infty(\nabla\dv(\zeta^0_j\CS_{0j}(\lambda)\tzeta^0_jF_1)
-\zeta^0_j\nabla\dv(\CS_{0j}(\lambda)\tzeta^0_jF_1))\\
&\quad+ \sum_{j=1}^\infty(\DV(\mu\bD(\zeta^1_j\CS_{1j}(\lambda)
\tzeta^1_jF))-\zeta^1_j\DV(\mu\bD(\CS_{1j}(\lambda)\tzeta^1_jF)) \\
&\quad-\sum_{j=1}^\infty(\nabla\dv(\zeta^1_j\CS_{1j}(\lambda)\tzeta^1_jF)
-\zeta^1_j\nabla\dv(\CS_{1j}(\lambda)\tzeta^1_jF))\\
&\quad+ 2\sum_{j=1}^\infty(\nabla\zeta^0_j)K_{0j}(
\CS_{0j}(\lambda)\tilde\zeta^0_j F_1)
+ 2\sum_{j=1}^\infty(\nabla\zeta^1_j)
K_{1j}(\CS_{1j}(\lambda)\tilde\zeta^1_jF, 
\CH_j(\lambda)\tilde\zeta^1_jF)\\
&\quad-\sum_{j=1}^\infty\mu^0_j\bD(\CS_{0j}(\lambda)\tzeta^0_jF_1)
(\nabla\zeta^0_j)
+\sum_{j=1}^\infty
\dv(\CS_{0j}(\lambda)\tzeta^0_jF_1)(\nabla\zeta^0_j)\\
&\quad-\sum_{j=1}^\infty\mu^1_j\bD(\CS_{1j}(\lambda)\tzeta^1_jF)
(\nabla\zeta^1_j)
 +\sum_{j=1}^\infty
\dv(\CS_{1j}(\lambda)\tzeta^1_jF_1)(\nabla\zeta^1_j); 
\end{align*}
\begin{align*}
&<<\CM(\lambda)F, \varphi>> 
= -\sum_{j=1}^\infty(\mu^0_j\bD(\CS_{0j}(\lambda)\tzeta^0_jF_1):
(\nabla^2\zeta^0_j), \varphi-e_j)_\Omega\\
&\quad +\sum_{j=1}^\infty
(\dv(\CS_{0j}(\lambda)\tzeta^0_jF_1)(\Delta\zeta^0_j),
\varphi-e_j)_\Omega \\
&\quad +\sum_{j=1}^\infty(K_{0j}(\CS_{0j}(\lambda)
\tzeta^0_j F_1)(\Delta\zeta^0_j), \varphi-e_j)_{\Omega}\\ 
&\quad -\sum_{j=1}^\infty(\mu^1_j\bD(\CS_{1j}(\lambda)\tzeta^1_jF):
(\nabla^2\zeta^1_j), \varphi)_\Omega 
+\sum_{j=1}^\infty
(\dv(\CS_{1j}(\lambda)\tzeta^1_jF)(\Delta\zeta^1_j),
\varphi)_\Omega \\
&\quad + \sum_{j=1}^\infty(K_{1j}(\CS_{1j}(\lambda)\tilde\zeta^1_jF, 
\CH_{j}(\lambda)\tilde\zeta^1_jF)(\Delta\zeta^1_j),
\varphi)_\Omega. 
\end{align*} 
}
Here and in the following, 
$\hat W^{-1}_{q,0}(\Omega)$ denotes the dual space of
 $\hat H^1_{q',0}(\Omega)$,
and $<<\cdot, \cdot>>$ denotes the duality
between $\hat W^{-1}_{q,0}(\Omega)$ and $\hat H^1_{q',0}(\Omega)$. 

Moreover, by \eqref{wd:3.1.2} and \eqref{132}, for $x \in \Gamma$ we have
\allowdisplaybreaks{
\begin{align*}
\CV^1_{21}(\lambda)F & = <\mu\bD(\CA_p(\lambda)F)\bn, \bn> -
\sigma\Delta_\Gamma\CB_p(\lambda)F - \dv \CA_p(\lambda)F \\
& - \sum_{j=1}^\infty \zeta^1_j\{<(\mu(x^1_j)\bD(\CS_{1j}(\lambda)\tzeta^1_jF)
\bn_j, \bn_j> - \sigma(x^1_j)\Delta_{\Gamma_j}\CH_j(\lambda)\tzeta^1_jF
\\
&\phantom{- \sum_{j=1}^\infty \zeta^1_j\{<(\mu(x^1_j)\bD(\CS_{1j}(\lambda)\tzeta^1_jF)
\bn_j, \bn_j>}
-\dv \CS_{1j}(\lambda)\tzeta^1_jF\}
\\
& = \sum_{j=1}^\infty<\mu\bD(\zeta^1_j\CS_{1j}(\lambda)\tzeta^1_j F)\bn, \bn> 
-\sum_{j=1}^\infty \sigma\Delta_\Gamma(\zeta^1_j
\CH_j(\lambda)\tzeta^1_jF) \\
&\phantom{
= \sum_{j=1}^\infty<\mu\bD(\zeta^1_j
\CS_{1j}(\lambda)\tzeta^1_j F)\bn, \bn>
}
-\sum_{j=1}^\infty\dv(\zeta^1_j\CS_{1j}(\lambda)\tzeta^1_j F) \\
& - \sum_{j=1}^\infty \zeta^1_j\{<\mu(x^1_j)\bD(\CS_{1j}(\lambda)\tzeta^1_jF)
\bn_j, \bn_j> - \sigma(x^1_j)\Delta_{\Gamma_j}\CH_j(\lambda)\tzeta^1_jF
\\
&\phantom{= \sum_{j=1}^\infty<\mu\bD
(\zeta^1_j\CS_{1j}(\lambda)\tzeta^1_j F)\bn, \bn>}
-\dv \CS_{1j}(\lambda)\tzeta^1_jF\}.
\end{align*}
}
Thus, we define an operator $\CL_b(\lambda)$  acting on $F \in \CY_q(\Omega)$
 by letting 
\allowdisplaybreaks{
\begin{align*}
\CL_b(\lambda)F  &= \sum_{j=1}^\infty[
<\mu(x)(\bD(\zeta^1_j\CS_{1j}(\lambda)\tzeta^1_j F) - 
\zeta^1_j\bD(\CS_{1j}(\lambda)\tzeta^1_jF))\bn, \bn> \\
&- \sigma(\Delta_\Gamma(\zeta^1_j\CH_j(\lambda)\tzeta^1_jF) 
- 
\zeta^1_j\Delta_\Gamma\CH_j(\lambda)\tzeta^1_jF) 
-(\nabla\zeta^1_j)\CS_{1j}(\lambda)\tzeta^1_jF], 
\end{align*}
}
and then, $\CV^1_{21}(\lambda)F = \CL_b(\lambda)F$ on $\Gamma$.

We now prove the $\CR$ boundedness of operator families  $\CL(\lambda)$, 
$\CM(\lambda)$ and $\CL_b(\lambda)$. We first prove that
\begin{equation}\label{63.49} 
\CR_{\CL(\CY_q(\Omega), \hat W^{-1}_{q,0}(\Omega))}(\{(\tau\pd_\tau)^\ell
\CM(\lambda) \mid \lambda \in \Lambda_{\kappa, \tilde\lambda_0}\})
\leq (\epsilon + C_{q, \epsilon}\tilde\lambda_0^{-1/2})r_b
\end{equation}
for $\ell=0,1$. In fact, if we set 
\begin{align*}
<<\CM^0_j(\lambda)F, \varphi>> 
&= -(\mu(x^0_j)\bD(\CS_{0j}(\lambda)\tzeta^0_jF_1):
(\nabla^2\zeta^0_j), \varphi-e_j)_\Omega\\
&\quad +
(\dv(\CS_{0j}(\lambda)\tzeta^0_jF_1)(\Delta\zeta^0_j),
\varphi-e_j)_\Omega\\
&\quad +(K_{0j}(\CS_{0j}(\lambda)
\tzeta^0_j F_1)(\Delta\zeta^0_j), \varphi-e_j)_{\Omega}, \\
<<\CM^1_j(\lambda)F, \varphi>> & = 
-(\mu(x^1_j)\bD(\CS_{1j}(\lambda)\tzeta^1_jF):
(\nabla^2\zeta^1_j), \varphi)_\Omega \\
&\quad +
(\dv(\CS_{1j}(\lambda)\tzeta^1_jF)(\Delta\zeta^1_j),
\varphi)_\Omega \\
&\quad
+ (K_{1j}(\CS_{1j}(\lambda)\tilde\zeta^1_jF, \CH_{j}(\lambda)\tilde\zeta^1_jF)
(\Delta\zeta^1_j),
\varphi)_\Omega, 
\end{align*}
then, 
by Lemma \ref{lem:5.7}, 
Lemma \ref{lem:63.3.6} and \eqref{young:1}, we have
\begin{align*}
&|<<\CM^0_j(\lambda)F, \varphi>>| \leq C_{M_2}\|\nabla\CS_{0j}(\lambda)
\tzeta^0_jF\|_{L_q(\BR^N)}\|\nabla\varphi\|_{L_{q'}(B^0_j)}, \\
&|<<\CM^1_j(\lambda)F, \varphi>>| \leq \{\epsilon(\|\CS_{1j}
\tzeta^1_jF\|_{H^2_q(\Omega_j)}
+ \|\CH_j(\lambda)\tzeta^1_jF\|_{H^3_q(\Omega_j)}) \\
&\quad
+ C_{\epsilon, M_2}(\|\CS_{1j}
\tzeta^1_jF\|_{H^1_q(\Omega_j)}
+ \|\CH_j(\lambda)\tzeta^1_jF\|_{H^2_q(\Omega_j)})\}
\|\nabla\varphi\|_{L_q(B^1_j\cap \Omega)}.
\end{align*}
By \eqref{D.1}, we have
$$\sum_{j=1}^\infty\|\nabla\varphi\|_{L_{q'}(B^0_j)}^{q'}
+ \sum_{j=1}^\infty\|\nabla\varphi\|_{L_{q'}(B^1_j \cap \Omega)}^{q'}
\leq C_{q'}\|\nabla\varphi\|_{L_{q'}(\Omega)}^{q'}
$$
for any $\varphi \in \hat H^1_{q',0}(\Omega)$.  By \eqref{136}, \eqref{D.1}, 
and Proposition \eqref{prop:5.4}, we have 
\begin{align*}
&\sum_{j=1}^\infty\|\nabla\CS_{0j}(\lambda)
\tzeta^0_jF_1\|_{L_q(\BR^N)}^q
+ \sum_{j=1}^\infty(\|\CS_{1j}(\lambda)\tzeta^1_jF\|_{H^2_q(\Omega_j)}^q
+ \|\CH_j(\lambda)\tzeta^1_jF\|_{H^3_q(\Omega_j)}^q) \\
&\quad\leq r_b^q(\sum_{j=1}^\infty\|\tzeta^0_jF_1\|_{L_q(\BR^N)}^q
+ \sum_{j=1}^\infty\|\tzeta^1_jF\|_{\CY_q(\Omega_j)}^q)
\leq r_b^qC_q\|F\|_{\CY_q(\Omega)}^q < \infty.
\end{align*}
Thus, by Proposition \ref{prop:5.1}, 
$\CM(\lambda)F = \sum_{j=1}^\infty \CM^0_j(\lambda)F 
+ \sum_{j=1}^\infty \CM^1_j(\lambda)F$ exists in 
$\hat W^{-1}_{q,0}(\Omega)$ for any $F \in \CY_q(\Omega)$
and 
\begin{align*}
&\|\CM(\lambda)F\|_{\hat W^{-1}_{q,0}(\Omega)}^q \\
&\quad 
\leq C_{M_2}^q\sum_{j=1}^\infty\|\nabla\CS_{0j}(\lambda)\tzeta^0_jF_1
\|_{L_q(\BR^N)}^q
+ \epsilon^q\sum_{j=1}^\infty\|\CS_{1j}(\lambda)
\tzeta^1_jF\|_{H^2_q(\Omega_j)}^q \\
&\quad 
+ \epsilon^q\sum_{j=1}^\infty\|\CH_j(\lambda)\tzeta^1_jF\|_{H^3_q(\Omega_j)}^q
+ C^q_{\epsilon, M_2}\sum_{j=1}^\infty\|\CS_{1j}(\lambda)\tzeta^1_jF
\|_{H^1_q(\Omega_j)}^q\\
&\quad
+ C^q_{\epsilon, M_2}\sum_{j=1}^\infty\|\CH_j(\lambda)\tzeta^1_jF
\|_{H^1_q(\Omega_j)}^q.
\end{align*} 
Analogously, by Proposition \ref{prop:5.1} we have 
\begin{align*}
&\|\sum_{\ell=1}^nr_\ell(u)\CM(\lambda_\ell)F_\ell
\|_{\hat W^{-1}_{q,0}(\Omega)}^q 
\leq C_{M_2}^q\sum_{j=1}^\infty\|\sum_{\ell=1}^nr_\ell(u)
\nabla\CS_{0j}(\lambda_\ell)\tzeta^0_jF_{1\ell}
\|_{L_q(\BR^N)}^q \\
&+ \epsilon^q\sum_{j=1}^\infty\{\|
\sum_{\ell=1}^nr_\ell(u)
\CS_{1j}(\lambda_\ell)\tzeta^1_jF_\ell\|_{H^2_q(\Omega_j)}^q
+
\sum_{\ell=1}^nr_\ell(u)
\CH_j(\lambda_\ell)\tzeta^1_jF_\ell\|_{H^3_q(\Omega_j)}^q
\}\\
& + C^q_{\epsilon, M_2}\sum_{j=1}^\infty\{\|
\sum_{\ell=1}^nr_\ell(u)\CS_{1j}(\lambda_\ell)\tzeta^1_jF_\ell
\|_{H^1_q(\Omega_j)}^q + \|
\sum_{\ell=1}^nr_\ell(u)\CH_j(\lambda_\ell)\tzeta^1_jF_\ell
\|_{H^1_q(\Omega_j)}^q\}.
\end{align*}
Noting that $\Omega \cap B^1_j = \Omega_j \cap B^1_j$,
by \eqref{136}, \eqref{D.1*}, Proposition \ref{prop:5.4}, 
and Proposition \eqref{prop:4.1}, we have
\begin{align*}
&\int^1_0\|\sum_{\ell=1}^nr_\ell(u)\CM(\lambda_\ell)F_\ell
\|_{\hat W^{-1}_{q,0}(\Omega)}^q\,du\\
&\quad
\leq C_{M_2}^q\tilde\lambda_0^{-q/2}r_b^q\int^1_0
\sum_{j=1}^\infty\|\sum_{\ell=1}^nr_\ell(u)
\tzeta^0_jF_{1\ell}
\|_{L_q(\BR^N)}^q\,du \\
&\quad\quad+ \epsilon^qr_b^q\int^1_0\sum_{j=1}^\infty\|
\sum_{\ell=1}^nr_\ell(u)\tzeta^1_jF_\ell\|_{\CY_q(\Omega_j)}^q\,du
\\
&\quad\quad
+ C^q_{\epsilon, M_2}\tilde\lambda_0^{-q/2}r_b^q \int^1_0 
\sum_{j=1}^\infty\|
\sum_{\ell=1}^nr_\ell(u)\tzeta^1_jF_\ell
\|_{\CY_q(\Omega_j)}^q\,du\\
&\quad \leq C_q(\epsilon^q + C_{\epsilon, M_2}^q\tilde\lambda_0^{-q/2})r_b^q
\int^1_0\|\sum_{\ell=1}^n r_\ell(u)F_\ell
\|_{\CY_q(\Omega)}^q\,du,
\end{align*}
which shows \eqref{63.49}.  Analogously, we can prove
\begin{equation}\label{63.50}\begin{split} 
\CR_{\CL(\CY_q(\Omega), L_q(\Omega)^N)}(\{(\tau\pd_\tau)^\ell
\CL(\lambda) \mid \lambda \in \Lambda_{\kappa, \tilde\lambda_0}\})
\leq C_{M_2}r_b\tilde\lambda_0^{-1/2}, \\
\CR_{\CL(\CY_q(\Omega), H^1_q(\Omega)^N)}(\{(\tau\pd_\tau)^\ell
\CL_b(\lambda) \mid \lambda \in \Lambda_{\kappa, \tilde\lambda_0}\})
\leq  C_{M_2}r_b\tilde\lambda_0^{-1/2}
\end{split}\end{equation}
for $\ell=0,1$.

We now use the following lemma.
\begin{lem} \label{lem:63.3.5} Let $1 < q < \infty$.  Then, 
there exists a linear map $\CE$ from $\hat W^{-1}_{q,0}(\Omega)$
into $L_q(\Omega)^N$ such that for any $F \in \hat W^{-1}_{q,0}(\Omega)$, 
$\|\CE(F)\|_{L_q(\Omega)} \leq C\|F\|_{\hat W^{-1}_{q,0}(\Omega)}$
and 
$$<F, \varphi> = (\CE(F), \nabla\varphi)_\Omega \quad
\text{for all $\varphi \in \hat H^1_{q',0}(\Omega)$}.
$$
\end{lem}
\pf The lemma follows from the Hahn-Banach theorem by indentifying
$\hat H^1_{q',0}(\Omega)$ with a closed subspace of $L_{q'}(\Omega)^N$ 
via the mapping: $\varphi \mapsto \nabla\varphi$. \qed \vskip0.5pc
Applying Lemma \ref{lem:63.3.5} and using \eqref{weak:1}
and \eqref{63.49}, we have 
\begin{equation}\label{63.51}
(\nabla\CV^1_{21}(\lambda)F, \nabla\varphi)_\Omega = (\CL(\lambda)F + 
\CE(\CM(\lambda)F), \nabla\varphi)_\Omega \quad\text{for all
$\varphi \in \hat H^1_{q',0}(\Omega)$}, 
\end{equation}
subject to $\CV^1_{21}(\lambda)F = \CL_b(\lambda)F$ on $\Gamma$, and 
\begin{equation}\label{63.52}
\CR_{\CL(\CY_q(\Omega), L_q(\Omega)^N)}
(\{(\tau\pd_\tau)^\ell \CE\circ\CM(\lambda) \mid 
\lambda \in \Sigma_{\sigma, \tilde\lambda_0}\}) \leq 
C(\epsilon + C_{q,\epsilon}\tilde\lambda_0^{-1/2})r_b,
\end{equation}
where $\CE\circ\CM(\lambda)$ denotes a bounded
 linear operator family acting on
$F$ by $\CE\circ\CM(\lambda)F = \CE(\CM(\lambda)F)$.  
By  Remark \ref{rem:3.2}, 
we have $\CV^1_{21}(\lambda)F = \CL_b(\lambda)F + \CK(\CL(\lambda)F + 
\CE(\CM(\lambda)F)-\nabla\CL_b(\lambda)F)$, and so by \eqref{63.50}
and \eqref{63.52}, we see that $\nabla \CV^1_{21}(\lambda) 
\in \Hol(\Sigma_{\sigma, \tilde\lambda_0}, 
\CL(\CY_q(\Omega), L_q(\Omega)^N))$ and
\begin{equation}\label{63.53}\begin{aligned}
&\CR_{\CL(\CY_q(\Omega), L_q(\Omega)^N)}
(\{(\tau\pd_\tau)^\ell \nabla \CV^1_{21}(\lambda) 
\mid \lambda \in \Sigma_{\sigma,
\tilde\lambda_0}\}) \\
&\quad \leq 
C_q(\epsilon  + C_{M_2, \epsilon}\tilde\lambda_0^{-1/2})r_b
\quad\text{for $\ell=0, 1$}.
\end{aligned}\end{equation}

Finally, by Lemma \ref{lem:63.3.6}, \eqref{young:1}, 
\eqref{D.1*}, and Proposition \ref{prop:5.1}, we have
\begin{gather*}
\CV^1_{22}(\lambda) \in \Hol(\Sigma_{\sigma, \tilde\lambda_0}, 
\CL(\CY_q(\Omega), L_q(\Omega)^N)), \\
\CR_{\CL(\CY_q(\Omega), L_q(\Omega)^N)}
(\{(\tau\pd_\tau)^\ell \CV^1_{21}(\lambda) \mid \lambda \in \Sigma_{\sigma,
\tilde\lambda_0}\}) \leq C_q(\epsilon +  
C_{\epsilon}\tilde\lambda_0^{-1/2})r_b
\end{gather*}
for $\ell=0,1$, which, combined with \eqref{63.53}
and the formula: $\CV^1_2(\lambda) = \nabla\CV^1_{21}(\lambda)
+ \CV^2_{22}(\lambda)$, leads to \eqref{est:6.2}. 


\subsection{Proof of Theorem \ref{main:thm3}, Existence part
for Eq. \eqref{eq:7.0}}\label{newsubsec:6.4}

Choosing  $\epsilon$ so small that $C_{q,r} r_b \epsilon \leq 1/4$, 
 and $\tilde\lambda_0$ so large
that 
\begin{equation}\label{147*}
C_qr_bC_{M_2, \epsilon}(\tilde\lambda_0^{-1}\gamma_\kappa
+ \tilde\lambda_0^{-1/2}) \leq 1/4
\end{equation}
in \eqref{146}, we have 
\begin{equation}\label{147} 
\CR_{\CL(\CY_q(\Omega))}(\{(\tau\pd_\tau)^\ell \bF_\lambda\CV(\lambda)
\mid \lambda \in \Lambda_{\kappa, \tilde\lambda_0}\}) \leq 1/2
\end{equation}
for $\ell=0,1$. Let $\lambda_*$ be a large number 
for which  
$\lambda_* \geq
(8C_qr_bC_{M_2,\epsilon})^2$, and then setting 
$\tilde\lambda_0 = \lambda_*\gamma_\kappa$,
we have \eqref{147*}. By \eqref{147*}, $(\bI-\bF_\lambda\CV(\lambda))^{-1}
= \sum_{j=1}^\infty (\bF_\lambda\CV(\lambda))^j$ exists in 
$\Hol(\Lambda_{\kappa,\tilde\lambda_0}, \CL(\CY_q(\Omega)))$ and 
\begin{equation}\label{148} \CR_{\CL(\CY_q(\Omega))}
(\{(\tau\pd_\tau)^\ell(\bI-\bF_\lambda\CV(\lambda))^{-1} \mid
\lambda \in \Lambda_{\kappa, \tilde\lambda_0}\}) \leq 4
\end{equation}
for $\ell=0,1$. Moreover, by \eqref{145} and \eqref{147} 
\begin{equation}\label{151}
\|\bF_\lambda V(\lambda)(\bff, d, \bh)\|_{\CY_q(\Omega)}
\leq (1/2)\|\bF_\lambda(\bff, d, \bh)\|_{\CY_q(\Omega)}.
\end{equation}
Since $\|\bF_\lambda(\bff, d, \bh)\|_{\CY_q(\Omega)}$ gives an
equivalent norm in $Y_q(\Omega)$ for $\lambda\not=0$, by \eqref{151} 
$(\bI-V(\lambda))^{-1} = \sum_{j=0}^\infty V(\lambda)^j$ exists in 
$\CL(Y_q(\Omega))$.  Since $\bu = \CA_p(\lambda)\bF_\lambda(\bff, d, \bh)$
and $h = \CB_p(\lambda)\bF_\lambda(\bff, d, \bh)$ satisfy 
Eq. \eqref{140}, 
setting
$$\bv=\CA_p(\lambda)\bF_\lambda(\bI-V(\lambda))^{-1}(\bff, d, \bh),
\quad \rho = \CB_p(\lambda)(\lambda)\bF_\lambda(\bI - V(\lambda))^{-1}
(\bff, d, \bh),
$$
we see that $\bv \in H^2_q(\Omega)^N$, 
$\rho \in H^3_q(\Omega)$ and $\bv$ and $\rho$ satisfy the equations:
\begin{equation}\label{150}\left\{\begin{aligned}
\lambda\bv - \DV(\mu\bD(\bv) - K(\bv, \rho)\bI) 
&= \bff &\quad&
\text{in $\Omega$}, \\
\lambda \rho + A_\kappa\cdot\nabla'_\Gamma \rho - \bv\cdot\bn & = d 
&\quad&\text{on $\Gamma$}, \\
(\mu\bD(\bv) - K(\bv, \rho)\bI -(\sigma\Delta_\Gamma\rho)\bI)\bn
&= \bh
&\quad&\text{on $\Gamma$},
\end{aligned}\right.\end{equation}
Moreover, by \eqref{145} we have 
$\bF_\lambda(\bI - V(\lambda))^{-1} = 
(\bI - \bF_\lambda\CV(\lambda))^{-1}\bF_\lambda$.  Thus, setting
$$\CA_r(\lambda) = \CA_p(\lambda)(\bI - \bF_\lambda \CV(\lambda))^{-1},
\quad
\CH_r(\lambda) = \CB_p(\lambda)(\bI - \bF_\lambda\CV(\lambda))^{-1}
$$
we see that $\bv = \CA_r(\lambda)\bF_\lambda(\bff, d, \bh)$ and 
$\rho = \CH_r(\lambda)\bF_\lambda(\bff, d, \bh)$ are solutions of 
Eq.\eqref{150}.  Since we may assume that 
$\lambda_*\gamma_\kappa \geq \lambda_1$ in \eqref{144}, 
by \eqref{144} and  \eqref{148}, 
we have
\begin{equation}\label{5.8.150}\begin{aligned}
\CR_{\CL(\CY_q(\Omega), H^{2-j}_q(\Omega)^N)}
(\{(\tau\pd_\tau)^\ell(\lambda^{j/2}\CA_r(\lambda)) \mid \lambda 
\in \Lambda_{\kappa, \lambda_*\gamma_\kappa}\}) & \leq C_qr_b, \\
\CR_{\CL(\CY_q(\Omega), H^{3-k}_q(\Omega))}
(\{(\tau\pd_\tau)^\ell(\lambda^{k}\CH_r(\lambda)) \mid \lambda 
\in \Lambda_{\kappa, \lambda_*\gamma_\kappa}\}) & \leq C_qr_b, 
\end{aligned}\end{equation}
for $\ell=0,1$, $j=0,1,2$ and $k=0,1$. 
This completes the 
proof of the existence part of Theorem \ref{main:thm3} for Eq. \eqref{eq:7.0}.
\subsection{A proof of the existence part of Theorem \ref{main:thm3}}
\label{subsec:new5}

We now prove the existence of $\CR$-bounded solution operators
of Eq. \eqref{rres:1.2}. Let $\CA_r(\lambda)$ and $\CH_r(\lambda)$ be the 
operators constructed in the previous subsection. 
Let $(\bff, d, \bh) \in Y_q(\Omega)$. 
Let $\bu = \CA_r(\lambda)\bF_\lambda(\bff, d, \bh)$ and 
$h = \CH_r(\lambda)\bF_\lambda(\bff, d, \bh)$, 
where, $\bF(\bff, d, \bh) 
= (\bff, d, \lambda^{1/2}\bh, \bh) \in \CY_q(\Omega)$ for 
$(\bff, d, \bh) \in Y_q(\Omega)$.
Then $\bu$ and $h$ satisfy the equations:
\begin{equation}\label{rres:5.1.2} 
\left\{\begin{aligned}
\lambda \bu - \DV(\mu\bD(\bu, h)-K(\bu, h)\bI) & = \bff 
&\quad&\text{in $\Omega$}, \\
\lambda h + A_\kappa\cdot\nabla'_\Gamma h
- \bn\cdot\bu + \CF_1\bu& = d  + \CF_1\bu&\quad&\text{on $\Gamma$}, \\
(\mu\bD(\bu) - K(\bu, h)\bI)\bn - (\CF_2h+\sigma\Delta_\Gamma h)\bn
&=\bh - (\CF_2h)\bn&\quad&\text{on $\Gamma$}.
\end{aligned}\right.\end{equation}
Let
\begin{align*}
&V_{R1}(\lambda)(\bff, d, \bh) = -\CF_1\CA_r(\lambda)\bF_\lambda(\bff, d, \bh),
 \quad \CV_{R1}(\lambda)F = -\CF_1\CA_r(\lambda)F, \\
&V_{R2}(\lambda)(\bff, d, \bh) 
= \CF_2\CH_r(\lambda)\bF_\lambda(\bff, d, \bh), \quad
\CV_{R2}(\lambda)F = \CF_2\CH_r(\lambda)F
\end{align*}
with $\bF_\lambda(\bff, d, \bh) = (\bff, d, \lambda^{1/2}\bh, \bh)$
and $F = (F_1, F_2, F_3, F_4) \in \CY_q(\Omega)$. 
Notice that  
\begin{align*}
\CF_1\bu &= -V_{R1}(\lambda)(\bff, d, \bh)
= -\CV_{R1}(\lambda)\bF_\lambda(\bff, d, \bh), \\
\CF_2h &=
 V_{R2}(\lambda)(\bff, d, \bh)
= \CV_{R2}(\lambda)\bF_\lambda(\bff, d, \bh). 
\end{align*}
Let $V_R(\lambda) = (0, V_{R1}(\lambda), V_{R2}(\lambda))$ and 
$\CV_R(\lambda) = (0, \CV_{R1}(\lambda), \CV_{R2}(\lambda))$, and then
\begin{equation}\label{5.8.150}
V_R(\lambda) = \CV_R(\lambda)\bF_\lambda.
\end{equation}  
By \eqref{5.8.150} and \eqref{cond-low}, we have
$$\CR_{\CL(\CY_q(\Omega))}(\{(\tau\pd_\tau)^\ell \bF_\lambda\CV_R(\lambda)
\mid \lambda \in \Lambda_{\kappa, \lambda_2}) \leq 
M_0(\lambda_2^{-1/2} + \lambda_2^{-1})r_b
\quad(\ell=0,1)
$$
for any $\lambda_2 > \lambda_*\gamma_\kappa$. 
Choosing $\lambda_1$ so large,  we have
\begin{equation}\label{5.8.151}
\CR_{\CL(\CY_q(\Omega))}(\{(\tau\pd_\tau)^\ell \bF_\lambda\CV_R(\lambda)
\mid \lambda \in \Lambda_{\kappa, \lambda_2}) \leq 1/2
\quad(\ell=0,1).
\end{equation}
By \eqref{5.8.150} and \eqref{5.8.151}, 
$$\|\bF_\lambda V_R(\lambda)(\bff, d, \bh)\|_{\CY_q(\Omega)}
\leq (1/2)\|\bF_\lambda(\bff, d, \bh)\|_{\CY_q(\Omega)}.
$$
Since $\|\bF_\lambda(\bff, d, \bh)\|_{\CY_q(\Omega)}$ gives an
equivalent norm in $Y_q(\Omega)$ for $\lambda\not=0$, by \eqref{151} 
$(\bI-V_R(\lambda))^{-1} = \sum_{j=0}^\infty V_R(\lambda)^j$ exists in 
$\CL(Y_q(\Omega))$.  Since $\bu = \CA_r(\lambda)\bF_\lambda(\bff, d, \bh)$
and $h = \CH_r(\lambda)\bF_\lambda(\bff, d, \bh)$ 
satisfy Eq. \eqref{rres:5.1.2}, 
setting
$$\bv=\CA_r(\lambda)\bF_\lambda(\bI-V_R(\lambda))^{-1}(\bff, d, \bh),
\quad \rho = \CH_r(\lambda)(\lambda)\bF_\lambda(\bI - V_R(\lambda))^{-1}
(\bff, d, \bh),
$$
we see that $\bv \in H^2_q(\Omega)^N$, 
$\rho \in H^3_q(\Omega)$ and $\bv$ and $\rho$ satisfy the equations:
\begin{equation}\label{5.8.152}\left\{\begin{aligned}
\lambda\bv - \DV(\mu\bD(\bv) - K(\bv, \rho)\bI) 
&= \bff &\quad&
\text{in $\Omega$}, \\
\lambda \rho + A_\kappa\cdot\nabla'_\Gamma \rho - \bv\cdot\bn 
+ \CF_1\bv& = d 
&\quad&\text{on $\Gamma$}, \\
(\mu\bD(\bv) - K(\bv, \rho)\bI -((\CF_2\rho+\sigma\Delta_\Gamma)\rho)\bI)\bn
&= \bh
&\quad&\text{on $\Gamma$},
\end{aligned}\right.\end{equation}
Moreover, by \eqref{5.8.150} we have 
$\bF_\lambda(\bI - V_R(\lambda))^{-1} = 
(\bI - \bF_\lambda\CV_R(\lambda))^{-1}\bF_\lambda$.  Thus, setting
$$\tilde\CA_r(\lambda) = \CA_r(\lambda)
(\bI - \bF_\lambda \CV(\lambda))^{-1},
\quad
\tilde\CH_r(\lambda) = \CH_r(\lambda)(\bI - \bF_\lambda\CV(\lambda))^{-1}
$$
we see that $\bv = \tilde\CA_r(\lambda)\bF_\lambda(\bff, d, \bh)$ and 
$\rho = \tilde\CH_r(\lambda)\bF_\lambda(\bff, d, \bh)$ are solutions of 
Eq.\eqref{5.8.152}.  And, by \eqref{5.8.150} 
and \eqref{5.8.151},  there exists a large number $\lambda_{**}
\geq \lambda_*$ for which 
\begin{align*}
\CR_{\CL(\CY_q(\Omega), H^{2-j}_q(\Omega)^N)}
(\{(\tau\pd_\tau)^\ell(\lambda^{j/2}\tilde\CA_r(\lambda)) \mid \lambda 
\in \Lambda_{\kappa, \lambda_{**}\gamma_\kappa}\}) & \leq C_qr_b, \\
\CR_{\CL(\CY_q(\Omega), H^{3-k}_q(\Omega))}
(\{(\tau\pd_\tau)^\ell(\lambda^{k}\tilde\CH_r(\lambda)) \mid \lambda 
\in \Lambda_{\kappa, \lambda_{**}\gamma_\kappa}\}) & \leq C_qr_b, 
\end{align*}
for $\ell=0,1$, $j=0,1,2$ and $k=0,1$. 
This completes the 
proof of the existence part of Theorem \ref{main:thm3}.

\subsection{Uniqueness}\label{subsec:5.9}

In this subsection, we shall prove the uniqueness part of  
Theorem \ref{main:thm0}.  The uniqueness part of Theorem \ref{main:thm3} 
 will be proved in the next subsection. 

Let $\bu \in H^2_q(\Omega)^N$ satisfy the homogeneous
equations:
\begin{equation}\label{uni:1}
\lambda\bu - \DV(\mu\bD(\bu) - K_0(\bu)\bI) = 0\,\,\,
\text{in $\Omega$}, \quad 
(\mu\bD(\bu) - K_0(\bu)\bI)\bn = 0\,\,\,
\text{on $\Gamma$}. 
\end{equation}
We shall prove that $\bu=0$ below. 
Let $\lambda_0$ be a large positive number such that
for any $\lambda \in \Sigma_{\epsilon, \lambda_0}$  
the existence part of Theorem \ref{main:thm0}
holds for $q' = q/(q-1)$. 
Let $J_{q}(\Omega)$ be a solenoidal space defined in 
\eqref{sol:1} and let $\bg$ be any element in $J_{q'}(\Omega)$.
Let $\bv \in H^2_{q'}(\Omega)^N$ be a solution of  
the equations:
\begin{equation}\label{uni:2}
\bar\lambda\bv-\DV(\mu\bD(\bv) - K_0(\bv)\bI)=\bg\,\,\,
\text{in $\Omega$}, \quad 
(\mu\bD(\bv) - K_0(\bv)\bI)\bn=0\,\,\,\text{on $\Gamma$}.
\end{equation}
We first observe that $\bv \in J_{q'}(\Omega)$.  In fact, for any
$\varphi \in \hat H^1_{q,0}(\Omega)$, we have
\begin{align}
0 &= (\bg, \nabla\varphi)_\Omega=\bar\lambda(\bv, \nabla\varphi)_\Omega
-(\DV(\mu\bD(\bv)), \nabla\varphi)_\Omega +(\nabla K_0(\bv),
\nabla\varphi)_\Omega \nonumber \\
&
= \bar\lambda(\bv, \nabla\varphi)_\Omega 
- (\nabla\dv\bv, \nabla\varphi)_\Omega. \label{uni:3}
\end{align}
Since $H^1_{q, 0}(\Omega) \subset \hat H^1_{q,0}(\Omega)$, for any 
$\varphi \in H^1_{q,0}(\Omega)$, we have 
\begin{equation}\label{weak:LD}
0 = \bar\lambda(\dv\bv, \varphi)_\Omega
 + (\nabla\dv\bv, \nabla\varphi)_\Omega.
\end{equation}
Choosing $\lambda_0 > 0$ larger if necessary, we have 
the uniquness of the resolvent problem 
\eqref{weak:LD} for the weak Dirichlet
operator, and so $\dv \bv=0$.  Putting this and \eqref{uni:3}
together gives $(\bv, \nabla\varphi)_\Omega=0$ for any $\varphi \in 
\hat H^1_{q,0}(\Omega)$, that is $\bv \in J_{q'}(\Omega)$. 
Analogously, we see that  $\bu \in J_q(\Omega)$, because
$\bu$ satisifies Eq. \eqref{uni:1}. 

Since $K_0(\bv) \in H^1_{q'}(\Omega) + \hat H^1_{q',0}(\Omega)$,
we write $K_0(\bv) = A_1 + A_2$ with $A_1 \in H^1_{q'}(\Omega)$
and $A_2 \in \hat H^1_{q',0}(\Omega)$.  
Since $\bu \in J_q(\Omega)$, we see that $\dv \bu=0$ in $\Omega$.
Thus,  by the definition of the solenoidal space $J_q(\Omega)$
and the divergence theorem of Gauss
\begin{align*}
(\bu, \nabla K_0(\bv))_\Omega &= (\bu, \nabla A_1)_\Omega
+ (\bu, \nabla A_2)_\Omega 
= (\bu\cdot\bn, A_1)_\Gamma -(\dv\bu, A_1)_\Omega \\
&= (\bu\cdot\bn, K_0(\bv))_\Gamma,  
\end{align*}
where we have used $K_0(\bv) = A_1$ on $\Gamma$.   Analogously, we have
$$(\nabla K_0(\bu), \bv)_\Omega 
= (K_0(\bu), \bv\cdot\bn)_\Gamma.
$$
Thus, by the divergence theorem of Gauss we have
\begin{align*}
(\bu, \bg)_\Omega& = (\bu, \bar \lambda \bv) 
-(\bu, \DV(\mu\bD(\bv) - K_0(\bv)\bI))_\Omega \\
& = \lambda(\bu, \bv) -(\bu, (\mu\bD(\bv)-K_0(\bv))\bn)_\Gamma
+ (\frac{\mu}{2}\bD(\bu), \bD(\bv))_\Omega\\
& = \lambda(\bu, \bv) + \frac12(\mu\bD(\bu), \bD(\bv))_\Omega.
\end{align*}
Analogously, we have
\begin{align*}
0& = (\lambda\bu - \DV(\mu\bD(\bu)-K_0(\bu)\bI), \bv)_\Omega 
-(\bu, \DV(\mu\bD(\bv) - K_0(\bv)\bI))_\Omega \\
& = \lambda(\bu, \bv) + \frac12(\mu\bD(\bu), \bD(\bv))_\Omega.
\end{align*}
Combining these two equalities yields that 
\begin{equation}\label{uni:6}
(\bu, \bg)_\Omega=0 \quad\text{for any $\bg \in J_{q'}(\Omega)$}.
\end{equation}
For any $\bff \in C^\infty_0(\Omega)^N$, let $\psi \in \hat H^1_{q',0}(\Omega)$
be a solution to the variational equation $(\bff, \nabla\varphi)_\Omega
=(\nabla\psi, \nabla\varphi)_\Omega$ for any $\varphi \in 
\hat H^1_{q,0}(\Omega)$.
Let $\bg = \bff - \nabla\psi$, and then $\bg \in J_{q'}(\Omega)$ and 
$(\bu, \nabla\psi)_\Omega = 0$, because $\bu \in J_q(\Omega)$.  
Thus, by \eqref{uni:6}, 
$(\bu, \bff)_\Omega = (\bu, \bg)_\Omega =0$, which, combined with
the arbitrariness of the choice of $\bff$, leads to  
$\bu = 0$.  This completes the proof of the uniqueness part of 
Theorem \ref{main:thm0}.


\subsection{A priori estimate} \label{sec:5.10}

In this subsection, we prove {\it a priori} estimates of solutions
of Eq. \eqref{rres:1.2}, from which the uniqueness part of 
Theorem \ref{main:thm3} follows immediately. 
In the previous subsection, we used the dual problem for Eq. \eqref{rres:1.1}
to prove 
the uniqueness, but in the present case, it is not clear what is
a suitable dual problem for Eq. \eqref{rres:1.2} 
to prove the uniqueness.  This is the 
reason why we consider the {\it a priori} estimates. 

To prove the {\it a priori} estimates, we need a slight
restriction of the domain $\Omega$.  Namely, 
in this subsection, we assume 
that $\Omega$ is a uniform $C^3$ domain whose inside has
a finite covering (cf. Definition \ref{dfn:2}), which is used to 
estimate the local norm of $K(\bu, h)$. 
The following theorem is the main result of this section.
\begin{thm}\label{thm:apri:1} 
Let $1 < q < \infty$. Let $\Omega$ be a uniformly $C^3$ domain 
whose inside has a finite covering. 
Then, there exists a $\lambda_0 > 0$ such that 
for any $\lambda \in \Lambda_{\sigma, \lambda_0\gamma_\kappa}$ 
and $(\bu, h)
\in H^2_q(\Omega)^N\times H^3_q(\Omega)$ satisfying Eq. \eqref{rres:1.2}, 
we have 
\begin{equation}\label{7.0}\begin{split}
&|\lambda|\|\bu\|_{L_q(\Omega)} + |\lambda|^{1/2}\|\bu\|_{H^1_q(\Omega)}
+ \|\bu\|_{H^2_q(\Omega)} + |\lambda|\|h\|_{H^2_q(\Omega)}
+ \|h\|_{H^3_q(\Omega)} \\
&\quad \leq C\{\|\bff\|_{L_q(\Omega)} + \|d\|_{W^{2-1/q}_q(\Gamma)} 
+ |\lambda|^{1/2}\|\bh\|_{L_q(\Omega)} + \|\bh\|_{H^1_q(\Omega)}\}.
\end{split}\end{equation}
\end{thm}
\begin{cor}\label{cor:unique} 
Let $1 < q < \infty$. Let $\Omega$ be a uniformly $C^3$ domain whose
inside has a finite covering. Then, there exists a $\lambda_0 > 0$ such that 
the uniqueness holds for Eq. \eqref{rres:1.2} 
for any $\lambda \in \Lambda_{\sigma, \lambda_0\gamma_\kappa}$.
\end{cor}
In what follows, we shall prove Theorem \ref{thm:apri:1}.  We first consider
Eq. \eqref{eq:7.0}. 
Let $\bu \in H^2_q(\Omega)^N$ and $h \in H^3_q(\Omega)$ satisfy 
Eq. \eqref{eq:7.0}.
We use the same notation as in Sect. \ref{subsec:5.6},  Proposition
\ref{prop:lap} and Definition \ref{dfn:2}.  Let $\psi^0$ be the function
given in Definition \ref{dfn:2}, that is 
$\psi_0 = \sum_{j=1}^\infty \zeta^0_j$, where $\zeta^0_j$ are the 
functions given in Proposition \ref{prop:lap}. Notice that 
${\rm supp}\, \psi^0 \subset \Omega$. 
Let 
$\bu^0 = \psi^0\bu$, and then $\bu^0$  satisfies the equations:
\begin{equation}\label{7.1}\begin{aligned}
\lambda \bu^0 - \DV(\mu\bD(\bu^0) - K_0(\bu^0)\bI)
&= \bff^0  &\quad&\text{in $\Omega$}, \\
(\mu\bD(\bu^0) - K_0(\bu_0)\bI)\bn &= 0 & \quad&\text{on $\Gamma$}.
\end{aligned}\end{equation}
Where, we have set
$$\bff^0 = \psi^0\bff + \psi^0\DV(\mu\bD(\bu))-\DV(\mu\bD(\psi^0\bu))
-(\psi^0\nabla K(\bu, h) - \nabla K_0(\psi^0\bu)),$$
and we have used the fact: 
$$K_0(\bu_0) = <\mu\bD(\bu^0)\bn, \bn> - \dv \bu^0 = 0 \quad
\text{on $\Gamma$}.
$$
We also let $\bu^1_j=\zeta^1_j\bu$ and $h_j=\zeta^1_j h$,
and then $\bu^1_j$ and $h_j$  satisfy the equations:
\begin{equation}\label{7.2}\left\{\begin{aligned}
\lambda \bu^1_j - \DV(\mu(x^1_j)\bD(\bu^1_j) - K_{1j}(\bu^1_j, h^1_j)\bI) 
&= \bff^1_j
&\quad&\text{in $\Omega_j$}, \\
\lambda h_j + A_\kappa(x^1_j)\cdot\nabla_{\Gamma_j}h_j - \bn_j\cdot\bu_j & = 
d_j &\quad&\text{on $\Gamma_j$}, \\
(\mu(x^1_j)\bD(\bu^1_j) - K_{1j}(\bu^1_j, h^1_j)\bI)\bn_j
-\sigma(x^1_j)(\Delta_{\Gamma_j}h_j)\bn_j & = \bh_j
&\quad&\text{on $\Gamma_j$}.
\end{aligned}\right.\end{equation}
Where, we have set 
\begin{align*}
&\bff^1_j  = \zeta^1_j\bff + \zeta^1_j\DV(\mu(x)\bD(\bu)) -
\DV(\mu(x)\bD(\zeta^1_j\bu)) \\
&\quad + \DV((\mu(x)-\mu(x^1_j))\bD(\zeta^1_j\bu))
-(\zeta^1_j\nabla K(\bu, h)-\nabla K_{1j}(\zeta^1_j\bu, \zeta^1_jh)); \\
&d_j  = \zeta^1_jd - \zeta^1_j(A_\kappa(x) - A_\kappa(x^1_j))
\cdot\nabla_{\Gamma_j}h_j \\
&\quad -A_\kappa(x^1_j)\cdot(\zeta^1_j\nabla_{\Gamma_j}h
-\nabla_{\Gamma_j}(\zeta^1_jh));\\
&\bh_j  = \zeta^1_j\bh - \{\zeta^1_j(\mu(x)-\mu(x^1_j))\bD(\bu)
+ \mu(x^1_j)(\zeta^1_j\bD(\bu) - \bD(\zeta^1_j\bu))\}\bn\\
&\quad + (\zeta^1_jK(\bu, h) - K_{1j}(\zeta^1_j\bu, \zeta^1_jh))\bn \\
&\quad + \zeta^1_j(\sigma(x) - \sigma(x^1_j))
(\Delta_{\Gamma}h)\bn
 + \sigma(x^1_j)(\zeta^1_j\Delta_\Gamma h
-\Delta_\Gamma(\zeta^1_jh))\bn.
\end{align*}
Set 
\begin{align*}
E_\lambda(\bu, h) &= |\lambda|^q\|\bu\|_{L_q(\Omega)}^q + 
|\lambda|^{q/2}\|\bu\|_{H^1_q(\Omega)}^q 
+ \|\bu\|_{H^2_q(\Omega)}^q \\
&\quad
+ |\lambda|^q\|h\|_{H^2_q(\Omega)}^q + \|h\|_{H^3_q(\Omega)}^q.
\end{align*}
Employing the similar argument to that in Subsec. \ref{subsec:n5.3}
and using Theorem \ref{main:thm0} to estimate $\bu^0$ in addition, 
for any positive number $\omega$ 
we have
\begin{align}
&E_\lambda(\bu, h)  \leq C\{
\|\bff^0_j\|_{L_q(\Omega)}^q +
\sum_{j=1}^\infty(
\|\bff^1_j\|_{L_q(\Omega_j)}^q + \|d_j\|_{W^{2-1/q}_q(\Gamma_j)}^q
\nonumber\\
&\phantom{E_\lambda(\bu, h)  \leq C\sum_{j=1}^\infty
aaaaaaaaaaaaaaa}
+|\lambda|^{q/2}\|\bh_j\|_{L_q(\Omega_j)}^q 
+ \|\bh_j\|_{H^1_q(\Omega_j)}^q)\} \nonumber \\
&\leq C\{\|\bff\|_{L_q(\Omega)}^q + \|d\|_{W^{2-1/q}_q(\Gamma)}^q 
+ |\lambda|^{q/2}\|\bh\|_{L_q(\Omega)} 
+ \|\bh\|_{H^1_q(\Omega)}^q +\gamma_\kappa^q\|\bh\|_{H^2_q(\Omega)}^q
\nonumber \\
&\quad +(\omega^q + M_1^q)
(\|\bu\|_{H^2_q(\Omega)}^q + \|\bh\|_{H^3_q(\Omega)}^q 
+ |\lambda|^{q/2}\|\bh\|_{H^1_q(\Omega))})\nonumber \\
&\quad 
+C_{\omega, M_2}(\|\bu\|_{H^1_q(\Omega)}
+ |\lambda|^{q/2}\|\bu\|_{L_q(\Omega)}
+\|h\|_{H^2_q(\Omega)}^q + |\lambda|^{q/2}\|h\|_{H^1_q(\Omega)}^q)
\nonumber \\
&\phantom{E_\lambda(\bu, h)  \leq C\sum_{j=1}^\infty
aaaaaaaaaaaaaaa}
+ \|K(\bu, h)\|_{L_q(\CO)}^q\}.\label{7.3}
\end{align}
Where,  $C_{\omega, M_2}$ is a constant depending on 
$\omega$ and $M_2$, and we have used the assumption that 
\begin{equation}\label{assump:dfn:2}
{\rm supp}\,\nabla \psi^0 
\cup \Bigl(\bigcup_{j=1}^\infty{\rm supp}\, \nabla\zeta^1_j
\Bigr)= \CO
\end{equation}
in Definition \ref{dfn:2}. 

To estimate $\|K(\bu, h)\|_{L_q(\CO)}$, we need the following 
Poincar\'es' type lemma.
\begin{lem}\label{lem:poin:1} Let $1 < q < \infty$ and let $\Omega$ be a 
uniformly $C^2$ domain whose inside has a  finite covering.  
Let $\CO$ be a set given in \eqref{assump:dfn:2}.  Then, we have
$$\|\varphi\|_{L_q(\CO)} \leq C_{q,\CO}\|\nabla\varphi\|_{L_q(\Omega)}
\quad\text{for any $\varphi \in \hat H^1_{q,0}(\Omega)$}$$
with some constant $C_{q, \CO}$ depending solely on $\CO$ and $q$.
\end{lem}
\pf ~Let $\CO_i$ ($i=1, \ldots, \iota$) be the sub-domains given in Definition 
\ref{dfn:2}, and then it is sufficient to prove that 
$$\|\varphi\|_{L_q(\CO_i)} \leq C\|\nabla\varphi\|_{L_q(\Omega)}
\quad\text{for any $\varphi \in \hat H_{q,0}(\Omega)$
and $i=1, \ldots, \iota$}.
$$
If $\CO_i \subset \Omega_R$ for some $R > 0$, since $\varphi|_{\Gamma}
= 0$, by the usual Poincar\'es' inequality we have 
$$\|\varphi\|_{L_q(\CO_i)} \leq \|\varphi\|_{L_q(\Omega_R)}
\leq C\|\nabla\varphi\|_{L_q(\Omega_R)}
\quad\text{for any $\varphi \in \hat H^1_{q,0}(\Omega)$}.
$$
Let $\CO_i$ be a subdomain for which the condition \thetag{b}
in Definition \ref{dfn:2} holds. 
Since the norms for $\varphi(\CA\circ\tau(y))$
and $\varphi(y)$ are equivalent, without loss of generality we
may assume that 
\begin{align*}
\CO_i \subset \{x = (x', x_N) \in \BR^N \mid a(x') < x_N < b, \quad
x' \in D\} &\subset \Omega \\
\{x = (x', x_N) \in \BR^N \mid x_N = a(x')\quad
x' \in D\} & \subset \Gamma.
\end{align*}
Since $\varphi \in \hat H^1_{q,0}(\Omega)$, we can write
$$\varphi(x', x_N) = \int^{x_N}_{a(x')}(\pd_s\varphi)(x', s)\,ds.
$$
because $\varphi(x', a(x')) = 0$. 
By Hardy inequality \eqref{ab:4}, we have
\begin{align*}
&\Bigl(\int^b_{a(x')}|\varphi(x', x_N)|^qx_N^{-q}\,dx_N\Bigr)^{1/q}
\\
&\quad
\leq \Bigl(\int^b_{a(x')}\Bigl(\int^{x_N}_{a(x')}
|(\pd_s\varphi)(x', s)|\,ds\Bigr)^q\,x_N^{-q}\,dx_N\Bigr)^{1/q}\\
&\quad
\leq \frac{q}{q-1}\Bigl(\int^b_{a(x')}
|s\pd_s\varphi(x', s)|^q\,s^{-q}
\,ds\Bigr)^{1/q},
\end{align*} 
and so, by Fubini's theorem we have 
\begin{align*}
\Bigl(\int_{\CO_i}|\varphi(x)|^q\,dx\Bigr)^{1/q}
& \leq 
\Bigl(\int_D\,dx'\int^b_{a(x')}|\varphi(x', x_N)|^q
\,dx_N\Bigr)^{1/q}\\
& \leq\Bigl(\int_D\,dx'\int^b_{a(x')}|\varphi(x', x_N)|^q\,
x_N^{-q}b^q\,dx_N\Bigr)^{1/q} \\
&\leq \frac{q b}{q-1}
\Bigl(\int_D\,dx'\int^b_{a(x')}|\pd_N\varphi(x)|^q\,dx_N
\Bigr)^{1/q}\\
& \leq bq'\|\nabla\varphi\|_{L_q(\Omega)}.
\end{align*}
This completes the proof of Lemma \ref{lem:poin:1}.  \qed
\vskip0.5pc
We now prove that for any $\omega > 0$ there exists a constant
$C_{\omega, M_2}$ depending on $\omega$ and $M_2$ such that 
\begin{equation}\label{7.4}\begin{aligned}
&\|K(\bu, h)\|_{L_q(\CO)} \\ 
&\quad
\leq \omega(\|\bu\|_{H^2_q(\Omega)} + \|h\|_{H^3_q(\Omega)}) 
+ C_{\omega, M_2}(\|\bu\|_{H^1_q(\Omega)}
+ \|h\|_{H^2_q(\Omega)}).
\end{aligned}\end{equation}
For this purpose, we estimate 
$|(K(\bu, h), \psi)_\Omega|$ for any $\psi \in C^\infty_0(\CO)$. 
 By Lemma \ref{lem:poin:1}, for any $\varphi \in \hat H^1_{q,0}(\Omega)$
we have 
$$|(\varphi, \psi)_\Omega| \leq \|\varphi\|_{L_q(\CO)}\|\psi\|_{L_{q'}(\CO)}
\leq C_{q,\CO}\|\nabla\varphi\|_{L_q(\CO)}\|\psi\|_{L_{q'}(\CO)}.
$$
Thus, by the Hahn-Banach theorem, there exists a $\bg \in L_{q'}(\Omega)^N$ such that $\|\bg\|_{L_{q'}(\Omega)} \leq C_{q', \CO}\|\psi\|_{L_{q'}(\CO)}$ and 
\begin{equation}\label{7.5} (\varphi, \psi)_\Omega
 = (\nabla\varphi, \bg)_\Omega
\end{equation}
for any $\varphi \in \hat H^1_{q,0}(\Omega)$. In particular, 
$\dv\bg = -\psi$, and therefore $\|\dv\bg\|_{L{q'}(\Omega)}
\leq \|\psi\|_{L_{q'}(\CO)}$.  By the assumption of the unique existence
of solutions of the weak Dirichlet problem and its regularity theorem,
Theorem \ref{thm:ap-r} in Subsec. \ref{sec:ap.2} below,
there exists a $\Psi \in \hat H^1_{q', 0}(\Omega)$ such 
that $\nabla^2\Psi \in L_q(\Omega)^N$, $\Psi$ satisfies the 
weak Dirichlet problem:
\begin{equation}\label{7.6}
(\nabla\Psi, \nabla\varphi)_\Omega = (\bg, \nabla\varphi)_\Omega
\quad\text{for any $\varphi \in \hat H^1_{q,0}(\Omega)$}
\end{equation}
and the estimate: 
\begin{equation}\label{7.7}
\|\nabla\Psi\|_{H^1_{q'}(\Omega)} \leq C_{q, \CO}\|\psi\|_{L_{q'}(\CO)}.
\end{equation}
Let $L = K(\bu, h) - \{<\mu\bD(\bu)\bn, \bn> - \sigma\Delta_\Gamma h
-\dv\bu\}$, where the Laplace-Beltrami operator $\Delta_\Gamma$ is suitably
extended in $\Omega$, 
and then $L \in \hat H^1_{q,0}(\Omega)$.  Thus, by \eqref{7.5},
\eqref{7.6} with $\varphi = L$ and the divergence theorem of 
Gau\ss, 
\begin{align*}
|(L, \psi)_\Omega| & = |(\nabla L , \bg)_\Omega| 
= |(\nabla\Psi, \nabla L)_\Omega| \\
&\leq |(\DV(\mu\bD(\bu)) - \nabla\dv\bu, \nabla\Psi)_\Omega|\\
& \qquad+ |(\nabla\{<\mu\bD(\bu)\bn, \bn> - \sigma\Delta_\Gamma h
-\dv\bu\}, \nabla\Psi)_\Omega| \\
& \leq C_{M_2}\{(\|\nabla\bu\|_{L_q(\Gamma)}
+ \|(h, \nabla h, \nabla^2h)\|_{L_q(\Gamma)})\|\nabla\Psi\|_{L_q(\Gamma)}
\\
&\qquad
+ (\|\nabla\bu\|_{L_q(\Omega)}+
\|h\|_{H^2_q(\Omega)})\|\nabla^2\Psi\|_{L_q(\Omega)}\}.
\end{align*}
Using the interpolation inequality: 
$\|v\|_{L_q(\Gamma)} \leq C\|\nabla v\|_{L_q(\Omega)}^{1/q}
\|v\|_{L_q(\Omega)}^{1-1/q}$ and \eqref{7.7}, we have
\begin{align*}
|(L, \psi)_\Omega| &\leq \{\omega(\|\nabla^2\bu\|_{L_q(\Omega)} 
+ \|h\|_{H^3_q(\Omega)}) \\
&\qquad + C_{\omega, M_2}
(\|\nabla\bu\|_{L_q(\Omega)} + \|h\|_{H^2_q(\Omega)})\}
\|\psi\|_{L_{q'}(\CO)},
\end{align*}
which leads to 
$$\|L\|_{L_q(\Omega)} \leq \omega
(\|\nabla^2\bu\|_{L_q(\Omega)} + \|h\|_{H^3_q(\Omega)})
+ C_{\omega, M_2}(\|\nabla\bu\|_{L_q(\Omega)}
+ \|h\|_{H^2_q(\Omega)}).$$
Thus, we have \eqref{7.4}.

Putting \eqref{7.3} and  \eqref{7.4} together and choosing
$\omega$ and $M_1$ small enough and $\lambda_0$ large enough, 
we have \eqref{7.0}.  This completes the proof of Theorem 
\ref{thm:apri:1} for Eq. \eqref{eq:7.0}. 

For Eq. \eqref{rres:1.2}, using the result we have proved just now,
we see that
$$
E_\lambda(\bu, h) \leq C(\|\bff\|_{L_q(\Omega)} + \|d\|_{W^{2-1/q}_q(\Gamma)}
+ |\lambda|^{1/2}\|\bh\|_{L_q(\Omega)} + \|\bh\|_{H^1_q(\Omega)} + R)
$$
with
$$R = \|\CF_1\bu\|_{W^{2-1/q}_q(\Gamma)} 
+ |\lambda|^{1/2}\|\CF_2h\|_{L_q(\Omega)} + \|\CF_2h\|_{H^1_q(\Omega)}).$$
Since 
$ R
\leq C(\|\bu\|_{H^1_q(\Omega)} + |\lambda|^{1/2}\|h\|_{H^2_q(\Omega)})$ 
as follows from \eqref{cond-low}, choosing $|\lambda|$ suitably large, 
$R$ can be absorbed by $E_\lambda(\bu, h)$, which completes the
proof of Theorem \ref{thm:apri:1}.
\section{Local well-posedness for arbitrary  data
in a uniform $C^3$ domain whose inside has
a finite covering} \label{sec:loc6}

In this section, we shall prove the unique existence of local 
in time solutions of Eq. \eqref{linear:6}.  
Namely, we shall prove the following theorem.
\begin{thm}\label{thm:loc.1}
Let $2 < p < \infty$, $N <  q < \infty$, $B > 0$
 and let $\Omega$ be a uniform $C^3$ domain whose
inside has a finite covering. Assume that $2/p + N/q < 1$ and 
that the weak Dirichlet problem is uniquely solvable for indices
$q$ and $q' = q/(q-1)$.  Assume that the mean curvature of $\Gamma$ satisfies
the condition \eqref{linear:2}.  Let $\bu_0\in B^{2(1-1/p)}_{q,p}(\Omega)^N$ and $\rho_0 \in W^{3-1/p-/q}_{q,p}(\Gamma)$ be initial data for which  
$\|\bu_0\|_{B^{2(1-1/p}_{q,p}(\Omega)} \leq B$ holds and 
the compatibility conditions: 
\begin{gather*}
\bu_0 - \bg(\bu_0, \Psi_\rho|_{t=0}) \in J_q(\Omega),
\quad \dv\bu_0 = g(\bu_0, \Psi_\rho|_{t=0}) \quad\text{in $\Omega$}, \\
(\mu\bD(\bu_0)\bn)_\tau = \bh'(\bu_0, \Psi_\rho|_{t=0})
\quad\text{on $\Gamma$}
\end{gather*}
are satisfied. 
Then, there exist a time $T > 0$ and a small number $\epsilon > 0$ 
depending on $R$ such that 
if $\|\rho_0\|_{B^{3-1/p-1/q}_{q,p}(\Gamma)} \leq \epsilon$, then problem
\eqref{linear:6} admits  unique solutions $\bu$, $\fq$ and $\rho$
with
\begin{equation}\label{sol.class:1}\begin{split}
\bu &\in L_p((0, T), H^2_q(\Omega)^N) \cap H^1_p((0, T), L_q(\Omega)^N), \\
\fq &\in L_p((0, T), H^1_q(\Omega) + \hat H^1_{q,0}(\Omega)), \nonumber \\
\rho &\in L_p((0, T), W^{3-1/q}_q(\Gamma)) 
\cap H^1_p((0, T), W^{2-1/q}_q(\Gamma))  
\end{split}\end{equation}
 possessing the estimate \eqref{trans:1} and 
$$E_{p,q,T}(\bu, \rho) \leq CB$$
for some constant $C$ independent of $B$.  Where, we have set
\begin{align*}
E_{p,q,T}(\bu, \rho)  &=  \|\bu\|_{L_p((0, T), H^2_q(\Omega))}
+ \|\pd_t\bu\|_{L_p((0, T), L_q(\Omega))} \\
&+ \|\pd_t\rho\|_{L_\infty((0, T), W^{1-1/q}_q(\Gamma))}
+ \|\pd_t\rho\|_{L_p((0, T), W^{2-1/q}_q(\Omega))} \\
&+ \|\rho\|_{L_p((0, T), W^{3-1/p}_q(\Omega))}.
\end{align*}
\end{thm}
In what follows, we shall prove Theorem \ref{thm:loc.1}
 by the Banach fixed point theorem.  Since $N < q < \infty$, by
Sobolev's inequality we have the following estimates:
\begin{gather}
\|f\|_{L_\infty(\Omega)} \leq C\|f\|_{H^1_q(\Omega)}, 
\quad \|fg\|_{H^1_q(\Omega)} \leq C\|f\|_{H^1_q(\Omega)}
\|g\|_{H^1_q(\Omega)}, \nonumber\\
\|fg\|_{W^{1-1/q}_q(\Gamma)} \leq C\|f\|_{W^{1-1/q}_q(\Gamma)}
\|g\|_{W^{1-1/q}_q(\Gamma)}. \label{sob:5.0}
\end{gather}
Moreover, since $2/p+N/q < 1$, we have
\begin{equation}\label{sol:5.3}
\|f\|_{H^1_\infty(\Omega)} \leq C\|f\|_{B^{2(1-1/p)}_{q,p}(\Omega)},
\quad \|f\|_{H^2_\infty(\Gamma)} 
\leq C\|f\|_{B^{3-1/p-1/q}_{q,p}(\Gamma)}.
\end{equation}
We define an  underlying space $\bU_T$ by letting
\begin{align*}
\bU_T  = \{(\bu, \rho) \mid \enskip 
&\bu \in H^1_p((0, T), L_q(\Omega)^N) \cap
L_p((0, T), H^2_q(\Omega)^N),  \\
& \rho \in H^1_p((0, T), W^{2-1/q}_q(\Gamma))
\cap L_p((0, T), W^{3-1/q}_q(\Gamma)),\\
&\bu|_{t=0} = \bu_0 \quad\text{in $\Omega$}, \quad 
\rho|_{t=0} = \rho_0
\quad\text{on $\Gamma$}, \\ 
&E_{p,q,T}(\bu, \rho) \leq L, 
\quad \sup_{t \in (0, T)} \|\rho(\cdot, t)\|_{H^1_\infty(\Gamma)}
\leq \delta\},
\end{align*}
where $L$ is a number determined later. Since $L$ is chosen 
large and $\epsilon$ small eventually, we may assume that 
$0 < \epsilon < 1 < L$ in the following. Thus, for example,
we will use the inequality: $1 + \epsilon + L < 3L$ below.

Let $(\bv, h) \in \bU_T$ and let $\bu$, $\fq$ and $\rho$ 
be solutions of the linear equations:
\begin{equation}\label{g.5.5}\left\{\begin{aligned}
\pd_t\bu - \DV(\mu\bD(\bu) - \fq\bI)  = \bff(\bv, \Psi_h)
+ \kappa_b\CH(\Gamma)&
&\quad&\text{in $\Omega^T$}, \\
\dv\bu = g(\bu, \Psi_h) = \dv\bg(\bv, \Psi_h)&
&\quad&\text{in $\Omega^T$}, \\
\pd_t\rho + <\bu_\kappa \mid \nabla'_\Gamma\rho> - \bu\cdot\bn 
= d_\kappa(\bv, \Psi_h)
&&\quad&\text{on $\Gamma^T$}, \\
(\mu\bD(\bu)\bn)_\tau = \bh'(\bv, \Psi_h)
&&\quad&\text{on $\Gamma^T$}, 
\\
<\mu\bD(\bu)\bn, \bn> -\fq
-\sigma (\Delta_\Gamma\rho + \BB\rho) 
= h_N(\bv, \Psi_h)+ \kappa_a\CH(\Gamma)
&&\quad&\text{on $\Gamma^T$}, \\
(\bu, \rho)|_{t=0} = (\bu_0, \rho_0)
&&\quad&\text{in $\Omega\times\Gamma$}.
\end{aligned}\right.\end{equation}
Where, the operator $\BB$ is defined in \eqref{bdyop:1} and we have set 
\begin{align*}
d_\kappa(\bv, \Psi_h)
&= d(\bv, \Psi_h) + <\bu_\kappa - \bv \mid \nabla'_\Gamma h>.
\end{align*} 
 Let $\Psi_h = \omega H_h\bn$, where 
$H_h$ is an extension of $h$ to $\Omega$ such that
\begin{equation}\label{f.6.0}\begin{aligned}
C_1\|H_h(\cdot, t)\|_{H^k_q(\Omega)} 
&\leq \|h(\cdot, t)\|_{W^{k-1/q}_q(\Gamma)}\leq
C_2\|H_h(\cdot, t)\|_{H^k_q(\Omega)} \\
C_1\|\pd_tH_h(\cdot, t)\|_{H^\ell_q(\Omega)} 
&\leq \|\pd_th(\cdot, t)\|_{W^{\ell-1/q}_q(\Gamma)}
\leq 
C_2\|\pd_tH_h(\cdot, t)\|_{H^\ell_q(\Omega)}
\end{aligned}\end{equation}
for $k=1,2,3$ and $\ell=1,2$ 
with some constants $C_1$ and $C_2$.  In particular, 
for $(\bv, h) \in \bU_T$, we may assume that 
\begin{equation}\label{f.6.1}
 \sup_{t \in (0, T)} \|\Psi_h(\cdot, t)\|_{H^1_\infty(\BR^N)}
\leq \delta.
\end{equation}
In fact, we can assume that 
$\sup_{t \in (0, T)} \|h(\cdot, t)\|_{H^1_\infty(\Gamma)}
\leq \delta$ with smaller $\delta$.

Since $(\bv, h) \in \bU_T$, we have  
\begin{equation}\label{nonloc:0}
E_{p,q,T}(\bv, h) \leq L,
\end{equation}
that is 
\begin{align}
&\|\pd_th\|_{L_\infty((0, T), W^{1-1/q}_q(\Gamma))}
+ \|\pd_th\|_{L_p((0, T), W^{2-1/q}_q(\Gamma))}
+ \|h\|_{L_p((0, T), W^{3-1/q}_q(\Gamma))}\nonumber \\
&\quad + \|\pd_t\bv\|_{L_p((, T), L_q(\Omega))}
+ \|\bv\|_{L_p((0, T), H^2_q(\Omega))}
\leq L. \label{nonloc:0*}
\end{align}
Moreover, we use the assumption:
\begin{equation}\label{nonloc:0**}
\|\bu_0\|_{B^{2(1-1/p)}_{q,p}(\Omega)}\leq B,
\quad \|\rho_0\|_{B^{3-1/p-1/q}_{q,p}(\Gamma)}
\leq \epsilon,
\end{equation}
where $\epsilon$ is a small constant determined later.

To obtain the estimates of $\bu$ and $\rho$, we shall use
Corollary \ref{cor:max.1}.  Thus, we shall estimate
the nonlinear functions appearing in the right hand side of Eq.\eqref{g.5.5}. 
We start with proving that  
\begin{equation}\label{nonloc:2}
\|\bff(\bv, \Psi_h)\|_{L_p((0, T), L_q(\Omega))} 
\leq C\{T^{1/p}(L+B)^2 + (\epsilon
+ T^{1/{p'}}L)L\}.
\end{equation}
The definition of $\bff(\bv, \Psi_h)$ is given by 
 replacing $\xi(t)$, $\rho$ and $\bu$ by 
$0$, $h$ and $\bv$ \eqref{form:f}. 
Since $|\bV_0(\bk)| \leq C|\bk|$ when $|\bk| \leq \delta$, 
by \eqref{f.6.1} we have
\begin{align*}
\|\bff(\bv, \Psi_h)\|_{L_q(\Omega)} 
&\leq C\{\|\bv\|_{L_\infty(\Omega)}\|\nabla\bv\|_{L_q(\Omega)}
+ \|\pd_t\Psi_h\|_{L_\infty(\Omega)}\|\nabla\bv\|_{L_q(\Omega)}
\\
&
+ \|\nabla\Psi_h\|_{L_\infty(\Omega)}\|\pd_t\bv\|_{L_q(\Omega)}
+ \|\nabla\Psi_h\|_{L_\infty(\Omega)}
\|\nabla^2\bv\|_{L_q(\Omega)}\\
&+ \|\nabla^2\Psi_h\|_{L_q(\Omega)}\|\nabla\bv\|_{L_\infty(\Omega)}\}.
\end{align*}
By \eqref{f.6.0} and \eqref{sob:5.0}, 
we have
\begin{equation}\label{f.6.4}\begin{aligned}
\|\bv(\cdot, t)\|_{L_\infty(\Omega)} & \leq 
C\|\bv(\cdot, t)\|_{H^1_q(\Omega)}, \\
\|\pd_t\Psi_h(\cdot, t)\|_{L_\infty(\Omega)} & \leq C\|\pd_th
(\cdot, t)\|_{W^{1-1/q}_q(\Gamma)}, \\
\|\nabla\Psi_h(\cdot, t)\|_{L_\infty(\Omega)} & \leq 
C\|h(\cdot, t)\|_{W^{2-1/q}_q(\Gamma)}, \\
\|\nabla\bv(\cdot, t)\|_{L_\infty(\Omega)}
& \leq C\|\bv(\cdot, t)\|_{H^2_q(\Omega)},
\end{aligned}\end{equation}
and so, 
\begin{align}
&\|\bff(\bv, \Psi_h)\|_{L_p((0, T), L_q(\Omega))}
\leq C\{\|\bv\|^2_{L_\infty((0, T), H^1_q(\Omega))}T^{1/p}
\nonumber \\
&\quad + \|\bv\|_{L_\infty((0, T), H^1_q(\Omega))}
\|\pd_th\|_{L_\infty((0, T), W^{1-1/q}_q(\Gamma))}T^{1/p}
\label{f.6.3}\\
&\quad+ \|h\|_{L_\infty((0, T), W^{2-1/q}_q(\Gamma))}
(\|\pd_t\bv\|_{L_p((0, T), L_q(\Omega))}
+ \|\bv\|_{L_p((0, T), H^2_q(\Omega))})\}.
\nonumber
\end{align}
Here and in the following, we use the estimate:
$$\|f\|_{L_p((0, T), X)} \leq T^{1/p}\|f\|_{L_\infty((0, T), X)}
$$
for the lower order term. 
In what follows, we shall use the following inequalities:
\begin{align}
&\|\bv\|_{L_\infty((0, T), B^{2(1-1/p)}_{q,p}(\Omega))}
\leq C\{\|\bu_0\|_{B^{2(1-1/p)}_{q,p}(\Omega)} \nonumber \\
&\quad + \|\bv\|_{L_p((0, T), H^2_q(\Omega))}
+ \|\pd_t\bv\|_{L_p((0, T), L_q(\Omega))}\}, \label{real:7.1.1}\\
&\|h\|_{L_\infty((0, T), B^{3-1/p-1/q}_{q,p}(\Gamma))} 
\leq C\{\|\rho_0\|_{B^{3-1/p-1/q}_{q,p}(\Gamma)}
\nonumber \\
&\quad
+ \|h\|_{L_p((0, T), W^{3-1/q}_q(\Gamma))}
 + \|\pd_th\|_{L_p((0, T), W^{2-1/q}_q(\Gamma))}\},  \label{real:7.1.2}
\\
&\|h(\cdot,t)\|_{L_\infty((0, T), W^{2-1/q}_q(\Gamma))} \nonumber\\
&\quad \leq \|\rho_0\|_{W^{2-1/q}_q(\Gamma)}
+ \int^T_0\|\pd_sh(\cdot, s)\|_{W^{2-1/q}_q(\Gamma)}\,ds \nonumber\\
&\quad \leq \|\rho_0\|_{W^{2-1/q}_q(\Gamma)}
+ T^{1/{p'}}L \label{nonloc:1}, 
\end{align}
for some constant $C$ independent of $T$. The inequalities: 
\eqref{real:7.1.1} and \eqref{real:7.1.2} will be proved later. 
Combining \eqref{f.6.3} with \eqref{real:7.1.1}, 
\eqref{nonloc:1}, \eqref{nonloc:0*}
and \eqref{nonloc:0**} gives \eqref{nonloc:2}.

We next consider $d_\kappa(\bu, \Psi_h)$. 
Since $\xi(t) = 0$, by \eqref{kinematic:4*} we have
$$d(\bv, \Phi_h) = \bv\cdot\bV_\Gamma(\bar\nabla\Psi_h)
\bar\nabla\Psi_h\otimes\bar\nabla\Psi_h 
-\pd_th<\bn, 
\bV_\Gamma(\bar\nabla\Psi_h)
\bar\nabla\Psi_h\otimes\bar\nabla\Psi_h>.
$$
Applying \eqref{repr:2.1***}, \eqref{f.6.1} and \eqref{sob:5.0}, 
we have 
\begin{align*}
&\|\bV_\Gamma(\bar\nabla\Psi_h)
\bar\nabla\Psi_h\otimes\bar\nabla\Psi_h 
\|_{W^{1-1/q}_q(\Gamma)}
\leq C\|h(\cdot,t)\|_{W^{2-1/q}_q(\Gamma)}, \\
&\|\bV_\Gamma(\bar\nabla\Psi_h)
\bar\nabla\Psi_h\otimes\bar\nabla\Psi_h
\|_{W^{2-1/q}_q(\Gamma)}\\
&\qquad \leq C\|h(\cdot, t)\|_{W^{2-1/q}_q(\Gamma)}
\|h(\cdot, t)\|_{W^{3-1/q}_q(\Gamma)},
\end{align*}
and so, by \eqref{f.6.4} and \eqref{sob:5.0} 
\begin{align*}
&\|d(\bv, \Psi_h)\|_{W^{1-1/q}_q(\Gamma)} \\
&
\leq C(\|\pd_th(\cdot, t)\|_{W^{1-1/q}_q(\Gamma)}
+ \|\bv(\cdot, t)\|_{H^1_q(\Omega)})
\|h(\cdot, t)\|_{W^{2-1/q}_q(\Gamma)}, \\
&\|d(\bv, \Psi_h)\|_{W^{2-1/q}_q(\Gamma)} \\
&
\leq C\{(\|\pd_th(\cdot, t)\|_{W^{2-1/q}_q(\Gamma)}
+ \|\bv(\cdot, t)\|_{H^2_q(\Omega)})
\|h(\cdot, t)\|_{W^{2-1/q}_q(\Gamma)}
\\
&+ (\|\pd_th(\cdot, t)\|_{W^{1-1/q}_q(\Gamma)}
+ \|\bv(\cdot, t)\|_{H^1_q(\Omega)})
\|h(\cdot, t)\|_{W^{2-1/q}_q(\Gamma)}
\|h(\cdot, t)\|_{W^{3-1/q}_q(\Gamma)}\}.
\end{align*}
Thus, by \eqref{nonloc:0*}, \eqref{nonloc:0**}, 
\eqref{real:7.1.1}, and \eqref{nonloc:1}, we have
\begin{equation}\label{nonloc:3}\begin{split}
\sup_{t\in(0, T)} \|d(\bv, \Psi_h)\|_{W^{1-1/q}_q(\Gamma)}
& \leq C(L+B)(\epsilon 
+T^{1/{p'}}L), \\
\|d(\bv, \Psi_h)\|_{L_p((0, T), W^{2-1/q}_q(\Gamma))}
& \leq C(L+B)L(\epsilon 
+T^{1/{p'}}L).
\end{split}\end{equation}

By \eqref{linear:5} and \eqref{sob:5.0}, we have
\begin{equation}\label{f.6.7}\begin{aligned}
&\|<\bu_\kappa-\bv\mid\nabla'_\Gamma h>\|_{W^{1-1/q}_q(\Gamma)}
\\
&\quad\leq C(\|\bu_0\|_{B^{2(1-1/p)}_{q,p}(\Omega)}
+ \|\bv(\cdot, t)\|_{H^1_q(\Omega)})\|h(\cdot, t)\|_{W^{2-1/q}_q(\Gamma)},
\\
&\|<\bu_\kappa-\bv\mid\nabla'_\Gamma h>\|_{W^{2-1/q}_q(\Gamma)}\\
&\quad \leq C(\|\bu_\kappa\|_{H^2_q(\Omega)}
+\|\bv(\cdot, t)\|_{H^2_q(\Omega)})
\|h(\cdot, t)\|_{W^{2-1/q}_q(\Gamma)} \\
&\qquad + \|\bu_\kappa-\bv(\cdot, t)\|_{H^1_q(\Omega)}
\|h(\cdot, t)\|_{W^{3-1/q}_q(\Gamma)}).
\end{aligned}\end{equation}
By \eqref{real:7.1.1}, \eqref{nonloc:1}, \eqref{nonloc:0*}, and 
\eqref{nonloc:0**}, we have
\begin{equation}\label{nonloc:8*}\begin{split}
&\|<\bu_\kappa-\bv\mid\nabla'_\Gamma h>
\|_{L_\infty((0, T), W^{1-1/q}_q(\Gamma))} \\
&\quad\leq C(B+L)(\epsilon + T^{1/{p'}}L).
\end{split}\end{equation}
By the definition of $\bu_\kappa$, we have
$$\bu_\kappa - \bu_0 = \frac{1}{\kappa}\int^\kappa_0(T(s)\tilde \bu_0 - \bu_0)\,ds,
$$
and so, we have
\begin{align*}
\|\bu_\kappa-\bu_0\|_{L_q(\Omega)} 
&\leq \frac{1}{\kappa}\int^\kappa_0\Bigl(\int^s_0\|\pd_s T(r)
\tilde\bu_0\|_{L_q(\Omega)}\,dr\Bigr)\,ds \\
&\leq C\kappa^{1/{p'}}\|\bu_0\|_{B^{2(1-1/p)}_{q,p}(\Omega)}.
\end{align*}
Let $s$ be a positive number such that $1 < s < 2(1-1/p)$, and then
by interpolation theory and \eqref{linear:5} we have
\begin{equation}\label{nonloc:5}\begin{split}
\|\bu_\kappa-\bu_0\|_{H^1_q(\Omega)}
&\leq C\|\bu_\kappa-\bu_0\|_{L_q(\Omega)}^{1-1/s}
\|\bu_\kappa - \bu_0\|_{W^s_q(\Omega)}^{1/s} \\
&\leq C\kappa^{(1-1/s)/p'}\|\bu_0\|_{B^{2(1-1/p)}_{q,p}(\Omega)}.
\end{split}\end{equation}
Writing $\bv(\cdot, t) - \bu_0 = \int^t_0\pd_r\bv(\cdot, r)\,dr$, we have
$$\|\bv(\cdot, t) - \bu_0\|_{L_q(\Omega)} \leq LT^{1/{p'}}.
$$
On the other hand, by \eqref{real:7.1.1} 
$$\|\bv(\cdot,t) - \bu_0\|_{B^{2(1-1/p)}_{q,p}(\Omega)} \leq L+\|\bu_0
\|_{B^{2(1-1/p)}_{q,p}(\Omega)},
$$
and so, 
\begin{equation}\label{nonloc:6}
\|\bv(\cdot, t)-\bu_0\|_{H^1_q(\Omega)}
\leq (L+B)T^{(1-1/s)/p'}.
\end{equation} 
Putting \eqref{nonloc:5} and \eqref{nonloc:6} together
gives
\begin{equation}\label{nonloc:7}
\sup_{t\in (0, T)} \|\bv(\cdot, t) - \bu_\kappa\|_{H^1_q(\Omega)}
\leq C(\kappa^{(1-1/s)/{p'}} + T^{(1-1/s)/p'})(L+B).
\end{equation}
By \eqref{f.6.7} 
\begin{multline*}
\|<\bu_\kappa-\bv\mid\nabla'_\Gamma h>\|_{L_p((0, T), W^{2-1/q}_q(\Gamma))}
\\
\leq C\{(\|\bu_\kappa\|_{H^2_q(\Omega)}T^{1/p} 
+ \|\bv\|_{L_p((0, T), H^2_q(\Omega))})
\|h\|_{L_\infty((0, T), W^{2-1/q}_q(\Gamma))}\\
+ \|\bu_\kappa-\bv\|_{L_\infty((0, T), H^1_q(\Omega))}
\|h\|_{L_p((0, T), W^{3-1/q}_q(\Gamma))}\}.
\end{multline*}
Thus, by \eqref{linear:5}, \eqref{nonloc:0*}, \eqref{nonloc:0**},
\eqref{nonloc:1} and \eqref{nonloc:7}, we have
\begin{equation}\label{nonloc:8}\begin{split}
&\|<\bu_\kappa-\bv \mid \nabla'_\Gamma h>\|_{L_p((0, T), 
W^{2-1/q}_q(\Gamma))} \\
&\quad\leq C\{(B\kappa^{-1/p}T^{1/p}+L)(\epsilon
+ T^{1/{p'}}L) \\
&\qquad 
+ L(L+B)(\kappa^{(1-1/s)/p'} + T^{(1-1/s)/p'})\}
\end{split}\end{equation}
for some constant $s \in (1, 2(1-1/p))$. 
Putting \eqref{nonloc:3}, \eqref{nonloc:8*}, 
and \eqref{nonloc:8} gives 
\begin{align}
&\|\tilde d_\kappa(\bv, \Psi_h)\|_{L_\infty((0, T), W^{1-1/q}_q(\Gamma))}
\leq C(L+B)(\epsilon + T^{1/{p'}}L)
\nonumber \\
&\|\tilde d_\kappa(\bv, \Psi_h)\|_{L_p((0, T), W^{2-1/q}_q(\Gamma))}
\leq C\{(L+B\kappa^{-1/p}T^{1/p})(\epsilon 
+ T^{1/{p'}}L) \nonumber \\
&\quad
+ L(L+B)(\epsilon + T^{1/{p'}}L
+ \kappa^{(1-1/s)/{p'}}
+ T^{(1-1/s)/p'})\}, \label{nonloc:8**}
\end{align}
where $s$ is a constant $ \in (1, 2(1-1/p))$.

We now estimate $g(\bv, \Psi_h)$ and $\bg(\bv, \Psi_h)$,  
which are defined in \eqref{form:g}.   In view of Corollary \ref{cor:max.1}, 
we have to extend $g(\bv, \Psi_h)$ and $\bg(\bv, \Psi_h)$ to the whole
time interval $\BR$.  Before turning to the extension of these functions,
we make a few definitions. 
Let $\tilde\bu_0 \in B^{2(1-1/p)}_{q,p}(\BR^N)^N$ be an extension of 
$\bu_0 \in B^{2(1-1/p)}_{q,p}(\Omega)^N$ to $\BR^N$ such that 
$$\bu_0 = \tilde\bu_0 \quad\text{in $\Omega$}, \quad  
\|\tilde\bu_0\|_{X(\BR^N)} 
\leq C\|\bu_0\|_{X(\Omega)}$$
for $X \in \{H^k_q, B^{2(1-1/p)}_{q,p}\}$. 
 Let  
\begin{equation}\label{tv:1}
T_v(t)\bu_0
= e^{-(2-\Delta)t}\tilde\bu_0
= \CF^{-1}[e^{-(|\xi|^2+2)t}\CF[\tilde\bu_0](\xi)], 
\end{equation}
 and then 
$T_v(0)\bu_0=\bu_0$ in $\Omega$ and 
\begin{align}
\|e^tT_v(t)\bu_0\|_{X(\Omega)} 
&\leq C\|\bu_0\|_{X(\Omega)}, \nonumber\\
\|e^tT_v(\cdot)\bu_0\|_{H^1_p((0, \infty), L_q(\Omega))}
&+ \|e^tT_v(\cdot)\bu_0\|_{L_p((0, \infty), H^2_q(\Omega))} 
\nonumber \\
&
\leq C\|\bu_0\|_{B^{2(1-1/p)}_{q,p}(\Omega)}.
\label{g.5.7} 
\end{align}
Let $\bW$, $P$, and $\Xi$ be solutions of the equations:
$$\left\{\begin{aligned}
\pd_t\bW + \lambda_0\bW- \DV(\mu\bD(\bW) - P\bI)  = 0,
\quad 
\dv\bW = 0&
&\quad&\text{in $\Omega\times(0, \infty)$}, \\
\pd_t\Xi + \lambda_0\Xi- \bW\cdot\bn 
= 0
&&\quad&\text{on $\Gamma\times(0, \infty)$}, \\
(\mu\bD(\bW)\bn - P\bI) \bn
-(\sigma\Delta_\Gamma \Xi) \bn
=0&&\quad&\text{on $\Gamma\times(0, \infty)$}, 
\\
(\bW, \Xi)|_{t=0} = (0, \rho_0)
&&\quad&\text{in $\Omega\times\Gamma$}.
\end{aligned}\right.
$$
Since 
\begin{align*}
&\|\bW\|_{L_\infty((0, \infty), B^{2(1-1/p)}_{q,p}(\Omega))}\\
&\quad 
\leq C(\|\bW\|_{L_p((0, \infty), H^2_q(\Omega))}
+ \|\pd_t\bW\|_{L_p((0, \infty), L_q(\Omega))}), \\
&\|\Xi\|_{L_\infty((0, \infty), B^{3-1/p-1/q}_{q,p}(\Gamma))}\\
&\quad \leq 
C(\|\Xi\|_{L_p((0, \infty), W^{3-1/q}_q(\Gamma))}
+ \|\pd_t\Xi\|_{L_p((0, \infty), W^{2-1/q}_q(\Gamma))}),
\end{align*}
as follow from real interpolation \eqref{tanabe:1} and \eqref{tanabe:2}.  
choosing  $\lambda_0$ large enough,
 by Theorem \ref{thm:max.2},
we know the existence of $\bW$ and $\Xi$ with
\begin{align*}
\bW &\in L_p((0, \infty), H^2_q(\Omega)^N) \cap 
H^1_p((0, \infty), L_q(\Omega)^N), \\
\Xi &\in L_p((0, \infty), W^{3-1/q}_q(\Gamma)) \cap 
H^1_p((0, \infty), W^{2-1/q}_q(\Gamma))
\end{align*}
possessing the estimates:
\begin{align*}
&\|e^t\bW\|_{L_\infty((0, \infty), B^{2(1-1/p)}_{q,p}(\Omega))}
+\|e^t\Xi\|_{L_\infty((0, \infty), B^{3-1/p-1/q}_{q,p}(\Gamma))}
\\
&\quad + \|e^t\bW\|_{L_p((0, \infty), H^2_q(\Omega))}
+ \|e^t\pd_t\bW\|_{L_p((0, \infty), L_q(\Omega))} \\
&\quad\quad +\|e^t\Xi\|_{L_p((0, \infty), W^{3-1/q}_q(\Gamma))}
 + \|e^t\pd_t\Xi\|_{L_p((0, \infty), W^{2-1/q}_q(\Gamma))}\\
&\qquad\quad \leq C\|\rho_0\|_{B^{3-1/p-1/q}_{q,p}(\Gamma)}.
\end{align*}
Moreover, by the trace theorem and the kinematic equation,
we have
\begin{align*}
&\|e^t\pd_t\Xi\|_{L_\infty((0, \infty), W^{1-1/q}_q(\Gamma))} \\
&\quad
\leq \lambda_0\|e^t\Xi\|_{L_\infty((0, \infty), W^{1-1/q}_q(\Gamma))}
+ \|e^t\bn\cdot\bW\|_{L_\infty((0, \infty), W^{1-1/q}_q(\Gamma))}\\
&\quad
\leq C\|\rho_0\|_{W^{3-1/p-1/q}_{q,p}(\Gamma)}.
\end{align*}
Setting $T_h(t)\rho_0 = \Psi_\Xi$, we have $T_h(0)\rho_0 = \Psi_{\rho_0}$
in $\Omega$, and 
\begin{equation}\label{g.5.8}\begin{split}
&\|e^tT_h(\cdot)\rho_0\|_{L_\infty(0, \infty), B^{3-1/p}_{q,p}(\Omega))}
+ \|e^t\pd_tT_h(\cdot)\rho_0\|_{L_\infty((0, \infty), H^1_q(\Omega))} \\
&\quad
+ \|e^tT_h(\cdot)\rho_0\|_{L_p((0, \infty), H^3_q(\Omega))}
+ \|e^t\pd_tT_h(\cdot)\rho_0\|_{L_p((0, \infty), H^2_q(\Omega))}\\
&\qquad
\leq C\|\rho_0\|_{B^{3-1/p-1/q}_{q,p}(\Gamma)}.
\end{split}\end{equation}
Let $\psi(t) \in C^\infty(\BR)$ equal one for $t > -1$ and 
zero for $t < -2$.  Given a function, $f(t)$, defined on $(0, T)$, 
the extension, $e_T[f]$, of $f$ is defined by letting
\begin{equation}\label{eT:1}
e_T[f](t) = \begin{cases} 0&\quad\text{for $t < 0$}, \\
f(t) & \quad \text{for $0 < t < T$}, \\
f(2T-t)&\quad\text{for $T < t < 2T$}, \\
0 & \quad\text{for $t > 2T$}.
\end{cases}
\end{equation}
Obviously, $e_T[f](t)=f(t)$ for $t \in (0, T)$ and $e_T[f](t) = 0$
for $t \not\in(0, 2T)$.  Moreover, if $f|_{t=0}= 0$, then  
\begin{equation}\label{eT:2}
\pd_te_T[f](t) = \begin{cases} 0&\quad\text{for $t < 0$}, \\
(\pd_tf)(t) & \quad \text{for $0 < t < T$}, \\
-(\pd_tf)(2T-t)&\quad\text{for $T < t < 2T$}, \\
0 & \quad\text{for $t > 2T$}.
\end{cases}
\end{equation}
We now define the extensions, $\CE_1[\bv]$, and  
$\CE_2[\Psi_h]$,  of  $\bv$ and $\Psi_h$ to 
$\BR$ by letting
\begin{equation}\label{eT:2*}\begin{split}
\CE_1[\bv] &= e_T[\bv-T_v(t)\bu_0] + \psi(t)T_v(|t|)\bu_0,
\\
\CE_2[\Psi_h] &= e_T[\Psi_h - T_h(t)\rho_0] + 
\psi(t)T_h(|t|)\rho_0.
\end{split}\end{equation}
Since $\bv|_{t=0} = T_v(0)\bu_0=\bu_0$ and $\Psi_h|_{t=0} = 
T_h(0)\rho_0$, we can differentiate $\CE_1[\bv]$ 
and $\CE_2[\Psi_h]$ once with respect to $t$ and we can use
the formula \eqref{eT:2}. Obviously, we have
\begin{equation}\label{eT:3}
\CE_1[\bv] = \bv, \quad \CE_2[\Psi_h] = \Psi_h 
\quad \text{in $\Omega^T$}.
\end{equation}
Let
\begin{equation}\label{g.5.10}\begin{split}
\tilde g(\bv, \Psi_h) & = -\{J_0(\nabla\CE_2[\Psi_h])\dv\CE_1[\bv]
\\
&\quad + (1+J_0(\nabla\CE_2[\Psi_h]))\bV_0(\nabla\CE_2[\Psi_h])
:\nabla\CE_1[\bv]\},\\
\tilde\bg(\bv, \Psi_h) & = -(1+J_0(\nabla\CE_2[\Psi_h]))
{}^\top\bV_0(\nabla\CE_2[\Psi_h])\CE_1[\bv],
\end{split}\end{equation}
and then, applying the Hanzawa transform: $x = y+ \CE_2[\Psi_h]$
instead of $x = y + \Psi_h$, 
by \eqref{form:g}, \eqref{div:3}, and \eqref{eT:3}, we have 
\begin{equation} \label{g.5.11} \begin{split}
\tilde g(\bv, \Psi_h) = g(\bv, \Psi_h),\quad
\tilde\bg(\bv, \Psi_h) = \bg(\bv, \Psi_h) \quad&\text{in 
$\Omega^T$}, \\
\dv \tilde\bg(\bv, \Psi_h) = \tilde g(\bv, \Psi_h)
\quad&\text{in $\Omega\times\BR$}.
\end{split}\end{equation}
By \eqref{g.5.8} and \eqref{nonloc:1}, we have
\begin{align*}
\sup_{t\in\BR}\|\CE_2[\Psi_h]\|_{H^1_\infty(\Omega)} 
&\leq C(\sup_{t\in(0, T)}\|\Psi_h\|_{H^2_q(\Omega)}
+ \|T(\cdot)\rho_0\|_{L_\infty, H^2_q(\Omega)}) \\
&\leq C(\|\rho_0\|_{W^{3-1/p1/q}_q(\Gamma)} + T^{1/{p'}}L)
\end{align*}
Thus, we choose $T$ and $\|\rho\|_{W^{3-1/p-1/q}_{q,p}(\Omega)}$ so small that
\begin{equation}\label{g.5.8*}
\sup_{t\in\BR}\|\CE_2[\Psi_h]\|_{H^1_\infty(\Omega)} \leq \delta.
\end{equation}
Since $\bV_0(0) = 0$, we may write $\tilde \bg(\bv, \Psi_h)$ as
\begin{equation}\label{g.5.10*}
\tilde \bg(\bv, \Psi_h) 
= \bV_\bg(\nabla\CE_2[\Psi_h])\nabla\CE_2[\Psi_h]\otimes\CE_1(\bv)
\end{equation}
with some matrix of $C^1$ functions, $\bV_\bg(\bk)$, defined on $|\bk|
< \delta$ such that $\|(\bV_g, \bV'_g)\|_{L_\infty(|\bk| \leq \delta)}
\leq C$, where $\bV'_\bg$ denotes the derivative of 
$\bV_\bg$ with respect to $\bk$. 
We may write the time derivative of $\tilde \bg(\bv, \Psi_h)$ as 
\begin{align*}
\pd_t\tilde\bg(\bv, \Psi_h) & = 
\bV_\bg(\nabla\CE_2[\Psi_h])\nabla\CE_2[\Psi_h]\otimes(\pd_t\CE_1(\bv))
\\
&+\bV_\bg(\nabla\CE_2[\Psi_h])(\pd_t\nabla\CE_2[\Psi_h])\otimes\CE_1(\bv)
\\
&+(\bV'_\bg(\nabla\CE_2[\Psi_h])\pd_t\nabla\CE_2[\Psi_h])
\nabla\CE_2[\Psi_h]\otimes\CE_1(\bv).
\end{align*}
and so by \eqref{g.5.8*} we have 
\begin{equation}\label{g.5.12*}\begin{split}
\|\pd_t\tilde\bg(\bv, \Psi_h)\|_{L_q(\Omega)}
&\leq C\{\|\nabla\CE_2[\Psi_h]\|_{H^1_q(\Omega)}
\|\pd_t\CE_1[\bv]\|_{L_q(\Omega)} \\
&\quad +\|\pd_t\nabla\CE_2[\Psi_h]\|_{L_q(\Omega)}
\|\CE_1[\bv]\|_{H^1_q(\Omega)}\}.
\end{split}\end{equation}
Thus, using \eqref{eT:2}, \eqref{sob:5.0},
\eqref{nonloc:1},  
\eqref{f.6.3}, \eqref{nonloc:0*}, \eqref{nonloc:0**}, 
 \eqref{g.5.7}, and \eqref{g.5.8}, for any $\gamma \geq 0$,  we have
\begin{align}
&\|e^{-\gamma t}
\pd_t\tilde\bg(\bv, \Psi_h)\|_{L_p(\BR, L_q(\Omega))}
\nonumber\\
&\quad
\leq C(\|\Psi_h\|_{L_\infty((0,T), H^2_q(\Omega))}
+ \|T_h(\cdot)\Psi_{\rho_0}\|_{L_\infty((0, \infty), H^2_q(\Omega))})
\nonumber\\
&\phantom{\leq C\{\|\Psi_h\|}
\times
(\|\pd_t\bv\|_{L_p((0, T), L_q(\Omega))}
+ \|T_v(\cdot)\bu_0\|_{H^1_p((0, \infty), L_q(\Omega))})
\nonumber \\
&\quad\quad +
T^{1/p}(\|\pd_t\Psi_h\|_{L_\infty((0, T), H^1_q(\Omega))}
+ \|\pd_tT_h(\cdot)\rho_0\|_{L_\infty((0, \infty), H^1_q(\Omega))})
\nonumber\\
&\phantom{\leq C\{\|\Psi_h\|}
\times
(\|\bv\|_{L_\infty((0, T), H^1_q(\Omega))} + 
\|T_v(\cdot)\bu_0\|_{L_\infty((0, \infty), H^1_q(\Omega))}) \nonumber\\
&\quad \leq C(LT^{1/{p'}} + \epsilon 
+ (L+\epsilon)T^{1/p} )(L+B). \nonumber \\
&\quad \leq C(\epsilon + L(T^{1/p'} + T^{1/p}))(L+B).
\label{g.5.12}
\end{align}

We next prove that
\begin{equation}\label{g.5.13}\begin{split}
&\|e^{-\gamma t}\tilde g(\bv, \Psi_h)\|_{H^{1/2}_p(\BR, L_q(\Omega))}
+ \|e^{-\gamma t}\tilde g(\bv, \Psi_h)\|_{L_p(\BR, H^1_q(\Omega))} \\
&\quad\leq C(\epsilon + T^{1/{p'}}L + T^{(q-N)/(pq)}L^{1+N/(2q)})(L+B)
\end{split}\end{equation}
for any $\gamma \geq 0$. 
In the sequel, to estimate $\|fg\|_{H^{1/2}_p(\BR, L_q(\Omega))}$, we use 
the following two lemmas. 
\begin{lem}\label{lem:g.5.0}
Let $1 < p < \infty$, $N < q < \infty$ and $0 < T \leq 1$. 
Let $f \in H^1_p(\BR, H^1_q(\Omega))$ and 
$g \in H^{1/q}_p(\BR, L_q(\Omega))
\cap L_p(\BR, H^1_q(\Omega))$.
If $f$ vanishes for $t \not\in (0, 2T)$,
then we have 
\begin{align*}
&\|fg\|_{H^{1/2}_p(\BR, L_q(\Omega))} + \|fg\|_{L_p(\BR, H^1_q(\Omega))}\\
&\leq C\{\|f\|_{L_\infty(\BR, H^1_q(\Omega))}
+T^{(q-N)/(pq)}\|\pd_tf\|_{L_\infty(\BR, L_q(\Omega))}^{1-N/(2q))}
\|\pd_tf\|_{L_p(\BR, H^1_q(\Omega))}^{N/(2q)}\}\\
&\quad\quad\times(
\|g\|_{H^{1/2}_p(\BR, L_q(\Omega))} + \|g\|_{L_p(\BR, H^1_q(\Omega))}).
\end{align*}
\end{lem}
\pf For a proof, see Shibata-Shimizu \cite[Lemma 2.6]{Shibata-Shimizu3}.
\qed
\begin{rem}\label{rem:g.5.0}
 We replace 
$\|f\|_{L_\infty(\BR, H^1_q(\Omega))}$  by 
$T^{1/{p'}}\|\pd_t f\|_{L_p(\BR, H^1_q(\Omega))}$
in Lemma \ref{lem:g.5.0}. In fact, 
since $f|_{t=0} =0$, representing $f(t) = \int^t_0\pd_sf(\cdot, s)\,ds$
and applying  H\"older's inequality, we have 
$$\|f\|_{L_\infty(\BR, H^1_q(\Omega))}
\leq T^{1/{p'}}\|\pd_t f\|_{L_p(\BR, H^1_q(\Omega))}.$$
\end{rem}
\begin{lem}\label{lem:g.5.1}
Let $1 < p < \infty$ and $N < q < \infty$. 
Let 
\begin{align*}
f &\in L_\infty(\BR, H^1_q(\Omega)) \cap
H^1_\infty(\BR, L_q(\Omega)), \\
g &\in
H^{1/2}_p(\BR, L_q(\Omega)) \cap L_p(\BR, H^1_q(\Omega)).
\end{align*}  Then, we have
\begin{align*}
&\|fg\|_{H^{1/2}_p(\BR, L_q(\Omega))} + \|fg\|_{L_p(\BR, H^1_q(\Omega))}
\\
&\leq C(\|f\|_{H^1_\infty(\BR, L_q(\Omega))}
+ \|f\|_{L_\infty(\BR, H^1_q(\Omega))})
(\|g\|_{H^{1/2}_p(\BR, L_q(\Omega))} + \|g\|_{L_p(\BR, H^1_q(\Omega))}).
\end{align*}
\end{lem}
\pf We can prove the lemma by using the complex interpolation:
\begin{align*}
&H^{1/2}_p(\BR, L_q(\Omega))\cap L_p(\BR, H^{1/2}_q(\Omega)) \\
&\quad = (L_p(\BR, L_q(\Omega)),
H^1_p(\BR, L_q(\Omega))\cap L_p(\BR, H^1_q(\Omega)))_{[1/2]}.
\end{align*}
\qed \vskip0.5pc 
To estimate $\|\nabla v\|_{H^{1/2}_p(\BR, L_q(\Omega))}$, 
we use the following lemma. 
\begin{lem}\label{lem:g.5.2} Let $1 <p, q < \infty$ and let $\Omega$
be a uniform $C^2$ domain.  Then, 
$$H^1_p(\BR, L_q(\Omega)) \cap L_p(\BR, H^2_q(\Omega))
\subset H^{1/2}_p(\BR, H^1_q(\Omega))$$
and 
$$\|u\|_{H^{1/2}_p(\BR, H^1_q(\Omega))}
\leq C\{\|u\|_{L_p(\BR, H^2_q(\Omega))}
+ \|\pd_tu\|_{L_p(\BR, L_q(\Omega))}\}.
$$
\end{lem}
\begin{rem}
This lemma was mentioned in Shibata-Shimizu \cite{SS2}
in the case that $\Omega$ is bounded and was proved
by Shibata \cite{S6} in the case that $\Omega$ is a uniform
$C^2$ domain. 
\end{rem}
Since $J_0(0) = 0$ and $\bV_0(0) = 0$, by \eqref{g.5.8*} we may write
\begin{equation}\label{g.5.10**}
\tilde g(\bv, \Psi_h) = V_g(\nabla\CE_2[\Psi_h])
\nabla \CE_2[\Psi_h]\otimes
\nabla \CE_1[\bv]
\end{equation}
with some matrix of $C^1$ functions, $V_g(\bk)$, defined 
on $|\bk| < \delta$ such that 
$\|(V_g, V'_g)\|_{L_\infty(|\bk| \leq \delta)} \leq C$, where 
$V'_g$ denotes the derivative of $V_g$ with respect to $\bk$. 
 By \eqref{sob:5.0}, \eqref{f.6.0},
\eqref{g.5.8*}, \eqref{eT:2*}, \eqref{g.5.8},
\eqref{nonloc:1},and \eqref{nonloc:0**}, 
we have
\begin{align*}
&\|V_g(\nabla \CE_2[\Psi_h])\nabla \CE_2[\Psi_h]
\|_{L_\infty(\BR, H^1_q(\Omega))}
\leq C\|\CE_2[\Psi_h]\|_{L_\infty(\BR, H^2_q(\Omega))}
\\
&\quad
\leq C(\|T_h(\cdot)\rho_0\|_{L_\infty(\BR, 
H^2_q(\Omega))}
+ \|\Psi_h\|_{L_\infty(\BR, H^2_q(\Omega))})
\leq C(\epsilon + T^{1/{p'}}L); \\
&\|\pd_t(V_g(\nabla \CE_2[\Psi_h])\nabla \CE_2[\Psi_h])
\|_{L_\infty(\BR, L_q(\Omega))}
\leq C\|\pd_t\CE_2[\Psi_h]\|_{L_\infty(\BR, H^1_q(\Omega))}
\\
&\quad \leq C(\epsilon + L ) \leq 2CL; \\
&\|\pd_t(V_g(\nabla \CE_2[\Psi_h])\nabla \CE_2[\Psi_h])
\|_{L_p(\BR, H^1_q(\Omega))}\\
&\quad 
\leq C\|(\pd_t\bar\nabla^2\CE_2[\Psi_h], \pd_t\bar\nabla\CE_2[\Psi_h]
\otimes\bar\nabla^2\CE_2[\Psi_h])\|_{L_p((\BR, L_q(\Omega))}
\\
&\quad\leq C(\|\pd_t\CE_2[\Psi_h]\|_{L_p(\BR, H^2_q(\Omega))}
+ \|\pd_t\CE_2[\Psi_h]\|_{L_\infty(\BR, H^1_q(\Omega))}
\|\CE_2[\Psi_h]\|_{L_p(\BR, H^3_q(\Omega))}\\
&\quad \leq C((\epsilon + L)+(\epsilon+L)^2) \leq 6CL^2.
\end{align*} 
Thus, by Lemma \ref{lem:g.5.1} and Lemma \ref{lem:g.5.2} we have 
\begin{align}
&\|e^{-\gamma t}\tilde g(\bv, \Psi_h)\|_{H^{1/2}_p(\BR, L_q(\Omega))}
+ 
\|e^{-\gamma t}\tilde g(\bv, \Psi_h)\|_{L_p(\BR, H^1_q(\Omega))}
\nonumber \\
&\quad \leq C(\epsilon + T^{1/{p'}}L 
+ T^{(q-N)/(pq)}L^{1+N/(2q)}) \label{g.5.13**}\\
&\phantom{\quad \leq C(\epsilon + T^{1/{p'}}L }\times
(\|e^{-\gamma t}\CE_1[\bv]\|_{H^1_p(\BR, L_q(\Omega))}
+ \|e^{-\gamma t}\CE_1[\bv]\|_{L_p(\BR, H^2_q(\Omega))})
\nonumber
\end{align}
for any $\gamma \geq 0$, which, combined with \eqref{nonloc:0*},
\eqref{nonloc:0**}, and \eqref{g.5.7}, leads to \eqref{g.5.13}.

We finally estimate $\bh'(\bv, \Psi_h)$ and 
$h_N(\bv, h)$ given in \eqref{58*} and \eqref{bdy:form}. 
In view of \eqref{58*}, we may  write 
$\bh'(\bv, \Psi_h)$ as
\begin{equation}\label{g.5.13***}
\bh'(\bv, \Psi_h) = 
\bV'_{\bh}(\bar\nabla\Psi_h)
\bar\nabla\Psi_h\otimes \nabla\bv
\end{equation}
with some matrix of $C^1$ functions, $\bV'_{\bh}(\bar\bk) 
= \bV'_{\bh}(y, \bar\bk)$, defined on 
$\overline{\Omega}\times\{\bar\bk \mid |\bar\bk| \leq \delta\}$
possessing the estimate : 
$$\sup_{|\bar\bk| \leq\delta}\|(\bV_\bh'(\cdot, \bar\bk), 
\pd_{\bar\bk}\bV_\bh'(\cdot, \bar\bk))\|_{H^1_\infty(\Omega)}
\leq C$$
with some constant $C$. Where, $\pd_{\bar\bk}\bV_\bh'$ 
denotes the derivative  of $\bV_\bh'$ with respect to $\bar\bk$.
We extend $\bh'(\bv, \Psi_h)$ to the whole 
time interval $\BR$ by letting
\begin{equation}\label{g.5.13*}
\tilde\bh'(\bv, \Psi_h) = 
\bV'_{\bh}(\bar\nabla\CE_2[\Psi_h]))
\bar\nabla\CE_2[\Psi_h]\otimes\nabla\CE_1[\bv].
\end{equation}
Employing the same argument as in the proof of \eqref{g.5.13}, 
for any $\gamma > 0$ 
we have
\begin{equation}\label{g.5.14}\begin{split}
&\|e^{-\gamma t}\tilde\bh'(\bv, \Psi_h)\|_{H^{1/2}_p(\BR, L_q(\Omega))}
+ \|e^{-\gamma t}\tilde\bh'(\bv,\Psi_\rho)\|_{L_p(\BR, H^1_q(\Omega))}\\
&\quad\leq 
 C(\epsilon + T^{1/{p'}}L + T^{(q-N)/(pq)}L^{1+N/(2q)})(L+B).
\end{split}\end{equation}


We finally consider  $h_N(\bv, \Psi_h)$. 
In \eqref{bdy:form} we may write
\begin{align*}
&-<\bn, \mu\bD(\bv)\bV_\Gamma(\bar\bk)\bar\bk> 
-<\bn, \mu(\CD_\bD(\bk)\nabla\bv)(\bn + \bV_\Gamma(\bar\bk)\bar\bk)>\\
&\quad = \bV_{\bh, N}(\bar\bk)\nabla\Psi_h\otimes\nabla\bv
\end{align*}
with some matrix of $C^1$ functions,
 $\bV_{\bh, N}(\bar\bk) = \bV_{\bh, N}(y, \bar\bk)$, 
defined on $\overline{\Omega}\times\{\bk \mid |\bk| \leq \delta\}$
such that  
$$\sup_{|\bar\bk| \leq\delta}\|(\bV_{\bh, N}(\cdot, \bar\bk),
\pd_{\bar\bk}\bV_{\bh, N}(\cdot, \bar\bk))\|_{H^1_\infty(\Omega)}
\leq C.
$$ 
Thus, we may write 
$h_N(\bv, \Psi_h)$ as
\begin{equation}\label{bdyfunc:2}h_N(\bv, \Psi_h) 
= \bV_{\bh, N}(\bar\nabla\Psi_h)\nabla\Psi_h\otimes\nabla\bv 
+ \sigma\bV'_\Gamma(\bar\nabla\Psi_h)\bar\nabla\Psi_h\otimes 
\bar\nabla^2\Psi_h,
\end{equation}
and so, we can define the extension of $h_N(\bv, \Psi_h)$ by letting
\begin{equation}\label{bdyfunc:2*}
\begin{aligned}
\tilde h_N(\bv, \Psi_h) &= 
\bV_{\bh, N}(\bar\nabla\CE_2[\Psi_h])
\bar\nabla\CE_2[\Psi_h]\otimes\nabla \CE_1[\bv] \\
&+ \sigma V^1_\Gamma(\bar\nabla\CE_2[\Psi_h])
\bar\nabla\CE_2[\Psi_h]\otimes\bar\nabla^2\CE_2[\Psi_h].
\end{aligned}\end{equation}
Using Lemma \ref{lem:g.5.0}, Lemma \ref{lem:g.5.1}, 
\eqref{g.5.7}, and \eqref{g.5.8}, we have
\begin{align*}
&\|\tilde h_N(\bv, \Psi_h)\|_{L_p(\BR, H^1_q(\Omega))}
+ \|\tilde h_N(\bv, \Psi_h)\|_{H^{1/2}_p(\BR, L_q(\Omega))}\\
&\quad\leq C(\epsilon + T^{1/{p'}}L 
+ T^{(q-N)/(pq)}L^{1+N/(2q)})\\
&\quad\times(\|e^{-\gamma t}\CE_1[\bv]\|_{L_p(\BR, H^2_q(\Omega))}
+ \|e^{-\gamma t}\CE_1[\bv]\|_{H^1_p(\BR, L_q(\Omega))}\\
&\qquad +
\|e^{-\gamma t}\bar\nabla^2\CE_2[\Psi_h]\|_{L_p(\BR, H^1_q(\Omega))}
+ \|e^{-\gamma t}\bar\nabla^2\CE_2[\Psi_h]\|_{H^{1/2}_p(\BR, L_q(\Omega))}).
\end{align*}
By  the fact that 
$H^1_p(\BR, L_q(\Omega)) \subset H^{1/2}_p(\BR, L_q(\Omega))$, we have
\begin{align*}
&\|e^{-\gamma t}\bar\nabla^2\CE_2[\Psi_h]\|_{H^{1/2}_p(\BR, 
L_q(\Omega))} 
\leq \|e^{-\gamma t}\bar\nabla^2\CE_2[\Psi_h]\|_{H^1_p(\BR, 
L_q(\Omega))} \leq C(\epsilon + L).
\end{align*}
Therefore, by \eqref{g.5.7}, \eqref{g.5.8},
\eqref{nonloc:0*}, and \eqref{nonloc:0**}, we have
\begin{equation}\label{g.5.18*}\begin{aligned}
&\|\tilde h_N(\bv, \Psi_h)\|_{L_p(\BR, H^1_q(\Omega))}
+ \|\tilde h_N(\bv, \Psi_h)\|_{H^{1/2}_p(\BR, L_q(\Omega))}\\
&\quad \leq C(\epsilon + T^{1/{p'}}L 
+ T^{(q-N)/(pq)}L^{1+N/(2q)})(B+L).
\end{aligned}\end{equation}

Let
\begin{align*}
\tilde E_{p,q,T}(\bu, \rho) &= \|\bu\|_{L_p((0, T), H^2_q(\Omega))}
+ \|\pd_t\bu\|_{L_p((0, T), L_q(\Omega))} \\
&+ \|\rho\|_{L_p((0, T), W^{3-1/q}_q(\Gamma))}
+ \|\pd_t\rho\|_{L_p((0, T), W^{2-1/q}_q(\Gamma))}.
\end{align*}
Notice that
$$E_{p,q,T}(\bu, \rho) = \tilde E_{p,q,T}(\bu, \rho) + \|\pd_t\rho\|_{L_\infty
((0, T), H^1_q(\Omega))}.
$$
Applying Corollary \ref{cor:max.1} and using the estimates 
\eqref{nonloc:2}, \eqref{nonloc:8**},
\eqref{g.5.12}, \eqref{g.5.14}, \eqref{g.5.18*}, 
and \eqref{nonloc:0**},  we have 
\begin{equation}\label{g.5.19}
\tilde E_{p,q,T}(\bu, \rho)
\leq Ce^{\gamma \kappa^{-b}T}D_\epsilon
\end{equation}
with
\begin{align*}
D_\epsilon &= B  
+ \kappa^{-b}\epsilon
+ T^{1/p}(L+B)^2 \\
&\quad + 
(\epsilon+ T^{1/{p'}}L)L 
+ (L+B\kappa^{-1/p}T^{1/p})(\epsilon + T^{1/{p'}}L)
\\
&\quad + L(L+B)(\epsilon + T^{1/{p'}}L
+ \kappa^{(1-1/s)/{p'}}
+ T^{(1-1/s)/{p'}}) \\
&\quad 
+(\epsilon + L(T^{1/{p'}} + T^{1/p})
+ T^{(q-N)/(pq)}L^{1+N/(2q)})(L+B)\}
\end{align*}
for some positive constants $\gamma$ and $C$ independent of
$B$, $\epsilon$ and $T$. Combining \eqref{real:7.1.1} and 
\eqref{real:7.1.2} with \eqref{g.5.19} yields that
\begin{equation}\label{g.5.19*}\begin{aligned}
\|\bu\|_{L_\infty((0, T), B^{2(1-1/p)}_{q,p}(\Omega))}
&\leq C(B + \tilde E_{p,q,T}(\bu, \rho)), \\
\|\rho\|_{L_\infty((0, T), B^{3-1/p-1/q}_{q,p}(\Gamma))}
&\leq C(\epsilon + \tilde E_{p,q,T}(\bu, \rho)).
\end{aligned}\end{equation}
Noting that $3-1/p-1/q > 2-1/q$, 
by the kinetic equation in Eq.\eqref{g.5.5} and \eqref{g.5.19*} 
 we have
\begin{align}
&\sup_{t \in (0, T)} \|\pd_t\rho(\cdot, t)\|_{W^{1-1/q}_q(\Gamma)}
\nonumber \\
&\quad
\leq \|\bu_\kappa\|_{H^1_q(\Omega)}\sup_{t\in (0, T)}
\|\rho(\cdot, t)\|_{W^{2-1/q}_q(\Gamma)} \nonumber \\
&\qquad +C\sup_{t\in(0, T)}\|\bu(\cdot,t)\|_{H^1_q(\Omega)}
+ \sup_{t \in(0,T)}\|d(\cdot, t)\|_{W^{1-1/q}_q(\Gamma)}
\nonumber\\
&\quad \leq CB(\epsilon + \tilde E_{p,q,T})
+ C(B + \tilde E_{p,q,T}(\bu, \rho)) + C(L+B)(\epsilon + T^{1/{p'}}L).
\label{g.5.20}
\end{align}
Since we may assume that $0 < \epsilon < 1$ and $B \geq 1$,  
combining \eqref{g.5.19} and \eqref{g.5.20} yields that 
\begin{equation}\label{g.5.21}
E_{p,q,T}(\bu, \rho) \leq A_1B + 
A_2(L+B)(\epsilon + T^{1/{p'}}L)
+ A_3Be^{\gamma \kappa^{-b}T}D_\epsilon.
\end{equation}
for some constants $A_1$, $A_2$ and $A_3$ independent of 
$B$, $\epsilon$, $\kappa$ and $T$. 
Let $\kappa = 
\epsilon =T$, and then 
we have $D_\epsilon = B + \epsilon^{1-b} + D'_\epsilon$ with 
\begin{align*}
D'_\epsilon &=  \epsilon^{1/p}(L+B)^2 
(\epsilon+ \epsilon^{1/{p'}}L)L 
+ (L+B)(\epsilon + \epsilon^{1/{p'}}L)
\\
&\quad + L(L+B)(\epsilon + \epsilon^{1/{p'}}L
+ 2\epsilon^{(1-1/s)/{p'}}) \\
&\quad+(\epsilon + L(\epsilon^{1/{p'}} + \epsilon^{1/p})
+ \epsilon^{(q-N)/(pq)}L^{1+N/(2q)})(L+B)\}.
\end{align*}
Choosing $\epsilon \in (0, 1)$
so small that 
\begin{gather*}
D'_\epsilon \leq B, \quad \epsilon^{1-b} \leq B, \quad 
\gamma\kappa^{-b}T = \gamma \epsilon^{1-b} \leq 1, \\
(L+B)(\epsilon + T^{1/{p'}}L) =(L+B)(\epsilon + \epsilon^{1/{p'}}L) \leq 2B
\end{gather*}
by \eqref{g.5.21}, we have
$E_{p,q,T}(\bu, \rho) \leq A_1B + 2A_2B + 3A_3eB^2$. 
Thus, setting $L = A_1B + 2A_2B + 3A_3eB^2$, we finally obtain 
\begin{equation}\label{g.5.22}
E_{p,q,T}(\bu, \rho) \leq L.
\end{equation}
Let $\CM$ be a map defined by letting 
$\CM(\bv, h) = (\bu, \rho)$, and then by 
\eqref{g.5.22} $\CM$ maps $\bU_T$ into itself. 
We can also prove that $\CM$ is a contraction map.
Namely, choosing $\kappa=\epsilon = T$ smaller
if necessary, we can show that for any $(\bv_i, h_i)
\in \bU_T$ ($i=1,2$), 
$$E_{p,q,T}(\CM(\bv_1, h_1) - \CM(\bv_2, h_2))
\leq (1/2)E_{p,q,T}((\bv_1, \rho_1) - (\bv_2, \rho_2)).
$$
Thus, by the contraction mapping principle, 
we have Theorem \ref{thm:loc.1}, which completes the 
proof of Theorem \ref{thm:loc.1}.

We finally prove the inequalities \eqref{real:7.1.1} and \eqref{real:7.1.2}. 
Let $\CE_1[\bv]$
be the function given in \eqref{eT:2*}.  By \eqref{eT:3} and 
\eqref{tanabe:2}, we have
\begin{align*}
&\|\bv\|_{L_\infty((0, T), B^{2(1-1/p)}_{q,p}(\Omega))} 
\leq \|\CE_1[\bv]\|_{L_\infty((0, T), B^{2(1-1/p)}_{q,p}(\Omega))} \\
& \leq C\{\|\CE_1[\bv]\|_{L_p((0, \infty), H^2_q(\Omega))} 
+ \|\pd_t\CE_1[\bv]\|_{L_p((0, \infty), L_q(\Omega))} \},
\end{align*}
which, combined with \eqref{g.5.7}, leads to the  inequality
\eqref{real:7.1.1}.  Analogously, using $\CE_2[\Psi_h]$ given
in \eqref{eT:2*} and \eqref{tanabe:2}, we have 
\begin{align*}
&\|\Psi_h\|_{L_\infty((0, T), B^{3-1/p-1/q}_{q,p}(\Gamma))}
\\
&\quad
\leq C\{\|\Psi_h\|_{L_p((0, T), H^3_q(\Omega))}
+ \|\pd_t\Psi_h\|_{L_p((0, T), H^2_q(\Omega))}\},
\end{align*}
which, combined with \eqref{g.5.8} and \eqref{sob:5.4}, leads to
the inequality in \eqref{real:7.1.2}.


\section{Global well-posedness in a bounded domain closed to a ball}
\label{sec:6}

In this section, we study the global well-posedness of 
Eq. \eqref{navier:1}.  As was stated in Subsec.\ref{subsec:2.4},
 we consider the problem
in the following setting.
Let $B_R$ be the ball of radius $R$ centered at the origin and let
$S_R$ be a sphere of radius $R$ centered at the origin.  We assume
that
\begin{itemize}
\item[\thetag{A.1}]\qquad
$|\Omega| = |B_R| = \dfrac{R^N\omega_N}{N}$, 
where $|D|$ denotes the Lebesgue measure of a Lebesgue measurable set 
$D$ in $\BR^N$ and $\omega_N$ is the area of $S_1$.
\item[\thetag{A.2}]\qquad
$\displaystyle{\int_\Omega x\,dx = 0}$. 
\item[\thetag{A.3}]\qquad
$\Gamma$ is a normal perturbation of $S_R$ given
by
$$\Gamma = \{x = y+\rho_0(y)\bn(y) \mid y \in S_R\}$$
with a given small function $\rho_0(y)$ defined on $S_R$.
\end{itemize}
Here, $\bn = y/|y|$ is the unit outer normal to $S_R$. $\bn$ is extended to
$\BR^N$ by $\bn = R^{-1}y$. 
Let $\Gamma_t$ be given by 
\begin{equation}\label{g.1.2}\begin{aligned}
\Gamma_t &= \{x = y + \rho(y, t)\bn + \xi(t) \mid y \in S_R\} \\
&= \{x = y + R^{-1}\rho(y, t)y + \xi(t) \mid y \in S_R\}
\end{aligned}\end{equation}
where $\rho(y, t)$ is an unknown function with $\rho(y,0) = \rho_0(y)$
for $y \in S_R$.  Let $\xi(t)$ be the barycenter point of the 
domain $\Omega_t$ defined by 
$$\xi(t) = \frac{1}{|\Omega|}\int_{\Omega_t} x\,dx.
$$
Here, by \eqref{mass:1}, we have $|\Omega_t| = |\Omega|$.  Notice that $\xi(t)$ is also unknown.  It follows from
the assumption \thetag{A.2} that $\xi(0) = 0$.  Moreover,
by \eqref{mom:3} we have 
\begin{equation}\label{g.1.3}
\xi'(t) = \frac{1}{|\Omega|}\int_{\Omega_t} \bv(x, t)\,dx.
\end{equation}
Given any function $\rho$ defined on $S_R$, let $H_\rho$ be a suitable extension
of $\rho$ such that $H_\rho=\rho$ on $S_R$ and \eqref{sob:5.4}
holds.
We define the Hanzawa transform centered at $\xi(t)$ by 
\eqref{g.1.5}. 
In the following, we set $\Psi_\rho = R^{-1}H_\rho y$ and we assume that
\begin{equation}\label{g.1.6}
\sup_{t \in (0, T)} \|\Psi_\rho(\cdot, t)\|_{H^1_\infty(B_R)} \leq \delta.
\end{equation}
We will choose $\delta$ so small that several  conditions
hold.  For example, since $|\Psi_\rho(y, t) - \Psi_\rho(y', t)|
\leq \|\nabla \Psi_\rho(\cdot, t)\|_{L_\infty(B_R)}|y-y'|$, 
if $0 < \delta < 1/4$, then 
$$|\bh_\rho(y, t) - \bh_\rho(y', t)| \geq (1-2\delta)|y-y'| 
\geq(1/2)|y-y'|
\quad \text{for any $y$, $y' \in B_R$},
$$
and so the Hanzawa transform is a bijective map from $B_R$ onto 
$\Omega_t$ with
\begin{equation}\label{g.1.7}
\Omega_t = \{x = \bh_\rho(y, t) \mid y\in B_R\}
\quad\text{for $t \in (0, T)$}.
\end{equation}
Let  $\bv_0(x)$ be an initial data for Eq. \eqref{navier:1}, 
and then we set $\bu_0(y) = \bv_0(\bh_{\rho_0}(y))$,
where $\bh_{\rho_0}(y) = y + R^{-1}H_{\rho_0}y$.  Notice that 
$\bh_{\rho_0}(y) = \bh_{\rho}(y, 0)$ if $\rho|_{t=0} = \rho_0$. 
Let $\bv$ and $\fp$ satisfy Eq. \eqref{navier:1} and we set  
\begin{equation}\label{g.1.9}
\bu(y, t) = \bv(\bh_\rho(y, t), t),
\quad \fq(y, t) = \fp(\bh_\rho(y, t), t) - \frac{\sigma(N-1)}{R}.
\end{equation}
We then see from the consideration in Sect. \ref{sec:2} that
$\bu$, $\fq$ and $\rho$ satisfy the following equations:
\begin{equation}\label{g.1.10}\left\{\begin{aligned}
\pd_t\bu - \DV(\mu\bD(\bu) - \fq\bI) & = \bff(\bu, \Psi_\rho)
&\quad&\text{in $B_R^T$}, \\
\dv\bu = g(\bu, H_\rho) &= \dv\bg(\bu, \Psi_\rho)
&\quad&\text{in $B_R^T$}, \\
\pd_t\rho - \bn\cdot P\bu & = \tilde d(\bu, \Psi_\rho)
&\quad&\text{on $S_R^T$}, \\
\BPi_0(\mu\bD(\bu)\bn) & = \bh'(\bu, \Psi_\rho)
&\quad&\text{on $S_R^T$}, \\
<\mu\bD(\bu)\bn, \bn> -\fq
-\sigma\CB\rho & = h_N(\bu, \Psi_\rho)
&\quad&\text{on $S_R^T$}, \\
(\bu, \rho)|_{t=0} & = (\bu_0, \rho_0)
&\quad&\text{in $B_R\times S_R$},
\end{aligned}\right.\end{equation}
where $B_R^T = B_R\times(0, T)$, 
$S_R^T = S_R\times(0, T)$, $\bn = y/|y|$ for $y \in S_R$, $\displaystyle{
P\bu = \bu - \frac{1}{|B_R|}\int_{B_R}\bu\,dy}$, $\CB = \Delta_{S_R}
+ \dfrac{N-1}{R^2}$, and $\Delta_{S_R}$ is the Laplace-Beltrami
operator on $S_R$.  Where, the functions in the right side of \eqref{g.1.10}, 
 $\bff(\bu, \Psi_\rho)$, 
$g(\bu, \Psi_\rho)$, $\bg(\bu, \Psi_\rho)$, and 
$\bh'(\bu, \Psi_\rho)$ 
are nonlinear terms given in \eqref{form:f}, \eqref{form:g}, and 
\eqref{58*}, respectively.  And, $h_N(\bu, \Psi_\rho)$  
is given in the formula \eqref{non:g.5} of Subsec. \ref{subsec.7.1} below. Since $dx = dy + J_0(\bk)dy$ with 
$$J_0(\bk) = \det \left(\begin{matrix} 
\frac{\pd \Psi_{\rho 1}(y, t)}{\pd y_1} & \cdots &
 \frac{\pd \Psi_{\rho 1}(y, t)}
{\pd y_N} \\
\vdots & \ddots &\vdots \\
\frac{\pd \Psi_{\rho N}(y, t)}{\pd y_1} &
\cdots &  \frac{\pd \Psi_{\rho N}(y, t)}
{\pd y_N}
\end{matrix}\right)
$$
for $\Psi_\rho= {}^\top(\Psi_{\rho1}, \ldots, \Psi_{\rho N})$, 
by \eqref{g.1.3}
we have
\begin{equation}\label{g.1.3*}
\xi'(t) = \frac{1}{|B_R|}\int_{B_R}\bu(y, t)\,dy 
+ \frac{1}{|B_R|}\int_{B_R}\bu(y, t)J_0(\bk)\,dy,
\end{equation}
where $\bk = \nabla \Psi_\rho$. 
Thus, by \eqref{kinematic:4} and \eqref{kinematic:4*}, we have
$\pd_t\rho - \bn\cdot P\bu = \tilde d(\bu, \Psi_\rho)$, 
where $\tilde d(\bu, \Psi_\rho)$ is given by letting 
\begin{equation}\label{kin:1}
\tilde d(\bu, \Psi_\rho) = d(\bu, \Psi_\rho) -<\bu  \mid \nabla'_\Gamma \rho>
+ \frac{1}{|\Omega|}\int_{B_R} \bu(y, t) J_0(\bk) \,dy,
\end{equation}
with $d(\bu, \Psi_\rho)$ given in \eqref{kinematic:4*}. 

We now state our main result of this section. For this 
purpose, we make several definitions.  From the assumptions
\thetag{A.1} and \thetag{A.2}  it follows that 
$\rho_0$ should satisfy the following conditions:
\begin{align*}
&\frac{R^N\omega_N}{N} = \int_\Omega\,dx 
= \int_{|\omega|=1}\int^{R+\rho_0(R\omega)}_0r^{N-1}\,dr\,d\omega
= \int_{|\omega|=1}\frac{(R+\rho_0(R\omega))^N}{N}\,d\omega, \\
&0 = \int_\Omega x_i\,dx = \int_{|\omega|=1}\int^{R+\rho_0(R\omega)}_0
\omega_ir^N\,dr\,d\omega = \int_{|\omega|}
\frac{(R+\rho_0(R\omega))^{N+1}}{N+1}\omega_i\,d\omega
\end{align*}
for $i=1, \ldots, N$, and so, we have the compatibility conditions
for $\rho_0$ as follows:
\begin{equation}\label{g.1.12}
\sum_{k=1}^N {}_NC_k\int_{S_R}(R^{-1}\rho_0(y))^k\,d\omega = 0, \quad
\sum_{k=1}^{N+1}{}_{N+1}C_k\int_{S_R} y_i(R^{-1}\rho_0(y))^k
\,d\omega=0
\end{equation}
where ${}_NC_k = \frac{N!}{k!(N-k)!}$
and $d\omega$ denotes the surface element of $S_R$, because 
$\int_{S_1}\,d\omega = \omega_N$ and $\int_{S_1}\omega_i\,d\omega = 0$.
Let $\CR_d = \{\bu \mid \bD(\bu) = 0\}$, and then $\CR_d$ is a 
finite dimensional vector space spanned by  
constant $N$-vectors and first order polynomials $x_i\be_j-x_j\be_i$
($i,j=1, \ldots, N$), where $\be_i$ are $N$-vector whose $i^{\rm th}$
component is $1$ and other components are zero. Let $M_d$ be the 
dimension of $\CR_d$ and  let 
$\bp_\ell = |B_R|^{-1/2}\be_\ell$ ($\ell=1, \ldots, N$), and 
$\bp_\ell$ ($\ell=N+1, \ldots, M_d$) be one of 
$c_R^1(x_i\be_j - x_j\be_i)$ with $c_R^1 =\sqrt{(N+2)/(2R^2|B_R|)}
$. Since $(\bp_\ell, \bp_m)_{B_R}
= \delta_{\ell, m}$ for $\ell$, $m = 1, \ldots, M_d$, the set
$\{\bp_\ell\}_{\ell=1}^{M_d}$ forms a orthogonal basis of 
$\CR_d$ with respect to the $L_2(B_R)$ inner-product
$(\cdot, \cdot)_{B_R}$.

We know that $\CB x_i = 0$ on $S_R$ for $i = 1, \ldots, N$.  
Setting $\varphi_1 = |S_R|^{-1/2}$ and $\varphi_\ell = 
c^2_Rx_\ell$ ($\ell=2, \ldots, N+1$) with
$c^2_R = \sqrt{N/(R^{N+1}\omega_N)}$,  
we see that  $(\varphi_i, \varphi_j)_{S_R} = \delta_{ij}$
$(i, j=1, \ldots, N+1)$, and therefore  
the set $\{\varphi_i\}_{i=1}^{N+1}$ forms an orthogonal basis 
of the space $\{\psi \mid \CB\psi = 0 \enskip\text{on $S_R$}\}
\cup \BC$ with respect to the $L_2(S_R)$ inner-product 
$(\cdot, \cdot)_{S_R}$.  In this section, we set  
\begin{align}
\CS_{p,q}((a, b)) &= \{(\bu, \fq, \rho) \mid \bu
\in H^1_p((a, b), L_q(B_R)^N) \cap L_p((a, b), 
H^2_q(B_R)^N), \nonumber  \\
&\fq \in L_p((a, b), H^1_q(S_R)), \nonumber\\
&\rho \in H^1_p((a, b), W^{2-1/q}_q(S_R))
\cap L_p((a, b), W^{3-1/q}_q(S_R))\}. \label{func:4.1}
\end{align}
Our main result of this section is the following.
\begin{thm}\label{thm:sec:5} Let $p$ and $q$ be real numbers
such that $2 < p < \infty$, $N < q < \infty$ and 
$2/p + N/q < 1$.  Assume that \thetag{A.1}, \thetag{A.2} and
\thetag{A.3} hold. Then, there exists a small number $\epsilon > 0$
such that if 
initial data $\bu_0 \in B^{2(1-1/p)}_{q,p}(B_R)$ and 
$\rho_0 \in B^{3-1/p-/q}_{q,p}(S_R)$ satisfy the smallness
condition: 
\begin{equation}\label{g.1.13}
\|\bu_0\|_{B^{2(1-1p)}_{q,p}(B_R)} + \|\rho_0\|_{B^{3-1/p-1/q}_{q,p}
(S_R)} \leq \epsilon, 
\end{equation}
the compatibility conditions \eqref{g.1.12} and 
\begin{equation}\label{g.1.14}\begin{split}
\dv\bu_0 = g(\bu_0, \rho_0) = \dv \bg(\bu_0, \rho_0)
\quad&\text{in $B_R$}, \\
(\mu\bD(\bu_0)\bn)_\tau = \bg'(\bu_0, \rho_0)
\quad&\text{on $S_R$},
\end{split}\end{equation}
and the orthogonal condition:
\begin{equation}\label{g.1.15}
(\bv_0, \be_i)_{\Omega} = 0, \quad
(\bv_0, x_i\be_j - x_j\be_i)_\Omega=0,
\end{equation}
for $i, j=1, \ldots, N$, 
then problem \eqref{g.1.10} with $T = \infty$ admits a unique solution
$(\bu, \fq, \rho) \in \CS_{p,q}((0, \infty))$ possessing the estimate:
\begin{align*}
&\|e^{\eta t}\pd_t\bu\|_{L_p((0, \infty), L_q(B_R))}
+ \|e^{\eta t}\bu\|_{L_p((0, \infty), H^2_q(B_R))}
+ \|e^{\eta t}\nabla \fq\|_{L_p((0, \infty), L_q(B_R))} \\
&\quad
+ \|e^{\eta t}\pd_t\rho\|_{L_p((0, \infty), W^{2-1/q}_q(S_R))}
+ \|e^{\eta t}\rho\|_{L_p((0, \infty), W^{3-1/q}_q(S_R))}
\leq C\epsilon
\end{align*}
for some positive constants $C$ and $\eta$ independent of $\epsilon$.
\end{thm}
\subsection{Derivation of 
nonlinear term $h_N(\bu, \Psi_\rho)$ in Eq. \eqref{g.1.10}}
\label{subsec.7.1}

In this subsection, we derive the nonlinear term 
$h_N(\bu, \Psi_\rho)$ in \eqref{g.1.10}. 
Let $\omega \in S_1$ be represented by 
$\omega =\omega(p_1, \ldots, p_{N-1})$ with a local coordinate
$(p_1, \ldots, p_{N-1})$, and then for $x = 
(R+\rho)\omega +\xi(t)\in \Gamma_t$, we have 
$$\frac{\pd x}{\pd p_j}= (R+\rho)\tau_j
+ \frac{\pd\rho}{\pd p_j}\omega
$$
where $\tau_j = \frac{\pd \omega}{\pd p_j}$. Since 
$\tau_j\cdot \omega = 0$, 
the $(i, j)^{\rm th}$ component of the first fundamental form 
$G_t$ of $\Gamma_t$ is given by 
$$g_{tij}=\frac{\pd x}{\pd p_i}\cdot\frac{\pd x}{\pd p_j}
= (R+\rho)^2g_{ij} + \frac{\pd \rho}{\pd p_i}\frac{\pd \rho}{\pd p_j},$$
where $g_{ij} = \tau_i\cdot\tau_j$ are the $(i, j)^{\rm th}$ elements of 
the  first fundamental form, G,  of $S_1$, and so 
\begin{align*}
G_t &= (R+\rho)^2(G+(R+\rho)^{-2}\nabla_p\rho\otimes
\nabla_p\rho) \\
&= (R+\rho)^2G(\bI + (R+\rho)^{-2}(G^{-1}\nabla_p\rho)\otimes \nabla_p\rho),
\end{align*}
where $\nabla_p\rho = {}^\top(\pd \rho/ \pd p_1, \ldots, \pd\rho/\pd p_{N-1})$. 
Since
\begin{equation}\label{38} 
\det(\bI + \ba'\otimes\bb') = 1 + \ba'\cdot\bb',
\quad
(\bI + \ba'\otimes\bb')^{-1} = \bI - 
\frac{\ba'\otimes\bb'}{1+\ba'\cdot\bb'}
\end{equation}
for any $N-1$ vectors $\ba'$ and $\bb' \in \BR^{N-1}$, we have
\begin{align*}
G_t^{-1}&=(R+\rho)^{-2}
\Bigl(\bI - \frac{(R+\rho)^{-2}(G^{-1}\nabla_p\rho)\otimes \nabla_p\rho}
{1 + (R+\rho)^{-2}<G^{-1}\nabla_p\rho, \nabla_p\rho>}
\Bigr)G^{-1} \\
&= (R+\rho)^{-2}G^{-1} + O_2.
\end{align*}
Here and in the following,  $O_2$ denote the same 
symbol as in \eqref{residue:1}.  Namely, 
$$g_t^{ij} = (R+\rho)^{-2}g^{ij} + O_2, 
$$
componentwise.

We next calculate the Christoffel symbols of $\Gamma_t$. Since
\begin{align*}
\tau_{ti} &= (R+\rho)\tau_i + \frac{\pd\rho}{\pd p_i}\omega, \\
\tau_{tij} &= (R+\rho)\tau_{ij} + 
\frac{\pd\rho}{\pd p_j}\tau_i + 
\frac{\pd\rho}{\pd p_i}\tau_j
+ \frac{\pd^2\rho}{\pd p_i\pd p_j}\omega,
\end{align*}
we have
\begin{align*}
<\tau_{tij}, \tau_{t\ell}> &= (R+\rho)^2<\tau_{ij}, \tau_\ell>
+ (R+\rho)(\frac{\pd\rho}{\pd p_\ell}\ell_{ij}
+ g_{i\ell}\frac{\pd\rho}{\pd p_j}
+ g_{j\ell}\frac{\pd\rho}{\pd p_i})\\
&+ \frac{\pd^2\rho}{\pd p_i\pd p_j}\frac{\pd \rho}{\pd p_\ell},
\end{align*}
and so
\begin{align*}
\Lambda^k_{tij}  &= g_t^{k\ell}<\tau_{tij}, \tau_{t\ell}> \\
&=<(R+\rho)^{-2}g^{k\ell} + O_2, (R+\rho)^2<\tau_{ij}, \tau_\ell>\\
&+ (R+\rho)(\frac{\pd\rho}{\pd p_\ell}\ell_{ij}
+ g_{i\ell}\frac{\pd\rho}{\pd p_j}
+ g_{j\ell}\frac{\pd\rho}{\pd p_i}) + \frac{\pd^2\rho}{\pd p_i\pd p_j}\frac{\pd \rho}{\pd p_\ell}>\\
& = \Lambda^k_{ij} + (R+\rho)^{-1}(g^{k\ell}\frac{\pd\rho}{\pd p_\ell}\ell_{ij}
+ \delta^k_i\frac{\pd\rho}{\pd p_j}
+ \delta^k_j\frac{\pd\rho}{\pd p_i}\\
&+((R+\rho)^{-2}g^{k\ell}\frac{\pd\rho}{\pd p_\ell}
+ O_2)\frac{\pd^2\rho}{\pd p_i \pd p_j}
+ O_2.
\end{align*}
Thus, 
\begin{align*}
&\Delta_{\Gamma_t}f  = g_t^{ij}(\pd_i\pd_jf- \Lambda^k_{tij}\pd_kf) \\
& =(R+\rho)^{-2}g^{ij}(\pd_i\pd_jf - \Lambda^k_{ij}\pd_kf)
+(A^k\nabla_p^2\rho +O_2)\pd_kf
\end{align*}
with
\begin{align*}
A^k\nabla_p\rho = 
((R+\rho)^{-4}g^{k\ell}g^{ij}\frac{\pd\rho}{\pd p_\ell} +O_2)\frac{\pd^2 \rho}
{\pd p_i\pd p_j},
\end{align*}
and so
\begin{align*}
&H(\Gamma_t)\bn_t = \Delta_{\Gamma_t}[(R+\rho)\omega + \xi(t)] \\
& =(R+\rho)^{-2}g^{ij}(\pd_i\pd_j - \Lambda^k_{ij}\pd_k)((R+\rho)\omega)
+ (A^k\nabla^2_p\rho + O_2)\pd_k((R+\rho)\omega)\\
& = (R+\rho)^{-1}g^{ij}(\pd_i\pd_j\omega-\Lambda^k_{ij}\pd_k\omega)
+ (R+\rho)^{-2}g^{ij}(\pd_i\rho\pd_j\omega + \pd_j\rho\pd_i\omega)\\
&\quad + (R+\rho)^{-2}g^{ij}(\pd_i\pd_j\rho-\Lambda^k_{ij}\pd_k\rho)\omega
+ ((A^k\nabla^2_p\rho + O_2)\pd_k\rho)\omega \\
&\quad +(A^k\nabla^2_p\rho + O_2)(R+\rho)\pd_k\omega.
\end{align*}
Combining this formula with \eqref{repr:2.1}, recalling $\bn=\omega$ in
this case and using $<\pd_i\omega, \omega> = 0$ gives 
\begin{align*}
&<H(\Gamma_t)\bn_t\, \bn_t> 
=<H(\Gamma_t)\bn_t, \omega- g^{k\ell}(\pd_\ell\rho)\tau_k +O_2>\\
& = (R+\rho)^{-1}<\Delta_{S1}\omega,  \omega> 
+(R+\rho)^{-2}\Delta_{S_1}\rho 
+(A^k\nabla^2_p\rho + O_2)\pd_k\rho\\
&+ (R+\rho)^{-1}g^{k\ell}\pd_\ell\rho<\Delta_{S1}\omega, \tau_\ell>
+(A^m\nabla^2_p\rho)\pd_m\rho + O_2\nabla^2_p\rho + O_2.
\end{align*}
Since $\Delta_{S_1}\omega= -(N-1)\omega$, we have
\begin{align*}
<H(\Gamma_t)\bn_t, \bn_t> 
&= -(R+\rho)^{-1}(N-1) + (R+\rho)^{-2}
\Delta_{S_1}\rho \\
&\quad +(A^m\nabla^2_p\rho)\pd_m\rho + O_2\nabla^2_p\rho + O_2.
\end{align*}
Since
\begin{align*}
(R+\rho)^{-1}=&R^{-1} - \rho R^{-2} + O(\rho^2), \\
(R+\rho)^{-2}\Delta_{S_1}\rho &= R^{-2}\Delta_{S_1}\rho 
+2R^{-3}\rho\Delta_{S_1}\rho + O_2\nabla^2_p\rho,
\end{align*}
we have 
\begin{equation}\label{44}\begin{aligned}
&<H(\Gamma_t)\bn_t, \bn_t> \\
&\quad = -\frac{N-1}{R} + \CB\rho 
+ 
2R^{-3}\rho \Delta_{S_1}\rho+(A^m\nabla^2_p\rho)\pd_m\rho + O_2\nabla^2_p\rho 
 + O_2.
\end{aligned}\end{equation}
Replacing the pressure term $\fq$ by $\fq + \sigma\frac{N-1}{R}$, 
we have
$$<\mu\bD(\bu)\bn, \bn> - \fq - \sigma\CB\rho
= h_N(\bu, \Psi_\rho)
$$
on $S^T_R$. In view of  \eqref{bdyfunc.1}, \eqref{2.5.2*}, 
and \eqref{bdyfunc:2},  $h_N(\bu, \Psi_\rho)$ may be defined by letting
\begin{equation}\label{non:g.5}
h_N(\bu, \Psi_\rho) = \bV_{h,N}(\bar\nabla\Psi_\rho)
\bar\nabla\Psi_\rho\otimes\nabla\bv 
+ \sigma \tilde \bV'_\Gamma(\bar\nabla\Psi_\rho)\bar\nabla\Psi_\rho\otimes
\bar\nabla^2\Psi_\rho,
\end{equation}
where $\bV_{h, N}(\bar\bk)$ and $\tilde\bV'_\Gamma(\bar\bk)$ are 
functions defined on $\Omega\times\{\bar\bk \mid |\bar\bk| \leq\delta\}$
possessing the estimate:
\begin{align*}
&\sup_{|\bar\bk|\leq\delta}\|(\bV_{h, N}(\cdot, \bar\bk), 
\pd_{\bar\bk}\bV_{h, N}(\cdot, \bar\bk))\|_{H^1_\infty(\Omega)}
\leq C, \\
&\sup_{|\bar\bk|\leq\delta}\|(\tilde\bV'_{\Gamma}(\cdot, \bar\bk), 
\pd_{\bar\bk}\tilde\bV'_{\Gamma}(\cdot, \bar\bk))\|_{H^1_\infty(\Omega)}
\leq C,
\end{align*}
for some constant $C$.

\subsection{Local well-posedness} \label{subsec:5.2}
In this subsection, we prove the local well-posedness of Eq. \eqref{g.1.10}.
\begin{thm}\label{thm:g.5.1} Let $N < q < \infty$, $2 < p < \infty$ and 
$T > 0$.  Assume that $2/p + N/q < 1$. Then, there exists a constant
$\epsilon > 0$ depending on $T$ such that if initial data $\bu_0
\in B^{2(1-1/p)}_{q,p}(B_R)$ and $\rho_0 \in B^{3-1p-1q}_{q,p}(S_R)$
satisfy the smallness condition:
\begin{equation}\label{g.5.2} \|\bu_0\|_{B^{2(1-1/p)}_{q,p}(B_R)}
+ \|\rho_0\|_{B^{3-1/p-1/q}_{q,p}(S_R)} \leq \epsilon^2
\end{equation}
and the compatibility condition:
\begin{equation}\label{g.5.3} \dv\bu_0 = g(\bu_0, \rho_0) 
\quad\text{in $B_R$}, \quad 
\BPi_0(\mu\bD(\bu_0)\bn) = \mu\bh'(\bu_0, \rho_0)
\quad\text{on $S_R$},
\end{equation}
then problem \eqref{g.1.10} admits unique solutions $(\bu, \fq, \rho)
\in \CS_{p,q}((0, T))$ possessing the estimates:
\begin{align*}
\sup_{0 < t < T}\|\Psi_\rho(\cdot, t)\|_{H^1_\infty(B_R)} \leq \delta,
\quad E_{p,q, T}(\bu, \rho) \leq \epsilon.
\end{align*}
Here and in the following,  for $\eta \in \BR$ and $0 \leq a < b \leq \infty$,
we set
\begin{align*}
E_{p,q, T}(\bu, \rho)
 &= \|\pd_t\bu\|_{L_p((0, T), L_q(B_R))} 
+ \|\bu\|_{L_p((0, T), H^2_q(B_R))}\\
&+
\|\pd_t\rho\|_{L_p((0, T), W^{2-1/q}_q(S_R))} 
+ \|\rho\|_{L_p((0, T), W^{3-1/q}_q(S_R))}\\
&
+ \|\pd_t\rho\|_{L_\infty((0, T), W^{1-1/q}_q(S_R)}).
\end{align*}
\end{thm}

To prove Theorem \ref{thm:g.5.1} the key tool is the maximal $L_p$-$L_q$
regularity for the following linear equations:
\begin{equation}\label{g.5.4} \left\{\begin{aligned}
\pd_t\bu - \DV(\mu\bD(\bu) - \fq\bI) & = \bff
&\quad&\text{in $B_R^T$}, \\
\dv\bu= g &=\dv\bg 
&\quad&\text{in $B_R^T$}, \\
\pd_t\rho - \bn\cdot P\bu & = d 
&\quad&\text{on $S_R^T$}, \\
\mu\bD(\bu) - \fq\bI)\bn - \sigma(\CB\rho)\bn & = \bh
&\quad&\text{on $S_R^T$}, \\
(\bu,\rho)|_{t=0} & = (\bu_0, \rho_0)
&\quad&\text{in $B_R\times S_R$}.
\end{aligned}\right.
\end{equation}
Setting $\CF_1\bu= \displaystyle{\frac{1}{|B_R|}\int_{B_R} \bu\,dx}$,
and $\CF_2\rho = \sigma R^{-2}(N-1)\rho$, 
by Theorem \ref{thm:max.1} and Theorem \ref{thm:max.2}, we have the following 
theorem.
\begin{thm}\label{thm:g.5.2} let $1 < p, q < \infty$ and 
$2/p + 1/q \not=0$. Let $T > 0$.  Then, there exists a $\gamma_0>0$
such that the following assertion holds:
Let $\bu_0 \in B^{2(1-1/p)}_{q,p}(B_R)^N$
and $\rho_0 \in B^{1-1/p-1/q}_{q,p}(S_R)$ be initial data for Eq. 
\eqref{g.5.4} and let $\bff$, $g$, $\bg$, $d$, $\bh$ be given functions in
the right side of Eq. \eqref{g.5.4} with 
\allowdisplaybreaks{
\begin{align*}
&\bff \in L_p((0, T), L_q(B_R)^N), \quad 
e^{-\gamma t}g \in H^1_p(\BR, H^1_q(B_R)) \cap H^{1/2}_p(\BR, L_q(B_R)), \\
&e^{-\gamma t}\bg \in H^1_p(\BR, L_q(B_R)), \quad 
d \in L_p((0, T), W^{2-1/q}_q(S_R)), \\
&e^{-\gamma t}\bh 
\in H^1_p(\BR, H^1_q(B_R)^N) \cap H^{1/2}_p(\BR, L_q(B_R)^N)
\end{align*}
}
for all $\gamma \geq \gamma_0$.  Assume that the compatibility condition:
$\dv\bu_0 = g|_{t=0}$ in $B_R$ holds. 
In addition, the compatibility condition:
$\BPi_0(\mu\bD(\bu_0)) = \BPi_0(\bh|_{t=0})$ on $\Gamma$
holds 
provided $2/p + 1/q < 1$.  Then, problem \eqref{g.5.4} admits
unique solutions $(\bu, \fq, \rho) \in \CS_{p,q}((0, \infty))$ 
possessing the estimate:
\begin{align}
&\|\pd_t\bu\|_{L_p((0, T) L_q(B_R))} 
+ \|\bu\|_{L_p((0, T), H^2_q(B_R))}
+
\|\pd_t\rho\|_{L_p((0, T), W^{2-1/q}_q(S_R))} \nonumber  \\
&+ \|\rho\|_{L_p((0, T), W^{3-1/q}_q(S_R))} 
\leq Ce^{\gamma T}(\|\bu_0\|_{B^{2(1-1/p)}_{q,p}(B_R)}
+ \|\rho_0\|_{B^{3-1/p-1/q}_{q,p}(S_R)}  \nonumber \\
&+ \|\bff\|_{L_p((0, T), L_q(B_R))}
+ \|e^{-\gamma t}(g, \bh)\|_{L_p(\BR, H^1_q(B_R))}
+ \|e^{-\gamma t}(g, \bh)\|_{H^{1/2}_p(\BR, L_q(B_R))}\nonumber \\
&+ \|e^{-\gamma t}\pd_t\bg\|_{L_p(\BR, L_q(B_R))}
+ \|d\|_{L_p((0, T), W^{2-1/q}_q(S_R)})
\label{max:est:7.1}
\end{align}
for any $\gamma \geq \gamma_0$ with some constant $C$ independent of 
$\gamma$.  
\end{thm}

{\bf Proof of Theorem \ref{thm:g.5.1}}.  In what follows, using the 
Banach fixed point argument, we prove Theorem \ref{thm:g.5.1}. 
Let $\CU_{\epsilon}$ be the underlying space defined by 
\begin{align}
\CU_{\epsilon} &= \{(\bu, \rho) \mid \bu \in H^1_p((0, T), L_q(B_R)^N)
\cap L_p((0, T), H^2_q(B_R)^N), 
\nonumber   \\
&\rho \in H^1_q((0, T), W^{2-1/q}_q(S_R)) \cap L_p((0, T), W^{3-1/q}_q(S_R)),
\nonumber \\
&\bu|_{t=0} = \bu_0 \quad\text{in $B_R$},\quad 
\rho|_{t=0} = \rho_0 \quad\text{on $S_R$},  \nonumber \\
&\sup_{0< t < T}\|\Psi_\rho(\cdot, t)\|_{H^1_\infty(B_R)} \leq \delta,
\quad E_{p,q, T}(\bu, \rho) \leq \epsilon\}.
\label{eq:under:1}
\end{align}
We recall that $\Psi_\rho = H_\rho\bn$ and $H_\rho$ is a suitable extension
of $\rho$ to $B_R$ such that $H_\rho = \rho$ on $S_R$ 
and \eqref{sob:5.4} holds. 
Since $N < q < \infty$ and $2/p + N/q < 1$, we have 
\begin{equation}\label{sob:5.3}\begin{split}
\sup_{t \in (0, T)}\|f(\cdot, t)\|_{H^1_\infty(B_R)} 
&\leq C_{p,q}\sup_{t \in (0, T)}\|f(\cdot, t)\|_{B^{2(1-1/p)}_{q,p}(B_R)},
\\
\sup_{t \in (0, T)}\|f(\cdot, t)\|_{H^2_\infty(B_R)} 
&\leq C_{p,q}\sup_{t \in (0, T)}\|f(\cdot, t)\|_{B^{3-1/p}_{q,p}(B_R)}.
\end{split}\end{equation} 
Let $(\bv, h) \in \CU_{\epsilon}$. We then consider the 
linear equations:
\begin{equation}\label{*g.5.5} \left\{\begin{aligned}
\pd_t\bu - \DV(\mu\bD(\bu) - \fq\bI) & = \bff(\bv, \Psi_h)
&\quad&\text{in $B_R^T$}, \\
\dv\bu= g(\bv, H_h) &=\dv\bg(\bv, \Psi_h) 
&\quad&\text{in $B_R^T$}, \\
\pd_t\rho - \bn\cdot P\bu & = \tilde d(\bv, \Psi_h)
&\quad&\text{on $S_R^T$}, \\
(\mu\bD(\bu) - \fq\bI)\bn - \sigma(\CB\rho)\bn & = \bh(\bv, \Psi_h)
&\quad&\text{on $S_R^T$}, \\
(\bu,\rho)|_{t=0} & = (\bu_0, \rho_0)
&\quad&\text{in $B_R\times S_R$}.
\end{aligned}\right.
\end{equation}Noting \eqref{real:7.1.1} and \eqref{real:7.1.2}, 
we have
\begin{equation}\label{small:7.1}\begin{aligned}
&\|\bv\|_{L_\infty((0, T), B^{2(1-1/p)}_{q,p}(B_R))} \leq C\epsilon^2, \\
&\|\rho\|_{L_\infty((0, T), B^{3-1/p-1/q}_{q,p}(S_R))} \leq C\epsilon^2, \\
& \sup_{t\in (0, T)}\|\Psi_h(\cdot, t)\|_{H^1_\infty(B_R)} \leq \delta, \\
&\|\pd_t\bv\|_{L_p((0, T), L_q(B_R))} + 
\|\bv\|_{L_p((0, T), H^2_q(B_R))} \\
&\quad+\|\pd_th\|_{L_p((0, T), W^{2-1/q}_q(S_R))} + 
\|h\|_{L_p((0, T), W^{3-1/q}_q(S_R))}\\ 
&\quad + \|\pd_th\|_{L_\infty((0, T), W^{1-1/q}_q(S_R))} 
\leq \epsilon^2.
\end{aligned}\end{equation}
Employing the same argument as in the proof of Theorem \ref{thm:loc.1}, 
we have
\begin{equation}\label{loc:7.1}\begin{split}
\|\bff(\bv, \Psi_h)\|_{L_q(B_R)} & \leq C\{(\|\bv(\cdot, t)\|_{H^1_q(B_R)}
+ \|\pd_t\Psi_h(\cdot, t)\|_{H^1_q(B_R)})
\|\bv(\cdot, t)\|_{H^1_q(B_R)} \\
&+ \|\Psi_h(\cdot, t)\|_{H^2_q(B_R)}
(\|\pd_t\bv(\cdot, t)\|_{L_q(B_R)}
+ \|\bv(\cdot, t)\|_{H^2_q(B_R)})\}.
\end{split}\end{equation}

Putting \eqref{loc:7.1}, \eqref{real:7.1.1}, \eqref{real:7.1.2},
\eqref{eq:under:1}, and \eqref{small:7.1},  we have
\begin{equation}\label{loc:7.2}
\|\bff(\bv, \Psi_h)\|_{L_p((0,m T), L_q(B_R))} \leq C\epsilon^2.
\end{equation}
Here and in the following, $C$ denotes a generic  constant independent of 
$\epsilon$.  Moreover, since we choose $\epsilon > 0$ small enough eventually,
we may assume that $0 < \epsilon <1$, and so $\epsilon^s \leq \epsilon^2$ for 
any $s \geq 2$.

We next consider $\tilde d(\bv, \Psi_h)$.  In view of \eqref{g.1.3*},
for $(\bv, h) \in \CU_\epsilon$ 
$$\xi'(t) = \frac{1}{|B_R|}\int_{B_R}\bv(y, t)\,dy + 
\frac{1}{|B_R|}\int_{B_R}\bv(y, t)\bJ_0(\bk)\,dy.
$$
Here and in the following, $\bk = \nabla\Psi_h$.  Thus, noting 
\eqref{small:7.1}, we have 
$$\sup_{t\in(0, T)}|\xi'(t)| \leq \sup_{t\in(0, T)}\|\bv(\cdot, t)
\|_{L_q(B_R)} \leq C\epsilon,
$$
and so by \eqref{small:7.1} and \eqref{sob:5.4.b} we have
\begin{align*}
&\|<\xi'(t)\mid\nabla'_\Gamma h>\|_{L_\infty((0, T), W^{1-1/q}_q(S_R))} \\
&\quad \leq C\epsilon\|h\|_{L_\infty((0, T), W^{2-1/q}_q(S_R))}
\leq C\epsilon^2, \\
&\|<\xi'(t)\mid V_\Gamma(\bar\bk)\bar\nabla\Psi_h\otimes
\bar\nabla\Psi_h>\|_{L_\infty((0, T), W^{1-1/q}_q(S_R))} \\
&\quad \leq C\epsilon\|h\|_{L_\infty((0, T), W^{2-1/q}_q(S_R))}^2
\leq C\epsilon^3,\\
&\|<\xi'(t)\mid\nabla'_\Gamma h>\|_{L_p((0, T), W^{2-1/q}_q(S_R))} \\
&\quad \leq C\epsilon\|h\|_{L_p((0, T), W^{3-1/q}_q(S_R))}
\leq C\epsilon^2, \\
&\|<\xi'(t)\mid V_\Gamma(\bar\bk)\bar\nabla\Psi_h\otimes
\bar\nabla\Psi_h>\|_{L_p((0, T), W^{2-1/q}_q(S_R))} \\
&\quad \leq C\epsilon(1+\|h\|_{L_\infty((0, T), W^{2-1/q}_q(S_R))})
\|h\|_{L_p((0, T), W^{3-1/q}_q(S_R))}
\leq C\epsilon^2(1+\epsilon). 
\end{align*}
Estimating the rest of $d(\bv, h)$ given in \eqref{kinematic:4*} 
defined by replacing $\bu$ and $\rho$ by $\bv$ and $h$ in the similar
manner to the proof of \eqref{nonloc:3}, we have
$$\|d(\bv, \Psi_h)\|_{L_\infty((0, T), W^{1-1/q}_q(S_R))} \leq C\epsilon^2,
\enskip \|d(\bv, \Psi_h)\|_{L_p((0, T), W^{2-1/q}_q(S_R))} \leq C\epsilon^2.
$$
 Moreover,
by \eqref{sob:5.0}, \eqref{real:7.1.1}, \eqref{real:7.1.2},  and
\eqref{small:7.1}, we have
\begin{align*}
&\|<\bv\mid\nabla'_\Gamma h> \|_{L_\infty((0, T), W^{1-1/q}_q(S_R))}
\\
&\quad \leq C\|\bv\|_{L_\infty((0, T), H^1_q(B_R))}
\|h\|_{L_\infty((0, T), W^{2-1/q}_q(S_R))}
\leq C\epsilon^2; \\
&\|<\bv\mid\nabla'_\Gamma h> \|_{L_p((0, T), W^{2-1/q}_q(S_R))}\\
&\quad \leq C\{\|\bv\|_{L_p((0, T), H^2_q(B_R))}
\|h\|_{L_\infty((0, T), W^{2-1/q}_q(S_R))} \\
&\quad+ \|\bv\|_{L_\infty((0, T), H^1_q(B_R))}
\|h\|_{L_p((0, T), W^{3-1/q}_q(S_R))}\}
\leq C\epsilon^2.
\end{align*}
And also, 
$$\int_{B_R}|\bv(y, t)J_0(\bk)|\,dy
\leq |B_R|^{1/q'}\|\bv(\cdot, t)\|_{L_q(B_R)}
\|\nabla\Psi_h(\cdot, t)\|_{L_\infty(B_R)},
$$
 and so by \eqref{small:7.1} we have
\begin{align*}
\Bigl\|\int_{B_R}\bv J_0\,dy\Bigr\|_{L_\infty((0, T), W^{1-1/q}_q(S_R))}
&\leq C\epsilon^2, \\ 
\Bigl\|\int_{B_R}\bv J_0\,dy\Bigr\|_{L_p((0, T), W^{2-1/q}_q(S_R))}
&\leq C\epsilon^2.
\end{align*}
Combining estimates obtained above gives
\begin{equation}\label{loc:7.3}\begin{split}
\|\tilde d(\bv, \Psi_h)\|_{L_\infty((0, T), W^{1-1/q}_q(S_R))} 
&\leq C\epsilon^2, \\
\|\tilde d(\bv, \Psi_h)\|_{L_p((0, T), W^{2-1/q}_q(S_R))} 
&\leq C\epsilon^2.
\end{split}\end{equation}

We now consider $g(\bv, \Psi_h)$ and $\bg(\bv, \Psi_h)$.
Let $\tilde g(\bv, \Psi_h)$ and $\tilde \bg(\bv, \Psi_h)$
be the extension of  $g(\bv, \Psi_h)$ and $\bg(\bv, \Psi_h)$
to the whole time line $\BR$ given in \eqref{g.5.10*} and
\eqref{g.5.10**}, respectively. 
By \eqref{g.5.12*} we have 
\begin{align*}
&\|e^{-\gamma t}\pd_t\tilde\bg(\bv, \Psi_h)\|_{L_p(\BR, L_q(B_R))}
\\
&\quad \leq C\{\|\CE_2[\Psi_h]\|_{L_\infty(\BR, H^2_q(B_R))}
\|e^{-\gamma t}\pd_t\CE_1[\bv]\|_{L_p(\BR, L_q(B_R))}\\
&\quad + \|e^{-\gamma t}\CE_2[\Psi_h]\|_{L_p(\BR, H^2_q(B_R))}
\|\CE_1[\bv]\|_{L_\infty(\BR, H^1_q(B_R))}\},
\end{align*}
and so by \eqref{g.5.7}, \eqref{g.5.8}, and \eqref{small:7.1}, we have
\begin{equation}\label{loc:7.4}
\|e^{-\gamma t}\pd_t\tilde \bg(\bv, \Psi_h)\|_{L_p(\BR, L_q(B_R))}
\leq C\epsilon^2.
\end{equation}
Analogously, by \eqref{g.5.7}, \eqref{g.5.8}, Lemma \ref{lem:g.5.1},
Lemma \ref{lem:g.5.2}, \eqref{eq:under:1}, \eqref{real:7.1.1}, 
and \eqref{real:7.1.2},
we have
\begin{equation}\label{loc:7.5}
\|e^{-\gamma t}\tilde g(\bv, \Psi_h)\|_{H^{1/2}_p(\BR, L_q(B_R))}
+ \|e^{-\gamma t}\tilde g(\bv, \Psi_h)\|_{L_p(\BR, H^1_q(B_R))}
\leq C\epsilon^2.
\end{equation}

We next consider $\bh'(\bv, \Psi_h)$. Let $\tilde\bh'(\bv, \Psi_h)$ 
be the extension of $\bh'(\bv, \Psi_h)$ to the whole time line
$\BR^N$ given by \eqref{g.5.13*}.  By Lemma \ref{lem:g.5.1} and 
Lemma \ref{lem:g.5.2}, we have
\begin{align*}
&\|e^{-\gamma t}\tilde \bh'(\bv, \Psi_h)\|_{H^{1/2}_p(\BR, L_q(B_R))}
+ \|e^{-\gamma t}\tilde\bh'(\bv, \Psi_h)\|_{L_p(\BR, H^1_q(\Omega))} \\
&\quad 
\leq C(\|\bar\nabla\CE_2[\Psi_h]\|_{L_\infty(\BR, H^1_q(B_R))}
+ \|\pd_t\bar\nabla\CE_2[\Psi_h]\|_{L_\infty(\BR, L_q(B_R))})\\
&\qquad \times(\|e^{-\gamma t}\CE_1[\bv]\|_{L_p(\BR, H^2_q(B_R))}
+ \|e^{-\gamma t}\pd_t\CE_1[\bv]\|_{L_p(\BR, L_q(B_R))}).
\end{align*}
Thus, by \eqref{g.5.7}, \eqref{g.5.8}, and \eqref{small:7.1}, we have
\begin{equation}\label{loc:7.6}
\|e^{-\gamma t}\tilde\bh'(\bv, \Psi_h))\|_{H^{1/2}_p(\BR, L_q(B_R))}
+ \|e^{-\gamma t}\tilde\bh'(\bv, \Psi_h)\|_{L_p(\BR, H^1_q(B_R))}
\leq C\epsilon^2.
\end{equation}

We finally consider $h_N(\bv, \Psi_h)$ given in \eqref{non:g.5}.
Let  $\tilde h_N(\bv, \Psi_h)$ be the extension 
of $h_N(\bv, \Psi_h)$ to the whole time interval $\BR$ given by 
\begin{equation}\label{non:g.5*}\begin{aligned}
\tilde h_N(\bv, \Psi_h) & = \bV_{h,N}(\bar\nabla\CE_2[\Psi_h]) 
(\bar\nabla\CE_2[\Psi_h], \nabla\CE_1[\bv]) \\
&+\sigma\tilde\bV'_\Gamma(\bar\nabla\CE_2[\Psi_h])
(\bar\nabla\CE_2[\Psi_h], \bar\nabla^2\CE_2[\Psi_h]).
\end{aligned}\end{equation}
By Lemma \ref{lem:g.5.1} and Lemma \ref{lem:g.5.2}, we have
\begin{align*}
&\|e^{-\gamma t}\tilde h_N(\bv, \Psi_h)\|_{H^{1/2}_p(\BR, L_q(B_R))}
+ \|e^{-\gamma t}\tilde h_N(\bv, \Psi_h)\|_{L_p(\BR, H^1_q(B_R))}
\\
&\quad \leq C(\|\pd_t\bar\nabla\CE_2[\Psi_h]\|_{L_\infty(\BR, L_q(B_R))}
+ \|\bar\nabla\CE_2[\Psi_h]\|_{L_\infty(\BR, H^1_q(B_R))})\\
&\qquad \times (\|e^{-\gamma t}\pd_t\CE_1[\bv]\|_{L_p(\BR, L_q(B_R))} 
+ \|e^{-\gamma t}\CE_1[\bv]\|_{L_p(\BR, H^2_q(B_R))} \\
&\qquad\qquad+ \|e^{-\gamma t}\pd_t\CE_2[\Psi_h]\|_{L_p(\BR, H^2_q(B_R))}
+ \|e^{-\gamma t}\CE_2[\Psi_h]\|_{L_p(\BR, H^3_q(B_R))}).
\end{align*}
Thus, by \eqref{g.5.7}, \eqref{g.5.8}, and \eqref{small:7.1}, we have
\begin{equation}\label{loc:7.7}
\|e^{-\gamma t}\tilde h(\bv, \Psi_h)\|_{H^{1/2}_p((\BR, L_q(B_R))}
+ \|e^{-\gamma t}\tilde h(\bv, \Psi_h)\|_{L_p((\BR, H^1_q(B_R))}
\leq C\epsilon^2.
\end{equation}
Applying Theorem \ref{thm:g.5.2} to Eq. \eqref{*g.5.5} and 
using \eqref{loc:7.2}, \eqref{loc:7.3}, \eqref{loc:7.4},
\eqref{loc:7.5}, \eqref{loc:7.6}, and  \eqref{loc:7.7}, we have
\begin{align}
&\|\bu\|_{L_p((0, T), H^2_q(B_R))} 
+ \|\pd_t\bu\|_{L_p((0, T), L_q(B_R))} \nonumber  \\
&+\|\rho\|_{L_p((0, T), W^{3-1/q}_q(S_R))} 
+ \|\pd_t\rho\|_{L_p((0, T), W^{2-1/q}_q(S_R))}
\leq Ce^{\gamma T}\epsilon^2\label{loc:7.8}
\end{align}
for some positive constants $C$ and $\gamma$.  Moreover, by the third equation
in \eqref{*g.5.5}, we have
\begin{equation}\label{loc:7.8.7}\begin{aligned}
&\|\pd_t\rho\|_{L_\infty((0, T), W^{1-1/q}_q(S_R))} \\
&\quad
\leq C(\|\bu\|_{L_\infty((0, T), H^1_q(B_R))}
+ \|\tilde d(\bv, \Psi_h)\|_{L_\infty((0, T), W^{1-1/q}_q(S_R))}),
\end{aligned}\end{equation}
which, combined with \eqref{real:7.1.1}, \eqref{real:7.1.2},
 and  \eqref{loc:7.3}, yields 
\begin{align*}
&\|\pd_t\rho\|_{L_\infty((0, T), W^{1-1/q}_q(S_R))} 
\leq C(\|\bu_0\|_
{B^{2(1-1/p)}_{q,p}(B_R))} \\
&\quad+ \|\bu\|_{L_p((0, T), H^2_q(B_R))} 
+ \|\pd_t\bu\|_{L_p((0, T), L_q(B_R))} + C\epsilon^2).
\end{align*}
Combining this inequality with \eqref{loc:7.8} gives
\begin{equation}\label{loc:7.9}
E_{p,q,T}(\bu, \rho) \leq Ce^{\gamma T}\epsilon^2
\end{equation}
for some positive constants $C$ and $\gamma$. 

Let $\CP$ be a map defined by $\CP(\bv, h) = (\bu, \rho)$, and then 
choosing $\epsilon$ so small that $Ce^{\gamma T}\epsilon\leq 1$, 
by \eqref{loc:7.9} we see that $E_{p,q,T}(\CP(\bv, h))
 \leq \epsilon$, which shows that $\CP$ maps $\CU_{\epsilon}$
into itself.  Let $(\bv_i, h_i)$ ($i = 1,2$) be any two elements 
of $\CU_{\epsilon}$, and then we have
$$E_{p,q,T}(\CP(\bv_1, h_1) - \CP(\bv_2, h_2))
\leq Ce^{\gamma T}\epsilon E_{p,q,T}((\bv_1, h_1)-(\bv_2, h_2))
$$
for some positive constants $C$ and $\gamma$, where we have used 
$E_{p,q,T}((\bv_i, \rho_i)) \leq \epsilon$ for $i=1, 2$. 
Choosing $\epsilon$ so small that 
$Ce^{\gamma T}\epsilon \leq 1/2$, we see that $\CP$ is a contraction
map from $\CU_{\epsilon}$ into itself, and so by the Banach fixed 
point theorem, there exists a unique $(\bu, \rho) \in 
\CU_{\epsilon}$ satisfying $\CP(\bu, \rho) = (\bu, \rho)$.
Thus, $(\bu, \rho)$ is a required unique soluton of Eq. \eqref{g.1.10},
which completes the proof of Theorem \ref{thm:g.5.1}. 

\subsection{Decay estimates of solutions for the linearized equations}
\label{subsec:5.3}

To prove Theorem \ref{thm:sec:5} the key tool is decay properties
of solutions of the Stokes equations:
\begin{equation}\label{*g.5.4} \left\{\begin{aligned}
\pd_t\bu - \DV(\mu\bD(\bu) - \fp\bI) & = \bff
&\quad&\text{in $B_R^T$}, \\
\dv\bu= g &=\dv\bg 
&\quad&\text{in $B_R^T$}, \\
\pd_t\rho - \bn\cdot P\bu & = d 
&\quad&\text{on $S_R^T$}, \\
(\mu\bD(\bu) - \fp\bI)\bn - \sigma(\CB\rho)\bn & = \bh
&\quad&\text{on $S_R^T$}, \\
(\bu,\rho)|_{t=0} & = (\bu_0, \rho_0)
&\quad&\text{in $B_R\times S_R$}.
\end{aligned}\right.
\end{equation}
Here and in the following, $\bn = y/|y| \in S_1$. 
We will prove the following theorem. 
\begin{thm}\label{thm:g.5.2*} Let $1 < p, q < \infty$ and 
$2/p + 1/q \not=0$. Let $\{\bp_\ell\}_{\ell=1}^M$ 
and $\{\varphi_j\}_{j=1}^{N+1}$
be orthogonal basis of $\CR_d = \{\bu \mid \bD(\bu) = 0\}$ and 
$\{\psi \mid \CB\psi=0 \quad\text{on $S_R$}\}\cup \BC$, which 
are given before Theorem \ref{thm:sec:5}.   Then, there exists an $\eta>0$
such that the following assertion holds:
Let $\bu_0 \in B^{2(1-1/p)}_{q,p}(B_R)^N$
and $\rho_0 \in B^{3-1/p-1/q}_{q,p}(S_R)$ be initial data for Eq. 
\eqref{*g.5.4} and let $\bff$, $g$, $\bg$, $d$, $\bh$ be given functions in
the right side of Eq. \eqref{*g.5.4}.  Assume that  
$$\bff \in L_p((0, T), L_q(B_R)^N), \quad 
d \in L_p((0, T), W^{2-1/q}_q(S_R)),
$$ and that there exist $\tilde g_\eta$, $\tilde \bg_\eta$ 
and $\tilde \bh_\eta$ such that $e^{\eta t}g = \tilde g_\eta$,
$e^{\eta t}\bg = \tilde \bg_\eta$, $e^{\eta t}\bh
= \tilde\bh_\eta$ for $t \in (0, T)$, $\dv \tilde g_\eta 
= \tilde g_\eta$ for $t \in \BR$, and 
\begin{align*}
&\tilde g_\eta  \in H^1_p(\BR, H^1_q(B_R)) \cap H^{1/2}_p(\BR, L_q(B_R)),
\quad
\tilde \bg_\eta \in H^1_p(\BR, L_q(B_R)), \\
&\tilde\bh_\eta  
\in H^1_p(\BR, H^1_q(B_R)^N) \cap H^{1/2}_p(\BR, L_q(B_R)^N).
\end{align*}
 Assume that the compatibility condition:
$\dv\bu_0 = g|_{t=0}$ in $B_R$ holds. 
In addition, the compatibility condition:
$(\mu\bD(\bu_0)\bn)_\tau = \bh_\tau|_{t=0}$ on $\Gamma$
holds 
provided $2/p + 1/q < 1$.  Then, problem \eqref{*g.5.4} admits
unique solutions $(\bu, \fp, \rho) \in \CS_{p,q}((0, T))$ 
possessing the estimate:
\begin{align}
\CI_{p,q, T}(\bu, \rho; \eta) \nonumber 
& \leq C\Bigl\{\CJ_{p,q, T}(\bu_0, \rho_0, \bff, g, \bg, d, \bh;
\eta) \nonumber \\
&\quad + \sum_{\ell=1}^M\Bigl(\int^T_0(e^{\eta s}
|(\bu(\cdot, s), \bp_\ell)_{B_R}|)^p\,ds\Bigr)^{1/p}
\nonumber\\
&+ \sum_{j=1}^{N+1}\Bigl(\int^T_0(e^{\eta s}
|(\rho(\cdot, s), \varphi_j)_{S_R}|)^p\,ds\Bigr)^{1/p}\,\Bigr\} 
\label{g.ineq:1}\end{align}
for some constant $C$ independent of 
$\eta$. 
Here and in the following, we set 
\begin{align*}
&\CI_{p,q, T}(\bu, \rho; \eta)  = \|e^{\eta t}\bu\|_{L_p((0, T), H^2_q(B_R))}
+ \|e^{\eta t}\pd_t\bu\|_{L_p((0, T), L_q(B_R))} \\
&\quad 
+ \|e^{\eta t}\rho\|_{L_p((0, T), W^{3-1/q}_q(S_R))}
+ \|e^{\eta t}\pd_t\rho\|_{L_p((0, T), W^{2-1/q}_q(S_R))}; \\
&\CJ_{p,q, T}(\bu_0, \rho_0, \bff, g, \bg, d, \bh; \eta) 
= \|\bu_0\|_{B^{2(1-1/p)}_{q,p}(B_R)}
+ \|\rho_0\|_{B^{3-1/p-1/q}_{q,p}(S_R)} \\
&+ \|e^{\eta t}\bff\|_{L_p((0, T), L_q(B_R))} 
+ \|e^{\eta t}d\|_{L_p((0, T), W^{2-1/q}_q(S_R))}
+ \|\pd_t\bg_\eta\|_{L_p(\BR, L_q(B_R))} \\
&+ \|(g_\eta, \bh_\eta)\|_{H^{1/2}_p(\BR, L_q(B_R))}
+ \|(g_\eta, \bh_\eta)\|_{L_p(\BR, H^1_q(B_R))}
\end{align*}
\end{thm}
%
To prove Theorem \ref{thm:g.5.2},  we first consider the
following shifted equations:
\begin{equation}\label{shift:1}\begin{cases*}
&\pd_t\bu_1 + \lambda_1\bu_1 - \DV(\mu\bD(\bu_1) - \fp_1\bI) = \bff
\quad&\text{in $B_R^T$}, \\
&\dv\bu_1 = g = \dv\bg
\quad&\text{in $B_R^T$},\\
&\pd_t\rho_1  + \lambda_1\rho_1- \bn\cdot P\bu_1
= d
\quad&\text{on $S_R^T$},\\
&(\mu\bD(\bu_1)- \fp_1\bI)\bn  -\sigma 
(\CB\rho_1)\bn = \bh 
\quad&\text{on $S_R^T$}, \\
&(\bu_1, \rho_1)|_{t=0} = (\bu_0, \rho_0)
\quad&\text{on $B_R\times S_R$}.
\end{cases*}\end{equation}
For the shifted equation \eqref{shift:1}, we have 
\begin{thm}\label{thm:g.3.2} Let $1 < p, q < \infty$ and $T > 0$.
Assume that $2/p + 1/q \not=1$. Let 
$\lambda_0$ be  
a constant given  in Theorem \ref{thm:rbdd:1.2}.
Let $\bu_0 \in B^{2(1-1/p)}_{q,p}(B_R)^N$
and $\rho_0 \in B^{1-1/p-1/q}_{q,p}(S_R)$ be initial data for Eq. 
\eqref{shift:1} and let $\bff$, $g$, $\bg$, $d$, $\bh$ be given functions in
the right side of Eq. \eqref{shift:1} satisfying the same condition 
as in Theorem \ref{thm:g.5.2*}. Moreover, there exist
$\tilde g_0$, $\tilde \bg_0$ and $\tilde \bh_0$ such that
$g= \tilde g_0$, $\bg = \tilde\bg_0$ and $\bh= \tilde \bh_0$ for 
$t \in (0, T)$, $\dv \tilde\bg_0 = \tilde g_0$ for $t \in \BR$,  and 
\begin{align*}
&\tilde g_0 \in H^1_p(\BR, H^1_q(B_R)) \cap H^{1/2}_p(\BR, L_q(B_R)),  
\quad \tilde \bg_0 \in H^1_p(\BR, L_q(B_R)), \\
&\tilde \bh_0 
\in H^1_p(\BR, H^1_q(B_R)^N) \cap H^{1/2}_p(\BR, L_q(B_R)^N).
\end{align*}
 Assume that the compatibility condition:
$\dv\bu_0 = g|_{t=0}$ in $B_R$ holds. 
In addition, the compatibility condition:
$(\mu\bD(\bu_0)\bn)_\tau = \bh_\tau|_{t=0}$ on $\Gamma$
holds 
provided $2/p + 1/q < 1$.  Then, 
for any $\lambda_1 > \lambda_0$, problem \eqref{shift:1} admits
unique solutions $(\bu_1, \fp_1, \rho_1) \in \CS_{p,q}((0, T))$ 
possessing the estimate:
$$
\CI_{p,q, T}(\bu_1, \rho_1; 0)
\leq C\CJ_{p,q,T}(\bu_0, \rho_0, \bff, g, \bg, d, \bh; 0)  
$$
for some constant $C$ independent of $T$.  
\end{thm}
\pf
Let $\CA_2(\lambda)$, $\CP_2(\lambda)$ and $\CH_2(\lambda)$ be 
operators given in Theorem \ref{thm:rbdd:1.2}.   Let 
 $\lambda_1$ and $\lambda_2$ be numbers for which  $\lambda_1 - \lambda_0
> \lambda_2 > 0$. 
Then, for any $\lambda \in -\lambda_2 + \Sigma_\epsilon$  
we have $\lambda + \lambda_1 \in 
\lambda_0 + \Sigma_\epsilon$, and so
 by Theorem \ref{thm:rbdd:1.2} we have 
\begin{align*}
\CR_{\CL(\CX_q(\Omega), H^{2-j}_q(\Omega)^N)}
(\{(\tau\pd_\tau)^\ell(\lambda^{j/2}\CA_2(\lambda+\lambda_1))\mid
\lambda \in -\lambda_2+\Sigma_{\epsilon_0}\}) &\leq r_b; \\
\CR_{\CL(\CX_q(\Omega), L_q(\Omega)^N)}
(\{(\tau\pd_\tau)^\ell(\nabla\CP_2(\lambda+\lambda_1))\mid
\lambda \in -\lambda_2+\Sigma_{\epsilon_0}\}) &\leq r_b; \\
\CR_{\CL(\CX_q(\Omega), W^{3-k}_q(\Gamma))}
(\{(\tau\pd_\tau)^\ell(\lambda^k\CH_2(\lambda+\lambda_1))\mid
\lambda \in -\lambda_2+\Sigma_{\epsilon_0}\}) &\leq r_b.
\end{align*}
for $\ell=0,1$, $j=0,1,2$, and $k=0,1$. 
Thus, we can choose $\gamma=0$ in the argument on 
Subsec. \ref{subsec:3.5}, and so we have the theorem. 
\qed \vskip0.5pc
For any $\eta > 0$, $e^{\eta t}\bu_1$, $e^{\eta t}\fp_1$ and 
$e^{\eta t}\rho_1$ satisfy the equations: 
\allowdisplaybreaks{
\begin{align*}
\left.\begin{aligned}
\pd_t(e^{\eta t}\bu_1) + (\lambda_1-\eta)e^{\eta t}\bu_1 - 
\DV(\mu\bD(e^{\eta t}\bu_1) - e^{\eta t}\fp_1\bI) = e^{\eta t}\bff& \\
\dv e^{\eta t}\bu_1 = e^{\eta t}g = \dv e^{\eta t}\bg&
\end{aligned}\right\}
\enskip\text{in $B_R^T$},&\\
\left. \begin{aligned}
\pd_t(e^{\eta t}\rho_1)  + (\lambda_1-\eta)e^{\eta t}\rho_1- 
\bn\cdot P(e^{\eta t}\bu_1)
&= e^{\eta t}d \\
(\mu\bD(e^{\eta t}\bu_1)- e^{\eta t}\fp_1\bI)\bn)  -\sigma 
(\CB(e^{\eta t}\rho_1))\bn &= e^{\eta t}\bh 
\end{aligned}\right\}\enskip
\text{on $S_R^T$},& \\
(e^{\eta t}\bu_1, e^{\eta t}\rho_1)|_{t=0} = (\bu_0, \rho_0)
\quad\text{on $B_R\times S_{R}$}.&
\end{align*}
} 
Given $\eta > 0$, we choose $\lambda_1>0$  in such a way that 
$\lambda_1 -\lambda_0> \eta > 0$, and then
 by Theorem \ref{thm:g.3.2}, we have the following
corollary.
\begin{cor}\label{cor:g.3.2} Let $1 < p, q < \infty$, $T > 0$
and $\eta > 0$. 
Assume that $2/p + 1/q \not=1$. 
Let $\bu_0 \in B^{2(1-1/p)}_{q,p}(B_R)^N$
and $\rho_0 \in B^{1-1/p-1/q}_{q,p}(S_R)$ be initial data for Eq. 
\eqref{shift:1} and let $\bff$, $g$, $\bg$, $d$, $\bh$ be given functions in
the right side of Eq. \eqref{shift:1} satisfying the same conditions
as in Theorem \ref{thm:g.5.2*}.  
 Assume that the compatibility condition:
$\dv\bu_0 = g|_{t=0}$ in $B_R$ holds. 
In addition, the compatibility condition:
$(\mu\bD(\bu_0)\bn)_\tau = \bh_\tau|_{t=0}$ on $\Gamma$
holds 
provided $2/p + 1/q < 1$.  Then, there exists a
$\lambda_1 > 0$ such that problem \eqref{shift:1} admits
unique solutions $(\bu_1, \fp_1, \rho_1) \in \CS_{p,q}((0, T))$ 
possessing the estimate:
\begin{equation}\label{ineq:g.1}
\CI_{p,q,T}(\bu_1, \rho_1; \eta)
\leq C\CJ_{p,q,T}(\bu_0, \rho_0, \bff, g, \bg, d, \bh; \eta)  
\end{equation}
for some constant $C$. 
\end{cor}

We consider  solutions $\bu$, $\fp$ and $\rho$ of problem
\eqref{*g.5.4} of the form:
$\bu = \bu_1 + \bv$, $\fp = \fp_1 + \fq$ and $\rho = \rho_1 + h$,
where $\bu_1$, $\fp_1$ and $\rho_1$ are solutions of the shifted 
equations \eqref{shift:1}, and then $\bv$, $\fq$ and $h$ should satisfy
the equations:
\begin{equation}\label{g.cor:1}\begin{cases*}
&\pd_t\bv  - \DV(\mu\bD(\bv) - \fq\bI) = -\lambda_1\bu_1,
\quad \dv \bv=0
\quad&\text{in $B_R^T$},\\
&\pd_t h- \bn\cdot P\bv
= -\lambda_1\rho_1
\quad&\text{on $S_R^T$},\\
&(\mu\bD(\bv) - \fq\bI)\bn  - \sigma 
(\CB h)\bn = 0 
\quad&\text{on $S_R^T$}, \\
&(\bv, h)|_{t=0} = (0, 0)
\quad&\text{on $B_R\times S_R$}.
\end{cases*}\end{equation}

Recall the definition of $J_q(\Omega)$ given in \eqref{sol:1} in Subsec.
\ref{subsec:3.4}, that is 
$$J_q(B_R) = \{\bff \in L_q(B_R)^N \mid (\bff, \nabla\varphi)_{B_R}
= 0\quad\text{for any $\varphi \in \hat H^1_{q',0}(B_R)$}\}.
$$
Recall that 
\begin{align*}
H^1_{q,0}(B_R) &= \{\varphi \in H^1_q(B_R) \mid \varphi|_{S_R} = 0\}, 
\\
\hat H^1_{q,0}(B_R) &=
\{\varphi \in L_{q, {\rm loc}}(B_R)
\mid \nabla\varphi \in L_q(B_R)^N, \enskip
\varphi|_{S_R} = 0\}.
\end{align*}
Since $C^\infty_0(B_R)$ is dense in $\hat H^1_{q',0}(B_R)$,
the necessary and sufficient condition in order that $\bu \in J_q(B_R)$ is 
that $\dv \bu = 0$ in $B_R$. Let $\psi \in H^1_{q,0}(B_R)$ be 
a solution of the variational equation:
\begin{equation}\label{g.sol:2}
(\nabla\psi, \nabla\varphi)_{B_R} = (\bu_1, \nabla\varphi)_{B_R}
\quad\text{for any $\varphi \in H^1_{q',0}(B_R)$},
\end{equation}
and let $\bw = \bu_1 - \nabla\psi$.  Then, $\bw \in J_q(B_R)$
and 
\begin{equation}\label{g.sol:3}
\|\bw\|_{L_q(B_R)} + \|\psi\|_{H^1_q(B_R)}
\leq C\|\bu_1\|_{L_q(B_R)}.
\end{equation}
Using $\bw$ and $\psi$ , we can rewrite the first equation in 
\eqref{g.cor:1} as follows:
$$\pd_t\bv-\DV(\mu\bD(\bv) - (\fq + \lambda_1\psi)\bI)
= -\lambda_1\bw, \quad \dv\bv=0
\quad\text{in $B_R^T$}.
$$
Thus, in what follows we may assume that 
\begin{equation}\label{g.sol:4}
\bu_1 \in H^1_p((0, T), J_q(B_R)) \cap 
L_p((0, T), H^2_q(B_R)^N).
\end{equation}

According to the argument in Subsec. \ref{subsec:3.3}, we introduce
a functional $P(\bv, h) \in H^1_q(B_R) + \hat H^1_{q,0}(B_R)$
that is a unique solution of the weak Dirichlet
problem 
\begin{equation}\label{g:wd:3}
(\nabla P(\bv, h), \nabla\varphi)_{B_R}
= (\DV(\mu\bD(\bv))-\nabla\dv\bv, \nabla\varphi)_{B_R}
\end{equation}
for any $\varphi \in \hat H^1_{q',0}(B_R)$, subject to 
\begin{equation}\label{g:wd:3*}
P(\bv, h) = \mu<\bD(\bv)\bn, \bn> - \sigma(\CB h) - \dv\bv
\quad\text{on $S_R$}.
\end{equation}
And then, to handle problem \eqref{g.cor:1} in the semigroup setting, 
we consider the initial value problem:
\begin{equation}\label{g.cor:2}\begin{cases*}
&\pd_t\bv  - \DV(\mu\bD(\bv) - P(\bv, h)\bI) = 0 
\quad&\text{in $B_R\times(0, \infty)$},\\
&\pd_t h- \bn\cdot P\bv
= 0
\quad&\text{on $S_R\times(0, \infty)$},\\
&(\mu\bD(\bv)-\fq\bI)\bn - \sigma 
(\CB h)\bn = 0 
\quad&\text{on $S_R\times(0, \infty)$}, \\
&(\bv, h)|_{t=0} = (\bv_0, \rho_0)
\quad&\text{on $B_R\times S_R$}.
\end{cases*}\end{equation}
Note that $(\mu\bD(\bv)-P(\bv, h)\bI)\bn - \sigma 
(\CB h)\bn = 0$ on $S_R\times(0, \infty)$ is equivalent to 
\begin{equation}\label{g.sol:5}
(\bD(\bv)\bn)_\tau = 0,
\quad \dv \bv = 0 \quad\text{ on  $S_R\times(0, \infty)$}.
\end{equation}
 
Defining 
$\CH_q(B_R)$, $\CD_q(B_R)$ and $\CA_q(\bv, h)$ by
\begin{equation}\label{space:0}\begin{split}
\CH_q(B_R) & = \{(\bv, h) \mid \bv \in J_q(B_R), \quad
h \in W^{2-1/q}_q(S_R)\}, \\
\CD_q(B_R) & = \{(\bv, h) \in \CH_q(B_R) \mid \bv \in H^2_q(B_R)^N,
\enskip h \in W^{3-1/q}_q(S_R), \\
&\phantom{ = \{(\bv, h) \in \CH_q(B_R) \mid }\,\,
(\bD(\bv)\bn)_\tau = 0 \enskip\text{on $S_R$}\}, \\
\CA_q(\bv, h) & = (\DV(\bD(\bv) - P(\bv, h)\bI), -\bn\cdot P\bv)
\quad \text{for $(\bv, h) \in \CD_q(B_R)$},
\end{split}\end{equation}
we see that Eq. \eqref{g.cor:2} is formulated by
\begin{equation}\label{g.sol:5*}
\pd_t U = \CA_q U \quad (t>0), \quad U|_{t=0} = U_0
\end{equation}
with $U = (\bv, h) \in \CD_q(B_R)$ for $t > 0$
and $U_0=(\bu_0, \rho_0) \in \CH_q(B_R)$. 
According to Theorem \ref{thm:semi:1.2}, we see that $\CA_q$ generates
a $C^0$ semigroup $\{T(t)\}_{t\geq 0}$ on $\CH_q(B_R)$. 
Moreover, if we define  
\begin{equation}\label{space:1}\begin{split}
\dot J_q(B_R) &= \{ \bff \in J_q(B_R) \mid 
(\bff, \bp_\ell)_{B_R} = 0 \enskip (\ell=1, \ldots, M)\}; \\
\dot W^\ell_q(S_R) & = \{ g \in W^\ell_q(S_R) \mid 
(g, \varphi_j)_{S_R} = 0 \enskip (j=1, \ldots, N+1)\}; \\
\dot\CH_q(B_R) & = \{(\bff, g) \mid 
\bff \in \dot J_q(B_R), \quad 
g \in \dot W^{2-1/q}_q(S_R)\}; \\
\|(\bff, g)\|_{\CH_q} &= \|\bff\|_{L_q(B_R)} 
+ \|g\|_{W^{2-1/q}_q(S_R)}; \\
\|(\bv, h)\|_{\CD_q} &= \|\bv\|_{H^2_q(B_R)} 
+ \|h\|_{W^{3-1/q}_q(S_R)},
\end{split}\end{equation}
then we have 
\begin{thm}\label{thm:g.3.3} Let $1 < q < \infty$.  Then,
 $\{T(t)\}_{t\geq0}$ is exponentially stable on $\dot\CH_q$,
that is 
\begin{equation}\label{decay:1}
\|T(t)(\bff, g)\|_{\CH_q}
\leq Ce^{-\eta_1 t}\|(\bff, g)\|_{\CH_q}
\end{equation}
for any $t >  0$ and $(\bff, g) \in \dot\CH_q(B_R)$
with some positive constants $C$ and $\eta_1$.
\end{thm}
Postponing the proof of Theorem \ref{thm:g.3.3} to the next section, 
we continue to prove Theorem \ref{thm:g.5.2}. 
Let 
$$\tilde\bu_1 = \bu_1 
- \sum_{\ell-1}^M(\bu_1(\cdot, t), \bp_\ell)_{B_R}\bp_\ell,\quad
\tilde\rho_1 = \rho_1 
- \sum_{j=1}^{N+1}(\rho_1(\cdot, t), \varphi_j)_{S_R}\varphi_j,
$$
and then $(\tilde\bu_1, \bp_\ell)_{B_R} = 0$
 ($\ell=1, \ldots, M$) and 
$(\tilde\rho_1, \varphi_j)_{S_R} = 0$ ($j=1, \ldots, N+1$). 
Moreover, since $\dv \bp_\ell=0$, by \eqref{g.sol:4} 
we have $\tilde\bu_1 \in L_p((0, T), \dot J_q(B_R))$. 
Let 
$$(\tilde\bv, \tilde  h)(\cdot, s) = \int^s_0T(s-r)(-\lambda_1
\tilde\bu_1(\cdot, r),
-\lambda_1\tilde\rho_1(\cdot, r))\,dr,$$
and then by the Duhamel principle, $\tilde\bv$ and $\tilde h$
satisfy the equations:
\begin{equation}\label{g.cor:4}\begin{cases*}
&\pd_t\tilde\bv  - \DV(\mu\bD(\tilde\bv) - 
P(\tilde\bv, \tilde h)\bI) 
= -\lambda_1\tilde\bu_1, \quad\dv\tilde\bv = 0
\quad&\text{in $B_R^T$},\\
&\pd_t \tilde h- \bn\cdot P\tilde\bv
= -\lambda_1\tilde\rho_1
\quad&\text{on $S_R^T$},\\
&(\mu\bD(\tilde \bv)-P(\tilde\bv, \tilde h)\bI)\bn 
-\sigma (\CB \tilde h)\bn = 0
\quad&\text{on $S_R^T$}, \\
&(\tilde\bv, \tilde h)|_{t=0} = (0, 0)
\quad\text{on $B_R\times S_R$}.
\end{cases*}\end{equation}
 By \eqref{decay:1}, 
\begin{align*}
&\|(\tilde\bv, \tilde h)(\cdot, s)
\|_{\CH_q}
\leq C\int^s_0e^{-\eta_1(s-r)}\|(\tilde\bu_1(\cdot, r),
\tilde\rho_1(\cdot, r))
\|_{\CH_q}\,dr \\
&\leq C\Bigl(\int^s_0e^{-\eta_1(s-r)}\,dr\Bigr)^{1/{p'}}
\Bigl(\int^s_0e^{-\eta_1(s-r)}
\|(\tilde\bu_1(\cdot, r),
\tilde\rho_1(\cdot, r))
\|_{\CH_q}^p\,dr
\Bigr)^{1/p}.
\end{align*}
Choosing $\eta > 0$ smaller if necessary, 
we may assume that $0 < \eta p < \eta_1$ without loss of generality.
Thus, by the inequality above we have 
\begin{align*}
&\int^t_0(e^{\eta s}\|(\tilde\bv, \tilde h)(\cdot, s)
\|_{\CH_q})^p\,ds \\
&\quad \leq C\int^t_0\Bigl(\int^s_0
e^{\eta sp}e^{-\eta_1(s-r)}\|(\bu_1(\cdot, r), 
\rho_1(\cdot, r))\|_{\CH_q}^p\,dr
\Bigr)\,ds \\
&\quad = 
\int^t_0\Bigl(\int^s_0
e^{-(\eta_1-p\eta)(s-r)}(e^{\eta r}\|(\bu_1(\cdot, r), 
\rho_1(\cdot, r))\|_{\CH_q})^p\,dr
\Bigr)\,ds\\
&\quad = 
\int^t_0(e^{\eta r}\|(\bu_1(\cdot, r), 
\rho_1(\cdot, r))\|_{\CH_q})^p
\Bigl(\int^t_r
e^{-(\eta_1-p\eta)(s-r)}\,ds
\Bigr)\,dr\\
&\quad
\leq (\eta_1-p\eta)^{-1}
\int^T_0(e^{\eta r}\|(\bu_1(\cdot, r), 
\rho_1(\cdot, r))\|_{\CH_q})^p\,dr,
\end{align*}
which, combined with \eqref{ineq:g.1}, leads to 
\begin{equation}\label{g.decay:2}\begin{split}
&\|e^{\eta s}(\tilde\bv, \tilde h)
\|_{L_p((0, t), \CH_q)} 
\leq C\CJ_{p,q,T}(\bu_0, \rho_0, \bff, g, \bg, d, \bh; \eta)
\end{split}\end{equation}
for any $t \in (0, T)$.
If $\tilde\bw$ and $\tilde k$ satisfy the shifted equations: 
$$\begin{cases*}
&\pd_t\tilde\bw + \lambda_1\tilde\bw - \DV(\mu\bD(\tilde\bw)
 - P(\tilde\bw, \tilde k)\bI) 
= \bff,  \quad\dv\tilde\bw= 0
\quad&\text{in $B_R^T$},\\
&\pd_t \tilde k + \lambda_1\tilde k- \bn\cdot P\tilde\bw
= d
\quad&\text{on $S_R^T$},\\
&(\mu\bD(\tilde \bw)-P(\tilde\bw, \tilde k)\bI)\bn 
-\sigma (\CB \tilde k)\bn = 0
\quad&\text{on $S_R^T$}, \\
&(\tilde\bw, \tilde k)|_{t=0} = (0, 0)
\quad\text{on $B_R\times S_R$},
\end{cases*}$$
where we have set  
$$\bff = -\lambda_1\tilde\bu_1 + \lambda_1\tilde\bv, \quad 
d =  -\lambda_1\tilde\rho_1 + \lambda_1\tilde h,
$$
by \eqref{ineq:g.1} and 
\eqref{g.decay:2}, we have 
$$
\CI_{p,q,T}(\tilde\bw, 
\tilde k; \eta) 
\leq C\CJ_{p,q,T}(\bu_0, \rho_0, \bff, f_d, 
\bff_d, g, \bh; \eta).
$$
But, noting that $\tilde \bv$ and $\tilde h$ satisfy Eq.\eqref{g.cor:4},
by the uniqueness, we see that
$\tilde \bw= \tilde \bv$ and $\tilde k = \tilde h$
for $t \in (0, T)$,
and so we have 
\begin{equation}\label{g.decay:3}
\CI_{p,q,T}(\tilde\bv,
\tilde h; \eta) 
\leq C\CJ_{p,q,T}(\bu_0, \rho_0, \bff, f_d, 
\bff_d, g, \bh; \eta).
\end{equation}
Let 
\begin{align*}
\bv &= \tilde\bv - \lambda_1\sum_{\ell=1}^M
\int^t_0(\bu_1(\cdot, s), \bp_\ell)_{B_R}\,ds\,\bp_\ell, \\
h  &= \tilde h - \lambda_1\sum_{j=1}^{N+1}
\int^t_0(\rho_1(\cdot, s), \varphi_j)_{S_R}\,ds\,\varphi_j.
\end{align*}
In this case, 
$$P(\bv, h) = P(\tilde\bv, \tilde h) 
+ \lambda_1(N-1)R^{-2}
\sigma\int^t_0(\rho_1(\cdot, s), \varphi_1)_{S_R}\,ds\,\varphi_1.$$
In fact, letting  
$$ \CC = \lambda_1(N-1)R^{-2}
\sigma\int^t_0(\rho_1(\cdot, s), \varphi_1)_{S_R}\,ds\,\varphi_1$$
for notational simplicity,  
we have 
$$
(\nabla(P(\bv, h) -(P(\tilde\bv, \tilde h)+\CC)), \nabla\psi)_{B_R}
= 0 
$$
for any $\psi \in \hat H^1_{q',0}(B_R)$, because $\nabla \varphi_1=0$. 
Moreover, on $S_R$ we have 
$$P(\bv, h) - (P(\tilde\bv, \tilde h) + \CC)
=0
$$ 
because 
$\bD(\bp_\ell) = 0$, $\dv \bp_\ell = 0$, and 
$$\CB h = \CB\tilde h -(N-1)R^{-2}\lambda_1\int^t_0(\rho_1(\cdot, s), \varphi_1)_{S_R}\,ds\,\varphi_1.$$  Thus, we have $P(\bv, h)
= P(\tilde\bv, \tilde h) + \CC.$

By \eqref{g.cor:4}, we have
\begin{alignat*}2
&\pd_t\bv  - \DV(\mu\bD(\bv) - P(\bv, h)\bI) 
= -\lambda_1\bu_1,  \quad\dv \bv = 0
\quad&\text{in $B_R^T$}, \\
&(\mu\bD(\bv)-P(\bv, h)\bI)\bn 
-\sigma (\CB h)\bn = 0
\quad&\text{on $S_R^T$}.
\end{alignat*}
Recall that $\bp_\ell\cdot\bn|_{S_R} = 0$ for $\ell=N+1, \ldots, M$.
Moreover, recalling that $\bp_\ell = |B_R|^{-1}\be_\ell$
($\ell=1, \ldots, N$), we have 
$$P\bp_\ell = |B_R|^{-1}
\Bigl(\be_\ell - |B_R|^{-1}\int_{B_R}\be_\ell\,dy\Bigr) =0,$$
and therefore, 
$$\pd_th-\bn\cdot P\bv = \pd_t\tilde h - \bn\cdot P\tilde\bv
-\lambda_1\sum_{j=1}^{N+1}(\rho_1(\cdot, t), \varphi_j)_{S_R}\,\varphi_j
= -\lambda_1\rho_1 \quad\text{on $S_R^T$}.
$$
Summing up, we have proved that $\bv$, 
$\fq = P(\bv, h)$, and $h$ satisfy the 
equations \eqref{g.cor:1}.

By \eqref{g.decay:3}, we have
\begin{equation}\label{g.decay:4}\begin{split}
&\|e^{\eta t}\pd_t(\bv, h)\|_{L_p((0, T), \CH_q)} 
\leq C\CJ_{p,q,T}(\bu_0, \rho_0, 
\bff, f_d, \bff_d, g, \bh; \eta)
\end{split}\end{equation}
for any $t \in (0, T)$.

To estimate 
$\|e^{\eta t}(\bv, h)\|_{L_p((0, T), \CD_q)}$, 
we use the following lemma.
\begin{lem}\label{lem:5.3.4} Let $1 < q < \infty$.  Let
$\bu \in H^2_q(B_R)^N\cap J_q(B_R)$
and $\rho \in W^{3-1/q}_q(S_R)$ satisfy the equations:
\begin{equation}\label{elliptic:1}\begin{cases*}
&-\DV(\mu\bD(\bu) - P(\bu, \rho)\bI)
 = \bff \quad&\text{in $B_R$}, \\
&\bn\cdot P\bu =g \quad&\text{on $S_R$}, \\
&(\mu\bD(\bu) - P(\bu, \rho)\bI)\bn - \sigma 
(\CB\rho)\bn = 0
\quad&\text{on $S_R$}.
\end{cases*}\end{equation}
Then, there exists a constant $C > 0$ such that 
\begin{equation}\label{elliptic:2}
\|(\bu, \rho)\|_{\CD_q}  \leq C\Bigl\{\|(\bff, g)\|_{\CH_q}
+ \sum_{\ell=1}^M|(\bu, \bp_\ell)_{B_R}|
+ \sum_{j=1}^{N+1}|(\rho, \varphi_j)_{S_R}|\Bigr\}.
\end{equation}
\end{lem}
Postponing the proof of Lemma \ref{lem:5.3.4} to 
the next section, we continue to prove Theorem \ref{thm:g.5.2}. 
By \eqref{g.cor:1}, $\bv$ and $h$ satisfy the elliptic equations:
\begin{align*}
&-\DV(\mu\bD(\bv) - P(\bv, h)\bI) 
= -\lambda_1\bu_1 -\pd_t\bv , \quad\dv \bv= 0
\quad&\text{in $B_R$}, \\
&\bn\cdot P\bv =\lambda_1\rho_1 +\pd_th  \quad&\text{on $S_R$}, \\
&(\mu\bD(\bv - P(\bv, h)\bI)\bn - \sigma 
(\CB h)\bn = 0
\quad&\text{on $S_R$},
\end{align*}
and therefore, applying Lemma \ref{lem:5.3.4} and using
\eqref{g.decay:4} yield that
\begin{equation}\label{g.decay:5}\begin{split}
&\|e^{\eta t}\bv\|_{L_p((0, T), H^2_q(B_R))}
+ \|e^{\eta t}h\|_{L_p((0, T), W^{3-1/q}_q(S_R))}\\
&\quad 
\leq C\Bigl\{\CJ_{p,q,T}(\bu_0, \rho_0, \bff, f_d, 
\bff_d, g, \bh; \eta) \\
&\qquad+ \sum_{\ell=1}^M\Bigl(\int^T_0(e^{\eta s}
|(\bv(\cdot, s), \bp_\ell)_{B_R}|)^p\,ds\Bigr)^{1/p}
\\
&\qquad\quad+ \sum_{j=1}^M\Bigl(\int^T_0(e^{\eta s}
|(h(\cdot, s), \varphi_j)_{S_R}|)^p\,ds\Bigr)^{1/p}\Bigr\}.
\end{split}\end{equation}
Let $\bu= \bu_1 + \bv$, $\fp = \fp_1 + \fq$ and $\rho = \rho_1 + h$.
By \eqref{shift:1} and \eqref{g.cor:1}, $\bu$, $\fp$ and 
$\rho$ satisfy the equations \eqref{*g.5.4}.  
Since
\begin{equation}\label{g.decay:6}\begin{split}
&\Bigl(\int^T_0e^{\eta s}|(\bv(\cdot, s), 
\bp_\ell)_{B_R}|)^p\,ds\Bigr)^{1/p} \\
&\quad \leq 
\Bigl(\int^T_0e^{\eta s}|(\bu(\cdot, s), 
\bp_\ell)_{B_R}|)^p\,ds\Bigr)^{1/p}
+ CJ_{\eta, T}, \\
&\Bigl(\int^T_0e^{\eta s}|(h(\cdot, s), 
\varphi_j)_{S_R}|)^p\,ds\Bigr)^{1/p} \\
&\quad \leq 
\Bigl(\int^T_0e^{\eta s}|(\rho(\cdot, s), 
\varphi_j)_{S_R}|)^p\,ds\Bigr)^{1/p}
+ CJ_{\eta, T},
\end{split}\end{equation}
 where we have
set  $J_{\eta, T} = \CJ_{p,q,T}(\bu_0, \rho_0, 
\bff, f_d, \bff_d, g, \bh; \eta)$,  
as follows from \eqref{ineq:g.1}, by \eqref{g.decay:4},
\eqref{g.decay:5}, and \eqref{g.decay:6}, we see that 
$\bu$, $\fp$ and $\rho$ satisfy the inequality \eqref{g.ineq:1}.

\subsection{Exponential stability
 of continuous analytic semigroup associated with
Eq. \eqref{g.cor:2}}\label{subsec:5.4}

In this subsection, we shall prove Theorem \ref{thm:g.3.3} stated
 in the previous subsection. For this purpose, we consider the equations:
\begin{equation}\label{eq:5.4.1}
(\lambda\bI - \CA_q)U = F
\end{equation}
for $F = (\bff, g) \in \dot J_q(B_R)$ 
and $U = (\bv, h) \in  \CD_q(B_R)\cap \dot J_q(B_R)$,
which is the resolvent problem corresponding to Eq. \eqref{g.sol:5*}. 
Here, $\bI$ is the 
identity operator, $\dot J_q(B_R)$ is the space defined in \eqref{space:1},
and $\CD_q(B_R)$ and $\CA_q$ are the domain and
 the operator defined in \eqref{space:0}. Since $\CR$ boundedness
implies the usual boundedness of operator families, by Theorem
\ref{thm:rbdd:1.2} and the observation in Subsec. \ref{subsec:3.3},
we have the following theorem.
\begin{thm}\label{thm:g.5.4}
Let $1 < q < \infty$ and $0 < \epsilon_0 < \pi/2$.  Then,
 there exists a $\lambda_0 > 0$ such that 
for any $\lambda \in \Sigma_{\epsilon_0, \lambda_0}$ and
$F \in \CH_q(B_R)$,  Eq. \eqref{eq:5.4.1} admits a unique
solution $U \in \CD_q(B_R)$ possessing the estimate:
\begin{equation}\label{eq:res.1}
|\lambda|\|U\|_{\CH_q} + \|U\|_{\CD_q} \leq 
C\|F\|_{\CH_q}
\end{equation}
for some constant $C > 0$.  Here, $\CH_q(B_R)$ is the space defined in 
\eqref{space:0}. 
\end{thm}
Our task in this subsection is to prove the 
following theorem.
\begin{thm}\label{thm:5.4} Let $1 < q < \infty$ and let
$\BC_+ = \{\lambda \in \BC \mid {\rm Re}\,\lambda \geq 0\}$. 
 Then, for any $\lambda \in \BC_+$ and $F \in \dot \CH_q(B_R)$,
Eq. \eqref{eq:5.4.1} admits a unique solution $U \in 
\CD_q(B_R) \cap \dot \CH_q(B_R)$ possessing the estimate \eqref{eq:res.1}.
\end{thm}
In the following, we shall prove Theorem \ref{thm:5.4}.  We start with
the following lemma. 
\begin{lem}\label{lem:5.4.1} Let $1 < q < \infty$,  let
$\lambda \in \BC\setminus(-\infty, 0)$ and $F =(\bff, g)
\in \dot \CH_q(B_R)$.
If  $U =(\bv, h)\in \CD_q(B_R)$ satisfies Eq. \eqref{eq:5.4.1}, then
$U=(\bv, h)$  
belongs to $\dot\CH_q(B_R)$. 
\end{lem}
\pf 
First we prove that $\bv \in \dot J_q(B_R)$.  By \eqref{eq:5.4.1}, 
we know that $\bv \in J_q(B_R) \cap H^2_q(B_R)^N$ satisfies the 
equations: 
\begin{alignat*}2
\lambda \bv - \DV(\mu\CD(\bv) - P(\bv, h)\bI) = \bff,
\quad \dv \bv &=0 
&\quad&\text{in $B_R$}, \\
(\mu\CD(\bv) - P(\bv, \fq)\bI)\bn -\sigma(\CB h)\bn & = 0
&\quad&\text{on $S_R$}.
\end{alignat*}
Since $F = (\bff, g) \in \dot \CH_q(B_R)$, we know that
$\bff \in J_q(B_R)$ and 
$(\bff, \bp_\ell)_{B_R} = 0$ for $\ell=1, \ldots, M$, and so 
by  the divergence theorem of Gau\ss\, we have 
\begin{align*}
&0 = (\bff, \bp_\ell)_{B_R} = (\lambda\bv - \DV(\mu(\bD(\bv)-
P(\bv, h)\bI), \bp_\ell)_{B_R} \\
&=\lambda(\bv, \bp_\ell)_{B_R} - \sigma (\CB h, 
\bn\cdot\bp_\ell)_{S_R}
+ \frac{\mu}{2}(\bD(\bv), \bD(\bp_\ell))_{B_R}-
(P(\bv, h), \dv \bp_\ell)_{B_R}.
\end{align*}
We see that 
\begin{equation}\label{intby:1}
(\CB h, \bn\cdot\bp_\ell)_{S_R}= 0 \quad(\ell=1, \ldots, M). 
\end{equation}
In fact, 
recalling that $\bp_\ell = |B_R|^{-1}\be_\ell$ ($\ell = 1, \ldots, N$)
and $\bn = y/|y| \in S_1$, 
we have 
\begin{align*}
(\CB h, \bn\cdot \bp_\ell)_{S_R} 
&= R^{-1}|B_R|^{-1}(h, \CB y_\ell)_{S_R} \\
&=R^{-3}|B_R|^{-1}(h, (N-1+\Delta_{S_1})y_\ell)_{S_R}
=0
\end{align*}
for $\ell=1, \ldots, N$.  Moreover, $\bp_\ell\cdot\bn 
=0$ for $\ell=N+1, \ldots, M$ because $\bp_\ell$ ($\ell=N+1
,\ldots, M$) is equal to  $c_{ij}(x_i\be_j - x_j\be_i)$
for some $i$ and $j$ and constant $c_{ij}$, and therefore 
$(Bh, \bn\cdot \bp_\ell)_{S_R}=0$ for $\ell=N+1, \ldots, M$. 

Since $\bD(\bp_\ell) = 0$ and $\dv \bp_\ell=0$, we have
$\lambda(\bv, \bp_\ell)_{B_R} = 0$, which, combined with
$\lambda\not=0$, leads to $(\bv, \bp_\ell)_{B_R}=0$, 
that is, $\bv \in \dot J_q(B_R)\cap \CD_q(B_R)$. 

Next, we prove that $h \in \dot W^{3-1/q}_q(S_R)$.  
We know that $h \in W^{3-1/q}_q(S_R)$ and $g$ satisfies the
equation:
$$\lambda h - \bn\cdot P \bv = g \quad\text{on $S_R$}.
$$
Since $(g, \varphi_j)_{S_R} = 0$, by the divergence theorem of 
Gau\ss\, we have   
\begin{align*}
0 &= (g, \varphi_j)_{S_R} = \lambda(h, \varphi_j)_{S_R}
-(P\bv\cdot\bn, \varphi_j)_{S_R} \\
&= \lambda(h, \varphi_j)_{S_R} - \int_{B_R}\dv((P\bv)\varphi_j)\,dx.
\end{align*}
Since  $\dv P\bv = \dv \bv = 0$ and since 
$\pd_\ell\varphi_j$ are constants, we have   
\begin{equation}\label{intby:2}\begin{split}
&\int_{B_R}\dv((P\bv)\varphi_j)\,dx \\
&= \int_{B_R}(\dv(P\bv))\varphi_j\,dx 
+ 
\sum_{\ell=1}^{N+1}\int_{B_R}\Bigl(v_\ell - 
|B_R|^{-1}\int_{B_R}v_\ell\,dy\Bigr)(\pd_\ell\varphi_j)\,dx
\\
&=0,
\end{split}\end{equation}
 Thus, 
$\lambda(h, \varphi_j)_{S_R} = 0$, which, combined with
$\lambda\not=0$, leads to $(h, \varphi_j)_{S_R} = 0$.
Therefore, we have $h \in \dot W^{3-1/q}_q(S_R)$. 
This completes the proof of Lemma \ref{lem:5.4.1}.
\qed \vskip0.5pc
Combining Theorem \ref{thm:g.5.4} and Lemma \ref{lem:5.4.1},
we have the following Corollary. 
\begin{cor}\label{cor:g.5.4} 
Let $1 < q < \infty$ and $0 < \epsilon_0 < \pi/2$. 
Then,
there exists a positive constant $\lambda_0$ such that
for any $\lambda \in \Sigma_{\epsilon_0, \lambda_0}$  
and $(\bff, g) \in \dot\CH_q(B_R)$,  Eq. \eqref{eq:5.4.1} 
 admits a unique solution $(\bv, h)
\in \CD_q\cap \dot\CH_q(B_R)$  
possessing the estimates \eqref{eq:res.1}.
\end{cor}

In view of Corollary \ref{cor:g.5.4}, in order to 
prove Theorem \ref{thm:5.4} it suffices to prove 
the following theorem. 
\begin{thm}\label{thm:5.4.2} Let $1 < q < \infty$ and 
let $\lambda_0$ be the same positive number as in
Corollary \ref{cor:g.5.4}.  Let 
$$Q_{\lambda_0} = 
\{\lambda \in \BC \mid {\rm Re}\,\lambda \geq 0, \enskip
|\lambda| \leq \lambda_0\}.$$
Then, for any $\lambda \in Q_{\lambda_0}$ and $(\bff, g) \in 
\dot\CH_q(B_R)$,  Eq. \eqref{eq:5.4.1} admits
a unique solution $(\bv, h) \in \CD_q(B_R) \cap \dot\CH_q(B_R)$ possessing
the estimate:
\begin{equation}\label{5.4.6}
\|(\bv, h)\|_{\CD_q}
\leq C\|(\bff, g)\|_{\CH_q}
\end{equation}
with some constant $C$ independent of $\lambda \in Q_{\lambda_0}$.
\end{thm}
\pf
We write  $\dot\CD_q = \CD_q(B_R) \cap \dot\CH_q(B_R)$ 
for the sake of simplicity. 
We first observe that
\begin{equation}\label{5.4.2.2}
\CA_q\dot\CD_q \subset \dot\CH_q(B_R).
\end{equation}
In fact, for $(\bv, h) \in \dot\CD_q$, we set $\CA_q(\bv, h) = (\bff, g)$,
that is $\DV(\mu(\bD(\bv) - P(\bv, h)\bI) = \bff$ in $B_R$ and 
$\bn\cdot P\bv = g$ on $S_R$. 
For any $\varphi \in \hat H^1_{q',0}(B_R)$, by \eqref{g:wd:3}
and the fact that $\dv \bv=0$, 
we have 
$$
(\bff, \nabla\varphi)_{B_R} 
= (\DV(\mu\bD(\bv) - P(\bv, h)\bI), \nabla\varphi)_{B_R} 
= (\nabla\dv\bv, \nabla\varphi)_{B_R}=0,
$$ 
which implies  $\bff \in J_q(B_R)$.

Next, we observe that
\begin{align*}
(\bff, \bp_\ell)_{B_R} &= (\DV(\mu\bD(\bv) - P(\bv, h)\bI),
\bp_\ell)_{B_R} \\
&= \sigma (\CB h, \bn\cdot\bp_\ell)_{S_R}
-\frac{\mu}{2}(\bD(\bv), \bD(\bp_\ell))_{B_R}
+ (P(\bv, h), \dv \bp_\ell)_{B_R}.
\end{align*}
Thus, by \eqref{intby:1} and the facts that $\bD(\bp_\ell)=0$
and $\dv\bp_\ell=0$, we have $(\bff, \bp_\ell)_{B_R}=0$,
and so,  $\bff \in \dot J_q(B_R)$.

Finally, by \eqref{intby:2} we have 
$$
(g, \varphi_j)_{S_R}  = (\bn\cdot P\bv, \varphi_j)_{S_R}
= \int_{B_R}\dv((P\bv)\varphi_j)\,dx =0, 
$$
which leads to $g \in \dot W^{3-1/q}_q(S_R)$. 
This completes the proof of  \eqref{5.4.2.2}.

In view of  Corollary \ref{cor:g.5.4},  $(\lambda_0\bI - \CA_q)^{-1}$
exists as a bounded linear operator from $\dot\CH_q(B_R)$ onto
 $\dot\CD_q$,
and then, the equation \eqref{eq:5.4.1} is rewritten as  
\begin{align*}
(\bff, g) &= (\lambda\bI - \CA_q)(\bv, h)
= (\lambda-\lambda_0)(\bv, h) + (\lambda_0\bI - \CA_q)(\bv, h)\\
&= (\bI + (\lambda-\lambda_0)(\lambda_0\bI - \CA_q)^{-1})
(\lambda_0\bI - \CA_q)(\bv, h).
\end{align*}
If $(\bI + (\lambda-\lambda_0)(\lambda_0\bI - \CA_q)^{-1})^{-1}$
exists as a bounded linear operator from $\dot\CH_q(B_R)$ into itself,
then we have 
\begin{equation}\label{solution:1}
(\bv, h) = (\lambda_0\bI-\CA_q)^{-1}
(\bI + (\lambda-\lambda_0)(\lambda_0\bI - \CA_q)^{-1})^{-1}(\bff, g).
\end{equation}
Thus, our task is to prove the existence of the
inverse operator 
$(\bI + (\lambda-\lambda_0)(\lambda_0\bI - \CA_q)^{-1})^{-1}$. 
Since $H^2_q(B_R)^N$ and $W^{3-1/q}_q(S_R)$ are compactly
embedded into $L_q(B_R)^N$ and $W^{2-1/q}_q(S_R)$,
respectively,  as follows
from the Rellich compact embedding theorem, 
$(\lambda_0\bI - \CA_q)^{-1}$ is a compact operator from  
$\dot\CH_q$ into itself. Thus, 
in view of  Riesz-Schauder theory, in order to 
prove the existence of the inverse operator 
$(\bI + (\lambda-\lambda_0)(\lambda_0\bI - \CA_q)^{-1})^{-1}$
it suffices to prove that the kernel of the map
$\bI + (\lambda-\lambda_0)(\lambda_0\bI - \CA_q)^{-1}$ is trivial. 
Thus, let $(\bff, g)$ be an element in $\dot\CH_q(B_R)$ such that
\begin{equation}\label{5.4.2.3}
(\bI + (\lambda-\lambda_0)(\lambda_0\bI-\CA_q)^{-1})(\bff, g) = (0, 0).
\end{equation}
Our task is to prove that $(\bff, g) = (0, 0)$. Since
$(\bff, g) = -(\lambda-\lambda_0)(\lambda_0\bI-\CA_q)^{-1}(\bff, g)
\in \dot\CD_q$, we have $(\lambda_0\bI- \CA_q)(\bff, g)
= -(\lambda-\lambda_0)(\bff, g)$, and so,  
$(\bff, g) \in \dot\CD_q$ satisfies the homogeneous equation:
\begin{equation}\label{5.4.2.4}
(\lambda\bI-\CA_q)(\bff, g) = (0, 0).
\end{equation}
Namely, 
$(\bff, g) \in \dot\CD_q$ satisfies the homogeneous equations:
\begin{equation}
\begin{cases*}
\lambda\bff - \DV(\mu\bD(\bff) - P(\bff, g)\bI) & = 0
\quad &\text{in $B_R$}, \\
\lambda g - \bn\cdot P\bff & = 0 
\quad&\text{on $S_R$}, \\
(\mu\bD(\bff) - P(\bff, g)\bI)\bn 
- \sigma (\CB g)\bn& = 0
\quad&\text{on $S_R$}.
\end{cases*} \label{5.4.7}
\end{equation}
First we consider the case where $2 \leq q < \infty$. 
Since $(\bff, g) \in \dot\CD_q \subset \dot\CD_2$,
by \eqref{5.4.7} and the divergence theorem of Gau\ss, we have 
\begin{align*}
0 &= (\lambda\bff - \DV(\bD(\bff) - P(\bff, g)\bI), \bff)_{B_R}
\nonumber \\
&=\lambda\|\bff\|_{L_2(B_R)}^2
- \sigma (\CB g, \bn\cdot\bff)_{S_R}
+ \frac{\mu}{2}\|\bD(\bff)\|_{L_2(B_R)}^2 
- (P(\bff, g), \dv\bff)_{B_R}.
\end{align*} 
For $h \in H^2_q(S_R)$ and $\bg = {}^\top(g_1, \ldots, g_N)$, 
we have 
\begin{equation}\label{intby:3}
(\CB h ,P\bg\cdot\bn)_{S_R} = (\CB h, \bg\cdot\bn)_{S_R},
\end{equation}
because recalling $\bn = y/|y| \in S_1$, we have 
\begin{align*}
&\sum_{j=1}^N|B_R|^{-1}\int_{B_R}g_\ell\,dy R^{-1}(\CB h, y_\ell)_{S_R}
\\
&\quad 
= \sum_{j=1}^N|B_R|^{-1}\int_{B_R}g_\ell\,dy\,R^{-1}(h, 
R^{-2}(N-1+\Delta_{S_1})y_\ell)_{S_R}=0.
\end{align*}
Moreover, $\dv \bff = 0$, because $\bff \in \dot J_q(B_R)$. Thus, 
noting that $\lambda g = P\bff\cdot\bn$ on $S_R$, we have 
\begin{equation} \lambda\|\bff\|_{L_2(B_R)}^2 -
\sigma \bar{\lambda}(\CB g, g)_{S_R} + \frac{\mu}{2}
\|\bD(\bff)\|_{L_2(B_R)}^2=0. \label{5.4.2.5}
\end{equation}
 To treat
$(\CB g, g)_{S_R}$, we use the following lemma. 
\begin{lem}\label{lem:5.4.4} Let
$$\dot H^2_2(S_R) = \{h \in H^2_2(S_R) \mid 
(h, 1)_{S_R} = 0, \quad (h, x_j)_{S_R} = 0
\enskip(j=1, \ldots, N)\}.$$
Then, 
\begin{equation}\label{5.4.8*}
-(\CB h, h) \geq c\|h\|_{L_2(S_R)}^2
\end{equation}
for any $h \in \dot H^2_2(S_R)$ with some constant
$c > 0$.
\end{lem}
Postponing the proof of Lemma \ref{lem:5.4.4}, 
we continue the proof of Theorem \ref{thm:5.4}. 
Since $g \in \dot W^{3-1/q}_q(S_R)
\subset \dot H^2_2(S_R)$, taking the real part of \eqref{5.4.2.5}
we have 
\begin{align*}
0 &= {\rm Re}\,\lambda(\|\bff\|_{L_2(B_R)}^2
-\sigma (\CB g, g)_{S_R})
+ \frac{\mu}{2}\|\bD(\bff)\|_{L_2(B_R)}^2 \\
& \geq {\rm Re}\,\lambda(\|\bff\|_{L_2(B_R)}^2
+c\sigma \|g\|_{L_2(S_R)}^2)
+ \frac{\mu}{2}\|\bD(\bff)\|_{L_2(B_R)}^2, 
\end{align*}
which, combined with ${\rm Re}\,\lambda \geq 0$, 
leads to $\bD(\bff) = 0$.   But, 
$(\bff, \bp_\ell)_{B_R} = 0$ for $\ell=1, \ldots, M$, 
and so, $\bff = 0$.  Thus, by the first equation in \eqref{5.4.7}, 
$\nabla P(\bff, g) = 0$, and so, 
$ P(\bff, g) = f_0$  with some constant $f_0$,
which, combined with the third equation in \eqref{5.4.7}, 
leads to  $\CB g = -\sigma^{-1} f_0$ on $S_R$. 
Since $(g, 1)_{S_R} = |S_R|(g, \varphi_1)_{S_R} = 0$, 
we have
$$-\sigma^{-1}|S_R|f_0 = (\CB g, 1)_{S_R}
= (g, \Delta_{S_R}1)_{S_R} + R^{-2}(N-1)(g, 1)_{S_R} = 0,$$ 
and so,  $f_0=0$, which implies that $\CB g=0$ on $S_R$.
Recalling that $g \in \dot W^{3-1/q}_q(S_R) \subset
\dot H^2_2(S_R)$, by Lemma \ref{lem:5.4.4} $g= 0$. 
Thus, we have $(\bff, g) = (0, 0)$, and so,  the formula
\eqref{solution:1} holds. Namely, problem \eqref{eq:5.4.1}
admits a unique solution$(\bv, h) \in \dot\CD_q$ possessing
the estimate: 
\begin{equation}\label{est:5.4.1}
\|(\bv, h)\|_{\CD_q} \leq C_\lambda\|(\bff, g)\|_{\CH_q(B_R)}
\end{equation}
with some constant $C_\lambda$ depending on $\lambda$,
when $2 \leq q < \infty$ and $\lambda \in Q_{\lambda_0}$.

Before considering the case where $1 < q < 2$, at this point 
we give a 
\vskip0.3pc\noindent
{\bf Proof of Lemma \ref{lem:5.4.4}}.~
Let $\{\lambda_j\}_{j=1}^\infty$ be the set of all eigen-values of 
the Laplace-Beltrami operator $\Delta_{S_R}$ on $S_R$. 
We may assume that $\lambda_1 > \lambda_2 > \lambda_3 > 
\cdots > \lambda_j >  \cdots \to -\infty$, and then $\lambda_1=0$
and $\lambda_2 = -(N-1)R^{-2}$. Let $E_j$ be the eigen-space
corresponding to $\lambda_j$, and then the dimension of $E_j$ 
is finite (cf. Neri \cite[Chapter III, Spherical Harmonics]{Neri}).  
Let $d_j = \dim E_j$, and then $d_1=1$ and 
$d_2 = N$.  Especially, $E_1 = \{a \mid a \in \BC\}$ and 
$E_2 = \{a_1x_1 + \cdots+ a_Nx_N 
\mid a_i \in \BC\,\, (i=1, \ldots, N)\}$.
Let $\{\varphi_{ij}\}_{j=1}^{d_i}$ be the orthogonal 
basis of $E_i$ in $L_2(S_R)$, and then for any $h \in \dot H^2_2(S_R)$
we have 
$$h = \sum_{i=3}^\infty\sum_{j=1}^{d_i} a_{ij}\varphi_{ij}
\quad(a_{ij} = (h, \varphi_{ij})_{S_R}),$$
because $(h, \varphi_{ij})_{S_R} = 0$ for $i = 1,2$.
Thus, we have 
$$
-(\CB h, h)_{S_R} = \sum_{i=3}^\infty\sum_{j=1}^{d_i}
|a_{ij}|^2(-\lambda_i-(N-1)R^{-1})\|\varphi_{ij}\|_{L_2(S_R)}^2.
$$
Since $-\lambda_i-(N-1)R^{-1} \geq c$ with some positive 
constant $c$ for any $i \geq 3$, we have \eqref{5.4.8*},
which completes the proof of Lemma \ref{lem:5.4.4}.
\qed \vskip0.5pc

Next, we consider the case where $1 < q < 2$. Let $(\bff, g) \in 
\dot\CD_q$ satisfy the homogeneous equations \eqref{5.4.7}.
First, we prove that 
\begin{equation}\label{null:1}
(\bff, \bg)_{B_R} = 0\quad 
\text{for any $\bg \in \dot J_{q'}(B_R)$}.
\end{equation} 
Let 
$(\bu, \rho) \in \dot\CD_{q'}$ be a solution of the 
equations:
\begin{equation}
\begin{cases*}
\bar{\lambda}\bu - \DV(\mu\bD(\bu) - 
P(\bu, \rho)\bI) & = \bg
\quad &\text{in $B_R$}, \\
\bar{\lambda}\rho - \bn\cdot P\bu & = 0 
\quad&\text{on $S_R$}, \\
(\mu\bD(\bu) - P(\bu, \rho)\bI)\bn
- \sigma (\CB \rho)\bn& = 0
\quad&\text{on $S_R$}.
\end{cases*} \label{5.4.8}
\end{equation}
Since $\lambda \in Q_{\lambda_0}$, $\bar{\lambda} \in Q_{\lambda_0}$,
and moreover $2 < q' < \infty$, and so,  by the fact proved above 
we know the unique existence of $(\bu, \rho) \in \dot\CD_{q'}$. 
By \eqref{5.4.7}, \eqref{5.4.8} and the divergence theorem
of Gau\ss, we have
\begin{align*}
&(\bff, \bg)_{B_R} 
= (\bff, \bar{\lambda}\bu -\DV(\mu\bD(\bu)
-P(\bu, \rho)\bI))_{B_R}\\
& =\lambda(\bff, \bu)_{B_R}
 -(\bff\cdot\bn, \sigma \CB\rho)_{S_R}
+ \frac{\mu}{2}(\bD(\bff), \bD(\bu))_{B_R} 
-(\dv\bff, P(\bu, \rho))_{B_R}.
\end{align*}
Noting that $P\bff\cdot\bn = \lambda g$ and $\dv \bff = 0$ and 
using \eqref{intby:3}, we have 
\begin{equation}\label{5.4.9}\begin{split}
(\bff, \bg)_{B_R} &= \lambda(\bff, \bu)_{B_R}
+\sigma \lambda\{(\nabla_{S_R}g, \nabla_{S_R}\rho)_{S_R}
\\
&-  R^{-2}(N-1)(g, \rho)_{S_R}\} + 
\frac{\mu}{2}(\bD(\bff), \bD(\bu))_{B_R}.
\end{split}\end{equation}
On the other hand, we have
\begin{align*}
0& = (\lambda\bff-\DV(\mu\bD(\bff)-P(\bff, g)\bI), \bu)_{B_R}\\
& = \lambda(\bff, \bu)_{B_R}
-\sigma (\CB g, \bn\cdot\bu)_{S_R}
+ \frac{\mu}{2}(\bD(\bff), \bD(\bu))_{B_R}
-(P(\bff, g), \dv \bu)_{B_R}.
\end{align*}
Noting that  $P\bu\cdot\bn = 
\bar{\lambda}\rho$ and $\dv \bu = 0$  and 
using \eqref{intby:3}, we have 
\begin{align*}
0 &= \lambda(\bff, \bu)_{B_R}
+\sigma \lambda\{(\nabla_{S_R}g, \nabla_{S_R}\rho)_{S_R}\\
&-  R^{-2}(N-1)(g, \rho)_{S_R}\} + 
\frac{\mu}{2}(\bD(\bff), \bD(\bu))_{B_R},
\end{align*}
which, combined with \eqref{5.4.9}, leads to  \eqref{null:1}.

Next, we prove that $(\bff, \bg)_{B_R}=0$ for any 
$\bg \in L_{q'}(B_R)^N$. Given $\bg 
\in L_{q'}(B_R)^N$, let $\psi \in \hat H^1_{q', 0}(B_R)$
be a solution to the variational equation: 
$$(\nabla \psi, \nabla\varphi)_{B_R} = (\bg, \nabla\varphi)_{B_R}
\quad\text{for any $\varphi \in \hat H^1_{q,0}(B_R)$}.
$$
Let $\bh = \bg - \nabla \psi$ and we decompose $\bg$
as 
$$\bg = \nabla\psi +\bh -\sum_{j=1}^M(\bh, \bp_\ell)_{B_R}
\bp_\ell + \sum_{j=1}^M(\bh, \bp_\ell)_{B_R}\bp_\ell.
$$
  Since 
$\bff \in \dot J_q(B_R)$, we have 
$(\bff, \bg)_{B_R} = (\bff, \bh - \sum_{j=1}^M
(\bh, \bp_\ell)\bp_\ell)_{B_R}$.  Since $\dv \bp_\ell = 0$, 
$\bh - \sum_{j=1}^M
(\bh, \bp_\ell)\bp_\ell \in \dot J_{q'}(B_R)$, 
and so,  by \eqref{null:1}
$(\bff, \bh - \sum_{j=1}^M(\bh, \bp_\ell)\bp_\ell
)_{B_R} = 0$, which implies that $(\bff, \bg)_{B_R} = 0$
for any $\bg \in L_{q'}(B_R)$.  Thus, we have 
$\bff=0$. By the first equation in \eqref{5.4.7}, 
$\nabla P(\bff, g)= 0$ in $B_R$, which leads to 
$P(\bff, g) = f_0$ with some constant $f_0$.
Thus, by the third equation in \eqref{5.4.7}, we have
$\CB g = -\sigma^{-1}f_0$ on $S_R$. 
Since $(g, 1)_{S_R} = 0$, we have 
$-\sigma^{-1}f_0|S_R| = (\CB g, 1)_{S_R} 
= R^{-2}(N-1)(g, 1)_{S_R} = 0$, 
which leads to $f_0=0$.  Thus, we have
$\CB g=0$ on $S_R$.  By the hypoellipticity of the
operator $\Delta_{S_R}$, we see that $g \in H^2_2(S_R)$,
and so, $g \in \dot H^2_2(S_R)$, which, combined with
Lemma \ref{lem:5.4.4},  leads to $g=0$.  Thus,
the formula \eqref{solution:1} holds, and therefore
problem \eqref{eq:5.4.1} admits a unique solution $(\bv, h)
\in \dot\CD_q$ possessing the estimate \eqref{est:5.4.1} 
when $1 < q < 2$ and $\lambda \in Q_{\lambda_0}$.

Finally, we prove that the constant in the estimate \eqref{est:5.4.1}
is independent of $\lambda \in Q_{\lambda_0}$. 
Let $\lambda \in Q_{\lambda_0}$ and 
$\mu \in \BC$, and we consider the equation:
\begin{equation}\label{5.4.2.10} (\mu\bI - \CA_q)(\bv, h) = (\bff, g).
\end{equation}
We write this equation as follows:
$$(\bff, g) = ((\mu-\lambda)\bI + (\lambda\bI - \CA_q))(\bv, h)
= (\bI + (\mu-\lambda)(\lambda\bI - \CA_q)^{-1})
(\lambda\bI - \CA_q)(\bv, h).
$$
Since $\|(\mu-\lambda)(\lambda\bI- \CA_q)^{-1}
\|_{\CL(\dot\CH_q(B_R))} \leq |\mu-\lambda|
C_\lambda$ as follows from \eqref{est:5.4.1}, 
where $\|\cdot\|_{\CL(\dot\CH_q(B_R))}$
denotes the operator norm of the bounded linear operator 
from $\dot\CH_q(B_R)$
into itself, choosing $\mu \in \BC$ in such a way that 
$|\mu-\lambda|C_\lambda \leq 1/2$, we see that 
the inverse operator $(\bI + (\mu-\lambda)(\lambda\bI - \CA_q)^{-1})^{-1}$
exists as a bounded linear operator from $\dot\CH_q(B_R)$ into itself and
$$\|(\bI + (\mu-\lambda)(\lambda\bI - \CA_q)^{-1})^{-1}\|_{\CL(\dot\CH_q(B_R))}
\leq 2.$$
Thus, $(\bv, h) = (\lambda\bI - \CA_q)^{-1}(\bI + 
(\mu-\lambda)(\lambda\bI - \CA_q)^{-1})^{-1}(\bff, g)$ belongs to
$\dot\CD_q$ and solves the equation \eqref{5.4.2.10}.  Moreover,
\begin{align*}
&\|(\bv, h)\|_{\CD_q}  \\
&
\leq \|(\lambda \bI - \CA_q)^{-1}\|_{\CL(\CH_q,
\CD_q)}\|(\bI + (\mu-\lambda)(\lambda\bI - \CA_q)^{-1}
\|_{\CL(\dot\CH_q(B_R))}
\|(\bff, g)\|_{\CH_q} \\
&\leq 2C_\lambda\|(\bff, g)\|_{\CH_q},
\end{align*}
provided that $|\mu-\lambda| \leq (2C_\lambda)^{-1}$, where
$\|\cdot\|_{\CL(\dot\CH_q, \CD_q)}$ denotes the operator norm of 
bounded linear operators from $\dot \CH_q(B_R)$ into $\CD_q(B_R)$. 
Since $Q_{\lambda_0}$ is a compact set, we have
\eqref{5.4.6}, which completes the proof of Theorem \ref{thm:5.4}.
\qed \vskip0.5pc

{\bf Proof of Lemma \ref{lem:5.3.4}}.~ Finally we prove Lemma \ref{lem:5.3.4}.
Let 
$$\tilde \bu = \bu - \sum_{\ell=1}^M(\bu, \bp_\ell)_{B_R}\,\bp_\ell,\quad
\tilde \rho = \rho - \sum_{j=1}^{N+1}(\rho, \varphi_j)_{S_R}\,\varphi_j.
$$
Since $\bD(\bp_\ell) = 0$, $\dv\bp_\ell=0$, and $\CB\varphi_j=0$ 
($j=2, \ldots, N+1$), we have
$$
(\nabla P(\tilde\bu, \tilde\rho), \nabla\psi)_{B_R} = (
\mu\DV \bD(\bu) - \nabla\dv\bu, \nabla\psi)_{B_R}
= (\nabla P(\bu, \rho), \nabla\psi)_{B_R}
$$
for any $\psi \in \hat H^1_{q',0}(B_R)$, subject to 
$$P(\tilde\bu, \tilde\rho) = <\mu\bD(\bu)\bn, \bn>-\sigma\CB\rho
+ \frac{\sigma(N-1)}{R^2}(\rho, \varphi_1)_{S_R}\varphi_1 - \dv\bu
\quad\text{on $\Gamma$},
$$
and so, $P(\tilde\bu, \tilde\varphi) = P(\bu, \rho) 
+ \frac{\sigma(N-1)}{R^2}(\rho, \varphi_1)_{S_R}\varphi_1$. 
Since $\nabla\varphi_1 = 0$,   
$\tilde\bu$ and $\tilde\rho$ satisfy the equations \eqref{elliptic:1}.
Moreover, $(\tilde \bu, \tilde\rho) \in \dot\CH_q(B_R) \cap \CD_q(B_R)$, and 
therefore by \eqref{5.4.2.2} and Theorem \ref{thm:5.4} with 
$\lambda=0$, we have
$\|(\tilde \bu, \tilde\rho)\|_{\CD_q(B_R)} \leq C\|(\bff, g)\|_{\CH_q(B_R)}$, 
which, combined with the estimate: 
$$\|(\bu, \rho)\|_{\CD_q(B_R)} \leq \|(\tilde\bu, \tilde\rho)\|_{\CD_q(B_R)}
+ C\Bigl\{\sum_{\ell=1}^M|(\bu, \bp_\ell)_{B_R}|
+ \sum_{j=1}^{N+1}|(\rho, \varphi_j)_{S_R}|\Bigr\},
$$
leads to  \eqref{elliptic:2}.
This completes the proof of Lemma \ref{lem:5.3.4}.

\subsection{Global Wellposedness, A proof of Theorem \ref{thm:sec:5}} 
\label{sec:5.6}
 In this section, we prove Theorem \ref{thm:sec:5}. 
Assume that the initial data $\bu_0 \in B^{2-2/p}_{q,p}(B_R)$ and 
$\rho_0 \in W^{3-1/p-1/q}_{q,p}(S_R)$ satisfy the smallness condition:
\begin{equation}\label{small:6.1}
\|\bu_0\|_{H^2_q(B_R)} + \|\rho_0\|_{W^{3-1/p-1/q}_{q,p}(S_R)}
\leq \epsilon
\end{equation}
with small constant $\epsilon >0$ as well as the compatibility condition
\eqref{compati:2}. For the notational simplicity, we write 
\begin{align*}
&\CI = \|\bu_0\|_{H^2_q(B_R)} + \|\rho_0\|_{W^{3-1/p-1/q}_{q,p}(S_R)}, \\
&\tilde E_{p,q,T}(\bu, \rho; \eta) 
= \|e^{\eta t}\bu\|_{L_p((0, T), H^2_q(B_R))} + \|e^{\eta t}
\pd_t\bu\|_{L_p((0, T), L_q(B_R))} \\
&\quad + \|e^{\eta t}\rho\|_{L_p((0, T), W^{3-1/q}_q(S_R))} 
 + \|e^{\eta t}\pd_t\rho\|_{L_p((0, T), W^{2-1/q}_q(S_R))} 
\\
&E_{p,q,T}(\bu, \rho; \eta) 
= \tilde E_{p,q,T}(\bu, \rho; \eta) 
+ \|e^{\eta t}\pd_t\rho\|_{L_\infty((0, T), W^{1-1/q}_q(S_R))}. 
\end{align*}
Notice that 
\begin{align*}
\|e^{\eta t}\bu\|_{L_\infty((0, T), B^{2(1-1/p)}_{q,p}(B_R))} 
&+ \|e^{\eta t}\rho\|_{L_\infty((0, T), B^{3-1/p-1/q}_{q,p}(S_R))}
\\
&\leq C(\CI + \tilde E_{p,q,T}(\bu, \rho; \eta)), 
\end{align*}
which follows from \eqref{real:7.1.1} and \eqref{real:7.1.2}. 
Since we choose $\epsilon$ small enough eventually, we may assume that
$0 < \CI \leq \epsilon < 1$. 
Let $T_0$ be a positive number $>2$.  In view of Theorem \ref{thm:g.5.1}, 
there exists a constant $\epsilon_1 > 0$
depending on $T_0$ such that if $\CI \leq \epsilon_1$, then 
for any $T \in (0, T_0]$ problem
\eqref{*g.5.4} admits  a unique solution $(\bu, \fq, \rho) 
\in \CS_{p,q}((0, T))$ satisfying the condition :
\begin{equation}\label{small:6.2}
\sup_{0 < t < T}\|\Psi_\rho(\cdot, t)\|_{H^1_\infty(B_R)} 
\leq \delta
\end{equation}
where $\delta \in (0,1/4)$ is the same constant as in \eqref{g.1.6}.   
We shall prove that $\bu$, $\fq$ and $\rho$
can be prolonged  beyond $T_0$ provided that $\epsilon > 0$ is small enough. 
Note that if we write
solutions of Eq. \eqref{g.1.10} by $U_T = (\bu_T, \fq_T, \rho_T) \in
\CS_{p,q}((0, T))$, then by the uniqueness of local in time solutions,
we see that $U_T = U_{T'}$ in $(0, T)$ for any $T$ and $T' \in (0, T_0]$ with
$T < T'$, and so in what follows we write $(\bu_T, \fq_T, \rho_T)$  simply by
$(\bu, \fq, \rho)$.
 Below, it is assumed that $0 < \epsilon \leq \epsilon_1$.

To prove Theorem \ref{thm:sec:5}, it suffices to prove that the inequality:
\begin{equation}\label{beyond:1}
E_{p,q,T}(\bu, \rho; \eta)
\leq M_3(\CI + E_{p,q,T}(\bu, \rho; \eta)^2 + E_{p,q,T}(\bu, \rho; \eta)^3)
\end{equation}
holds for any $T \in (0, T_0]$ with some constant $M_3 > 0$ 
independent of $\epsilon$, $T$, and $T_0$, where $\eta$ is 
the same positive constant as in Theorem \ref{thm:g.5.2}.

In fact, let $r_0(\epsilon)$ and  $r_\pm(\epsilon)$ be three different 
solutions of the algebraic equation: $x^3 + x^2 + \epsilon - M_3^{-1}x = 0$
with 
$r_0(\epsilon) = M_3\epsilon + O(\epsilon^2)$, 
$r_+(\epsilon) = M_3^{-1} + O(M_3^{-2}) + O(\epsilon)$, and 
$r_-(\epsilon) = -1 - M_3^{-1} + O(M_3^{-2}) + O(\epsilon)$ 
as $M_3 \to \infty$ and $\epsilon \to 0$.  Since $E_{p,q,T}(\bu, \rho; \eta)
\geq 0 > r_-(\epsilon)$, by \eqref{beyond:1} one of the 
following cases holds: 
$$E_{p,q,T}(\bu, \rho; \eta)\leq r_0(\epsilon), 
\quad 
E_{p,q,T}(\bu, \rho; \eta) \geq r_+(\epsilon).$$
Since we may change $M_3$ larger for \eqref{beyond:1} to hold if necessary,
and since we choose $\epsilon_1 > 0$ relatively small, we may assume that
\begin{equation}\label{small:6.2*}
e\epsilon_1 < M_3/2.
\end{equation}   
By \eqref{small:6.2*}, 
$E_{p,q,1/\eta}(\bu, \rho; \eta) \leq e\epsilon_1 <M_3/2 < r_+(\epsilon)$, 
and therefore, 
$$E_{p,q, T}(\bu, \rho; \eta) \leq r_0(\epsilon) \quad\text{
for $T \in (0, 1/\eta)$}.
$$  
But, 
$E_{p,q,T}(\bu, \rho; \eta)$ is a continuous function 
with respect to $T \in (0, T_0)$, which yields that
\begin{equation}\label{beyond:0} 
E_{p,q,T}(\bu, \rho; \eta) \leq r_0(\epsilon)
\quad\text{for any $T \in (0, T_0)$.}
\end{equation}

Let $T_0' = T_0-1/2$. 
Thus,  choosing $\epsilon > 0$
small enough and employing the same argument as that in proving
Theorem \ref{thm:g.5.1}, we see that there exist   unique 
solutions  
$(\bv, \fp, h) \in \CS_{p,q}((T'_0, T'_0+1))$  
of the equations:
\begin{equation}\label{np*}\begin{cases*}
&\partial_t \bv -\DV(\mu \bD(\bv) - \fp\bI)= \bff(\bv, \Psi_h)
\quad &\text{in $B_R\times(T'_0, T'_0+1)$},\\
&\dv \bv = g(\bv, \Psi_h) = \dv \bg(\bv, \Psi_h)
\quad &\text{in $B_R\times(T'_0, T'_0+1)$},\\
&\pd_th - \omega\cdot P\bv = 
\tilde d(\bv, \Psi_h) 
\quad &\text{on $S_R\times(T'_0, T'_0+1)$}, \\
&(\mu\bD(\bv)\omega)_\tau = \bh'(\bv, \Psi_h) 
\quad &\text{on $S_R\times(T'_0, T'_0+1)$}, \\
&<\mu\bD(\bv)\omega, \omega> - \fp - \sigma 
\CB\rho = h_N(\bv,\Psi_h)
\quad &\text{on $S_R\times(T'_0, T'_0+1)$}, \\
& (\bv, h)|_{t=T'_0} = (\bu(\cdot, T'_0), \rho(\cdot, T'_0))\quad 
&\text{in $B_R\times S_R$},
\end{cases*}\end{equation}
which satisfies the condition:
$$\sup_{T'_0 < t < T'_0+1}\|\Psi_h(\cdot, t)\|_{H^1_\infty(B_R)} \leq
\delta.
$$ 
Let
\begin{align*}
\bu_1 &= \begin{cases} \bu &\enskip (0 < t \leq T'_0), \\
\bv & \enskip(T'_0< t < T'_0+1),
\end{cases}
\quad 
\fq_1 = \begin{cases} \fq &\enskip (0 < t \leq T'_0), \\
\fp & \enskip(T'_0 < t < T'_0+1),
\end{cases}
\\
\rho_1 &= \begin{cases} \rho &\enskip (0 < t \leq T'_0), \\
h & \enskip(T'_0 < t < T'_0+1),
\end{cases}
\end{align*}
and then $(\bu_1,\fq_1, \rho_1)$ belongs to $\CS_{p,q}((0, T'_0+1))$ and 
 satisfies  
the condition:
$$\sup_{0 < t < T'_0+1}\|\Psi_{\rho_1}(\cdot, t)\|_{H^1_\infty(B_R)} \leq
\delta,
$$ 
and Eq. \eqref{*g.5.4} in  $(0, T'_0+1)$. Since 
$T'_0+ 1 =T_0+1/2$, we can prolong the solutions.  Repeating this 
argument,  we can prolong 
$\bu$, $\fq$ and $\rho$ to the time interval 
$(0, \infty)$.  

Below,  we prove \eqref{beyond:1}. Notice that $(\bu, \fq,\rho) \in 
\CS_{p,q}((0, T))$ satisfies Eq.\eqref{g.1.10}.
We shall estimate the right side of Eq. \eqref{g.1.10}. 
In view of \eqref{loc:7.1}, we have
\begin{equation}\label{eq:7.5.1}
\|e^{\eta t}\bff(\bu, \Psi_\rho)\|_{L_p((0, T), L_q(B_R))}
\leq C(\CI + E_{p,q,T}(\bu, \rho, \eta))E_{p,q,T}(\bu, \rho, \eta).
\end{equation}
Here and in the following, we use the inequality: 
$$E_{p,q,T}(\bu, \rho, 0) \leq E_{p,q,T}(\bu, \rho, \eta).$$
Employing the same argument as in the proof of \eqref{loc:7.3},
we have
\begin{equation}\label{eq:7.5.2}\begin{aligned}
\|e^{\eta t}\tilde d(\bu, \Psi_\rho)\|_{L_\infty((0, T), W^{1-1/q}_q(S_R))}
&\leq C (\CI+E_{p,q,T}(\bu, \rho, \eta)^2, \\
\|e^{\eta t}\tilde d(\bu, \Psi_\rho)\|_{L_p((0, T), W^{2-1/q}_q(S_R))}
&\leq C(\CI + E_{p,q,T}(\bu, \rho, \eta))^2E_{p,q,T}(\bu, \rho, \eta).
\end{aligned}\end{equation}
We now consider  $g(\bu, \Psi_\rho)$ and $\bg(\bu, \Psi_\rho)$.
Let $\tilde g_\eta = \tilde g(e^{\eta t}\bu, \Psi_\rho)$ 
and $\tilde \bg_\eta = \tilde \bg(e^{\eta t}\bu, \Psi_\rho)$,
where $\tilde g(\bv, \psi_h)$ and $\tilde \bg(\bv, \Psi_h)$ are the
 functions defined in \eqref{g.5.10}. 
Since $e^{\eta t}\bu|_{t=0} = \bu_0$, we note that 
in \eqref{eT:2*} 
$$\CE_1[e^{\eta t}\bu] = e_T[e^{\eta t}\bu - T_v(t)\bu_0]
+ \psi(t)T_v(|t|)\bu_0.
$$
Moreover, we may assume that $0 < \eta < 1$, and so the way  of estimating 
 $\tilde g_\eta$ 
and $\tilde \bg_\eta$ is the same as in 
Sect. \ref{sec:loc6}.  Moreover, by the same reason as in \eqref{g.5.11},
we have
\begin{gather}
\tilde g_\eta = e^{\eta t}g(\bu, \Psi_\rho),
\quad \tilde \bg_\eta = e^{\eta t}\bg(\bu, \Psi_\rho)
\quad\text{for $t \in (0, T)$}, \nonumber\\
\dv \tilde \bg_\eta = 
\tilde g_\eta \quad\text{for $t \in \BR$}.
\label{g.5.11.7}
\end{gather}
Using  \eqref{g.5.12*}, \eqref{g.5.7}, and \eqref{g.5.8}, 
we have
\begin{equation}\label{eq:7.5.3}
\|\tilde \bg_\eta \|_{L_p(\BR, L_q(B_R))}
\leq C(\CI + E_{p,q,T}(\bu, \rho, \eta))E_{p,q,T}(\bu, \rho, \eta).
\end{equation}
Employing the same argument as in proving \eqref{g.5.13**}
and using Lemma \ref{lem:g.5.1}, Lemma \ref{lem:g.5.2}, 
\eqref{g.5.7}, and \eqref{g.5.8}, we have 
\begin{equation}\label{eq:7.5.4}\begin{aligned}
&\|\tilde g_\eta \|_{H^{1/2}_p(\BR, L_q(B_R))}
+ \|\tilde g_\eta\|_{L_p(\BR, H^1_q(B_R))} \\
&\quad\leq C(\CI + E_{p,q,T}(\bu, \rho, \eta))E_{p,q,T}(\bu, \rho, \eta).
\end{aligned}\end{equation}

We now consider $\bh'(\bu, \Psi_\rho)$.  Since $\bh'(\bv, \Psi_h)$
is written like \eqref{g.5.13***}, 
following \eqref{g.5.13*} we define $\tilde \bh'_\eta$
by setting 
$$\bh'_\eta = \bv'_\bh(\bar\nabla\CE_2[\Psi_\rho])\bar\nabla
\CE_2[\Psi_\rho]\otimes\nabla \CE_1[e^{\eta t}\bu].
$$
And then, we have $e^{\eta t}\bh'(\bu, \Psi_\rho) = \tilde \bh_\eta'$
for $t \in (0, T)$.  Moreover, 
using Lemma \ref{lem:g.5.1}, Lemma \ref{lem:g.5.2}, 
\eqref{g.5.7}, and \eqref{g.5.8}, we have 
\begin{equation}\label{eq:7.5.5}\begin{aligned}
&\|\tilde \bh'_\eta \|_{H^{1/2}_p(\BR, L_q(B_R))}
+ \|\tilde \bh'_\eta\|_{L_p(\BR, H^1_q(B_R))} \\
&\quad\leq C(\CI + E_{p,q,T}(\bu, \rho, \eta))E_{p,q,T}(\bu, \rho, \eta).
\end{aligned}\end{equation}

We finally consider $h_N(\bu, \Psi_\rho)$, which is given in
\eqref{non:g.5}.  According to the formula \eqref{non:g.5},
we define $\tilde h_{N, \eta}$ by setting 
\begin{align*}
\tilde h_{N, \eta} &= \bV_{h,N}(\bar\nabla\CE_2[\Psi_\rho])
\bar\nabla\CE_2[\Psi_\rho]\otimes\nabla\CE_1[e^{\eta t}\bv] \\
&+ \sigma \tilde \bV'_\Gamma(\bar\nabla\CE_2[\Psi_\rho])
\bar\nabla\CE_2[\Psi_\rho]\otimes
\bar\nabla^2\CE_2[e^{\eta t}\Psi_{\rho}].
\end{align*}
We have $e^{\eta t}h_N(\bu, \Psi_\rho) = \tilde h_{N, \eta}$ 
for $t \in (0, T)$.
Since  $T_h(t)\rho_0|_{t=0} = \Psi_{\rho_0} =  e^{\eta t}\Psi_\rho|_{t=0}$, 
we note that in \eqref{eT:2*} 
$$\CE_2[e^{\eta t}\Psi_{\rho}]
=e_T[e^{\eta t} \Psi_\rho - T_h(t)\rho_0] + \psi(t)T_h(|t|)\rho_0.
$$
Thus, the way of estimating $\tilde h_N$ is the same as in  Sect. 6, 
and so, noting that
$$\|\bar\nabla^2\CE_2[e^{\eta t}\Psi_{\rho}]\|_{H^{1/2}_p(\BR, L_q(B_R))}
\leq 
\|\bar\nabla^2\CE_2[e^{\eta t}\Psi_{\rho}]\|_{H^1_p(\BR, L_q(B_R))},
$$
by Lemma \ref{lem:g.5.1}, Lemma \ref{lem:g.5.2}, 
\eqref{g.5.7}, and \eqref{g.5.8}, we have 
\begin{equation}\label{eq:7.5.6}\begin{aligned}
&\|\tilde h_{N, \eta} \|_{H^{1/2}_p(\BR, L_q(B_R))}
+ \|\tilde h_{N, \eta}\|_{L_p(\BR, H^1_q(B_R))} \\
&\quad\leq C(\CI + E_{p,q,T}(\bu, \rho, \eta))E_{p,q,T}(\bu, \rho, \eta).
\end{aligned}\end{equation}
Thus, applying Theorem \ref{thm:g.5.2*}
and using \eqref{eq:7.5.1}, \eqref{eq:7.5.2}, \eqref{eq:7.5.3},
\eqref{eq:7.5.4}, \eqref{eq:7.5.5}, and \eqref{eq:7.5.6}, we have 
\begin{equation}\label{5.6.3}\begin{split}
\tilde E_{p,q,T}(\bu, \rho; \eta) &\leq C\{\CI + 
E_{p,q,T}(\bu, \rho; \eta)^2 +
E_{p,q,T}(\bu, \rho; \eta)^3\\
&+ \sum_{\ell=1}^M
\Bigl(\int^T_0(e^{\eta s}|(\bu(\cdot, s), \bp_\ell)_{B_R}|)^p\,ds\Bigr)^{1/p}
\\
&+ \sum_{j=1}^{N+1}\Bigl(\int^T_0(e^{\eta s}
|(\rho(\cdot, s), \varphi_j)_{S_R}|)^p\,ds\Bigr)^{1/p}\}.
\end{split}\end{equation}
Here, we have used the inequalities:
$$(\CI + x)^2 \leq \CI + x^2, \quad (\CI + x)^2x \leq \CI + 3x^2 + x^3$$
for any $0 < \CI < 1$ and $x > 0$. 

By \eqref{loc:7.8.7}, we have 
\begin{align*}
&\|e^{\eta t}\pd_t\rho\|_{L_\infty((0, T), W^{1-1/q}_q(S_R))} \\
&\quad
\leq C(\|e^{\eta t}\bu\|_{L_\infty((0, T), H^1_q(B_R))}
+ \|e^{\eta t}\tilde d(\bv, \Psi_h)\|_{L_\infty((0, T), W^{1-1/q}_q(S_R))}),
\end{align*}
and so, by \eqref{eq:7.5.2} and \eqref{real:7.1.1}, we have
$$\|e^{\eta t}\pd_t\rho\|_{L_\infty((0, T), W^{1-1/q}_q(S_R))}
\leq C(\CI + \tilde E_{p,q,T}(\bu, \rho; \eta) + 
E_{p,q,T}(\bu, \rho; \eta)^2),$$
which, combined with \eqref{5.6.3}, leads to 
\begin{equation}\label{5.6.3.7}\begin{split}
E_{p,q,T}(\bu, \rho; \eta) &\leq C\{\CI + 
E_{p,q,T}(\bu, \rho; \eta)^2 +
E_{p,q,T}(\bu, \rho; \eta)^3\\
&+ \sum_{\ell=1}^M
\Bigl(\int^T_0(e^{\eta s}|(\bu(\cdot, s), \bp_\ell)_{B_R}|)^p\,ds\Bigr)^{1/p}
\\
&+ \sum_{j=1}^{N+1}\Bigl(\int^T_0(e^{\eta s}
|(\rho(\cdot, s), \varphi_j)_{S_R}|)^p\,ds\Bigr)^{1/p}\}.
\end{split}\end{equation}
Our final task is to prove that
\begin{align}
\sum_{\ell=1}^M\Bigl(\int^t_0
(e^{\eta s}|(\bu(\cdot, s), \bp_\ell)_{B_R}|)^p\,ds
\Bigr)^{1/p}
&+ \sum_{j=1}^{N+1}
\Bigl(\int^t_0
(e^{\eta s}|(\rho(\cdot, s), \varphi_j)_{S_R}|)^p\,ds
\Bigr)^{1/p} \nonumber \\
\leq C(\CI+E_{p,q,T}(\bu, \rho; \eta))E_{p,q,T}(\bu, \rho; \eta)
\label{5.6.4}
\end{align}
with some constant $C > 0$ independent of $\epsilon$ and 
$T > 0$.  If we have \eqref{5.6.4}, then putting  \eqref{5.6.3.7} and 
\eqref{5.6.4} together gives \eqref{beyond:1}. 


From now on, we prove \eqref{5.6.4}. Let $\bh_z^{-1}(x, t)$ be the inverse map
of the Hanzawa transform $x = \bh_z(y,t)$ for each $t \in 
[0, T]$, which is the diffeomorphism from $B_R$ onto $\Omega_t$ of 
$H^3_q$ class for each $t \in [0, T]$. Let $\bv = \bu\circ\bh_z^{-1}$ and $\fp
=(\fq + \sigma(N-1)R^{-1})\circ\bh_z^{-1}$, and then $\bv$ and $\fp$
satisfy the equations \eqref{navier:1}. Let $J$ be the Jacobian of the 
Hanzawa transform.  By \eqref{mass:1}, we have 
$|\Omega_t| = |\Omega|$, 
which, combined with \thetag{A.1}, leads to 
$|\Omega_t| = |B_R|$. 
Since $\Gamma_t$ is given by \eqref{g.1.2}, using the 
polar coordinates : $x-\xi(t) = s\omega$ for $\omega \in S_1$ and 
$s \in (0, R + \rho(R\omega, t))$, we have
\begin{align*}
&|B_R| = |\Omega_t| = \int_{\Omega_t}\,dx
= \int_{|\omega|=1}\,d\omega\int_0^{R+\rho(R\omega, t)}s^{N-1}\,ds\\
&= \frac1N\int_{|\omega|=1}(R+\rho(R\omega, t))^N\,d\omega 
 = |B_R| + \frac1N\sum_{j=1}^N {}_NC_j R^{1-j}
\int_{S_R}\rho(y, t)^j\,d\tau,
\end{align*}
which leads to 
\begin{equation}\label{orth:3}
(\rho,1)_{S_R} = -\sum_{j=2}^N{}_NC_j\,R^{1-j}
\int_{S_R}\rho(y, t)^j\,d\tau.
\end{equation}
Since $\|\rho(\cdot, t)\|_{L_\infty(S_R)} \leq \|H_\rho(\cdot, t)
\|_{L_\infty(B_R)} \leq \delta$ as follows from 
$H_\rho|_{S_R} = \rho$ on $S_R$, 
by \eqref{orth:3} and \eqref{real:7.1.2} we have
\begin{equation}\label{orth:4*}\begin{split}
&
\Bigl(\int^T_0(e^{\eta t}|(\rho(\cdot, t), \varphi_1)_{S_R}|)^p\,dt\Bigr)^{1/p}
\\
&\leq C\sup_{0 < t < T}\|\rho(\cdot, t)\|_{W^{2-1/q}_q(S_R)}
\|e^{\eta t}\rho\|_{L_p((0, T), L_q(S_R))} \\
&\leq C(\CI + E_{p,q,T}(\bu, \rho; \eta))E_{p,q,T}(\bu, \rho; \eta).
\end{split}\end{equation}
Since $\xi(t) = |\Omega|^{-1}\int_{\Omega_t}x\,dx$ and since $|\Omega_t|
= |\Omega|$, we have
$$0 = \int_{\Omega_t}x\,dx - |\Omega|\xi(t) = \int_{\Omega_t}
(x- \xi(t))\,dx.$$
Using the polar coordinates again, we have 
\begin{align*}
0  &= \int_{\Omega_t}(x_i - \xi_i(t))\,dx 
= \int_{|\omega| = 1}\,d\omega\int^{R + \rho(R\omega, t)}_0
(s\omega_i)s^{N-1}\,ds \\
& = \frac{1}{N+1}
\int_{|\omega|=1}\omega_i(R+\rho(R\omega, t))^{N+1}\,d\omega\\
& = \frac{R^{N+1}}{N+1}\int_{|\omega|=1}\omega_i\,d\omega
+R^N\int_{|\omega|=1}\rho(R\omega, t)\omega_i\,d\omega \\
&\quad + \frac{1}{N+1}\sum_{k=2}^{N+1}{}_{N+1}C_k\,R^{N+1-k}
\int_{|\omega|=1}\rho(R\omega, t)^k\omega_i\,d\omega \\
&= (\rho, x_i)_{S_R} + \frac{1}{N+1}\sum_{k=2}^{N+1}{}_{N+1}C_k
R^{1-k}(\rho^k, x_i)_{S_R}.
\end{align*}
Recalling that $\varphi_j$ equals constant $\times x_j$ 
($j=2, \ldots, N+1$),  we have
\begin{equation}\label{orth:4}\begin{split}
&
\Bigl(\int^T_0(e^{\eta t}|(\rho(\cdot, t), \varphi_i)_{S_R}|)^p\,dt\Bigr)^{1/p}
\\
&\leq C\sup_{0 < s < T}\|\rho(\cdot, t)\|_{W^{2-1/q}_q(S_R)}
\|e^{\eta t}\rho\|_{L_p((0, T), L_q(S_R))} \\
&\leq C(\CI+E_{p,q,T}(\bu, \rho; \eta))E_{p,q,T}(\bu, \rho; \eta).
\end{split}\end{equation}

By \eqref{mom:1} and \eqref{mom:2}, we have 
\begin{align}
\int_{\Omega_t}\bv(x, t)\,dx &
= \int_\Omega \bv_0(x)\,dx, \label{conserve:2} \\
\int_{\Omega_t}(x_iv_j(x, t) - x_jv_i(x, t))\,dx & 
= \int_\Omega (x_iv_{0j}(x) - x_jv_{0i}(x))\,dx.
\label{conserve:3}
\end{align}
Putting  \eqref{conserve:2}, 
\eqref{conserve:3} and \eqref{g.1.15} together gives 
\begin{equation}\label{conserve:4}
(\bv, \be_i)_{\Omega_t} = 0, \quad (\bv, x_i\be_j - x_j\be_i)_{\Omega_t}
= 0
\end{equation}
for $i, j=1, \ldots, N$. 
Since the Jacobian $J$ of the Hanzawa transform has the form:
$J = 1 + J_0(\bk)$, 
it then follows from \eqref{conserve:4} that 
\begin{align*}
0 &= \int_{B_R}\bu(y, t)\cdot\bp_\ell(y + H_\rho(y, t)\bn(y) + \xi(t))J\,dy  \\
&= (\bu, \bp_\ell)_{B_R} + 
\int_{B_R}\bu(y, t)\cdot\bp_\ell(y)J_0(\bk)\,dy  \\
&+ \int_{B_R}\bu(y, t)\tilde\bp_\ell(H_\rho(y, t)y + \xi(t))
(1 + J_0(\bk))\,dy
\end{align*}
where $\tilde\bp_\ell = 0$ for $\ell=1, \ldots, N$ and
$\tilde\bp_\ell$ are some matrices of first order polynomials for $\ell=N+1,
\ldots, M$. Since $\xi(0) = 0$ as follows from \thetag{A.2}
and since  $|J_0(\bk)|
\leq C|\nabla \Psi_\rho(\cdot,t)|$,  
using \eqref{g.1.3} and \eqref{g.1.6}, we have
\begin{align*}
|\xi(t)| &\leq \frac{1}{|B_R|}\int^t_0\int_{B_R}|\bu(y, t)|(1+|J_0(\bk)|)\,dy
\\
&\leq C\Big(\int^T_0e^{-\eta p's}\,ds\Bigr)^{1/{p'}}
\|e^{\eta t}\bu\|_{L_p((0, T), L_q(B_R))},
\end{align*}
and so, we have 
\begin{align*}
&\Bigl(\int^T_0(e^{\eta t}|(\bu(\cdot, t),\, \bp_\ell)_{B_R}|)^p\,dt
\Bigr)^{1/p}\\
&\leq \Bigl(\int^T_0\Bigl(e^{\eta t}
|\int_{B_R}\bu(y, t)\cdot\bp_\ell(y)J_0(\bk)\,dy|\Bigr)^p\,dt\Bigr)^{1/p}\\
&\phantom{aaaaaaaaa}+
\Bigl(\int^T_0\Bigl(e^{\eta t}
|\int_{B_R}\bu(y, t)\tilde\bp_\ell(\Psi_\rho(y, t) + \xi(t))\,dy|
\Bigr)^p\,dt\Bigr)^{1/p}\\
&\leq C\Bigl(\int^T_0(e^{\eta t}\|\bu(\cdot, t)\|_{L_q(B_R)})^p\,dt\Bigr)^{1/p}
(\|\Psi_\rho\|_{L_\infty((0, T), L_\infty(B_R))}
+ \sup_{t \in (0, T)}|\xi(t)|) \\
&\leq C(\CI+E_{p,q,T}(\bu, \rho; \eta))E_{p,q,T}(\bu, \rho; \eta).
\end{align*}
Putting this and \eqref{orth:4} together gives \eqref{5.6.4}.
This completes the proof of Theorem \ref{thm:sec:5}.

\section{Global well-posedness in an exterior domain}

In this section, we consider Eq. \eqref{navier:1} in the case where
$\sigma = 0$ and $\Omega$ is an exterior domain in $\BR^N$
whose boundary $\Gamma$ is a compact $C^2$ hypersurface.  
The problem we consider in this section is the following:
\begin{equation}\label{navier:2}\left\{\begin{aligned}
\pd_t\bv + (\bv\cdot\nabla)\bv - 
\DV(\mu\bD(\bv) - \fp\bI)  =0
&\quad&\text{in $\bigcup_{0 < t < T}\Omega_t\times\{t\}$}, \\
\dv\bv  = 0&\quad&\text{in $\bigcup_{0 < t < T}\Omega_t\times\{t\}$}, \\
(\mu\bD(\bv) - \fp\bI)\bn_t  = 0&
\quad&\text{on $\bigcup_{0 < t < T}\Gamma_t\times\{t\}$}, \\
V_n  = \bv\cdot\bn_t
&\quad&\text{on $\bigcup_{0 < t < T}\Gamma_t\times\{t\}$}, \\
\bv|_{t=0}  = \bv_0 \quad\text{in $\Omega_0$}, 
\quad \Omega_t|_{t=0} = \Omega_0=\Omega&
\end{aligned}\right.\end{equation}
In this section, we assume that $\mu$ is a positive constant. 
If we use the Hanazawa transform to transform $\Omega_t$ to some fixed 
domain, as was discussed in Sect. \ref{sec:loc6}, we have to require
$W^{3-1/q}_q$ regularity of the height function $\rho$  
representing $\Gamma_t$.
However, this regularity is obtained by surface tension, that is
the Laplace-Beltrami operator on $\Gamma_t$. We now consider the case
where the surface tension is not taken into account.  Thus, 
we can not obtain $W^{3-1/q}_q$ regularity of the height function.
Thus, we can not use the Hanzawa transform in the present case.
Another method is to use the Lagrange transform.  However, 
we can not expect the exponential decay unlike Sect. \ref{sec:6} for the
solutions of the Stokes equations with free boundary condition,
because $\Omega$ is unbounded domain. The decay of 
the solutions of the Stokes equations
with free boundary condition is only polynomial order, 
which is not sufficient to 
controle the term: first derivatives of $\int^t_0\bv(y, s)\,ds$
times the second derivatives of $\bv$. 

To overcome this difficulty, the idea here is to 
use the Lagrange transform only near the boundary.
Let $R$ be a
positive number for which $\CO = \BR^N\setminus\Omega 
\subset B_{R/2}$. Let $\kappa \in C^\infty_0(B_{2R})$ equal 
one in $B_R$. Let $\bu(y, t)$ be the velocity field in the Lagrange 
coordinates $\{\xi\}$.  We consider the partial Lagrange transform:
\begin{equation}\label{plag:1}
x = X_\bu(y, t) = y + \int^t_0\kappa(y)\bu(y, s)\,ds.
\end{equation} 
Assume that
\begin{equation}\label{assump:2}
\int^T_0\|\kappa(\cdot)\bu(\cdot, s)\|_{H^1_\infty(\Omega)}\,
ds \leq \delta,
\end{equation}
where $\delta > 0$ is a small number that will be chosen 
in such a way that several conditions hold. For example, 
if $\delta < 1/2$, then the map $x = X_\bu$ is injective
for each $t \in (0, T)$. Let 
$$\Psi(y, t) = \int^t_0\kappa(y)\bu(y, s)\,ds$$
and so $X_\bu(y, t) = y + \Psi(y, t)$. Let  
$$\Omega_t = \{x = X_\bu(y, t) \mid y \in \Omega \},
\quad \Gamma_t = \{x = X_\bu(y, t) \mid y \in \Gamma\}.$$
Let $y = X^{-1}_\bu(x, t)$ be the inverse 
of the transformation: $x = X_\bu(y, t)$ given in \eqref{plag:1}
and set 
$$\bv(x, t) = \bu(X^{-1}_\bu(x, t), t), \quad 
\fp(x, t) = \fq(X^{-1}_\bu(x, t), t).$$
We observe that 
\begin{align*}
V_n=\frac{\pd x}{\pd t}\cdot\bn_t = \bu(\xi, t)\cdot\bn_t
= \bv(x, t)\cdot\bn_t
\end{align*}
on $\Gamma_t$, because $\kappa=1$ on $\Gamma_t$, 
and so  the kinematic equation is autmatically satisfied. 
If $\bv$ and $\fp$ satisfy Eq.\eqref{navier:2}, then  
employing the same argument as in Subsec.\ref{sub:1}
and using the formula:
\begin{equation}\label{normal:8.1}
\bn_t = \frac{{}^\top A\bn}{|{}^\top A\bn|}
\quad\text{with}\enskip A = \Bigl(\frac{\pd x}{\pd y}\Bigr)^{-1}
= \bI + {}^\top\bV_0(\nabla\Psi(y,t)),
\end{equation}
which will be seen in the begining of Subsec \ref{subsec:6.3} 
below, 
and \eqref{change:3*}, 
we see that $\bu(y, t)$ and $\fq(y,t)$ satisfy the following 
equations:
\begin{equation}\label{navier:3}\left\{
\begin{aligned}
\rho\pd_t \bu
- \DV(\mu \bD(\bu) - \fq\bI)&= \bff(\bu), 
&\quad&\text{in $\Omega^T$}, \\
\dv\bu = g(\bu) &= \dv\bg(\bu)
&\quad&\text{in $\Omega^T$}, \\
(\mu\bD(\bu) - \fq\bI)\bn & = \bh(\bu)
&\quad &\text{on $\Gamma^T$}, \\
\bu|_{t=0} &= \bu_0 &\quad &\text{in $\Omega$}. 
\end{aligned}\right.\end{equation}
Here, $\bu_0 = \bv_0$, $\bn$ is the unit outer normal to $\Gamma$,
and the nonlinear terms in the right side of Eq. \eqref{navier:3}
are given as follows:
\begin{align}
\bff(\bu)|_i & = -\rho\sum_{j,k=1}^N(1-\kappa)u_j
(\delta_{jk} + V_{0jk}(\bk))\frac{\pd u_i}{\pd y_k} \nonumber\\
&\quad-\rho\sum_{\ell=1}^N \frac{\pd\Psi_\ell}{\pd y_i}
\Bigl(\frac{\pd u_\ell}{\pd t} 
+ \sum_{j,k=1}^N(1-\kappa)u_j
(\delta_{jk} + V_{0jk}(\bk))\frac{\pd u_\ell}{\pd y_k}\Bigr) \nonumber \\
&\quad + \mu\Bigl(\sum_{j=1}^N \frac{\pd}{\pd y_j}(\CD_\bD(\bk)\nabla\bu)_{ij}
\nonumber \\
&\quad + \sum_{j,k=1}^NV_{0jk}(\bk)\frac{\pd}{\pd y_k}(\bD(\bu)_{ij}
+ (\CD_\bD(\bk)\nabla\bu)_{ij})  \nonumber \\ 
&\quad + \sum_{j,k,\ell=1}^N\frac{\pd\Psi_\ell}{\pd y_i}
(\delta_{jk} + V_{0jk}(\bk))
\frac{\pd}{\pd y_k}(\bD(\bu)_{\ell j} + (\CD_\bD(\bk)\nabla\bu)_{\ell j})
\Bigr), \label{non-term:1}\\
g(\bu) &= -(J_0(\bk)\dv\bu + (1+J_0(\bk))\bV_0(\bk):\nabla\bu),
\nonumber \\
\bg(\bu) &= -(1+J_0(\bk)){}^\top\bV_0(\bk)\bu,
\label{non-term:2} \\
\bh(\bu) &= -\mu\{\bD(\bu){}^\top\bV_0(\bk)\bn
+ (\CD_\bD(\bk)\nabla\bu)(\bI + {}^\top\bV_0(\bk))\bn \nonumber \\
&\quad + {}^\top(\nabla\Psi)(\bD(\bu) + \CD_\bD(\bk)\nabla\bu)
(\bI + {}^\top\bV_0(\bk))\bn\}. \label{non-term:3}
\end{align}
where $\bk = \nabla \int^t_0(\kappa(y)\bu(y, s)\,ds$,  
$\Psi = \int^t_0\kappa(y)\bu(y, s)\,ds$, and 
 we have used the fact that $\pd_t\Psi_j = \kappa u_j$.

The main result of this section is the following theorem
that showes the unique existence theorem of global in
time solutions of  Eq.\eqref{navier:3} and asymptotics
as $t\to\infty$. 
\begin{thm} \label{6.thm:main} 
Let $N \geq 3$ and let $q_1$ and $q_2$ be exponents such that 
$$max(N, \frac{2N}{N-2}) < q_2 < \infty, \quad 
1/q_1 = 1/q_2 + 1/N. $$
Let $b$ and $p$ be numbers defined by
\begin{equation}\label{6.cond:0}\begin{split}
b = \frac{3N}{2q_2} + \frac12, \enskip 
p = \frac{2q_2(1+\sigma)}{q_2-N}
\end{split}\end{equation}
with some very small positive number $\sigma$. 
Then, there exists an $\epsilon > 0$ such that 
if initial data $\bu_0 \in B^{2(1-1/p)}_{q_2,p}(\Omega)^N
\cap B^{2(1-1/p)}_{q_{1/2}, p}(\Omega)^N$ 
satisfies the compatibility condition:
$$
\dv \bu_0= 0 \quad\text{in $\Omega$}, 
\quad \bD(\bu_0)\bn - <\bD(\bu_0)\bn, \bn>\bn = 0 
\quad\text{on $\Gamma$}, 
$$
and the smallness condition:
\begin{equation}\label{6.small.1}\|\bu_0\|_{B^{2(1-1/p)}_{q_2,p}} 
+ \|\bu_0\|_{B^{2(1-1/p)}_{q_1/2,p}} \leq \epsilon,
\end{equation}
then Eq. \eqref{navier:3} admits  unique solutions $\bu$  
and $\fq$ with
\begin{equation}\label{6.reg.1}
\begin{split}
\bu &\in L_p((0, \infty), H^2_{q_2}(\Omega)^N) 
\cap H^1_p((0, \infty), L_{q_2}(\Omega)^N),
\\
\fq &\in L_p((0, \infty), H^1_{q_2}(\Omega)
+  \hat H^1_{q_2, 0}(\Omega)),
\end{split}\end{equation}
possessing the estimate $[\bu]_\infty \leq C\epsilon$ with
\begin{align}
&[\bu]_T  = 
\Bigl\{
\int^T_0(<s>^b\|\bu(\cdot, s)\|_{H^1_\infty(\Omega)})^p\,ds 
\nonumber \\ 
&+\int^T_0 (<s>^{(b-\frac{N}{2q_1})}
\|\bu(\cdot, s)\|_{H^1_{q_1}(\Omega)})^p\,ds 
+ (\sup_{0 < s < T}<s>^{\frac{N}{2q_1}}\|\bu(\cdot, s)
\|_{L_{q_1}(\Omega)})^p
\nonumber  \\
&+ \int^T_0(<s>^{(b-\frac{N}{2q_2})}(\|\bu(\cdot, s)\|_{H^2_{q_2}(\Omega)}
+ \|\pd_t\bu(\cdot, s)\|_{L_{q_2}(\Omega)}))^p\,ds\Bigr\}^{1/p}. 
\label{6.norm.1}\end{align}
Here, $<s> = (1+s^2)^{\frac12}$ and 
$C$ is a constant that is independent of $\epsilon$.
\end{thm}
\begin{rem} Let $p' = p/(p-1)$, that is 
$1/{p'} = 1-1/p$. And then, 
$$\frac{1}{p'} = \frac{(1+2\sigma)q_2 + N}{2q_2(1+\sigma)}.
$$
We choose $\sigma > 0$ small enough in such a way that 
the following relations hold:, 
\begin{align}
&1 < q_1 < 2, \enskip \frac{N}{q_1} > b > \frac{1}{p'}, 
\enskip \Bigl(\frac{N}{q_1}-b\Bigr)p > 1,
\enskip
\Bigl(b-\frac{N}{2q_2}\Bigr)p > 1, 
\enskip b \geq \frac{N}{2q_1}, \nonumber \\
&b \geq \frac{N}{q_2},
\enskip \Bigl(\frac{N}{2q_2} + \frac12
\Bigr)p' < 1, \enskip bp' > 1, \enskip 
\Bigl(b - \frac{N}{2q_2}\Bigr)p' > 1, 
\enskip \frac{N}{q_2} + \frac{2}{p} < 1.
\label{6.cond:1}\end{align}
\end{rem}
\begin{rem} The exponent $q_2$ is used to control 
the nonlinear terms, and so $q_2$ is chosen in
such a way that $N < q_2 < \infty$.  Let 
\begin{equation}\label{6.cond:1*}
\quad \frac{1}{q_1} = \frac{1}{N} + \frac{1}{q_2},
\quad \frac{1}{q_3} = \frac{1}{q_1} + \frac{1}{q_2}.
\end{equation}
The exponent $q_2$ is used to control the nonlinear term,
and so $N < q_2 < \infty$ is required. And  
the condition: 
$q_2 > \frac{2N}{N-2}$ 
implies that  $q_1 > 2$ and $q_3 > 1$, which is necessary to prove 
Theorem \ref{6.thm:main}.  Thus,  we assume that 
$$\max(N, \frac{2N}{N-2}) < q_2 < \infty.$$
\end{rem}
\begin{rem} If we choose $\delta > 0$ in \eqref{assump:2}, 
then $x = X_\bu(y, t)$ becomes
a diffeomorphism with suitable regularity from $\Omega$ onto $\Omega_t$,
and so the original problem \eqref{navier:2} is globally well-posed. 
\end{rem}

\subsection{Maximal $L_p$-$L_q$ regularity in an exterior domain}
\label{subsec:6.1}
In this subsection, we study the maximal $L_p$-$L_q$ regularity of 
solutions to the Stokes equations with free boundary condition:
\begin{equation}\label{6.max.1}
\left\{\begin{aligned}
\pd_t\bu - \DV(\mu\bD(\bu) - \fq\bI) = \bff, \quad
\dv \bu = g &= \dv\bg &\quad&\text{in $\Omega^T$}, \\
(\mu\bD(\bu)- \fq\bI)\bn & = \bh
&\quad&\text{on $\Gamma^T$}, \\
\bu|_{t=0} & = \bu_0 &\quad&\text{in $\Omega$}.
\end{aligned}\right.\end{equation}
Let 
\begin{align*}
\hat H^1_{q,0}(\Omega) & = \{\varphi \in L_{q, {\rm loc}}(\Omega) \mid
\nabla\varphi \in L_q(\Omega)^N, \quad \varphi|_\Gamma=0\}, \\
J_q(\Omega) & = \{\bu \in L_q(\Omega)^N \mid 
(\bu, \nabla\varphi)_\Omega=0 
\quad\text{for any $\varphi \in \hat H^1_{q', 0}(\Omega)$}\}.
\end{align*}
We start with the following proposition.
\begin{prop}\label{prop:6.1} Let $1 < q < \infty$.  If 
$\bu \in H^1_q(\Omega)$ satisfies $\dv \bu = 0$ in $\Omega$,
then $\bu \in J_q(\Omega)$.
\end{prop}
To prove Proposition \ref{prop:6.1}, we need the following lemma.
\begin{lem}\label{lem:bog} Let $1 < q < \infty$, $m \in \BN_0$  and let
$G$ be a bounded domain whose boundary $\pd G$ is a hypersurface of
$C^{m+1}$ class. Let $H^0_{q,0}(G) = L_q(G)$ and for $m \geq 1$ let 
$$H^m_{q, 0}(G) = \{f \in H^m(G) \mid \pd_x^\alpha f|_{\pd G} = 0
\quad\text{for $|\alpha| \leq m-1$}\}.
$$
Then, there exists a linear operator 
$\BB:  H^m_{q,0}(G) \to  H^{m+1}_{q,0}(G)^N$
having the following properties:
\begin{itemize}
\item[\thetag1]
There exists a $\rho \in C^\infty_0(G)$ such that
$\rho \geq 0$, $\displaystyle{ \int_G \rho \,dx = 1}$, and \\
$\displaystyle{\dv \BB[f] = f-\rho\int_\Omega f\,dx}$. 
In particular, if $\displaystyle{\int_G f\,dx = 0}$, then 
$\dv\BB[f] = f$. 
\item[\thetag2] We have the estimate:
$$\|\BB[f]\|_{H^{k+1}_q(G)} \leq C_{q,k,G}\|f\|_{H^k_q(G)}
\quad(k=0, \ldots, m).$$
\item[\thetag3]~ If $f = \pd g/\pd x_i$ with some 
$g \in H^{m+1}_{q,0}(G)$, then 
$$\|\BB[f]\|_{H^{k}_q(G)} \leq C_{q,k,G}\|g\|_{H^k_q(G)}
\quad(k=0, \ldots, m).$$
\end{itemize}
\end{lem}
\begin{rem}
\thetag1 Since $\BB[f] \in H^{m+1}_{q,0}(G)$, if the $0$
extension of $\BB[f]$ to $\BR^N$ is written simply by
$\BB[f]$, then $\BB[f] \in H^{m+1}_q(\BR^N)^N$
and ${\rm supp}\, \BB[f] \subset G$.
\\
\thetag2 Lemma \ref{lem:bog} was proved by Bogovski
\cite{Bog1, Bog2} (cf. also Galdi \cite{Gal})
\end{rem} 
To apply Lemma \ref{lem:bog}, we use the following lemma.
\begin{lem}\label{lem:8.4} Let $1 < q < \infty$ and let 
$2R < L_1 < L_2 < L_3 < L_4 < 5R$. 
Let $\chi$ be a function in $C^\infty(\BR^N)$ such that 
$\chi(x) = 1$ for $x \in B_{L_2}$ and 
$\chi(x) = 0$ for $x \not\in B_{L_3}$.  
If $\bv \in H^2_q(G)^N$, $G \in \{\BR^N, \Omega, \Omega_{5R}\}$,
 satisfies 
$\dv \bv=0$ in $D_{L_1, L_4}$, then $(\nabla\chi)\cdot\bv \in
H^3_{q,0,a}(D_{L_2, L_3})$. Here, we have set 
$$H^3_{q,0,a}(D_{L_2, L_3})
= \{f \in \dot H^3_{q,0}(D_{L_2, L_3}) \mid \int_{D_{L_2, L_3}} f(x)\,dx = 0\}.
$$
Here and in the following, we set $D_{L,M}=
\{x \in \BR^N \mid L < |x| < M\}$ for $0 < L < M$.
\end{lem}
To prove Lemma \ref{lem:8.4}, we need the following lemma.
\begin{lem}\label{lem:8.3}
Let $1 < q < \infty$ and let $2R < L_1 < L_2 < L_3 < L_4 < L_5 < L_6
< 5R$.  Let $\chi$ be a function in 
$C^\infty(\BR^N)$ such that ${\rm supp}\,
\nabla\chi \subset D_{L_3, L_4}$.  If $\bu 
\in H^2_q(G)$, $G \in \{\BR^N, \Omega, \Omega_{5R}\}$,
 satisfies $\dv\bu = 0$ in 
$D_{L_1, L_6}$, then there exists a $\bv 
\in H^2_q(\BR^N)^N$ such that ${\rm supp}\, \bv 
\subset D_{L_2, L_5}$, $\dv \bv=0$ in 
$\BR^N$ and $(\nabla\chi)\cdot\bv = (\nabla\chi)\cdot\bu$
in $\BR^N$. 
\end{lem}
\pf
Let $A_0$, $A_1$, $\ldots$, $A_5$ and $B_0$, $B_1$, 
$\ldots$, $B_5$ be numbers such that
$L_2 = A_5 < A_4 < A_3 < A_2 < A_1 <A_0< L_3 < L_4 < 
B_0 < B_1 < B_2 < B_3 < B_4 < B_5 = L_5$. 
Let $\varphi \in C^\infty_0(\BR^N)$ such that 
$\varphi(x) = 1$ for $A_2 < |x| < B_2$ and $\varphi(x) = 0$
for $|x| < A_3$ or $|x| > B_3$. Note that $\varphi(x) = 1$ on
${\rm supp}\,\nabla\chi$.  Let $E = D_{A_4, A_1} \cup D_{B_1, B_4}$.
Since $\dv \bu = 0$ in $D_{L_1, L_6}$ and ${\rm supp}\, \varphi
\subset D_{A_3, B_3}$, we have $\dv(\varphi\bu) = 
(\nabla\varphi)\cdot\bu$ in $D_{A_4, B_4}$ and 
$\dv(\varphi\bu) = (\nabla\varphi)\cdot \bu \in H^2_{q,0}(E)$. 
Moreover, we have 
\begin{align*}
\int_E (\nabla\varphi)\cdot\bu\,dx& = \int_{D_{A_4, B_4}}(\nabla\varphi)\cdot
\bu\,dx = \int_{D_{A_4, B_4}}\dv(\varphi\bu)\,dx \\
&= -\int_{S_{A_4}}\frac{x}{|x|}\cdot(\varphi\bu)\,d\tau
+ \int_{S_{B_4}}\frac{x}{|x|}\cdot(\varphi\bu)\,d\tau= 0,
\end{align*}
which leads to  $(\nabla\varphi)\cdot\bu \in H^2_{q,0,a}(E)$. 
By Lemm \ref{lem:bog}, $\bv = \varphi\bu - \BB[(\nabla\varphi)\cdot\bu]$
has the properties: $\bv \in H^2_q(\BR^N)$, 
${\rm supp}\,\bv \subset D_{A_3, B_3}$, and 
$\dv \bv = 0$ in $\BR^N$.  Moreover, $(\nabla\chi)\cdot\bu
= (\nabla\chi)\cdot\bv$ in $\BR^N$, because 
$\varphi = 1$ on ${\rm supp}\,\nabla\chi$ and $\BB[(\nabla\varphi)\cdot
\bu] = 0$ on ${\rm supp}\, \nabla\chi$.
\qed \vskip0.5pc
{\bf Proof of Lemma \ref{lem:8.4}.}~ 
By Lemma \ref{lem:8.3}, there exists a 
$\bw \in H^2_q(\BR^N)^N$ such that
$(\nabla\chi)\cdot\bv = (\nabla\chi)\cdot\bw$ in $\BR^N$, ${\rm supp}\,\bw 
\subset B_{L_4}$ and 
$\dv \bw=0$ in $\BR^N$. Then, the  assertion follows from
the following observation:
\begin{align*}
\int_{D_{L_2, L_3}}(\nabla\chi)\cdot\bv\,dx &= \int_{D_{L_2, L_3}}
(\nabla\chi)\cdot\bw\,dx = \int_{B_{L_4}}(\nabla\chi)\cdot\bw\,dx\\
&= \int_{B_{L_4}}\dv(\chi\bw)\,dx = 
\int_{S_{L_4}}\frac{x}{|x|}\cdot(\chi\bw)\,d\tau=0,
\end{align*}
because $\chi|_{S_{L_4}} = 0$. This completes the proof of 
Lemma \ref{lem:8.4}. \qed
\vskip0.5pc
{\bf Proof of Proposition \ref{prop:6.1}}~ We use the Sobolev cut off
function.  Let $\chi(t) \in C^\infty_0((-1, 1))$ equal one for $|t| \leq 1/2$,
and set 
$$\chi_L(x) = \chi\Bigl(\frac{\ln\ln|x|}{\ln\ln L}\Bigr).
$$
Notice that
\begin{equation}\label{6.sob.1}
|\nabla \chi_L(x)| \leq \frac{c}{\ln\ln L}\frac{1}{|x|\ln|x|}
\end{equation}
and $\nabla\chi_L(x)$ vanishes outside of $\dot\Omega_L$, where 
$\dot\Omega_L= \{x \in \Omega \mid \exp\sqrt{\ln L} < |x| < L\}$. 
Since $\bu \in H^1_q(\Omega)$ and $\dv \bu = 0$ in $\Omega$, 
for any $\varphi \in H^1_{q', 0}(\Omega)$ we have
$$(\bu, \nabla\varphi)_\Omega = \lim_{L\to\infty}
(\chi_L\bu, \nabla\varphi)_\Omega 
= -\lim_{L\to\infty} (\nabla\chi_L\cdot\bu, \varphi)_\Omega.
$$
Let $\varphi_0$ be the zero extension of $\varphi$ to $\BR^N$, that is
$\varphi_0(x) = \varphi(x)$ for $x \in \Omega$ and $\varphi_0(x) = 0$ for 
$x \not\in \Omega$.  Since $\varphi \in H^1_{q',0}(\Omega)$, $\varphi_0
\in \hat H^1_{q'}(\BR^N)$, where 
$$\hat H^1_q(\BR^N)=
\{\psi \in L_{q', {\rm loc}}(\BR^N) \mid \nabla\psi \in L_{q'}(\BR^N)\}.$$
Then, we know (cf. \cite{Gal}) 
that there exist  constants $c\not=0$ and $C$ for which 
$$\Bigl\|\frac{\varphi_0-c}{d}\Bigr\|_{L_{q'}(\BR^N)}
\leq C\|\nabla\varphi_0\|_{L_{q'}(\BR^N)} 
= C\|\nabla\varphi\|_{L_{q'}(\Omega)},
$$
where $d(x) = (1+|x|)\log(2+|x|)$. Noting that
${\rm supp}\, \nabla \chi_L \subset \Omega_R$, we have
\begin{align*}
(\nabla\chi_L\cdot\bu, \varphi)_\Omega
&= c\int_{\dot\Omega_L} (\nabla\chi_L(x))\cdot\bu(x)\,dx \\
&+ \int_{\dot\Omega_L}
(d(x)(\nabla\chi_L)(x)\cdot\bu(x)) \frac{\varphi_0(x)-c}{d(x)}\,dx.
\end{align*}
By Lemma \ref{lem:8.4}, we have
$$\int_{\dot\Omega_L} (\nabla\chi_L(x))\cdot\bu(x)\,dx = 0.
$$
Moreover, by \eqref{6.sob.1} and H\"older's inequality, we have
\begin{align*}
&\Bigl|\int_{\dot\Omega_L}
(d(x)(\nabla\chi_L)(x)\cdot\bu(x)) \frac{\varphi_0(x)-c}{d(x)}\,dx
\Bigr| \\
&\quad \leq \Bigl(\int_{\dot\Omega_L}
(|d(x)(\nabla\chi_L)(x)||\bu(x)|)^q\,dx\Bigr)^{1/q}
\Bigl\|\frac{\varphi_0-c}{d}\Bigr\|_{L_{q'}(\BR^N)} \\
&\quad \leq\frac{C}{\ln \ln L}\Bigl(\int_{\dot\Omega_L}|\bu(x)|^q\,dx
\Bigr)^{1/q}\|\nabla\varphi\|_{L_q(\Omega)} \to 0
\quad\text{as $L\to \infty$}.
\end{align*}
Therefore, we have $(\bu, \nabla\varphi)_\Omega=0$ for any $\varphi \in 
\hat H_{q',0}(\Omega)$, that is $\bu \in J_q(\Omega)$. 
This completes the proof of Proposition \ref{prop:6.1}. 
\qed
\vskip0.5pc
We next consider the weak Dirichlet problem:
\begin{equation}\label{6.weak.1}
(\nabla u, \nabla\varphi)_\Omega = (\bff, \nabla\varphi)_\Omega
\quad\text{for any $\varphi \in \hat H^1_{q',0}(\Omega)$}.
\end{equation}
Then, we know the following fact.
\begin{prop}\label{prop:6.2} Let $1 < q < \infty$ and let $\Omega$ be an 
exterior domain in $\BR^N$ $(N \geq 2)$.  Then, 
the weak Dirichlet problem is uniquely solvable.  Namely, 
for any $\bff \in L_q(\Omega)^N$, problem \eqref{6.weak.1}
admits a unique solution $u \in \hat H^1_{q,0}(\Omega)$
possessing the estimate: $\|\nabla u\|_{L_q(\Omega)} \leq 
C\|\bff\|_{L_q(\Omega)}$. 
\end{prop}
\begin{rem}
\thetag{1}~
This proposition was proved by Pruess and Simonett \cite[Section 7.4]{PS17}
and by Shibata \cite[Theorem 18]{S6} independently. \\
\thetag2~Let $\Omega = \BR^N\setminus S_1$ and $\Gamma=S_1$, where $S_1$
denotes the unit sphere in $\BR^N$.  Let 
$$f(x) = \begin{cases} \ln |x| &\quad N=2, \\
|x|^{-(N-2)}-1 &\quad N \geq 3.
\end{cases}$$
Then, $f(x)$ satisfies the strong Dirichlet problem:
$\Delta f=0$ in $\Omega$ and $f|_{\Gamma} = 0$.  Moreover, 
$f \in H^1_{q, 0}(\Omega)$ provided that $q > N/(N-1)$.  However,
$f$ does not satisfy the weak Dirichlet problem:
$$(\nabla f, \nabla\varphi)_\Omega= 0 \quad\text{for any
$\varphi \in \hat H^1_{q',0}(\Omega)$}.
$$
In fact, $C^\infty_0(\Omega)$ is not dense in $\hat H^1_{q',0}(\Omega)$
when $1 < q' < N$. The detailed is discussed in Shibata \cite[Appendix A]{S6}.
\end{rem}
By Theorem \ref{thm:max.1} in Sect. \ref{sec:3},
we have the following theorem.
\begin{thm}\label{linearthm:6.1} Let $1 < p, q < \infty$ with 
$2/p + 1/q \not=1$ and $0 < T < \infty$. Let 
\begin{equation}\label{compati:0}\begin{split}
\bu_0 &\in B^{2(1-1/p)}_{q,p}(\Omega)^N, \enskip  
\bff \in L_p((0, T), L_q(\Omega)^N), \\ 
g &\in L_p(\BR, H^1_q(\Omega)) \cap H^{1/2}_p(\BR, L_q(\Omega)), 
\quad 
\bg
\in H^1_p(\BR, L_q(\Omega)^N), \\  
\bh &\in H^{1/2}_p(\BR, L_q(\Omega)^N) \cap L_p(\BR, H^1_q(\Omega)^N)
\end{split}\end{equation}
which satisfy the compatibility condition:
$$
\dv\bu_0 = g|_{t=0} \quad\text{in $\Omega$}
$$
and, in addition, 
\begin{equation}\label{6.compati:2}
(\mu\bD(\bu_0)\bn - \bh|_{t=0})_\tau = 0 
\quad\text{on $\Gamma$}
\end{equation}
provided that $2/p + 1/q < 1$. Then, 
problem \eqref{6.max.1} admits unique solutions 
$\bu$ and $\fp$ with
\begin{equation}\label{6.compati:3}\begin{split}
\bu & \in L_p((0, T), H^2_q(\Omega)^N) \cap 
H^1_p((0, T), L_q(\Omega)^N), \\
\fp & \in L_p((0, T), H^1_q(\Omega) + \hat H^1_{q,0}(\Omega))
\end{split}\end{equation} 
satisfying the estimates
\begin{align}
&\|\bu\|_{L_p((0, T), H^2_q(\Omega))}+
\|\pd_t\bu\|_{L_p((0, T), L_q(\Omega))} \nonumber \\
&\quad 
\leq C_\gamma e^{\gamma T}\bigl[
\|\bu_0\|_{B^{2(1-1/p)}_{q,p}(\Omega)} 
 + \|\bff\|_{L_p((0, T), L_q(\Omega)} 
+ \|(g, \bh)\|_{L_p(\BR, H^1_q(\Omega))} \nonumber\\
&\qquad +
\|(g, \bh)\|_{H^{1/2}_p(\BR, L_q(\Omega))} 
+ \|\bg\|_{H^1_p(\BR, L_q(\Omega))} \bigr]
\label{6.est:2.4.0}\end{align}
for some positive constants $C$ and $\gamma$. 
\end{thm}
\begin{rem} 
In the case where $2/p + 1/q < 1$,  $\mu\bD(\bu_0) 
\in B^{1-2/p}_{q,p}(\Omega)$ and $1-2/p > 1/q$, and so  
$\mu\bD(\bu_0)|_\Gamma$ exists. Since 
$\bh \in H^{1/2}_p(\BR, L_q(\Omega)^N) \cap L_p(\BR, H^1_q(\Omega)^N)$,
by complex interpolation theory, $\bh 
\in H^{\theta/2}_p(\BR, H^{1-\theta}_q(\Omega)^N)$ for any 
$\theta \in (0, 1)$.  Since $2/p + 1/q < 1$, we can choose 
$\theta$ in such a way that $1-\theta > 1/q$ and $1/p < \theta/2$.
Thus, $\bh|_{t=0} \in H^{1-\theta}_q(\Omega)^N$, and so 
the trace of $\bh|_{t=0}$ to $\Gamma$ exists. 
\end{rem}
\subsection{Local Well-posedness of Eq. \eqref{navier:3}}
\label{subsec:6.2}

In this subsection, we prove the local well-posedness of 
Eq. \eqref{navier:3}.  The following theorem is the main result of this
subsection.
\begin{thm}\label{6.thm:loc} Let $2 < p < \infty$, $N < q < \infty$ and 
$S > 0$.  Let $\Omega$ be an exterior domain in $\BR^N$ $(N \geq 2)$ whose
boundary $\Gamma$ is a $C^2$ compact hypersurface.  Assume that 
$2/p + N/q < 1$.  Then, there exists a time $T > 0$ depending on
$S$ such that if initial data $\bu_0 \in B^{2(1-1/p)}_{q,p}(\Omega)^N$
satisfies $\|\bu_0\|_{B^{2(1-1/p)}_{q,p}(\Omega)} \leq S$ and 
the compatibility condition:
\begin{equation}\label{comp:4.1}
\dv \bu_0 = 0 \enskip\text{in $\Omega$}, \enskip 
(\bD(\bu_0)\bn)_\tau = 0 
\enskip\text{on $\Gamma$},
\end{equation}
then
problem \eqref{navier:2} admits a unique solution $(\bu, \fq)$ with
\begin{align*}
\bu& \in L_p((0, T), H^2_q(\Omega)^N) \cap H^1_p((0, T), L_q(\Omega)^N), 
\\
\fq  &\in L_p((0, T), H^1_q(\Omega) + \hat H^1_{q,0}(\Omega))
\end{align*}
possessing the estimate:
\begin{gather*}
\|\bu\|_{L_p((0, T), H^2_q(\Omega))}
+ \|\pd_t\bu\|_{L_p((0, T), L_q(\Omega))} \leq CS, \\
\int^T_0\|\kappa(\cdot)\bv(\cdot, s)\|_{H^1_\infty(\Omega)}\,ds \leq\delta
\end{gather*}
for some constant $C > 0$ independent of $T$ and $S$. 
Here, $\delta$ is the constant appearing in \eqref{assump:2}.
\end{thm}
\pf
Let $T$ and $L$ be a 
positive numbers determined late and let 
\begin{align*}
\CI_T = \{&\bu \in L_p((0, T), H^2_q(\Omega)^N)
\cap H^1_p((0, T), L_q(\Omega)^N) \mid \bu|_{t=0} = \bu_0, \\
&<\bu>_T := \|\bu\|_{L_p((0, T), H^2_q(\Omega))} 
+ \|\pd_t\bu\|_{L_p((0, T), L_q(\Omega))}\leq L, \\
&\int^T_0\|\kappa(\cdot)\bu(\cdot, s)\|_{H^1_\infty(\Omega)}\,ds \leq
\delta\}.
\end{align*}
Since $T > 0$ is chosen small enough eventually, we may assume that 
$0 < T \leq 1$.   Given $\bv \in \CI_T$, let $\bu$ be a solution of 
linear equations: 
\begin{equation}\label{6.g.0}
\begin{cases*}
\pd_t\bu - \DV(\mu\bD(\bu) - \fq\bI) & = \bff(\bv)
\quad&\text{in $\Omega^T$}, \\
\dv\bu = g(\bv) & = \dv \bg(\bv)
\quad&\text{in $\Omega^T$}, \\
(\mu\bD(\bu) - \fq\bI)\bn & = \bh(\bv) 
\quad&\text{on $\Gamma^T$}, \\
\bu|_{t=0} & = \bu_0
\quad&\text{in $\Omega$}.
\end{cases*}
\end{equation}
Notice that
\begin{equation}\label{6.g:1}
<\bv>_T \leq L, \quad \int^T_0\|\kappa(\cdot)\bv(\cdot, s)
\|_{H^1_\infty(\Omega)}\,ds \leq \delta.
\end{equation}
To solve \eqref{6.g.0}, we use Theorem \ref{linearthm:6.1}. To this end, 
we introduce $\CE[\bv] = \CE_1[\bv]$, which is the functions 
defined in \eqref{eT:2*} of Sect. \ref{sec:loc6}, and $e_T$, 
which is the extension map defined by \eqref{eT:1}.
Recall that the following formulas hold:
\begin{align}
<\CE_1[\bv]>_\infty \leq C(\|\bu_0\|_{B^{2(1-1/p)}_{q,p}(\Omega)}
+<\bv>_T) \leq C(S+L)&. \label{6.eg.2}
\end{align}
For the sake of simplicity, we may write $g(\bv)$, $\bg(\bv)$ 
and $\bh(\bv)$ given
in Eq. \eqref{non-term:2} and \eqref{non-term:3} as
\begin{align*}
g(\bv) &= v_1(\int^t_0\nabla(\kappa\bv)\,ds)
\int^t_0\nabla(\kappa\bv)\,ds\otimes\nabla\bv, \\
\bg(\bv)  &= \bv_2(\int^t_0\nabla(\kappa\bv)\,ds)
\int^t_0\nabla(\kappa\bv)\,ds\otimes\bv, \\
\bh(\bv) &= 
\bv_3(\int^t_0\nabla(\kappa\bv)\,ds)
\int^t_0\nabla(\kappa\bv)\,ds\otimes \nabla\bv, 
\end{align*}
with some matrices of $C^1$ functions 
$v_1(\bk)$, $\bv_2(\bk)$, and  $\bv_3(\bk)$ defined for
$|\bk| < \delta$.  Note that $\bv_3(\bk)$ depends also on
$x \in \BR^N$ with 
$$\sup_{|\bk| \leq \delta}\|(\bv_3(\cdot, \bk), 
\pd_\bk\bv_3(\cdot, \bk)\|_{H^1_\infty(\BR^N)} \leq C
$$
with some constant $C$, 
because $\bn$ is defined on $\BR^N$ with $\bn\|_{H^1_\infty(\BR^N)}
< \infty$.

Let $e_T$ be the operator defined in \eqref{eT:1}, and 
set 
\begin{alignat*}2
g_0 = e_T[g(\bv)], \quad
\bg_0 = e_T[\bg(\bv)], \quad
\bh_0=e_T[\bh(\bv)].
\end{alignat*}
Since 
\begin{align*}
g_0(\cdot, t) &= \begin{cases} 0 \quad & t < 0, \\
[g(\bv)](\cdot, t) \quad & 0 < t < T, \\
[g(\bv)](\cdot, 2T-t) \quad & T < t < 2T, \\
0 \quad & t > 2T, 
\end{cases}
\\
\bg_0(\cdot, t) &= \begin{cases} 0 \quad & t < 0, \\
[\bg(\bv)](\cdot, t) \quad & 0 < t < T, \\
[\bg(\bv)](\cdot, 2T-t) \quad & T < t < 2T, \\
0 \quad & t > 2T, 
\end{cases}
\end{align*}
as follows from $g(\bv)=0$ and 
$\bg(\bv) = 0$ at $t=0$, 
and since $\dv \bg(\bv) = g(\bv)$ for $0 < t < T$,  we have 
$\dv\bg_0 = g_0$ for any $t \in \BR$.  Moreover,  we have
\begin{equation}\label{6.r.1}\begin{split}
g_0 &= e_T[v_1(\int^t_0\nabla(\kappa\bv)\,ds)\int^t_0
\nabla(\kappa\bv)\,ds\otimes \nabla\CE_1[\bv]], 
\\
\bg_0  &= e_T[\bv_2(\int^t_0\nabla(\kappa\bv)\,ds)\int^t_0
\nabla(\kappa\bv)\,ds\otimes \CE_1[\bv]]\\
\bh_0  &= e_T[\bv_3(\int^t_0\nabla(\kappa\bv)\,ds)
\int^t_0\nabla(\kappa\bv)\,ds\otimes
\nabla\CE_1[\bv]],
\end{split}\end{equation}
because $\CE_1[\bv] = \bv$ in $(0, T)$. 
Where, $\CE_1[\bv]$ has been defined in \eqref{eT:2*}.

Let 
$\bu$ and $\fq$ be solutions of the linear equations:
\begin{equation}\label{6.g.0*}
\begin{cases*}
\pd_t\bu - \DV(\mu\bD(\bu) - \fq\bI) & = \bff
\quad&\text{in $\Omega^T$}, \\
\dv\bu = g_0 & = \dv \bg_0
\quad&\text{in $\Omega^T$}, \\
(\mu\bD(\bu) - \fq\bI)\bn & = \bh_0 
\quad&\text{on $\Gamma^T$}, \\
\bu|_{t=0} & = \bu_0
\quad&\text{in $\Omega$},
\end{cases*}
\end{equation}
and then $\bu$ and $\fq$ are also solutions of the equations \eqref{6.g.0},
because  $g_0 = g(\bv)$, $\bg_0 = \bg(\bv)$ and 
$\bh_0 = \bh(\bv)$ for $t \in (0, T)$. 
 In  the following, 
using the Banach fixed point theorem, we prove
that there exists a unique $\bu \in \CI_T$ such that $\bu = \bv$,
which is a required solution of Eq. \eqref{navier:3}.

Applying Theorem \ref{linearthm:6.1} gives that
\begin{align}
<\bu>_T & \leq C\{\|\bu_0\|_{B^{2(1-1/p)}_{q,p}(\Omega)} + 
\|\bff_0\|_{L_p(\BR, L_q(\Omega))} 
+ \|(g_0, \bh_0)\|_{L_p(\BR, H^1_q(\Omega))} 
\nonumber \\
&\quad + \|(g_0, \bh_0)\|_{H^{1/2}_p(\BR, L_q(\Omega))}
+ \|\pd_t\bg_0\|_{L_p(\BR, L_q(\Omega))}\},
\label{6.g:1*}
\end{align}
provided that the right hand side in \eqref{6.g:1*} is finite. 
In the following, $C$ denotes generic constants independent of 
$T$ and $L$. 
Recalling $\bk = \nabla\int^t_0(\kappa\bv)\,ds$
and the definition of $\bff(\bu)|_i$ given in
\eqref{non-term:1},  we have  
\begin{align}
\|\bff(\bv)\|_{L_q(\Omega)} 
&\leq C\Bigl\{\|\bv(\cdot, t)\|_{L_\infty(\Omega)}
\|\nabla\bv(\cdot, t)\|_{L_q(\Omega)} \nonumber  \\
&+\int^T_0
\|\bv(\cdot, s)\|_{H^1_\infty(\Omega)}
\,ds\|(\nabla^2\bv, \pd_t\bv)\|_{L_q(\Omega)} \nonumber \\
&+ \int^T_0\|\bv(\cdot, s)\|_{H^2_q(\Omega)}\,ds
\|\nabla\bv(\cdot, t)\|_{L_\infty(\Omega)} 
\Bigr\}. \label{6.g:2}
\end{align}
By  H\"older's inequality, 
the Sobolev inequality and the assumption: $2/p+N/q < 1$, 
\begin{align}
\int^T_0\|\bv(\cdot, s)\|_{H^1_\infty(\Omega)}\,ds
&\leq C\Bigl(\int^T_0\|\bv(\cdot, s)\|_{H^2_q(\Omega)}^p\,ds\Bigr)^{1/p}
T^{1/{p'}} \leq CT^{1/{p'}}L, \nonumber \\
\int^T_0\|\bv(\cdot, s)\|_{H^2_q(\Omega)}\,ds
&\leq C\Bigl(\int^T_0\|\bv(\cdot, s)\|_{H^2_q(\Omega)}^p\,ds\Bigr)^{1/p}
T^{1/{p'}}
\leq CT^{1/{p'}}L, \nonumber \\
\|\bv(\cdot, t)\|_{L_\infty(\Omega)}
&\leq C\|\bv(\cdot, t)\|_{H^1_q(\Omega)},
\nonumber \\
\|\nabla\bv(\cdot, t)\|_{L_\infty(\Omega)}
&\leq C\|\bv(\cdot, t)\|_{B^{2(1-1/p)}_{q,p}(\Omega)}.
\label{6.g:3}
\end{align}
Moreover, by \eqref{tanabe:1} and \eqref{g.5.2}, we have 
\begin{equation}\label{6.g:4}\begin{split}
\sup_{0 \leq t \leq T}\|\bv(\cdot, t)\|_{B^{2(1-1/p)}_{q,p}(\Omega)}
&\leq \sup_{t \in (0, \infty)}\|\CE_1[\bv]\|_{B^{2(1-1/p)}_{q,p}(\Omega)}
\leq C(S+L).
\end{split}\end{equation}
 By \eqref{6.g:3} and \eqref{6.g:4}  
\begin{equation}\label{6.g:8}\begin{split}
&\sup_{0\leq t \leq T}(\|\bv(\cdot,t)\|_{L_\infty(\Omega)} 
\|\bv(\cdot,t)\|_{H^1_q(\Omega)})
\leq C\sup_{0 < t < T} \|\bv(\cdot, t)\|_{H^1_q(\Omega)}^2 \\
&\quad \leq C\sup_{0 < t < T}\|\bv(\cdot, t)\|_{B^{2(1-1/p)}_{q,p}(\Omega)}^2
\leq C(S+L)^2,
\end{split}\end{equation}
because $2(1-1/p) > 1$ as follows from $p > 2$
and $N < q < \infty$.  Combining \eqref{6.g:2}, 
\eqref{6.g:3} and \eqref{6.g:8}, we have
\begin{equation}\label{el:1}
\|\bff(\bv)\|_{L_p((0, T), L_q(\Omega))} \leq C(S+L)^2(T^{1/{p'}} + T^{1/p}).
\end{equation}

We next consider the estimate of $g_0$ and $\bh_0$. 
By \eqref{eT:1}, \eqref{6.g:1}, 
\eqref{6.eg.2}, the Sobolev inequality and the assumption: $N < q
< \infty$, we have
\begin{align*}
\|g_0(\cdot, t)\|_{H^1_q(\Omega)}
&\leq C\int^T_0\|\bv(\cdot, s)\|_{H^2_q(\Omega)}\,ds
\|\bv(\cdot, t)\|_{H^2_q(\Omega)} 
\quad\text{for $t \in (0, T)$}, \\ 
\|g_0(\cdot, t)\|_{H^1_q(\Omega)}
&\leq C\int^{2T-t}_0\|\bv(\cdot, s)\|_{H^2_q(\Omega)}\,ds
\|\bv(\cdot, 2T-t)\|_{H^2_q(\Omega)}
\end{align*}
for $t \in (T, 2T)$, and 
$\|g_0(\cdot, t)\|_{H^1_q(\Omega)}= 0$ for $t \not\in [0, 2T]$. 
Thus, 
$$
\|g_0(\cdot, t)\|_{H^1_q(\Omega)}
\leq C\begin{cases} 0 \quad & t < 0, \\
T^{1/p'}<\bv>_T
\|\bv(\cdot, t)\|_{H^2_q(\Omega)} 
\quad& 0 < t < T, \\
T^{1/p'}<\bv>_T
\|\bv(\cdot, 2T-t)\|_{H^2_q(\Omega)} 
\quad& T < t < 2T, \\
0 \quad & t > 2T,
\end{cases}$$
which, combined with \eqref{6.g:1},  gives that 
\begin{equation}\label{el:2}
\|g_0\|_{L_p(\BR, H^1_q(\Omega))} \leq CL^2T^{1/{p'}}.
\end{equation}
Analogously, 
\begin{equation}\label{el:3}
\|\bh_0\|_{L_p(\BR, H^1_q(\Omega))} \leq CL^2T^{1/{p'}}.
\end{equation}

Since 
\begin{align*}
\|\pd_t\int^t_0\nabla(\kappa\bv)\,ds\|_{L_p((0, T), H^1_q(\Omega))}
& \leq C<\bv>_T\leq CL, \\
\|\pd_t\int^t_0\nabla(\kappa\bv)\,ds\|_{L_\infty((0, T),L_q(\Omega))}
&\leq C\|\bv\|_{L_\infty((0, T), H^1_q(\Omega))}
\leq C(S+L)
\end{align*}
as follows from \eqref{6.g:4},  applying Lemma \ref{lem:g.5.0}
and Lemma \ref{lem:g.5.1} 
to $g_0$ and $\bh_0$, and using the formula of $g_0$ and $\bh_0$ given 
in \eqref{6.r.1} and the estimates \eqref{6.g:1}
we have 
\begin{align}
\|g_0\|_{H^{1/2}_p(\BR, L_q(\Omega))}
&\leq C(L+S)^2(T^{1/{p'}} + T^{\frac{q-N}{pq}}),
\label{el:6.g4} \\
\|\bh_0\|_{H^{1/2}_p(\BR, L_q(\Omega))}
&\leq C(L+S)L(T^{1/{p'}} + T^{\frac{q-N}{pq}}).
\label{el:6.g5}
\end{align}

We finally estimate $\pd_t\bg_0$.  To this end,
we write
\begin{align*}
\pd_t\bg_0 &= 
\bv_2(\int^t_0\nabla(\kappa\bv)\,ds)
\int^t_0\nabla(\kappa\bv)\,ds\otimes\pd_t\CE_1[\bv](\cdot, t) \\
&+ 
\bv_2(\int^t_0\nabla(\kappa\bv)\,ds)
\nabla(\kappa\bv)\otimes\CE_1[\bv](\cdot, t) \\
&
+ \bv'_2(\int^t_0\nabla(\kappa\bv)\,ds)\nabla(\kappa\bv)
\int^t_0\nabla(\kappa\bv)\,ds\otimes\CE_1[\bv](\cdot, t))
 \quad\text{for $t \in (0, T)$}, \\
\pd_t\bg_0 
& = \bv_2(\int^{2T-t}_0\nabla(\kappa\bv)\,ds)
\int^{2T-t}_0\nabla(\kappa\bv)\,ds\otimes\pd_t\CE_1[\bv](\cdot, t) \\
&+ 
\bv_2(\int^{2T-t}_0\nabla(\kappa\bv)\,ds)
\nabla(\kappa\bv)\otimes\CE_1[\bv](\cdot, t) \\
&
+ \bv'_2(\int^{2T-t}_0\nabla(\kappa\bv)\,ds)\nabla(\kappa\bv)
\int^{2T-t}_0\nabla(\kappa\bv)\,ds\otimes\CE_1[\bv](\cdot, t))
\end{align*}
for $t \in (T, 2T)$, 
and $\pd_t \bg_0 = 0$ for $t \not\in [0, 2T]$, 
where $\bv'_2(\bk) = \nabla_\bk \bv_2(\bk)$. By \eqref{6.g:4}, 
$$\|\pd_t\bg_0\|_{L_p(\BR, L_q(\Omega))}
\leq C(T^{\frac{1}{p'}} + T^{\frac{1}{p}})(L+S)^2,
$$
which, combined with \eqref{6.g:1*}, \eqref{el:1},
\eqref{el:2}, \eqref{el:3}, \eqref{el:6.g4}, and \eqref{el:6.g5},
leads to  
$$<\bu>_T \leq C(L+S)^2(T^{\frac{1}{p'}} + T^{\frac{1}{p}}
+ T^{\frac{q-N}{pq}}).
$$
Choosing $T > 0$ so small that
$$C(L+S)(T^{\frac{1}{p'}} + T^{\frac{1}{p}}
+ T^{\frac{q-N}{pq}}) \leq 1/2,$$
we have $<\bu>_T \leq (C+1/2)S+L/2$. 
Thus, choosing $L = (2C+1)S$, we have 
\begin{equation}\label{6.g.main}<\bu>_T \leq L.\end{equation}
Moreover, we have 
\begin{align*}
\int^T_0\|\kappa(\cdot)\bu(\cdot, s)\|_{H^1_\infty(\Omega)}\,ds
&\leq C_q\|\kappa\|_{H^1_\infty(\Omega)}T^{1/{p'}}
\Bigl(\int^T_0\|\bu(\cdot, s)\|_{H^2_q(\Omega)}^p\,ds\Bigr)^{1/p}\\
&\leq C_q\|\kappa\|_{H^1_\infty(\Omega)}LT^{1/{p'}},
\end{align*}
and so choosing $T > 0$ in such a way that 
$C_q\|\kappa\|_{H^1_\infty(\Omega)}LT^{1/{p'}} \leq \delta$,
we have 
$$\int^T_0\|\kappa(\cdot)\bu(\cdot, s)\|_{H^1_\infty(\Omega)}\,ds\leq \delta.$$
Thus,  $\bu \in \CI_T$.  Let $\bQ$ be a map defined by
$\bQ\bv = \bu$, and then $\bQ$ maps $\CI_T$ into
itself. 

Given $\bv_i \in \CI_T$ ($i=1,2$), considering the equations satisfied
by $\bu_2 - \bu_1 = \bQ\bv_2- \bQ\bv_1$ and employing the same argument
as that in proving \eqref{6.g.main}, 
 we can show that 
$$<\bQ\bv_1- \bQ\bv_2>_T \leq 
C(L+S)(T^{\frac1{p'}} + T^{\frac1p} + 
T^{\frac{q-N}{pq}})<\bv_1- \bv_2>_T
$$
holds. 
Choosing $T$ smaller if necessary, we may assume that
$C(L+S)(T^{\frac{1}{p'}} + T^{\frac1p} + 
T^{\frac{q-N}{pq}}) \leq 1/2$, and so 
$\bQ$ is a contration map on $\CI_T$.  By
the Banach fixed point theorem, there exists a 
unique $\bu \in \CI_T$ such that $\bQ\bu = \bu$, 
which is a required unique solution of Eq.
\eqref{navier:2}.  This completes the proof of Theorem \ref{6.thm:loc}.
\qed \vskip0.5pc

Employing the similar argument to that in Subsec. \ref{subsec:5.2}, 
 we can prove the following theorem, which is used to 
prove the global well-posedness. 
\begin{thm}\label{thm:6.gloc*} Let $2 < p < \infty$, $N < q < \infty$ and 
$T > 0$.  Let $\Omega$ be an exterior domain in $\BR^N$ $(N \geq 2)$, whose
boundary $\Gamma$ is a $C^2$ compact hypersurface.  Assume that 
$2/p + N/q < 1$.  Then, there exists an $\epsilon_1 > 0$ depending on
$T$ such that if initial data $\bu_0 \in B^{2(1-1/p)}_{q,p}(\Omega)^N$
satisfies $\|\bu_0\|_{B^{2(1-1/p)}_{q,p}(\Omega)} \leq \epsilon_1$ and 
the compatibility condition \eqref{comp:4.1}, 
then
problem \eqref{navier:2} admits a unique solution $(\bu, \fq)$ with
\begin{align*}
\bu& \in L_p((0, T), H^2_q(\Omega)^N) \cap H^1_p((0, T), L_q(\Omega)^N), \\
\fq & \in L_p((0, T), H^1_q(\Omega) + \hat H^1_{q,0}(\Omega))
\end{align*}
possessing the estimate:
$$\|\bu\|_{L_p((0, T), H^2_q(\Omega))}
+ \|\pd_t\bu\|_{L_p((0, T), L_q(\Omega))} \leq \epsilon_1,\quad
\int^T_0\|\kappa(\cdot)\bu(\cdot, s)\|_{H^1_\infty(\Omega)} \leq\delta.
$$
\end{thm} 
\subsection{A new formulation of Eq. \eqref{navier:3}}\label{subsec:6.3}
Let $T > 0$ and let 
\begin{equation}\label{6.3:1}\begin{split}
\bu &\in H^1_p((0, T), L_q(\Omega)^N) \cap 
L_p((0, T), H^2_q(\Omega)^N), \\ \fq
&\in L_p((0, T), H^1_q(\Omega) + \hat H^1_{q,0}(\Omega))
\end{split}\end{equation}
 be solutions of Eq. \eqref{navier:3} satisfying the condition
\eqref{assump:2}. In what follows, we rewrite Eq. \eqref{navier:3}
in order that the nonlinear terms have suitable decay properties. 

In the following, we repeat the argument in Subsec: \ref{sub:1} and
Subsec: \ref{subsec:2.3}. Let $\bV_0(\bk) = (V_{0ij}(\bk))$ be the 
$N\times N$ matrix defined in \eqref{trans:2} in Subsec. \ref{sub:1}
and set  $\dfrac{\pd y}{\pd x} = \bI + \bV_0(\bk)
: = (a_{ij}(t))$, 
 where $\bk=\{k_{ij} \mid i,j=1, \ldots, N\}$ are
 the variables corresponding to 
$\int^t_0\nabla(\kappa(y)\bv(y, s))\,ds$. 
 And also, let  $\bn_t = {}^\top(n_{t1},
\ldots, n_{tN})$ and $\bn = {}^\top(n_1, \ldots, n_N)$
be the unit outer normal to $\Gamma_t$ and $\Gamma$, respectively.
Since
\begin{align*}
0 = \bn_t\cdot dx = \sum_{j=1}^N n_{tj}dx_j 
= \sum_{j,k=1}^N n_{tj}\frac{\pd x_j}{\pd y_k}dy_k,
\end{align*}
we see that ${}^\top\dfrac{\pd x}{\pd y}\bn_t$ is parallel to $\bn$, 
that is ${}^\top\dfrac{\pd x}{\pd y}\bn_t = c\bn$ for some 
$c \in \BR\setminus\{0\}$, and so we have
$\bn_t = c{}^\top\dfrac{\pd y}{\pd x}\bn$.  Since $|\bn_t| = 1$, 
we have \eqref{normal:8.1}. 
Thus, by \eqref{trans:2} and \eqref{change:1*} we have
\begin{equation}\label{6.trans:1}
\frac{\pd}{\pd x_i} = \sum_{j=1}^N
a_{ji}(t)\frac{\pd}{\pd y_j}, \quad
n_{ti} = d(t)\sum_{j=1}^N a_{ji}(t)n_j
\end{equation}
where $d(t)
=|(\bI + {}^\top\bV_0(\bk))\bn|$. Let 
$J$ be the Jacobian of the partial Lagrange transform \eqref{plag:1}
and set 
$$\ell_{ij} = \delta_{ij} + \int^t_0 \frac{\pd}{\pd y_j}(\kappa(y)
u_i(y, s))\,ds
$$
where $\bu={}^\top(u_1, \ldots, u_N)$. Let 
\begin{equation}\label{6.trans:0}
a_{ij}(t) = \delta_{ij} + \tilde a_{ij}(t), \enskip
J(t) = 1+ \tilde J(t), \enskip 
\ell_{ij}(t) = \delta_{ij} + \tilde \ell_{ij}(t)
\end{equation}
with  
\begin{equation}\label{6.trans:4}\begin{split}
\tilde a_{ij}(t) &= b_{ij}(\int^t_0\nabla( \kappa(y) \bv(y, s))\,ds), 
\quad
\tilde J(t) = K(\int^t_0\nabla( \kappa(y) \bv(y, s))\,ds), 
\\ 
\tilde \ell_{ij}(t) &= m_{ij}(\int^t_0
\nabla( \kappa(y) \bv(y, s))\,ds):= \int^t_0\frac{\pd}{\pd y_j}
(\kappa(y)u_i(y, t))\,ds.
\end{split}\end{equation}
Here, if we use the symbols defined in  Subsec. \ref{sub:1},  then
 $b_{ij} = V_{0ij}$ and $K=J_0$, and so  $b_{ij}$ and $K$ are 
smooth functions defined 
on $\{\bk \mid |\bk| \leq \delta\}$ such that $b_{ij}(0) = K(0)
 = 0$. 

Let $\bv(x, t) = \bu(y,t)$ and $\fp(x, t) = \fq(y, t)$, 
and then $\bv$ and $\fp$ are solutions of Eq. \eqref{navier:1} with
\begin{align*}
\Omega_t &= \{x = y + \int^t_0\kappa(y)\bu(y, s)\,ds
\mid y \in \Omega\}, 
\\
\Gamma_t &=\{x = y + \int^t_0\kappa(y)\bu(y, s)\,ds
\mid y \in \Gamma\}.
\end{align*}
 By \eqref{6.trans:1}, 
$$
\frac{\pd  v_i}{\pd x_j} + \frac{\pd v_j}{\pd x_i}
= D_{ij, t}(\bu) : = D_{ij}(\bu) + 
\tilde D_{ij}(t)\nabla\bu
$$
with  
\begin{equation}\label{6.trans:2}
D_{ij}(\bu) = \frac{\pd u_i}{\pd y_j}
+ \frac{\pd u_j}{\pd y_j}, \quad 
\tilde D_{ij}(t)\nabla\bu = \sum_{k=1}^N(\tilde a_{kj}(t)
\frac{\pd u_i}{\pd y_k} + \tilde a_{ki}(t)
\frac{\pd u_j}{\pd y_k}).
\end{equation}
By \eqref{div:2} in Subsec. \ref{sub:1}
we also have an important formula:
\begin{equation}\label{6.trans:5}
\dv\bv = \sum_{j=1}^N\frac{\pd v_j}{\pd x_j}
= \sum_{j,k=1}^N J(t)a_{kj}(t)\frac{\pd u_j}{\pd y_k}
= \sum_{j,k=1}^N \frac{\pd}{\pd y_k}(J(t)a_{kj}(t)u_j),
\end{equation}
which implies that 
\begin{equation}\label{6.trans:6}
\sum_{j,k=1}^N(\tilde a_{kj}(t) + \tilde J(t)a_{kj}(t))
\frac{\pd u_j}{\pd y_k}
= \sum_{j,k=1}^N\frac{\pd}{\pd y_k}\{(\tilde a_{kj}(t)
+ \tilde J(t)a_{kj}(t))u_j\}.
\end{equation} 
And then, Eq. \eqref{navier:3} is written as follows:
\begin{equation}\label{navier:2*}\begin{aligned}
\sum_{i=1}^N \ell_{is}(t)(\pd_t u_i 
+ (1-\kappa)\sum_{j,k=1}^Nu_ja_{kj}(t)\frac{\pd u_i}{\pd y_k})
\\
\qquad - \mu\sum_{i,j.k=1}^N \ell_{is}(t)a_{kj}(t)\frac{\pd}{\pd y_k}
D_{ij,t}(\bu) - \frac{\pd \fq}{\pd y_s} = 0&
&\quad&\text{in $\Omega^T$}, \\
\sum_{j,k=1}^N J(t)a_{kj}(t)\frac{\pd u_j}{\pd y_k}
= \sum_{j,k=1}^N \frac{\pd}{\pd y_k}(J(t)a_{kj}(t)u_j)
= 0& &\quad&\text{in $\Omega^T$}, \\
\mu\sum_{i,j,k=1}^N \ell_{is}(t)a_{kj}(t)D_{ij, t}(\bu)n_k
- \fq n_s = 0& &\quad &\text{on $\Gamma^T$}, \nonumber \\
\bu|_{t=0} = \bu_0& &\quad &\text{in $\Omega$},
\end{aligned}\end{equation}
where $s$ runs from $1$ through $N$. Here, we have used the fact that
$(\ell_{ij}) = \bA^{-1}$.


In order to get some decay properties of the nonlinear
terms, we write
$$\int^t_0 \nabla(\kappa(y)\bu(y, s))\,ds
=\int^T_0 \nabla(\kappa(y)\bu(y, s))\,ds -
\int^T_t \nabla(\kappa(y)\bu(y, s))\,ds.
$$ 
In \eqref{6.trans:4}, by the Taylor formula we write
\begin{alignat}2
 a_{ij}(t) &=  a_{ij}(T) + \CA_{ij}(t), &\quad 
 \ell_{ij}(t) &=  \ell_{ij}(T) + \CL_{ij}(t), \nonumber \\
D_{ij, t}(\bu) &= D_{ij, T}(\bu) + \CD_{ij}(t)\nabla\bu,
&\quad 
J(t) &=  J(T) + \CJ(t) \label{6.trans:7}
\end{alignat}
with
\begin{align*}
&\CA_{ij}(t) = -\int^1_0b_{ij}'(\int^T_0\nabla(\kappa(y)\bu(y, s))\,ds
-\theta\int^T_t\nabla(\kappa(y)\bu(y, s))\,ds)\,d\theta\, \\
&\phantom{= -\int^1_0b_{ij}'(\int^T_0\nabla(\kappa(y)\bu(y, s))\,ds}
\times
\int^T_t\nabla(\kappa(y)\bu(y, s))\,ds \nonumber\\
&\CL_{ij}(t) = -\int^T_t\frac{\pd}{\pd y_j}(\kappa(y)u_i(y, s))\,ds,
\\
&\CD_{ij}(t)\nabla\bu
= \sum_{k=1}^N(\CA_{kj}(t)\frac{\pd u_i}{\pd y_k}
+ \CA_{ki}(t)\frac{\pd u_j}{\pd y_k}), \\
&\CJ(t)  = -\int^1_0 K'(\int^T_0\nabla(\kappa(y)\bu(y, s))\,ds
-\theta\int^T_t\nabla(\kappa(y)\bu(y, s))\,ds)\,d\theta\, \\
&\phantom{= -\int^1_0b_{ij}'(\int^T_0\nabla(\kappa(y)\bu(y, s))\,ds}
\times
\int^T_t\nabla(\kappa(y)\bu(y, s))\,ds, 
\end{align*}
where $b'_{ij}$ and $K'$ are derivatives of $b_{ij}$ and $K$ with respect
to $\bk$. By  the relation: 
\begin{equation}\label{inverse:1}
\sum_{s=1}^N\ell_{is}(T)a_{sm}(T) = \delta_{si},
\end{equation}
 the first equation in \eqref{navier:2*} 
is rewritten as follows: 
$$\pd_tu_m - \sum_{j,k=1}^Na_{kj}(T)
\frac{\pd}{\pd y_k}(\mu D_{mj, T}(\bu) -\delta_{mj}\fq) = \tilde f_m(\bu)$$
with
\begin{align}
&\tilde f_m(\bu)= -\sum_{s=1}^Na_{sm}(T)\{\sum_{i=1}^N\CL_{is}(t)\pd_tu_i + 
\sum_{i,j,k=1}^N(1-\kappa)\ell_{is}(t)a_{kj}(t)u_j
\frac{\pd u_i}{\pd y_k}\}\nonumber \\ 
&\quad+\mu\sum_{s=1}^Na_{sm}(T)\Bigl\{\sum_{i,j,k=1}^N\ell_{is}(T)a_{kj}(T)
\frac{\pd}{\pd y_k}(\CD_{ij}(t)\nabla\bu) \label{non:1}\\
&\quad+ \sum_{i,j,k=1}^N \ell_{is}(T)\CA_{kj}(t)\frac{\pd}{\pd y_k}
D_{ij, t}(\bu) 
 +\sum_{i,j,k=1}^N \CL_{is}(t)a_{kj}(t)\frac{\pd}{\pd y_k}D_{ij,t}
(\bu)\Bigr\}. \nonumber
\end{align}
Next, by \eqref{6.trans:5}
$$\widetilde{\dv} \bu = \tilde g(\bu) = \dv\tilde \bg(\bu)$$
with
\begin{align}
\widetilde{\dv}\bu &= \sum_{j,k=1}^NJ(T)a_{kj}(T)\frac{\pd u_j}{\pd y_k}
= \sum_{j,k=1}^N\frac{\pd}{\pd y_k}(J(T)a_{kj}(T)u_j), 
\nonumber \\
\tilde g(\bu) & = \sum_{j,k=1}^N(J(T)\CA_{kj}(t) + \CJ(t)a_{kj}(t))
\frac{\pd u_j}{\pd y_k}, \label{non:1*} \\
\tilde g_k(\bu) & = 
\sum_{j=1}^N(J(T)\CA_{kj}(t) + \CJ(t)a_{kj}(t))u_j,
\quad \tilde\bg(\bu) = {}^\top(\tilde g_1(\bu), \ldots, \tilde g_N(\bu)).
\nonumber
\end{align}

Finally, we consider the boundary condition. 
Let $\tilde\bn$ be an $N$-vector defined on $\BR^N$ such 
that $\tilde\bn = \bn$ on $\Gamma$ and $\|\tilde\bn\|_{H^2_\infty(\BR^N)}
\leq C$.  In what follows, $\tilde \bn$ is simply written
by $\bn = {}^\top(n_1, \ldots, n_N)$. 
By \eqref{6.trans:1} and \eqref{inverse:1}
$$\sum_{j,k=1}^Na_{kj}(T)(\mu D_{mj,T}(\bu)
- \delta_{mj}\fq)n_k= \tilde h_m(\bu)$$
with
\begin{equation}\label{non:1**}\begin{split}
\tilde h_m(\bu) &= -\mu\sum_{j,k=1}^N(a_{kj}(T)\CD_{mj}(t)\nabla\bu
+ \CA_{kj}(t)D_{mj,t}(\bu))n_k \\
&-\mu\sum_{i,j,k,s=1}^Na_{sm}(T)\CL_{is}(t)a_{kj}(t)D_{ij,t}(\bu)n_k.
\end{split}\end{equation}
By \eqref{6.trans:5}, 
\begin{align*}
&\sum_{j,k=1}^Na_{kj}(T)
\frac{\pd}{\pd y_k}(\mu D_{mj, T}(\bu) -\delta_{mj}\fq) \\
&\quad
= J(T)^{-1}\sum_{k=1}^N\frac{\pd}{\pd y_k}
[\sum_{j=1}^N\{J(T)a_{kj}(T)(\mu D_{mj, T}(\bu) - \delta_{mj}\fq)\}].
\end{align*}
Thus, letting 
\begin{align*}
S_{mk}(\bu, \fq) &= \sum_{j=1}^NJ(T)a_{kj}(T)
(D_{mj,T}(\bu) - \delta_{mj}\fq), \quad 
\tilde \bS(\bu, \fq) = (S_{ij}(\bu, \fq)), 
\\
\tilde\bff(\bu) &= {}^\top(\tilde f_1(\bu), \ldots, \tilde f_N(\bu)), \quad 
\tilde \bh(\bu) = {}^\top(\tilde h_1(\bu), \ldots, \tilde h_N(\bu)),
\end{align*} 
and using \eqref{6.trans:5}, 
we see that $\bu$ and $\fq$ satisfy the following equations:
\begin{equation}\label{neweq:4}\left\{\begin{aligned}
\pd_t\bu - J(T)^{-1}\DV\tilde\bS(\bu, \fq) = \tilde\bff(\bu)&
&\quad&\text{in $\Omega^T$},\\
\widetilde{\dv}\bu = \tilde g(\bu) = \dv\tilde\bg(\bu)&
&\quad&\text{in $\Omega^T$}, 
\\
\tilde\bS(\bu, \fq)\bn  = J(T)\tilde\bh(\bu)&
&\quad&\text{on $\Gamma^T$}, \\
\bu|_{t=0}  = \bu_0&
&\quad &\text{in $\Omega$}.
\end{aligned}\right.
\end{equation}
This is a new formula of Eq. \eqref{navier:3} which is  satisfied
by local in time solutions
$\bu$ and $\fq$ of Eq. \eqref{navier:3}. 
The corresponding linear equations to Eq. \eqref{neweq:4} 
is the followings:
\begin{equation}\label{neweq:7}\left\{\begin{aligned}
\pd_t\bu - J(T)^{-1}\DV\tilde\bS(\bu, \fq) = \bff&
&\quad&\text{in $\Omega^T$},\\
\widetilde{\dv}\bu = g = \dv\bg&
&\quad&\text{in $\Omega^T$}, 
\\
\tilde\bS(\bu, \fq)\bn  =\bh&
&\quad&\text{on $\Gamma^T$}, \\
\bu|_{t=0}  = \bu_0&
&\quad &\text{in $\Omega$}.
\end{aligned}\right.
\end{equation}
We call
Eq. \eqref{neweq:7} the slightly perturbed Stokes equations. 
\subsection{Slightly perturbed Stokes equations} \label{subsec:6.4}
In this subsection we summarize some results obtained by
Shibata \cite{S7} concerning
the slightly perturbed Stokes equations.  Let $r$ be an exponent 
such that $N < r < \infty$.  Let $a_{ij}(T)$, $\tilde a_{ij}(T)$, 
$J(T)$ and $\tilde J(T)$
be functions defined in \eqref{6.trans:0} with $t=T$.  We assume that
\begin{equation}\label{sassump:1}
\|(\tilde a_{ij}(T), \tilde J(T))\|_{L_\infty(\Omega)}
+ \|\nabla(\tilde a_{ij}(T), \tilde J(T))\|_{L_r(\Omega)} \leq
\sigma
\end{equation}
with some small constant $\sigma > 0$.  In view of Theorem \ref{thm:6.gloc*},
we can choose $\sigma > 0$ small as much as we want if we choose the initial
data small. Since ${\rm supp}\,\kappa \subset B_{2R}$, 
$\tilde a_{ij}(T)$ and $\tilde J(T)$ vanish for $x \not\in B_{2R}$.
In the following, we write $a_{ij}(T)$, $\tilde a_{ij}(T)$, 
$J(T)$ and $\tilde J(T)$ simply by
$a_{ij}$, $\tilde a_{ij}$, 
$J$ and $\tilde J$, respectively. And also, we write $\bA = (a_{ij}(T))$.

To state the compatibility condition for Eq. \eqref{neweq:7}, 
we modify the equation slightly. 
 Notice that
it follows from \eqref{6.trans:5} with $u_j =\delta_{mj}\fq$ that 
\begin{align*}
\sum_{j,k=1}^NJ^{-1}\frac{\pd}{\pd y_k}(Ja_{kj}\delta_{mj}\fq)
&= \sum_{j,k=1}^Na_{kj}\frac{\pd}{\pd y_k}(\delta_{mj}\fq)
= \sum_{k=1}^Na_{km}\frac{\pd\fq }{\pd y_k} \\
&= {}^\top\bA\nabla \fq|_m.
\end{align*}
In the following, we set 
\begin{gather*}
\tilde\nabla\fq = {}^\top\bA\nabla, 
\enskip
\tilde \bD(\bu) = (D_{ij, T}(\bu)), \enskip
\tilde \bn = d_\bn^{-1}{}^\top\bA\bn, 
\enskip \tilde n_i = d_\bn^{-1}\sum_{j=1}^Na_{ji}(T)n_j,
\end{gather*}
with $d_\bn = (\sum_{i,j=1}^Na_{ij}(T)n_in_j)^{1/2}$,
where we have set $\tilde\bn={}^\top(\tilde n_1, \ldots, \tilde n_N)$ and 
$\bn = {}^\top(n_1, \ldots, n_N)$. Note that 
$|\tilde\bn| = 1$.  And then, 
we can write
\begin{equation}\label{sform:1}\begin{split}
J^{-1}\DV\tilde \bS(\bu, \fq) &=  J^{-1}\DV(J\bA\mu\tilde \bD(\bu))
-\tilde\nabla\fq, \\
\tilde \bS(\bu, \fq)\bn  &= Jd_\bn (\mu\tilde\bD(\bu) - \fq\bI)\tilde \bn.
\end{split}\end{equation}
And then, Eq. \eqref{neweq:7} is rewritten as
\begin{equation}\label{neweq:7.7}\left\{\begin{aligned}
\pd_t\bu -  J^{-1}\DV(J\bA\mu\tilde \bD(\bu))
+\tilde\nabla\fq = \bff&
&\quad&\text{in $\Omega^T$},\\
\widetilde{\dv}\bu = g = \dv\bg&
&\quad&\text{in $\Omega^T$}, 
\\
\mu\tilde\bD(\bu)\tilde\bn  - \fq\tilde \bn =(Jd_\bn)^{-1}\bh&
&\quad&\text{on $\Gamma^T$}, \\
\bu|_{t=0}  = \bu_0&
&\quad &\text{in $\Omega$}.
\end{aligned}\right.
\end{equation}

To obtain the maximal $L_p$-$L_q$ regularity and 
some decay properties of solutions of Eq. \eqref{neweq:7.7}, 
we first 
consider the following generalized resolvent problem corresponding
to Eq. \eqref{neweq:4}:
\begin{equation}\label{spres:1}\left\{\begin{aligned}
\lambda \bu -   J^{-1}\DV(J\bA\mu\tilde \bD(\bu))
+\tilde\nabla\fq = \bff&
&\quad&\text{in $\Omega^T$},\\
\widetilde{\dv}\bu = g = \dv\bg&
&\quad&\text{in $\Omega^T$}, 
\\
\mu\tilde\bD(\bu)\tilde\bn  - \fq\tilde \bn =(Jd_\bn)^{-1}\bh&
&\quad&\text{on $\Gamma^T$}, \\
\bu|_{t=0}  = \bu_0&
&\quad &\text{in $\Omega$}.
\end{aligned}\right.
\end{equation}

The following theorem was proved by Shibata \cite{S7}.
\begin{thm}\label{spthm:1} Let $1 < q \leq r$ and $0 < \epsilon_0 < 
\pi/2$.  Assume that $\Omega$ is an exterior domain whose boundary
$\Gamma$ is a compact $C^2$ hypersurface.  Let 
\begin{align*}
X''_q(\Omega) & = \{(\bff, g, \bg, \bh) \mid \bff, \bg \in L_q(\Omega)^N,
\quad g \in H^1_q(\Omega), \quad \bh \in H^1_q(\Omega)^N\},\\
\CX''_q(\Omega) & = \{(F_1, \ldots, F_N) \mid F_1, F_4, F_6 \in L_q(\Omega)^N,
\quad F_2 \in H^1_q(\Omega), \\
&\phantom{= \{(F_1, \ldots, F_N) \mid }
 F_3 \in L_q(\Omega),\quad
F_5 \in H^1_q(\Omega)^N\}.
\end{align*}
Then, there exist a constant $\lambda_0 > 0$ 
and operator families $\CA_s(\lambda)$
and $\CP_s(\lambda)$ with
\begin{align*}
\CA_s(\lambda) &\in \Hol(\Sigma_{\epsilon_0, \lambda_0},
\CL(\CX''_q(\Omega), H^2_q(\Omega)^N)), 
\\
\CP_s(\lambda) &\in \Hol(\Sigma_{\epsilon_0, \lambda_0}, 
\CL(\CX''_q(\Omega), H^1_q(\Omega) + \hat H^1_{q,0}(\Omega)))
\end{align*}
such that for any $\lambda \in \Sigma_{\epsilon_0, \lambda_0}$ and 
$(\bff, g, \bg, \bh) \in X''_q(\Omega)$, 
$\bu= \CA_s(\lambda)\bF_\lambda$ and $\fq = \CP_s(\lambda)\bF_\lambda$
are unique solutions of Eq. \eqref{spres:1}, where
$$\bF_\lambda = (\bff, g, \lambda^{1/2}g, \lambda\bg, \bh, \lambda^{1/2}\bh).$$
Moreover, we have
\begin{align*}
\CR_{\CL(\CX''_q(\Omega), H^{2-j}_q(\Omega)^N)}
(\{(\tau\pd_\tau)^\ell(\lambda^{j/2} \CA_s(\lambda)) \mid 
\lambda \in \Sigma_{\epsilon_0, \lambda_0}\})& \leq r_b, \\
\CR_{\CL(\CX''_q(\Omega), L_q(\Omega)^N)}
(\{(\tau\pd_\tau)^\ell( \nabla\CP_s(\lambda)) \mid 
\lambda \in \Sigma_{\epsilon_0, \lambda_0}\})& \leq r_b
\end{align*}
for $\ell=0,1$ and $j=0,1,2$ with 
 some constant $r_b$.
\end{thm}
To obtain the decay properties of solutions to Eq. \eqref{neweq:7.7},
we first consider the time shifted equations:
\begin{equation}\label{s-neweq:7.7}\left\{\begin{aligned}
\pd_t\bu + \lambda_0\bu -   J^{-1}\DV(J\bA\mu\tilde \bD(\bu))
+\tilde\nabla\fq = \bff&
&\quad&\text{in $\Omega^T$},\\
\widetilde{\dv}\bu = g = \dv\bg&
&\quad&\text{in $\Omega^T$}, 
\\
\mu\tilde\bD(\bu)\tilde\bn  - \fq\tilde \bn =(Jd_\bn)^{-1}\bh&
&\quad&\text{on $\Gamma^T$}, \\
\bu|_{t=0}  = \bu_0&
&\quad &\text{in $\Omega$}.
\end{aligned}\right.
\end{equation}
Employing the same argument as in Subsec. \ref{subsec:3.5},
we have the following maximal $L_p$-$L_q$ regularity theorem
for Eq. \eqref{s-neweq:7.7} with large $\lambda_0 > 0$. 
\begin{thm}\label{sthm:max} Let $1<p, q < \infty$ and assume that
$2/p + N/q \not=1$.  Then, there exist constants $\sigma > 0$ 
and $\lambda_0 > 0$ such that
if  \eqref{sassump:1} holds, then the following
assertion holds:   
Let $\bu_0 \in B^{2(1-1/p)}_{q,p}(\Omega)^N$ be initial data for Eq. 
\eqref{s-neweq:7.7} and let $\bff$, $g$, $\bg$, $d$, $\bh$ be given functions 
for Eq. \eqref{s-neweq:7.7} with 
\begin{alignat*}2
\bff &\in L_p(\BR, L_q(\Omega)^N), &\quad 
g &\in H^1_p(\BR, H^1_q(\Omega)) \cap H^{1/2}_p(\BR, L_q(B_R)), \\
\bg &\in H^1_p(\BR, L_q(\Omega)^N), &\quad 
\bh 
&
\in H^1_p(\BR, H^1_q(\Omega)^N) \cap H^{1/2}_p(\BR, L_q(\Omega)^N).
\end{alignat*}
 Assume that the compatibility conditions:
$$
\widetilde{\dv}\bu_0 = g|_{t=0} \quad\text{in $\Omega$}.$$
When $2/p + 1/q < 1$, in addition we assume that 
$$Jd_\bn\tilde\BPi_0(\mu\tilde \bD(\bu_0)\tilde\bn) = \tilde\BPi_0(\bh|_{t=0})
\quad\text{on $\Gamma$}, 
$$
where for any $N$-vector $\bd$ we set 
$\tilde\BPi_0\bd = \bd - <\bd, \tilde\bn>\tilde\bn$. 
Then, problem \eqref{s-neweq:7.7} admits
unique solutions $\bu$ and $\fq$ with
\begin{align*}
\bu &\in H^1_p((0, T), L_q(\Omega)^N) \cap L_p((0, T), H^2_q(\Omega)^N),
\\
\fq &\in L_p((0, T), H^1_q(\Omega) + \hat H^1_{q,0}(\Omega))
\end{align*}
possessing the estimate:
\begin{align*}
&\|\bu\|_{L_p((0, T), H^2_q(\Omega))}
+\|\pd_t\bu\|_{L_p((0, T), L_q(\Omega))} \\
&\quad \leq C(\|\bu_0\|_{B^{2(1-1/p)}_{q,p}(\Omega)}
+ \|\bff\|_{L_p(\BR, L_q(\Omega))} 
+ \|(g, \bh)\|_{H^{1/2}_p(\BR, L_q(\Omega))} \\
&\phantom{\quad \leq C(\|\bu_0\|_{B^{2(1-1/p)}_{q,p}(\Omega)}}\,\,
 + \|(g, \bh)\|_{L_p(\BR, H^1_q(\Omega))}
+ \|\pd_t\bg\|_{L_p(\BR, L_q(\Omega))})
\end{align*}
for some constant $C$. 
\end{thm}

Since $\pd_t(<t>^b\bu) = <t>^b\pd_t\bu + b<t>^{b-1}\bu$, if $\bu$ and 
$\fq$ satisfy Eq. \eqref{s-neweq:7.7}, then 
$<t>^b\bu$ and $<t>^b\fq$ satisfy the equations:
\begin{align*}
\pd_t(<t>^b\bu)  +\lambda_0(<t>^b\bu)
- J(T)^{-1}\DV(J\bA\mu\tilde\bD(<t>^b\bu)) &\\
+ \tilde\nabla(<t>^b\fq)
= 
<t>^b\tilde\bff + b<t>^{b-1}\bu&
\quad\text{in $\Omega^T$},\\
\quad\widetilde{\dv}<t>^b\bu = <t>^b\tilde g = \dv
(<t>^b\tilde\bg)&
\quad\text{in $\Omega^T$}, 
\\
\mu\tilde\bD(<t>^b\bu)\tilde\bn - <t>^b\fq\tilde\bn  
= <t>^b(Jd_\bn)^{-1}\tilde\bh&
\quad\text{on $\Gamma^T$}, \\
<t>^b\bu|_{t=0}  = \bu_0&
\quad \text{in $\Omega$}. 
\end{align*} 
Thus, repeated use of Theorem \ref{sthm:max} yields that
\begin{align}
&\|<t>^b\bu\|_{L_p((0, T), H^2_q(\Omega))}
+\|<t>^b\pd_t\bu\|_{L_p((0, T), L_q(\Omega))}
\nonumber  \\
& \leq C(\|\bu_0\|_{B^{2(1-1/p)}_{q,p}(\Omega)}
+ \|<t>^b\bff\|_{L_p(\BR, L_q(\Omega))} 
\nonumber \\
&+ \|(<t>^bg, <t>^b\bh)\|_{H^{1/2}_p(\BR, L_q(\Omega))} 
 + \|(<t>^bg, <t>^b\bh)\|_{L_p(\BR, H^1_q(\Omega))} 
\nonumber \\
&+ \|\pd_t(<t>^b\bg)\|_{L_p(\BR, L_q(\Omega))}),
\label{ineq:g.2}
\end{align}
provided that the right hand side is finite.  
Where, $g_b$, $\bg_b$ and $\bh_b$ are suitable extension of 
$<t>^bg$, $<t>^b\bg$, and $<t>^b\bh$ to the whole time 
interval $\BR$ such that 
$<t>g = g_b$, $<t>^b\bg = \bg_b$, $<t>^b\bh=\bh_b$ for 
$t \in (0, T)$ and $\dv \bg_b = g_b$ for $t \in \BR$. 

We next consider so called $L_p$-$L_q$ decay estimate of
semigroup associated with Eq. \eqref{neweq:7.7}. 
Thus, we formulate Eq. \eqref{neweq:7.7} 
in the semigroup setting.  For this purpose,  we have to
eliminate $\fq$.  We start with the weak Dirichlet problem:
\begin{equation}\label{pwd:1}
(\tilde\nabla u, J\tilde\nabla\varphi)_\Omega 
= (\bff, J\tilde\nabla\varphi)_\Omega
\quad \text{for any $\varphi \in \hat H^1_{q',0}(\Omega)$}.
\end{equation}
Here, 
$$\tilde \nabla\varphi= {}^\top(\sum_{k-1}^Na_{k1}\frac{\pd\varphi}{\pd x_k},
\ldots, \sum_{k-1}^Na_{kN}\frac{\pd\varphi}{\pd x_k})
= {}^\top\bA\nabla\varphi.
$$
Since $\tilde a_{ij}$ and $\tilde J$ vanish outside of $B_{2R}$, 
$\widetilde{\dv}\bu = \dv\bu$ and $\tilde \nabla\varphi= \nabla\varphi$
in $\BR^N\setminus B_{2R}$.  Thus, by Proposition \ref{prop:6.2}, we
have the following result.
\begin{prop}\label{prop:6.3} Let $1 < q \leq r$.  Then, for 
any $\bff \in L_q(\Omega)^N$ problem \eqref{pwd:1} admits a unique solution
$u \in \hat H^1_{q,0}(\Omega)$ possessing the estimate:
$$\|\nabla u\|_{L_q(\Omega)} \leq C\|\bff\|_{L_q(\Omega)}.$$
\end{prop} 

We next consider the following slightly perturned Dirichlet problem. 
Given $\bu \in H^2_q(\Omega)^N$, let $K(\bu)$ be a unique solution
of the weak Dirichlet problem:
\begin{equation}\label{pwd:2}
(\tilde\nabla K(\bu), J\tilde\nabla\varphi)_\Omega
= (\mu J^{-1}\DV(J\bA\tilde \bD(\bu))-\tilde\nabla\widetilde{\dv}\,\bu,
J\tilde\nabla\varphi)_\Omega
\end{equation}
for any $\varphi \in \hat H^1_{q',0}(\Omega)$, subject to
$$K(\bu) = <\tilde \bD(\bu)\tilde\bn, \tilde\bn> - \widetilde{\dv}\,\bu.
$$
This problem also can be treated by small perturbation of the weak Dirichlet
problem. 

We now consider the evolution equations: 
\begin{equation}\label{neweq:5*}
\left\{\begin{aligned}
\pd_t\bu -  J^{-1}\DV(J\bA\mu\tilde \bD(\bu))
+\tilde\nabla K(\bu)= 0&
&\quad&\text{in $\Omega\times(0, \infty)$}, \\
\mu\tilde\bD(\bu)\tilde\bn -  K(\bu)\tilde\bn  = 0&
&\quad&\text{on $\Gamma\times(0, \infty)$}, \\
\bu|_{t=0} = \bu_0&,
\end{aligned}\right. \end{equation}
which is corresponding to Eq. \eqref{neweq:7} with 
$\lambda_0=0$, $\bff=g=\bg=\bh=0$, and 
$\bu_0 \in \tilde J_q(\Omega)$. Where, we have set 
$$\tilde J_q(\Omega) = \{\bff \in L_q(\Omega) \mid
(\bff, J\tilde\nabla\varphi)_\Omega = 0
\quad\text{for any $\varphi \in \hat H^1_{q',0}(\Omega)$}
\}.
$$
Given $\bff \in \tilde J_q(\Omega)$, let $\bu \in H^2_q(\Omega)^N$ and 
$\fq \in H^1_q(\Omega) + \hat H^1_{q,0}(\Omega)$ be
solutions of the resolvent equations: 
\begin{equation}\label{spres:2}\left\{\begin{aligned}
\lambda\bu -  J^{-1}\DV(J\bA\mu\tilde \bD(\bu))
+\tilde\nabla K(\bu) = \bff,
\quad\widetilde{\dv}\,\bu=0&
\quad&\text{in $\Omega$}, \\
\mu\tilde\bD(\bu)\tilde\bn -  K(\bu)\tilde\bn   = 0&
\quad&\text{on $\Gamma$}.
\end{aligned}\right.
\end{equation}
Employing the same argument as in Subsec. \ref{subsec:3.3}, we have
$\fq = K(\bu)$, and so by Theorem \ref{spthm:1} we have
$\bu = \CA_s(\lambda)\bff$ and 
$$|\lambda|\|\bu\|_{L_q(\Omega)} + \|\bu\|_{H^2_q(\Omega)}
\leq r_b\|\bff\|_{L_q(\Omega)}$$
for any $\lambda \in \Sigma_{\epsilon_0, \lambda_0}$.
 Let 
\begin{align*}
\CD_{s,q}(\Omega) & = \{\bu \in \tilde J_q(\Omega) \cap H^2_q(\Omega)^N \mid 
\mu\tilde\bD(\bu)\tilde\bn -  K(\bu)\tilde\bn |_{\Gamma} = 0\}, \\
\CA_s\bu &= J^{-1}\DV(J\bA\mu\tilde \bD(\bu))
-\tilde\nabla K(\bu) \quad\text{for
$\bu \in\CD_{s,q}(\Omega)$}.
\end{align*}
Then, Eq. \eqref{neweq:5*} is written by
$$\dot \bu - \CA_s\bu = 0 \quad t>0, \quad \bu|_{t=0} = \bu_0
$$
for $\bu_0 \in \tilde J_q(\Omega)$. 
Notice that $\mu\tilde\bD(\bu)\tilde\bn -  K(\bu)\tilde\bn|_{\Gamma} = 0$
 for $\bu \in \CD_{s,q}(\Omega)$ is equivalent to
\begin{equation}\label{scompati:1}
\tilde \BPi_0[\mu\tilde\bD(\bu)\tilde\bn] =0 \quad\text{on $\Gamma$}.
\end{equation}
We then see that the operator $\CA_s$ generates a $C^0$
analytic semigroup $\{T_s(t)\}_{t\geq 0}$ on $\tilde J_q(\Omega)$.
Moreover, we have the following theorem.
\begin{thm}\label{thm:Lp-Lq}
Assume that $N \geq 3$ and that $\mu$ is a positive 
constant. Then,  there exists a $\sigma > 0$ such that
if the assumption \eqref{sassump:1} holds, then for any $q \in (1, r]$, 
there exists a $C^0$ analytic semigroup $\{T_s(t)\}_{t\geq 0}$ 
associated with Eq. \eqref{neweq:5*} such that
for any $\bu_0 \in \tilde J_q(\Omega)$,
 $\bu$ is a unique solution of Eq. \eqref{neweq:5*} with
\begin{align*}\bu &= T_s(t)\bu_0 \\ 
&\in C^0([0, \infty), \tilde J_q(\Omega))
\cap C^1((0, \infty), L_q(\Omega)^N)
\cap C^0((0, \infty), \CD_{s,q}(\Omega)).
\end{align*}

Moreover, for any $p \in [q, \infty)$ or $p=\infty$ and for any 
$\bff \in \tilde J_q(\Omega)$ and $t > 0$  
we have the following estimates:
\begin{equation}\label{lp-lq:1}\begin{split}
\|T_s(t)\bff\|_{L_p(\Omega)} & \leq C_{q,p}t^{-\frac12
\left(\frac1q-\frac1p\right)}\|\bff\|_{L_q(\Omega)}, \\
\|\nabla T_s(t)\bff\|_{L_p(\Omega)} & \leq C_{q,p}t^{-\frac12-\frac12
\left(\frac1q-\frac1p\right)}\|\bff\|_{L_q(\Omega)}.
\end{split}\end{equation}
\end{thm}
\begin{rem}
This theorem was proved by Shibata \cite{S7}. 
\end{rem}

We finally consider the following equations:
\begin{equation}\label{neweq:6*}\left\{
\begin{aligned}
\pd_t\bu -  J^{-1}\DV(J\bA\mu\tilde \bD(\bu))
+\tilde\nabla \fq  = \bff,
\quad\widetilde{\dv}\,\bu&=0 
&\quad&\text{in $\Omega^T$}, \\
\mu\tilde\bD(\bu)\tilde\bn -  \fq\tilde\bn  &= 0
&\quad&\text{on $\Gamma^T$},\\
\bu|_{t=0} &=0
&\quad&\text{in $\Omega$}.
\end{aligned}\right.
\end{equation}
Let $\psi \in \hat H^1_{q,0}(\Omega)$ be a solution of the weak
Dirichlet problem:
$$(\tilde\nabla\psi, J\tilde\nabla \varphi)_\Omega
= (\bff, J\tilde\nabla\varphi)_\Omega
\quad\text{for any $\varphi \in 
\hat H^1_{q',0}(\Omega)$}.
$$
Let $\bg = \bff - \tilde\nabla\psi$, and then 
$\bg \in \tilde J_q(\Omega)$ and 
$$\|\bg\|_{L_q(\Omega)} + \|\nabla \psi\|_{L_q(\Omega)}
\leq C\|\bff\|_{L_q(\Omega)}.
$$
Using this decomposition, we can rewrite Eq.
\eqref{neweq:6*} as 
$$\left\{
\begin{aligned}
\pd_t\bu -  J^{-1}\DV(J\bA\mu\tilde \bD(\bu))
+\tilde\nabla (\fq-\psi) = \bg,
\quad\widetilde{\dv}\,\bu&=0 
&\quad&\text{in $\Omega^T$}, \\
\mu\tilde\bD(\bu)\tilde\bn -  (\fq-\psi)\tilde\bn   &= 0
&\quad&\text{on $\Gamma^T$},\\
\bu|_{t=0} &=0
&\quad&\text{in $\Omega$},
\end{aligned}\right.
$$
where we have used the formula \eqref{sform:1} and
$\psi|_\Gamma=0$. Since $\bg \in \tilde J_q(\Omega)$ for 
any $t \in (0, T)$, we see that by Duhamel's principle, we have
\begin{equation}\label{sduhamel:1}
\bu = \int^t_0 T(t-s)\bg(s)\,ds = \int^t_0T(t-s)(\bff(s)-\nabla\psi(s))\,ds.
\end{equation}
Moreover, we have  $\fq = K(\bu) + \psi$. 
This is a solution formula of Eq. \eqref{neweq:6*}.

\subsection{Estimates for the nonlinear terms}\label{subsec:6.5}

Let $\tilde\bff(\bu)$, $\tilde g(\bu)$, $\tilde\bg(\bu)$ 
and $\tilde\bh(\bu)$ are functions defined
in Sect. \ref{subsec:6.3}. 
 In this subsection, we estimate these functions.
In the following we write
\begin{alignat*}2
\|<t>^\alpha \bw\|_{L_p((0, T), X)}
&= \Bigl\{\int^T_0(<t>^\alpha\|\bw(\cdot, t)\|_X)^p\,dt\Bigr\}^{1/p}
&\quad &1 \leq p < \infty, \\ 
\|<t>^\alpha \bw\|_{L_\infty((0, T), X)}
&=\underset{0 < t < T}{\rm esssup}\,<t>^\alpha\|\bw(\cdot, t)\|_X
&\quad &p = \infty.
\end{alignat*}
Let $\tilde\bff ={}^\top(\tilde f_1(\bu), \ldots, f_N(\bu))$ 
be the vector of functions given in
\eqref{non:1}. We first prove that
\begin{equation}\label{6.cond:2}\begin{split}
\|<t>^b\tilde\bff\|_{L_p((0, T), L_{q_1/2}(\Omega))}
&+ \|<t>^b\tilde\bff\|_{L_p((0, T), L_{q_2}(\Omega))} \\
&\leq C(\CI + [\bu]_T^2),
\end{split}\end{equation}
where $[\bu]_T$ is the norm defined in Theorem \ref{6.thm:main}
and  $\CI= \|\bu_0\|_{B^{2(1-1/p)}_{q_2, p}(\Omega)} 
+ \|\bu_0\|_{B^{2(1-1/p)}_{q_1/2,p}(\Omega)}$. 
Here and in the following, $C$ denotes generic constants independent
of $\CI$, $[\bu]_T$, $\delta$, and $T$. The value of $C$ may
change from line to line. 
Since we choose 
$\CI$ small enough eventually, we may assume that $0 < \CI \leq 1$. 
Especially, we use the estimates:
$$\CI^2 \leq \CI, \quad
 \CI [\bu]_T \leq \frac12(\CI^2 + [\bu]_T^2) \leq \CI + [\bu]_T^2$$
below.

Since 
\begin{align*}
&\int^\beta_\alpha \|\nabla(\kappa\bu(\cdot, s))\|_{L_\infty(\Omega)}\,ds \\
&\leq C(1+\alpha)^{-b+\frac{1}{p'}}\Bigl(\int^\beta_\alpha(<s>^b\|\bu(\cdot,
s)\|_{H^1_\infty(\Omega)})^p\,ds\Bigr)^{1/p},
\\
&\int^\beta_\alpha \|\nabla^2(\kappa\bu(\cdot, s))\|_{L_q(\Omega)}\,ds\\
&\leq C(1+\alpha)^{-b+\frac{N}{2q_2}+\frac{1}{p'}} 
\Bigl(\int^\beta_\alpha(<s>^{b-\frac{N}{2q_2}}
\|\bu(\cdot, s)\|_{H^2_{q_2}(\Omega)})^p\,ds\Bigr)^{1/p}
\end{align*}
for any $0 \leq \alpha < \beta \leq T$, where $q \in [1, q_2]$, 
we have 
\begin{equation}\label{6.est:6}\begin{split}
\int^\beta_\alpha \|\nabla(\kappa\bu(\cdot, s))\|_{L_\infty(\Omega)}\,ds
&\leq C[\bu]_T(1+\alpha)^{-b+\frac{1}{p'}}, \\
\int^\beta_\alpha \|\nabla^2(\kappa\bu(\cdot, s))\|_{L_q(\Omega)}\,ds
&\leq C[\bu]_T
\end{split}\end{equation}
for any $0 \leq \alpha < \beta \leq T$, 
where $q \in [1, q_2]$, because $b > \frac{N}{2q_2} + \frac{1}{p'}$
as follows from \eqref{6.cond:1}. By \eqref{real:7.1.1}, we have
\begin{multline}\label{6.est:6*}
\sup_{t \in (0, T)}<t>^{b-\frac{N}{2q_2}}
\|\bu(\cdot, t)\|_{B^{2(1-1/p)}_{q_2,p}(\Omega)}
\leq C(\|\bu_0\|_{B^{2(1-1/p)}_{q_2,p}(\Omega)} \\
+ \|<t>^{b-\frac{N}{2q_2}}\bu\|_{L_p((0, T), H^2_{q_2}(\Omega))}
+ \|<t>^{b-\frac{N}{2q_2}}\pd_t\bu\|_{L_p((0, T), L_{q_2}(\Omega))}\}.
\end{multline}
Since $2/p + N/q_2 < 1$, $B_{q_2, p}^{2(1-1/p)}(\Omega)$ is
continuously imbedded into $H^1_\infty(\Omega)$, and so 
by \eqref{6.est:6*} 
\begin{equation}\label{6.est:6**}
\|<t>^{b-\frac{N}{q_2}}\bu\|_{L_\infty((0, T), H^1_\infty(\Omega))}
\leq C(\CI + [\bu]_T).
\end{equation}
Applying  \eqref{assump:2}, \eqref{6.est:6} and \eqref{6.est:6*} 
to the formulas in \eqref{6.trans:0} and
\eqref{6.trans:4} and using the fact
 that $-b+ \frac{1}{p'} < -\frac{N}{2q_2}$ and 
$-b + \frac{N}{2q_2}\leq -\frac{N}{2q_2}$, which follows 
from \eqref{6.cond:1},  give  
\begin{align}
&\|(a_{ij}(t),  J(t), \ell_{ij}(t), \CA_{ij}(t), 
\CJ(t), \CL_{ij}(t))\|_{L_\infty(\Omega)} \leq C,  \nonumber \\
&\|(\CA_{ij}(t), \CJ(t), \CL_{ij}(t))\|_{L_\infty(\Omega)} \nonumber \\
&\leq C\int^T_t\|\nabla(\kappa\bu(\cdot, s))\|_{L_\infty(\Omega)}\,ds
\leq C[\bu]_T<t>^{-b+\frac{1}{p'}}
 \leq C[\bu]_T<t>^{-\frac{N}{2q_2}}, \nonumber\\
&\|\nabla( a_{ij}(t), J(t), \ell_{ij}(t), \CA_{ij}(t), \CJ(t), \CL_{ij}(t)) 
\|_{L_q(\Omega)} \nonumber\\
&\leq C\int^T_0\|\nabla^2(\kappa\bu(\cdot, s))\|_{L_q(\Omega)}
\leq C[\bu]_T,  \nonumber\\
&\|\pd_t(a_{ij}(t), J(t), \ell_{ij}(t), \CA_{ij}(t), \CJ(t), \CL_{ij}(t)
)\|_{L_\infty(\Omega)} 
\leq C\|\nabla(\kappa\bu(\cdot, t))\|_{L_\infty(\Omega)} \nonumber\\
&\leq C(\CI+[\bu]_T)<t>^{-b+\frac{N}{2q_2}}
\leq C(\CI+[\bu]_T)<t>^{-\frac{N}{2q_2}}
\label{6.est:4}\end{align}
for any $t \in (0, T]$, where $q \in [1, q_2]$. 
Moreover, we have 
\begin{equation}\label{6.est:4*}
(\tilde a_{ij}, \tilde J, \tilde \ell_{ij}, 
\CA_{ij},  \CJ, \CL_{ij})(x, t) = 0
\quad\text{for $x \not\in B_{2R}$ and $t \in [0, T]$}.
\end{equation}
By \eqref{6.est:4} and \eqref{6.est:4*},
$$\|a_{sm}(T)\CL_{is}(t)\pd_tu_i(t)\|_{L_q(\Omega)}
\leq C[\bu]_T<t>^{-b+\frac{1}{p'}}
\|\pd_tu_i(t)\|_{L_{q_2}(\Omega)}
$$
for any $q \in [1, q_2]$.  Since $\frac{1}{p'} < b - \frac{N}{2q_2}$
as follows from \eqref{6.cond:1}, we have
$$\|<t>^ba_{sm}(T)\CL_{is}\pd_tu_i\|_{L_p((0, T), L_q(\Omega))}
\leq C(\CI + [\bu]_T^2)
$$
for any $q \in [1, q_2]$. 

Next, by H\"older's inequality,   
\begin{align*}
&<t>^b\|\bu(\cdot, t)\cdot\nabla\bu(\cdot, t)\|_{L_{q_1/2}(\Omega)}
\\
&\quad
\leq <t>^{\frac{N}{2q_1}}\|\bu(\cdot, t)\|_{L_{q_1}(\Omega)}
<t>^{b-\frac{N}{2q_1}}\|\nabla\bu(\cdot, t)\|_{L_{q_1}(\Omega)}
\end{align*}
and so, by \eqref{6.est:4}, we have
$$\|<t>^ba_{sm}(T)\ell_{is}a_{kj}u_j\frac{\pd u_i}{\pd\xi_k}
\|_{L_p((0, T), L_{q_1/2}(\Omega))} \leq C[\bu]_T^2.
$$
Since  
\begin{align*}
&<t>^b\|\bu\cdot\nabla\bu(\cdot, t)\|_{L_{q_2}(\Omega)} \\
&\quad
\leq <t>^{\frac{N}{2q_2}}\|\bu(\cdot, t)\|_{L_\infty(\Omega)}
<t>^{b-\frac{N}{2q_2}}\|\nabla\bu(\cdot, t)\|_{L_{q_2}(\Omega)},
\end{align*}
 by \eqref{6.est:4} 
$$\|<t>^ba_{sm}(T)\ell_{is}a_{kj}u_j\frac{\pd u_i}{\pd\xi_k}
\|_{L_p((0, T), L_{q_2}(\Omega))} \leq C(\CI + [\bu]_T)[\bu]_T
\leq C(\CI + [\bu]_T^2). 
$$
Since
\begin{align*}
\frac{\pd}{\pd\xi_k}(\CD_{ij}(t)\nabla \bu)
&=\sum_{m=1}^N(\CA_{m j}(t)\frac{\pd^2 u_i}{\pd\xi_k\pd\xi_m}
+ \CA_{m i}(t)\frac{\pd^2 u_j}{\pd\xi_k\pd\xi_m})\\
&+
\sum_{m=1}^N((\frac{\pd}{\pd \xi_m}
 \CA_{m j}(t))\frac{\pd u_i}{\pd\xi_m}
+ (\frac{\pd}{\pd \xi_k}
\CA_{m i}(t))\frac{\pd u_j}{\pd\xi_m}),
\end{align*}
by \eqref{6.est:4} and \eqref{6.est:4*},
\begin{equation}\label{6.est:3.1}\begin{split}
&<t>^b\|\frac{\pd}{\pd\xi_k}(\CD_{ij}(\cdot)\nabla\bu)\|_{L_q(\Omega)}\\
&\leq C[\bu]_T\{<t>^{b-\frac{N}{2q_2}}
\|\nabla^2\bu(\cdot t)\|_{L_{q_2}(\Omega)}
+ <t>^b\|\nabla\bu(\cdot, t)\|_{L_\infty(\Omega)}\},
\end{split}\end{equation}
for any $q \in [1, q_2]$, and therefore  
$$\|<t>^ba_{sm}(T)\ell_{is}(T)a_{kj}(T)\frac{\pd}{\pd\xi_k}
(\CD_{ij}(\cdot)\nabla\bu)\|_{L_p((0, T), L_q(\Omega))}
\leq C[\bu]_T^2
$$
for any $q \in [1, q_2]$. Since 
\begin{align*}\frac{\pd}{\pd\xi_k}D_{ij, T}(\bu)
&=\sum_{m=1}^N(a_{m j}(T)\frac{\pd^2 u_i}{\pd\xi_k\pd\xi_m}
+ a_{m i}(T)\frac{\pd^2 u_j}{\pd\xi_k\pd\xi_m}) \\
&+
\sum_{m=1}^N((\frac{\pd}{\pd \xi_m}
 a_{m j}(T))\frac{\pd u_i}{\pd\xi_m}
+ (\frac{\pd}{\pd \xi_k}
a_{m i}(T))\frac{\pd u_j}{\pd\xi_m}),
\end{align*}
by \eqref{6.est:4} and \eqref{6.est:4*}, 
\begin{align*}
&<t>^b\|a_{sm}(T)\ell_{is}(T)\CA_{kj}(t)\frac{\pd}{\pd\xi_k}
D_{ij, T}(\bu)\|_{L_q(\Omega)} \\
&\quad
\leq C[\bu]_T\{<t>^{b-\frac{N}{2q_2}}\|\nabla^2\bu(\cdot t)\|_{L_{q_2}(\Omega)}
+ <t>^b\|\nabla\bu(\cdot, t)\|_{L_\infty(\Omega)}\},
\end{align*}
and so  
\begin{align*}
\|<t>^ba_{sm}(T)\ell_{is}(T)\CA_{kj}\frac{\pd}{\pd\xi_k}
D_{ij, T}(\bu)\|_{L_p((0, T), L_q(\Omega))} 
\leq C[\bu]_T^2
\end{align*}
for any $q \in [1, q_2]$. Analogously, we have
\begin{align*}
\|<t>^ba_{sm}(T)\CL_{is}a_{kj}\frac{\pd}{\pd\xi_k}
D_{ij, T}(\bu)\|_{L_p((0, T), L_q(\Omega))} 
\leq C[\bu]_T^2
\end{align*}
for any $q \in [1, q_2]$. Summing up, we have obtained \eqref{6.cond:2}.


We now consider  $\tilde g$  and $\tilde \bh={}^\top(\tilde h_1(\bu), 
\ldots, \tilde h_N(\bu))$, which have been defined in \eqref{non:1*}
and \eqref{non:1**}, respectively. 
To estimate the $H^{\frac12}_p$ norm,
we use the following lemma.
\begin{lem}\label{6.half:1}
Let $f \in H^1_\infty(\BR, L_\infty(\Omega))$ and 
$g \in H^{\frac12}_p(\BR, L_{q_2}(\Omega))$.  Assume that 
$f(x, t) = 0$ for $(x, t) \not\in B_R\times\BR$.  
 Then,  
\begin{equation}\label{6.half:2}
\|fg\|_{H^{\frac12}_p(\BR, L_q(\Omega))}
\leq C_q\|f\|_{H^1_\infty(\BR, L_\infty(\Omega))}
\|g\|_{H^{\frac12}_p(\BR, L_{q_2}(\Omega))}.
\end{equation}
for any $q \in [1, q_2]$ with some constant $C_q$ depending on
$q$ and $q_2$. 
\end{lem}
\pf
To prove the lemma, we use the fact that
\begin{equation}\label{6.half:3}
H^{\frac12}_p(\BR, L_q(\Omega)) = (L_p(\BR, L_q(\Omega)),
H^1_p(\BR, L_q(\Omega)))_{[\frac12]},
\end{equation}
where $(\cdot, \cdot)_{[1/2]}$ denotes a complex interpolation functor. 
Let $q \in [1, q_2]$.  Noting that $f(x, t) = 0$ for 
$(x, t) \not\in B_R\times\BR$, we have 
\begin{align*}
\|\pd_t(fg)\|_{L_q(\Omega)}
\leq \|\pd_t f\|_{L_\infty(\Omega)}\|g\|_{L_{q_2}(\Omega)} 
+ \|f\|_{L_\infty(\Omega)}\|\pd_tg\|_{L_{q_2}(\Omega)},
\end{align*}
and therefore
$$\|\pd_t(fg)\|_{L_p(\BR, L_q(\Omega))}
\leq C\|f\|_{H^1_\infty(\BR, L_\infty(\Omega))}
\|g\|_{H^1_p(\BR, L_{q_2}(\Omega))}
$$
for any $q \in [1, q_2]$. 
Moreover, we easily see that 
$$\|fg\|_{L_p(\BR, L_q(\Omega))}
\leq C\|f\|_{L_\infty(\BR, L_\infty(\Omega))}
\|g\|_{L_p(\BR, L_{q_2}(\Omega))}.
$$
Thus, by \eqref{6.half:3}, we have \eqref{6.half:2}, which completes the 
proof of Lemma \ref{6.half:1}.
\qed \vskip0.5pc

To use the maximal $L_p$-$L_q$ estimate, we have to extend $\tilde g$, 
$\tilde \bg$ and 
$\tilde \bh$ to $\BR$ with respect to time variable.  
For this purpose, we introduce an extension operator 
$\tilde e_T$.
Let $f$ be a function defined on $(0, T)$ such that $f|_{t=T} = 0$, and then
$\tilde e_T$ is an operator acting on $f$ defined by 
\begin{equation}\label{6.ref:1}[\tilde e_Tf](t)
= \begin{cases} 0 &\quad (t > T), \\
f(t) &\quad (0 < t < T), \\
f(-t) &\quad (-T < t < 0), \\
0 &\quad (t < -T).
\end{cases}
\end{equation}
\begin{lem}\label{6.half:4} Let $1 < p < \infty$, $1 \leq q \leq q_2$,
 and $0 \leq a \leq b$. Let $f \in H^1_\infty((0, T), L_\infty(\Omega))$ and 
$g \in H^1_p((0, T), L_{q_2}(\Omega)) \cap L_p((0, T), H^2_{q_2}(\Omega))$.
Assume that $f|_{t=T} = 0$ and $f = 0$ for $(x, t) \not\in B_R^T$. 
Let $<t> = (1+t^2)^{1/2}$.  Then, we have
\begin{equation}\label{half:5}\begin{split}
&\|\tilde e_T(<t>^af\nabla g)\|_{H^{\frac12}_p(\BR, L_q(\Omega))} 
\leq C\|<t>^{\frac{N}{2q_2}}f\|_{H^1_\infty((0, T), L_\infty(\Omega))}\\
&\times(\|<t>^{b-\frac{N}{2q_2}}g\|_{L_p((0, T), H^2_{q_2}(\Omega))}
+ \|<t>^{b-\frac{N}{2q_2}}\pd_tg\|_{L_p((0, T), L_{q_2}(\Omega))}
\\
&\phantom{\times(\|<t>^{b-\frac{N}{2q_2}}g\|_{L_p((0, T), H^2_{q_2}(\Omega))}}
\,\,\,
+ \|g|_{t=0}\|_{B^{2(1-1/p)}_{q_2, p}(\Omega)}).
\end{split}\end{equation}
\end{lem}
\pf
Let $f_0(t) = <t>^{a-b+\frac{N}{2q_2}}f(t)$
 and $g_0(t) = <t>^{b-\frac{N}{2q_2}}
g(t)$, and then $<t>^af\nabla g = f_0\nabla g_0$. Let $h$ be a function in 
$B^{2(1-1/p)}_{q_2, p}(\BR^N)$ such that $h = g|_{t=0}$ in $\Omega$ and 
$\|h\|_{B^{2(1-1/p)}_{q_2, p}(\Omega)} \leq 
C\|g|_{t=0}\|_{B^{2(1-1/p)}_{q_2, p}(\Omega)}$. 
Similarly to \eqref{tv:1}, we define $T_v(t)h$ by letting  
$T_v(t)h= e^{-(2-\Delta)}h$. 

Recall the operator
$e_T$ defined in \eqref{eT:1} and note that 
$g_0|_{t=0} = g|_{t=0} = T_v(t)h|_{t=0}$ in $\Omega$.  Let  
$$G(t) = e_T[g_0 - T_v(\cdot)h](t) + T_v(t)h
$$
for $t > 0$ and let 
$$[\iota g](t) = \begin{cases} G(t) &\quad (t > 0), \\
G(-t) & \quad (t < 0), \end{cases}
\quad 
[\iota f](t) = \begin{cases}0&\quad ( t > T), \\
 f_0(t) &\quad (0 < t < T), \\
f_0(-t) & \quad (-T < t < 0), \\
0 &\quad (t < -T).\end{cases}
$$
Since $G(t) = g_0(t)$ for $0 < t < T$, we have
\begin{align*}
\tilde e_T[<t>^af\nabla g](t)
&= \begin{cases} 0 &\, (t > T) \\
f_0(t)\nabla g_0(t) &\,(0 < t < T) \\
f_0(-t)\nabla g_0(-t) & \,(-T < t < 0) \\
0 &\, (t < -T) 
\end{cases} \\
&= 
\begin{cases} 0 &\, (t > T) \\
f_0(t)\nabla G(t) &\,(0 < t < T) \\
f_0(-t)\nabla G(-t) & \,(-T < t < 0) \\
0 &\,(t < -T) 
\end{cases} 
= [\iota f](t)\nabla[\iota g](t).
\end{align*}
By Lemma \ref{6.half:1}, 
\begin{align*}
\|\tilde e_T[<t>^af\nabla g]\|_{H^{\frac12}_p(\BR, L_q(\Omega))}
&= \|[\iota f]\nabla[\iota g]\|_{H^{\frac12}_p(\BR, L_q(\Omega))}
\\
&\leq C\|\iota f\|_{H^1_\infty(\BR, L_q(\Omega))}
\|\nabla(\iota g)\|_{H^{\frac12}_p(\BR, L_{q_2}(\Omega))}.
\end{align*}
Since $f_0(t)|_{t=T} = 0$, we have
$$\|\iota f\|_{H^1_\infty(\BR, L_\infty(\Omega))}
=2\|f_0\|_{H^1_\infty((0, T), L_\infty(\Omega))}
\leq \|<t>^{\frac{N}{2q_2}}f\|_{H^1_\infty((0, T), L_\infty(\Omega))},
$$
because $a - b \leq 0$. 

To estimate $\|\nabla(\iota g)\|_{H^{\frac12}_p(\BR, L_{q_2}(\Omega))}$,
we use Lemma \ref{lem:g.5.2}. And then, by Lemma \ref{lem:g.5.2}
and  \eqref{g.5.7}, we have
\begin{align*}
&\|\nabla(\iota g)\|_{H^{\frac12}_p(\BR, L_{q_2}(\Omega))}
\leq C(\|\iota g\|_{H^1_p(\BR, L_{q_2}(\Omega))}
+ \|\iota g\|_{L_p(\BR, H^2_{q_2}(\Omega))}) \\
& \leq C(\|G\|_{H^1_p((0, \infty), L_{q_2}(\Omega))}
+ \|G\|_{L_p((0, \infty), H^2_{q_2}(\Omega))})\\
&  \leq C(\|g_0-T_v(\cdot)h\|_{H^1_p((0, T), L_{q_2}(\Omega))}
+ \|g_0-T_v(\cdot)h\|_{L_p((0, T), H^2_{q_2}(\Omega))}\\
&\qquad + \|T_v(\cdot)h\|_{H^1_p((0, T), L_{q_2}(\Omega))}
+ \|T_v(\cdot)h\|_{L_p((0, \infty), H^2_{q_2}(\Omega))})\\
& \leq C(\|<t>^{b-\frac{N}{2q_2}}\pd_tg\|_{L_p((0, T), L_{q_2}(\Omega))}
+ \|<t>^{b-\frac{N}{2q_2}}g\|_{L_p((0, T), H^2_{q_2}(\Omega))}
\\
&\phantom{\leq C(\|<t>^{b-\frac{N}{2q_2}}\pd_tg
\|_{L_p((0, T), L_{q_2}(\Omega))}}\,\,
+ \|g|_{t=0}\|_{B^{2(1-1/p)}_{q_2, p}(\Omega)}).
\end{align*}
This completes the proof of Lemma \ref{6.half:4}.
\qed \vskip0.5pc

Recall the definitions of $\tilde g(\bu)$ 
and $\tilde h_m(\bu)$ given in \eqref{non:1*}
and \eqref{non:1**}.  By Lemma \ref{6.half:4} and \eqref{6.est:4} 
\begin{align}
&\|\tilde e_T[<t>^a\tilde g(\bu)]\|_{H^{\frac12}_p(\BR, L_q(\Omega))}
\nonumber \\
&\quad \leq \sum_{j,k=1}^N\|<t>^{\frac{N}{2q_2}}(J(T)\CA_{kj}(\cdot)
+ T_v(\cdot)a_{kj}(\cdot))\|_{H^1_\infty((0, T), L_\infty(\Omega))}
\nonumber \\
&\quad\quad \times
(\|<t>^{b-\frac{N}{2q_2}}\bu\|_{L_p((0, T), H^2_{q_2}(\Omega))}
+ \|<t>^{b-\frac{N}{2q_2}}\pd_t\bu\|_{L_p((0, T), L_{q_2}(\Omega))}
\nonumber \\
&\phantom{\quad\quad \times
(\|<t>^{b-\frac{N}{2q_2}}\bu\|_{L_p((0, T), H^2_{q_2}(\Omega))}}\,\,
+ \|\bu_0\|_{B^{2(1-1/p)}_{q_2, p}(\Omega)}) \nonumber \\
&\quad \leq C(\CI+[\bu]_T^2)
\label{6.cond:3}\end{align}
for any $a \in [0, b]$ and $q \in [1, q_2]$.
Analogously, we have
\begin{equation}\label{6.cond:4}
\|\tilde e_T[<t>^a \tilde\bh(\bu)]\|_{H^{\frac12}_p(\BR, L_q(\Omega))}
\leq C(\CI + [\bu]_T^2)
\end{equation}
for any $a \in [0, b]$ and $q \in [1, q_2]$. 

Next, by \eqref{6.est:4}, \eqref{6.est:4*} and \eqref{6.ref:1}, 
\begin{align}
&\|\tilde e_T[<t>^a\tilde g(\bu)]\|_{L_p(\BR, H^1_q(\Omega))} \nonumber \\
&\quad
\leq \sum_{j,k=1}^N\|<t>^{\frac{N}{2q_2}}(J(T)\CA_{kj}(\cdot)
+ \CJ(\cdot)a_{kj}(\cdot))\|_{L_\infty((0, T), L_\infty(\Omega))}
\nonumber \\
&\phantom{\quad
\leq \sum_{j,k=1}^N\|<t>^{\frac{N}{2q_2}}(J(T)\CA_{kj}(\cdot)}
\times
\|<t>^{b-\frac{N}{2q_2}}\bu\|_{L_p((0, T), H^2_{q_2}(\Omega))}\nonumber \\
&\quad +\sum_{j,k=1}^N\|\nabla(J(T)\CA_{kj}(\cdot)
+ \CJ(\cdot)a_{kj}(\cdot))\|_{L_\infty((0, T), L_q(\Omega))}
\nonumber \\
&\phantom{\quad +\sum_{j,k=1}^N\|\nabla(J(T)\CA_{kj}(\cdot)}
\times\|<t>^b\bu\|_{L_p((0, T), H^1_\infty(\Omega))} \nonumber \\
&\quad \leq C[\bu]_T^2
\label{6.cond:5}\end{align}
for any $a \in [0, b]$ and $q \in [1, q_2]$. Analogously, we have  
\begin{equation}\label{6.cond:6}
\|\tilde e_T[<t>^b\tilde\bh(\bu)]\|_{L_p((0, T), H^1_q(\Omega))}
\leq C[\bu]_T^2
\end{equation}
for any $a \in [0, b]$ and $q \in [1, q_2]$. 

We finally consider $\tilde \bg={}^\top(\tilde g_1(\bu), 
\ldots, \tilde g_N(\bu))$. 
Let $\tilde g_k(\bu)$ be functions given in \eqref{non:1*}.
Since 
\begin{align*}
\pd_t \tilde g_k(\bu) & = \sum_{j=1}^N(J(T)\pd_t \CA_{kj}(t)
+ (\pd_t\CJ(t))a_{kj}(t)
+ \CJ(t)(\pd_t a_{kj}(t))u_j \\
&+ \sum_{j=1}^N(J(T)\CA_{kj}(t) + \CJ(t)a_{kj}(t))\pd_t u_j
\end{align*}
and since $\|(J(T), a_{kj}(t), \CJ(t))\|_{L_\infty(\Omega)} \leq C$
as follows from \eqref{6.est:4}, by \eqref{6.est:4*} we have
\begin{align*}
&\|\tilde e_T[<t>^a\pd_tg_k(\bu)]\|_{L_p(\BR, L_q(\Omega))} \\
&\quad
\leq \sum_{j=1}^N(\|<t>^{\frac{N}{2q_2}}
\pd_t(\CA_{kj}, \CJ, a_{kj})\|_{L_\infty((0, T), 
L_\infty(\Omega))}\\
&\phantom{\quad
\leq \sum_{j=1}^N(\|<t>^{\frac{N}{2q_2}}
\pd_t(\CA_{kj},}\times
\|<t>^{b-\frac{N}{2q_2}}\bu\|_{L_p((0, T), L_{q_2}(\Omega))} \\
&\quad + \|<t>^\frac{N}{2q_2}(\CA_{kj}, \CJ)
\|_{L_\infty((0, T), L_\infty(\Omega))}
\|<t>^{b-\frac{N}{2q_2}}\pd_t\bu
\|_{L_p((0, T), L_{q_2}(\Omega))}),
\end{align*}
which, combined with \eqref{6.est:4}, leads to 
\begin{equation}\label{6.cond:7}
\|\tilde e_T[<t>^a\tilde \bg(\bu)]\|_{L_p(\BR, L_q(\Omega))}
\leq C(\CI + [\bu]_T^2)
\end{equation}
for any $a \in [0, b]$ and $q \in [1, q_2]$. 

\subsection{A Proof of Theorem \ref{6.thm:main}}\label{subsec:6.6}

The stragety of proving  Theorem \ref{6.thm:main} is
to prolong local in time solutions to any time interval,
which is  the same idea as that in Subsec. \ref{sec:5.6}.
Let $T$ be a positive number $> 2$.  
Let $\bu$ and $\fq$ be solutions of Eq. \eqref{navier:3},
which satisfy the regularity condition \eqref{6.reg.1} and the 
condition \eqref{assump:2}. Let $[\cdot]_T$ be the norm 
defined in \eqref{6.norm.1}.  To prove Theorem \ref{6.thm:main}
it suffices to prove that
\begin{equation}\label{6.prolong:1}
[\bu]_T \leq C(\CI + [\bu]_T^2)
\end{equation}
for some constant $C > 0$, where 
$$\CI = \|\bu_0\|_{B^{2(1-1/p)}_{q_2, p}(\Omega)} +
\|\bu_0\|_{B^{2(1-1/p)}_{q_1/2, p}(\Omega)}.
$$
If we show \eqref{6.prolong:1}, employing the same argument as that 
in Subsec. \ref{sec:5.6} and using Theorem \ref{thm:6.gloc*} 
concerning the almost global unique existence theorem, we can 
show that there exists a small constant $\epsilon > 0$ such that
if $\CI \leq \epsilon$ then $[\bu]_T \leq C\epsilon$ for some constant
$C > 0$ independent of $\epsilon$, and so we can 
prolong $\bu$ to any time interval
beyond $(0, T)$.  Thus,  we have Theorem \ref{6.thm:main}. 
In view of Theorem \ref{thm:6.gloc*}, there exists an 
$\epsilon_1 > 0$ such that 
if $\|\bu_0\|_{B^{2(1-1/p)}_{q_2, p}(\Omega)} \leq \epsilon_1$,
then $\bu$ and $\fq$ mentioned above surely exist. We assume that
$0 < \epsilon \leq \epsilon_1$.
Thus, our task below is to prove \eqref{6.prolong:1}.
In the following, we use the results stated in Subsec. \ref{subsec:6.4}
with $r=q_2$ and Subsec. \ref{subsec:6.5}.


As was seen in Subsec. \ref{subsec:6.2}, $\bu$ and $\fq$ 
satisfy Eq. \eqref{neweq:4}.  To estimate 
$\bu$, we divide $\bu$ and $\fq$ into two parts as $\bu = \bw + \bv$, 
and $\fq = \fr + \fp$, 
where $\bw$ and $\fr$ are  solutions of the equations: 
\begin{equation}\label{neweq:5}
\begin{cases*}
\pd_t\bw + \lambda_0\bw - J(T)^{-1}\DV(J\bA\mu\tilde\bD(\bw))
+ \tilde\nabla \fq &
= \tilde\bff(\bu)
\quad&\text{in $\Omega^T$},\\
\quad\widetilde{\dv}\bw = \tilde g(\bu) &
= \dv\tilde \bg(\bu)
\quad&\text{in $\Omega^T$}, 
\\
Jd_\bn(\mu\tilde\bD(\bw)\tilde\bn - \fq\tilde\bn)& = \tilde \bh(\bu)
\quad&\text{on $\Gamma^T$}, \\
\bw|_{t=0} & = \bu_0
\quad &\text{in $\Omega$},
\end{cases*}
\end{equation}
and $\bv$ and $\fp$ are  solutions of  the 
equations:
\begin{equation}\label{neweq:6}
\begin{cases*}
\pd_t\bv - J(T)^{-1}\DV(J\bA\mu\tilde\bD(\bv))
+ \tilde\nabla \fr &= -\lambda_0\bw
\quad&\text{in $\Omega^T$},\\
\quad\widetilde{\dv}\bv &= 0
\quad&\text{in $\Omega^T$}, 
\\
\mu\tilde\bD(\bv)\tilde\bn - \fr\tilde\bn & = 0
\quad&\text{on $\Gamma^T$}, \\
\bv|_{t=0} & = 0
\quad &\text{in $\Omega$}.
\end{cases*}
\end{equation}
Concerning the estimate of $\bw$,
applying \eqref{ineq:g.2} and using estimations:
\eqref{6.cond:2}, \eqref{6.cond:3}, 
\eqref{6.cond:4}, \eqref{6.cond:5}, 
\eqref{6.cond:6}, and \eqref{6.cond:7}, we have
\begin{equation}\label{6.7.1}\|<t>^b\pd_t\bw\|_{L_p((0, T), L_q(\Omega))}
+ \|<t>^b\bw\|_{L_p((0, T), H^2_q(\Omega))}
\leq C(\CI + [\bu]_T^2)
\end{equation}
for $q=q_1/2$, $q_1$, $q_2$. 
By Sobolev's inequality, we have
\begin{equation}\label{6.7.3}
\int^T_0(<s>^b\|\bw(\cdot, s)\|_{H^1_\infty(\Omega)})^p\,ds
\leq 
\int^T_0(<s>^b\|\bw(\cdot, s)\|_{H^2_{q_2}(\Omega)})^p\,ds,
\end{equation}
because $N < q_2 < \infty$. By \eqref{real:7.1.1}, 
\begin{equation}\label{6.7.2}\begin{aligned}
&\sup_{0<t<T}<s>^{\frac{N}{2q_1}}\|\bw(\cdot, t)\|_{L_{q_1}(\Omega)}
\leq C(\|\bu_0\|_{B^{2(1-1/q_1)}_{q_1,p}(\Omega)} \\
&\quad + \|<t>^{\frac{N}{2q_1}}\pd_t\bw\|_{L_p((0, T), L_{q_1}(\Omega))}
+ \|<t>^{\frac{N}{2q_1}}\bw\|_{L_p((0, T), H^2_{q_1}(\Omega))}).
\end{aligned}\end{equation}
Since $q_1/2 < q_1 < q_2$, we have 
$\|\bu_0\|_{B^{2(1-1/q_1)}_{q_1,p}(\Omega)} \leq \CI$.  Thus,
putting \eqref{6.7.1}, \eqref{6.7.3}, and \eqref{6.7.2} together yields that
\begin{equation}\label{6.5.2}
[\bw]_T \leq C(\CI + [\bu]_T^2).
\end{equation}

We next consider $\bv$. Let $\psi$ be a solution of the weak Dirichlet
problem
$$(\tilde\nabla\psi, J\tilde\nabla\varphi)_\Omega
= (-\lambda_0\bw, J\tilde\nabla\varphi)_\Omega
\quad\text{for any $\varphi \in \hat H^1_{q',0}(\Omega)$}.
$$
Let $\bP\bw = -\lambda_0\bw - \nabla\psi$. Since 
$$
\|\nabla\psi\|_{L_q(\Omega)} \leq C_q\|\bw\|_{L_q(\Omega)}
$$
for any $q \in (1, q_2]$,  we have  
\begin{equation}\label{proj:2}
\|\bP\bw\|_{L_q(\Omega)} \leq C_q\|\bw\|_{L_q(\Omega)}
\end{equation}
for any $q \in (1, q_2]$. 
Moreover, by \eqref{sduhamel:1}, we have
\begin{equation}\label{proj:3}
\bv(\cdot, t) = \int^t_0T(t-s)(\bP\bw)(\cdot, s)\,ds.
\end{equation}
Using the estimates \eqref{lp-lq:1} and \eqref{proj:2} yields that
\begin{equation}\label{6.5.4}\begin{split}
\|\nabla^j\bv(\cdot, t)\|_{L_r(\Omega)} 
&\leq C_{r, \tilde q_1}\int^{t-1}_0(t-s)^{-\frac{j}{2}
-\frac{N}{2}\left(\frac{1}{\tilde q_1}-\frac1r\right)}
\|\bw(\cdot, s)\|_{L_{\tilde q_1}(\Omega)}\,ds \\
&+ C_{r, \tilde q_2}\int^t_{t-1}(t-s)^{-\frac{j}{2}
-\frac{N}{2}\left(\frac{1}{\tilde q_2}-\frac1r\right)}
\|\bw(\cdot, s)\|_{L_{\tilde q_2}(\Omega)}\,ds
\end{split}\end{equation}
for $j=0,1$, for any $t > 1$ and for any indices $r$, $\tilde q_1$ and 
$\tilde q_2$ such that $1 < \tilde q_1, \tilde q_2
\leq r\leq\infty$ and $\tilde q_1$, $\tilde q_2 \leq q_2$, 
where $\nabla^0\bv = \bv$ and $\nabla^1\bv = \nabla\bv$.

Recall that $T > 2$.  In what follows,  we prove that 
\begin{align}
\Bigl(\int^T_2(<t>^b\|\bv(\cdot, t)\|_{H^1_\infty(\Omega)})^p\,dt\Bigr)^{1/p}
& \leq C(\CI + [\bu]_T^2), \label{6.5.6}\\
\sup_{2 \leq t \leq T}(<t>^{\frac{N}{2q_1}}
\|\bv(\cdot, t)\|_{L_{q_1}(\Omega)} & \leq C(\CI +[\bu]_T^2), \label{6.5.7}
\\
\Bigl(\int^T_2(<t>^{b-\frac{N}{2q_1}}
\|\bv(\cdot, t)\|_{H^1_{q_1}(\Omega)})^p\,dt\Bigr)^{1/p}
& \leq C(\CI + [\bu]_T^2), \label{6.5.7*} \\ 
\Bigl(\int^T_2(<t>^{b-\frac{N}{2q_2}}\|\bv(\cdot, t)
\|_{L_{q_2}(\Omega)})^p\,dt\Bigr)^{1/p}
& \leq C(\CI + [\bu]_T^2). \label{6.5.8}
\end{align}
By \eqref{6.5.4} with $r = \infty$, $\tilde q_1= q_1/2$ and 
$\tilde q_2= q_2$, 
\begin{align*}
\|\bv(\cdot, t)\|_{H^1_\infty(\Omega)}
& \leq C\int^t_0\|T(t-s)\bP\bw(\cdot, s)\|_{H^1_\infty(\Omega)}\,ds\\
&= C(I_{\infty}(t) + II_\infty(t) + III_\infty(t))
\end{align*}
with
\begin{align*}
I_\infty(t) &=
\int^{t/2}_0(t-s)^{-\frac{N}{q_1}}\|\bw(\cdot, s)\|_{L_{q_1/2}(\Omega)}\,ds, \\
II_\infty(t) &= 
\int^{t-1}_{t/2}
(t-s)^{-\frac{N}{q_1}}\|\bw(\cdot, s)\|_{L_{q_1/2}(\Omega)}\,ds, \\
III_\infty(t) &= \int^t_{t-1}(t-s)^{-\frac{N}{2q_2}-\frac12}
\|\bw(\cdot, s)\|_{L_{q_2}(\Omega)}\,ds.
\end{align*}
Since
\begin{align*}
&I_\infty(t) \\
&\leq (t/2)^{-\frac{N}{q_1}}
\Bigl(\int^{t/2}_0<s>^{-bp'}\,ds\Bigr)^{1/{p'}}
\Bigl(\int^{t/2}_0(<s>^b\|\bw(\cdot, s)\|_{L_{q_1/2}(\Omega)}
)^p\,ds\Bigr)^{1/p} \\
&\leq C(bp'-1)^{-1/{p'}}(\CI + [\bu]_T^2)t^{-\frac{N}{q_1}}
\end{align*}
as follows from the condition: $bp' > 1$ in \eqref{6.cond:1},  
by the condition: $(\frac{N}{q_1} - b)p > 1$ in \eqref{6.cond:1}, 
we have 
\begin{align*}
\int^T_2(<t>^bI_\infty(t))^p\,dt
&\leq C\int^T_2<t>^{-\left(\frac{N}{q_1}-b\right)p}\,dt
(\CI +[\bu]_T^2)^p \\
&\leq C((\frac{N}{q_1}-b)p-1)^{-1}(\CI + [\bu]_T^2)^p.
\end{align*}
By H\"older's inequality, 
\begin{align*}
&<t>^bII_\infty(t) \\
& \leq C\int^{t-1}_{t/2}
(t-s)^{-\frac{N}{q_1}}<s>^b\|\bw(\cdot, s)\|_{L_{q_1/2}(\Omega)}\,ds\\
& \leq C\Bigl(\int^{t-1}_{t/2}(t-s)^{-\frac{N}{q_1}}\,ds
\Bigr)^{1/{p'}}\\
&\qquad\times
\Bigl(\int^{t-1}_{t/2}(t-s)^{-\frac{N}{q_1}}
(<s>^b\|\bw(\cdot, s)\|_{L_{q_1/2}(\Omega)})^p\,ds\Bigr)^{1/p}\\
& \leq C\Bigl(\frac{N}{q_1}-1\Bigr)^{-1/{p'}}
\Bigl(\int^{t-1}_{t/2}(t-s)^{-\frac{N}{q_1}}(<s>^b\|\bw(\cdot, s)
\|_{L_{q_1/2}(\Omega)})^p\,ds\Bigr)^{1/p}
\end{align*}
because $N/q_1 = N/q_2 + 1 > 1$. By the Fubini-Tonelli theorem and 
\eqref{6.5.2}, 
\begin{align*}
&\int^T_2(<t>^bII_\infty(t))^p\,dt \\ 
&\leq C\Bigl(\frac{N}{q_1}-1\Bigr)^{-\frac{p}{p'}}
\int^T_2\,dt\int^{t-1}_{t/2}(t-s)^{-\frac{N}{q_1}}
(<s>^b\|\bw(\cdot, s)\|_{L_{q_1/2}(\Omega)})^p\,ds
\\
&\leq C\Bigl(\frac{N}{q_1}-1\Bigr)^{-\frac{p}{p'}}
\int^{T-1}_1(<s>^b\|\bw(\cdot, s)\|_{L_{q_1/2}(\Omega)}
)^p\,ds\int^{2s}_{s+1}(t-s)^{-\frac{N}{q_1}}\,dt\\
&\leq C\Bigl(\frac{N}{q_1}-1\Bigr)^{-p}(\CI +[\bu]_T^2)^p.
\end{align*}
Since $\frac{N}{2q_2} + \frac12 < 1$ as follows from 
$q_2 > N$,  by H\"older's inequality, 
\begin{align*}
&<t>^bIII_\infty(t)\\
& \leq C\int^t_{t-1}(t-s)^{-\frac{N}{2q_2}-\frac12}
<s>^b\|\bw(\cdot, s)\|_{L_{q_2}(\Omega)}\,ds \\
& \leq C\Bigl(\int^t_{t-1}(t-s)^{-\frac{N}{2q_2}-\frac12}\,ds
\Bigr)^{1/{p'}}\\
&\qquad\times\Bigl(\int^t_{t-1}(t-s)^{-\frac{N}{2q_2}-\frac12}
(<s>^b\|\bw(\cdot, s)\|_{L_{q_2}(\Omega)})^p\,ds
\Bigr)^{1/p}\\
&\leq C\Bigl(\frac{N}{2q_2}-\frac12\Bigr)^{-1/{p'}}
\Bigl(\int^t_{t-1}(t-s)^{-\frac{N}{2q_2}-\frac12}
(<s>^b\|\bw(\cdot, s)\|_{L_{q_2}(\Omega)})^p\,ds
\Bigr)^{1/p}.
\end{align*}
 By the Fubini-Tonelli theorem, we have 
\begin{align*}
&\int^T_2(<t>^bIII_\infty(t))^p\,dt \\
&\leq C\Bigl(1-\frac{N}{2q_2}\Bigr)^{-\frac{p}{p'}} \\
&\quad\times\int^T_2\,dt\int^t_{t-1}(t-s)^{-\frac{N}{2q_2}-\frac12}
(<s>^b\|\bw(\cdot, s)\|_{L_{q_2}(\Omega)})^p\,ds \\
&\leq C\Bigl(1-\frac{N}{2q_2}\Bigr)^{-\frac{p}{p'}} \\
&\quad\times\int^T_1(<s>^b\|\bw(\cdot, s)\|_{L_{q_2}(\Omega)}
)^p\,ds \int^{s+1}_s(t-s)^{-\frac{N}{2q_2}-\frac12}\,dt\\
& = C\Bigl(1-\frac{N}{2q_2}\Bigr)^{-p}(\CI + [\bu]_T^2)^p.
\end{align*}
Summing up, we have obtained \eqref{6.5.6}. 

Next, we prove \eqref{6.5.7}. By \eqref{6.5.4} with $r=q_1$, $\tilde q_1 = 
q_1/2$ and $\tilde q_2= q_1$,  
$$\|\bv(\cdot, t)\|_{L_{q_1}(\Omega)} \leq C(I_{q_1,\infty}(t)
+ II_{q_1, 1}(t) + III_{q_1, 1}(t))
$$
with
\begin{align*}
I_{q_1, 1}(t) & = \int^{t/2}_0(t-s)^{-\frac{N}{2q_1}}\|\bw(\cdot, s)
\|_{L_{q_1/2}(\Omega)}\,ds, \\
II_{q_1, 1}(t) & = \int^{t-1}_{t/2}(t-s)^{-\frac{N}{2q_1}}\|\bw(\cdot, s)
\|_{L_{q_1/2}(\Omega)}\,ds, \\
III_{q_1, 1}(t) & = \int^t_{t-1}\|\bw(\cdot, s)
\|_{L_{q_1}(\Omega)}\,ds.
\end{align*}
By \eqref{6.5.2} 
\begin{align*}
I_{q_1, 1}(t) 
& \leq (t/2)^{-\frac{N}{2q_1}}
\Bigl(\int^{t/2}_0<s>^{-bp'}\,ds
\Bigr)^{1/{p'}}\\
&\quad\times
\Bigl(\int^T_0(<s>^b\|\bw(\cdot, s)\|_{L_{q_1/2}(\Omega)})^p\,
ds\Bigr)^{1/p}\\
&\leq Ct^{-\frac{N}{2q_1}}(\CI +[\bu]_T^2).
\end{align*}
Analogously, by H\"older's inequality and \eqref{6.5.2},  
\allowdisplaybreaks{
\begin{align*}
II_{q_1, 1}(t) & \leq C\int^{t-1}_{t/2}(t-s)^{-\frac{N}{2q_1}}
<s>^{-b}<s>^b\|\bw(\cdot, s)\|_{L_{q_1/2}(\Omega)}\,ds \\
& \leq C<t>^{-b}\Bigl(\int^{t-1}_{t/2} (t-s)^{-\frac{Np'}{2q_1}}\,ds
\Bigr)^{1/{p'}}\\
&\quad\times
\Bigl(\int^T_0(<s>^b\|\bw(\cdot, s)\|_{L_{q_1/2}(\Omega)}
)^p\,ds\Bigr)^{1/p}\\
& = C\Bigl(1-\frac{Np'}{2q_1}\Bigr)^{1/{p'}}
<t>^{-b-\frac{N}{2q_1} + \frac{1}{p'}}(\CI +[\bu]_T^2)\\
& \leq C\Bigl(1-\frac{Np'}{2q_1}\Bigr)^{1/{p'}}
<t>^{-\frac{N}{2q_1}}(\CI +[\bu]_T^2),
\end{align*}
}
because $b > \frac{1}{p'}$. Finally, by \eqref{6.5.2}, 
\begin{align*}
III_{q_1, 1}(t) & \leq Ct^{-b}\int^t_{t-1}
<s>^b\|\bw(\cdot, s)\|_{L_{q_1/2}(\Omega)}\,ds \\
& \leq Ct^{-b}\Bigl(\int^t_{t-1}\,ds\Bigr)^{1/{p'}}
\Bigl(\int^T_0(<s>^b\|\bw(\cdot, s)\|_{L_{q_1/2}(\Omega)})^p
\,ds\Bigr)^{1/p} \\
&\leq Ct^{-b}(\CI + [\bu]_T^2).
\end{align*}
Summing up, we have obtained \eqref{6.5.7}. 

Next, we prove \eqref{6.5.7*}.  By \eqref{6.5.4} with 
$r=q_1$, $\bar q_1 = q_1/2$ and $\bar q_2=q_1$, 
$$\|\bv(\cdot, t)\|_{H^1_{q_1}(\Omega)}
\leq C(I_{q_1, 2}(t) + II_{q_1, 2}(t) + III_{q_1, 2}(t))
$$
with 
\begin{align*}
I_{q_1, 2}(t) &= \int^{t/2}_0 (t-s)^{-\frac{N}{2q_1}}
\|\bw(\cdot, s)\|_{L_{q_1/2}(\Omega)}\,ds,\\
II_{q_1, 2}(t) & = \int^{t-1}_{t/2} (t-s)^{-\frac{N}{2q_1}}
\|\bw(\cdot, s)\|_{L_{q_1/2}(\Omega)}\,ds, \\
III_{q_1, 2}(t) & = \int^{t}_{t-1} (t-s)^{-\frac12}
\|\bw(\cdot, s)\|_{L_{q_1}(\Omega)}\,ds.
\end{align*}
By \eqref{6.5.2}, 
\begin{align*}
I_{q_1, 2}(t)  
& \leq (t/2)^{-\frac{N}{2q_1}}
\Bigl(\int^{t/2}_0<s>^{-bp'}\,ds\Bigr)^{1/{p'}} \\
&\quad\times
\Bigl(\int^{t/2}_0(<s>^b\|\bw(\cdot, s)\|_{L_{q_1/2}(\Omega)})^p\,ds
\Bigr)^{1/p}\\
& \leq Ct^{-\frac{N}{2q_1}}(\CI + [\bu]_T^2),
\end{align*}
and so,  by the condition: $(\frac{N}{q_1} - b) p > 1$ in \eqref{6.cond:1}
$$\Bigl(\int^T_2(<t>^{b-\frac{N}{2q_1}}I_{q_1, 2}(t))^p\,dt
\Bigr)^{1/p} \leq C\Bigl((\frac{N}{q_1}-b)p - 1\Bigr)^{-1/p}
(\CI + [\bu]_T^2).
$$
By  H\"older's inequality, 
\begin{align*}
&<t>^{b-\frac{N}{2q_1}}II_{q_1, 2}(t) \\
&\leq C<t>^{-\frac{N}{2q_1}}\int^{t-1}_{t/2}(t-s)^{-\frac{N}{2q_1}}
<s>^b\|\bw(\cdot, s)\|_{L_{q_1/2}(\Omega)}\,ds \\
& \leq C<t>^{-\frac{N}{2q_1}}\Bigl(\int^{t-1}_{t/2}(t-s)^{-\frac{Np'}{2q_1}}
\,ds\Bigr)^{1/p'}\\
&\qquad\times
\Bigl(\int^T_0(<s>^b\|\bw(\cdot, s)\|_{L_{q_1/2}(\Omega)}
)^p\,ds\Bigr)^{1/p}\\
& \leq C(1+ t)^{-\left(\frac{N}{q_1} -\frac{1}{p'}\right)}(\CI + [\bu]_T^2).
\end{align*}
Since $(\frac{N}{q_1} - \frac{1}{p'})p  > 1$ as follows from
$\frac{N}{q_1} = 1 + \frac{N}{q_2} > 1 = \frac{1}{p} + \frac{1}{p'}$,  we have
$$\Bigl(\int^T_2(<t>^{b-\frac{N}{2q_1}}II_{q_1, 2}(t))^p\,dt
\Bigr)^{1/p}\leq C\Bigl((\frac{N}{q_1}-b)p-1\Bigr)^{-1/p}
(\CI + [\bu]_T^2). 
$$
Since 
\begin{align*}
<t>^{b-\frac{N}{2q_1}}III_{q_1, 2}(t) &\leq \int^t_{t-1}(t-s)^{-\frac12}
<s>^{b-\frac{N}{2q_1}}\|\bw(\cdot, s)\|_{L_{q_1}(\Omega)}\,ds
\\
&\leq \Bigl(\int^t_{t-1}(t-s)^{-\frac12}\,ds\Bigr)^{1/p'}
\\
&\quad\times
\Bigl(\int^t_{t-1}(t-s)^{-\frac12}
(<s>^b\|\bw(\cdot, s)\|_{L_{q_1}(\Omega)})^p\,ds
\Bigr)^{1/p},
\end{align*}
by the Fubini-Tonelli theorem, we have
\begin{align*}
&\int^T_2(<t>^{b-\frac{N}{2q_1}}III_{q_1, 2}(t))^p\,dt\\
&\quad\leq 2^{\frac{p}{p'}}\int^T_2\,dt\int^t_{t-1}
(t-s)^{-\frac12}(<s>^b\|\bw(\cdot, s)\|_{L_{q_1}(\Omega)})
^p\,ds \\
&\quad  \leq 2^{\frac{p}{p'}}\int^T_0
(<s>^b\|\bw(\cdot, s)\|_{L_{q_1}(\Omega)})
^p\,ds\\
&\qquad\times
\int^{s+1}_s(t-s)^{-\frac12}\,dt
 = 2^p\|<t>^b\bw\|_{L_p((0, T), L_{q_1}(\Omega))},
\end{align*}
which, combined with \eqref{6.7.1} with $q=q_1$, leads to 
$$
\Bigl(\int^T_2(<t>^{b-\frac{N}{2q_1}}III_{q_1, 2}(t))^p\,dt
\Bigr)^{1/p} \leq C(\CI +[\bu]_T^2).
$$
Summing up, we have obtained \eqref{6.5.7*}.

Finally, we prove \eqref{6.5.8}. By \eqref{6.5.4} with $r=q_2$, 
$\tilde q_1=q_1/2$ and $\tilde q_2=q_2$,
$$\|\bv(\cdot, t)\|_{L_{q_2}(\Omega)}
\leq C(I_{q_2}(t) + II_{q_2}(t) + III_{q_2}(t))
$$
with 
\begin{align*}
I_{q_2}(t) & = \int^{t/2}_0(t-s)^{-\frac{N}{2}
\left(\frac{2}{q_1} - \frac{1}{q_2}\right)}\|\bw(\cdot, s)
\|_{L_{q_1/2}(\Omega)}\,ds, \\
II_{q_2}(t) & = \int^{t-1}_{t/2}(t-s)^{-\frac{N}{2}
\left(\frac{2}{q_1} - \frac{1}{q_2}\right)}\|\bw(\cdot, s)
\|_{L_{q_1/2}(\Omega)}\,ds, \\
III_{q_2}(t) & = \int^{t}_{t-1}\|\bw(\cdot, s)
\|_{L_{q_2}(\Omega)}\,ds.
\end{align*}
By H\"older's inequality, 
\begin{align*}
I_{q_2}(t) &\leq (t/2)^{-\frac{N}{2}\left(\frac{2}{q_1}
-\frac{1}{q_2}\right)}\Bigl(\int^{t/2}_0<s>^{-bp'}\,ds\Bigr)^{1/{p'}}
\\
&\qquad\times 
\Bigl(\int^{t/2}_0(<s>^b\|\bw(\cdot, s)\|_{L_{q_1/2}(\Omega)})^p\,ds
\Bigr)^{1/p}\\
& \leq C<t>^{-\frac{N}{2}\left(\frac{2}{q_1}
-\frac{1}{q_2}\right)}(\CI + [\bu]_T^2)
\end{align*}
for $t \geq 2$. Since 
$$\frac{N}{2}\Bigl(\frac{2}{q_1}- \frac{1}{q_2}\Bigr) -
\Bigl(b-\frac{N}{2q_2}\Bigr)
= \frac{N}{q_1} - b,$$
by the condition: $(\frac{N}{q_1}-b)p > 1$
in \eqref{6.cond:1}, 
\begin{align*}
&\Bigl(\int^T_2(<t>^{b-\frac{N}{2q_2}}I_{q_2}(t))^p\,dt
\Bigr)^{1/p} \\
&\quad
\leq C\Bigl(\int^T_2t^{-\left(\frac{N}{q_1} - b\right)p}\,dt\Bigr)^{1/p}
(\CI +[\bu]_T^2)\\
&\quad
\leq C\Bigl((\frac{N}{q_1}-b)p-1\Bigr)^{-1/p}(\CI + [\bu]_T^2).
\end{align*}
Since 
$$\frac{N}{2}
\Bigl(\frac{2}{q_1} - \frac{1}{q_2}\Bigr)=\frac{N}{2}\Bigl(\frac{1}{q_2}
+\frac{2}{N}\Bigr) = \frac{N}{2q_2} + 1 > 1,$$
by H\"older's inequality
\begin{align*}
&<t>^{b-\frac{N}{2q_2}}
II_{q_2}(t)  \\
& \leq C\int^{t-1}_{t/2}
(t-s)^{-\left(\frac{N}{2q_2} + 1\right)}
<s>^{b-\frac{N}{2q_2}}\|\bw(\cdot, s)
\|_{L_{q_1/2}(\Omega)}\,ds \\
&\leq C\Bigl(\int^{t-1}_{t/2}
(t-s)^{-\left(\frac{N}{2q_2} + 1\right)}
\,ds\Bigr)^{1/{p'}}\\
&\qquad\times 
\Bigl(\int^{t-1}_{t/2}(t-s)^{-\left(\frac{N}{2q_2} + 1\right)}
(<s>^b\|\bw(\cdot, s)\|_{L_{q_1/2}(\Omega)})^p\,ds
\Bigr)^{1/p}\\
&\leq C\Bigl(\frac{N}{2q_2}\Bigr)^{-1/{p'}}
\Bigl(\int^{t-1}_{t/2}(t-s)^{-\left(\frac{N}{2q_2} + 1\right)}
(<s>^b\|\bw(\cdot, s)\|_{L_{q_1/2}(\Omega)})^p\,ds
\Bigr)^{1/p}, 
\end{align*}
and so,  by the Fubini-Tonelli theorem 
and \eqref{6.5.2}
\begin{align*}
&\int^T_2(<t>^{b-\frac{N}{2q_2}}
II_{q_2}(t))^p\,dt \\
&\leq C\Bigl(\frac{N}{2q_2}\Bigr)^{-p/{p'}}
\int^T_2\,dt
\int^{t-1}_{t/2}(t-s)^{-\left(\frac{N}{2q_2} + 1\right)}
(<s>^b\|\bw(\cdot, s)\|_{L_{q_1/2}(\Omega)})^p\,ds\\
&\leq C\Bigl(\frac{N}{2q_2}\Bigr)^{-p/{p'}}
\int^T_0(<s>^b\|\bw(\cdot, s)\|_{L_{q_1/2}(\Omega)})^p\,ds
\\
&\phantom{\leq C\Bigl(\frac{N}{2q_2}\Bigr)^{-p/{p'}}
\int^T_0(}\times
\int^{2s}_{s+1}(t-s)^{-\left(\frac{N}{2q_2} + 1\right)}\,dt
 \leq C\Bigl(\frac{N}{2q_2}\Bigr)^{-p}(\CI + [\bu]_T^2)^p.
\end{align*}
Analogously, by H\"older's inequality
\begin{align*}
<t>^{b-\frac{N}{2q_2}}III_{q_2}(t)
& \leq C\int^t_{t-1}<s>^{b-\frac{N}{2q_2}}
\|\bw(\cdot, s)\|_{L_{q_2}(\Omega)}\,ds, \\
& \leq C\Bigl(\int^t_{t-1}\,ds\Bigr)^{1/{p'}}
\Bigl(\int^t_{t-1}(<s>^b\|\bw(\cdot, s)\|_{L_{q_2}(\Omega)})^p
\,ds\Bigr)^{1/p} \\
& = C
\Bigl(\int^t_{t-1}(<s>^b\|\bw(\cdot, s)\|_{L_{q_2}(\Omega)})^p
\,ds\Bigr)^{1/p},
\end{align*}
and so,  by the Fubini-Tonelli theorem  and \eqref{6.5.2} 
\begin{align*}
&\int^T_2 (<t>^{b-\frac{N}{2q_2}}III_{q_2}(t))^p\,dt
\leq C\int^T_2\,dt
\int^t_{t-1}(<s>^b\|\bw(\cdot, s)\|_{L_{q_2}(\Omega)})^p
\,ds\\
&\quad \leq C\int^T_0(<s>^b\|\bw(\cdot, s)\|_{L_{q_2}(\Omega)})^p\,ds
\int^{s+1}_s \,dt
\leq C(\CI + [\bu]_T^2)^p.
\end{align*}
Summing up, we have obtained \eqref{6.5.8}.

We next consider the case where $t \in (0, 2)$. 
In view of Theorem \ref{spthm:1},
the usual maximal $L_p$-$L_q$ theorem holds for Eq. \eqref{neweq:6},
and so  we have 
\begin{equation}\label{6.5.9}\begin{split}
&\|\bv\|_{L_p((0, 2), H^2_q(\Omega))}
+ \|\pd_t\bv\|_{L_p((0, 2), L_q(\Omega))} \\
&\quad\leq C_q\|\lambda_0\bw\|_{L_p((0, 2), L_q(\Omega))}
\leq C(\CI + [\bu]_T^2)
\end{split}\end{equation}
for any $q \in [q_1/2, q_2]$.  Since $2(1-1/p) > N/q_2 + 1$ as follows
from $2/p + N/q_2 < 1$, by \eqref{real:7.1.7}, Sobolev's inequality, and 
\eqref{6.5.9}, we have
\begin{equation}\label{6.5.10.1}\begin{aligned}
&\Bigl(\int^2_0\|\bv(\cdot, t)\|_{H^1_\infty(\Omega)}\,dt\Bigr)^{1/p}
\leq  C\sup_{0 < t < 2}\|\bv(\cdot, t)\|_{B^{2(1-1/p)}_{q_2,p}(\Omega)}
\\
&\leq C(\|\bu_0\|_{B^{2(1-1/p)}_{q_2, p}(\Omega)}
 + \|\bv\|_{L_p((0, 2), H^2_{q_2}(\Omega))}
+ \|\pd_t\bv\|_{L_p((0, 2), L_{q_2}(\Omega))})\\
&\leq C(\CI + [\bu]_T^2). 
\end{aligned}\end{equation}
Moreover, by \eqref{real:7.1.1} and \eqref{6.5.9}
\begin{equation}\label{6.5.10.2} \begin{aligned}
&\sup_{0 < t < 2}\|\bv(\cdot, t)\|_{L_{q_1}(\Omega)}\\
&\leq C(\|\bu_0\|_{B^{2(1-1/p)}_{q_1, p}(\Omega)}  
+ \|\bv\|_{L_p((0, 2), H^2_{q_1}(\Omega))} 
 + \|\pd_t\bv\|_{L_p((0, 2), L_{q_1}(\Omega))})\\
&\leq C(\CI + [\bu]_T^2). 
\end{aligned}\end{equation}
Combining \eqref{6.5.6},
\eqref{6.5.7}, \eqref{6.5.7*}, \eqref{6.5.8}, \eqref{6.5.9},
\eqref{6.5.10.1} and \eqref{6.5.10.2},
we have
\begin{equation}\label{6.5.11}\begin{split}
&\|<t>^b\bv\|_{L_p((0, T), H^1_\infty(\Omega))}
+ \|<t>^{\frac{N}{2q_1}}\bv\|_{L_\infty((0, T), L_{q_1}(\Omega))} \\
&\quad
+ \|<t>^{b-\frac{N}{2q_1}}\bv\|_{L_p((0, T), H^1_{q_1}(\Omega))}
+ \|<t>^{b-\frac{N}{2q_2}}\bv\|_{L_p((0, T), L_{q_2}(\Omega))}
\\
&\qquad
\leq C(\CI + [\bu]_T^2).
\end{split}\end{equation}
From \eqref{neweq:6}, $\bv$ satisfies the equations:
\begin{alignat*}2
\pd_t\bv + \lambda_0\bv - J(T)^{-1}\DV(J\bA\mu\tilde\bD(\bv))
+ \tilde\nabla \fr &= -\lambda_0\bw + \lambda_0\bv
&\quad&\text{in $\Omega^T$},\\
\widetilde{\dv}\bv &= 0
&\quad&\text{in $\Omega^T$},\\
\mu\tilde\bD(\bv)\tilde\bn - \fr\tilde\bn  &= 0
&\quad&\text{on $\Gamma^T$}, \\
\bv|_{t=0} & = 0
&\quad &\text{in $\Omega$}, 
\end{alignat*}
and so by \eqref{ineq:g.2} we have  
\begin{align*}
&\|<t>^{b-\frac{N}{2q_2}}\bv
\|_{L_p((0, T), H^2_{q_2}(\Omega))}
+ \|<t>^{b-\frac{N}{2q_2}}\pd_t\bv
\|_{L_p((0, T), L_{q_2}(\Omega))}\\
&\quad
\leq C\|<t>^{b-\frac{N}{2q_2}}(\bv, \bw)\|_{L_p((0, T), L_{q_2}(\Omega))},
\end{align*}
which, combined with  \eqref{6.5.11}, leads to
\begin{equation}\label{6.5.12}\begin{split}
[\bv]_T 
\leq C(\CI + [\bu]_T^2).
\end{split}\end{equation}
Since $\bu = \bw + \bv$, by \eqref{6.5.2} and \eqref{6.5.12}, we see that 
$\bu$ satisfies the inequality \eqref{6.prolong:1}, which completes the proof
of Theorem \ref{6.thm:main}. 


\end{document}